\documentclass[10pt,a4paper]{amsart}
\usepackage{amssymb,amsxtra}
\usepackage[english,french]{babel}
\newtheorem{thm}{Theorem}[section]
\newtheorem{lem}[thm]{Lemma}

\theoremstyle{definition}
\newtheorem{defn}[thm]{Definition}

\theoremstyle{remark}

\theoremstyle{plain}

\theoremstyle{remark}

\newtheorem*{example}{Example}
\numberwithin{equation}{section}

\begin{document}
\title{Exponential sums and rank of triple persymmetric matrices over $\mathbb{F}_{2}$ }
\author{Jorgen~Cherly}
\address{D\'epartement de Math\'ematiques, Universit\'e de    
    Brest, 29238 Brest cedex~3, France}
\email{Jorgen.Cherly@univ-brest.fr}
\thanks{}
 \maketitle 
 \begin{abstract}
 
 Notre travail concerne  une g\' en\' eralisation des r\' esultats obtenus dans :
 {Exponential sums and rank of  double  persymmetric  matrices over  $\mathbf{F}_2 $  }\\
{arXiv : 0711.1937}. \\[0.1 cm]
 Soit $\mathbb{K}^{3} $  le  $\mathbb{K}$ - espace vectoriel de dimension 3  o\`{u} $\mathbb{K}$  d\'{e}note 
    le corps des s\'eries  de Laurent formelles  $ \mathbb{F}_{2}((T^{-1})). $ Nous calculons en particulier des sommes exponentielles 
    (dans $\mathbb{K}^{3}$) de la forme \\
     $  \sum_{\deg Y \leq  k-1}\sum_{\deg Z \leq s-1}E(tYZ)\sum_{\deg U \leq s+m-1}E(\eta YU)\sum_{\deg V \leq s+m+l-1}E(\xi YV) $
   o\`{u} $ (t,\eta,\xi ) $ est dans la boule  unit\'{e} de $\mathbb{K}^{3}.$\\ 
    Nous d\'{e}montrons qu'elles d\'{e}pendent uniquement
   du rang de matrices triples persym\'{e}triques  avec des entr\'{e}es dans $\mathbb{F}_{2},$
 c'est-\`{a}-dire des matrices de la forme    $\left[{A\over{B \over C}}\right] $ o\`u  A  est une matrice  $ s \times k $
    persym\'{e}trique , B une matrice  $ (s+m) \times k $  persym\'{e}trique  et  C une  matrice  $ (s+m +l) \times k $  persym\'{e}trique  (une matrice \;$ [\alpha _{i,j}]  $ est 
      persym\'{e}trique  si  $ \alpha _{i,j} = \alpha _{r,s} $ \; pour  \; i+j = r+s).
        En outre, nous \'{e}tablissons plusieurs 
      formules concernant des propri\'{e}t\'{e}s de rang de partitions de matrices triples persym\'{e}triques, ce qui nous 
      conduit \`{a} une formule r\'{e}currente du nombre   $    \Gamma_{i}^{\left[s\atop{ s+m \atop s+m+l} \right]\times k} $
  des matrices de rang i de la forme  $ \left[{A\over{B \over C}}\right] $
    Nous d\'{e}duisons de cette formule r\'{e}currente 
   que si  $  0\leq i\leq\inf(s-1, k-1),$ le nombre    $  \Gamma_{i}^{\left[s\atop{ s+m \atop s+m+l} \right]\times k}  $ d\'{e}pend uniquement de i.
    D'autre part, si   $ i\geq 2s+m+1, k\geq i,\;   \Gamma_{i}^{\left[s\atop{ s+m \atop s+m+l} \right]\times k}  $ peut \^{e}tre  calcul\'{e}
   \`{a} partir du nombre  $   \Gamma_{2s' +m' +1}^{\left[s'\atop{ s'+m' \atop s'+m'+l'} \right]\times k'}   $de matrices  de rang (2s'+m'+1) de la forme  
       $ \left[{A' \over{B' \over C'}}\right]$ o\`{u}  A' est une matrice   $ s' \times k' $ persym\'{e}trique,\; B' une matrice  $ (s'+m') \times k' $ persym\'{e}trique et C'  une matrice  $ (s'+m'+l') \times k' $
        o\`{u}  s', m',l' et k' d\' {e}pendent  de  i, s, m,l et k.
          La preuve de ce r\'{e}sultat est bas\'{e}e sur une formule
        du nombre de matrices de rang i de la forme     $\left[A\over b_{-}\right] $   o\`{u}  A est une matrice double  persym\'{e}trique  et  $ b_{-} $
 une matrice ligne avec entr\'{e}es dans $ \mathbb{F}_{2}.$ Nous montrons \'{e}galement que le nombre R de repr\'{e}sentations dans 
  $ \mathbb{F}_{2}[T] $ 
 des \'{e}quations polynomiales
   \[\left\{\begin{array}{cc}
YZ + Y_{1}Z_{1} + \ldots + Y_{q-1}Z_{q-1}= 0 \\
YU + Y_{1}U_{1} + \ldots + Y_{q-1}U_{q-1}= 0\\
YV + Y_{1}V_{1} + \ldots + Y_{q-1}V_{q-1}= 0
    \end{array}\right.\]\\ 
  associ\'{e}es aux sommes exponentielles \\
   $\sum_{degY\leq  k-1}\sum_{degZ\leq s-1}E(tYZ)\sum_{\deg U\leq s+m-1}E(\eta YU) \sum_{\deg V \leq s+m+l-1}E(\xi YV)  $
  est donn\'{e} par une int\'{e}grale sur la boule unit\'{e} de $\mathbb{K}^{3}$ et est une combinaison lin\'{e}aire de
  $ \Gamma_{i}^{\left[s\atop{ s+m \atop s+m+l} \right]\times k} $ pour $ i\geq 0. $ Nous pouvons alors calculer explicitement le nombre R.
  Notre article est,  pour des raisons de longueur,  limit\' e au cas $m \geqslant 0,\; l=0.$
 \end{abstract}
  \selectlanguage{english}
  \begin{abstract}
   
   Our work concerns a generalization of the results obtained in :
  {Exponential sums and rank of  double  persymmetric  matrices over  $\mathbf{F}_2 $  }\\
{arXiv : 0711.1937}. \\[0.1 cm]
Let $ \mathbb{K}^{3}  $ be the 3-dimensional vectorspace over $ \mathbb{K}$  where $ \mathbb{K}$ denotes the field of
 Laurent Series  $ \mathbb{F}_{2}((T^{-1})). $ We compute in particular   exponential sums, (in $\mathbb{K}^{3} $)  of the form \\
  $  \sum_{\deg Y \leq  k-1}\sum_{\deg Z \leq s-1}E(tYZ)\sum_{\deg U \leq s+m-1}E(\eta YU) \sum_{\deg V \leq s+m+l-1}E(\xi YV) $ where $ (t,\eta,\xi ) $ is in the unit interval of  $\mathbb{K}^{3}.$  We  show that they only depend  on the rank  of some associated  triple persymmetric matrices  
with entries in  $\mathbb{F}_{2}, $ that is matrices of the form  $\left[{A\over{B \over C}}\right] $ 
 where A is  a  $ s \times k $ persymmetric matrix,  B a  $ (s+m) \times k $ persymmetric matrix and C is a  $ (s+m+l) \times k $ persymmetric matrix
  (A matrix $\; [\alpha _{i,j}]  $ is  persymmetric if $ \alpha _{i,j} = \alpha _{r,s}  \; for  \; i+j = r+s ).$
   Besides, we establish several formulas concerning  rank properties of  partitions of  triple  persymmetric matrices, which leads to a  recurrent formula for the number $  \Gamma_{i}^{\left[s\atop{ s+m \atop s+m+l} \right]\times k} $  of rank i  matrices of the form 
   $\left[{A\over{B \over C}}\right] $  
    We deduce from the recurrent formula that  if  $  0\leq i\leq\inf(s-1, k-1) $ then  $  \Gamma_{i}^{\left[s\atop{ s+m \atop s+m+l} \right]\times k} $  depends only on  i.
  On the other hand,  if   $ i\geq 2s+m+1, k\geq i,\;   \Gamma_{i}^{\left[s\atop{ s+m \atop s+m+l} \right]\times k}  $ 
 can be computed from the number   $   \Gamma_{2s' +m' +1}^{\left[s'\atop{ s'+m' \atop s'+m'+l'} \right]\times k'}   $    of rank (2s'+m'+1) matrices of the form  $ \left[{A' \over{B' \over C'}}\right]$ 
 where A' is  a  $ s' \times k' $ persymmetric matrix,  B' a  $ (s'+m') \times k' $ persymmetric matrix 
 and C' a  $ (s'+m'+l') \times k' $ persymmetric matrix , where s', m',l' and k' depend on i, s, m,l and k.
  The proof of this result is based on a formula  of the number of rank i matrices of the form  $\left[A\over b_{-}\right] $ 
  where A is double persymmetric and $ b_{-} $ a one-row matrix with entries in  $ \mathbb{F}_{2}.$
   We also  prove that the number R of representations in  $ \mathbb{F}_{2}[T] $ 
 of the polynomial equations
   \[\left\{\begin{array}{cc}
YZ + Y_{1}Z_{1} + \ldots + Y_{q-1}Z_{q-1}= 0 \\
YU + Y_{1}U_{1} + \ldots + Y_{q-1}U_{q-1}= 0 \\
YV + Y_{1}V_{1} + \ldots + Y_{q-1}V_{q-1}= 0
    \end{array}\right.\]\\ 
  associated to the exponential sums \\
   $\sum_{degY\leq  k-1}\sum_{degZ\leq s-1}E(tYZ)\sum_{\deg U\leq s+m-1}E(\eta YU) \sum_{\deg U\leq s+m+l-1}E(\xi YV)  $
   is given by an integral over the unit interval of $\mathbb{K}^{3}$, and is a linear combination of the
    $ \Gamma_{i}^{\left[s\atop{ s+m \atop s+m+l} \right]\times k} $ pour $ i\geq 0. $
    We can  then compute  explicitly the number R.Our article is for reasons of length limited to the case $m \geqslant 0,\; l=0$
     \end{abstract}
 \tableofcontents 
 \allowdisplaybreaks  
  \newpage
 \textbf{CROQUIS DE LA PREUVE}\\[0.01 cm]
 \begin{enumerate}
\item Nous avons calcul\' e  explicitement le nombre de matrices de rang i de
   la forme  $\left[A\over B\right] $ o\`{u} B  est une matrice double persym\'{e}trique.
 \item Nous avons obtenu  une formule r\' ecurrente du nombre $   \Gamma_{i}^{\left[s\atop{ s+m\atop s+m+l} \right]\times k}$
 des matrices de rang i de la forme $\left[{A\over{B \over C}}\right] $ o\`u  A  est une matrice  $ s \times k $
    persym\'{e}trique , B une matrice  $ (s+m) \times k $  persym\'{e}trique  et  C une  matrice  $ (s+m +l) \times k $  persym\'{e}trique 
  \item  Nous ne disposons au d\' ebut que d'un  syst\`eme de 2 \' equations  en $3s+2m+l$ inconnues :
   \begin{align*}
  \sum_{ i = 0}^{\inf(3s+2m+l,k)}  \Gamma  _{i}^{\left[s\atop{ s+m\atop s+m+l} \right]\times k} & = 2^{3k+3s+2m+l-3}  \\
 &  \text{et} \\
  \sum_{i = 0}^{\inf(3s+2m+l,k)}  \Gamma  _{i}^{\left[s\atop{ s+m\atop s+m+l} \right]\times k}\cdot2^{-i} & =  2^{2k+3s+2m+l-3} + 2^{3k-3} - 2^{2k-3}. 
\end{align*}
 \item  On commencera par calculer  explicitement le nombre  $   \Gamma_{i}^{\left[s\atop{ s+m\atop s+m+l} \right]\times k}$ pour $i \leq s-1,\;
  i \leq k, \;m\geq 0,\; l\geq 0,$ en utilisant  successivement les  2 \' equations  dans (3).
  \item Nous  calculerons ensuite  explicitement le nombre  $   \Gamma_{i}^{\left[s\atop{ s+m\atop s+m+l} \right]\times k}$  pour 
$  s\leqslant i \leqslant\  inf(2s+m,k) $ en utilisant  successivement  les r\' esultats obtenus dans (4) et (3).
\item A partir de ces r\' esultats nous calculerons le nombre  $   \Gamma_{2s+m+1}^{\left[s\atop{ s+m\atop s+m+l} \right]\times k}$.
\item Guid\' e par des r\' esultats obtenus dans (2), j'\' emets les hypoth\`eses suivantes :
  \begin{align*}
\Gamma_{2s+1+m+j}^{\left[s\atop{ s+m\atop s+m+l} \right]\times k} & = 16^{j}\Gamma_{2s+1+m}^{\left[s\atop{ s+m\atop s+m+(l-j)} \right]\times (k-j)}   &  \text{si  }   0\leq j\leq l, \\
\Gamma_{2s+1+m+l+j}^{\left[s\atop{ s+m\atop s+m+l} \right]\times k} & = 16^{2j+l}\Gamma_{2s+1+m-j}^{\left[s\atop{ s+(m-j)\atop s+(m-j)} \right]\times (k-(2j+l))} & \text{si }  0\leq j\leq m, \\
\Gamma_{2s+1+2m+l+j}^{\left[s\atop{ s+m\atop s+m+l} \right]\times k} & = 16^{2m+l+3j}\Gamma_{2(s-j)+1}^{\left[s -j\atop{ s -j\atop s -j} \right]\times (k-(2m+l) -3j)} & \text{si }  0\leq j\leq s-1.
\end{align*}
\item  Une fois ces hypoth\`es v\' erifi\' ees, nous pouvons utiliser (6) et (7) pour obtenir le nombre  $   \Gamma_{i}^{\left[s\atop{ s+m\atop s+m+l} \right]\times k}$
pour $  2s+m+2 \leqslant i \leqslant \inf(3s+2m+l,k) .$ 
\end{enumerate}
\newpage

 \textbf{SKETCH OF THE PROOF}\\[0.01 cm]
 \begin{enumerate}
\item  We have explicitly calculated  the number of rank i matrices of the form  $\left[A\over B\right] $ where B is a double persymmetric matrix.
\item We have obtained a recurrent formula for the number $   \Gamma_{i}^{\left[s\atop{ s+m\atop s+m+l} \right]\times k}$ of rank i matrices 
of the form  $\left[{A\over{B \over C}}\right] $  where A is  a  $ s \times k $ persymmetric matrix,  B a  $ (s+m) \times k $ persymmetric matrix and C is a  $ (s+m+l) \times k $ persymmetric matrix.
\item In the beginning we only have a system of two equations in $3s+2m+l$ unknowns :
  \begin{align*}
  \sum_{ i = 0}^{\inf(3s+2m+l,k)}  \Gamma  _{i}^{\left[s\atop{ s+m\atop s+m+l} \right]\times k} & = 2^{3k+3s+2m+l-3}  \\
 &  \text{et} \\
  \sum_{i = 0}^{\inf(3s+2m+l,k)}  \Gamma  _{i}^{\left[s\atop{ s+m\atop s+m+l} \right]\times k}\cdot2^{-i} & =  2^{2k+3s+2m+l-3} + 2^{3k-3} - 2^{2k-3}. 
  \end{align*}
  \item We start by calculating  explicitly the number  $   \Gamma_{i}^{\left[s\atop{ s+m\atop s+m+l} \right]\times k}$ pour $i \leq s-1,\;
  i \leq k, \;m\geq 0,\; l\geq 0$ by using successively the two equations in (3).
  \item  We then  compute explicitly the number   $   \Gamma_{i}^{\left[s\atop{ s+m\atop s+m+l} \right]\times k}$  for
$  s\leqslant i \leqslant\  inf(2s+m,k) $ using successively the results obtained in (4) and (3).
  \item From these results we calculate the number   $   \Gamma_{2s+m+1}^{\left[s\atop{ s+m\atop s+m+l} \right]\times k}$ 
  of triple persymmetric $(3s+2m+l)\times k$ rank ($2s+m+1$) matrices.
  \item Leaded by the results obtained in (2), I suggest the following hypotheses : 
   \begin{align*}
\Gamma_{2s+1+m+j}^{\left[s\atop{ s+m\atop s+m+l} \right]\times k} & = 16^{j}\Gamma_{2s+1+m}^{\left[s\atop{ s+m\atop s+m+(l-j)} \right]\times (k-j)}   &  \text{si  }   0\leq j\leq l, \\
\Gamma_{2s+1+m+l+j}^{\left[s\atop{ s+m\atop s+m+l} \right]\times k} & = 16^{2j+l}\Gamma_{2s+1+m-j}^{\left[s\atop{ s+(m-j)\atop s+(m-j)} \right]\times (k-(2j+l))} & \text{si }  0\leq j\leq m, \\
\Gamma_{2s+1+2m+l+j}^{\left[s\atop{ s+m\atop s+m+l} \right]\times k} & = 16^{2m+l+3j}\Gamma_{2(s-j)+1}^{\left[s -j\atop{ s -j\atop s -j} \right]\times (k-(2m+l) -3j)} & \text{si }  0\leq j\leq s-1.
\end{align*}
\item Having proved  these hypotheses, we can use (6) and (7) to obtain the number  $   \Gamma_{i}^{\left[s\atop{ s+m\atop s+m+l} \right]\times k}$
for $  2s+m+2 \leqslant i \leqslant \inf(3s+2m+l,k) .$

\end{enumerate}

  \newpage

 \section{\textbf{Exponential sums and rank of matrices of the form  $\left[{A\over B}\right],$ where B is a double persymmetric matrix with entries in $\mathbb{F}_{2}$} }
 \label{sec 1}
  \subsection{Notation}
  \label{subsec 1}
  
   \subsubsection{Analyses on  $\mathbb{K}$}
  \label{subsubsec 1}
  We denote by $ \mathbb{F}_{2}\big(\big({\frac{1}{T}}\big) \big)
 = \mathbb{K} $ the completion
 of the field $\mathbb{F}_{2}(T), $ the field of  rational fonctions over the
 finite field\; $\mathbb{F}_{2}$, for the  infinity  valuation \;
 $ \mathfrak{v}=\mathfrak{v}_{\infty }$ \;defined by \;
 $ \mathfrak{v}\big(\frac{A}{B}\big) = degB -degA $ \;
 for each pair (A,B) of non-zero polynomials.
 Then every element non-zero t in  $\mathbb{F}_{2}\big(\big({\frac{1}{T}}\big) \big) $
 can be expanded in a unique way in a convergent Laurent series
                              $  t = \sum_{j= -\infty }^{-\mathfrak{v}(t)}t_{j}T^j
                                 \; where\; t_{j}\in \mathbb{F}_{2}. $\\
  We associate to the infinity valuation\; $\mathfrak{v}= \mathfrak{v}_{\infty }$
   the absolute value \; $\vert \cdot \vert_{\infty} $\; defined by \;
  \begin{equation*}
  \vert t \vert_{\infty} =  \vert t \vert = 2^{-\mathfrak{v}(t)}. \\
\end{equation*}
    We denote  E the  Character of the additive locally compact group
$  \mathbb{F}_{2}\big(\big({\frac{1}{T}}\big) \big) $ defined by \\
\begin{equation*}
 E\big( \sum_{j= -\infty }^{-\mathfrak{v}(t)}t_{j}T^j\big)= \begin{cases}
 1 & \text{if      }   t_{-1}= 0, \\
  -1 & \text{if      }   t_{-1}= 1.
    \end{cases}
\end{equation*}
  We denote $\mathbb{P}$ the valuation ideal in $ \mathbb{K},$ also denoted the unit interval of  $\mathbb{K},$ i.e.
  the open ball of radius 1 about 0 or, alternatively, the set of all Laurent series 
   $$ \sum_{i\geq 1}\alpha _{i}T^{-i}\quad (\alpha _{i}\in  \mathbb{F}_{2} ) $$ and, for every rational
    integer j,  we denote by $\mathbb{P}_{j} $
     the  ideal $\left\{t \in \mathbb{K}|\; \mathfrak{v}(t) > j \right\}. $
     The sets\; $ \mathbb{P}_{j}$\; are compact subgroups  of the additive
     locally compact group  $ \mathbb{K}. $\\
      All $ t \in \mathbb{F}_{2}\Big(\Big(\frac{1}{T}\Big)\Big) $ may be written in a unique way as
$ t = [t] + \left\{t\right\}, $   $  [t] \in \mathbb{F}_{2}[T] ,
 \; \left\{t\right\}\in \mathbb{P}  ( =\mathbb{P}_{0}). $\\
 We denote by dt the Haar measure on  $ \mathbb{K} $\; chosen so that \\
  $$ \int_{\mathbb{P}}dt = 1. $$\\
   \begin{defn}
\label{defn 1.1}We introduce the following definitions in  $ \mathbb{K}: $
\begin{itemize}
\item Let s, m and k denote rational integers such that $ s\geq 2,\; m\geq 0 \; and \; k\geq 1. $\\

\item  A matrix\;$ D = [\alpha _{i,j}]  $ is said to be persymmetric if $ \alpha _{i,j} = \alpha _{r,s} $ \; whenever \; i+j = r+s. \vspace{0.5 cm}
\item Set $ t = \sum_{i\geq 1}\alpha _{i}T^{-i}\in \mathbb{P}, $ 
we denote by $ D_{  s \times k}(t) $   the following $ s \times k $ persymmetric matrix 
 $$   \left ( \begin{array} {cccccc}
\alpha _{1} & \alpha _{2} & \alpha _{3} &  \ldots & \alpha _{k-1}  &  \alpha _{k} \\
\alpha _{2 } & \alpha _{3} & \alpha _{4}&  \ldots  &  \alpha _{k} &  \alpha _{k+1} \\
\vdots & \vdots & \vdots    &  \vdots & \vdots  &  \vdots \\
\alpha _{s-1} & \alpha _{s} & \alpha _{s +1} & \ldots  &  \alpha _{s+k-3} &  \alpha _{s+k-2}  \\
\alpha _{s} & \alpha _{s+1} & \alpha _{s +2} & \ldots  &  \alpha _{s+k-2} &  \alpha _{s+k-1} 
 \end{array}  \right). $$ \vspace{0.5 cm}
 \item  Set $ \eta  = \sum_{i\geq 1}\beta _{i}T^{-i}\in \mathbb{P}, $ 
we denote by $ D_{ (s +m )\times k}(\eta ) $   the following $ (s+m)\times k $ persymmetric matrix 
 $$   \left ( \begin{array} {cccccc}
\beta  _{1} & \beta  _{2} & \beta  _{3} & \ldots  &  \beta_{k-1} &  \beta _{k}  \\
\beta  _{2} & \beta  _{3} & \beta  _{4} & \ldots  &  \beta_{k} &  \beta _{k+1}  \\
\vdots & \vdots & \vdots    &  \vdots & \vdots  &  \vdots \\
\beta  _{m+1} & \beta  _{m+2} & \beta  _{m+3} & \ldots  &  \beta_{k+m-1} &  \beta _{k+m}  \\
\vdots & \vdots & \vdots    &  \vdots & \vdots  &  \vdots \\
\beta  _{s+m-1} & \beta  _{s+m} & \beta  _{s+m+1} & \ldots  &  \beta_{s+m+k-3} &  \beta _{s+m+k-2}  \\
\beta  _{s+m} & \beta  _{s+m+1} & \beta  _{s+m+2} & \ldots  &  \beta_{s+m+k-2} &  \beta _{s+m+k-1} 
\end{array}  \right). $$ \vspace{0.5 cm}
\end{itemize}
\end{defn}

   \subsubsection{Analyses on the two-dimensional  $\mathbb{K}$-vectorspace}
  \label{subsubsec 2}
    Let $\mathbb{K}\times \mathbb{K}  = \mathbb{K}^2 $
 be the 2-dimensional vector space 
  over $ \mathbb{K}. $ Let $ (t,\eta )\in \mathbb{K}^2  $  and  $ \vert (t,\eta )\vert =
  sup \left\{\vert t \vert ,\vert \eta \vert \right\} =  2^{-inf(\mathfrak{v}(t),\mathfrak{v}(\eta ))}. $\\
 It is easy to see that $ (t,\eta )\longrightarrow   \vert (t,\eta )\vert $
  is an ultrametric valuation on $ \mathbb{K}^2, $ 
  that is,   $ (t,\eta )\longrightarrow   \vert (t,\eta )\vert $ is a norm and
  $  \vert((t,\eta ) + (t',\eta '))\vert \leq max\left\{\vert(t,\eta ) \vert ,\vert (t',\eta ')\vert \right\}. $ \\
   We denote by  $ d (t,\eta )=dtd\eta $ the  Haar measure on   $ \mathbb{K}^2 $
  chosen so that the measure on the unit interval of $ \mathbb{K}^2  $ is equal to one, thus 
 $$ \iint_{\mathbb{P}\times \mathbb{P}} d (t,\eta )=
   \int_{\mathbb{P}}dt\int_{\mathbb{P}}d\eta = 1\cdot 1 =1. $$\\
   
 Let  $(t,\eta ) = \big( \sum_{i= -\infty }^{-\mathfrak{v}(t)}t_{i}T^i ,
\sum_{i = -\infty }^{-\mathfrak{v}(\eta )}\eta _{i}T^i   \big) \in  \mathbb{K}^2, $ 
 we denote $\chi $  the  Character on  $(\mathbb{K}^2, +) $ defined by \\
\begin{equation*}
 \chi\big( \sum_{i= -\infty }^{-\mathfrak{v}(t)}t_{i}T^i ,
\sum_{ i = -\infty }^{-\mathfrak{v}(\eta )}\eta _{i}T^i   \big) =
  E\big( \sum_{i = -\infty }^{-\mathfrak{v}(t)}t_{i}T^i\big)\cdot
  E\big( \sum_{i = -\infty }^{-\mathfrak{v}(\eta )}\eta _{i}T^i\big) =
\begin{cases}
 1 & \text{if      }   t_{-1} +  \eta _{-1}= 0, \\
  -1 & \text{if      }   t_{-1} + \eta _{-1}=1.
    \end{cases}
\end{equation*}

\begin{defn}
\label{defn 1.2}We introduce the following definitions in the two-dimensional  $ \mathbb{K}$ -vectorspace.
\begin{itemize}
\item Let k, s and m denote rational integers such that $ k\geq 1,\; s\geq 2 \; and \; m\geq 0. $\\

\item We denote by $\mathbb{P}/\mathbb{P}_{i}\times \mathbb{P}/\mathbb{P}_{j} $
a complete set of coset representatives of $\mathbb{P}_{i}\times\mathbb{P}_{j} $
in  $\mathbb{P}\times\mathbb{P},$ for instance $\mathbb{P}/\mathbb{P}_{s+k-1}\times \mathbb{P}/\mathbb{P}_{s+m +k-1} $
denotes  a complete set of coset representatives of $\mathbb{P}_{s+k-1}\times\mathbb{P}_{s+m +k-1} $ in  $\mathbb{P}\times\mathbb{P}.$\\

\item $Set \;(t,\eta )= (\sum_{i\geq 1}\alpha _{i}T^{-i},\sum_{i\geq 1}\beta  _{i}T^{-i})
\in \mathbb{P}\times\mathbb{P}. $\\
  We denote by  $  D^{\left[\stackrel{s}{s+m}\right] \times k }(t,\eta  ) $ any  $(2s+m)\times k $   matrix, 
such that  after a rearrangement of the rows, if necessary,  we can  obtain  the following double persymmetric matrix 
$ \left[{D_{s  \times k}(t)\over D_{(s+m )\times k}(\eta )}\right] $ \\

 $$   \left ( \begin{array} {cccccc}
\alpha _{1} & \alpha _{2} & \alpha _{3} &  \ldots & \alpha _{k-1}  &  \alpha _{k} \\
\alpha _{2 } & \alpha _{3} & \alpha _{4}&  \ldots  &  \alpha _{k} &  \alpha _{k+1} \\
\vdots & \vdots & \vdots    &  \vdots & \vdots  &  \vdots \\
\alpha _{s-1} & \alpha _{s} & \alpha _{s +1} & \ldots  &  \alpha _{s+k-3} &  \alpha _{s+k-2}  \\
\alpha _{s} & \alpha _{s+1} & \alpha _{s +2} & \ldots  &  \alpha _{s+k-2} &  \alpha _{s+k-1}  \\
\hline \\
\beta  _{1} & \beta  _{2} & \beta  _{3} & \ldots  &  \beta_{k-1} &  \beta _{k}  \\
\beta  _{2} & \beta  _{3} & \beta  _{4} & \ldots  &  \beta_{k} &  \beta _{k+1}  \\
\vdots & \vdots & \vdots    &  \vdots & \vdots  &  \vdots \\
\beta  _{m+1} & \beta  _{m+2} & \beta  _{m+3} & \ldots  &  \beta_{k+m-1} &  \beta _{k+m}  \\
\vdots & \vdots & \vdots    &  \vdots & \vdots  &  \vdots \\
\beta  _{s+m-1} & \beta  _{s+m} & \beta  _{s+m+1} & \ldots  &  \beta_{s+m+k-3} &  \beta _{s+m+k-2}  \\
\beta  _{s+m} & \beta  _{s+m+1} & \beta  _{s+m+2} & \ldots  &  \beta_{s+m+k-2} &  \beta _{s+m+k-1} 
\end{array}  \right). $$ \\
\item $Set \;(t,\eta )= (\sum_{i\geq 1}\alpha _{i}T^{-i},\sum_{i\geq 1}\beta  _{i}T^{-i})
\in \mathbb{P}\times\mathbb{P}. $\\
Let $  \Gamma _{i}^{\Big[\substack{s \\ s+m }\Big] \times k} $ denote the number of 
double persymmetric $ (2s+m)\times k $ rank i  matrices  of the form $ \left[{D_{s  \times k}(t)\over D_{(s+m )\times k}(\eta )}\right], $that is \vspace{0.1 cm}\\
 $  \Gamma _{i}^{\Big[\substack{s \\ s+m }\Big] \times k}  =   
Card  \left\{(t,\eta )\in \mathbb{P}/\mathbb{P}_{k+s-1}\times \mathbb{P}/\mathbb{P}_{k+s+m-1}
\mid   r(D^{\big[\stackrel{s}{s+m}\big] \times k }(t,\eta )) = i
 \right\}. $ 
\end{itemize}
\end{defn}

  \subsubsection{Analyses on the three-dimensional  $\mathbb{K}$-vectorspace}
  \label{subsubsec 3}
    Let $\mathbb{K}\times \mathbb{K}\times \mathbb{K}  = \mathbb{K}^3 $
 be the 3-dimensional vector space over $ \mathbb{K}. $  
 Let $ (t,\eta,\xi  )\in \mathbb{K}^3  $  and  $ \vert (t,\eta ,\xi )\vert =
  sup \left\{\vert t \vert ,\vert \eta \vert ,\vert \xi  \vert\right\} =  2^{-inf(\mathfrak{v}(t),\mathfrak{v}(\eta ),\mathfrak{v}(\xi  ))}. $\\
 It is easy to see that $ (t,\eta ,\xi )\longrightarrow   \vert (t,\eta,\xi  )\vert $
  is an ultrametric valuation on $ \mathbb{K}^3. $ \\
  That is, $ (t,\eta ,\xi )\longrightarrow   \vert (t,\eta,\xi  )\vert $ is a norm and
  $  \vert((t,\eta ,\xi ) + (t',\eta ',\xi' ))\vert \leq max\left\{\vert(t,\eta ,\xi ) \vert ,\vert (t',\eta ',\xi ')\vert \right\}. $ \\
   We denote by  $ d (t,\eta ,\xi )=dtd\eta d\xi $ the  Haar measure on   $ \mathbb{K}^3 $
  chosen so that the measure on the unit interval of  $ \mathbb{K}^3 $ is equal to one, that is 
 $$ \iiint_{\mathbb{P}\times \mathbb{P}\times \mathbb{P}} d (t,\eta ,\xi ) =
   \int_{\mathbb{P}}dt\int_{\mathbb{P}}d\eta\int_{\mathbb{P}}d\xi  = 1\cdot 1\cdot 1 =1. $$
 Let  $(t,\eta ,\xi ) = \big( \sum_{i= -\infty }^{-\mathfrak{v}(t)}t_{i}T^i ,
\sum_{ i = -\infty }^{-\mathfrak{v}(\eta )}\eta _{i}T^i ,
 \sum_{i = -\infty }^{-\mathfrak{v}(\xi )}\xi  _{i}T^i \big) \in  \mathbb{K}^3, $ 
 we denote $\psi  $  the  Character on  $(\mathbb{K}^3, +) $ defined by \\
\begin{align*}
 \psi\big( \sum_{i = -\infty }^{-\mathfrak{v}(t)}t_{i}T^i ,
\sum_{i = -\infty }^{-\mathfrak{v}(\eta )}\eta _{i}T^i ,
 \sum_{i = -\infty }^{-\mathfrak{v}(\xi )}\xi  _{i}T^i   \big) & =
  E\big( \sum_{i = -\infty }^{-\mathfrak{v}(t)}t_{i}T^i\big)\cdot
  E\big( \sum_{i = -\infty }^{-\mathfrak{v}(\eta )}\eta _{i}T^i\big)\cdot
   E\big( \sum_{i = -\infty }^{-\mathfrak{v}(\xi )}\xi _{i}T^i\big) \\
    & =
   \begin{cases}
 1 & \text{if      }   t_{-1} +  \eta _{-1} +\xi _{-1}     = 0 \\
  -1 & \text{if      }   t_{-1} + \eta _{-1}  +\xi _{-1}     =1
    \end{cases}
\end{align*}
 \begin{defn}
\label{defn 1.3}We introduce the following definitions in the three-dimensional  $ \mathbb{K} $-vectorspace.
 \begin{itemize}
\item We denote by $ \mathbb{P}/\mathbb{P}_{k+s-1} \times  \mathbb{P}/\mathbb{P}_{k+s+m-1} \times  \mathbb{P}/\mathbb{P}_{k}  $
  a complete set of coset representatives of \\
  $\mathbb{P}_{k+s-1} \times \mathbb{P}_{k+s+m-1}\times \mathbb{P}_{k}\; in \; 
   \mathbb{P}\times \mathbb{P}\times \mathbb{P}. $
\item  $  Set \;(t,\eta ,\xi ) =
 \Big(\sum_{i\geq 1}\alpha _{i}T^{-i},\sum_{i\geq  1}\beta _{i}T^{-i},\sum_{i\geq  1}\gamma  _{i}T^{-i}\Big) \in
 \mathbb{P}\times \mathbb{P}\times \mathbb{P}$  
  We denote by  $ D^{\left[s\atop{ s+m\atop 1} \right]\times k}(t,\eta ,\xi )  $ any  $(2s+m +1)\times k $   matrix, 
such that  after a rearrangement of the rows, if necessary,  we can  obtain  the following triple persymmetric matrix 
$ \left[{D_{s  \times k}(t)\over{ D_{(s+m )\times k}(\eta )\over D_{1\times k}(\xi )}}\right] $ 
\begin{equation*}
  \left ( \begin{array} {cccccc}
\alpha _{1} & \alpha _{2} & \alpha _{3} &  \ldots & \alpha _{k-1}  &  \alpha _{k} \\
\alpha _{2 } & \alpha _{3} & \alpha _{4}&  \ldots  &  \alpha _{k} &  \alpha _{k+1} \\
\vdots & \vdots & \vdots    &  \vdots & \vdots  &  \vdots \\
\alpha _{s-1} & \alpha _{s} & \alpha _{s +1} & \ldots  &  \alpha _{s+k-3} &  \alpha _{s+k-2}  \\
\alpha _{s} & \alpha _{s+1} & \alpha _{s +2} & \ldots  &  \alpha _{s+k-2} &  \alpha _{s+k-1}  \\
\hline \\
\beta  _{1} & \beta  _{2} & \beta  _{3} & \ldots  &  \beta_{k-1} &  \beta _{k}  \\
\beta  _{2} & \beta  _{3} & \beta  _{4} & \ldots  &  \beta_{k} &  \beta _{k+1}  \\
\vdots & \vdots & \vdots    &  \vdots & \vdots  &  \vdots \\
\beta  _{m+1} & \beta  _{m+2} & \beta  _{m+3} & \ldots  &  \beta_{k+m-1} &  \beta _{k+m}  \\
\vdots & \vdots & \vdots    &  \vdots & \vdots  &  \vdots \\
\beta  _{s+m-1} & \beta  _{s+m} & \beta  _{s+m+1} & \ldots  &  \beta_{s+m+k-3} &  \beta _{s+m+k-2}  \\
\beta  _{s+m} & \beta  _{s+m+1} & \beta  _{s+m+2} & \ldots  &  \beta_{s+m+k-2} &  \beta _{s+m+k-1} \\
\hline \\
\gamma  _{1} & \gamma  _{2} & \gamma  _{3} & \ldots  &  \gamma _{k-1} &  \gamma  _{k}  
\end{array}  \right).
\end{equation*}  
 \item   Set $$   \;(\xi_{1},\xi _{2},\ldots,\xi _{n}, t,\eta ) 
  =  \Big(\sum_{i\geq 1}\gamma  _{1i}T^{-i},\sum_{i\geq 1}\gamma  _{2i}T^{-i},\ldots,\sum_{i\geq 1}\gamma  _{ni}T^{-i},\sum_{i\geq  1}\alpha  _{i}T^{-i},\sum_{i\geq  1}\beta  _{i}T^{-i}\Big) \in
  \mathbb{P}^{n+2}$$ \\
 We denote by $  D^{\left[(n)\atop{ 1+m\atop 1+m+l} \right]\times k}(\xi, t,\eta  )  $
    the following $(n+2m+l+2)\times k $ matrix, where the n first rows form a  $n\times k $ matrix
  over the finite field  $ \mathbb{F}_{2} $, the following (1+m)  rows  form  a $ (1+m)\times k $ persymmetric matrix 
   and the $(1+m+l) $ last rows form a $ (1+m+l)\times k $ persymmetric matrix.  
  $$  \left ( \begin{array} {cccccc}
\gamma _{11} & \gamma   _{12} & \gamma  _{13} & \ldots  & \gamma  _{1,k-1} &  \gamma  _{1,k}  \\
 \gamma  _{21} & \gamma  _{22} & \gamma  _{23} & \ldots  & \gamma  _{2,k-1} &   \gamma _{2,k}\\
\vdots & \vdots & \vdots   &  \ldots  & \vdots  &  \vdots \\
\vdots & \vdots & \vdots    &  \ldots & \vdots  &  \vdots \\
 \gamma _{n1} & \gamma  _{n2} &  \gamma  _{n3} & \ldots  & \gamma  _{n,k-1} &  \gamma  _{n,k}\\
 \hline
 \alpha _{1} & \alpha _{2} & \alpha _{3} &  \ldots & \alpha _{k-1}  &  \alpha _{k} \\
\alpha _{2 } & \alpha _{3} & \alpha _{4}&  \ldots  &  \alpha _{k} &  \alpha _{k+1} \\
\vdots & \vdots & \vdots   &  \ldots  & \vdots  &  \vdots \\
\vdots & \vdots & \vdots    &  \ldots & \vdots  &  \vdots \\
\alpha _{1+m} & \alpha _{2+m} & \alpha _{3+m} & \ldots  &  \alpha _{k+m-1} &  \alpha _{k+m}  \\
\hline \\
\beta  _{1} & \beta  _{2} & \beta  _{3} & \ldots  &  \beta_{k-1} &  \beta _{k}  \\
\beta  _{2} & \beta  _{3} & \beta  _{4} & \ldots  &  \beta_{k} &  \beta _{k+1}  \\
\vdots & \vdots & \vdots    &  \vdots & \vdots  &  \vdots \\
\beta  _{1+m+l} & \beta  _{2+m+l} & \beta  _{3+m+l} & \ldots  &  \beta_{k+m+l-1} &  \beta _{k+m+l}  
  \end{array}  \right) $$
\item Consider the following partition of the  matrix $  D^{\left[s\atop{ s+m\atop 1} \right]\times k}(t,\eta ,\xi )   $\\
    $$   \left ( \begin{array} {cccccc}
\alpha _{1} & \alpha _{2} & \alpha _{3} &  \ldots & \alpha _{k-1}  &  \alpha _{k} \\
\alpha _{2 } & \alpha _{3} & \alpha _{4}&  \ldots  &  \alpha _{k} &  \alpha _{k+1} \\
\vdots & \vdots & \vdots    &  \vdots & \vdots  &  \vdots \\
\alpha _{s-1} & \alpha _{s} & \alpha _{s +1} & \ldots  &  \alpha _{s+k-3} &  \alpha _{s+k-2}  \\
\alpha _{s} & \alpha _{s+1} & \alpha _{s +2} & \ldots  &  \alpha _{s+k-2} &  \alpha _{s+k-1}  \\
\beta  _{1} & \beta  _{2} & \beta  _{3} & \ldots  &  \beta_{k-1} &  \beta _{k}  \\
\beta  _{2} & \beta  _{3} & \beta  _{4} & \ldots  &  \beta_{k} &  \beta _{k+1}  \\
\vdots & \vdots & \vdots    &  \vdots & \vdots  &  \vdots \\
\beta  _{m+1} & \beta  _{m+2} & \beta  _{m+3} & \ldots  &  \beta_{k+m-1} &  \beta _{k+m}  \\
\vdots & \vdots & \vdots    &  \vdots & \vdots  &  \vdots \\
\beta  _{s+m-1} & \beta  _{s+m} & \beta  _{s+m+1} & \ldots  &  \beta_{s+m+k-3} &  \beta _{s+m+k-2}  \\
\beta  _{s+m} & \beta  _{s+m+1} & \beta  _{s+m+2} & \ldots  &  \beta_{s+m+k-2} &  \beta _{s+m+k-1} \\
\hline \\
\gamma  _{1} & \gamma  _{2} & \gamma  _{3} & \ldots  &  \gamma _{k-1} &  \gamma  _{k}  
\end{array}  \right). $$ 
 \item  We define $$  \sigma _{i,i}^{\left[\stackrel{s}{\stackrel{s+m }{\overline{\gamma _{1 -}}}}\right] \times k } $$
 to be the cardinality of the following set 
$$ \left\{(t,\eta,\xi  )\in \mathbb{P}/\mathbb{P}_{k+s-1}\times \mathbb{P}/\mathbb{P}_{k+s+m-1}\times \mathbb{P}/\mathbb{P}_{k}
\mid  r( D^{\left[s\atop s+m\right]\times k}(t,\eta  ) ) 
 = r( D^{\left[s\atop{ s+m\atop 1} \right]\times k}(t,\eta ,\xi ) ) = i \right\}.$$
 \item   We define $$  \sigma _{i-1,i}^{\left[\stackrel{s}{\stackrel{s+m }{\overline{\gamma _{1 -}}}}\right] \times k } $$
 to be the cardinality of the following set 
$$ \left\{(t,\eta,\xi  )\in \mathbb{P}/\mathbb{P}_{k+s-1}\times \mathbb{P}/\mathbb{P}_{k+s+m-1}\times \mathbb{P}/\mathbb{P}_{k}
\mid  r( D^{\left[s\atop s+m\right]\times k}(t,\eta  ) ) = i-1,\quad
  r( D^{\left[s\atop{ s+m\atop 1} \right]\times k}(t,\eta ,\xi ) ) = i \right\}.$$
\item $Set \;(t,\eta,\xi  )= (\sum_{i\geq 1}\alpha _{i}T^{-i},\sum_{i\geq 1}\beta  _{i}T^{-i},\sum_{i\geq 1}\gamma   _{i}T^{-i})
\in \mathbb{P}\times\mathbb{P}\times\mathbb{P} $
Let $   \Gamma_{i}^{\left[s\atop{ s+m\atop 1} \right]\times k} $ denote the number of 
triple persymmetric $ (2s+m+1)\times k $ rank i  matrices  of the form
 $ \left[{D_{s  \times k}(t)\over{ D_{(s+m )\times k}(\eta )\over D_{1\times k}(\xi )}}\right],$ that is 
 $$   \Gamma_{i}^{\left[s\atop{ s+m\atop 1} \right]\times k} =   
Card  \left\{(t,\eta,\xi  )\in \mathbb{P}/\mathbb{P}_{k+s-1}\times \mathbb{P}/\mathbb{P}_{k+s+m-1}\times \mathbb{P}/\mathbb{P}_{k}
\mid   r( D^{\left[s\atop{ s+m\atop 1} \right]\times k}(t,\eta ,\xi ) ) = i \right\}.$$ 
\item $Set \;(t,\eta,\xi  )= (\sum_{i\geq 1}\alpha _{i}T^{-i},\sum_{i\geq 1}\beta  _{i}T^{-i},\sum_{i\geq 1}\gamma   _{i}T^{-i})
\in \mathbb{P}\times\mathbb{P}\times\mathbb{P} $\\
Let $   \Gamma_{i}^{\left[1\atop{ 1+m\atop 1+m+l} \right]\times k} $ denote the number of 
triple persymmetric $ (2m+l+3)\times k $ rank i  matrices  of the form
 $ \left[{D_{1  \times k}(t)\over{ D_{(1+m )\times k}(\eta )\over D_{(1+m+l)\times k}(\xi )}}\right],$ that is \\
 $$   \Gamma_{i}^{\left[1\atop{ 1+m\atop 1+m+l} \right]\times k} =   
Card  \left\{(t,\eta,\xi  )\in \mathbb{P}/\mathbb{P}_{k}\times \mathbb{P}/\mathbb{P}_{k+m}\times \mathbb{P}/\mathbb{P}_{k+m+l}
\mid   r( D^{\left[1\atop{ 1+m\atop 1+m+l} \right]\times k}(t,\eta ,\xi ) ) = i \right\}.$$
\end{itemize}
\end{defn}

    \subsection{Introduction}
  \label{subsec 2}
  We generalize the results obtained in [1] \\[0.01 cm]
  
     \subsection{Statement of main results}
  \label{subsec 3}

  \begin{thm}
\label{thm 1.4}
For all  $1\leq i\leq inf(k,2s+m+1) $ we have \\
\begin{equation*}
  \Gamma_{i}^{\left[s\atop{ s+m\atop 1} \right]\times k}= ( 2^{k} -  2^{i-1})\cdot \Gamma _{i-1}^{\left[\stackrel{s}{s+m}\right]\times k}
 +  2^{i}\cdot\Gamma _{i}^{\left[\stackrel{s}{s+m}\right]\times k}
\end{equation*}
 \end{thm}
 
 \begin{thm}
\label{thm 1.5} 
  We have the following reduction formulas \\
  \begin{align*}
\Gamma_{m+3+j}^{\left[1\atop{ 1+m\atop 1+m+l} \right]\times k} & = 2^{4j}\Gamma_{m+3}^{\left[1\atop{ 1+m\atop 1+m+(l-j)} \right]\times (k-j)}   &  \text{if  }   0\leq j\leq l   \\
\Gamma_{m+l+3+j}^{\left[1\atop{ 1+m\atop 1+m+l} \right]\times k} & = 2^{8j+4l}\Gamma_{m+3-j}^{\left[1\atop{ 1+(m-j)\atop 1+(m-j)} \right]\times (k-(2j+l))} & \text{if  }  0\leq j\leq m 
\end{align*}
 \end{thm}
\begin{thm}
\label{thm 1.6}
We have for $ m\geq 2,\;  l\geq 3 $ and  $k \geq 2m+l $\\
\small
$ \Gamma_{i}^{\left[1\atop{ 1+m\atop 1+m+l} \right]\times k}  $  is equal to \\
\[ =
\begin{cases}
1 & \text{if    } i = 0\\
2^{k} +17 & \text{if    } i = 1 \\
 (21\cdot2^{3i-7} - 3\cdot2^{2i-5})\cdot2^{k} + 315\cdot2^{4i-8} -21\cdot2^{3i-6}& \text{if    }2\leq  i\leq m \\
 21\cdot2^{k+3m-4} + 13\cdot2^{k+2m-3} +315\cdot2^{4m-4} -85\cdot2^{3m-3} & \text{if    } i = m+1 \\
 2^{2k+m} + 21\cdot2^{k+3m-1} + 9\cdot2^{k+2m-1} +315\cdot2^{4m} - 157\cdot2^{3m} &\text{if    } i = m+2 \\
 3\cdot2^{2k-m+2i-6} + 21\cdot2^{k+3i-7} + 21\cdot2^{k-m+3i-8} +315\cdot2^{4i-8} - 315\cdot2^{4i-m-9}& \text{if    } m+3\leq i\leq m+l \\
  3\cdot2^{2k+m+2l-4} + 21\cdot2^{k+3m +3l-4} + 149\cdot2^{k+2m+3l-5} +315\cdot2^{4m+4l-4} - 827\cdot2^{3m+4l-5} &\text{if    } i = m+l+1 \\
 11\cdot2^{2k+m+2l-2} + 21\cdot2^{k+3m +3l-1} + 261\cdot2^{k+2m+3l-2} +315\cdot2^{4m+4l} - 1627\cdot2^{3m+4l-1}& \text{if    } i = m+l+2 \\ 
21\cdot2^{2k-2m-l+3i-9} + 21\cdot2^{k+3i-7} + 525\cdot2^{k-2m-l+4i-11} +315\cdot2^{4i-8} - 3255\cdot2^{5i-2m-l-12}& \text{if    } m+l+3\leq i\leq 2m+l+1 \\ 
53\cdot2^{2k+4m+2l-3} - 159\cdot2^{k+6m+3l-3} +53\cdot2^{8m+4l-2} & \text{if    } i = 2m+l+2 \\
 2^{3k+2m+l} -7\cdot2^{2k+4m+2l} +7\cdot2^{k+6m+3l+1} -2^{8m+4l+3}& \text{if    } i = 2m+l+3
 \end{cases}\]
\end{thm}
 \begin{thm}
\label{thm 1.7}Let $ \Gamma_{i}^{\left[ (n)\atop{ 1+m\atop 1+m+l} \right]\times k}  $ denote the number of matrices of the form
 $\left[{A\over{B \over C}}\right] $of rank i such that  A is a  $ n\times k $ matrix  over $ \mathbb{F}_{2},$ B a 
$(1+m)\times k$ persymmetric matrix and C  a $ (1+m+l)\times k $ persymmetric  matrix  
 and where $\Gamma_{i}^{\left[ 1+m\atop 1+m+l \right]\times k}$ denote
   the number of double persymmetric $ (2+2m+l)\times k $  rank i matrices of the form  $\left[{B \over C}\right] $\\
     Then  $\Gamma_{i}^{\left[ (n)\atop{ 1+m\atop 1+m+l} \right]\times k} $ expressed as a linear combination of the $\Gamma_{i-j}^{\left[ 1+m\atop 1+m+l \right]\times k}$ is equal to
   \begin{equation*}
 \sum_{j= 0}^{n}2^{(n-j)\cdot(i-j)} a_{j}^{(n)}\prod_{l=1}^{j}(2^{k}- 2^{i-l})\cdot\Gamma_{i-j}^{\left[ 1+m\atop 1+m+l \right]\times k}
\quad for \quad 0\leq i\leq inf(k,n+2m+l+2),
 \end{equation*}
  where $ a_{j}^{(n)} $ satisfies the linear recurrence relation 
 \begin{equation*}
 a_{j}^{(n)} = 2^{j}\cdot a_{j}^{(n-1)} + a_{j-1}^{(n-1)},\quad n = 2,3,4,\ldots       \quad for\quad 1\leq j\leq n-1. 
\end{equation*}
 The explicit value of $ a_{j}^{(n)}$ is given by the following formula\\
  \begin{equation*}
 a_{j}^{(n)} =    (-1)^{j}\cdot2^{jn - {j(j-1)\over 2}} + \sum_{s =0}^{j-1} (-1)^{s}\prod_{l=0}^{j-(s+1)}{2^{n+1}- 2^{l}\over 2^{j-s}-2^{l}}\cdot2^{s(n-j) +{s(s+1)\over 2}}
  \quad for \quad  1\leq j\leq n-1.
   \end{equation*}
   We set  
   \begin{align*}
   a_{0}^{(n)} & =   a_{n}^{(n)} = 1    \\ and\quad
  \Gamma_{i-j}^{\left[ 1+m\atop 1+m+l \right]\times k}  & = 0 \quad if \quad  i-j \notin \{0,1,2,\ldots, inf(k,2m+l+2)\}.
   \end{align*}
  \end{thm}
  
\begin{thm}
\label{thm 1.8}

 Let $ f_{m,l,k}(t,\eta,\xi  _{1}, \xi _{2},\ldots, \xi _{n} ) $  be the exponential sum  in $ \mathbb{P}^{n+2} $ defined by\\
 \small
$ (t,\eta,\xi  _{1}, \xi _{2},\ldots, \xi _{n} ) \in \mathbb{P}^{n+2}\longrightarrow \\
  \sum_{deg Y\leq k-1}\sum_{deg Z\leq1+ m}E(tYZ)\sum_{deg U\leq1+ m+l}E(\eta YU)\sum_{deg V_{1}\leq  0}E(\xi _{1} YV_{1})
  \sum_{deg V_{2} \leq 0}E(\xi _{2} YV_{2}) \ldots \sum_{deg V_{n} \leq 0} E(\xi  _{n} YV_{n}). $\vspace{0.5 cm}\\
 \normalsize 
  Set $$(t,\eta,\xi  _{1}, \xi _{2},\ldots, \xi _{n} ) =
  \big(\sum_{i\geq 1}\alpha _{i}T^{-i}, \sum_{i\geq 1}\beta  _{i}T^{-i},\sum_{i\geq 1}\gamma  _{1i}T^{-i}), \ldots, \sum_{i\geq 1}\gamma  _{ni}T^{-i}) \in\mathbb{P}^{n+2}.   $$     
        
   Then
  $$ f_{m,l,k}(t,\eta,\xi _{1}, \xi _{2},\ldots, \xi _{n} ) = 
  2^{k+2m +l+n+4-r(D^{\left[(n)\atop{ 1+m\atop 1+m+l} \right]\times k}(\xi_{1},\xi _{2},\ldots,\xi _{n}, t,\eta  ))} $$
 where

$$D^{\left[(n)\atop{ 1+m\atop 1+m+l} \right]\times k}(\xi_{1},\xi _{2},\ldots,\xi _{n}, t,\eta  ) $$
 denotes  the following  $(2m+l+n+2) \times k $ matrix
   $$  \left ( \begin{array} {cccccc}
\gamma _{11} & \gamma   _{12} & \gamma  _{13} & \ldots  & \gamma  _{1,k-1} &  \gamma  _{1,k}  \\
 \gamma  _{21} & \gamma  _{22} & \gamma  _{23} & \ldots  & \gamma  _{2,k-1} &   \gamma _{2,k}\\
\vdots & \vdots & \vdots   &  \ldots  & \vdots  &  \vdots \\
\vdots & \vdots & \vdots    &  \ldots & \vdots  &  \vdots \\
 \gamma _{n1} & \gamma  _{n2} &  \gamma  _{n3} & \ldots  & \gamma  _{n,k-1} &  \gamma  _{n,k}\\
 \hline
 \alpha _{1} & \alpha _{2} & \alpha _{3} &  \ldots & \alpha _{k-1}  &  \alpha _{k} \\
\alpha _{2 } & \alpha _{3} & \alpha _{4}&  \ldots  &  \alpha _{k} &  \alpha _{k+1} \\
\vdots & \vdots & \vdots   &  \ldots  & \vdots  &  \vdots \\
\vdots & \vdots & \vdots    &  \ldots & \vdots  &  \vdots \\
\alpha _{1+m} & \alpha _{2+m} & \alpha _{3+m} & \ldots  &  \alpha _{k+m-1} &  \alpha _{k+m}  \\
\hline \\
\beta  _{1} & \beta  _{2} & \beta  _{3} & \ldots  &  \beta_{k-1} &  \beta _{k}  \\
\beta  _{2} & \beta  _{3} & \beta  _{4} & \ldots  &  \beta_{k} &  \beta _{k+1}  \\
\vdots & \vdots & \vdots    &  \vdots & \vdots  &  \vdots \\
\beta  _{1+m+l} & \beta  _{2+m+l} & \beta  _{3+m+l} & \ldots  &  \beta_{k+m+l-1} &  \beta _{k+m+l}  
  \end{array}  \right) $$
  
Then the number of solutions \\
 $(Y_1,Z_1,U_{1},V_{1}^{(1)},V_{2}^{(1)}, \ldots,V_{n}^{(1)}, Y_2,Z_2,U_{2},V_{1}^{(2)},V_{2}^{(2)}, 
\ldots,V_{n}^{(2)},\ldots  Y_q,Z_q,U_{q},V_{1}^{(q)},V_{2}^{(q)}, \ldots,V_{n}^{(q)}   ) $ \vspace{0.5 cm}\\
 of the polynomial equations  \vspace{0.5 cm}
  \[\left\{\begin{array}{c}
 Y_{1}Z_{1} +Y_{2}Z_{2}+ \ldots + Y_{q}Z_{q} = 0  \\
  Y_{1}U_{1} +Y_{2}U_{2}+ \ldots + Y_{q}U_{q} = 0 \\
   Y_{1}V_{1}^{(1)} + Y_{2}V_{1}^{(2)} + \ldots  + Y_{q}V_{1}^{(q)} = 0  \\
    Y_{1}V_{2}^{(1)} + Y_{2}V_{2}^{(2)} + \ldots  + Y_{q}V_{2}^{(q)} = 0\\
    \vdots \\
   Y_{1}V_{n}^{(1)} + Y_{2}V_{n}^{(2)} + \ldots  + Y_{q}V_{n}^{(q)} = 0 
 \end{array}\right.\]
    satisfying the degree conditions \\
                   $$  \deg Y_i \leq k-1 , \quad \deg Z_i \leq m +1 ,\quad \deg U_i \leq m +l+1 ,
                   \quad \deg V_{j}^{i} \leq 0 , \quad  for \quad 1\leq j\leq n  \quad 1\leq i \leq q $$ \\
                   
  is equal to the following integral over the unit interval in $ \mathbb{K}^{n+2} $
 $$ \int_{\mathbb{P}^{n+2}} f_{m,l,k}^{q}(t,\eta ,\xi _{1},\xi _{2},\ldots,\xi _{n} )dt d \eta d \xi _{1}d \xi _{2}\ldots d\xi_{n}. $$
  Observing that $ f_{m,l,k}^{q}(t,\eta ,\xi _{1},\xi _{2},\ldots,\xi _{n} )$ is constant on cosets of $ \mathbb{P}_{k+m}\times \mathbb{P}_{k+m+l}\times\mathbb{P}_{k}^{n}, $\;
  the above integral is equal to 
$$  2^{q(k+2m+l+n+4) - (2m+l+k(n+2))}\sum_{i = 0}^{inf(k,n+2m+l+2)}
  \Gamma_{i}^{\left[ (n)\atop{ 1+m\atop 1+m+l} \right]\times k} 2^{-iq}.$$
\end{thm}

\begin{example} n=2, m=1, $ l=3.$\;$ k\geq 5, $\; q = 3\vspace{0.1 cm}\\

\begin{align*}
& \Gamma_{i}^{\left[ (2)\atop{ 2\atop 5} \right]\times k}\\
& =2^{2i}\Gamma_{i}^{\left[ 2\atop 5 \right]\times k}+                                                                             
     3\cdot2^{i-1}(2^{k}-2^{i-1})\cdot\Gamma_{i-1}^{\left[ 2\atop 5 \right]\times k} 
     +(2^{k}-2^{i-1})(2^{k}-2^{i-2})\cdot \Gamma_{i-2}^{\left[ 2\atop 5 \right]\times k}
       \quad for\quad 0\leq i\leq \inf(k,9) 
\end{align*}

\[\Gamma_{i}^{\left[(2)\atop{ 1+1\atop 1+1+3} \right]\times k}  =
\begin{cases}
1 & \text{if    } i = 0\\
3\cdot2^{k} + 33 & \text{if    } i = 1 \\
 2^{2k} + 83\cdot2^{k} +886 & \text{if    } i = 2 \\
 33\cdot2^{2k} + 978\cdot2^{k} + 29352  &\text{if    } i = 3 \\
  2^{3k+1} + 182\cdot2^{2k} +16408\cdot2^{k} + 937408 & \text{if    } i = 4 \\
 3\cdot2^{3k+1} +189\cdot2^{2k+3} +8191\cdot2^{k+6}+12911\cdot2^{10}  & \text{if    } i = 5 \\ 
  3\cdot2^{3k+3} +36672\cdot2^{2k} +11271168\cdot2^{k} -399015936  & \text{if    } i = 6 \\
 11\cdot2^{3k+5} +1022464\cdot2^{2k} - 100679680\cdot2^{k} + 2163212288 & \text{if    } i = 7 \\  
 117\cdot2^{3k+7} - 3354624\cdot2^{2k} +214695936\cdot2^{k} - 3925868544 & \text{if    } i = 8 \\
  2^{4k+5} - 480\cdot2^{3k+5} +2240\cdot2^{2k+10} -1920\cdot2^{k+16} +2^{31} & \text{if    } i = 9
 \end{cases}\]
 If k = 5 we obtain \\
 
  \[\Gamma_{i}^{\left[(2)\atop{ 1+1\atop 1+1+3} \right]\times 5}  =
\begin{cases}
1 & \text{if    } i = 0\\
129 & \text{if    } i = 1 \\
 4566 & \text{if    } i = 2 \\
 94440 &\text{if    } i = 3 \\
  1714368& \text{if    } i = 4 \\
 31740928  & \text{if    } i = 5 \\ 
 \end{cases}\]
Then the number of solutions \\

 $(Y_1,Z_1,U_{1},V_{1}^{(1)},V_{2}^{(1)}, Y_2,Z_2,U_{2},V_{1}^{(2)},V_{2}^{(2)}, 
 Y_3,Z_3,U_{3},V_{1}^{(3)},V_{2}^{(3)}  ) $\\
 
 of the polynomial equations  \vspace{0.5 cm}
  \[\left\{\begin{array}{c}
 Y_{1}Z_{1} +Y_{2}Z_{2}+ Y_{3}Z_{3} = 0  \\
  Y_{1}U_{1} +Y_{2}U_{2}+ Y_{3}U_{3} = 0 \\
   Y_{1}V_{1}^{(1)} + Y_{2}V_{1}^{(2)} + Y_{3}V_{1}^{(3)} = 0  \\
    Y_{1}V_{2}^{(1)} + Y_{2}V_{2}^{(2)} + Y_{3}V_{2}^{(3)} = 0\\
    \end{array}\right.\]
    satisfying the degree conditions 
                   $$  \deg Y_i \leq 4 , \quad \deg Z_i \leq 2 ,\quad \deg U_i \leq 5 ,
                   \quad \deg V_{j}^{i} \leq 0 , \quad  for \quad 1\leq j\leq 2  \quad 1\leq i \leq 3 $$ 
 is equal to 
 $$  \int_{\mathbb{P}^{4}} f_{1,3,5}^{3}(t,\eta ,\xi _{1},\xi _{2} )dt d \eta d \xi _{1}d \xi _{2}= 2^{23}\sum_{i = 0}^{5}
  \Gamma_{i}^{\left[ (2)\atop{ 2\atop 5} \right]\times k} 2^{- 3i}= 13281\cdot2^{20}.$$
  
$$ \text{where for} \;(t,\eta,\xi  _{1}, \xi _{2} ) =
  \big(\sum_{i\geq 1}\alpha _{i}T^{-i}, \sum_{i\geq 1}\beta  _{i}T^{-i},\sum_{i\geq 1}\gamma  _{1i}T^{-i}), \sum_{i\geq 1}\gamma  _{2i}T^{-i}) \in\mathbb{P}^{4}.$$   
 \begin{align*}
&  f_{1,3,5}(t,\eta ,\xi _{1},\xi _{2} ) = 
   \sum_{deg Y\leq 4}\sum_{deg Z\leq 2}E(tYZ)\sum_{deg U\leq 5}E(\eta YU)\sum_{deg V_{1}\leq  0}E(\xi _{1} YV_{1})
  \sum_{deg V_{2} \leq 0}E(\xi _{2} YV_{2})\\
  &  =  2^{16-r(D^{\left[(2)\atop{ 2\atop 5} \right]\times k}(\xi_{1},\xi _{2} t,\eta  ))}
  \end{align*}
and \\

   $$ D^{\left[(2)\atop{ 2\atop 5} \right]\times k}(\xi_{1},\xi _{2}, t,\eta  ) = \left ( \begin{array} {cccccc}
\gamma _{11} & \gamma   _{12} & \gamma  _{13}  & \gamma  _{1,4} &  \gamma  _{1,5}  \\
 \gamma  _{21} & \gamma  _{22} & \gamma  _{23} &  \gamma  _{2,4} &   \gamma _{2,5}\\
 \hline
 \alpha _{1} & \alpha _{2} & \alpha _{3} &  \alpha _{4}  &  \alpha _{5} \\
\alpha _{2 } & \alpha _{3} & \alpha _{4} &  \alpha _{5} &  \alpha _{6} \\
\hline \\
\beta  _{1} & \beta  _{2} & \beta  _{3}  &  \beta_{4} &  \beta _{5}  \\
\beta  _{2} & \beta  _{3} & \beta  _{4}  &  \beta_{5} &  \beta _{6}  \\
\beta  _{3} & \beta  _{4} & \beta  _{5}  &  \beta_{6} &  \beta _{7}  \\
\beta  _{4} & \beta  _{5} & \beta  _{6}  &  \beta_{7} &  \beta _{8}  \\
\beta  _{5} & \beta  _{6} & \beta  _{7}  &  \beta_{8} &  \beta _{9}  \\
 \end{array}  \right) $$
 \end{example}
  \subsection{Some Lemmas}
  \label{subsec 4}
  \begin{lem}
\label{lem 1.9}
Let $ (t,\eta ,\xi ) \in  \mathbb{P}\times \mathbb{P}\times \mathbb{P}  $ and
$$ g_{s,m,k}(t,\eta,\xi  ) = g(t,\eta,\xi  ) = \sum_{deg Y\leq k-1}\sum_{deg Z\leq s-1}E(tYZ)\sum_{deg U\leq  s+m-1}E(\eta YU)\sum_{deg V = 0}E(\xi YV). $$  \\
Then \\
  \begin{equation*}
 g(t,\eta,\xi  ) = \begin{cases}
 2^{k+2s+m-  r(D^{\big[\stackrel{s}{s+m}\big] \times k }(t,\eta ) ) }  & \text{if }
 r(D^{\big[\stackrel{s}{s+m}\big] \times k }(t,\eta ) ) =  r( D^{\left[\stackrel{s}{\stackrel{s+m}{1}}\right] \times k }(t,\eta ,\xi ) ) \\
     0  & \text{otherwise }.
    \end{cases}
\end{equation*}
\end{lem}

\begin{proof}We generalize the  argument given in the proof of Lemma 6.2 in [1]\\
 Consider the matrix \\
$$ D^{\left[s\atop{ s+m\atop 1} \right]\times k}(t,\eta ,\xi ) =  
       \left ( \begin{array} {cccccc}
\alpha _{1} & \alpha _{2} & \alpha _{3} &  \ldots & \alpha _{k-1}  &  \alpha _{k} \\
\alpha _{2 } & \alpha _{3} & \alpha _{4}&  \ldots  &  \alpha _{k} &  \alpha _{k+1} \\
\vdots & \vdots & \vdots    &  \vdots & \vdots  &  \vdots \\
\alpha _{s-1} & \alpha _{s} & \alpha _{s +1} & \ldots  &  \alpha _{s+k-3} &  \alpha _{s+k-2}  \\
\alpha _{s} & \alpha _{s+1} & \alpha _{s +2} & \ldots  &  \alpha _{s+k-2} &  \alpha _{s+k-1}  \\
\beta  _{1} & \beta  _{2} & \beta  _{3} & \ldots  &  \beta_{k-1} &  \beta _{k}  \\
\beta  _{2} & \beta  _{3} & \beta  _{4} & \ldots  &  \beta_{k} &  \beta _{k+1}  \\
\vdots & \vdots & \vdots    &  \vdots & \vdots  &  \vdots \\
\beta  _{m+1} & \beta  _{m+2} & \beta  _{m+3} & \ldots  &  \beta_{k+m-1} &  \beta _{k+m}  \\
\vdots & \vdots & \vdots    &  \vdots & \vdots  &  \vdots \\
\beta  _{s+m-1} & \beta  _{s+m} & \beta  _{s+m+1} & \ldots  &  \beta_{s+m+k-3} &  \beta _{s+m+k-2}  \\
\beta  _{s+m} & \beta  _{s+m+1} & \beta  _{s+m+2} & \ldots  &  \beta_{s+m+k-2} &  \beta _{s+m+k-1} \\
\hline \\
\gamma  _{1} & \gamma  _{2} & \gamma  _{3} & \ldots  &  \gamma _{k-1} &  \gamma  _{k}  
\end{array}  \right). $$ \vspace{0.5 cm}\\
It follows\\
\begin{align*}
 g(t,\eta ,\xi ) &  = 2^{2s+m}\sum_{deg Y\leq k-1\atop {Y\in ker D^{\big[\stackrel{s}{s+m}\big] \times k }(t,\eta ) }}E(\xi Y) \\
& = 2^{2s+m}\cdot\Big[\sum_{deg Y\leq k-1\atop {Y\in ker D^{\left[\stackrel{s}{\stackrel{s+m}{1}}\right] \times k }(t,\eta ,\xi )}} 1
 -   \sum_{\substack{deg Y\leq k-1 \\
 Y\in kerD^{\big[\stackrel{s}{s+m}\big] \times k }(t,\eta ) \\
 \sum_{j = 0}^{k -1}\gamma  _{j}\delta _{j}= 1}}1 \Big] \\
 & = 2^{2s+m}\cdot\Big[2\cdot 2^{k - r( D^{\left[\stackrel{s}{\stackrel{s+m}{1}}\right] \times k }(t,\eta ,\xi ))}
 - 2^{k - r( D^{\big[\stackrel{s}{s+m}\big] \times k }(t,\eta ) ) }\Big].
\end{align*}
\end{proof}

\begin{lem}
\label{lem 1.10}
Set  $  h(t,\eta )  = \sum_{deg \leq k-1}\sum_{deg Z \leq s-1}E(tYZ)\sum_{deg U \leq s+m-1}E(tYU) $ and let $ q\geq 2 $ be an integer, then \\
            $$ g^q(t,\eta,\xi  ) = g(t,\eta,\xi  )\cdot h^{q-1}(t,\eta ). $$
\end{lem}
\begin{proof}
Obviously we get 
$$ g^2(t,\eta,\xi  ) =\big[ 2^{2s+m}\sum_{deg Y_{1}\leq k-1\atop {Y_{1}\in \ker D^{\big[\stackrel{s}{s+m}\big] \times k }(t,\eta ) }}E(\xi Y_{1}) \big]\cdot
 \big[ 2^{2s+m}\sum_{deg Y_{2}\leq k-1\atop {Y_{2}\in \ker D^{\big[\stackrel{s}{s+m}\big] \times k }(t,\eta ) }} E(\xi Y_{2}) \big]. $$
 Define
\[
\left\{\begin{array}{cc}
Y_{1} + Y_{2} = Y_{3}  &  deg Y_{3}\leq k-1 \\
              Y_{1} = Y_{4}   &   deg Y_{4}\leq k-1.
\end{array}\right.\]
Then we obtain  \\
\begin{align*}
 g^2(t,\eta,\xi  ) & = 2^{2(2s+m)}\sum_{deg Y_{1}\leq k-1\atop {Y_{1}\in \ker D^{\big[\stackrel{s}{s+m}\big] \times k }(t,\eta ) }}
 \sum_{deg Y_{2}\leq k-1\atop {Y_{2}\in \ker D^{\big[\stackrel{s}{s+m}\big] \times k }(t,\eta ) }}E(\xi  (Y_{1}+ Y_{2}))  \\
 & = \big[ 2^{2s+m}\sum_{deg Y_{4}\leq k-1\atop {Y_{4}\in \ker D^{\big[\stackrel{s}{s+m}\big] \times k }(t,\eta ) }} 1 \big]\cdot
 \big[ 2^{1+m}\sum_{deg Y_{3}\leq k-1\atop {Y_{3}\in \ker D^{\big[\stackrel{s}{s+m}\big] \times k }(t,\eta ) }} E(\xi  Y_{ 3})\big] \\
 & = h(t,\eta )\cdot g(t,\eta,\xi  ).
\end{align*}
By recurrence on q we get   $g^q(t,\eta,\xi  ) = g(t,\eta,\xi  )\cdot h^{q-1}(t,\eta ). $
\end{proof}
\begin{lem}
\label{lem 1.11} We have
$$\int_{\mathbb{P}}\int_{\mathbb{P}}\int_{\mathbb{P}}g^q(t,\eta,\xi  ) dt d\eta d\xi  =
2^{q(k+2s+m ) - (3k+2s+m-2)}\sum_{i = 0}^{\inf(k,2s+m)}
\sigma _{i,i}^{\left[\stackrel{s}{\stackrel{s+m }{\overline{\gamma _{1 -}}}}\right] \times k } 2^{-iq} $$
 \end{lem}
 \begin{proof}  Recall  Definition \ref{defn 1.3},  then Lemma \ref{lem 1.9} gives ,  by observing that $ g(t,\eta,\xi  ) $ is
 constant on cosets of   $\mathbb{P}_{k+s-1} \times \mathbb{P}_{k+s+m-1}\times \mathbb{P}_{k}$\\
  \begin{align*}
\int_{\mathbb{P}}\int_{\mathbb{P}}\int_{\mathbb{P}}g^q(t,\eta,\xi  ) dt d\eta d\xi & =
 \sum_{(t,\eta,\xi  )\mathbb{P}/\mathbb{P}_{k+s-1}\times \mathbb{P}/\mathbb{P}_{k+s+m-1}\times \mathbb{P}/\mathbb{P}_{k}\atop
 {  r(D^{\big[\stackrel{s}{s+m}\big] \times k }(t,\eta ) ) =  r( D^{\left[\stackrel{s}{\stackrel{s+m}{1}}\right] \times k }(t,\eta ,\xi ) ) }}
  2^{k+2s+m-  r(D^{\big[\stackrel{s}{s+m}\big] \times k }(t,\eta ) ) }  \int_{\mathbb{P}_{k+s-1}}dt \int_{\mathbb{P}_{k+s+m-1}}d\eta \int_{\mathbb{P}_{k}}d\xi \\
 & = \sum_{i = 0}^{\inf(k,2s+m)}\sigma _{i,i}^{\left[\stackrel{s}{\stackrel{s+m }{\overline{\gamma _{1 -}}}}\right] \times k }  2^{q(k+2s+m- i)}2^{-(k+s-1)}2^{-(k+s+m-1)}2^{-k}.
  \end{align*}
 \end{proof}
 \begin{lem} Equally
\label{lem 1.12}
$$\int_{\mathbb{P}}\int_{\mathbb{P}}\int_{\mathbb{P}}g^q(t,\eta,\xi  ) dt d\eta d\xi  =
2^{q(k+2s+m ) - (3k+2s+m-2)}\sum_{i = 0}^{\inf(k,2s+m)}
 2^{i}\cdot\Gamma_{i}^{\left[s\atop s+m\right]\times k}2^{-iq} $$
\end{lem}
 \begin{proof}
From  Lemma  \ref{lem 1.10} recalling  Definition \ref{defn 1.2} we get by Fubini's theorem  \\ 
 \begin{align*}
& \int_{\mathbb{P}}\int_{\mathbb{P}}\int_{\mathbb{P}}g^q(t,\eta,\xi  ) dt d\eta d\xi   =
\int_{\mathbb{P}}\int_{\mathbb{P}}\int_{\mathbb{P}}g(t,\eta,\xi  )h^{q-1}(t,\eta ) dt d\eta d\xi \\
&  =\int_{\mathbb{P}}\int_{\mathbb{P}}h^{q-1}(t,\eta )\big(\int_{\mathbb{P}}g(t,\eta,\xi  )d\xi \big)dtd\eta  
 =  2^{2s+m}\int_{\mathbb{P}}\int_{\mathbb{P}}h^{q-1}(t,\eta )dtd\eta  \\
 & = 2^{2s+m}\cdot  \sum_{(t,\eta )\in \mathbb{P}/\mathbb{P}_{k+s-1}\times \mathbb{P}/\mathbb{P}_{k+s+m-1}}
 2^{(q-1)(2s+m+k- r(D^{\big[\stackrel{s}{s+m}\big] \times k }(t,\eta ) )}\int_{\mathbb{P}_{k+s-1}}dt \int_{\mathbb{P}_{k+s+m-1}}dt d\eta \\
 & =  2^{2s+m}\cdot \sum_{i = 0}^{\inf(2s+m,k)} \sum_{(t,\eta )\in \mathbb{P}/\mathbb{P}_{k+s-1}\times \mathbb{P}/\mathbb{P}_{k+s+m-1}}
  2^{(q-1)(2s+m+k- i) }\int_{\mathbb{P}_{k+s-1}}dt \int_{\mathbb{P}_{k+s+m-1}}dt d\eta \\
  & =  2^{2s+m}\cdot \sum_{i = 0}^{\inf(2s+m,k)}\Gamma _{i}^{\Big[\substack{s  \\ s+m }\Big] \times k} 
  \cdot2^{-(q-1)i}\cdot2^{(q-1)(2s+m+k)}\cdot2^{-(k+s-1)}\cdot2^{-(k+s+m-1)}.
\end{align*}
\end{proof}

\begin{lem}
\label{lem 1.13}  For all $ 0\leq i\leq \inf(k,2s+m) $ we have $\sigma _{i,i}^{\left[\stackrel{s}{\stackrel{s+m }{\overline{\gamma _{1 -}}}}\right] \times k }$ =
$ 2^{i}\cdot\Gamma_{i}^{\left[s\atop s+m\right]\times k} $ 
\end{lem}
\begin{proof}
Lemma \ref{lem 1.11}  and Lemma \ref{lem 1.12} give \\
$$\sum_{i=0}^{\inf(k,2s+m)}\big( \sigma _{i,i}^{\left[\stackrel{s}{\stackrel{s+m }{\overline{\gamma _{1 -}}}}\right] \times k }  -  2^{i}\cdot
\Gamma_{i}^{\left[s\atop s+m\right]\times k}\big)\cdot2^{-iq} = 0 \quad  for\; all\; q \geq 2. $$\\
\end{proof}
\begin{lem}
\label{lem 1.14} For all  $1\leq i\leq inf(k,2s+m+1) $  we have\\
 $$ \sigma _{i-1,i}^{\left[\stackrel{s}{\stackrel{s+m }{\overline{\gamma _{1 -}}}}\right] \times k } 
  = (2^k - 2^{i-1})\cdot\Gamma_{i-1}^{\left[s\atop s+m\right]\times k}$$
\end{lem}
\begin{proof}
Consider the matrix \\
 $$  \left[{D_{s  \times k}(t)\over{ D_{(s+m )\times k}(\eta )\over D_{1\times k}(\xi )}}\right] = \left ( \begin{array} {cccccc}
\alpha _{1} & \alpha _{2} & \alpha _{3} &  \ldots & \alpha _{k-1}  &  \alpha _{k} \\
\alpha _{2 } & \alpha _{3} & \alpha _{4}&  \ldots  &  \alpha _{k} &  \alpha _{k+1} \\
\vdots & \vdots & \vdots    &  \vdots & \vdots  &  \vdots \\
\alpha _{s-1} & \alpha _{s} & \alpha _{s +1} & \ldots  &  \alpha _{s+k-3} &  \alpha _{s+k-2}  \\
\alpha _{s} & \alpha _{s+1} & \alpha _{s +2} & \ldots  &  \alpha _{s+k-2} &  \alpha _{s+k-1}  \\
\hline \\
\beta  _{1} & \beta  _{2} & \beta  _{3} & \ldots  &  \beta_{k-1} &  \beta _{k}  \\
\beta  _{2} & \beta  _{3} & \beta  _{4} & \ldots  &  \beta_{k} &  \beta _{k+1}  \\
\vdots & \vdots & \vdots    &  \vdots & \vdots  &  \vdots \\
\beta  _{m+1} & \beta  _{m+2} & \beta  _{m+3} & \ldots  &  \beta_{k+m-1} &  \beta _{k+m}  \\
\vdots & \vdots & \vdots    &  \vdots & \vdots  &  \vdots \\
\beta  _{s+m-1} & \beta  _{s+m} & \beta  _{s+m+1} & \ldots  &  \beta_{s+m+k-3} &  \beta _{s+m+k-2}  \\
\beta  _{s+m} & \beta  _{s+m+1} & \beta  _{s+m+2} & \ldots  &  \beta_{s+m+k-2} &  \beta _{s+m+k-1} \\
\hline \\
\gamma  _{1} & \gamma  _{2} & \gamma  _{3} & \ldots  &  \gamma _{k-1} &  \gamma  _{k}  
\end{array}  \right). $$ \vspace{0.5 cm}\\
Obviously by Lemma \ref{lem 1.13}
\begin{align*}
2^{k}\cdot\Gamma _{i-1}^{\left[\stackrel{s}{s+m}\right]\times k}  & =
 \sigma _{i-1,i-1}^{\left[\stackrel{s}{\stackrel{s+m }{\overline{\gamma _{1 -}}}}\right] \times k } 
 +  \sigma _{i-1,i}^{\left[\stackrel{s}{\stackrel{s+m }{\overline{\gamma _{1 -}}}}\right] \times k } \\
& \Leftrightarrow  \sigma _{i-1,i}^{\left[\stackrel{s}{\stackrel{s+m }{\overline{\gamma _{1 -}}}}\right] \times k }
  = 2^{k}\cdot\Gamma _{i-1}^{\left[\stackrel{s}{s+m}\right]\times k}  - 2^{i-1}\cdot\Gamma _{i-1}^{\left[\stackrel{s}{s+m}\right]\times k}\\
&  \Leftrightarrow \sigma _{i-1,i}^{\left[\stackrel{s}{\stackrel{s+m }{\overline{\gamma _{1 -}}}}\right] \times k }  =
 ( 2^{k} -  2^{i-1})\cdot \Gamma _{i-1}^{\left[\stackrel{s}{s+m}\right]\times k} && \text{if } 1\leq i\leq \inf(k, 2s+m+1).
 \end{align*}
\end{proof}
\begin{lem}
\label{lem 1.15}
For all  $1\leq i\leq inf(k,2s+m+1) $ we have \\
\begin{equation}
\label{eq 1.1}
  \Gamma_{i}^{\left[s\atop{ s+m\atop 1} \right]\times k}= ( 2^{k} -  2^{i-1})\cdot \Gamma _{i-1}^{\left[\stackrel{s}{s+m}\right]\times k}
 +  2^{i}\cdot\Gamma _{i}^{\left[\stackrel{s}{s+m}\right]\times k}
\end{equation}
 \end{lem}
\begin{proof} From Lemmas \ref{lem 1.13}, \ref{lem 1.14} we get \\
 $$  \Gamma_{i}^{\left[s\atop{ s+m\atop 1} \right]\times k}
=  \sigma _{i,i}^{\left[\stackrel{s}{\stackrel{s+m }{\overline{\gamma _{1 -}}}}\right] \times k }   + \sigma _{i-1,i}^{\left[\stackrel{s}{\stackrel{s+m }{\overline{\gamma _{1 -}}}}\right] \times k } 
 = 2^{i}\cdot\Gamma_{i}^{\left[s\atop s+m\right]\times k} + ( 2^{k} -  2^{i-1})\cdot\Gamma_{i-1}^{\left[s\atop s+m\right]\times k}. $$
  \end{proof}
 \begin{lem}
\label{lem 1.16} 
  We have the following reduction formulas \\
  \begin{align}
\Gamma_{m+3+j}^{\left[1\atop{ 1+m\atop 1+m+l} \right]\times k} & = 2^{4j}\Gamma_{m+3}^{\left[1\atop{ 1+m\atop 1+m+(l-j)} \right]\times (k-j)}   &  \text{if  }   0\leq j\leq l  \label{eq 1.2} \\
\Gamma_{m+l+3+j}^{\left[1\atop{ 1+m\atop 1+m+l} \right]\times k} & = 2^{8j+4l}\Gamma_{m+3-j}^{\left[1\atop{ 1+(m-j)\atop 1+(m-j)} \right]\times (k-(2j+l))} & \text{if  }  0\leq j\leq m \label{eq 1.3}
\end{align}
 \end{lem}
 \begin{proof} 
  From  \eqref{eq 1.1} we get with $s \rightarrow 1+m,\;m\rightarrow l$ \\
   \begin{equation}
\label{eq 1.4}
    \Gamma_{i}^{\left[1\atop{ 1+m\atop 1+m+l} \right]\times k}= ( 2^{k} -  2^{i-1})\cdot \Gamma _{i-1}^{\left[\stackrel{1+m}{1+m+l}\right]\times k}
 +  2^{i}\cdot\Gamma _{i}^{\left[\stackrel{1+m}{1+m+l}\right]\times k} 
\end{equation}
Using the reduction formulas proved in [2, section 3] we get from \eqref{eq 1.4}with $i\rightarrow m+3+j. $\\
\begin{align*}  
 & \Gamma_{m+3+j}^{\left[1\atop{ 1+m\atop 1+m+l} \right]\times k}= ( 2^{k} -  2^{m+2+j})\cdot \Gamma _{m+2+j}^{\left[\stackrel{1+m}{1+m+l}\right]\times k}
 +  2^{m+3+j}\cdot \Gamma _{m+3+j}^{\left[\stackrel{1+m}{1+m+l}\right]\times k} \\
 & = ( 2^{k} -  2^{m+2+j})\cdot2^{3j}\cdot \Gamma _{m+2}^{\left[\stackrel{1+m}{1+m+(l-j)}\right]\times (k-j)}
 +  2^{m+3+j}\cdot2^{3j}\cdot\Gamma _{m+3}^{\left[\stackrel{1+m}{1+m+(l-j)}\right]\times (k-j)} \\
 & =  ( 2^{k+3j} -  2^{m+2+4j})\cdot \Gamma _{m+2}^{\left[\stackrel{1+m}{1+m+(l-j)}\right]\times (k-j)}
 +  2^{m+3+4j}\cdot \Gamma _{m+3}^{\left[\stackrel{1+m}{1+m+(l-j)}\right]\times (k-j)} \\
 & \\
  & \text{From \eqref{eq 1.4}with $i\rightarrow m+3,\;l\rightarrow l-j,\;k\rightarrow k-j$ we obtain}\\
 & \Gamma_{m+3}^{\left[1\atop{ 1+m\atop 1+m+(l-j)} \right]\times (k-j)}=  ( 2^{k-j} -  2^{m+2})\cdot \Gamma _{m+2}^{\left[\stackrel{1+m}{1+m+(l-j)}\right]\times (k-j)}
 +  2^{m+3}\cdot\Gamma _{m+3}^{\left[\stackrel{1+m}{1+m+(l-j)}\right]\times (k-j)} \\
 & \text{which proves \eqref{eq 1.2}}
\end{align*}

We prove \eqref{eq 1.3} in the same way, that is\\
\begin{align*}  
 & \Gamma_{m+l+3+j}^{\left[1\atop{ 1+m\atop 1+m+l} \right]\times k}= ( 2^{k} -  2^{m+l+2+j})\cdot \Gamma _{m+l+2+j}^{\left[\stackrel{1+m}{1+m+l}\right]\times k}
 +  2^{m+l+3+j}\cdot\Gamma _{m+l+3+j}^{\left[\stackrel{1+m}{1+m+l}\right]\times k}\\
  & = ( 2^{k} -  2^{m+l+2+j})\cdot2^{6j+3l}\cdot \Gamma _{m+2-j}^{\left[\stackrel{m+1-j}{m+1-j}\right]\times (k-l-2j)}
 +  2^{m+l+3+j}\cdot2^{6j+3l}\cdot\Gamma _{m+3-j}^{\left[\stackrel{m+1-j}{m+1-j}\right]\times (k-l-2j)}\\
 & =  ( 2^{k+6j+3l+3} -  2^{m+4l+7j+2})\cdot \Gamma _{m+2}^{\left[\stackrel{m+1-j}{m+1-j}\right]\times (k-l-2j)}
 +  2^{m+4l+7j+3}\cdot\Gamma _{m+3}^{\left[\stackrel{m+1-j}{m+1-j}\right]\times (k-l-2j)} \\
 & \\
  & \text{From \eqref{eq 1.4}with $i\rightarrow m+3-j,\;m\rightarrow m-j,\;l\rightarrow  0,\;k\rightarrow k-(2j+l)$ we obtain}\\
 & \Gamma_{m+3-j}^{\left[1\atop{ 1+(m-j)\atop 1+(m-j)} \right]\times (k-(2j+l))} =  ( 2^{k-(2j+l)} -  2^{m+2-j})\cdot \Gamma _{m+2-j}^{\left[\stackrel{1+(m-j)}{1+(m-j)}\right]\times (k-(2j+l))}
 +  2^{m+3-j}\cdot \Gamma _{m+3-j}^{\left[\stackrel{1+(m-j)}{1+(m-j)}\right]\times (k-(2j+l))}  \\
 & \text{which proves \eqref{eq 1.3}}
\end{align*}
  \end{proof}
 \begin{lem}
\label{lem 1.17}
We have for $m\geq 2,\;k>m+3$\\
\[ \Gamma_{m+3}^{\left[1\atop{ 1+m\atop 1+m+l} \right]\times k}  =
\begin{cases}
 3\cdot2^{2k+m} + 21\cdot2^{k+3m+2} + 21\cdot2^{k+2m+1} +315\cdot2^{4m+4} - 315\cdot2^{3m+3}& \text{if    } l\geq 3 \\
 3\cdot2^{2k+m} + 21\cdot2^{k+3m +2} + 149\cdot2^{k+2m+1} +315\cdot2^{4m+4} - 827\cdot2^{3m+3} &\text{if    } l =2 \\
 11\cdot2^{2k+m} + 21\cdot2^{k+3m+2} + 261\cdot2^{k+2m+1} +315\cdot2^{4m+4} - 1627\cdot2^{3m+3}& \text{if    } l=1 \\
 21\cdot2^{2k+m} + 21\cdot2^{k+3m+2} + 525\cdot2^{k+2m+1} +315\cdot2^{4m+4} - 3255\cdot2^{3m+3}& \text{if    } l=0 
 \end{cases}\]
\end{lem}
\begin{proof}
\underline{The case $l\geq 3$}\\
From \eqref{eq 1.4} we get with $i\rightarrow m+3$\\
\begin{align*}
& \Gamma_{m+3}^{\left[1\atop{ 1+m\atop 1+m+l} \right]\times k}= ( 2^{k} -  2^{m+2})\cdot \Gamma _{m+2}^{\left[\stackrel{1+m}{1+m+l}\right]\times k}
 +  2^{m+3}\cdot\Gamma _{m+3}^{\left[\stackrel{1+m}{1+m+l}\right]\times k}\\
 & \text{Applying Theorem 3.13 in [1, see section 3] with $s\rightarrow 1+m,\;m\rightarrow l$ we obtain}\\
&  \Gamma_{m+3}^{\left[1\atop{ 1+m\atop 1+m+l} \right]\times k}= ( 2^{k} -  2^{m+2})\cdot[3\cdot2^{k+m} +21\cdot(2^{3m+2} -2^{2m+1})]
 +  2^{m+3}\cdot[3\cdot2^{k+m+2} +21\cdot(2^{3m+5} -2^{2m+4})]\\
 & =  3\cdot2^{2k+m} + 21\cdot2^{k+3m+2} + 21\cdot2^{k+2m+1} +315\cdot2^{4m+4} - 315\cdot2^{3m+3}
\end{align*}
The other cases are proved in a similar way by applying  results proved in [2, see section 3]

\end{proof}
  
\begin{lem}
\label{lem 1.18}
We have for $ m\geq 2 $ and $l\geq 3 $ and  $k \geq 2m+l $\\
\small
$ \Gamma_{i}^{\left[1\atop{ 1+m\atop 1+m+l} \right]\times k}  $  is equal to \\
\[ =
\begin{cases}
1 & \text{if    } i = 0\\
2^{k} +17 & \text{if    } i = 1 \\
 (21\cdot2^{3i-7} - 3\cdot2^{2i-5})\cdot2^{k} + 315\cdot2^{4i-8} -21\cdot2^{3i-6}& \text{if    }2\leq  i\leq m \\
 21\cdot2^{k+3m-4} + 13\cdot2^{k+2m-3} +315\cdot2^{4m-4} -85\cdot2^{3m-3} & \text{if    } i = m+1 \\
 2^{2k+m} + 21\cdot2^{k+3m-1} + 9\cdot2^{k+2m-1} +315\cdot2^{4m} - 157\cdot2^{3m} &\text{if    } i = m+2 \\
 3\cdot2^{2k-m+2i-6} + 21\cdot2^{k+3i-7} + 21\cdot2^{k-m+3i-8} +315\cdot2^{4i-8} - 315\cdot2^{4i-m-9}& \text{if    } m+3\leq i\leq m+l \\
  3\cdot2^{2k+m+2l-4} + 21\cdot2^{k+3m +3l-4} + 149\cdot2^{k+2m+3l-5} +315\cdot2^{4m+4l-4} - 827\cdot2^{3m+4l-5} &\text{if    } i = m+l+1 \\
 11\cdot2^{2k+m+2l-2} + 21\cdot2^{k+3m +3l-1} + 261\cdot2^{k+2m+3l-2} +315\cdot2^{4m+4l} - 1627\cdot2^{3m+4l-1}& \text{if    } i = m+l+2 \\ 
21\cdot2^{2k-2m-l+3i-9} + 21\cdot2^{k+3i-7} + 525\cdot2^{k-2m-l+4i-11} +315\cdot2^{4i-8} - 3255\cdot2^{5i-2m-l-12}& \text{if    } m+l+3\leq i\leq 2m+l+1 \\ 
53\cdot2^{2k+4m+2l-3} - 159\cdot2^{k+6m+3l-3} +53\cdot2^{8m+4l-2} & \text{if    } i = 2m+l+2 \\
 2^{3k+2m+l} -7\cdot2^{2k+4m+2l} +7\cdot2^{k+6m+3l+1} -2^{8m+4l+3}& \text{if    } i = 2m+l+3
 \end{cases}\]
\end{lem}
\begin{proof}
Using \eqref{eq 1.4} and the results obtained in [2, see section 3] together  with the Lemmas \ref{lem 1.16} and \ref{lem 1.17}
we can easily  compute $ \Gamma_{i}^{\left[1\atop{ 1+m\atop 1+m+l} \right]\times k}  $ for $0\leq i\leq 2m+l+3$

\end{proof}

\begin{example}We have for $ m=2,\; l=3, \;k\geq 7$\\
 
\[\Gamma_{i}^{\left[1\atop{ 1+2\atop 1+2+3} \right]\times k}  =
\begin{cases}
1 & \text{if    } i = 0\\
2^{k} +17 & \text{if    } i = 1 \\
 9\cdot2^{k} + 294  & \text{if    } i =2 \\
 55\cdot2^{k+1} +4360 & \text{if    } i = 3 \\
 2^{2k+2} + 93\cdot2^{k+3} + 70592 &\text{if    } i = 4 \\
  3\cdot2^{2k+2} + 189\cdot2^{k+5} + 1128960  & \text{if    } i =5 \\
  3\cdot2^{2k+4} + 317\cdot2^{k+5} + 13869056 &\text{if    } i = 6 \\
  11\cdot2^{2k+6} + 429\cdot2^{k+11} + 893\cdot2^{17}  & \text{if    } i = 7 \\ 
 21\cdot2^{2k+8} + 693\cdot2^{k+14} - 1541406720  & \text{if    } i = 8 \\ 
53\cdot2^{2k+11} -159\cdot2^{k+18} +53\cdot2^{26}  & \text{if    } i = 9 \\ 
 2^{3k+7} -7\cdot2^{2k+14} +7\cdot2^{k+22} -2^{31}& \text{if    } i = 10
 \end{cases}\]
\end{example}
\begin{lem}
\label{lem 1.19}
We have for $ m\geq 2,\; l=2 $ and  $k \geq 2m+2 $\\
\small
\[\Gamma_{i}^{\left[1\atop{ 1+m\atop 1+m+2} \right]\times k} =
\begin{cases}
1 & \text{if    } i = 0\\
2^{k} +17 & \text{if    } i = 1 \\
 (21\cdot2^{3i-7} - 3\cdot2^{2i-5})\cdot2^{k} + 315\cdot2^{4i-8} -21\cdot2^{3i-6}& \text{if    }2\leq  i\leq m \\
 21\cdot2^{k+3m-4} + 13\cdot2^{k+2m-3} +315\cdot2^{4m-4} -85\cdot2^{3m-3} & \text{if    } i = m+1 \\
 2^{2k+m} + 21\cdot2^{k+3m-1} + 9\cdot2^{k+2m-1} +315\cdot2^{4m} - 157\cdot2^{3m} &\text{if    } i = m+2 \\
3\cdot2^{2k+m} + 21\cdot2^{k+3m+2} + 149\cdot2^{k+2m+1} +315\cdot2^{4m+4} - 827\cdot2^{3m+3} &\text{if    } i = m+3 \\ 11\cdot2^{2k+m+2} + 21\cdot2^{k+3m +5} + 261\cdot2^{k+2m+4} +315\cdot2^{4m+8} - 1627\cdot2^{3m+7}& \text{if    } i = m+4 \\ 
 21\cdot2^{2k-2m+3i-11} + 21\cdot2^{k+3i-7} + 525\cdot2^{k-2m+4i-13} +315\cdot2^{4i-8} - 3255\cdot2^{5i-2m-14}& \text{if    } m+5\leq i\leq 2m+3 \\ 
53\cdot2^{2k+4m+1} - 159\cdot2^{k+6m+3} +53\cdot2^{8m+6} & \text{if    } i = 2m+4 \\
 2^{3k+2m+2} -7\cdot2^{2k+4m+4} +7\cdot2^{k+6m+7} -2^{8m+11}& \text{if    } i = 2m+5
 \end{cases}\]
\end{lem}

\begin{example}We have for $ m=2,\; l=2, \;k\geq 6$\\
 
\[\Gamma_{i}^{\left[1\atop{ 1+2\atop 1+2+2} \right]\times k}  =
\begin{cases}
1 & \text{if    } i = 0\\
2^{k} +17 & \text{if    } i = 1 \\
 9\cdot2^{k} + 294  & \text{if    } i =2 \\
 55\cdot2^{k+1} +4360 & \text{if    } i = 3 \\
 2^{2k+2} + 93\cdot2^{k+3} + 70592 &\text{if    } i = 4 \\
  3\cdot2^{2k+2} + 317\cdot2^{k+5} + 866816 &\text{if    } i = 5 \\
 11\cdot2^{2k+4} + 429\cdot2^{k+8} + 7315456  & \text{if    } i = 6 \\ 
21\cdot2^{2k+6} + 693\cdot2^{k+11} - 735\cdot2^{17}  & \text{if    } i = 7 \\ 
53\cdot2^{2k+9} -159\cdot2^{k+15} +53\cdot2^{22}  & \text{if    } i = 8 \\ 
 2^{3k+4} -7\cdot2^{2k+12} +7\cdot2^{k+19} -2^{27}& \text{if    } i = 9
 \end{cases}\]
\end{example}

\begin{lem}
\label{lem 1.20}
We have for $ m\geq 2 $ and $l=1 $ and  $k \geq 2m+1 $\\
\small
$ \Gamma_{i}^{\left[1\atop{ 1+m\atop 1+m+1} \right]\times k}  $  is equal to \\
\[ =
\begin{cases}
1 & \text{if    } i = 0\\
2^{k} +17 & \text{if    } i = 1 \\
 (21\cdot2^{3i-7} - 3\cdot2^{2i-5})\cdot2^{k} + 315\cdot2^{4i-8} -21\cdot2^{3i-6}& \text{if    }2\leq  i\leq m \\
 21\cdot2^{k+3m-4} + 13\cdot2^{k+2m-3} +315\cdot2^{4m-4} -85\cdot2^{3m-3} & \text{if    } i = m+1 \\
 2^{2k+m} + 21\cdot2^{k+3m-1} + 73\cdot2^{k+2m-1} +315\cdot2^{4m} - 413\cdot2^{3m} &\text{if    } i = m+2 \\
 11\cdot2^{2k+m} + 21\cdot2^{k+3m+2} + 261\cdot2^{k+2m+1} +315\cdot2^{4m+4} - 1627\cdot2^{3m+3}& \text{if    } i = m+3 \\ 
 21\cdot2^{2k-2m+3i-10} + 21\cdot2^{k+3i-7} + 525\cdot2^{k-2m+4i-12} +315\cdot2^{4i-8} - 3255\cdot2^{5i-2m-13}& \text{if    } m+4\leq i\leq 2m+2 \\ 
53\cdot2^{2k+4m-1} - 159\cdot2^{k+6m} +53\cdot2^{8m+2} & \text{if    } i = 2m+3 \\
 2^{3k+2m+1} -7\cdot2^{2k+4m+2} +7\cdot2^{k+6m+4} -2^{8m+7}& \text{if    } i = 2m+4
 \end{cases}\]
\end{lem}

\begin{example}We have for $ m=2,\; l=1, \;k\geq 5 $\\
 
\[\Gamma_{i}^{\left[1\atop{ 1+2\atop 1+2+1} \right]\times k}  =
\begin{cases}
1 & \text{if    } i = 0\\
2^{k} +17 & \text{if    } i = 1 \\
 9\cdot2^{k} + 294  & \text{if    } i =2 \\
 55\cdot2^{k+1} +4360 & \text{if    } i = 3 \\
  2^{2k+2} + 1256\cdot2^{k} + 54208 &\text{if    } i = 4 \\
 11\cdot2^{2k+2} + 13728\cdot2^{k} +457216  & \text{if    } i = 5 \\ 
21\cdot2^{2k+4} + 693\cdot2^{k+8} - 735\cdot2^{13}  & \text{if    } i = 6 \\
  53\cdot2^{2k+7} -159\cdot2^{k+12} +53\cdot2^{18}  & \text{if    } i = 7 \\ 
 2^{3k+5} -7\cdot2^{2k+10} +7\cdot2^{k+16} -2^{23}& \text{if    } i = 8
 \end{cases}\]
\end{example}

\begin{lem}
\label{lem 1.21}
We have for $ m\geq 2 $ and $l=0 $ and  $k \geq 2m $\\
\small
$ \Gamma_{i}^{\left[1\atop{ 1+m\atop 1+m} \right]\times k}  $  is equal to \\
\[ =
\begin{cases}
1 & \text{if    } i = 0\\
2^{k} +17 & \text{if    } i = 1 \\
 (21\cdot2^{3i-7} - 3\cdot2^{2i-5})\cdot2^{k} + 315\cdot2^{4i-8} -21\cdot2^{3i-6}& \text{if    }2\leq  i\leq m \\
 21\cdot2^{k+3m-4} + 45\cdot2^{k+2m-3} +315\cdot2^{4m-4} -213\cdot2^{3m-3} & \text{if    } i = m+1 \\
 3\cdot2^{2k+m} + 21\cdot2^{k+3m-1} + 129\cdot2^{k+2m-1} +315\cdot2^{4m} - 813\cdot2^{3m} &\text{if    } i = m+2 \\
 21\cdot2^{2k-2m+3i-9} + 21\cdot2^{k+3i-7} + 525\cdot2^{k-2m+4i-11} +315\cdot2^{4i-8} - 3255\cdot2^{5i-2m-12}& \text{if    } m+3\leq i\leq 2m+1 \\ 
53\cdot2^{2k+4m-3} - 159\cdot2^{k+6m-3} +53\cdot2^{8m-2} & \text{if    } i = 2m+2 \\
 2^{3k+2m} -7\cdot2^{2k+4m} +7\cdot2^{k+6m+1} -2^{8m+3}& \text{if    } i = 2m+3
 \end{cases}\]
\end{lem}

\begin{example}We have for $ m=2,\; l=0, \;k\geq 4 $\\
 
\[\Gamma_{i}^{\left[1\atop{ 1+2\atop 1+2} \right]\times k}  =
\begin{cases}
1 & \text{if    } i = 0\\
2^{k} +17 & \text{if    } i = 1 \\
 9\cdot2^{k} + 294  & \text{if    } i =2 \\
 87\cdot2^{k+1} +3336 & \text{if    } i = 3 \\
 3\cdot2^{2k+2} + 213\cdot2^{k+3} + 28608 &\text{if    } i = 4 \\
 21\cdot2^{2k+2} + 693\cdot2^{k+5} - 735\cdot2^{9}  & \text{if    } i = 5 \\
 53\cdot2^{2k+5} -159\cdot2^{k+9} +53\cdot2^{14}  & \text{if    } i = 6 \\ 
 2^{3k+4} -7\cdot2^{2k+8} +7\cdot2^{k+13} -2^{19}& \text{if    } i = 7
 \end{cases}\]
\end{example}

\begin{lem}
\label{lem 1.22}
We have for $ m=1 $ and $l\geq 3 $ and  $k \geq 2+l $\\
\small
$ \Gamma_{i}^{\left[1\atop{ 1+1\atop 1+1+l} \right]\times k}  $  is equal to \\
\[ =
\begin{cases}
1 & \text{if    } i = 0\\
2^{k} +17 & \text{if    } i = 1 \\
17\cdot2^{k} +230  & \text{if    } i = 2 \\                   
 2^{2k+1} + 51\cdot2^{k+1} + 3784  & \text{if    } i = 3 \\
 3\cdot2^{2k+2i-7} + 105\cdot2^{k+3i-9} + 945\cdot2^{4i-10} &\text{if    } 4\leq i\leq l+1 \\
 3\cdot2^{2k+2l-3} + 233\cdot2^{k+3l-3} + 433\cdot2^{4l-2} &\text{if    } i = l+2\\
11\cdot2^{2k+2l-1} + 345\cdot2^{k+3l}  -367\cdot2^{4l+2} &\text{if    } i = l+3 \\
53\cdot2^{2k+2l+1} -159\cdot2^{k+3l+3} + 3392\cdot2^{4l} &\text{if    } i = l+4 \\ 
    2^{3k+l+2} -7\cdot2^{2k+2l+4} +7\cdot2^{k+3l+7} -2^{4l+11}& \text{if    } i = l+5
 \end{cases}\]

 \[\Gamma_{i}^{\left[2\atop 5 \right]\times k}  =
\begin{cases}
1 & \text{if    } i = 0\\
9 & \text{if    } i = 1 \\
 2^{k+1} + 62  & \text{if    } i =2 \\
 3\cdot2^{k+1} + 504 & \text{if    } i = 3 \\
 3\cdot2^{k+3} + 4032  &\text{if    } i = 4 \\
 11\cdot2^{k+5} + 15872 & \text{if    } i = 5 \\
  21\cdot2^{k+7} -86016  & \text{if    } i = 6 \\ 
 2^{2k+5} -3\cdot2^{k+10} +2^{16}& \text{if    } i = 7
 \end{cases}\]

\[\Gamma_{i}^{\left[1\atop{ 1+1\atop 1+1+3} \right]\times k}  =
\begin{cases}
1 & \text{if    } i = 0\\
2^{k} +17 & \text{if    } i = 1 \\
 17\cdot2^{k} + 230  & \text{if    } i =2 \\
 2^{2k+1} + 51\cdot2^{k+1} +3784 & \text{if    } i = 3 \\
 3\cdot2^{2k+1} + 105\cdot2^{k+3} + 60480  &\text{if    } i = 4 \\
 3\cdot2^{2k+3} + 233\cdot2^{k+6} +433\cdot2^{10}  & \text{if    } i = 5 \\
  11\cdot2^{2k+5} +345\cdot2^{k+9} -367\cdot2^{14}  & \text{if    } i = 6 \\ 
   53\cdot2^{2k+7} -651264\cdot2^{k} +13893632  & \text{if    } i = 7 \\
 2^{3k+5} -7\cdot2^{2k+10} +7\cdot2^{k+16} -2^{23}& \text{if    } i = 8
 \end{cases}\]

\[\Gamma_{i}^{\left[1\atop{ 1+1\atop 1+1+2} \right]\times k}  =
\begin{cases}
1 & \text{if    } i = 0\\
2^{k} +17 & \text{if    } i = 1 \\
 17\cdot2^{k} + 230  & \text{if    } i =2 \\
 2^{2k+1} + 51\cdot2^{k+1} +3784 & \text{if    } i = 3 \\
3\cdot2^{2k+1} + 233\cdot2^{k+3} +433\cdot2^{6}  & \text{if    } i = 4 \\
 11\cdot2^{2k+3} +345\cdot2^{k+6} -367\cdot2^{10}  & \text{if    } i = 5 \\
 53\cdot2^{2k+5} - 159\cdot2^{k+9} +3392\cdot2^{8} & \text{if    } i = 6 \\
 2^{3k+4} -7\cdot2^{2k+8} +7\cdot2^{k+13} -2^{19}& \text{if    } i = 7
 \end{cases}\]

\[\Gamma_{i}^{\left[1\atop{ 1+1\atop 1+1+1} \right]\times k}  =
\begin{cases}
1 & \text{if    } i = 0\\
2^{k} +17 & \text{if    } i = 1 \\
 17\cdot2^{k} + 230  & \text{if    } i =2 \\
 2^{2k+1} + 115\cdot2^{k+1} +1736 & \text{if    } i = 3 \\
 11\cdot2^{2k+1} +345\cdot2^{k+3} -367\cdot2^{6}  & \text{if    } i = 4 \\
 53\cdot2^{2k+3} - 159\cdot2^{k+6} +3392\cdot2^{4} & \text{if    } i = 5 \\
 2^{3k+3} -7\cdot2^{2k+6} +7\cdot2^{k+10} -2^{15}& \text{if    } i = 6
 \end{cases}\]

\[\Gamma_{i}^{\left[1\atop{ 1+1\atop 1+1} \right]\times k}  =
\begin{cases}
1 & \text{if    } i = 0\\
2^{k} +17 & \text{if    } i = 1 \\
 33\cdot2^{k} + 102  & \text{if    } i =2 \\
3\cdot2^{2k+1} + 171\cdot2^{k+1} - 1464 & \text{if    } i = 3 \\
53\cdot2^{2k+1} - 159\cdot2^{k+3} +3392 & \text{if    } i = 4 \\
 2^{3k+2} -7\cdot2^{2k+4} +7\cdot2^{k+7} -2^{11}& \text{if    } i = 5
 \end{cases}\]
 \end{lem}
    \begin{lem} 
     \label{lem 1.23}
\underline {The case  m = 0 , $ k\geq 3 $}\\
\begin{equation*}
\Gamma _{i}^{\left[\stackrel{1}{\stackrel{1}{1}}\right]\times k} = \begin{cases}
1 & \text{if  } i = 0 \\
 7(2^k-1) &  \text{if  }    i=1 \\
7\cdot2^{2k} - 21\cdot2^k + 14  & \text{if   } i = 2 \\
2^{3k} - 7\cdot2^{2k}  + 7\cdot2^{k+1} - 8  & \text{if   } i = 3
\end{cases}
\end{equation*}
\underline {The case  m = 1 , $ k\geq 4 $}\\
\begin{equation*}
   \Gamma _{i}^{\left[\stackrel{1}{\stackrel{1}{1+1}}\right]\times k}    = \begin{cases}
 1  & \text{if  } i = 0 \\
 3\cdot2^k + 9   &  \text{if  }    i=1 \\
2^{2k}+47\cdot2^k -54     & \text{if   } i = 2 \\
27\cdot2^{2k}-151\cdot2^k +172 & \text{if   } i = 3 \\
2^{3k+1} -7\cdot2^{2k+2} + 7\cdot2^{k+4} -128 & \text{if   } i = 4
\end{cases}
\end{equation*}
\underline {The case  m = 2 , $ k\geq 5 $}\\
\begin{equation*}
 \Gamma _{i}^{\left[\stackrel{1}{\stackrel{1}{1+2}}\right]\times k}  = \begin{cases}
 1  & \text{if  } i = 0 \\
3\cdot2^k + 9  &  \text{if  }    i=1 \\
2^{2k} + 15\cdot2^{k} +158  & \text{if   } i = 2 \\
3\cdot2^{2k}+ 191\cdot2^{k+1} - 1576   & \text{if   } i = 3 \\
27\cdot2^{2k} -81\cdot2^{k+4} +3456   & \text{if   } i = 4 \\
2^{3k+2} - 7\cdot2^{2k+4}  +7\cdot2^{k+7} - 2048  & \text{if   } i = 5
\end{cases}
\end{equation*}
\underline {The case $4\leq  k \leq 2+m $}\\
\begin{equation*}
 \Gamma _{i}^{\left[\stackrel{1}{\stackrel{1}{1+m}}\right]\times k}= \begin{cases}
1 & \text{if  } i = 0 \\
3\cdot2^k + 9  &  \text{if  }    i=1 \\
2^{2k}+15\cdot2^k +158  &  \text{if  }    i= 2  \\
3\cdot2^{2k+ 2i -6 } + 63\cdot[2^{k +3i -8} +5\cdot2^{4i-9}] & \text{if   } 3\leq i\leq k-1 \\
2^{3k +m} - 47\cdot2^{4k -9} & \text{if   } i = k
\end{cases}
\end{equation*}
 \underline {The case $ 3\leq m\leq k-3 $}
\begin{equation*}
 \Gamma _{i}^{\left[\stackrel{1}{\stackrel{1}{1+m}}\right]\times k}= \begin{cases}
1 & \text{if  } i = 0 \\
3\cdot2^k + 9   &  \text{if  }    i=1 \\
2^{2k}+15\cdot2^k +158  &  \text{if  }    i= 2  \\
3\cdot2^{2k+ 2i -6 } + 63\cdot[2^{k +3i -8} +5\cdot2^{4i-9}] & \text{if   } 3\leq i\leq m \\
 3\cdot2^{2k+ 2m - 4 }  + 191\cdot 2^{k +3m - 5}  - 197\cdot2^{4m-5}    & \text{if   } i = m +1  \\
  27\cdot2^{2k+ 2m - 2}  - 81\cdot 2^{k +3m - 2}  + 27\cdot2^{4m-1}     & \text{if   } i = m +2 \\
 2^{3k+m} -7\cdot2^{2k+2m} + 7\cdot2^{k+3m+1} - 2^{4m +3}    & \text{if   } i = m +3
\end{cases}
\end{equation*}
\end{lem}

\begin{proof}
The proofs of   Lemma \ref{lem 1.19}, Lemma \ref{lem 1.20}. Lemma \ref{lem 1.21},  Lemma \ref{lem 1.22} and Lemma \ref{lem 1.23} are similar to the proof of Lemma \ref{lem 1.18}.
\end{proof}

 \begin{lem}
\label{lem 1.24}Let $ \Gamma_{i}^{\left[ (n)\atop{ 1+m\atop 1+m+l} \right]\times k}  $ denote the number of matrices of the form
 $\left[{A\over{B \over C}}\right] $of rank i such that  A is a  $ n\times k $ matrix  over $ \mathbb{F}_{2},$ B a 
$(1+m)\times k$ persymmetric matrix and C  a $ (1+m+l)\times k $ persymmetric  matrix  
 and where $\Gamma_{i}^{\left[ 1+m\atop 1+m+l \right]\times k}$ denotes
   the number of double persymmetric $ (2+2m+l)\times k $  rank i matrices of the form  $\left[{B \over C}\right] $\\
     Then  $\Gamma_{i}^{\left[ (n)\atop{ 1+m\atop 1+m+l} \right]\times k} $ expressed as a linear combination of the $\Gamma_{i-j}^{\left[ 1+m\atop 1+m+l \right]\times k}$ is equal to
   \begin{equation}
   \label{eq 1.5}
 \sum_{j= 0}^{n}2^{(n-j)\cdot(i-j)} a_{j}^{(n)}\prod_{l=1}^{j}(2^{k}- 2^{i-l})\cdot\Gamma_{i-j}^{\left[ 1+m\atop 1+m+l \right]\times k}
\quad for \quad 0\leq i\leq inf(k,n+2m+l+2),
 \end{equation}
  where $ a_{j}^{(n)} $ satisfies the linear recurrence relation 
 \begin{equation}
 \label{eq 1.6}
 a_{j}^{(n)} = 2^{j}\cdot a_{j}^{(n-1)} + a_{j-1}^{(n-1)},\quad n = 2,3,4,\ldots       \quad for\quad 1\leq j\leq n-1. 
\end{equation}
 The explicit value of $ a_{j}^{(n)}$ is given by the following formula\\
  \begin{equation}
  \label{eq 1.6}
 a_{j}^{(n)} =    (-1)^{j}\cdot2^{jn - {j(j-1)\over 2}} + \sum_{s =0}^{j-1} (-1)^{s}\prod_{l=0}^{j-(s+1)}{2^{n+1}- 2^{l}\over 2^{j-s}-2^{l}}\cdot2^{s(n-j) +{s(s+1)\over 2}}
  \quad for \quad  1\leq j\leq n-1.
   \end{equation}
   We set  
   \begin{align*}
   a_{0}^{(n)} & =   a_{n}^{(n)} = 1    \\ and\quad
  \Gamma_{i-j}^{\left[ 1+m\atop 1+m+l \right]\times k}  & = 0 \quad if \quad  i-j \notin \{0,1,2,\ldots, inf(k,2m+l+2)\}.
   \end{align*}
  \end{lem}
  
  \begin{proof}
  The proof of Lemma \ref{lem 1.23} is somewhat similar to the proof of Lemma 11.1 in [1], indeed 
  a analogous proof  of  Lemma \ref{lem 1.15} gives the following generalization \\
  
\begin{equation}
\label{eq 1.7}
 \Gamma_{i}^{\left[ (n)\atop{ 1+m\atop 1+m+l} \right]\times k}  = 2^{i}\cdot \Gamma_{i}^{\left[ (n-1)\atop{ 1+m\atop 1+m+l} \right]\times k} 
    +      (2^{k}-2^{i-1})\cdot  \Gamma_{i-1}^{\left[ (n-1)\atop{ 1+m\atop 1+m+l} \right]\times k}\quad for \quad 0\leq i\leq inf(k,n+2m+l+2).
    \end{equation}
We prove \eqref{eq 1.5} by induction on n. \vspace{0.5 cm}\\
 From \eqref{eq 1.7} we get
respectively for n =1 and n = 2 \\
\begin{align}
   \Gamma_{i}^{\left[1\atop{ 1+m\atop 1+m+l} \right]\times k} & = ( 2^{k} -  2^{i-1})\cdot \Gamma _{i-1}^{\left[\stackrel{1+m}{1+m+l}\right]\times k}
 +  2^{i}\cdot\Gamma _{i}^{\left[\stackrel{1+m}{1+m+l}\right]\times k}  \quad  for \quad  0\leq i\leq inf(k,2m+l+3) \label{eq 1.8} \\
 \Gamma_{i}^{\left[(2)\atop{ 1+m\atop 1+m+l} \right]\times k} & = ( 2^{k} -  2^{i-1})\cdot  \Gamma_{i-1}^{\left[1\atop{ 1+m\atop 1+m+l} \right]\times k}
 +  2^{i}\cdot  \Gamma_{i}^{\left[1\atop{ 1+m\atop 1+m+l} \right]\times k}   \quad  for \quad  0\leq i\leq inf(k,2m+l+4) \label{eq 1.9} 
 \end{align}
From \eqref{eq 1.8} and \eqref{eq 1.9} we deduce \\
\begin{align}
 \Gamma_{i}^{\left[(2)\atop{ 1+m\atop 1+m+l} \right]\times k} 
& = 2^{2i}\Gamma _{i}^{\left[\stackrel{1+m}{1+m+l}\right]\times k}+
     3\cdot2^{i-1}(2^{k}-2^{i-1})\cdot\Gamma _{i-1}^{\left[\stackrel{1+m}{1+m+l}\right]\times k} \label{eq 1.10}\\
       &  +(2^{k}-2^{i-1})(2^{k}-2^{i-2})\cdot\Gamma _{i-2}^{\left[\stackrel{1+m}{1+m+l}\right]\times k}
      \quad for\quad 0\leq i\leq inf(k,2m+l+4). \nonumber
\end{align}

Hence by \eqref{eq 1.8} and \eqref{eq 1.10} the formula \eqref{eq 1.5} holds for n = 1 and n = 2. \vspace{0.5 cm}\\
Assume now that \eqref{eq 1.5} holds for the number n-1, that is
 \begin{align}    
   2^{i}\cdot \Gamma_{i}^{\left[ (n-1)\atop{ 1+m\atop 1+m+l} \right]\times k}   =    
 \sum_{j= 0}^{n-1}2^{(n-1-j)\cdot(i-j)}2^{i} a_{j}^{(n-1)}\prod_{l=1}^{j}(2^{k}- 2^{i-l})\cdot
\Gamma _{i-j}^{\left[\stackrel{1+m}{1+m+l}\right]\times k}\label{eq 1.11}
 \end{align}
 
\begin{align}
 (2^{k}-2^{i-1})\cdot \Gamma_{i-1}^{\left[ (n-1)\atop{ 1+m\atop 1+m+l} \right]\times k} & =
  (2^{k}-2^{i-1})\big[ \sum_{j= 0}^{n-1}2^{(n-1-j)\cdot(i-1-j)} a_{j}^{(n-1)}\prod_{l=1}^{j}(2^{k}- 2^{i-1-l})\cdot
\Gamma _{i-1-j}^{\left[\stackrel{1+m}{1+m+l}\right]\times k} \big] \label{eq 1.12} \\ & =
  \sum_{j= 0}^{n-1}2^{(n-1-j)\cdot(i-1-j)} a_{j}^{(n-1)}\prod_{l=1}^{j+1}(2^{k}- 2^{i-l})\cdot
 \Gamma _{i-(j+1)}^{\left[\stackrel{1+m}{1+m+l}\right]\times k}   \nonumber \\
 & = \sum_{j= 1}^{n}2^{(n-j)\cdot(i-j)} a_{j-1}^{(n-1)}\prod_{l=1}^{j}(2^{k}- 2^{i-l})\cdot\Gamma _{i- j}^{\left[\stackrel{1+m}{1+m+l}\right]\times k} . \nonumber
\end{align}

From \eqref{eq 1.11} ,  \eqref{eq 1.12} and \eqref{eq 1.6} it follows
 \begin{multline}
 \Gamma_{i}^{\left[ (n)\atop{ 1+m\atop 1+m+l} \right]\times k}    = 2^{ni} \Gamma _{i}^{\left[\stackrel{1+m}{1+m+l}\right]\times k}
+ \sum_{j= 1}^{n-1}2^{(n-j)\cdot(i-j)}(2^{j}a_{j}^{(n-1)}  +
 a_{j-1}^{(n-1)})\prod_{l=1}^{j}(2^{k}- 2^{i-l})\cdot\Gamma _{i-j}^{\left[\stackrel{1+m}{1+m+l}\right]\times k} \nonumber \\
 + a_{n-1}^{(n-1)}\prod_{l=1}^{n}(2^{k}- 2^{i-l}) \Gamma _{i- n}^{\left[\stackrel{1+m}{1+m+l}\right]\times k}   \hspace{6 cm}\\
  =  \sum_{j= 0}^{n}2^{(n-j)\cdot(i-j)} a_{j}^{(n)}\prod_{l=1}^{j}(2^{k}- 2^{i-l})\cdot
\Gamma _{i-j}^{\left[\stackrel{1+m}{1+m+l}\right]\times k} \hspace{6 cm}
\end{multline}
which proves \eqref{eq 1.5}.\\
By formula \eqref{eq 1.5} with $ m = l = 0,$ we obtain\\
\begin{equation}
\label{eq 1.13}
 \Gamma_{i}^{\left[ (n)\atop{ 1\atop 1} \right]\times k} =  \sum_{j= 0}^{n}2^{(n-j)\cdot(i-j)} a_{j}^{(n)}\prod_{l=1}^{j}(2^{k}- 2^{i-l})\cdot
\Gamma _{i-j}^{\left[\stackrel{1}{1}\right]\times k} \quad for \quad 0\leq i\leq inf(k,n+2)
\end{equation}
We have obviously \begin{equation}
\label{eq 1.14}
\Gamma _{i-j}^{\left[\stackrel{1}{1}\right]\times k}   = \begin{cases}
 1 & \text{if      }   i-j = 0, \\
 3\cdot(2^{k} -1) & \text{if      }   i-j = 1,\\
 2^{2k} -3\cdot2^{k} +2 & \text{if      }   i-j = 2.
    \end{cases}
\end{equation}
From \eqref{eq 1.13} and  \eqref{eq 1.14} we deduce
\begin{align}
 \Gamma_{i}^{\left[ (n)\atop{ 1\atop 1} \right]\times k}& =  \sum_{j= i-2}^{i}2^{(n-j)\cdot(i-j)} a_{j}^{(n)}\prod_{l=1}^{j}(2^{k}- 2^{i-l})\cdot
\Gamma _{i-j}^{\left[\stackrel{1}{1}\right]\times k} \label{eq 1.15}\\
& = a_{i-2}^{(n)}\cdot2^{2n-2i+4}\prod_{l=1}^{i-2}(2^{k}- 2^{i-l})\cdot( 2^{2k} -3\cdot2^{k} +2) 
+  a_{i-1}^{(n)}\cdot2^{n-i+1}\cdot\prod_{l=1}^{i-1}(2^{k}- 2^{i-l})\cdot3\cdot(2^{k} -1) \nonumber \\
& +  a_{i}^{(n)}\cdot\prod_{l=1}^{i}(2^{k}- 2^{i-l}) \nonumber  \\
 & = [2^{2n-2i+4} \cdot a_{i-2}^{(n)}  +  3\cdot2^{n-i+1}\cdot a_{i-1}^{(n)}  +   a_{i}^{(n)}] \cdot \prod_{l=1}^{i}(2^{k}- 2^{i-l}) \nonumber \\
   & = [2^{2n-2i+4} \cdot a_{i-2}^{(n)}  +  3\cdot2^{n-i+1}\cdot a_{i-1}^{(n)}  +   a_{i}^{(n)}] \cdot \prod_{l=0}^{i-1}(2^{k}- 2^{l}) \nonumber                                                                                                        
\end{align}
On the other hand by George Landsberg [3] we have
\begin{equation}
\label{eq 1.16}
 \Gamma_{i}^{\left[ (n)\atop{ 1\atop 1} \right]\times k} =
 \prod_{l = 0}^{i-1}{ (2^{n+2} -2^{l})(2^{k}-2^{l}) \over (2^{i}- 2^{l})}. 
\end{equation} 

 Hence by \eqref{eq 1.15} and \eqref{eq 1.16} we have the formula
 \begin{equation}
 \label{eq 1.17}
  2^{2n-2i+4} \cdot a_{i-2}^{(n)}  +  3\cdot2^{n-i+1}\cdot a_{i-1}^{(n)}  +   a_{i}^{(n)}  =  \prod_{l = 0}^{i-1}{ 2^{n+2} -2^{l} \over 2^{i}- 2^{l}}. 
 \end{equation}
 We can then from \eqref{eq 1.17} compute successively $ a_{i}^{(n)}\quad i = 1,2,\ldots\quad for\; n\geq i $\\
  We get $$ a_{1}^{(n)}  = 2^n -1\quad for\quad n\geq 1$$
 $$ a_{2}^{(n)} = {2^{2n-1}-3\cdot2^{n-1} + 1\over 3} \quad  for \quad n\geq 2 $$
To establish that  $ a_{j}^{(n)}$ is given by the  formula \eqref{eq 1.7}
we compute the sum  $$ 2^{2n-2i+4} \cdot a_{i-2}^{(n)}  +  3\cdot2^{n-i+1}\cdot a_{i-1}^{(n)}  +   a_{i}^{(n)}.$$
We get: \\
\begin{align*}
& 2^{2n-2i+4} \cdot\left[ (-1)^{i-2}\cdot2^{(i-2)n - {(i-2)(i-3)\over 2}} + \sum_{s =0}^{i-3} (-1)^{s}\prod_{l=0}^{i-s-3}{2^{n+1}- 2^{l}\over 2^{i-s-2}-2^{l}}\cdot2^{s(n-i+2) +{s(s+1)\over 2}} \right]\\
& + 3\cdot2^{n-i+1}\cdot\left[ (-1)^{i-1}\cdot2^{(i-1)n - {(i-1)(i-2)\over 2}} + \sum_{s =0}^{i-2} (-1)^{s}\prod_{l=0}^{i-s-2}{2^{n+1}- 2^{l}\over 2^{i-s-1}-2^{l}}\cdot2^{s(n-i+1) +{s(s+1)\over 2}} \right]\\
& +   (-1)^{i}\cdot2^{in - {i(i-1)\over 2}} + \sum_{s =0}^{i-1} (-1)^{s}\prod_{l=0}^{i-(s+1)}{2^{n+1}- 2^{l}\over 2^{i-s}-2^{l}}\cdot2^{s(n-i) +{s(s+1)\over 2}}
 =  \prod_{l = 0}^{i-1}{ 2^{n+2} -2^{l} \over 2^{i}- 2^{l}}. 
\end{align*}
\end{proof}

\begin{lem}
\label{lem 1.25}

 Let $ f_{m,l,k}(t,\eta,\xi  _{1}, \xi _{2},\ldots, \xi _{n} ) $  be the exponential sum  in $ \mathbb{P}^{n+2} $ defined by\\
 \small
$ (t,\eta,\xi  _{1}, \xi _{2},\ldots, \xi _{n} ) \in \mathbb{P}^{n+2}\longrightarrow \\
  \sum_{deg Y\leq k-1}\sum_{deg Z\leq1+ m}E(tYZ)\sum_{deg U\leq1+ m+l}E(\eta YU)\sum_{deg V_{1}\leq  0}E(\xi _{1} YV_{1})
  \sum_{deg V_{2} \leq 0}E(\xi _{2} YV_{2}) \ldots \sum_{deg V_{n} \leq 0} E(\xi  _{n} YV_{n}). $\vspace{0.5 cm}\\
 \normalsize 
  Set $$(t,\eta,\xi  _{1}, \xi _{2},\ldots, \xi _{n} ) =
  \big(\sum_{i\geq 1}\alpha _{i}T^{-i}, \sum_{i\geq 1}\beta  _{i}T^{-i},\sum_{i\geq 1}\gamma  _{1i}T^{-i}), \ldots, \sum_{i\geq 1}\gamma  _{ni}T^{-i}) \in\mathbb{P}^{n+2}.   $$     
        
   Then
  $$ f_{m,l,k}(t,\eta,\xi _{1}, \xi _{2},\ldots, \xi _{n} ) = 
  2^{k+2m +l+n+4-r(D^{\left[(n)\atop{ 1+m\atop 1+m+l} \right]\times k}(\xi_{1},\xi _{2},\ldots,\xi _{n}, t,\eta  ))} $$
 where

$$D^{\left[(n)\atop{ 1+m\atop 1+m+l} \right]\times k}(\xi_{1},\xi _{2},\ldots,\xi _{n}, t,\eta  ) $$
 denotes  the following  $(2m+l+n+2) \times k $ matrix
   $$  \left ( \begin{array} {cccccc}
\gamma _{11} & \gamma   _{12} & \gamma  _{13} & \ldots  & \gamma  _{1,k-1} &  \gamma  _{1,k}  \\
 \gamma  _{21} & \gamma  _{22} & \gamma  _{23} & \ldots  & \gamma  _{2,k-1} &   \gamma _{2,k}\\
\vdots & \vdots & \vdots   &  \ldots  & \vdots  &  \vdots \\
\vdots & \vdots & \vdots    &  \ldots & \vdots  &  \vdots \\
 \gamma _{n1} & \gamma  _{n2} &  \gamma  _{n3} & \ldots  & \gamma  _{n,k-1} &  \gamma  _{n,k}\\
 \hline
 \alpha _{1} & \alpha _{2} & \alpha _{3} &  \ldots & \alpha _{k-1}  &  \alpha _{k} \\
\alpha _{2 } & \alpha _{3} & \alpha _{4}&  \ldots  &  \alpha _{k} &  \alpha _{k+1} \\
\vdots & \vdots & \vdots   &  \ldots  & \vdots  &  \vdots \\
\vdots & \vdots & \vdots    &  \ldots & \vdots  &  \vdots \\
\alpha _{1+m} & \alpha _{2+m} & \alpha _{3+m} & \ldots  &  \alpha _{k+m-1} &  \alpha _{k+m}  \\
\hline \\
\beta  _{1} & \beta  _{2} & \beta  _{3} & \ldots  &  \beta_{k-1} &  \beta _{k}  \\
\beta  _{2} & \beta  _{3} & \beta  _{4} & \ldots  &  \beta_{k} &  \beta _{k+1}  \\
\vdots & \vdots & \vdots    &  \vdots & \vdots  &  \vdots \\
\beta  _{1+m+l} & \beta  _{2+m+l} & \beta  _{3+m+l} & \ldots  &  \beta_{k+m+l-1} &  \beta _{k+m+l}  
  \end{array}  \right) $$
  
Then the number of solutions \\
 $(Y_1,Z_1,U_{1},V_{1}^{(1)},V_{2}^{(1)}, \ldots,V_{n}^{(1)}, Y_2,Z_2,U_{2},V_{1}^{(2)},V_{2}^{(2)}, 
\ldots,V_{n}^{(2)},\ldots  Y_q,Z_q,U_{q},V_{1}^{(q)},V_{2}^{(q)}, \ldots,V_{n}^{(q)}   ) $ \vspace{0.5 cm}\\
 of the polynomial equations  \vspace{0.5 cm}
  \[\left\{\begin{array}{c}
 Y_{1}Z_{1} +Y_{2}Z_{2}+ \ldots + Y_{q}Z_{q} = 0  \\
  Y_{1}U_{1} +Y_{2}U_{2}+ \ldots + Y_{q}U_{q} = 0 \\
   Y_{1}V_{1}^{(1)} + Y_{2}V_{1}^{(2)} + \ldots  + Y_{q}V_{1}^{(q)} = 0  \\
    Y_{1}V_{2}^{(1)} + Y_{2}V_{2}^{(2)} + \ldots  + Y_{q}V_{2}^{(q)} = 0\\
    \vdots \\
   Y_{1}V_{n}^{(1)} + Y_{2}V_{n}^{(2)} + \ldots  + Y_{q}V_{n}^{(q)} = 0 
 \end{array}\right.\]
    satisfying the degree conditions \\
                   $$  \deg Y_i \leq k-1 , \quad \deg Z_i \leq m +1 ,\quad \deg U_i \leq m +l+1 ,
                   \quad \deg V_{j}^{i} \leq 0 , \quad  for \quad 1\leq j\leq n  \quad 1\leq i \leq q $$ \\
                   
  is equal to the following integral over the unit interval in $ \mathbb{K}^{n+2} $
 $$ \int_{\mathbb{P}^{n+2}} f_{m,l,k}^{q}(t,\eta ,\xi _{1},\xi _{2},\ldots,\xi _{n} )dt d \eta d \xi _{1}d \xi _{2}\ldots d\xi_{n}. $$
  Observing that $ f_{m,l,k}^{q}(t,\eta ,\xi _{1},\xi _{2},\ldots,\xi _{n} )$ is constant on cosets of $ \mathbb{P}_{k+m}\times \mathbb{P}_{k+m+l}\times\mathbb{P}_{k}^{n}, $\;
  the above integral is equal to 
$$  2^{q(k+2m+l+n+4) - (2m+l+k(n+2))}\sum_{i = 0}^{inf(k,n+2m+l+2)}
  \Gamma_{i}^{\left[ (n)\atop{ 1+m\atop 1+m+l} \right]\times k} 2^{-iq}.$$
\end{lem}
\begin{proof}
 The proof of Lemma \ref{lem 1.24} is  just  a generalization of the proof of  Corollary 11.5. (see [1]) 
\end{proof}

\begin{example} n=2, m=1, $ l=3.$\;$ k\geq 5, $\; q = 3\vspace{0.1 cm}\\

\begin{align*}
& \Gamma_{i}^{\left[ (2)\atop{ 2\atop 5} \right]\times k}\\
& =2^{2i}\Gamma_{i}^{\left[ 2\atop 5 \right]\times k}+                                                                             
     3\cdot2^{i-1}(2^{k}-2^{i-1})\cdot\Gamma_{i-1}^{\left[ 2\atop 5 \right]\times k} 
     +(2^{k}-2^{i-1})(2^{k}-2^{i-2})\cdot \Gamma_{i-2}^{\left[ 2\atop 5 \right]\times k}
       \quad for\quad 0\leq i\leq \inf(k,9) 
\end{align*}

\[\Gamma_{i}^{\left[(2)\atop{ 1+1\atop 1+1+3} \right]\times k}  =
\begin{cases}
1 & \text{if    } i = 0\\
3\cdot2^{k} + 33 & \text{if    } i = 1 \\
 2^{2k} + 83\cdot2^{k} +886 & \text{if    } i = 2 \\
 33\cdot2^{2k} + 978\cdot2^{k} + 29352  &\text{if    } i = 3 \\
  2^{3k+1} + 182\cdot2^{2k} +16408\cdot2^{k} + 937408 & \text{if    } i = 4 \\
 3\cdot2^{3k+1} +189\cdot2^{2k+3} +8191\cdot2^{k+6}+12911\cdot2^{10}  & \text{if    } i = 5 \\ 
  3\cdot2^{3k+3} +36672\cdot2^{2k} +11271168\cdot2^{k} -399015936  & \text{if    } i = 6 \\
 11\cdot2^{3k+5} +1022464\cdot2^{2k} - 100679680\cdot2^{k} + 2163212288 & \text{if    } i = 7 \\  
 117\cdot2^{3k+7} - 3354624\cdot2^{2k} +214695936\cdot2^{k} - 3925868544 & \text{if    } i = 8 \\
  2^{4k+5} - 480\cdot2^{3k+5} +2240\cdot2^{2k+10} -1920\cdot2^{k+16} +2^{31} & \text{if    } i = 9
 \end{cases}\]
 If k = 5 we obtain \\
 
  \[\Gamma_{i}^{\left[(2)\atop{ 1+1\atop 1+1+3} \right]\times 5}  =
\begin{cases}
1 & \text{if    } i = 0\\
129 & \text{if    } i = 1 \\
 4566 & \text{if    } i = 2 \\
 94440 &\text{if    } i = 3 \\
  1714368& \text{if    } i = 4 \\
 31740928  & \text{if    } i = 5 \\ 
 \end{cases}\]
Then the number of solutions \\

 $(Y_1,Z_1,U_{1},V_{1}^{(1)},V_{2}^{(1)}, Y_2,Z_2,U_{2},V_{1}^{(2)},V_{2}^{(2)}, 
 Y_3,Z_3,U_{3},V_{1}^{(3)},V_{2}^{(3)}  ) $\\
 
 of the polynomial equations  \vspace{0.5 cm}
  \[\left\{\begin{array}{c}
 Y_{1}Z_{1} +Y_{2}Z_{2}+ Y_{3}Z_{3} = 0  \\
  Y_{1}U_{1} +Y_{2}U_{2}+ Y_{3}U_{3} = 0 \\
   Y_{1}V_{1}^{(1)} + Y_{2}V_{1}^{(2)} + Y_{3}V_{1}^{(3)} = 0  \\
    Y_{1}V_{2}^{(1)} + Y_{2}V_{2}^{(2)} + Y_{3}V_{2}^{(3)} = 0\\
    \end{array}\right.\]
    satisfying the degree conditions 
                   $$  \deg Y_i \leq 4 , \quad \deg Z_i \leq 2 ,\quad \deg U_i \leq 5 ,
                   \quad \deg V_{j}^{i} \leq 0 , \quad  for \quad 1\leq j\leq 2  \quad 1\leq i \leq 3 $$ 
 is equal to 
 $$  \int_{\mathbb{P}^{4}} f_{1,3,5}^{3}(t,\eta ,\xi _{1},\xi _{2} )dt d \eta d \xi _{1}d \xi _{2}= 2^{23}\sum_{i = 0}^{5}
  \Gamma_{i}^{\left[ (2)\atop{ 2\atop 5} \right]\times k} 2^{- 3i}= 13281\cdot2^{20}.$$
  
$$ \text{where for} \;(t,\eta,\xi  _{1}, \xi _{2}, \xi _{2} ) =
  \big(\sum_{i\geq 1}\alpha _{i}T^{-i}, \sum_{i\geq 1}\beta  _{i}T^{-i},\sum_{i\geq 1}\gamma  _{1i}T^{-i}), \sum_{i\geq 1}\gamma  _{2i}T^{-i}) \in\mathbb{P}^{4}.$$   
 \begin{align*}
&  f_{1,3,5}(t,\eta ,\xi _{1},\xi _{2} ) = 
   \sum_{deg Y\leq 4}\sum_{deg Z\leq 2}E(tYZ)\sum_{deg U\leq 5}E(\eta YU)\sum_{deg V_{1}\leq  0}E(\xi _{1} YV_{1})
  \sum_{deg V_{2} \leq 0}E(\xi _{2} YV_{2})\\
  &  =  2^{16-r(D^{\left[(2)\atop{ 2\atop 5} \right]\times k}(\xi_{1},\xi _{2} t,\eta  ))}
  \end{align*}
and \\

   $$ D^{\left[(2)\atop{ 2\atop 5} \right]\times k}(\xi_{1},\xi _{2}, t,\eta  ) = \left ( \begin{array} {cccccc}
\gamma _{11} & \gamma   _{12} & \gamma  _{13}  & \gamma  _{1,4} &  \gamma  _{1,5}  \\
 \gamma  _{21} & \gamma  _{22} & \gamma  _{23} &  \gamma  _{2,4} &   \gamma _{2,5}\\
 \hline
 \alpha _{1} & \alpha _{2} & \alpha _{3} &  \alpha _{4}  &  \alpha _{5} \\
\alpha _{2 } & \alpha _{3} & \alpha _{4} &  \alpha _{5} &  \alpha _{6} \\
\hline \\
\beta  _{1} & \beta  _{2} & \beta  _{3}  &  \beta_{4} &  \beta _{5}  \\
\beta  _{2} & \beta  _{3} & \beta  _{4}  &  \beta_{5} &  \beta _{6}  \\
\beta  _{3} & \beta  _{4} & \beta  _{5}  &  \beta_{6} &  \beta _{7}  \\
\beta  _{4} & \beta  _{5} & \beta  _{6}  &  \beta_{7} &  \beta _{8}  \\
\beta  _{5} & \beta  _{6} & \beta  _{7}  &  \beta_{8} &  \beta _{9}  \\
 \end{array}  \right) $$

\end{example}
   \subsection{Proofs of Theorems \ref{thm 1.4},\ref{thm 1.5},\ref{thm 1.6},\ref{thm 1.7} and \ref{thm 1.8}.}
  \label{subsec 5}
  \subsubsection{Proof of Theorem \ref{thm 1.4}}
\label{subsubsec 5.1}.
   Follows from Lemma \ref{lem 1.15}
 \subsubsection{Proof of Theorem \ref{thm 1.5}}
\label{subsubsec 5.2}.
   Follows from Lemma \ref{lem 1.16}
  \subsubsection{Proof of Theorem \ref{thm 1.6}}
\label{subsubsec 5.3}.
   Follows from Lemma \ref{lem 1.18} 
  \subsubsection{Proof of Theorem \ref{thm 1.7}}
\label{subsubsec 5.4}.
   Follows from Lemma \ref{lem 1.23}  
   \subsubsection{Proof of Theorem \ref{thm 1.8}}
\label{subsubsec 5.5}.
   Follows from Lemma \ref{lem 1.24}

   \section{\textbf{A recurrent formula for the number $   \Gamma_{i}^{\left[s\atop{ s+m\atop s+m+l} \right]\times k}$ of rank i matrices
   of the form  $\left[{A\over{B \over C}}\right] $ such that  A is a  $ s\times k $ persymmetric matrix over $ \mathbb{F}_{2},$ B a 
$(s+m)\times k$ persymmetric matrix and C  a $ (s+m+l)\times k $ persymmetric  matrix } }
  \label{sec 2}
  \subsection{Notation}
  \label{subsec 1}
   \begin{defn}
\label{defn 2.1}We introduce the following definitions in the three - dimensional  $ \mathbb{K} $- vectorspace.\\
\small
\begin{itemize}
\item Let k,s,m and $l$ denote rational integers such that $ k\geq 1,\; s\geq 2 \; and \; m\geq 0, \; l\geq 0 $
\item We denote by $\mathbb{P}/\mathbb{P}_{i}\times \mathbb{P}/\mathbb{P}_{j}\times \mathbb{P}/\mathbb{P}_{r} $
a complete set of coset representatives of $\mathbb{P}_{i}\times\mathbb{P}_{j}\times\mathbb{P}_{r} $
in  $\mathbb{P}\times\mathbb{P}\times\mathbb{P},$ for instance $\mathbb{P}/\mathbb{P}_{s+k-1}\times \mathbb{P}/\mathbb{P}_{s+m +k-1}\times \mathbb{P}/\mathbb{P}_{s+m+l +k-1} $
denotes  a complete set of coset representatives of $\mathbb{P}_{s+k-1}\times\mathbb{P}_{s+m +k-1}\times\mathbb{P}_{s+m+l +k-1} $ in  $\mathbb{P}\times\mathbb{P}\times\mathbb{P}.$
\item $Set \;(t,\eta,\xi  )= (\sum_{i\geq 1}\alpha _{i}T^{-i},\sum_{i\geq 1}\beta  _{i}T^{-i},\sum_{i\geq 1}\gamma  _{i}T^{-i})
\in \mathbb{P}\times\mathbb{P}\times\mathbb{P}. $\\
\item  We denote by  $  D^{\left[s\atop{ s+m\atop s+m+l} \right]\times k}(t,\eta,\xi )  $ any  $(2s+m +l)\times k $   matrix, 
such that  after a rearrangement of the rows, if necessary,  we can  obtain  the following triple persymmetric matrix 
$\left[{D_{s  \times k}(t)\over{ D_{(s+m )\times k}(\eta )\over D_{(s+m +l)\times k}(\xi  )}}\right] $ 
 $$   \left ( \begin{array} {cccccc}
\alpha _{1} & \alpha _{2} & \alpha _{3} &  \ldots & \alpha _{k-1}  &  \alpha _{k} \\
\alpha _{2 } & \alpha _{3} & \alpha _{4}&  \ldots  &  \alpha _{k} &  \alpha _{k+1} \\
\vdots & \vdots & \vdots    &  \vdots & \vdots  &  \vdots \\
\alpha _{s-1} & \alpha _{s} & \alpha _{s +1} & \ldots  &  \alpha _{s+k-3} &  \alpha _{s+k-2}  \\
\alpha _{s} & \alpha _{s+1} & \alpha _{s +2} & \ldots  &  \alpha _{s+k-2} &  \alpha _{s+k-1}  \\
\hline \\
\beta  _{1} & \beta  _{2} & \beta  _{3} & \ldots  &  \beta_{k-1} &  \beta _{k}  \\
\beta  _{2} & \beta  _{3} & \beta  _{4} & \ldots  &  \beta_{k} &  \beta _{k+1}  \\
\vdots & \vdots & \vdots    &  \vdots & \vdots  &  \vdots \\
\beta  _{m+1} & \beta  _{m+2} & \beta  _{m+3} & \ldots  &  \beta_{k+m-1} &  \beta _{k+m}  \\
\vdots & \vdots & \vdots    &  \vdots & \vdots  &  \vdots \\
\beta  _{s+m-1} & \beta  _{s+m} & \beta  _{s+m+1} & \ldots  &  \beta_{s+m+k-3} &  \beta _{s+m+k-2}  \\
\beta  _{s+m} & \beta  _{s+m+1} & \beta  _{s+m+2} & \ldots  &  \beta_{s+m+k-2} &  \beta _{s+m+k-1} \\
\hline
\gamma  _{1} & \gamma   _{2} & \gamma  _{3} & \ldots  &  \gamma _{k-1} &  \gamma  _{k}  \\
\gamma   _{2} & \gamma  _{3} & \gamma   _{4} & \ldots  &  \gamma _{k} &  \gamma  _{k+1}  \\
\vdots & \vdots & \vdots    &  \vdots & \vdots  &  \vdots \\
\gamma  _{m +l +1} & \gamma   _{m+l+2} & \gamma _{m+l+3} & \ldots  &  \gamma _{k+m+l-1} &  \beta _{k+m +l}  \\
\vdots & \vdots & \vdots    &  \vdots & \vdots  &  \vdots \\
\gamma  _{s+m+l -1} & \gamma  _{s+m +l} & \gamma  _{s+m+l+1} & \ldots  &  \gamma _{s+m+l+k-3} &  \gamma  _{s+m+l+k-2}  \\
\gamma  _{s+m+l} & \gamma  _{s+m+l+1} & \gamma  _{s+m+l+2} & \ldots  &  \gamma _{s+m+l+k-2} &  \gamma  _{s+m+l+k-1} 
\end{array}  \right) $$ 
\item We denote by  $  D^{\left[\stackrel{s-1}{\stackrel{s+m -1}{\stackrel{s+m+l-1}{\overline{\alpha  _{s -}+\beta _{s+m-}}}}}\right] \times k } (t,\eta,\xi   ) $
  any  $(3s+2m+l-2)\times k $   matrix, 
such that  after a rearrangement of the rows, if necessary,  we can  obtain a matrix where the first $(3s+2m+l-3)$  rows form  the following triple persymmetric matrix 
$ \left[{D_{s -1 \times k}(t)\over {D_{(s+m-1 )\times k}(\eta )\over D_{(s+m+l-1)\times k}(\xi )}}\right] $ and the last row is equal to\\
$(\alpha  _{s }+\beta _{s+m},\alpha  _{s +1}+\beta _{s+m+1},\ldots,\alpha  _{s +k-1}+\beta _{s+m+k-1}) $\\
 $$   \left ( \begin{array} {cccccc}
\alpha _{1} & \alpha _{2} & \alpha _{3} &  \ldots & \alpha _{k-1}  &  \alpha _{k} \\
\alpha _{2 } & \alpha _{3} & \alpha _{4}&  \ldots  &  \alpha _{k} &  \alpha _{k+1} \\
\vdots & \vdots & \vdots    &  \vdots & \vdots  &  \vdots \\
\alpha _{s-1} & \alpha _{s} & \alpha _{s +1} & \ldots  &  \alpha _{s+k-3} &  \alpha _{s+k-2}  \\
\hline \\
\beta  _{1} & \beta  _{2} & \beta  _{3} & \ldots  &  \beta_{k-1} &  \beta _{k}  \\
\beta  _{2} & \beta  _{3} & \beta  _{4} & \ldots  &  \beta_{k} &  \beta _{k+1}  \\
\vdots & \vdots & \vdots    &  \vdots & \vdots  &  \vdots \\
\beta  _{m+1} & \beta  _{m+2} & \beta  _{m+3} & \ldots  &  \beta_{k+m-1} &  \beta _{k+m}  \\
\vdots & \vdots & \vdots    &  \vdots & \vdots  &  \vdots \\
\beta  _{s+m-1} & \beta  _{s+m} & \beta  _{s+m+1} & \ldots  &  \beta_{s+m+k-3} &  \beta _{s+m+k-2}  \\
\hline \\
\gamma  _{1} & \gamma   _{2} &  \gamma _{3} & \ldots  & \gamma  _{k-1} &  \gamma  _{k}  \\
\gamma  _{2} & \gamma  _{3} & \gamma  _{4} & \ldots  & \gamma  _{k} &  \gamma  _{k+1}  \\
\vdots & \vdots & \vdots    &  \vdots & \vdots  &  \vdots \\
 \gamma  _{m+l+1} &  \gamma _{m+l+2} &  \gamma _{m+l+3} & \ldots  & \gamma _{k+m+l-1} &  \gamma  _{k+m+l}  \\
\vdots & \vdots & \vdots    &  \vdots & \vdots  &  \vdots \\
\gamma  _{s+m+l-1} & \gamma  _{s+m+l} &  \gamma  _{s+m+l+1} & \ldots  & \gamma  _{s+m+l+k-3} &  \gamma  _{s+m+l+k-2}  \\
 \hline\\
  \alpha  _{s }+\beta _{s+m} & \alpha  _{s +1}+\beta _{s+m+1} & \alpha  _{s +2}+\beta _{s+m+2}& \ldots & \alpha  _{s +k-2}+\beta _{s+m+k-2} & \alpha  _{s +k-1}+\beta _{s+m+k-1}
\end{array}  \right). $$ \\
\item We denote by  $  D^{\left[\stackrel{s-1}{\stackrel{s+m -1}{\stackrel{s+m+l-1}{\overline{\alpha  _{s -}+\gamma  _{s+m+l-}}}}}\right] \times k } (t,\eta,\xi   ) $
  any  $(3s+2m+l-2)\times k $   matrix, 
such that  after a rearrangement of the rows, if necessary,  we can  obtain a matrix where the first $(3s+2m+l-3)$  rows form  the following triple persymmetric matrix 
$ \left[{D_{s -1 \times k}(t)\over {D_{(s+m-1 )\times k}(\eta )\over D_{(s+m+l-1)\times k}(\xi )}}\right] $ and the last row is equal to\\
$(\alpha  _{s }+\gamma  _{s+m+l},\alpha  _{s +1}+\gamma  _{s+m+l+1},\ldots,\alpha  _{s +k-1}+\gamma _{s+m+l+k-1}) $\\
 $$   \left ( \begin{array} {cccccc}
\alpha _{1} & \alpha _{2} & \alpha _{3} &  \ldots & \alpha _{k-1}  &  \alpha _{k} \\
\alpha _{2 } & \alpha _{3} & \alpha _{4}&  \ldots  &  \alpha _{k} &  \alpha _{k+1} \\
\vdots & \vdots & \vdots    &  \vdots & \vdots  &  \vdots \\
\alpha _{s-1} & \alpha _{s} & \alpha _{s +1} & \ldots  &  \alpha _{s+k-3} &  \alpha _{s+k-2}  \\
\hline \\
\beta  _{1} & \beta  _{2} & \beta  _{3} & \ldots  &  \beta_{k-1} &  \beta _{k}  \\
\beta  _{2} & \beta  _{3} & \beta  _{4} & \ldots  &  \beta_{k} &  \beta _{k+1}  \\
\vdots & \vdots & \vdots    &  \vdots & \vdots  &  \vdots \\
\beta  _{m+1} & \beta  _{m+2} & \beta  _{m+3} & \ldots  &  \beta_{k+m-1} &  \beta _{k+m}  \\
\vdots & \vdots & \vdots    &  \vdots & \vdots  &  \vdots \\
\beta  _{s+m-1} & \beta  _{s+m} & \beta  _{s+m+1} & \ldots  &  \beta_{s+m+k-3} &  \beta _{s+m+k-2}  \\
\hline \\
\gamma  _{1} & \gamma   _{2} &  \gamma _{3} & \ldots  & \gamma  _{k-1} &  \gamma  _{k}  \\
\gamma  _{2} & \gamma  _{3} & \gamma  _{4} & \ldots  & \gamma  _{k} &  \gamma  _{k+1}  \\
\vdots & \vdots & \vdots    &  \vdots & \vdots  &  \vdots \\
 \gamma  _{m+l+1} &  \gamma _{m+l+2} &  \gamma _{m+l+3} & \ldots  & \gamma _{k+m+l-1} &  \gamma  _{k+m+l}  \\
\vdots & \vdots & \vdots    &  \vdots & \vdots  &  \vdots \\
\gamma  _{s+m+l-1} & \gamma  _{s+m+l} &  \gamma  _{s+m+l+1} & \ldots  & \gamma  _{s+m+l+k-3} &  \gamma  _{s+m+l+k-2}  \\
 \hline\\
  \alpha  _{s }+\gamma  _{s+m+l} & \alpha  _{s +1}+\gamma  _{s+m+l+1} & \alpha  _{s +2}+\gamma  _{s+m+l+2}& \ldots & \alpha  _{s +k-2}+\gamma  _{s+m+l+k-2} & \alpha  _{s +k-1}+\gamma  _{s+m+l+k-1}
\end{array}  \right). $$ \\

We denote by  $  D^{\left[\stackrel{s-1}{\stackrel{s+m -1}{\stackrel{s+m+l-1}{\overline{\alpha  _{s -}+\beta _{s+m-}+\gamma _{s+m+l-}}}}}\right] \times k } (t,\eta,\xi   ) $
  any  $(3s+2m+l-2)\times k $   matrix, 
such that  after a rearrangement of the rows, if necessary,  we can  obtain a matrix where the first $(3s+2m+l-3)$  rows form  the following triple persymmetric matrix 
$ \left[{D_{s -1 \times k}(t)\over {D_{(s+m-1 )\times k}(\eta )\over D_{(s+m+l-1)\times k}(\xi )}}\right] $ and the last row is equal to\\
$(\alpha  _{s }+\beta _{s+m}+\gamma _{s+m+l},\alpha  _{s +1}+\beta  _{s+m+1}+\gamma _{s+m+l+1},\ldots,\alpha  _{s +k-1}+\beta _{s+m+k-1}+\gamma _{s+m+l+k-1}) $\\

\SMALL
 $$   \left ( \begin{array} {cccccc}
\alpha _{1} & \alpha _{2}  &  \ldots & \alpha _{k-1}  &  \alpha _{k} \\
\alpha _{2 } & \alpha _{3} &  \ldots  &  \alpha _{k} &  \alpha _{k+1} \\
\vdots & \vdots & \vdots    & \vdots  &  \vdots \\
\alpha _{s-1} & \alpha _{s} & \ldots  &  \alpha _{s+k-3} &  \alpha _{s+k-2}  \\
\hline \\
\beta  _{1} & \beta  _{2}  & \ldots  &  \beta_{k-1} &  \beta _{k}  \\
\beta  _{2} & \beta  _{3}  & \ldots  &  \beta_{k} &  \beta _{k+1}  \\
\vdots & \vdots    &  \vdots & \vdots  &  \vdots \\
\beta  _{m+1} & \beta  _{m+2}  & \ldots  &  \beta_{k+m-1} &  \beta _{k+m}  \\
\vdots & \vdots    &  \vdots & \vdots  &  \vdots \\
\beta  _{s+m-1} & \beta  _{s+m}  & \ldots  &  \beta_{s+m+k-3} &  \beta _{s+m+k-2}  \\
\hline \\
\gamma  _{1} & \gamma   _{2}  & \ldots  & \gamma  _{k-1} &  \gamma  _{k}  \\
\gamma  _{2} & \gamma  _{3}  & \ldots  & \gamma  _{k} &  \gamma  _{k+1}  \\
\vdots & \vdots    &  \vdots & \vdots  &  \vdots \\
 \gamma  _{m+l+1} &  \gamma _{m+l+2}  & \ldots  & \gamma _{k+m+l-1} &  \gamma  _{k+m+l}  \\
\vdots & \vdots   &  \vdots & \vdots  &  \vdots \\
\gamma  _{s+m+l-1} & \gamma  _{s+m+l}  & \ldots  & \gamma  _{s+m+l+k-3} &  \gamma  _{s+m+l+k-2}  \\
 \hline\\
  \alpha  _{s }+\beta _{s+m}+\gamma _{s+m+l} & \alpha  _{s +1}+\beta _{s+m+1}+\gamma _{s+m+l+1} & \ldots & \alpha  _{s +k-2}+\beta _{s+m+k-2} +\gamma _{s+m+l+k-2}& \alpha  _{s +k-1}+\beta _{s+m+k-1}+\gamma _{s+m+l+k-1}
\end{array}  \right). $$ \\
\normalsize

We denote by  $  D^{\left[s-1\atop{ s+m-1\atop {s+m+l-1 \atop{\overline{\alpha_{s -} \atop{\beta_{s+m -}} }}}}\right]\times k} (t,\eta,\xi   ) $
  any  $(3s+2m+l-1)\times k $   matrix, 
such that  after a rearrangement of the rows, if necessary,  we can  obtain a matrix where the first $(3s+2m+l-3)$  rows form  the following triple persymmetric matrix 
$ \left[{D_{s -1 \times k}(t)\over {D_{(s+m-1 )\times k}(\eta )\over D_{(s+m+l-1)\times k}(\xi )}}\right] $ and the last two rows is equal to\\
the following $2\times k$ matrix 
 $\begin{pmatrix}
\alpha  _{s } & \alpha  _{s +1} & \ldots &\alpha  _{s +k-1} \\
\beta _{s+m} & \beta  _{s+m+1} & \ldots & \beta _{s+m+k-1}
\end{pmatrix},$ that is we get a matrix of the form\\

 $$   \left ( \begin{array} {cccccc}
\alpha _{1} & \alpha _{2}  &  \ldots & \alpha _{k-1}  &  \alpha _{k} \\
\alpha _{2 } & \alpha _{3} &  \ldots  &  \alpha _{k} &  \alpha _{k+1} \\
\vdots & \vdots & \vdots    & \vdots  &  \vdots \\
\alpha _{s-1} & \alpha _{s} & \ldots  &  \alpha _{s+k-3} &  \alpha _{s+k-2}  \\
\hline \\
\beta  _{1} & \beta  _{2}  & \ldots  &  \beta_{k-1} &  \beta _{k}  \\
\beta  _{2} & \beta  _{3}  & \ldots  &  \beta_{k} &  \beta _{k+1}  \\
\vdots & \vdots    &  \vdots & \vdots  &  \vdots \\
\beta  _{m+1} & \beta  _{m+2}  & \ldots  &  \beta_{k+m-1} &  \beta _{k+m}  \\
\vdots & \vdots    &  \vdots & \vdots  &  \vdots \\
\beta  _{s+m-1} & \beta  _{s+m}  & \ldots  &  \beta_{s+m+k-3} &  \beta _{s+m+k-2}  \\
\hline \\
\gamma  _{1} & \gamma   _{2}  & \ldots  & \gamma  _{k-1} &  \gamma  _{k}  \\
\gamma  _{2} & \gamma  _{3}  & \ldots  & \gamma  _{k} &  \gamma  _{k+1}  \\
\vdots & \vdots    &  \vdots & \vdots  &  \vdots \\
 \gamma  _{m+l+1} &  \gamma _{m+l+2}  & \ldots  & \gamma _{k+m+l-1} &  \gamma  _{k+m+l}  \\
\vdots & \vdots   &  \vdots & \vdots  &  \vdots \\
\gamma  _{s+m+l-1} & \gamma  _{s+m+l}  & \ldots  & \gamma  _{s+m+l+k-3} &  \gamma  _{s+m+l+k-2}  \\
 \hline\\
  \alpha  _{s } & \alpha  _{s +1} & \ldots & \alpha  _{s +k-2}& \alpha  _{s +k-1}\\
   \beta _{s+m} & \beta _{s+m+1} & \ldots & \beta _{s+m+k-2} & \beta _{s+m+k-1}
\end{array}  \right). $$ \\ 
We introduce in a similar way the following notations
$$  D^{\left[s-1\atop{ s+m-1\atop {s+m+l-1 \atop{\overline{\beta_{s+m -} \atop{\alpha_{s -}} }}}}\right]\times k} (t,\eta,\xi   ),\quad  
     D^{\left[s-1\atop{ s+m-1\atop {s+m+l-1 \atop{\overline{\alpha_{s -} +\beta_{s+m -}\atop{\gamma _{s+m+l  -}} }}}}\right]\times k} (t,\eta,\xi   ),\quad
   D^{\left[s-1\atop{ s+m-1\atop {s+m+l-1 \atop{\overline{\alpha_{s -} +\beta_{s+m -}\atop{\beta_{s+m -} + \gamma _{s+m+l  -}} }}}}\right]\times k}, \quad
  D^{\left[s-1\atop{ s+m-1\atop {s+m+l-1 \atop{\overline{\alpha_{s -} +\gamma _{s+m+l  -}\atop{\beta_{s+m -}} }}}}\right]\times k},$$
    $$  D^{\left[s-1\atop{ s+m-1\atop {s+m+l-1 \atop{\overline{\alpha_{s -} \atop{\gamma _{s+m+l -}} }}}}\right]\times k} (t,\eta,\xi   ),\quad
     D^{\left[s-1\atop{ s+m-1\atop {s+m+l-1 \atop{\overline{\gamma _{s+m+l  -} \atop{\alpha _{s -}} }}}}\right]\times k} (t,\eta,\xi   )\quad \text{and}\quad
      D^{\left[s-1\atop{ s+m-1\atop {s+m+l-1 \atop{\overline{\alpha_{s -} \atop{\beta_{s+m -}+\gamma _{s+m+l -}} }}}}\right]\times k} (t,\eta,\xi   ) $$\\
  \item 
  We define $$  \sigma _{i,i}^{\left[\stackrel{s-1}{\stackrel{s+m -1}{\stackrel{s+m+l-1}{\overline{\alpha  _{s -}}}}}\right] \times k } $$
 to be the cardinality of the following set 
$$ \left\{(t,\eta,\xi  )\in \mathbb{P}/\mathbb{P}_{k+s-1}\times \mathbb{P}/\mathbb{P}_{k+s+m-2}\times \mathbb{P}/\mathbb{P}_{k+s+m+l-2}
\mid  r( D^{\left[s-1\atop{ s+m-1\atop s+m+l-1} \right]\times k}(t,\eta ,\xi ) )
 = r( D^{\left[s\atop{ s+m-1\atop s+m+l-1} \right]\times k}(t,\eta ,\xi ) ) = i \right\}.$$\\
 \item  
  We define $$  \sigma _{i,i}^{\left[\stackrel{s-1}{\stackrel{s+m -1}{\stackrel{s+m+l-1}{\overline{\alpha  _{s -}+\beta _{s+m-}}}}}\right] \times k } $$
 to be the cardinality of the following set 
$$ \left\{(t,\eta,\xi  )\in \mathbb{P}/\mathbb{P}_{k+s-1}\times \mathbb{P}/\mathbb{P}_{k+s+m-1}\times \mathbb{P}/\mathbb{P}_{k+s+m+l-2}
\mid  r( D^{\left[s-1\atop{ s+m-1\atop s+m+l-1} \right]\times k}(t,\eta ,\xi )) 
 = r( D^{\left[s-1\atop{ s+m-1\atop{ s+m+l-1\atop \alpha _{s-}+\beta _{s+m-}}} \right]\times k}(t,\eta ,\xi ) ) = i \right\}.$$
 \item   
   We define $$  \sigma _{i,i}^{\left[\stackrel{s-1}{\stackrel{s+m -1}{\stackrel{s+m+l-1}{\overline{\alpha  _{s -}+\gamma  _{s+m+l-}}}}}\right] \times k } $$
 to be the cardinality of the following set 
$$ \left\{(t,\eta,\xi  )\in \mathbb{P}/\mathbb{P}_{k+s-1}\times \mathbb{P}/\mathbb{P}_{k+s+m-2}\times \mathbb{P}/\mathbb{P}_{k+s+m+l-1}
\mid  r( D^{\left[s-1\atop{ s+m-1\atop s+m+l-1} \right]\times k}(t,\eta ,\xi ) )
 = r( D^{\left[s-1\atop{ s+m-1\atop{ s+m+l-1\atop \alpha _{s-}+\gamma  _{s+m+l -}}} \right]\times k}(t,\eta ,\xi ) ) = i \right\}.$$
\item   
   We define $$  \sigma _{i,i,i}^{\left[s-1\atop{ s+m-1\atop {s+m+l-1 \atop{\overline{\alpha_{s -} \atop{\beta_{s+m -}} }}}}\right]\times k}=
     \sigma _{i,i,i}^{\left[s-1\atop{ s+m-1\atop {s+m+l-1 \atop{\overline{\beta _{s+m -} \atop{\alpha _{s -}} }}}}\right]\times k} $$
 to be the cardinality of the following set 
  $$\begin{array}{l}\Bigg\{(t,\eta,\xi  )\in \mathbb{P}/\mathbb{P}_{k+s-1}\times \mathbb{P}/\mathbb{P}_{k+s+m-1}\times \mathbb{P}/\mathbb{P}_{k+s+m+l-2}
\mid  r( D^{\left[s-1\atop{ s+m-1\atop s+m+l-1} \right]\times k}(t,\eta ,\xi ) )
 =   r( D^{\left[s\atop{ s+m-1\atop s+m+l-1} \right]\times k}(t,\eta ,\xi ) )\\
   = r( D^{\left[s\atop{ s+m\atop s+m+l-1} \right]\times k}(t,\eta ,\xi ) )= i  \Bigg\}. \end{array} $$\\
 \item 
  We define $$    \sigma _{i,i,i}^{\left[s-1\atop{ s+m-1\atop {s+m+l-1 \atop{\overline{\alpha_{s -} \atop{\beta_{s+m -}+ \gamma _{s+m+l-}} }}}}\right]\times k}   $$
 to be the cardinality of the following set 
  $$\begin{array}{l}\Bigg\{(t,\eta,\xi  )\in \mathbb{P}/\mathbb{P}_{k+s-1}\times \mathbb{P}/\mathbb{P}_{k+s+m-1}\times \mathbb{P}/\mathbb{P}_{k+s+m+l-1}
\mid  r( D^{\left[s-1\atop{ s+m-1\atop s+m+l-1} \right]\times k}(t,\eta ,\xi )) 
 =   r( D^{\left[s\atop{ s+m-1\atop s+m+l-1} \right]\times k}(t,\eta ,\xi )) \\
 = r(D^{\left[\stackrel{s}{\stackrel{s+m -1}{\stackrel{s+m+l-1}{\overline{\beta  _{s+m -}+\gamma  _{s+m+l-}}}}}\right] \times k }(t,\eta ,\xi )) = i   \Bigg\}. \end{array} $$\\
  \item 
  We introduce in a similar way the following definitions \\
  
$ \sigma _{i,i,i}^{\left[s-1\atop{ s+m-1\atop {s+m+l-1 \atop{\overline{\alpha_{s -}+\beta _{s+m-} \atop{ \gamma _{s+m+l-}} }}}}\right]\times k}, \quad
\sigma _{i,i,i}^{\left[s-1\atop{ s+m-1\atop {s+m+l-1 \atop{\overline{\alpha_{s -}+\beta _{s+m-} \atop{\beta _{s+m-} + \gamma _{s+m+l-}} }}}}\right]\times k}, \quad
 \sigma _{i,i,i}^{\left[s-1\atop{ s+m-1\atop {s+m+l-1 \atop{\overline{\alpha_{s -} \atop{ \gamma  _{s+m+l-}} }}}}\right]\times k}, \quad  \text{and} \quad
\sigma _{i,i,i}^{\left[s-1\atop{ s+m-1\atop {s+m+l-1 \atop{\overline{\alpha   _{s -}+\gamma _{s+m+l-} \atop{\beta _{s+m -}} }}}}\right]\times k}. $\\
 \item   Let  $ ( j_{1}, j_{2},   j_{3},j_{4})  \in \mathbb{N}^{4}$ we define
$$ \sigma _{j_{1},j_{2},j_{3},j_{4}}^{\left[s-1\atop{ s+m-1\atop {s+m+l-1 \atop{\overline{\alpha_{s -} \atop{\beta_{s+m -} \atop \gamma_{s+m+l -}} }}}}\right]\times k}
  =  \sigma _{j_{1},j_{2},j_{3},j_{4}}^{\left[\alpha \atop{ \beta \atop \gamma } \right]\times k}$$
 to be the cardinality of the following set 
 $$\begin{array}{l}\Big\{ (t,\eta,\xi  ) \in \mathbb{P}/\mathbb{P}_{k+s -1}\times
           \mathbb{P}/\mathbb{P}_{k+s+m-1}\times \mathbb{P}/\mathbb{P}_{k+s+m+l-1}
\mid r(  D^{\left[s-1\atop{ s+m-1\atop s+m+l-1} \right]\times k}(t,\eta,\xi ) ) = j_{1} \\
 r(  D^{\left[s\atop{ s+m-1\atop s+m+l-1} \right]\times k}(t,\eta,\xi ) )= j_{2},\quad  
 r(  D^{\left[s\atop{ s+m\atop s+m+l-1} \right]\times k}(t,\eta,\xi ) ) = j_{3},  \quad
  r(  D^{\left[s\atop{ s+m\atop s+m+l} \right]\times k}(t,\eta,\xi ) ) = j_{4} \Big\}
    \end{array}$$ \\
    For any permutation of $ \left\{\alpha ,\beta ,\gamma \right\} $ we define in a similar manner
  $$ \sigma _{j_{1},j_{2},j_{3},j_{4}}^{\left[s-1\atop{ s+m-1\atop {s+m+l-1 \atop{\overline{\beta_{s+m -} \atop{\alpha_{s -} \atop \gamma_{s+m+l -}} }}}}\right]\times k}
  =  \sigma _{j_{1},j_{2},j_{3},j_{4}}^{\left[\beta  \atop{ \alpha  \atop \gamma } \right]\times k},$$
   $$ \sigma _{j_{1},j_{2},j_{3},j_{4}}^{\left[s-1\atop{ s+m-1\atop {s+m+l-1 \atop{\overline{\beta_{s+m -} \atop{\gamma _{s +m+l -} \atop \alpha _{s -}} }}}}\right]\times k}
  =  \sigma _{j_{1},j_{2},j_{3},j_{4}}^{\left[\beta  \atop{ \gamma  \atop \alpha  } \right]\times k},$$
    $$ \sigma _{j_{1},j_{2},j_{3},j_{4}}^{\left[s-1\atop{ s+m-1\atop {s+m+l-1 \atop{\overline{\alpha _{s -} \atop{\gamma _{s +m+l -} \atop \beta  _{s+m  -}} }}}}\right]\times k}
  =  \sigma _{j_{1},j_{2},j_{3},j_{4}}^{\left[\alpha  \atop{ \gamma  \atop \beta  } \right]\times k},$$
     $$ \sigma _{j_{1},j_{2},j_{3},j_{4}}^{\left[s-1\atop{ s+m-1\atop {s+m+l-1 \atop{\overline{\gamma  _{s+m+l -} \atop{\alpha  _{s  -} \atop \beta  _{s+m  -}} }}}}\right]\times k}
  =  \sigma _{j_{1},j_{2},j_{3},j_{4}}^{\left[\gamma  \atop{ \alpha  \atop \beta  } \right]\times k},$$
 $$ \sigma _{j_{1},j_{2},j_{3},j_{4}}^{\left[s-1\atop{ s+m-1\atop {s+m+l-1 \atop{\overline{\gamma  _{s+m+l -} \atop{\beta  _{s+m  -} \atop \alpha   _{s -}} }}}}\right]\times k}
  =  \sigma _{j_{1},j_{2},j_{3},j_{4}}^{\left[\gamma  \atop{ \beta   \atop \alpha   } \right]\times k},$$
  \newpage
\item 
\begin{align*}
& \text{The sum}\quad
  \sigma _{i,i,i,i}^{\left[\beta | \gamma  \atop{ \gamma | \beta \atop \alpha | \alpha  } \right]\times k}
+  \sigma _{i-1,i,i,i}^{\left[\beta | \gamma  \atop{ \gamma | \beta \atop \alpha | \alpha  } \right]\times k}
+  \sigma _{i-1,i-1,i,i}^{\left[\beta | \gamma  \atop{ \gamma | \beta \atop \alpha | \alpha  } \right]\times k}
+  \sigma _{i-2,i-1,i,i}^{\left[\beta | \gamma  \atop{ \gamma | \beta \atop \alpha | \alpha  } \right]\times k} \\
& \text{is equal to the sum}\quad
  \sigma _{i,i,i,i}^{\left[\beta  \atop{ \gamma  \atop \alpha } \right]\times k}
+  \sigma _{i-1,i,i,i}^{\left[\beta  \atop{ \gamma \atop \alpha  } \right]\times k}
+  \sigma _{i-1,i-1,i,i}^{\left[\beta  \atop{ \gamma \atop \alpha   } \right]\times k}
+  \sigma _{i-2,i-1,i,i}^{\left[\beta   \atop{ \gamma \atop  \alpha  } \right]\times k}\\
& \text{and is also equal to the sum}\quad
   \sigma _{i,i,i,i}^{\left[ \gamma  \atop{ \beta  \atop \alpha  } \right]\times k}
+  \sigma _{i-1,i,i,i}^{\left[ \gamma  \atop{  \beta \atop  \alpha  } \right]\times k}
+  \sigma _{i-1,i-1,i,i}^{\left[ \gamma  \atop{  \beta \atop  \alpha  } \right]\times k}
+  \sigma _{i-2,i-1,i,i}^{\left[ \gamma  \atop{  \beta \atop  \alpha  } \right]\times k}.\\
& \\
& \text{Analogous sums are defined in a similar way}
\end{align*}
  \item  $$   \Gamma_{i}^{\left[s\atop{ s+m\atop s+m+l} \right]\times k} =   
Card  \left\{(t,\eta,\xi  )\in \mathbb{P}/\mathbb{P}_{k+s-1}\times \mathbb{P}/\mathbb{P}_{k+s+m-1}\times \mathbb{P}/\mathbb{P}_{k+s+m+l-1}
\mid   r( D^{\left[s\atop{ s+m\atop s+m+l} \right]\times k}(t,\eta ,\xi ) ) = i \right\}.$$

  We recall that the rank of a matrix does not change under elementary row operations.\\
\item We denote by $\ker $  D the nullspace of the matrix D and r(D) the rank of the matrix D.
\item To simplify the notations concerning the exponential sums used in the proofs, we introduce the following definitions.\\
 \item  Let $  g_{1}(t,\eta,\xi  ) $ be the quadratic  exponential sum in $\mathbb{P}\times\mathbb{P}\times\mathbb{P}$ defined by
$$ (t,\eta,\xi  ) \in  \mathbb{P}\times \mathbb{P}\times\mathbb{P}\longmapsto  
  \sum_{deg Y\leq k-1}\sum_{deg Z = s-1}E(tYZ)\sum_{deg U \leq s+m-1}E(\eta YU) \sum_{deg V \leq s+m+l-1}E(\xi  YV) \in \mathbb{Z}  $$\\
\item    Let $  h_{1}(t,\eta,\xi  ) $ be the quadratic  exponential sum in $\mathbb{P}\times\mathbb{P}\times\mathbb{P}$ defined by
$$ (t,\eta,\xi  ) \in  \mathbb{P}\times \mathbb{P}\times\mathbb{P}\longmapsto \\ 
 \sum_{deg Y\leq k-1}E(tYT^{s-1}) \sum_{deg Z\leq  s-2}E(tYZ)\sum_{deg U \leq s+m-2}E(\eta YU) \sum_{deg V \leq s+m+l-2}E(\xi YV) \in \mathbb{Z}  $$ \\  
 \item    Let $  h_{2}(t,\eta,\xi  ) $ be the quadratic  exponential sum in $\mathbb{P}\times\mathbb{P}\times\mathbb{P}$ defined by
$$ (t,\eta,\xi  ) \in  \mathbb{P}\times \mathbb{P}\times\mathbb{P}\longmapsto $$ \\
$$ \sum_{deg Y\leq k-1}E(tYT^{s-1})E(\eta YT^{s+m-1}) \sum_{deg Z\leq  s-2}E(tYZ)\sum_{deg U \leq s+m-2}E(\eta YU) \sum_{deg V \leq s+m+l-2}E(\xi YV) \in \mathbb{Z}  $$ \\  
  \item    Let $  h_{3}(t,\eta,\xi  ) $ be the quadratic  exponential sum in $\mathbb{P}\times\mathbb{P}\times\mathbb{P}$ defined by
$$ (t,\eta,\xi  ) \in  \mathbb{P}\times \mathbb{P}\times\mathbb{P}\longmapsto$$ \\ 
$$ \sum_{deg Y\leq k-1}E(tYT^{s-1})E(\xi YT^{s+m+l-1})    \sum_{deg Z\leq  s-2}E(tYZ)\sum_{deg U \leq s+m-2}E(\eta YU) \sum_{deg V \leq s+m+l-2}E(\xi YV) \in \mathbb{Z}  $$ \\  
   \item    Let $  h_{4}(t,\eta,\xi  ) $ be the quadratic  exponential sum in $\mathbb{P}\times\mathbb{P}\times\mathbb{P}$ defined by \\ 
   
    $ (t,\eta,\xi  ) \in  \mathbb{P}\times \mathbb{P}\times\mathbb{P}\longmapsto  $\\
 $$ \sum_{deg Y\leq k-1}E(tYT^{s-1})E(\eta YT^{s+m-1}) E(\xi YT^{s+m+l-1})    \sum_{deg Z\leq  s-2}E(tYZ)\sum_{deg U \leq s+m-2}E(\eta YU) \sum_{deg V \leq s+m+l-2}E(\xi YV) \in \mathbb{Z}  $$ \\  
 
  \item  Let $  \varphi  (t,\eta,\xi  ) $ be the quadratic  exponential sum in $\mathbb{P}\times\mathbb{P}\times\mathbb{P}$ defined by
$$ (t,\eta,\xi  ) \in  \mathbb{P}\times \mathbb{P}\times\mathbb{P}\longmapsto  
  \sum_{deg Y\leq k-1}\sum_{deg Z \leq  s-2}E(tYZ)\sum_{deg U \leq s+m-2}E(\eta YU) \sum_{deg V \leq s+m+l-2}E(\xi YV) \in \mathbb{Z}  $$\\
    \item  Let $  \varphi_{1}  (t,\eta,\xi  ) $ be the quadratic  exponential sum in $\mathbb{P}\times\mathbb{P}\times\mathbb{P}$ defined by
$$ (t,\eta,\xi  ) \in  \mathbb{P}\times \mathbb{P}\times\mathbb{P}\longmapsto  
  \sum_{deg Y\leq k-1}\sum_{deg Z \leq  s-2}E(tYZ)\sum_{deg U \leq s+m-1}E(\eta YU) \sum_{deg V \leq s+m+l-1}E(\xi YV) \in \mathbb{Z}  $$\\
    \item  Let $  \varphi_{2}  (t,\eta,\xi  ) $ be the quadratic  exponential sum in $\mathbb{P}\times\mathbb{P}\times\mathbb{P}$ defined by
$$ (t,\eta,\xi  ) \in  \mathbb{P}\times \mathbb{P}\times\mathbb{P}\longmapsto  
  \sum_{deg Y\leq k-1}\sum_{deg Z \leq  s-2}E(tYZ)\sum_{deg U \leq s+m-2}E(\eta YU) \sum_{deg V \leq s+m+l-1}E(\xi YV) \in \mathbb{Z}  $$\\
    \item  Let $  \varphi_{3}  (t,\eta,\xi  ) $ be the quadratic  exponential sum in $\mathbb{P}\times\mathbb{P}\times\mathbb{P}$ defined by
$$ (t,\eta,\xi  ) \in  \mathbb{P}\times \mathbb{P}\times\mathbb{P}\longmapsto  
  \sum_{deg Y\leq k-1}\sum_{deg Z \leq  s-2}E(tYZ)\sum_{deg U \leq s+m-1}E(\eta YU) \sum_{deg V \leq s+m+l-2}E(\xi YV) \in \mathbb{Z}  $$\\
 \item    Let $  \kappa _{1}(t,\eta,\xi  ) $ be the quadratic  exponential sum in $\mathbb{P}\times\mathbb{P}\times\mathbb{P}$ defined by
$$ (t,\eta,\xi  ) \in  \mathbb{P}\times \mathbb{P}\times\mathbb{P}\longmapsto $$\\ 
$$ \sum_{deg Y\leq k-1}E(tYT^{s-1})E(\eta YT^{s+m-1}) \sum_{deg Z\leq  s-2}E(tYZ)\sum_{deg U \leq s+m-2}E(\eta YU) \sum_{deg V \leq s+m+l-1}E(\xi YV) \in \mathbb{Z}  $$ \\    
 \item    Let $  \kappa _{2}(t,\eta,\xi  ) $ be the quadratic  exponential sum in $\mathbb{P}\times\mathbb{P}\times\mathbb{P}$ defined by
$$ (t,\eta,\xi  ) \in  \mathbb{P}\times \mathbb{P}\times\mathbb{P}\longmapsto  
 \sum_{deg Y\leq k-1}E(tYT^{s-1}) \sum_{deg Z\leq  s-2}E(tYZ)\sum_{deg U \leq s+m-2}E(\eta YU) \sum_{deg V \leq s+m+l-1}E(\xi YV) \in \mathbb{Z}  $$ \\   
   \item    Let $  g _{2}(t,\eta,\xi  ) $ be the quadratic  exponential sum in $\mathbb{P}\times\mathbb{P}\times\mathbb{P}$ defined by
$$ (t,\eta,\xi  ) \in  \mathbb{P}\times \mathbb{P}\times\mathbb{P}\longmapsto  
 \sum_{deg Y\leq k-1}E(tYT^{s-1}) \sum_{deg Z\leq  s-2}E(tYZ)\sum_{deg U \leq s+m-1}E(\eta YU) \sum_{deg V \leq s+m+l-2}E(\xi YV) \in \mathbb{Z}  $$  \\
   \end{itemize}
\end{defn}

   \subsection{Introduction}
  \label{subsec 2}
Adapting the method used in Section 5 of [2]
we establish in this section a recurrent formula for the number of rank i matrices of the form  $\left[{A\over{B \over C}}\right], $ 
where A,B and C are persymmetric. \\[0.01 cm]

   \subsection{Computation of exponential sums in $ \mathbb{K}^3 $}
  \label{subsec 3}
  \small
  \begin{lem}
\label{lem 2.1}
Let $ (t,\eta ,\xi ) \in  \mathbb{P}\times \mathbb{P}\times \mathbb{P}  $ and
\begin{align*}
g_{1}(t,\eta,\xi ) & = \sum_{deg Y\leq k-1}\sum_{deg Z = s-1}E(tYZ)\sum_{deg U \leq s+m-1}E(\eta YU) \sum_{deg V \leq s+m+l-1}E(\xi YV)  \\
 h_{1}(t,\eta,\xi) & = \sum_{deg Y\leq k-1}E(tYT^{s-1}) \sum_{deg Z\leq  s-2}E(tYZ)\sum_{deg U \leq s+m-2}E(\eta YU) \sum_{deg V \leq s+m+l-2}E(\xi YV) \\
  h_{2}(t,\eta,\xi) & = \sum_{deg Y\leq k-1}E(tYT^{s-1})E(\eta YT^{s+m-1}) \sum_{deg Z\leq  s-2}E(tYZ)\sum_{deg U \leq s+m-2}E(\eta YU) \sum_{deg V \leq s+m+l-2}E(\xi YV)   \\
 h_{3}(t,\eta,\xi) & = \sum_{deg Y\leq k-1}E(tYT^{s-1})E(\xi YT^{s+m+l-1}) \sum_{deg Z\leq  s-2}E(tYZ)\sum_{deg U \leq s+m-2}E(\eta YU) \sum_{deg V \leq s+m+l-2}E(\xi YV) \\
 h_{4}(t,\eta,\xi) & = \sum_{deg Y\leq k-1}E(tYT^{s-1})E(\eta YT^{s+m-1}) E(\xi YT^{s+m+l-1})\sum_{deg Z\leq  s-2}E(tYZ)\sum_{deg U \leq s+m-2}E(\eta YU) \sum_{deg V \leq s+m+l-2}E(\xi  YV ) \\
 \varphi(t,\eta,\xi) & = \sum_{deg Y\leq k-1}\sum_{deg Z\leq  s-2}E(tYZ)\sum_{deg U \leq s+m-2}E(\eta YU) \sum_{deg V \leq s+m+l-2}E(\xi  YV ) \\
 \varphi_{1}(t,\eta,\xi) & = \sum_{deg Y\leq k-1}\sum_{deg Z\leq  s-2}E(tYZ)\sum_{deg U \leq s+m-1}E(\eta YU) \sum_{deg V \leq s+m+l-1}E(\xi  YV ) \\
 \varphi_{2}(t,\eta,\xi) & = \sum_{deg Y\leq k-1}\sum_{deg Z\leq  s-2}E(tYZ)\sum_{deg U \leq s+m-2}E(\eta YU) \sum_{deg V \leq s+m+l-1}E(\xi  YV ) \\
 \varphi_{3}(t,\eta,\xi) & = \sum_{deg Y\leq k-1}\sum_{deg Z\leq  s-2}E(tYZ)\sum_{deg U \leq s+m-1}E(\eta YU) \sum_{deg V \leq s+m+l-2}E(\xi  YV ) \\
  \kappa _{1}(t,\eta,\xi) & = \sum_{deg Y\leq k-1}E(tYT^{s-1})E(\eta YT^{s+m-1}) \sum_{deg Z\leq  s-2}E(tYZ)\sum_{deg U \leq s+m-2}E(\eta YU) \sum_{deg V \leq s+m+l-1}E(\xi YV)   \\
 \kappa _{2}(t,\eta,\xi) & = \sum_{deg Y\leq k-1}E(tYT^{s-1}) \sum_{deg Z\leq  s-2}E(tYZ)\sum_{deg U \leq s+m-2}E(\eta YU) \sum_{deg V \leq s+m+l-1}E(\xi YV)   \\
 g _{2}(t,\eta,\xi) & = \sum_{deg Y\leq k-1}E(tYT^{s-1})\sum_{deg Z\leq  s-2}E(tYZ)\sum_{deg U \leq s+m-1}E(\eta YU) \sum_{deg V \leq s+m+l-2}E(\xi YV)   
\end{align*}
Then \\
  \begin{align}
  \label{eq 2.1}
& g_{1}(t,\eta,\xi  ) = 2^{3s+2m+l-1}\cdot
 \sum_{deg Y\leq k-1\atop {Y\in ker D^{\left[\stackrel{s-1}{\stackrel{s+m}{s+m+l}}\right] \times k }(t,\eta ,\xi )}}E(tYT^{s-1}) \\
& =  \begin{cases}
 2^{k+3s+2m+l-1-  r( D^{\left[\stackrel{s-1}{\stackrel{s+m}{s+m+l}}\right] \times k }(t,\eta ,\xi ) )  }  & \text{if }
  r( D^{\left[\stackrel{s-1}{\stackrel{s+m}{s+m+l}}\right] \times k }(t,\eta ,\xi ) ) =  r( D^{\left[\stackrel{s}{\stackrel{s+m}{s+m+l}}\right] \times k }(t,\eta ,\xi ) ) \nonumber \\
     0  & \text{otherwise }.
    \end{cases}
 \end{align}
 
 \begin{align}
   \label{eq 2.2}
& h_{1}(t,\eta,\xi  ) = 2^{3s+2m+l-3}\cdot
 \sum_{deg Y\leq k-1\atop {Y\in ker D^{\left[\stackrel{s-1}{\stackrel{s+m-1}{s+m+l-1}}\right] \times k }(t,\eta ,\xi )}}E(tYT^{s-1}) \\
& =  \begin{cases}
 2^{k+3s+2m+l-3-  r( D^{\left[\stackrel{s-1}{\stackrel{s+m-1}{s+m+l-1}}\right] \times k }(t,\eta ,\xi ) )  }  & \text{if }
  r( D^{\left[\stackrel{s-1}{\stackrel{s+m-1}{s+m+l-1}}\right] \times k }(t,\eta ,\xi ) ) =  r( D^{\left[\stackrel{s}{\stackrel{s+m-1}{s+m+l-1}}\right] \times k }(t,\eta ,\xi ) ) \nonumber \\
     0  & \text{otherwise }.
    \end{cases}
 \end{align}
 
  \begin{align}
    \label{eq 2.3}
& h_{2}(t,\eta,\xi  ) = 2^{3s+2m+l-3}\cdot
 \sum_{deg Y\leq k-1\atop {Y\in ker D^{\left[\stackrel{s-1}{\stackrel{s+m-1}{s+m+l-1}}\right] \times k }(t,\eta ,\xi )}}E(tYT^{s-1})E(\eta YT^{s+m-1}) \\
& =  \begin{cases}
 2^{k+3s+2m+l-3-  r( D^{\left[\stackrel{s-1}{\stackrel{s+m-1}{s+m+l-1}}\right] \times k }(t,\eta ,\xi ) )  }  & \text{if }
  r( D^{\left[\stackrel{s-1}{\stackrel{s+m-1}{s+m+l-1}}\right] \times k }(t,\eta ,\xi ) ) =  r( D^{\left[s-1\atop{ s+m-1\atop {s+m+l-1 \atop{\overline{\alpha_{s -} +  \beta_{s+m -} } }}}\right]\times k} (t,\eta ,\xi ) ) \nonumber \\
     0  & \text{otherwise }.
    \end{cases}
 \end{align}
 
  \begin{align}
    \label{eq 2.4}
& h_{3}(t,\eta,\xi  ) = 2^{3s+2m+l-3}\cdot
 \sum_{deg Y\leq k-1\atop {Y\in ker D^{\left[\stackrel{s-1}{\stackrel{s+m-1}{s+m+l-1}}\right] \times k }(t,\eta ,\xi )}}E(tYT^{s-1})E(\xi  YT^{s+m+l-1}) \\
& =  \begin{cases}
 2^{k+3s+2m+l-3-  r( D^{\left[\stackrel{s-1}{\stackrel{s+m-1}{s+m+l-1}}\right] \times k }(t,\eta ,\xi ) )  }  & \text{if }
  r( D^{\left[\stackrel{s-1}{\stackrel{s+m-1}{s+m+l-1}}\right] \times k }(t,\eta ,\xi ) ) =  r( D^{\left[s-1\atop{ s+m-1\atop {s+m+l-1 \atop{\overline{\alpha_{s -} +  \gamma _{s+m+l -} } }}}\right]\times k} (t,\eta ,\xi ) ) \nonumber \\
     0  & \text{otherwise }.
    \end{cases}
 \end{align}
 
  \begin{align}
    \label{eq 2.5}
& h_{4}(t,\eta,\xi  ) = 2^{3s+2m+l-3}\cdot
 \sum_{deg Y\leq k-1\atop {Y\in ker D^{\left[\stackrel{s-1}{\stackrel{s+m-1}{s+m+l-1}}\right] \times k }(t,\eta ,\xi )}}E(tYT^{s-1}) E(\eta YT^{s+m-1})E(\xi  YT^{s+m+l-1}) \\
& =  \begin{cases}
 2^{k+3s+2m+l-3-  r( D^{\left[\stackrel{s-1}{\stackrel{s+m-1}{s+m+l-1}}\right] \times k }(t,\eta ,\xi ) )  }  & \text{if }
  r( D^{\left[\stackrel{s-1}{\stackrel{s+m-1}{s+m+l-1}}\right] \times k }(t,\eta ,\xi ) ) =  r( D^{\left[s-1\atop{ s+m-1\atop {s+m+l-1 \atop{\overline{\alpha_{s -} + \beta _{s+m -} + \gamma _{s+m+l -} } }}}\right]\times k} (t,\eta ,\xi ) ) \nonumber  \\
     0  & \text{otherwise }.
    \end{cases}
 \end{align}
\begin{equation}
  \label{eq 2.6}
 \varphi(t,\eta,\xi) = 2^{3s+2m+l-3}\cdot
 \sum_{deg Y\leq k-1\atop {Y\in ker D^{\left[\stackrel{s-1}{\stackrel{s+m-1}{s+m+l-1}}\right] \times k }(t,\eta ,\xi )}}1
  =  2^{k+3s+2m+l-3-  r( D^{\left[\stackrel{s-1}{\stackrel{s+m-1}{s+m+l-1}}\right] \times k }(t,\eta ,\xi ) )  } 
  \end{equation}
  \begin{equation}
  \label{eq 2.7}
 \varphi_{1}(t,\eta,\xi) = 2^{3s+2m+l-1}\cdot
 \sum_{deg Y\leq k-1\atop {Y\in ker D^{\left[\stackrel{s-1}{\stackrel{s+m}{s+m+l}}\right] \times k }(t,\eta ,\xi )}}1
  =  2^{k+3s+2m+l-1-  r( D^{\left[\stackrel{s-1}{\stackrel{s+m}{s+m+l}}\right] \times k }(t,\eta ,\xi ) )  } 
  \end{equation}
    \begin{equation}
  \label{eq 2.8}
 \varphi_{2}(t,\eta,\xi) = 2^{3s+2m+l-2}\cdot
 \sum_{deg Y\leq k-1\atop {Y\in ker D^{\left[\stackrel{s-1}{\stackrel{s+m-1}{s+m+l}}\right] \times k }(t,\eta ,\xi )}}1
  =  2^{k+3s+2m+l-2-  r( D^{\left[\stackrel{s-1}{\stackrel{s+m-1}{s+m+l}}\right] \times k }(t,\eta ,\xi ) )  } 
  \end{equation}
    \begin{equation}
  \label{eq 2.9}
 \varphi_{3}(t,\eta,\xi) = 2^{3s+2m+l-2}\cdot
 \sum_{deg Y\leq k-1\atop {Y\in ker D^{\left[\stackrel{s-1}{\stackrel{s+m}{s+m+l-1}}\right] \times k }(t,\eta ,\xi )}}1
  =  2^{k+3s+2m+l-2-  r( D^{\left[\stackrel{s-1}{\stackrel{s+m}{s+m+l-1}}\right] \times k }(t,\eta ,\xi ) )  } 
  \end{equation}
     \begin{align}
    \label{eq 2.10}
& \kappa _{1}(t,\eta,\xi  ) = 2^{3s+2m+l-2}\cdot
 \sum_{deg Y\leq k-1\atop {Y\in ker D^{\left[\stackrel{s-1}{\stackrel{s+m-1}{s+m+l}}\right] \times k }(t,\eta ,\xi )}}E(tYT^{s-1})E(\eta YT^{s+m-1}) \\
& =  \begin{cases}
 2^{k+3s+2m+l-2-  r( D^{\left[\stackrel{s-1}{\stackrel{s+m-1}{s+m+l}}\right] \times k }(t,\eta ,\xi ) )  }  & \text{if }
  r( D^{\left[\stackrel{s-1}{\stackrel{s+m-1}{s+m+l}}\right] \times k }(t,\eta ,\xi ) ) =  r( D^{\left[s-1\atop{ s+m-1\atop {s+m+l \atop{\overline{\alpha_{s -} +  \beta_{s+m -} } }}}\right]\times k} (t,\eta ,\xi ) ) \nonumber \\
     0  & \text{otherwise }.
    \end{cases}
 \end{align}
     \begin{align}
    \label{eq 2.11}
& \kappa _{2}(t,\eta,\xi  ) = 2^{3s+2m+l-2}\cdot
 \sum_{deg Y\leq k-1\atop {Y\in ker D^{\left[\stackrel{s-1}{\stackrel{s+m-1}{s+m+l}}\right] \times k }(t,\eta ,\xi )}}E(tYT^{s-1}) \\
& =  \begin{cases}
 2^{k+3s+2m+l-2-  r( D^{\left[\stackrel{s-1}{\stackrel{s+m-1}{s+m+l}}\right] \times k }(t,\eta ,\xi ) )  }  & \text{if }
  r( D^{\left[\stackrel{s-1}{\stackrel{s+m-1}{s+m+l}}\right] \times k }(t,\eta ,\xi ) ) =  r( D^{\left[s-1\atop{ s+m-1\atop {s+m+l \atop{\overline{\alpha_{s -}  } }}}\right]\times k} (t,\eta ,\xi ) ) \nonumber \\
     0  & \text{otherwise }.
    \end{cases}
 \end{align}
      \begin{align}
    \label{eq 2.12}
& g_{2}(t,\eta,\xi  ) = 2^{3s+2m+l-2}\cdot
 \sum_{deg Y\leq k-1\atop {Y\in ker D^{\left[\stackrel{s-1}{\stackrel{s+m}{s+m+l-1}}\right] \times k }(t,\eta ,\xi )}}E(tYT^{s-1}) \\
& =  \begin{cases}
 2^{k+3s+2m+l-2-  r( D^{\left[\stackrel{s-1}{\stackrel{s+m}{s+m+l-1}}\right] \times k }(t,\eta ,\xi ) )  }  & \text{if }
  r( D^{\left[\stackrel{s-1}{\stackrel{s+m}{s+m+l-1}}\right] \times k }(t,\eta ,\xi ) ) =  r( D^{\left[s-1\atop{ s+m\atop {s+m+l-1 \atop{\overline{\alpha_{s -}  } }}}\right]\times k} (t,\eta ,\xi ) ) \nonumber \\
     0  & \text{otherwise }.
    \end{cases}
 \end{align}
    \end{lem}
    \begin{proof}
 The proof of Lemma \ref{lem 2.1} is somewhat similar to the proofs of the  results obtained in [2, see section 4].
 \end{proof}
 \begin{lem}
\label{lem 2.2}
Let $ (t,\eta ,\xi ) \in  \mathbb{P}\times \mathbb{P}\times \mathbb{P}  $ and q a rational integer $\geq 2,$ then we have 
\begin{align}
\displaystyle
 &  g_{1}^{q}(t,\eta,\xi ) = g_{1}(t,\eta,\xi )\cdot \varphi_{1}^{q-1}(t,\eta,\xi)  \label{eq 2.13} \\
 &  g_{1}(t,\eta,\xi ) =  \varphi_{1}(t,\eta,\xi)  \quad \text{if  } \quad  g_{1}(t,\eta,\xi )  \neq 0   \label{eq 2.14} \\ 
  & \nonumber \\
   &  h_{1}^{q}(t,\eta,\xi ) = h_{1}(t,\eta,\xi )\cdot \varphi^{q-1}(t,\eta,\xi)  \label{eq 2.15} \\
 &  h_{1}(t,\eta,\xi ) =  \varphi(t,\eta,\xi)  \quad \text{if  } \quad  h_{1}(t,\eta,\xi )  \neq 0   \label{eq 2.16} \\
  & \nonumber \\
  &  h_{2}^{q}(t,\eta,\xi ) = h_{2}(t,\eta,\xi )\cdot \varphi^{q-1}(t,\eta,\xi)  \label{eq 2.17} \\
 &  h_{2}(t,\eta,\xi ) =  \varphi(t,\eta,\xi)  \quad \text{if  } \quad  h_{2}(t,\eta,\xi )  \neq 0   \label{eq 2.18} \\
  & \nonumber \\
 &  h_{3}^{q}(t,\eta,\xi ) = h_{3}(t,\eta,\xi )\cdot \varphi^{q-1}(t,\eta,\xi)  \label{eq 2.19} \\
 &  h_{3}(t,\eta,\xi ) =  \varphi(t,\eta,\xi)  \quad \text{if  } \quad  h_{3}(t,\eta,\xi )  \neq 0   \label{eq 2.20} \\
 & \nonumber \\ 
 &  h_{4}^{q}(t,\eta,\xi ) = h_{4}(t,\eta,\xi )\cdot \varphi^{q-1}(t,\eta,\xi)  \label{eq 2.21} \\
 &  h_{4}(t,\eta,\xi ) =  \varphi(t,\eta,\xi)  \quad \text{if  } \quad  h_{4}(t,\eta,\xi )  \neq 0   \label{eq 2.22} \\
  & \nonumber \\
 &  \kappa _{1}^{q}(t,\eta,\xi ) = \kappa _{1}(t,\eta,\xi )\cdot \varphi_{2}^{q-1}(t,\eta,\xi)  \label{eq 2.23} \\
 &  \kappa _{1}(t,\eta,\xi ) =  \varphi_{2}(t,\eta,\xi)  \quad \text{if  } \quad  \kappa _{1}(t,\eta,\xi )  \neq 0   \label{eq 2.24} \\
  & \nonumber \\
  &  \kappa _{2}^{q}(t,\eta,\xi ) = \kappa _{2}(t,\eta,\xi )\cdot \varphi_{2}^{q-1}(t,\eta,\xi)  \label{eq 2.25} \\
 &  \kappa _{2}(t,\eta,\xi ) =  \varphi_{2}(t,\eta,\xi)  \quad \text{if  } \quad  \kappa _{2}(t,\eta,\xi )  \neq 0   \label{eq 2.26} \\
  & \nonumber \\
  &  g _{2}^{q}(t,\eta,\xi ) = g _{2}(t,\eta,\xi )\cdot \varphi_{3}^{q-1}(t,\eta,\xi)  \label{eq 2.27} \\
 &  g _{2}(t,\eta,\xi ) =  \varphi_{3}(t,\eta,\xi)  \quad \text{if  } \quad  g_{2}(t,\eta,\xi )  \neq 0   \label{eq 2.28} 
  \end{align}
 \end{lem}
  \begin{proof}
To prove \eqref{eq 2.13} and \eqref{eq 2.14} we proceed as in [1 section 4], that is :\\
From  \eqref{eq 2.1} we get
\begin{align*}
 &  g_{1}^{2}(t,\eta,\xi  ) 
  =  \big[2^{3s+2m+l-1}\cdot
 \sum_{deg Y_{1}\leq k-1\atop {Y_{1}\in ker D^{\left[\stackrel{s-1}{\stackrel{s+m}{s+m+l}}\right] \times k }(t,\eta ,\xi )}}E(tY_{1}T^{s-1}) \big]
\big[2^{3s+2m+l-1}\cdot
 \sum_{deg Y_{2}\leq k-1\atop {Y_{2}\in ker D^{\left[\stackrel{s-1}{\stackrel{s+m}{s+m+l}}\right] \times k }(t,\eta ,\xi )}}E(tY_{2}T^{s-1}) \big] \\
 & \\
 & \text{ We set }  \left\{\begin{array}{cc}
Y_{1} + Y_{2} = Y_{3}  &  deg Y_{3}\leq k-1, \\
              Y_{2} = Y_{4}   &   deg Y_{4}\leq k-1.
\end{array}\right.\\
& \text{Then we obtain}\\
&  g_{1}^{2}(t,\eta,\xi  ) = 2^{(3s+2m+l-1)2}\sum_{deg Y_{1}\leq k-1\atop {Y_{1}\in ker D^{\left[\stackrel{s-1}{\stackrel{s+m}{s+m+l}}\right] \times k }(t,\eta ,\xi )}} 
 \sum_{deg Y_{2}\leq k-1\atop {Y_{2}\in ker D^{\left[\stackrel{s-1}{\stackrel{s+m}{s+m+l}}\right] \times k }(t,\eta ,\xi )}}E(t(Y_{1}+Y_{2})T^{s-1}) \\
 & =   =  \big[2^{3s+2m+l-1}\cdot
 \sum_{deg Y_{3}\leq k-1\atop {Y_{3}\in ker D^{\left[\stackrel{s-1}{\stackrel{s+m}{s+m+l}}\right] \times k }(t,\eta ,\xi )}}E(tY_{3}T^{s-1}) \big]
\big[2^{3s+2m+l-1}\cdot
 \sum_{deg Y_{4}\leq k-1\atop {Y_{4}\in ker D^{\left[\stackrel{s-1}{\stackrel{s+m}{s+m+l}}\right] \times k }(t,\eta ,\xi )}}1 \big] \\
 & =  g_{1}(t,\eta,\xi)\cdot\varphi _{1}(t,\eta,\xi) 
 \end{align*}
 The other assertions are proved in a similar way.
 \end{proof}
 \begin{lem}
\label{lem 2.3}
Let $ (t,\eta ,\xi ) \in  \mathbb{P}\times \mathbb{P}\times \mathbb{P},  $ and q a rational integer $\geq 2,$ then we have for $1\leq i\leq q-1$\\
\begin{align}
\displaystyle
 &  h_{1}^{i}(t,\eta,\xi )\cdot h_{2}^{q-i}(t,\eta,\xi ) = \label{eq 2.29} \\
& =  \begin{cases}
 \varphi^{q} (t,\eta ,\xi )  & \text{if }
  r( D^{\left[\stackrel{s-1}{\stackrel{s+m-1}{s+m+l-1}}\right] \times k }(t,\eta ,\xi ) ) =  r( D^{\left[s-1\atop{ s+m-1\atop {s+m+l-1 \atop{\overline{\alpha_{s -}  } }}}\right]\times k} (t,\eta ,\xi ) )
=  r( D^{\left[s-1\atop{ s+m-1\atop {s+m+l-1 \atop{\overline{ \beta_{s+m -} } }}}\right]\times k} (t,\eta ,\xi ) )   \nonumber \\
     0  & \text{otherwise }.
    \end{cases}
\end{align} 
 \begin{align}
\displaystyle
 &  h_{1}^{i}(t,\eta,\xi )\cdot h_{3}^{q-i}(t,\eta,\xi ) = \label{eq 2.30} \\
& =  \begin{cases}
 \varphi^{q} (t,\eta ,\xi )  & \text{if }
  r( D^{\left[\stackrel{s-1}{\stackrel{s+m-1}{s+m+l-1}}\right] \times k }(t,\eta ,\xi ) ) =  r( D^{\left[s-1\atop{ s+m-1\atop {s+m+l-1 \atop{\overline{\alpha_{s -}  } }}}\right]\times k} (t,\eta ,\xi ) )
=  r( D^{\left[s-1\atop{ s+m-1\atop {s+m+l-1 \atop{\overline{ \gamma _{s+m+l -} } }}}\right]\times k} (t,\eta ,\xi ) )   \nonumber \\
     0  & \text{otherwise }.
    \end{cases}
\end{align} 
  \begin{align}
\displaystyle
 &  h_{1}^{i}(t,\eta,\xi )\cdot h_{4}^{q-i}(t,\eta,\xi ) = \label{eq 2.31} \\
& =  \begin{cases}
 \varphi^{q} (t,\eta ,\xi )  & \text{if }
  r( D^{\left[\stackrel{s-1}{\stackrel{s+m-1}{s+m+l-1}}\right] \times k }(t,\eta ,\xi ) ) =  r( D^{\left[s-1\atop{ s+m-1\atop {s+m+l-1 \atop{\overline{\alpha_{s -}  } }}}\right]\times k} (t,\eta ,\xi ) )
=  r( D^{\left[s-1\atop{ s+m-1\atop {s+m+l-1 \atop{\overline{\beta _{s+m-} + \gamma _{s+m+l -} } }}}\right]\times k} (t,\eta ,\xi ) )   \nonumber \\
     0  & \text{otherwise }.
    \end{cases}
\end{align} 
 \begin{align}
\displaystyle
 &  h_{2}^{i}(t,\eta,\xi )\cdot h_{3}^{q-i}(t,\eta,\xi ) = \label{eq 2.32} \\
& =  \begin{cases}
 \varphi^{q} (t,\eta ,\xi )  & \text{if }
  r( D^{\left[\stackrel{s-1}{\stackrel{s+m-1}{s+m+l-1}}\right] \times k }(t,\eta ,\xi ) ) =  r( D^{\left[s-1\atop{ s+m-1\atop {s+m+l-1 \atop{\overline{\alpha_{s -} +\beta _{s+m-} } }}}\right]\times k} (t,\eta ,\xi ) )
=  r( D^{\left[s-1\atop{ s+m-1\atop {s+m+l-1 \atop{\overline{\beta _{s+m-} + \gamma _{s+m+l -} } }}}\right]\times k} (t,\eta ,\xi ) )   \nonumber \\
     0  & \text{otherwise }.
    \end{cases}
\end{align} 
 \begin{align}
\displaystyle
 &  h_{2}^{i}(t,\eta,\xi )\cdot h_{4}^{q-i}(t,\eta,\xi ) = \label{eq 2.33} \\
& =  \begin{cases}
 \varphi^{q} (t,\eta ,\xi )  & \text{if }
  r( D^{\left[\stackrel{s-1}{\stackrel{s+m-1}{s+m+l-1}}\right] \times k }(t,\eta ,\xi ) ) =  r( D^{\left[s-1\atop{ s+m-1\atop {s+m+l-1 \atop{\overline{\alpha_{s -} +\beta _{s+m-} } }}}\right]\times k} (t,\eta ,\xi ) )
=  r( D^{\left[s-1\atop{ s+m-1\atop {s+m+l-1 \atop{\overline{ \gamma _{s+m+l -} } }}}\right]\times k} (t,\eta ,\xi ) )   \nonumber \\
     0  & \text{otherwise }.
    \end{cases}
\end{align} 
 \begin{align}
\displaystyle
 &  h_{3}^{i}(t,\eta,\xi )\cdot h_{4}^{q-i}(t,\eta,\xi ) = \label{eq 2.34} \\
& =  \begin{cases}
 \varphi^{q} (t,\eta ,\xi )  & \text{if }
  r( D^{\left[\stackrel{s-1}{\stackrel{s+m-1}{s+m+l-1}}\right] \times k }(t,\eta ,\xi ) ) =  r( D^{\left[s-1\atop{ s+m-1\atop {s+m+l-1 \atop{\overline{\alpha_{s -} +\gamma  _{s+m+l-} } }}}\right]\times k} (t,\eta ,\xi ) )
=  r( D^{\left[s-1\atop{ s+m-1\atop {s+m+l-1 \atop{\overline{ \beta  _{s+m -} } }}}\right]\times k} (t,\eta ,\xi ) )   \nonumber \\
     0  & \text{otherwise }.
    \end{cases}
\end{align} 
 \begin{align}
\displaystyle
 &  \kappa _{1}^{i}(t,\eta,\xi )\cdot \kappa _{2}^{q-i}(t,\eta,\xi ) = \label{eq 2.35} \\
& =  \begin{cases}
 \varphi_{2}^{q} (t,\eta ,\xi )  & \text{if }
  r( D^{\left[\stackrel{s-1}{\stackrel{s+m-1}{s+m+l}}\right] \times k }(t,\eta ,\xi ) ) =  r( D^{\left[s-1\atop{ s+m-1\atop {s+m+l \atop{\overline{\alpha_{s -} +\beta  _{s+m-} } }}}\right]\times k} (t,\eta ,\xi ) )
=  r( D^{\left[s-1\atop{ s+m-1\atop {s+m+l \atop{\overline{ \alpha  _{s -} } }}}\right]\times k} (t,\eta ,\xi ) )   \nonumber \\
     0  & \text{otherwise }.
    \end{cases}
\end{align} 
\end{lem}
 \begin{proof}
 Follows from Lemma \ref{lem 2.1} and Lemma \ref{lem 2.2} using elementary rank considerations.
 \end{proof} 
  \begin{lem}
\label{lem 2.4}
Let $ (t,\eta ,\xi ) \in  \mathbb{P}\times \mathbb{P}\times \mathbb{P},  $ and q a rational integer $\geq 3,$ then we have for any
integers i,j,r such that
\begin{align*}
\displaystyle
   \begin{cases}
   i+j+r = q   \\
   1\leq i,j,r\leq q-2
 \end{cases}
\end{align*} 
\Small
\begin{align}
\displaystyle
\label{eq 2.36}
 &  h_{1}^{i}(t,\eta,\xi )\cdot h_{2}^{j}(t,\eta,\xi )\cdot h_{3}^{r}(t,\eta,\xi )  =  h_{1}^{i}(t,\eta,\xi )\cdot h_{3}^{j}(t,\eta,\xi )\cdot h_{4}^{r}(t,\eta,\xi )  \\
 & =  h_{2}^{i}(t,\eta,\xi )\cdot h_{3}^{j}(t,\eta,\xi )\cdot h_{4}^{r}(t,\eta,\xi ) =  h_{1}^{i}(t,\eta,\xi )\cdot h_{2}^{j}(t,\eta,\xi )\cdot h_{4}^{r}(t,\eta,\xi ) \nonumber \\
 & =  \begin{cases}
 \varphi^{q} (t,\eta ,\xi )  & \text{if }
  r( D^{\left[\stackrel{s-1}{\stackrel{s+m-1}{s+m+l-1}}\right] \times k }(t,\eta ,\xi ) ) =  r( D^{\left[s-1\atop{ s+m-1\atop {s+m+l-1 \atop{\overline{\alpha_{s -}  } }}}\right]\times k} (t,\eta ,\xi ) ) 
  =  r( D^{\left[s-1\atop{ s+m-1\atop {s+m+l-1 \atop{\overline{ \beta_{s+m -} } }}}\right]\times k} (t,\eta ,\xi ) )  
  =  r( D^{\left[s-1\atop{ s+m-1\atop {s+m+l-1 \atop{\overline{ \gamma _{s+m+l -} } }}}\right]\times k} (t,\eta ,\xi ) ) \nonumber \\
    0  & \text{otherwise }.\nonumber
    \end{cases}
\end{align} 
 \end{lem}
  \begin{proof}
 Follows from Lemma \ref{lem 2.1} and Lemma \ref{lem 2.2} using elementary row operations.
 \end{proof} 
  \begin{lem}
\label{lem 2.5}
Let $ (t,\eta ,\xi ) \in  \mathbb{P}\times \mathbb{P}\times \mathbb{P},  $ and q a rational integer $\geq 4,$ then we have for any
integers i,j,r,p such that
\begin{align*}
\displaystyle
   \begin{cases}
   i+j+r+p = q   \\
   1\leq i,j,r,p\leq q-3
 \end{cases}
\end{align*} 
\Small
\begin{align}
\displaystyle
\label{eq 2.37}
 &  h_{1}^{i}(t,\eta,\xi )\cdot h_{2}^{j}(t,\eta,\xi )\cdot h_{3}^{r}(t,\eta,\xi )\cdot h_{4}^{p}(t,\eta,\xi ) \\
  & =  \begin{cases}
 \varphi^{q} (t,\eta ,\xi )  & \text{if }
  r( D^{\left[\stackrel{s-1}{\stackrel{s+m-1}{s+m+l-1}}\right] \times k }(t,\eta ,\xi ) ) =  r( D^{\left[s-1\atop{ s+m-1\atop {s+m+l-1 \atop{\overline{\alpha_{s -}  } }}}\right]\times k} (t,\eta ,\xi ) ) 
  =  r( D^{\left[s-1\atop{ s+m-1\atop {s+m+l-1 \atop{\overline{ \beta_{s+m -} } }}}\right]\times k} (t,\eta ,\xi ) )  
  =  r( D^{\left[s-1\atop{ s+m-1\atop {s+m+l-1 \atop{\overline{ \gamma _{s+m+l -} } }}}\right]\times k} (t,\eta ,\xi ) ) \nonumber \\
    0  & \text{otherwise }.\nonumber
    \end{cases}
\end{align} 
 \end{lem}
\begin{proof}
Similarly to the proof of Lemma \ref{lem 2.4}.
\end{proof}  
   \subsection{Rank properties of partitions on triple persymmetric matrices}
  \label{subsec 4}
    \begin{lem}
\label{lem 2.6}We have :\\
\begin{align}
& \sigma _{i,i}^{\left[\stackrel{s-1}{\stackrel{s+m -1}{\stackrel{s+m+l-1}{\overline{\alpha  _{s -}+\beta _{s+m-}+\gamma _{s+m+l-}}}}}\right] \times k } 
= 2\cdot \sigma _{i,i}^{\left[\stackrel{s-1}{\stackrel{s+m -1}{\stackrel{s+m+l-1}{\overline{\alpha  _{s -}+\beta _{s+m-}}}}}\right] \times k } 
= 2\cdot \sigma _{i,i}^{\left[\stackrel{s-1}{\stackrel{s+m -1}{\stackrel{s+m+l-1}{\overline{\alpha  _{s -}+\gamma  _{s+m+l-}}}}}\right] \times k } 
= 4\cdot \sigma _{i,i}^{\left[\stackrel{s-1}{\stackrel{s+m -1}{\stackrel{s+m+l-1}{\overline{\alpha  _{s -}}}}}\right] \times k } \label{eq 2.38}\\
& \nonumber    \\
& \sigma _{i,i}^{\left[\stackrel{s-1}{\stackrel{s+m -1}{\stackrel{s+m+l}{\overline{\alpha  _{s -}+\beta _{s+m-}}}}}\right] \times k } 
= 2\cdot \sigma _{i,i}^{\left[\stackrel{s-1}{\stackrel{s+m -1}{\stackrel{s+m+l}{\overline{\alpha  _{s -}}}}}\right] \times k } \label{eq 2.39}\\
& \nonumber    \\
& \sigma _{i,i,i}^{\left[s-1\atop{ s+m-1\atop {s+m+l-1 \atop{\overline{\alpha_{s -} \atop{\beta_{s+m -}+ \gamma _{s+m+l-}} }}}}\right]\times k}
= 2\cdot\sigma _{i,i,i}^{\left[s-1\atop{ s+m-1\atop {s+m+l-1 \atop{\overline{\alpha_{s -} \atop{\beta_{s+m -}} }}}}\right]\times k}
= 2\cdot\sigma _{i,i,i}^{\left[s-1\atop{ s+m-1\atop {s+m+l-1 \atop{\overline{\beta _{s+m -} \atop{\alpha _{s -}} }}}}\right]\times k} \label{eq 2.40}\\
&  \nonumber   \\
& \sigma _{i,i,i}^{\left[s-1\atop{ s+m-1\atop {s+m+l-1 \atop{\overline{\alpha_{s -}+\beta _{s+m-} \atop{\beta_{s+m -}+ \gamma _{s+m+l-}} }}}}\right]\times k}
= 2\cdot\sigma _{i,i,i}^{\left[s-1\atop{ s+m-1\atop {s+m+l-1 \atop{\overline{\alpha_{s -} \atop{\beta_{s+m -}} }}}}\right]\times k}
= 2\cdot\sigma _{i,i,i}^{\left[s-1\atop{ s+m-1\atop {s+m+l-1 \atop{\overline{\beta _{s+m -} \atop{\alpha _{s -}} }}}}\right]\times k}\label{eq 2.41}   \\
&  \nonumber  \\
& \sigma _{i,i,i}^{\left[s-1\atop{ s+m-1\atop {s+m+l-1 \atop{\overline{\alpha_{s -}+\beta _{s+m-} \atop{ \gamma _{s+m+l-}} }}}}\right]\times k}
= 2\cdot\sigma _{i,i,i}^{\left[s-1\atop{ s+m-1\atop {s+m+l-1 \atop{\overline{\alpha_{s -} \atop{\gamma _{s+m+l -}} }}}}\right]\times k}
= 2\cdot\sigma _{i,i,i}^{\left[s-1\atop{ s+m-1\atop {s+m+l-1 \atop{\overline{\gamma  _{s+m+l -} \atop{\alpha _{s -}} }}}}\right]\times k} \label{eq 2.42} \\
&   \nonumber  \\
& \sigma _{i,i,i}^{\left[s-1\atop{ s+m-1\atop {s+m+l-1 \atop{\overline{\alpha_{s -}+\gamma  _{s+m+l-} \atop{ \beta  _{s+m-}} }}}}\right]\times k}
= 2\cdot\sigma _{i,i,i}^{\left[s-1\atop{ s+m-1\atop {s+m+l-1 \atop{\overline{\alpha_{s -} \atop{\beta  _{s+m -}} }}}}\right]\times k}
= 2\cdot\sigma _{i,i,i}^{\left[s-1\atop{ s+m-1\atop {s+m+l-1 \atop{\overline{\beta  _{s+m -} \atop{\alpha _{s -}} }}}}\right]\times k} \label{eq 2.43} 
 \end{align}
 \end{lem}
  \begin{proof}
  By way of example we establish the following equality in \eqref{eq 2.38}
  $$ \sigma _{i,i}^{\left[\stackrel{s-1}{\stackrel{s+m -1}{\stackrel{s+m+l-1}{\overline{\alpha  _{s -}+\beta _{s+m-}+\gamma _{s+m+l-}}}}}\right] \times k } 
= 4\cdot \sigma _{i,i}^{\left[\stackrel{s-1}{\stackrel{s+m -1}{\stackrel{s+m+l-1}{\overline{\alpha  _{s -}}}}}\right] \times k }$$
  Consider the following $(3s+2m+l-2)\times k $ matrix  denoted by 
   $  D^{\left[\stackrel{s-1}{\stackrel{s+m -1}{\stackrel{s+m+l-1}{\overline{\alpha  _{s -}+\beta _{s+m-}+\gamma _{s+m+l-}}}}}\right] \times k } (t,\eta,\xi   ) $
\tiny  
 $  \bordermatrix{%
                    &                   &                                           &                            &    \cr
   r_{1} &    \alpha _{1} & \alpha _{2}  &  \ldots & \alpha _{k-1}  &  \alpha _{k} \cr
 r_{2} &   \alpha _{2 } & \alpha _{3} &  \ldots  &  \alpha _{k} &  \alpha _{k+1} \cr
  \vdots &   \vdots  & \vdots    &  \vdots & \vdots  &  \vdots  \cr
 r_{s-1}   &  \alpha _{s-1} & \alpha _{s} & \ldots  &  \alpha _{s+k-3} &  \alpha _{s+k-2}  \cr
r_{s}   & \beta  _{1} & \beta  _{2}  & \ldots  &  \beta_{k-1} &  \beta _{k}  \cr
 r_{s+1} & \beta  _{2} & \beta  _{3}  & \ldots  &  \beta_{k} &  \beta _{k+1}  \cr
\vdots  &  \vdots  & \vdots    &  \vdots & \vdots  &  \vdots \cr
r_{s+m} &      \beta  _{m+1} & \beta  _{m+2}  & \ldots  &  \beta_{k+m-1} &  \beta _{k+m}  \cr
\vdots  & \vdots & \vdots    &  \vdots & \vdots  &  \vdots \cr
   r_{2s+m-2} & \beta  _{s+m-1} & \beta  _{s+m}  & \ldots  &  \beta_{s+m+k-3} &  \beta _{s+m+k-2}  \cr
 r_{2s+m-1} &  \gamma  _{1} & \gamma   _{2}  & \ldots  & \gamma  _{k-1} &  \gamma  _{k}  \cr
 r_{2s+m} &\gamma  _{2} & \gamma  _{3}  & \ldots  & \gamma  _{k} &  \gamma  _{k+1}  \cr
\vdots & \vdots   &  \vdots & \vdots  &  \vdots \cr
 r_{2s+2m+l-1} & \gamma  _{m+l+1} &  \gamma _{m+l+2}  & \ldots  & \gamma _{k+m+l-1} &  \gamma  _{k+m+l}  \cr
\vdots & \vdots   &  \vdots  & \vdots  &  \vdots \cr
 r_{3s+2m-3} &\gamma  _{s+m+l-1} & \gamma  _{s+m+l} & \ldots  & \gamma  _{s+m+l+k-3} &  \gamma  _{s+m+l+k-2}\cr 
   \hline
r_{3s+2m-2} & \alpha _{s} + \beta  _{s+m} +\gamma _{s+m+l} & \alpha _{s+1} + \beta  _{s+m+1} +\gamma _{s+m+l+1}  & \ldots 
& \alpha _{s+k-2} + \beta  _{s+k+m-2}  + \gamma _{s+m+l+k-2}&  \alpha _{s+k-1} + \beta  _{s+k+m-1}  +\gamma _{s+m+l+k-1}\cr
} $.  \vspace{0.5 cm}\\
\normalsize

 We recall  that the rank of a matrix does not change under elementary row operations.\\
   On the above  matrix, we add to the j-th row  the s+m+j-th row and the $ 2s+2m+l-1+j-th $ row for $ 0\leq j\leq s-2 $ obtaining 
   $$  \left( \begin{smallmatrix}
  \alpha _{1} + \beta  _{m+1}+\gamma _{m+l+1} & \alpha _{2} + \beta  _{m+2} +\gamma _{m+l+2}  & \ldots & 
\alpha _{k-1}  + \beta  _{m+k-1}+\gamma _{m+l+k-1} &  \alpha _{k}  + \beta  _{m+k}+\gamma _{m+l+k} \\
 \alpha _{2 } + \beta  _{m+2}+\gamma _{m+l+2} & \alpha _{3} + \beta  _{m+3}+\gamma _{m+l+3} &  \ldots  & 
 \alpha _{k} + \beta  _{m+k}+\gamma _{m+l+k} &  \alpha _{k+1} + \beta  _{m+k+1}+\gamma _{m+l+k+1}\\
  \vdots &  \vdots  & \vdots    &  \vdots & \vdots  \\
   \alpha _{s-1} + \beta  _{s-1+m}+\gamma _{s-1+m+l} & \alpha _{s} + \beta  _{s+m} +\gamma _{s+m+l} & \ldots  & 
 \alpha _{s+k-3} + \beta  _{s+k+m-3}+\gamma _{s+k+m+l-3} &  \alpha _{s+k-2}  + \beta  _{s+k+m-2} +\gamma _{s+k+m+l-2}\\
 \beta  _{1} & \beta  _{2}  & \ldots  &  \beta_{k-1} &  \beta _{k}  \\
\beta  _{2} & \beta  _{3}  & \ldots  &  \beta_{k} &  \beta _{k+1}  \\
\vdots & \vdots   &  \vdots & \vdots  &  \vdots \\
\beta  _{m+1} & \beta  _{m+2}  & \ldots  &  \beta_{k+m-1} &  \beta _{k+m}  \\
\vdots & \vdots    &  \vdots & \vdots  &  \vdots \\
\beta  _{s+m-1} & \beta  _{s+m}  & \ldots  &  \beta_{s+m+k-3} &  \beta _{s+m+k-2}  \\
\gamma  _{1} & \gamma   _{2} & \ldots  & \gamma  _{k-1} &  \gamma  _{k}  \\
\gamma  _{2} & \gamma  _{3}  & \ldots  & \gamma  _{k} &  \gamma  _{k+1}  \\
\vdots & \vdots   &  \vdots & \vdots  &  \vdots \\
 \gamma  _{m+l+1} &  \gamma _{m+l+2}  & \ldots  & \gamma _{k+m+l-1} &  \gamma  _{k+m+l}  \\
\vdots & \vdots     &  \vdots  & \vdots  &  \vdots \\
\gamma  _{s+m+l-1} & \gamma  _{s+m+l}  & \ldots  & \gamma  _{s+m+l+k-3} &  \gamma  _{s+m+l+k-2}  \\  
\hline \\
\alpha _{s}+ \beta  _{s+m}+\gamma _{s+m+l}  & \alpha _{s+1}+ \beta  _{s+m+1} + \gamma _{s+m+l+1} & \ldots  &  \alpha _{s+k-2}+ \beta_{s+m+k-2} +\gamma _{s+m+l+k-2} &  \alpha _{s+k-1}+   \beta _{s+m+k-1}+\gamma _{s+m+l+k-1}
  \end{smallmatrix}\right) $$\\
  $$\big\uparrow $$
In the above  matrix we set 
$\begin{cases}
 \alpha _{i}+\beta _{i+m}+\gamma _{i+m+l} & = \alpha  _{i}'\quad \text{for} \quad 1\leq i\leq k+s-2, \\
  \beta _{j} & = \beta  _{j}\quad \text{for}\quad  1\leq j \leq k+s+m-2.\\
  \gamma _{r}& =\gamma _{r}\quad \text{for}\quad  1\leq r \leq k+s+m+l-2.\\
\end{cases} $ 
We remark that the map $ \kappa:  \mathbb{F}_{2}^{3k+3s+2m+l-6}\longmapsto \mathbb{F}_{2}^{3k+3s+2m+l-6}$
defined by \vspace{0.1 cm} \\
 $ (\alpha _{1},\alpha _{2}, \ldots,\alpha _{k+s-2},\beta _{1},
 \beta _{2},\ldots,\beta _{k+s+m-2},\gamma _{1},\gamma _{2},\ldots,\gamma _{k+s+m+l-2})\\
 \longmapsto 
  (\alpha _{1}',\alpha _{2}', \ldots,\alpha _{k+s-2}',\beta _{1},
 \beta _{2},\ldots,\beta _{k+s+m-2},\gamma _{1},\gamma _{2},\ldots,\gamma _{k+s+m+l-2})$\\
   is an isomorphisme. We then obtain  the below  matrix \vspace{0.1 cm} \\
 \centering{$\downarrow $}
 $$ \left( \begin{smallmatrix}
  \alpha _{1}'  & \alpha _{2}'  & \ldots & 
\alpha _{k-1}'   &  \alpha _{k}'   \\
 \alpha _{2 }'  & \alpha _{3}'  &  \ldots  & 
 \alpha _{k}'  &  \alpha _{k+1}' \\
  \vdots &  \vdots  & \vdots    &  \vdots & \vdots  \\
   \alpha _{s-1}'  & \alpha _{s}'  & \ldots  & 
 \alpha _{s+k-3}'  &  \alpha _{s+k-2}'  \\
 \beta  _{1} & \beta  _{2}  & \ldots  &  \beta_{k-1} &  \beta _{k}  \\
\beta  _{2} & \beta  _{3}  & \ldots  &  \beta_{k} &  \beta _{k+1}  \\
\vdots & \vdots   &  \vdots & \vdots  &  \vdots \\
\beta  _{m+1} & \beta  _{m+2}  & \ldots  &  \beta_{k+m-1} &  \beta _{k+m}  \\
\vdots & \vdots    &  \vdots & \vdots  &  \vdots \\
\beta  _{s+m-1} & \beta  _{s+m}  & \ldots  &  \beta_{s+m+k-3} &  \beta _{s+m+k-2}  \\
\gamma  _{1} & \gamma   _{2} & \ldots  & \gamma  _{k-1} &  \gamma  _{k}  \\
\gamma  _{2} & \gamma  _{3}  & \ldots  & \gamma  _{k} &  \gamma  _{k+1}  \\
\vdots & \vdots   &  \vdots & \vdots  &  \vdots \\
 \gamma  _{m+l+1} &  \gamma _{m+l+2}  & \ldots  & \gamma _{k+m+l-1} &  \gamma  _{k+m+l}  \\
\vdots & \vdots     &  \vdots  & \vdots  &  \vdots \\
\gamma  _{s+m+l-1} & \gamma  _{s+m+l}  & \ldots  & \gamma  _{s+m+l+k-3} &  \gamma  _{s+m+l+k-2}  \\  
\hline \\
\alpha _{s}' & \alpha _{s+1}' & \ldots  &  \alpha _{s+k-2}' &  \alpha _{s+k-1}+   \beta _{s+m+k-1}+\gamma _{s+m+l+k-1}
  \end{smallmatrix}\right) $$\\
 
  We then get  the following equivalences : 
  \begin{align*}
&  rank \left( \begin{smallmatrix}
 \alpha _{1} & \alpha _{2}  &  \ldots & \alpha _{k-1}  &  \alpha _{k} \\
\alpha _{2 } & \alpha _{3} &  \ldots  &  \alpha _{k} &  \alpha _{k+1} \\
\vdots & \vdots    &  \vdots & \vdots  &  \vdots \\
\alpha _{s-1} & \alpha _{s}  & \ldots  &  \alpha _{s+k-3} &  \alpha _{s+k-2}  \\
\beta  _{1} & \beta  _{2}  & \ldots  &  \beta_{k-1} &  \beta _{k}  \\
\beta  _{2} & \beta  _{3}  & \ldots  &  \beta_{k} &  \beta _{k+1}  \\
\vdots & \vdots    &  \vdots & \vdots  &  \vdots \\
\beta  _{m+1} & \beta  _{m+2}  & \ldots  &  \beta_{k+m-1} &  \beta _{k+m}  \\
\vdots & \vdots    &  \vdots & \vdots  &  \vdots \\
\beta  _{s+m-1} & \beta  _{s+m}  & \ldots  &  \beta_{s+m+k-3} &  \beta _{s+m+k-2}  \\
\gamma  _{1} & \gamma   _{2}  & \ldots  & \gamma  _{k-1} &  \gamma  _{k}  \\
\gamma  _{2} & \gamma  _{3}  & \ldots  & \gamma  _{k} &  \gamma  _{k+1}  \\
\vdots & \vdots   &  \vdots & \vdots  &  \vdots \\
 \gamma  _{m+l+1} &  \gamma _{m+l+2}  & \ldots  & \gamma _{k+m+l-1} &  \gamma  _{k+m+l}  \\
\vdots & \vdots    &  \vdots & \vdots  &  \vdots \\
\gamma  _{s+m+l-1} & \gamma  _{s+m+l}  & \ldots  & \gamma  _{s+m+l+k-3} &  \gamma  _{s+m+l+k-2}  \\
\hline \\
 \alpha _{s}+ \beta  _{s+m}+\gamma _{s+m+l}  & \alpha _{s+1}+ \beta  _{s+m+1} + \gamma _{s+m+l+1}  & \ldots  &  \alpha _{s+k-2}+ \beta_{s+m+k-2} +\gamma _{s+m+l+k-2} &  \alpha _{s+k-1}+   \beta _{s+m+k-1}+\gamma _{s+m+l+k-1}
  \end{smallmatrix}\right) \\
 & \\
  &  =  rank \left( \begin{smallmatrix}
 \alpha _{1} & \alpha _{2} & \alpha _{3} &  \ldots & \alpha _{k-1}  &  \alpha _{k} \\
\alpha _{2 } & \alpha _{3} & \alpha _{4}&  \ldots  &  \alpha _{k} &  \alpha _{k+1} \\
\vdots & \vdots & \vdots    &  \vdots & \vdots  &  \vdots \\
\alpha _{s-1} & \alpha _{s} & \alpha _{s +1} & \ldots  &  \alpha _{s+k-3} &  \alpha _{s+k-2}  \\
\beta  _{1} & \beta  _{2} & \beta  _{3} & \ldots  &  \beta_{k-1} &  \beta _{k}  \\
\beta  _{2} & \beta  _{3} & \beta  _{4} & \ldots  &  \beta_{k} &  \beta _{k+1}  \\
\vdots & \vdots & \vdots    &  \vdots & \vdots  &  \vdots \\
\beta  _{m+1} & \beta  _{m+2} & \beta  _{m+3} & \ldots  &  \beta_{k+m-1} &  \beta _{k+m}  \\
\vdots & \vdots & \vdots    &  \vdots & \vdots  &  \vdots \\
\beta  _{s+m-1} & \beta  _{s+m} & \beta  _{s+m+1} & \ldots  &  \beta_{s+m+k-3} &  \beta _{s+m+k-2}\\
\gamma  _{1} & \gamma   _{2}& \gamma   _{3}  & \ldots  & \gamma  _{k-1} &  \gamma  _{k}  \\
\gamma  _{2} & \gamma  _{3} & \gamma   _{4} & \ldots  & \gamma  _{k} &  \gamma  _{k+1}  \\
\vdots & \vdots   &  \vdots &  \vdots & \vdots  &  \vdots \\
 \gamma  _{m+l+1} &  \gamma _{m+l+2} & \gamma _{m+l+3} & \ldots  & \gamma _{k+m+l-1} &  \gamma  _{k+m+l}  \\
\vdots & \vdots    &  \vdots &  \vdots  & \vdots  &  \vdots \\
\gamma  _{s+m+l-1} & \gamma  _{s+m+l} & \gamma  _{s+m+l+1} & \ldots  & \gamma  _{s+m+l+k-3} &  \gamma  _{s+m+l+k-2}  \\  
   \end{smallmatrix}\right) \vspace{2 cm}\\
&  \Longleftrightarrow  \\
&   rank \left( \begin{smallmatrix}
  \alpha _{1} + \beta  _{m+1}+\gamma _{m+l+1} & \alpha _{2} + \beta  _{m+2} +\gamma _{m+l+2}  & \ldots & 
\alpha _{k-1}  + \beta  _{m+k-1}+\gamma _{m+l+k-1} &  \alpha _{k}  + \beta  _{m+k}+\gamma _{m+l+k} \\
 \alpha _{2 } + \beta  _{m+2}+\gamma _{m+l+2} & \alpha _{3} + \beta  _{m+3}+\gamma _{m+l+3} &  \ldots  & 
 \alpha _{k} + \beta  _{m+k}+\gamma _{m+l+k} &  \alpha _{k+1} + \beta  _{m+k+1}+\gamma _{m+l+k+1}\\
  \vdots &  \vdots  & \vdots    &  \vdots & \vdots  \\
   \alpha _{s-1} + \beta  _{s-1+m}+\gamma _{s-1+m+l} & \alpha _{s} + \beta  _{s+m} +\gamma _{s+m+l} & \ldots  & 
 \alpha _{s+k-3} + \beta  _{s+k+m-3}+\gamma _{s+k+m+l-3} &  \alpha _{s+k-2}  + \beta  _{s+k+m-2} +\gamma _{s+k+m+l-2}\\
 \beta  _{1} & \beta  _{2}  & \ldots  &  \beta_{k-1} &  \beta _{k}  \\
\beta  _{2} & \beta  _{3}  & \ldots  &  \beta_{k} &  \beta _{k+1}  \\
\vdots & \vdots   &  \vdots & \vdots  &  \vdots \\
\beta  _{m+1} & \beta  _{m+2}  & \ldots  &  \beta_{k+m-1} &  \beta _{k+m}  \\
\vdots & \vdots    &  \vdots & \vdots  &  \vdots \\
\beta  _{s+m-1} & \beta  _{s+m}  & \ldots  &  \beta_{s+m+k-3} &  \beta _{s+m+k-2}  \\
\gamma  _{1} & \gamma   _{2} & \ldots  & \gamma  _{k-1} &  \gamma  _{k}  \\
\gamma  _{2} & \gamma  _{3}  & \ldots  & \gamma  _{k} &  \gamma  _{k+1}  \\
\vdots & \vdots   &  \vdots & \vdots  &  \vdots \\
 \gamma  _{m+l+1} &  \gamma _{m+l+2}  & \ldots  & \gamma _{k+m+l-1} &  \gamma  _{k+m+l}  \\
\vdots & \vdots     &  \vdots  & \vdots  &  \vdots \\
\gamma  _{s+m+l-1} & \gamma  _{s+m+l}  & \ldots  & \gamma  _{s+m+l+k-3} &  \gamma  _{s+m+l+k-2}  \\  
\hline \\
\alpha _{s}+ \beta  _{s+m}+\gamma _{s+m+l}  & \alpha _{s+1}+ \beta  _{s+m+1} + \gamma _{s+m+l+1} & \ldots  &  \alpha _{s+k-2}+ \beta_{s+m+k-2} +\gamma _{s+m+l+k-2} &  \alpha _{s+k-1}+   \beta _{s+m+k-1}+\gamma _{s+m+l+k-1}
  \end{smallmatrix}\right) \\
  & \\
  &  =  rank \left( \begin{smallmatrix}
 \alpha _{1} & \alpha _{2} & \alpha _{3} &  \ldots & \alpha _{k-1}  &  \alpha _{k} \\
\alpha _{2 } & \alpha _{3} & \alpha _{4}&  \ldots  &  \alpha _{k} &  \alpha _{k+1} \\
\vdots & \vdots & \vdots    &  \vdots & \vdots  &  \vdots \\
\alpha _{s-1} & \alpha _{s} & \alpha _{s +1} & \ldots  &  \alpha _{s+k-3} &  \alpha _{s+k-2}  \\
\beta  _{1} & \beta  _{2} & \beta  _{3} & \ldots  &  \beta_{k-1} &  \beta _{k}  \\
\beta  _{2} & \beta  _{3} & \beta  _{4} & \ldots  &  \beta_{k} &  \beta _{k+1}  \\
\vdots & \vdots & \vdots    &  \vdots & \vdots  &  \vdots \\
\beta  _{m+1} & \beta  _{m+2} & \beta  _{m+3} & \ldots  &  \beta_{k+m-1} &  \beta _{k+m}  \\
\vdots & \vdots & \vdots    &  \vdots & \vdots  &  \vdots \\
\beta  _{s+m-1} & \beta  _{s+m} & \beta  _{s+m+1} & \ldots  &  \beta_{s+m+k-3} &  \beta _{s+m+k-2}\\
\gamma  _{1} & \gamma   _{2}& \gamma   _{3}  & \ldots  & \gamma  _{k-1} &  \gamma  _{k}  \\
\gamma  _{2} & \gamma  _{3} & \gamma   _{4} & \ldots  & \gamma  _{k} &  \gamma  _{k+1}  \\
\vdots & \vdots   &  \vdots &  \vdots & \vdots  &  \vdots \\
 \gamma  _{m+l+1} &  \gamma _{m+l+2} & \gamma _{m+l+3} & \ldots  & \gamma _{k+m+l-1} &  \gamma  _{k+m+l}  \\
\vdots & \vdots    &  \vdots &  \vdots  & \vdots  &  \vdots \\
\gamma  _{s+m+l-1} & \gamma  _{s+m+l} & \gamma  _{s+m+l+1} & \ldots  & \gamma  _{s+m+l+k-3} &  \gamma  _{s+m+l+k-2}  \\  
   \end{smallmatrix}\right) \vspace{2 cm}\\
& \Longleftrightarrow   \vspace{2 cm}\\ 
&   rank \left( \begin{smallmatrix}
  \alpha _{1}'  & \alpha _{2}'  & \ldots & 
\alpha _{k-1}'   &  \alpha _{k}'   \\
 \alpha _{2 }'  & \alpha _{3}'  &  \ldots  & 
 \alpha _{k}'  &  \alpha _{k+1}' \\
  \vdots &  \vdots  & \vdots    &  \vdots & \vdots  \\
   \alpha _{s-1}'  & \alpha _{s}'  & \ldots  & 
 \alpha _{s+k-3}'  &  \alpha _{s+k-2}'  \\
 \beta  _{1} & \beta  _{2}  & \ldots  &  \beta_{k-1} &  \beta _{k}  \\
\beta  _{2} & \beta  _{3}  & \ldots  &  \beta_{k} &  \beta _{k+1}  \\
\vdots & \vdots   &  \vdots & \vdots  &  \vdots \\
\beta  _{m+1} & \beta  _{m+2}  & \ldots  &  \beta_{k+m-1} &  \beta _{k+m}  \\
\vdots & \vdots    &  \vdots & \vdots  &  \vdots \\
\beta  _{s+m-1} & \beta  _{s+m}  & \ldots  &  \beta_{s+m+k-3} &  \beta _{s+m+k-2}  \\
\gamma  _{1} & \gamma   _{2} & \ldots  & \gamma  _{k-1} &  \gamma  _{k}  \\
\gamma  _{2} & \gamma  _{3}  & \ldots  & \gamma  _{k} &  \gamma  _{k+1}  \\
\vdots & \vdots   &  \vdots & \vdots  &  \vdots \\
 \gamma  _{m+l+1} &  \gamma _{m+l+2}  & \ldots  & \gamma _{k+m+l-1} &  \gamma  _{k+m+l}  \\
\vdots & \vdots     &  \vdots  & \vdots  &  \vdots \\
\gamma  _{s+m+l-1} & \gamma  _{s+m+l}  & \ldots  & \gamma  _{s+m+l+k-3} &  \gamma  _{s+m+l+k-2}  \\  
\hline \\
\alpha _{s}' & \alpha _{s+1}' & \ldots  &  \alpha _{s+k-2}' &  \alpha _{s+k-1}+   \beta _{s+m+k-1}+\gamma _{s+m+l+k-1}
  \end{smallmatrix}\right) \\
  & \\
  &  =  rank \left( \begin{smallmatrix}
 \alpha _{1}' & \alpha _{2}' & \alpha _{3}' &  \ldots & \alpha _{k-1}'  &  \alpha _{k}' \\
\alpha _{2 }' & \alpha _{3}' & \alpha _{4}'&  \ldots  &  \alpha _{k}' &  \alpha _{k+1}' \\
\vdots & \vdots & \vdots    &  \vdots & \vdots  &  \vdots \\
\alpha _{s-1}' & \alpha _{s}' & \alpha _{s +1}' & \ldots  &  \alpha _{s+k-3}' &  \alpha _{s+k-2}'  \\
\beta  _{1} & \beta  _{2} & \beta  _{3} & \ldots  &  \beta_{k-1} &  \beta _{k}  \\
\beta  _{2} & \beta  _{3} & \beta  _{4} & \ldots  &  \beta_{k} &  \beta _{k+1}  \\
\vdots & \vdots & \vdots    &  \vdots & \vdots  &  \vdots \\
\beta  _{m+1} & \beta  _{m+2} & \beta  _{m+3} & \ldots  &  \beta_{k+m-1} &  \beta _{k+m}  \\
\vdots & \vdots & \vdots    &  \vdots & \vdots  &  \vdots \\
\beta  _{s+m-1} & \beta  _{s+m} & \beta  _{s+m+1} & \ldots  &  \beta_{s+m+k-3} &  \beta _{s+m+k-2}\\
\gamma  _{1} & \gamma   _{2}& \gamma   _{3}  & \ldots  & \gamma  _{k-1} &  \gamma  _{k}  \\
\gamma  _{2} & \gamma  _{3} & \gamma   _{4} & \ldots  & \gamma  _{k} &  \gamma  _{k+1}  \\
\vdots & \vdots   &  \vdots &  \vdots & \vdots  &  \vdots \\
 \gamma  _{m+l+1} &  \gamma _{m+l+2} & \gamma _{m+l+3} & \ldots  & \gamma _{k+m+l-1} &  \gamma  _{k+m+l}  \\
\vdots & \vdots    &  \vdots &  \vdots  & \vdots  &  \vdots \\
\gamma  _{s+m+l-1} & \gamma  _{s+m+l} & \gamma  _{s+m+l+1} & \ldots  & \gamma  _{s+m+l+k-3} &  \gamma  _{s+m+l+k-2}  \\  
   \end{smallmatrix}\right) \vspace{2 cm}\\
  \end{align*}
 Observing that for all $ \alpha  _{s+k-1}'\in \mathbb{F}_{2} $\\
 $$ Card\left\{(\alpha _{s+k-1},\beta _{s+k+m-1},\gamma _{s+k+m+l-1})\in \mathbb{F}_{2}^{3} \vert \quad
\alpha _{s+k-1} +  \beta _{s+k+m-1}+\gamma _{s+k+m+l-1} = \alpha  _{s+k-1}'\right\} = 4, $$ \\
     we obtain  by combining the above results  the desired equality \\
 The other equalities are obtained in a somewhat similarly way.
 \end{proof}
  \subsection{Integration formulas for exponential sums over the unit interval of $ \mathbb{K}^3 $ }
  \label{subsec 5}
    \begin{lem}
\label{lem 2.7}
Let $ (t,\eta ,\xi ) \in  \mathbb{P}\times \mathbb{P}\times \mathbb{P}  $ and q a rational integer $\geq 2,$ then we have \\
\begin{align}
\displaystyle
 & \int_{\mathbb{P}^{3}}  g_{1}^{q}(t,\eta,\xi )dtd\eta d\xi
  =\int_{\left\{(t,\eta ,\xi )\in\mathbb{P}^{3}\mid  g_{1}(t,\eta,\xi )\neq 0 \right\}} \varphi _{1}^{q}(t,\eta,\xi  )dtd\eta d\xi \label{eq 2.44}\\
   & =  2^{(k+3s+2m+l-1)q -(3k +3s+2m+l-3)}
 \sum_{i=0}^{\inf(k,3s+2m+l-1)}\sigma _{i,i}^{\left[s-1\atop{ s+m\atop {s+m+l \atop{\overline{\alpha_{s -}  } }}}\right]\times k}\cdot2^{-iq} \nonumber\\
 & = 2^{(k+3s+2m+l-1)q -(3k +3s+2m+l-3)}
 \sum_{i=0}^{\inf(k,3s+2m+l-1)}\big[\sigma _{i,i,i,i}^{\left[\beta | \gamma  \atop{ \gamma | \beta \atop \alpha | \alpha  } \right]\times k}
+  \sigma _{i-1,i,i,i}^{\left[\beta | \gamma  \atop{ \gamma | \beta \atop \alpha | \alpha  } \right]\times k}
+  \sigma _{i-1,i-1,i,i}^{\left[\beta | \gamma  \atop{ \gamma | \beta \atop \alpha | \alpha  } \right]\times k}
+  \sigma _{i-2,i-1,i,i}^{\left[\beta | \gamma  \atop{ \gamma | \beta \atop \alpha | \alpha  } \right]\times k}\big]\cdot2^{-iq} \nonumber  \\
   & \int_{\mathbb{P}^{3}}  h_{1}^{q}(t,\eta,\xi )dtd\eta d\xi =\int_{\left\{(t,\eta ,\xi )\in\mathbb{P}^{3}\mid  h_{1}(t,\eta,\xi )\neq 0 \right\}} \varphi^{q}(t,\eta,\xi  )dtd\eta d\xi \label{eq 2.45}\\
   & =  2^{(k+3s+2m+l-3)q -(3k +3s+2m+l-5)}
 \sum_{i=0}^{\inf(k,3s+2m+l-3)}\sigma _{i,i}^{\left[s-1\atop{ s+m-1\atop {s+m+l-1 \atop{\overline{\alpha_{s -}  } }}}\right]\times k}\cdot2^{-iq} \nonumber\\
  & = 2^{(k+3s+2m+l-3)q -(3k +3s+2m+l-3)}
 \sum_{i=0}^{\inf(k,3s+2m+l-3)}\big[\sigma _{i,i,i,i}^{\left[\alpha  | \alpha  \atop{ \beta  | \gamma \atop \gamma  | \beta   } \right]\times k}
+ \sigma _{i,i,i,i+1}^{\left[\alpha  | \alpha  \atop{ \beta  | \gamma \atop \gamma  | \beta   } \right]\times k}
+ \sigma _{i,i,i+1,i+1}^{\left[\alpha  | \alpha  \atop{ \beta  | \gamma \atop \gamma  | \beta   } \right]\times k}
+ \sigma _{i,i,i+1,i+2}^{\left[\alpha  | \alpha  \atop{ \beta  | \gamma \atop \gamma  | \beta   } \right]\times k} \big]\cdot2^{-iq} \nonumber \\
  & \int_{\mathbb{P}^{3}}  h_{2}^{q}(t,\eta,\xi )dtd\eta d\xi =\int_{\left\{(t,\eta ,\xi )\in\mathbb{P}^{3}\mid  h_{2}(t,\eta,\xi )\neq 0 \right\}} \varphi^{q}(t,\eta,\xi  )dtd\eta d\xi \label{eq 2.46}\\
   & =  2^{(k+3s+2m+l-3)q -(3k +3s+2m+l-4)}
 \sum_{i=0}^{\inf(k,3s+2m+l-3)}\sigma _{i,i}^{\left[s-1\atop{ s+m-1\atop {s+m+l-1 \atop{\overline{\alpha_{s -} +\beta _{s+m-} } }}}\right]\times k}\cdot2^{-iq} \nonumber\\
   & =  2^{(k+3s+2m+l-3)q -(3k +3s+2m+l-4)}
 \sum_{i=0}^{\inf(k,3s+2m+l-3)}2\cdot\sigma _{i,i}^{\left[s-1\atop{ s+m-1\atop {s+m+l-1 \atop{\overline{\alpha_{s -}  } }}}\right]\times k}\cdot2^{-iq} \nonumber\\
  & = 2^{(k+3s+2m+l-3)q -(3k +3s+2m+l-3)}
 \sum_{i=0}^{\inf(k,3s+2m+l-3)}\big[\sigma _{i,i,i,i}^{\left[\alpha  | \alpha  \atop{ \beta  | \gamma \atop \gamma  | \beta   } \right]\times k}
+ \sigma _{i,i,i,i+1}^{\left[\alpha  | \alpha  \atop{ \beta  | \gamma \atop \gamma  | \beta   } \right]\times k}
+ \sigma _{i,i,i+1,i+1}^{\left[\alpha  | \alpha  \atop{ \beta  | \gamma \atop \gamma  | \beta   } \right]\times k}
+ \sigma _{i,i,i+1,i+2}^{\left[\alpha  | \alpha  \atop{ \beta  | \gamma \atop \gamma  | \beta   } \right]\times k} \big]\cdot2^{-iq} \nonumber \\
 & \int_{\mathbb{P}^{3}}  h_{3}^{q}(t,\eta,\xi )dtd\eta d\xi =\int_{\left\{(t,\eta ,\xi )\in\mathbb{P}^{3}\mid  h_{3}(t,\eta,\xi )\neq 0 \right\}} \varphi^{q}(t,\eta,\xi  )dtd\eta d\xi \label{eq 2.47}\\
   & =  2^{(k+3s+2m+l-3)q -(3k +3s+2m+l-4)}
 \sum_{i=0}^{\inf(k,3s+2m+l-3)}\sigma _{i,i}^{\left[s-1\atop{ s+m-1\atop {s+m+l-1 \atop{\overline{\alpha_{s -} +\gamma  _{s+m+l-} } }}}\right]\times k}\cdot2^{-iq} \nonumber\\
   & =  2^{(k+3s+2m+l-3)q -(3k +3s+2m+l-4)}
 \sum_{i=0}^{\inf(k,3s+2m+l-3)}2\cdot\sigma _{i,i}^{\left[s-1\atop{ s+m-1\atop {s+m+l-1 \atop{\overline{\alpha_{s -}  } }}}\right]\times k}\cdot2^{-iq} \nonumber\\
   & = 2^{(k+3s+2m+l-3)q -(3k +3s+2m+l-3)}
 \sum_{i=0}^{\inf(k,3s+2m+l-3)}\big[\sigma _{i,i,i,i}^{\left[\alpha  | \alpha  \atop{ \beta  | \gamma \atop \gamma  | \beta   } \right]\times k}
+ \sigma _{i,i,i,i+1}^{\left[\alpha  | \alpha  \atop{ \beta  | \gamma \atop \gamma  | \beta   } \right]\times k}
+ \sigma _{i,i,i+1,i+1}^{\left[\alpha  | \alpha  \atop{ \beta  | \gamma \atop \gamma  | \beta   } \right]\times k}
+ \sigma _{i,i,i+1,i+2}^{\left[\alpha  | \alpha  \atop{ \beta  | \gamma \atop \gamma  | \beta   } \right]\times k} \big]\cdot2^{-iq} \nonumber \\
  & \int_{\mathbb{P}^{3}}  h_{4}^{q}(t,\eta,\xi )dtd\eta d\xi =\int_{\left\{(t,\eta ,\xi )\in\mathbb{P}^{3}\mid  h_{4}(t,\eta,\xi )\neq 0 \right\}} \varphi^{q}(t,\eta,\xi  )dtd\eta d\xi \label{eq 2.48}\\
   & =  2^{(k+3s+2m+l-3)q -(3k +3s+2m+l-3)}
 \sum_{i=0}^{\inf(k,3s+2m+l-3)}\sigma _{i,i}^{\left[s-1\atop{ s+m-1\atop {s+m+l-1 \atop{\overline{\alpha_{s -} +\beta_{s+m-}+ \gamma  _{s+m+l-} } }}}\right]\times k}\cdot2^{-iq} \nonumber\\
   & =  2^{(k+3s+2m+l-3)q -(3k +3s+2m+l-3)}
 \sum_{i=0}^{\inf(k,3s+2m+l-3)}4\cdot\sigma _{i,i}^{\left[s-1\atop{ s+m-1\atop {s+m+l-1 \atop{\overline{\alpha_{s -}  } }}}\right]\times k}\cdot2^{-iq} \nonumber\\
    & = 2^{(k+3s+2m+l-3)q -(3k +3s+2m+l-3)}
 \sum_{i=0}^{\inf(k,3s+2m+l-3)}\big[\sigma _{i,i,i,i}^{\left[\alpha  | \alpha  \atop{ \beta  | \gamma \atop \gamma  | \beta   } \right]\times k}
+ \sigma _{i,i,i,i+1}^{\left[\alpha  | \alpha  \atop{ \beta  | \gamma \atop \gamma  | \beta   } \right]\times k}
+ \sigma _{i,i,i+1,i+1}^{\left[\alpha  | \alpha  \atop{ \beta  | \gamma \atop \gamma  | \beta   } \right]\times k}
+ \sigma _{i,i,i+1,i+2}^{\left[\alpha  | \alpha  \atop{ \beta  | \gamma \atop \gamma  | \beta   } \right]\times k} \big]\cdot2^{-iq} \nonumber \\
 & \int_{\mathbb{P}^{3}}  \kappa _{1}^{q}(t,\eta,\xi )dtd\eta d\xi
  =\int_{\left\{(t,\eta ,\xi )\in\mathbb{P}^{3}\mid  \kappa _{1}(t,\eta,\xi )\neq 0 \right\}} \varphi_{2}^{q}(t,\eta,\xi  )dtd\eta d\xi \label{eq 2.49}\\
   & =  2^{(k+3s+2m+l-2)q -(3k +3s+2m+l-3)}
 \sum_{i=0}^{\inf(k,3s+2m+l-2)}\sigma _{i,i}^{\left[s-1\atop{ s+m-1\atop {s+m+l \atop{\overline{\alpha_{s -} +\beta_{s+m-} } }}}\right]\times k}\cdot2^{-iq} \nonumber\\
   & =  2^{(k+3s+2m+l-2)q -(3k +3s+2m+l-3)}
 \sum_{i=0}^{\inf(k,3s+2m+l-2)}2\cdot\sigma _{i,i}^{\left[s-1\atop{ s+m-1\atop {s+m+l \atop{\overline{\alpha_{s -}  } }}}\right]\times k}\cdot2^{-iq} \nonumber\\
    & = 2^{(k+3s+2m+l-2)q -(3k +3s+2m+l-3)}
 \sum_{i=0}^{\inf(k,3s+2m+l-2)}\big[\sigma _{i,i,i,i}^{\left[\gamma  \atop{ \alpha \atop  \beta   } \right]\times k}
+ \sigma _{i,i,i,i+1}^{\left[\gamma  \atop{ \alpha  \atop  \beta   } \right]\times k}
+ \sigma _{i-1,i,i,i}^{\left[\gamma  \atop{ \alpha  \atop \beta   } \right]\times k}
+ \sigma _{i-1,i,i,i+1}^{\left[\gamma  \atop{ \alpha \atop  \beta   } \right]\times k} \big]\cdot2^{-iq} \nonumber \\
 & \int_{\mathbb{P}^{3}}  \kappa _{2}^{q}(t,\eta,\xi )dtd\eta d\xi
  =\int_{\left\{(t,\eta ,\xi )\in\mathbb{P}^{3}\mid  \kappa _{2}(t,\eta,\xi )\neq 0 \right\}} \varphi_{2}^{q}(t,\eta,\xi  )dtd\eta d\xi \label{eq 2.50}\\
   & =  2^{(k+3s+2m+l-2)q -(3k +3s+2m+l-4)}
 \sum_{i=0}^{\inf(k,3s+2m+l-2)}\sigma _{i,i}^{\left[s-1\atop{ s+m-1\atop {s+m+l \atop{\overline{\alpha_{s -}  } }}}\right]\times k}\cdot2^{-iq} \nonumber\\
   & =  2^{(k+3s+2m+l-2)q -(3k +3s+2m+l-3)}
 \sum_{i=0}^{\inf(k,3s+2m+l-2)}2\cdot\sigma _{i,i}^{\left[s-1\atop{ s+m-1\atop {s+m+l \atop{\overline{\alpha_{s -}  } }}}\right]\times k}\cdot2^{-iq} \nonumber\\
    & = 2^{(k+3s+2m+l-2)q -(3k +3s+2m+l-3)}
 \sum_{i=0}^{\inf(k,3s+2m+l-2)}\big[\sigma _{i,i,i,i}^{\left[\gamma  \atop{ \alpha \atop  \beta   } \right]\times k}
+ \sigma _{i,i,i,i+1}^{\left[\gamma  \atop{ \alpha  \atop  \beta   } \right]\times k}
+ \sigma _{i-1,i,i,i}^{\left[\gamma  \atop{ \alpha  \atop \beta   } \right]\times k}
+ \sigma _{i-1,i,i,i+1}^{\left[\gamma  \atop{ \alpha \atop  \beta   } \right]\times k} \big]\cdot2^{-iq} \nonumber \\
 & \int_{\mathbb{P}^{3}}  g _{2}^{q}(t,\eta,\xi )dtd\eta d\xi
  =\int_{\left\{(t,\eta ,\xi )\in\mathbb{P}^{3}\mid  g _{2}(t,\eta,\xi )\neq 0 \right\}} \varphi_{3}^{q}(t,\eta,\xi  )dtd\eta d\xi \label{eq 2.51}\\
   & =  2^{(k+3s+2m+l-2)q -(3k +3s+2m+l-4)}
 \sum_{i=0}^{\inf(k,3s+2m+l-2)}\sigma _{i,i}^{\left[s-1\atop{ s+m\atop {s+m+l-1 \atop{\overline{\alpha_{s -}  } }}}\right]\times k}\cdot2^{-iq} \nonumber\\
   & =  2^{(k+3s+2m+l-2)q -(3k +3s+2m+l-3)}
 \sum_{i=0}^{\inf(k,3s+2m+l-2)}2\cdot\sigma _{i,i}^{\left[s-1\atop{ s+m\atop {s+m+l-1 \atop{\overline{\alpha_{s -}  } }}}\right]\times k}\cdot2^{-iq} \nonumber\\
    & = 2^{(k+3s+2m+l-2)q -(3k +3s+2m+l-3)}
 \sum_{i=0}^{\inf(k,3s+2m+l-2)}\big[\sigma _{i,i,i,i}^{\left[\beta  \atop{ \alpha \atop  \gamma   } \right]\times k}
+ \sigma _{i,i,i,i+1}^{\left[\beta  \atop{ \alpha  \atop  \gamma   } \right]\times k}
+ \sigma _{i-1,i,i,i}^{\left[\beta  \atop{ \alpha  \atop \gamma   } \right]\times k}
+ \sigma _{i-1,i,i,i+1}^{\left[\beta  \atop{ \alpha \atop  \gamma  } \right]\times k} \big]\cdot2^{-iq} \nonumber 
  \end{align}
  \end{lem}
  \begin{proof}
To prove \eqref{eq 2.44} we have by \eqref{eq 2.1} observing that $ g_{1}^{q}(t,\eta,\xi ) $ is constant on cosets of
$ \mathbb{P}_{k+s-1}\times\mathbb{P}_{k+s+m-1}\times\mathbb{P}_{k+s+m+l-1}$\\
 \begin{align*}
&\int_{\mathbb{P}^{3}} g_{1}^{q}(t,\eta,\xi  )dtd\eta d\xi = \int_{\left\{(t,\eta ,\xi )\in\mathbb{P}^{3}\mid  g_{1}(t,\eta,\xi )\neq 0 \right\}} \varphi _{1}^{q}(t,\eta,\xi  )dtd\eta d\xi \\
& = \sum_{(t,\eta,\xi  )\in \mathbb{P}/\mathbb{P}_{k+s-1}\times \mathbb{P}/\mathbb{P}_{k+s+m-1}\times \mathbb{P}/\mathbb{P}_{k+s+m+l-1}\atop
 {  r( D^{\left[\stackrel{s-1}{\stackrel{s+m}{s+m+l}}\right] \times k }(t,\eta ,\xi ) )
  =  r( D^{\left[\stackrel{s}{\stackrel{s+m}{s+m+l}}\right] \times k }(t,\eta ,\xi ) ) }}
  2^{\big(k+3s+2m+l-1-  r( D^{\left[\stackrel{s-1}{\stackrel{s+m}{s+m+l}}\right] \times k }(t,\eta ,\xi ) )\big)q  } \int_{\mathbb{P}_{k+s-1}}dt \int_{\mathbb{P}_{k+s+m-1}}d\eta \int_{\mathbb{P}_{k+s+m+l-1}}d\xi \\
 & = \sum_{i = 0}^{\inf(k,3s+2m+l-1)}\sum_{(t,\eta,\xi  )\in \mathbb{P}/\mathbb{P}_{k+s-1}\times \mathbb{P}/\mathbb{P}_{k+s+m-1}\times \mathbb{P}/\mathbb{P}_{k+s+m+l-1}\atop
 {  r( D^{\left[\stackrel{s-1}{\stackrel{s+m}{s+m+l}}\right] \times k }(t,\eta ,\xi ) )
  =  r( D^{\left[\stackrel{s}{\stackrel{s+m}{s+m+l}}\right] \times k }(t,\eta ,\xi ) ) }= i}
  2^{(k+3s+2m+l-1-i)q  } \int_{\mathbb{P}_{k+s-1}}dt \int_{\mathbb{P}_{k+s+m-1}}d\eta \int_{\mathbb{P}_{k+s+m+l-1}}d\xi \\
  & = 2^{(k+3s+2m+l-1)q -(3k +3s+2m+l-3)}
 \sum_{i=0}^{\inf(k,3s+2m+l-1)}\sigma _{i,i}^{\left[s-1\atop{ s+m\atop {s+m+l \atop{\overline{\alpha_{s -}  } }}}\right]\times k}\cdot2^{-iq} \\
 &  = 2^{(k+3s+2m+l-1)q -(3k +3s+2m+l-3)}
 \sum_{i=0}^{\inf(k,3s+2m+l-1)}\big[\sigma _{i,i,i,i}^{\left[\beta | \gamma  \atop{ \gamma | \beta \atop \alpha | \alpha  } \right]\times k}
+  \sigma _{i-1,i,i,i}^{\left[\beta | \gamma  \atop{ \gamma | \beta \atop \alpha | \alpha  } \right]\times k}
+  \sigma _{i-1,i-1,i,i}^{\left[\beta | \gamma  \atop{ \gamma | \beta \atop \alpha | \alpha  } \right]\times k}
+  \sigma _{i-2,i-1,i,i}^{\left[\beta | \gamma  \atop{ \gamma | \beta \atop \alpha | \alpha  } \right]\times k}\big]\cdot2^{-iq} \nonumber  
 \end{align*}  
The other integral equalities are obtained in a similar way.  
  \end{proof}
    \begin{lem}
\label{lem 2.8}
Let $ (t,\eta ,\xi ) \in  \mathbb{P}\times \mathbb{P}\times \mathbb{P},  $ and q a rational integer $\geq 2,$ then we have for $1\leq i\leq q-1$\\
 \begin{align}
\displaystyle
  & \int_{\mathbb{P}^{3}}  h_{1}^{i}(t,\eta,\xi )\cdot h_{2}^{q-i}(t,\eta,\xi ) dtd\eta d\xi
   =\int_{\left\{(t,\eta ,\xi )\in\mathbb{P}^{3}\mid  h_{1}(t,\eta,\xi )\cdot h_{2}(t,\eta,\xi )\neq 0 \right\}} \varphi^{q}(t,\eta,\xi  )dtd\eta d\xi \label{eq 2.52}\\
   & =  2^{(k+3s+2m+l-3)q -(3k +3s+2m+l-4)}
 \sum_{i=0}^{\inf(k,3s+2m+l-3)}\sigma _{i,i,i}^{\left[s-1\atop{ s+m-1\atop {s+m+l-1 \atop{\overline{\alpha_{s -} \atop{\beta_{s+m -} } }}}}\right]\times k} \cdot2^{-iq} \nonumber\\
  & = 2^{(k+3s+2m+l-3)q -(3k +3s+2m+l-3)} \sum_{i=0}^{\inf(k,3s+2m+l-3)}  \big[\sigma _{i,i,i,i}^{\left[\alpha  | \beta  \atop{ \beta  | \alpha  \atop \gamma  | \gamma  } \right]\times k}
 + \sigma _{i,i,i,i+1}^{\left[\alpha  | \beta  \atop{ \beta | \alpha \atop \gamma  | \gamma  } \right]\times k} \big]\cdot2^{-iq} \nonumber
  \end{align}

  \begin{align}
\displaystyle
  & \int_{\mathbb{P}^{3}}  h_{1}^{i}(t,\eta,\xi )\cdot h_{3}^{q-i}(t,\eta,\xi ) dtd\eta d\xi
   =\int_{\left\{(t,\eta ,\xi )\in\mathbb{P}^{3}\mid  h_{1}(t,\eta,\xi )\cdot h_{3}(t,\eta,\xi )\neq 0 \right\}} \varphi^{q}(t,\eta,\xi  )dtd\eta d\xi \label{eq 2.53}\\
   & =  2^{(k+3s+2m+l-3)q -(3k +3s+2m+l-4)}
 \sum_{i=0}^{\inf(k,3s+2m+l-3)}\sigma _{i,i,i}^{\left[s-1\atop{ s+m-1\atop {s+m+l-1 \atop{\overline{\alpha_{s -} \atop{\gamma _{s+m+l -} } }}}}\right]\times k} \cdot2^{-iq} \nonumber\\
 & =  2^{(k+3s+2m+l-3)q -(3k +3s+2m+l-3)}
 \sum_{i=0}^{\inf(k,3s+2m+l-3)}\big[\sigma _{i,i,i,i}^{\left[\alpha  | \gamma  \atop{ \gamma  | \alpha  \atop \beta  | \beta  } \right]\times k}
 + \sigma _{i,i,i,i+1}^{\left[\alpha  | \gamma  \atop{ \gamma  | \alpha \atop \beta  | \beta  } \right]\times k} \big]\cdot2^{-iq} \nonumber  
 \end{align}
 
  \begin{align}
\displaystyle
  & \int_{\mathbb{P}^{3}}  h_{1}^{i}(t,\eta,\xi )\cdot h_{4}^{q-i}(t,\eta,\xi ) dtd\eta d\xi
   =\int_{\left\{(t,\eta ,\xi )\in\mathbb{P}^{3}\mid  h_{1}(t,\eta,\xi )\cdot h_{4}(t,\eta,\xi )\neq 0 \right\}} \varphi^{q}(t,\eta,\xi  )dtd\eta d\xi \label{eq 2.54}\\
   & =  2^{(k+3s+2m+l-3)q -(3k +3s+2m+l-3)}
 \sum_{i=0}^{\inf(k,3s+2m+l-3)}\sigma _{i,i,i}^{\left[s-1\atop{ s+m-1\atop {s+m+l-1 \atop{\overline{\alpha_{s -} \atop{\beta _{s+m-} + \gamma _{s+m+l -} } }}}}\right]\times k} \cdot2^{-iq} \nonumber\\
 & = 2^{(k+3s+2m+l-3)q -(3k +3s+2m+l-3)}
 \sum_{i=0}^{\inf(k,3s+2m+l-3)}2\cdot\sigma _{i,i,i}^{\left[s-1\atop{ s+m-1\atop {s+m+l-1 \atop{\overline{\alpha_{s -} \atop{\beta _{s+m-} } }}}}\right]\times k} \cdot2^{-iq} \nonumber\\
  & = 2^{(k+3s+2m+l-3)q -(3k +3s+2m+l-3)} \sum_{i=0}^{\inf(k,3s+2m+l-3)}  \big[\sigma _{i,i,i,i}^{\left[\alpha  | \beta  \atop{ \beta  | \alpha  \atop \gamma  | \gamma  } \right]\times k}
 + \sigma _{i,i,i,i+1}^{\left[\alpha  | \beta  \atop{ \beta | \alpha \atop \gamma  | \gamma  } \right]\times k} \big]\cdot2^{-iq} \nonumber
 \end{align}

  \begin{align}
\displaystyle
  & \int_{\mathbb{P}^{3}}  h_{2}^{i}(t,\eta,\xi )\cdot h_{3}^{q-i}(t,\eta,\xi ) dtd\eta d\xi
   =\int_{\left\{(t,\eta ,\xi )\in\mathbb{P}^{3}\mid  h_{2}(t,\eta,\xi )\cdot h_{3}(t,\eta,\xi )\neq 0 \right\}} \varphi^{q}(t,\eta,\xi  )dtd\eta d\xi \label{eq 2.55}\\
   & =  2^{(k+3s+2m+l-3)q -(3k +3s+2m+l-3)}
 \sum_{i=0}^{\inf(k,3s+2m+l-3)}\sigma _{i,i,i}^{\left[s-1\atop{ s+m-1\atop {s+m+l-1 \atop{\overline{\alpha_{s -}+\beta _{s+m-} \atop{\beta _{s+m-} + \gamma _{s+m+l -} } }}}}\right]\times k} \cdot2^{-iq} \nonumber\\
 & = 2^{(k+3s+2m+l-3)q -(3k +3s+2m+l-3)}
 \sum_{i=0}^{\inf(k,3s+2m+l-3)}2\cdot\sigma _{i,i,i}^{\left[s-1\atop{ s+m-1\atop {s+m+l-1 \atop{\overline{\alpha_{s -} \atop{\beta _{s+m-} } }}}}\right]\times k} \cdot2^{-iq} \nonumber\\
  & = 2^{(k+3s+2m+l-3)q -(3k +3s+2m+l-3)} \sum_{i=0}^{\inf(k,3s+2m+l-3)}  \big[\sigma _{i,i,i,i}^{\left[\alpha  | \beta  \atop{ \beta  | \alpha  \atop \gamma  | \gamma  } \right]\times k}
 + \sigma _{i,i,i,i+1}^{\left[\alpha  | \beta  \atop{ \beta | \alpha \atop \gamma  | \gamma  } \right]\times k} \big]\cdot2^{-iq} \nonumber
 \end{align}

   \begin{align}
\displaystyle
  & \int_{\mathbb{P}^{3}}  h_{2}^{i}(t,\eta,\xi )\cdot h_{4}^{q-i}(t,\eta,\xi ) dtd\eta d\xi
   =\int_{\left\{(t,\eta ,\xi )\in\mathbb{P}^{3}\mid  h_{2}(t,\eta,\xi )\cdot h_{4}(t,\eta,\xi )\neq 0 \right\}} \varphi^{q}(t,\eta,\xi  )dtd\eta d\xi \label{eq 2.56}\\
   & =  2^{(k+3s+2m+l-3)q -(3k +3s+2m+l-3)}
 \sum_{i=0}^{\inf(k,3s+2m+l-3)}\sigma _{i,i,i}^{\left[s-1\atop{ s+m-1\atop {s+m+l-1 \atop{\overline{\alpha_{s -}+\beta _{s+m-} \atop{ \gamma _{s+m+l -} } }}}}\right]\times k} \cdot2^{-iq} \nonumber\\
 & = 2^{(k+3s+2m+l-3)q -(3k +3s+2m+l-3)}
 \sum_{i=0}^{\inf(k,3s+2m+l-3)}2\cdot\sigma _{i,i,i}^{\left[s-1\atop{ s+m-1\atop {s+m+l-1 \atop{\overline{\alpha_{s -} \atop{\gamma  _{s+m+l-} } }}}}\right]\times k} \cdot2^{-iq} \nonumber\\
  & =  2^{(k+3s+2m+l-3)q -(3k +3s+2m+l-3)}
 \sum_{i=0}^{\inf(k,3s+2m+l-3)}\big[\sigma _{i,i,i,i}^{\left[\alpha  | \gamma  \atop{ \gamma  | \alpha  \atop \beta  | \beta  } \right]\times k}
 + \sigma _{i,i,i,i+1}^{\left[\alpha  | \gamma  \atop{ \gamma  | \alpha \atop \beta  | \beta  } \right]\times k} \big]\cdot2^{-iq} \nonumber  
 \end{align}
   \begin{align}
   \displaystyle
  & \int_{\mathbb{P}^{3}}  h_{3}^{i}(t,\eta,\xi )\cdot h_{4}^{q-i}(t,\eta,\xi ) dtd\eta d\xi
   =\int_{\left\{(t,\eta ,\xi )\in\mathbb{P}^{3}\mid  h_{3}(t,\eta,\xi )\cdot h_{4}(t,\eta,\xi )\neq 0 \right\}} \varphi^{q}(t,\eta,\xi  )dtd\eta d\xi \label{eq 2.57}\\
   & =  2^{(k+3s+2m+l-3)q -(3k +3s+2m+l-3)}
 \sum_{i=0}^{\inf(k,3s+2m+l-3)}\sigma _{i,i,i}^{\left[s-1\atop{ s+m-1\atop {s+m+l-1 \atop{\overline{\alpha_{s -}+\gamma  _{s+m+l-} \atop{ \beta  _{s+m -} } }}}}\right]\times k} \cdot2^{-iq} \nonumber\\
  & = 2^{(k+3s+2m+l-3)q -(3k +3s+2m+l-3)}
 \sum_{i=0}^{\inf(k,3s+2m+l-3)}2\cdot\sigma _{i,i,i}^{\left[s-1\atop{ s+m-1\atop {s+m+l-1 \atop{\overline{\alpha_{s -} \atop{\beta  _{s+m-} } }}}}\right]\times k} \cdot2^{-iq} \nonumber\\
 & = 2^{(k+3s+2m+l-3)q -(3k +3s+2m+l-3)} \sum_{i=0}^{\inf(k,3s+2m+l-3)}  \big[\sigma _{i,i,i,i}^{\left[\alpha  | \beta  \atop{ \beta  | \alpha  \atop \gamma  | \gamma  } \right]\times k}
 + \sigma _{i,i,i,i+1}^{\left[\alpha  | \beta  \atop{ \beta | \alpha \atop \gamma  | \gamma  } \right]\times k} \big]\cdot2^{-iq} \nonumber
  \end{align}

   \begin{align}
   \displaystyle
  & \int_{\mathbb{P}^{3}}  \kappa _{1}^{i}(t,\eta,\xi )\cdot \kappa _{2}^{q-i}(t,\eta,\xi ) dtd\eta d\xi
   =\int_{\left\{(t,\eta ,\xi )\in\mathbb{P}^{3}\mid  \kappa _{1}(t,\eta,\xi )\cdot \kappa _{2}(t,\eta,\xi )\neq 0 \right\}} \varphi_{2}^{q}(t,\eta,\xi  )dtd\eta d\xi \label{eq 2.58}\\
   & =  2^{(k+3s+2m+l-2)q -(3k +3s+2m+l-3)}
 \sum_{i=0}^{\inf(k,3s+2m+l-2)}\sigma _{i,i,i}^{\left[s-1\atop{ s+m-1\atop {s+m+l \atop{\overline{\alpha_{s -}+\beta   _{s+m-} \atop{ \alpha  _{s -} } }}}}\right]\times k} \cdot2^{-iq} \nonumber\\
 & = 2^{(k+3s+2m+l-2)q -(3k +3s+2m+l-3)} \sum_{i=0}^{\inf(k,3s+2m+l-2)}  \big[\sigma _{i,i,i,i}^{\left[\gamma  | \gamma  \atop{ \alpha  | \beta   \atop \beta  | \alpha  } \right]\times k}
 + \sigma _{i-1,i,i,i}^{\left[\gamma  | \gamma   \atop{ \alpha  | \beta  \atop \beta  | \alpha  } \right]\times k} \big]\cdot2^{-iq} \nonumber
  \end{align}
   \end{lem}
 \begin{proof} 
Proceed as in the proof of \eqref{eq 2.44} using the results in Lemma \ref{lem 2.3} and Lemma \ref{lem 2.6}.
 \end{proof} 
  \begin{lem}
\label{lem 2.9}
Let $ (t,\eta ,\xi ) \in  \mathbb{P}\times \mathbb{P}\times \mathbb{P}  $ and q a rational integer $\geq 3,$ then we have for any
integers i,j,r such that
\begin{align*}
\displaystyle
   \begin{cases}
   i+j+r = q   \\
   1\leq i,j,r\leq q-2
 \end{cases}
\end{align*}
the following integral identities : 
\begin{align}
\displaystyle
  & \int_{\mathbb{P}^{3}}  h_{1}^{i}(t,\eta,\xi )\cdot h_{2}^{j}(t,\eta,\xi )\cdot h_{3}^{r}(t,\eta,\xi )  dtd\eta d\xi
   =\int_{\left\{(t,\eta ,\xi )\in\mathbb{P}^{3}\mid  h_{1}(t,\eta,\xi )\cdot h_{2}(t,\eta,\xi )\cdot h_{3}(t,\eta,\xi )\neq 0 \right\}} \varphi^{q}(t,\eta,\xi  )dtd\eta d\xi \label{eq 2.59}\\
&  = \int_{\mathbb{P}^{3}}  h_{1}^{i}(t,\eta,\xi )\cdot h_{3}^{j}(t,\eta,\xi )\cdot h_{4}^{r}(t,\eta,\xi )  dtd\eta d\xi
   =\int_{\left\{(t,\eta ,\xi )\in\mathbb{P}^{3}\mid  h_{1}(t,\eta,\xi )\cdot h_{3}(t,\eta,\xi )\cdot h_{4}(t,\eta,\xi )\neq 0 \right\}} \varphi^{q}(t,\eta,\xi  )dtd\eta d\xi \nonumber  \\
  & =  \int_{\mathbb{P}^{3}}  h_{2}^{i}(t,\eta,\xi )\cdot h_{3}^{j}(t,\eta,\xi )\cdot h_{4}^{r}(t,\eta,\xi )  dtd\eta d\xi
   =\int_{\left\{(t,\eta ,\xi )\in\mathbb{P}^{3}\mid  h_{2}(t,\eta,\xi )\cdot h_{3}(t,\eta,\xi )\cdot h_{4}(t,\eta,\xi )\neq 0 \right\}} \varphi^{q}(t,\eta,\xi  )dtd\eta d\xi \nonumber  \\
 &  = \int_{\mathbb{P}^{3}}  h_{1}^{i}(t,\eta,\xi )\cdot h_{2}^{j}(t,\eta,\xi )\cdot h_{4}^{r}(t,\eta,\xi )  dtd\eta d\xi
   =\int_{\left\{(t,\eta ,\xi )\in\mathbb{P}^{3}\mid  h_{1}(t,\eta,\xi )\cdot h_{2}(t,\eta,\xi )\cdot h_{4}(t,\eta,\xi )\neq 0 \right\}} \varphi^{q}(t,\eta,\xi  )dtd\eta d\xi \nonumber  \\
   & =  2^{(k+3s+2m+l-3)q -(3k +3s+2m+l-3)}
 \sum_{i=0}^{\inf(k,3s+2m+l-3)} \sigma _{i,i,i,i}^{\left[s-1\atop{ s+m-1\atop {s+m+l-1 \atop{\overline{\alpha_{s -} \atop{\beta_{s+m -} \atop \gamma_{s+m+l -}} }}}}\right]\times k} \cdot2^{-iq} \nonumber\\
  & =  2^{(k+3s+2m+l-3)q -(3k +3s+2m+l-3)} \sum_{i=0}^{\inf(k,3s+2m+l-3)} \sigma _{i,i,i,i}^{\left[\alpha \atop{ \beta \atop \gamma } \right]\times k}\cdot2^{-iq} \nonumber
  \end{align}
  \end{lem}
\begin{proof} 
Proceed as in the proof of \eqref{eq 2.44} using the results in Lemma \ref{lem 2.4}.
 \end{proof} 
  \begin{lem}
\label{lem 2.10}
Let $ (t,\eta ,\xi ) \in  \mathbb{P}\times \mathbb{P}\times \mathbb{P},  $ and q a rational integer $\geq 4,$ then we have for any
integers i,j,r,p such that
\begin{align*}
\displaystyle
   \begin{cases}
   i+j+r+p = q   \\
   1\leq i,j,r,p\leq q-3
 \end{cases}
\end{align*} 
the following integral identities : 
\begin{align}
\displaystyle
  & \int_{\mathbb{P}^{3}}  h_{1}^{i}(t,\eta,\xi )\cdot h_{2}^{j}(t,\eta,\xi )\cdot h_{3}^{r}(t,\eta,\xi ) \cdot h_{4}^{p}(t,\eta,\xi ) dtd\eta d\xi
   =\int_{\left\{(t,\eta ,\xi )\in\mathbb{P}^{3}\mid  h_{1}(t,\eta,\xi )\cdot h_{2}(t,\eta,\xi )\cdot h_{3}(t,\eta,\xi )\cdot h_{4}(t,\eta,\xi )\neq 0 \right\}} \varphi^{q}(t,\eta,\xi  )dtd\eta d\xi \label{eq 2.60}\\
 & =  2^{(k+3s+2m+l-3)q -(3k +3s+2m+l-3)}
 \sum_{i=0}^{\inf(k,3s+2m+l-3)} \sigma _{i,i,i,i}^{\left[s-1\atop{ s+m-1\atop {s+m+l-1 \atop{\overline{\alpha_{s -} \atop{\beta_{s+m -} \atop \gamma_{s+m+l -}} }}}}\right]\times k} \cdot2^{-iq} \nonumber\\
  & =  2^{(k+3s+2m+l-3)q -(3k +3s+2m+l-3)} \sum_{i=0}^{\inf(k,3s+2m+l-3)} \sigma _{i,i,i,i}^{\left[\alpha \atop{ \beta \atop \gamma } \right]\times k}\cdot2^{-iq} \nonumber
  \end{align} 
\end{lem}
\begin{proof} 
Proceed as in the proof of \eqref{eq 2.44} using the results in Lemma \ref{lem 2.5}.
 \end{proof} 
 \subsection{Some rank formulas for partitions of triple persymmetric matrices}
  \label{subsec 6}
  \begin{lem}
\label{lem 2.11}
 We can write : \\
 \begin{equation}
 \label{eq 2.61}
   g_{1}(t,\eta,\xi  ) =  \kappa _{1}(t,\eta,\xi  ) + \kappa _{2}(t,\eta,\xi  )
\end{equation}
Let q be a rational integer $\geq 2, $ then we obtain  from the binomial formula   \\
\begin{align}
\displaystyle
 g_{1}^q = ( \kappa _{1}+\kappa _{2})^q   =  \kappa _{1}^q  + \kappa _{2}^q    +   \sum_{i=1}^{q-1}\binom{q}{i} \kappa _{1}^i  \kappa  _{2}^{q-i}  \label{eq 2.62}  
\end{align}
  
By integrating  \eqref{eq 2.62}  on the unit interval of $\mathbb{K}^3 $ we get \\
\begin{align}
\displaystyle
& \int_{\mathbb{P}^3} g_{1}^q(t,\eta,\xi )dtd\eta d\xi
 =  2\int_{\{(t,\eta ,\xi )\in \mathbb{P}^3 \mid \kappa _{1}(t,\eta,\xi ) \neq 0 \}}\varphi_{2} ^q(t,\eta ,\xi )dtd\eta d\xi 
 + (2^q -2)\int_{\{(t,\eta ,\xi )\in \mathbb{P}^3 \mid \kappa _{1}(t,\eta,\xi )\cdot \kappa _{2}(t,\eta,\xi ) \neq 0 \}}\varphi_{2} ^q(t,\eta ,\xi )dtd\eta d\xi  \label{eq 2.63} 
\end{align}
\end{lem}
\begin{proof}
From Lemma \ref{lem 2.1} we get easily \eqref{eq 2.61}, indeed we have \\
\begin{align*}
&  g_{1}(t,\eta,\xi  ) =  \sum_{deg Y\leq k-1}\sum_{deg Z = s-1}E(tYZ)\sum_{deg U \leq s+m-1}E(\eta YU) \sum_{deg V \leq s+m+l-1}E(\xi YV)\\
& =   \sum_{deg Y\leq k-1}E(tYT^{s-1})E(\eta YT^{s+m-1}) \sum_{deg Z\leq  s-2}E(tYZ)\sum_{deg U \leq s+m-2}E(\eta YU) \sum_{deg V \leq s+m+l-1}E(\xi YV)   \\
& +   \sum_{deg Y\leq k-1}E(tYT^{s-1}) \sum_{deg Z\leq  s-2}E(tYZ)\sum_{deg U \leq s+m-2}E(\eta YU) \sum_{deg V \leq s+m+l-1}E(\xi  YV)\\
& =  \kappa _{1}(t,\eta,\xi  ) + \kappa _{2}(t,\eta,\xi).
\end{align*}
To prove \eqref{eq 2.63} we proceed as follows  \\

 From \eqref{eq 2.49},\eqref{eq 2.50} we get \\
  \begin{align}
\label{eq 2.64}
\displaystyle
\int_{\mathbb{P}^{3}}\big( \kappa _{1}^q(t,\eta,\xi  ) + \kappa _{2}^q(t,\eta,\xi  ) )\big)dtd\eta d\xi  =  2\int_{\{(t,\eta ,\xi )\in \mathbb{P}^3 \mid \kappa _{1}(t,\eta,\xi ) \neq 0 \}}\varphi_{2} ^q(t,\eta ,\xi )dtd\eta d\xi \\
& \nonumber
\end{align}

By \eqref{eq 2.58} we get \\

   \begin{align}
\label{eq 2.65}
\displaystyle
& \int_{\mathbb{P}^{3}}\sum_{i=1}^{q-1}\binom{q}{i} \kappa _{1}^i  \kappa _{2}^{q-i}dtd\eta d\xi
 = \sum_{i=1}^{q-1}\binom{q}{i}\int_{\mathbb{P}^{3}} \kappa _{1}^i  \kappa _{2}^{q-i}dtd\eta d\xi
 = (2^q -2)\int_{\{(t,\eta ,\xi )\in \mathbb{P}^3 \mid \kappa _{1}(t,\eta,\xi )\cdot \kappa _{2}(t,\eta,\xi ) \neq 0 \}}\varphi_{2} ^q(t,\eta ,\xi )dtd\eta d\xi, \\
  & \nonumber
 \end{align} 
 \eqref{eq 2.63} follows now from \eqref{eq 2.64} and \eqref{eq 2.65}.
\end{proof} 
  
   \begin{lem}
\label{lem 2.12}
We have the following rank formula :\\
\begin{align}
 \label{eq 2.66}
 \displaystyle 
\boxed{ \sigma _{i-1,i-1,i,i}^{\left[ \gamma  \atop{  \beta \atop  \alpha  } \right]\times k}
+  \sigma _{i-2,i-1,i,i}^{\left[ \gamma  \atop{  \beta \atop  \alpha  } \right]\times k}=
2\cdot\big[  \sigma _{i-1,i-1,i-1,i}^{\left[ \gamma  \atop{  \alpha  \atop  \beta  } \right]\times k}
+  \sigma _{i-2,i-1,i-1,i}^{\left[ \gamma  \atop{  \alpha \atop  \beta  } \right]\times k}  \big] \quad \text{for $0\leq i\leq \inf(k,3s+2m+l-1)$}}
\end{align}
 \end{lem}
 
 \begin{proof}
We obtain by  combining  \eqref{eq 2.44}, \eqref{eq 2.49},\eqref{eq 2.58} and \eqref{eq 2.62}.     \\
\begin{align*}
\displaystyle
& \int_{\mathbb{P}^3} g_{1}^q(t,\eta,\xi )dtd\eta d\xi
 =  2\int_{\{(t,\eta ,\xi )\in \mathbb{P}^3 \mid \kappa _{1}(t,\eta,\xi ) \neq 0 \}}\varphi_{2} ^q(t,\eta ,\xi )dtd\eta d\xi 
 + (2^q -2)\int_{\{(t,\eta ,\xi )\in \mathbb{P}^3 \mid \kappa _{1}(t,\eta,\xi )\cdot \kappa _{2}(t,\eta,\xi ) \neq 0 \}}\varphi_{2} ^q(t,\eta ,\xi )dtd\eta d\xi \\
 & \Longleftrightarrow \\
  &  2^{(k+3s+2m+l-1)q -(3k +3s+2m+l-3)}
 \sum_{i=0}^{\inf(k,3s+2m+l-1)}\big[\sigma _{i,i,i,i}^{\left[\beta | \gamma  \atop{ \gamma | \beta \atop \alpha | \alpha  } \right]\times k}
+  \sigma _{i-1,i,i,i}^{\left[\beta | \gamma  \atop{ \gamma | \beta \atop \alpha | \alpha  } \right]\times k}
+  \sigma _{i-1,i-1,i,i}^{\left[\beta | \gamma  \atop{ \gamma | \beta \atop \alpha | \alpha  } \right]\times k}
+  \sigma _{i-2,i-1,i,i}^{\left[\beta | \gamma  \atop{ \gamma | \beta \atop \alpha | \alpha  } \right]\times k}\big]\cdot2^{-iq}  \\
& = 2\cdot2^{(k+3s+2m+l-2)q -(3k +3s+2m+l-3)}
 \sum_{i=0}^{\inf(k,3s+2m+l-2)}\big[\sigma _{i,i,i,i}^{\left[\gamma  \atop{ \alpha \atop  \beta   } \right]\times k}
+ \sigma _{i,i,i,i+1}^{\left[\gamma  \atop{ \alpha  \atop  \beta   } \right]\times k}
+ \sigma _{i-1,i,i,i}^{\left[\gamma  \atop{ \alpha  \atop \beta   } \right]\times k}
+ \sigma _{i-1,i,i,i+1}^{\left[\gamma  \atop{ \alpha \atop  \beta   } \right]\times k} \big]\cdot2^{-iq} \\
& + (2^q -2)\cdot 2^{(k+3s+2m+l-2)q -(3k +3s+2m+l-3)} \sum_{i=0}^{\inf(k,3s+2m+l-2)}  \big[\sigma _{i,i,i,i}^{\left[\gamma  | \gamma  \atop{ \alpha  | \beta   \atop \beta  | \alpha  } \right]\times k}
 + \sigma _{i-1,i,i,i}^{\left[\gamma  | \gamma   \atop{ \alpha  | \beta  \atop \beta  | \alpha  } \right]\times k} \big]\cdot2^{-iq} 
\end{align*}
\newpage

\begin{align*} 
 & \text{Thus we get} \quad  \eqref{eq 2.62} \Longleftrightarrow  \\
  & \sum_{i=0}^{\inf(k,3s+2m+l-1)}\big[\sigma _{i,i,i,i}^{\left[\beta | \gamma  \atop{ \gamma | \beta \atop \alpha | \alpha  } \right]\times k}
+  \sigma _{i-1,i,i,i}^{\left[\beta | \gamma  \atop{ \gamma | \beta \atop \alpha | \alpha  } \right]\times k}
+  \sigma _{i-1,i-1,i,i}^{\left[\beta | \gamma  \atop{ \gamma | \beta \atop \alpha | \alpha  } \right]\times k}
+  \sigma _{i-2,i-1,i,i}^{\left[\beta | \gamma  \atop{ \gamma | \beta \atop \alpha | \alpha  } \right]\times k}\big]\cdot2^{-iq}  \\
& = 2\cdot2^{-iq}
 \sum_{i=0}^{\inf(k,3s+2m+l-2)}\big[\sigma _{i,i,i,i}^{\left[\gamma  \atop{ \alpha \atop  \beta   } \right]\times k}
+ \sigma _{i,i,i,i+1}^{\left[\gamma  \atop{ \alpha  \atop  \beta   } \right]\times k}
+ \sigma _{i-1,i,i,i}^{\left[\gamma  \atop{ \alpha  \atop \beta   } \right]\times k}
+ \sigma _{i-1,i,i,i+1}^{\left[\gamma  \atop{ \alpha \atop  \beta   } \right]\times k} \big]\cdot2^{-iq} \\
& + (2^q -2)\cdot 2^{-q} \sum_{i=0}^{\inf(k,3s+2m+l-2)}  \big[\sigma _{i,i,i,i}^{\left[\gamma  | \gamma  \atop{ \alpha  | \beta   \atop \beta  | \alpha  } \right]\times k}
 + \sigma _{i-1,i,i,i}^{\left[\gamma  | \gamma   \atop{ \alpha  | \beta  \atop \beta  | \alpha  } \right]\times k} \big]\cdot2^{-iq} \\
\end{align*}
\underline{The case $ k\leq 3s+2m+l-2 $}
\begin{align*} 
 &  \eqref{eq 2.62} \Longleftrightarrow  \\
  & \sum_{i=0}^{k}\big[\sigma _{i,i,i,i}^{\left[\beta | \gamma  \atop{ \gamma | \beta \atop \alpha | \alpha  } \right]\times k}
+  \sigma _{i-1,i,i,i}^{\left[\beta | \gamma  \atop{ \gamma | \beta \atop \alpha | \alpha  } \right]\times k}
+  \sigma _{i-1,i-1,i,i}^{\left[\beta | \gamma  \atop{ \gamma | \beta \atop \alpha | \alpha  } \right]\times k}
+  \sigma _{i-2,i-1,i,i}^{\left[\beta | \gamma  \atop{ \gamma | \beta \atop \alpha | \alpha  } \right]\times k}\big]\cdot2^{-iq}  \\
& = 2\cdot2^{-iq}
 \sum_{i=0}^{k}\big[\sigma _{i,i,i,i}^{\left[\gamma  \atop{ \alpha \atop  \beta   } \right]\times k}
+ \sigma _{i,i,i,i+1}^{\left[\gamma  \atop{ \alpha  \atop  \beta   } \right]\times k}
+ \sigma _{i-1,i,i,i}^{\left[\gamma  \atop{ \alpha  \atop \beta   } \right]\times k}
+ \sigma _{i-1,i,i,i+1}^{\left[\gamma  \atop{ \alpha \atop  \beta   } \right]\times k} \big]\cdot2^{-iq} 
 + (2^q -2)\cdot 2^{-q} \sum_{i=0}^{k}  \big[\sigma _{i,i,i,i}^{\left[\gamma  | \gamma  \atop{ \alpha  | \beta   \atop \beta  | \alpha  } \right]\times k}
 + \sigma _{i-1,i,i,i}^{\left[\gamma  | \gamma   \atop{ \alpha  | \beta  \atop \beta  | \alpha  } \right]\times k} \big]\cdot2^{-iq} \\
 & \Longleftrightarrow \\
   & \sum_{i=0}^{k}\big[\sigma _{i,i,i,i}^{\left[\beta | \gamma  \atop{ \gamma | \beta \atop \alpha | \alpha  } \right]\times k}
+  \sigma _{i-1,i,i,i}^{\left[\beta | \gamma  \atop{ \gamma | \beta \atop \alpha | \alpha  } \right]\times k}
+  \sigma _{i-1,i-1,i,i}^{\left[\beta | \gamma  \atop{ \gamma | \beta \atop \alpha | \alpha  } \right]\times k}
+  \sigma _{i-2,i-1,i,i}^{\left[\beta | \gamma  \atop{ \gamma | \beta \atop \alpha | \alpha  } \right]\times k}\big]\cdot2^{-iq}  \\
& =  \sum_{i=1}^{k+1}2\cdot\big[\sigma _{i-1,i-1,i-1,i-1}^{\left[\gamma  \atop{ \alpha \atop  \beta   } \right]\times k}
+ \sigma _{i-1,i-1,i-1,i}^{\left[\gamma  \atop{ \alpha  \atop  \beta   } \right]\times k}
+ \sigma _{i-2,i-1,i-1,i-1}^{\left[\gamma  \atop{ \alpha  \atop \beta   } \right]\times k}
+ \sigma _{i-2,i-1,i-1,i}^{\left[\gamma  \atop{ \alpha \atop  \beta   } \right]\times k} \big]\cdot2^{-iq} \\
& +  \sum_{i=0}^{k}  \big[\sigma _{i,i,i,i}^{\left[\gamma  | \gamma  \atop{ \alpha  | \beta   \atop \beta  | \alpha  } \right]\times k}
 + \sigma _{i-1,i,i,i}^{\left[\gamma  | \gamma   \atop{ \alpha  | \beta  \atop \beta  | \alpha  } \right]\times k} \big]\cdot2^{-iq} 
  -  \sum_{i=1}^{k+1}2\cdot \big[\sigma _{i-1,i-1,i-1,i-1}^{\left[\gamma  | \gamma  \atop{ \alpha  | \beta   \atop \beta  | \alpha  } \right]\times k}
 + \sigma _{i-2,i-1,i-1,i-1}^{\left[\gamma  | \gamma   \atop{ \alpha  | \beta  \atop \beta  | \alpha  } \right]\times k} \big]\cdot2^{-iq} \\
 & \Longrightarrow \\
& \sum_{i=2}^{k}\Big( \sigma _{i-1,i-1,i,i}^{\left[ \gamma  \atop{  \beta \atop  \alpha  } \right]\times k}
+  \sigma _{i-2,i-1,i,i}^{\left[ \gamma  \atop{  \beta \atop  \alpha  } \right]\times k} -
2\cdot\big[  \sigma _{i-1,i-1,i-1,i}^{\left[ \gamma  \atop{  \alpha  \atop  \beta  } \right]\times k}
+  \sigma _{i-2,i-1,i-1,i}^{\left[ \gamma  \atop{  \alpha \atop  \beta  } \right]\times k}  \big] \Big)\cdot2^{-iq} = 0 \quad \text{for all $ q \geq 2 $}\\
 & \Longrightarrow \\
&  \sigma _{i-1,i-1,i,i}^{\left[ \gamma  \atop{  \beta \atop  \alpha  } \right]\times k}
+  \sigma _{i-2,i-1,i,i}^{\left[ \gamma  \atop{  \beta \atop  \alpha  } \right]\times k} -
2\cdot\big[  \sigma _{i-1,i-1,i-1,i}^{\left[ \gamma  \atop{  \alpha  \atop  \beta  } \right]\times k}
+  \sigma _{i-2,i-1,i-1,i}^{\left[ \gamma  \atop{  \alpha \atop  \beta  } \right]\times k}  \big] = 0 \quad \text{for $2\leq i\leq k $}
 \end{align*}
 \newpage
  \underline{The case $ k\geq  3s+2m+l-1 $}
\begin{align*} 
 &  \eqref{eq 2.62} \Longleftrightarrow  \\
  & \sum_{i=0}^{3s+2m+l-1}\big[\sigma _{i,i,i,i}^{\left[\beta | \gamma  \atop{ \gamma | \beta \atop \alpha | \alpha  } \right]\times k}
+  \sigma _{i-1,i,i,i}^{\left[\beta | \gamma  \atop{ \gamma | \beta \atop \alpha | \alpha  } \right]\times k}
+  \sigma _{i-1,i-1,i,i}^{\left[\beta | \gamma  \atop{ \gamma | \beta \atop \alpha | \alpha  } \right]\times k}
+  \sigma _{i-2,i-1,i,i}^{\left[\beta | \gamma  \atop{ \gamma | \beta \atop \alpha | \alpha  } \right]\times k}\big]\cdot2^{-iq}  \\
& = 2\cdot2^{-iq}
 \sum_{i=0}^{3s+2m+l-2}\big[\sigma _{i,i,i,i}^{\left[\gamma  \atop{ \alpha \atop  \beta   } \right]\times k}
+ \sigma _{i,i,i,i+1}^{\left[\gamma  \atop{ \alpha  \atop  \beta   } \right]\times k}
+ \sigma _{i-1,i,i,i}^{\left[\gamma  \atop{ \alpha  \atop \beta   } \right]\times k}
+ \sigma _{i-1,i,i,i+1}^{\left[\gamma  \atop{ \alpha \atop  \beta   } \right]\times k} \big]\cdot2^{-iq} 
 + (2^q -2)\cdot 2^{-q} \sum_{i=0}^{3s+2m+l-2}  \big[\sigma _{i,i,i,i}^{\left[\gamma  | \gamma  \atop{ \alpha  | \beta   \atop \beta  | \alpha  } \right]\times k}
 + \sigma _{i-1,i,i,i}^{\left[\gamma  | \gamma   \atop{ \alpha  | \beta  \atop \beta  | \alpha  } \right]\times k} \big]\cdot2^{-iq} \\
 & \Longleftrightarrow \\
   & \sum_{i=0}^{3s+2m+l-1}\big[\sigma _{i,i,i,i}^{\left[\beta | \gamma  \atop{ \gamma | \beta \atop \alpha | \alpha  } \right]\times k}
+  \sigma _{i-1,i,i,i}^{\left[\beta | \gamma  \atop{ \gamma | \beta \atop \alpha | \alpha  } \right]\times k}
+  \sigma _{i-1,i-1,i,i}^{\left[\beta | \gamma  \atop{ \gamma | \beta \atop \alpha | \alpha  } \right]\times k}
+  \sigma _{i-2,i-1,i,i}^{\left[\beta | \gamma  \atop{ \gamma | \beta \atop \alpha | \alpha  } \right]\times k}\big]\cdot2^{-iq}  \\
& =  \sum_{i=1}^{3s+2m+l-1}2\cdot\big[\sigma _{i-1,i-1,i-1,i-1}^{\left[\gamma  \atop{ \alpha \atop  \beta   } \right]\times k}
+ \sigma _{i-1,i-1,i-1,i}^{\left[\gamma  \atop{ \alpha  \atop  \beta   } \right]\times k}
+ \sigma _{i-2,i-1,i-1,i-1}^{\left[\gamma  \atop{ \alpha  \atop \beta   } \right]\times k}
+ \sigma _{i-2,i-1,i-1,i}^{\left[\gamma  \atop{ \alpha \atop  \beta   } \right]\times k} \big]\cdot2^{-iq} \\
& +  \sum_{i=0}^{3s+2m+l-2}  \big[\sigma _{i,i,i,i}^{\left[\gamma  | \gamma  \atop{ \alpha  | \beta   \atop \beta  | \alpha  } \right]\times k}
 + \sigma _{i-1,i,i,i}^{\left[\gamma  | \gamma   \atop{ \alpha  | \beta  \atop \beta  | \alpha  } \right]\times k} \big]\cdot2^{-iq} 
  -  \sum_{i=1}^{3s+2m+l-1}2\cdot \big[\sigma _{i-1,i-1,i-1,i-1}^{\left[\gamma  | \gamma  \atop{ \alpha  | \beta   \atop \beta  | \alpha  } \right]\times k}
 + \sigma _{i-2,i-1,i-1,i-1}^{\left[\gamma  | \gamma   \atop{ \alpha  | \beta  \atop \beta  | \alpha  } \right]\times k} \big]\cdot2^{-iq} \\
 & \Longrightarrow \\
& \sum_{i=2}^{3s+2m+l-2}\Big( \sigma _{i-1,i-1,i,i}^{\left[ \gamma  \atop{  \beta \atop  \alpha  } \right]\times k}
+  \sigma _{i-2,i-1,i,i}^{\left[ \gamma  \atop{  \beta \atop  \alpha  } \right]\times k} -
2\cdot\big[  \sigma _{i-1,i-1,i-1,i}^{\left[ \gamma  \atop{  \alpha  \atop  \beta  } \right]\times k}
+  \sigma _{i-2,i-1,i-1,i}^{\left[ \gamma  \atop{  \alpha \atop  \beta  } \right]\times k}  \big] \Big)\cdot2^{-iq} = 0 \quad \text{for all $ q \geq 2 $}\\
& \Big( \sigma _{3s+2m+l-3,3s+2m+l-2,3s+2m+l-1,3s+2m+l-1}^{\left[ \gamma  \atop{  \beta \atop  \alpha  } \right]\times k}
- 2\cdot \sigma _{3s+2m+l-3,3s+2m+l-2,3s+2m+l-2,3s+2m+l-1}^{\left[ \gamma  \atop{  \alpha \atop  \beta  } \right]\times k}\Big)\cdot2^{-(3s+2m+l-1)q} = 0 \; \text{for all $ q \geq 2 $}\\ 
 & \Longrightarrow \\
&  \sigma _{i-1,i-1,i,i}^{\left[ \gamma  \atop{  \beta \atop  \alpha  } \right]\times k}
+  \sigma _{i-2,i-1,i,i}^{\left[ \gamma  \atop{  \beta \atop  \alpha  } \right]\times k} -
2\cdot\big[  \sigma _{i-1,i-1,i-1,i}^{\left[ \gamma  \atop{  \alpha  \atop  \beta  } \right]\times k}
+  \sigma _{i-2,i-1,i-1,i}^{\left[ \gamma  \atop{  \alpha \atop  \beta  } \right]\times k}  \big] = 0 \quad \text{for $2\leq i\leq 3s+2m+l-2 $}\\
  & \sigma _{3s+2m+l-3,3s+2m+l-2,3s+2m+l-1,3s+2m+l-1}^{\left[ \gamma  \atop{  \beta \atop  \alpha  } \right]\times k}
- 2\cdot \sigma _{3s+2m+l-3,3s+2m+l-2,3s+2m+l-2,3s+2m+l-1}^{\left[ \gamma  \atop{  \alpha \atop  \beta  } \right]\times k}= 0
 \end{align*}
 
\end{proof} 
   \begin{lem}
\label{lem 2.13}
 We can write : \\
 \begin{equation}
 \label{eq 2.67}
   \kappa _{2}(t,\eta,\xi  ) =  h _{1}(t,\eta,\xi  ) + h _{3}(t,\eta,\xi  )
\end{equation}
Let q be a rational integer $\geq 2, $ then we obtain  from the binomial formula   \\
\begin{align}
\displaystyle
 \kappa _{2}^q = ( h _{1}+ h _{3})^q   =  h _{1}^q  + h _{3}^q    +   \sum_{i=1}^{q-1}\binom{q}{i} h _{1}^i  h _{3}^{q-i}  \label{eq 2.68}  
\end{align}
  
By integrating  \eqref{eq 2.68}  on the unit interval of $\mathbb{K}^3 $ we get \\
\begin{align}
\displaystyle
& \int_{\mathbb{P}^3} \kappa _{2}^q(t,\eta,\xi )dtd\eta d\xi
 =  2\int_{\{(t,\eta ,\xi )\in \mathbb{P}^3 \mid h _{1}(t,\eta,\xi ) \neq 0 \}}\varphi ^q(t,\eta ,\xi )dtd\eta d\xi 
 + (2^q -2)\int_{\{(t,\eta ,\xi )\in \mathbb{P}^3 \mid h _{1}(t,\eta,\xi )\cdot h _{3}(t,\eta,\xi ) \neq 0 \}}\varphi ^q(t,\eta ,\xi )dtd\eta d\xi  \label{eq 2.69} 
\end{align}
\end{lem}
\begin{proof} 
Similar to the proof of Lemma \ref{lem 2.11}, using Lemma \ref{lem 2.1}, \eqref{eq 2.45},\eqref{eq 2.47} and \eqref{eq 2.53}.
\end{proof}
   \begin{lem}
\label{lem 2.14}
We have the following rank formula :\\
\begin{align}
 \label{eq 2.70}
 \displaystyle 
\boxed{ \sigma _{i-1,i,i,i}^{\left[ \gamma  \atop{  \alpha  \atop  \beta   } \right]\times k}
+  \sigma _{i-1,i,i,i+1}^{\left[ \gamma  \atop{  \alpha \atop  \beta  } \right]\times k}=
2\cdot\big[  \sigma _{i-1,i-1,i,i}^{\left[ \alpha  \atop{  \gamma  \atop  \beta  } \right]\times k}
+  \sigma _{i-1,i-1,i,i+1}^{\left[ \alpha  \atop{  \gamma  \atop  \beta  } \right]\times k}  \big] \quad \text{for $0\leq i\leq \inf(k,3s+2m+l-2)$}}
\end{align}
 \end{lem}
 \begin{proof}
We proceed as in Lemma \ref{lem 2.12} by  combining  \eqref{eq 2.50}, \eqref{eq 2.45}, \eqref{eq 2.47}  and \eqref{eq 2.53}.     \\
\begin{align*}
\displaystyle
& \int_{\mathbb{P}^3} \kappa _{2}^q(t,\eta,\xi )dtd\eta d\xi
 =  2\int_{\{(t,\eta ,\xi )\in \mathbb{P}^3 \mid h _{1}(t,\eta,\xi ) \neq 0 \}}\varphi ^q(t,\eta ,\xi )dtd\eta d\xi 
 + (2^q -2)\int_{\{(t,\eta ,\xi )\in \mathbb{P}^3 \mid h _{1}(t,\eta,\xi )\cdot h _{3}(t,\eta,\xi ) \neq 0 \}}\varphi ^q(t,\eta ,\xi )dtd\eta d\xi \\
 & \Longleftrightarrow \\
  &  2^{(k+3s+2m+l-2)q -(3k +3s+2m+l-3)}
 \sum_{i=0}^{\inf(k,3s+2m+l-2)}\big[\sigma _{i,i,i,i}^{\left[ \gamma  \atop{ \alpha \atop \beta   } \right]\times k}
+  \sigma _{i,i,i,i+1}^{\left[ \gamma  \atop{ \alpha \atop \beta   } \right]\times k}
+  \sigma _{i-1,i,i,i}^{\left[ \gamma  \atop{ \alpha  \atop \beta  } \right]\times k}
+  \sigma _{i-1,i,i,i+1}^{\left[ \gamma  \atop{ \alpha \atop \beta  } \right]\times k}\big]\cdot2^{-iq}  \\
  & = 2\cdot 2^{(k+3s+2m+l-3)q -(3k +3s+2m+l-3)}
 \sum_{i=0}^{\inf(k,3s+2m+l-3)}\big[\sigma _{i,i,i,i}^{\left[\alpha  | \alpha  \atop{ \beta  | \gamma \atop \gamma  | \beta   } \right]\times k}
+ \sigma _{i,i,i,i+1}^{\left[\alpha  | \alpha  \atop{ \beta  | \gamma \atop \gamma  | \beta   } \right]\times k}
+ \sigma _{i,i,i+1,i+1}^{\left[\alpha  | \alpha  \atop{ \beta  | \gamma \atop \gamma  | \beta   } \right]\times k}
+ \sigma _{i,i,i+1,i+2}^{\left[\alpha  | \alpha  \atop{ \beta  | \gamma \atop \gamma  | \beta   } \right]\times k} \big]\cdot2^{-iq}  \\
 & + (2^q-2)\cdot 2^{(k+3s+2m+l-3)q -(3k +3s+2m+l-3)}
 \sum_{i=0}^{\inf(k,3s+2m+l-3)}\big[\sigma _{i,i,i,i}^{\left[\alpha  | \gamma  \atop{ \gamma  | \alpha  \atop \beta  | \beta  } \right]\times k}
 + \sigma _{i,i,i,i+1}^{\left[\alpha  | \gamma  \atop{ \gamma  | \alpha \atop \beta  | \beta  } \right]\times k} \big]\cdot2^{-iq} \\
 & \Longleftrightarrow  \\
 &  \sum_{i=0}^{\inf(k,3s+2m+l-2)}\big[\sigma _{i,i,i,i}^{\left[ \gamma  \atop{ \alpha \atop \beta   } \right]\times k}
+  \sigma _{i,i,i,i+1}^{\left[ \gamma  \atop{ \alpha \atop \beta   } \right]\times k}
+  \sigma _{i-1,i,i,i}^{\left[ \gamma  \atop{ \alpha  \atop \beta  } \right]\times k}
+  \sigma _{i-1,i,i,i+1}^{\left[ \gamma  \atop{ \alpha \atop \beta  } \right]\times k}\big]\cdot2^{-iq}  \\
  & = 
 \sum_{i=1}^{\inf(k,3s+2m+l-3)+1}2\cdot\big[\sigma _{i-1,i-1,i-1,i-1}^{\left[\alpha  | \alpha  \atop{ \beta  | \gamma \atop \gamma  | \beta   } \right]\times k}
+ \sigma _{i-1,i-1,i-1,i}^{\left[\alpha  | \alpha  \atop{ \beta  | \gamma \atop \gamma  | \beta   } \right]\times k}
+ \sigma _{i-1,i-1,i,i}^{\left[\alpha  | \alpha  \atop{ \beta  | \gamma \atop \gamma  | \beta   } \right]\times k}
+ \sigma _{i-1,i-1,i,i+1}^{\left[\alpha  | \alpha  \atop{ \beta  | \gamma \atop \gamma  | \beta   } \right]\times k} \big]\cdot2^{-iq}  \\
 & 
 + \sum_{i=0}^{\inf(k,3s+2m+l-3)}\big[\sigma _{i,i,i,i}^{\left[\alpha  | \gamma  \atop{ \gamma  | \alpha  \atop \beta  | \beta  } \right]\times k}
 + \sigma _{i,i,i,i+1}^{\left[\alpha  | \gamma  \atop{ \gamma  | \alpha \atop \beta  | \beta  } \right]\times k} \big]\cdot2^{-iq} 
  -  \sum_{i =1}^{\inf(k,3s+2m+l-3)+1}2\cdot\big[\sigma _{i-1,i-1,i-1,i-1}^{\left[\alpha  | \gamma  \atop{ \gamma  | \alpha  \atop \beta  | \beta  } \right]\times k}
 + \sigma _{i-1,i-1,i-1,i}^{\left[\alpha  | \gamma  \atop{ \gamma  | \alpha \atop \beta  | \beta  } \right]\times k} \big]\cdot2^{-iq} 
\end{align*}

\end{proof} 
   \begin{lem}
\label{lem 2.15}
 We can write : \\
 \begin{equation}
 \label{eq 2.71}
   g _{2}(t,\eta,\xi  ) =  h _{1}(t,\eta,\xi  ) + h _{2}(t,\eta,\xi  )
\end{equation}
Let q be a rational integer $\geq 2, $ then we obtain  from the binomial formula   \\
\begin{align}
\displaystyle
 g _{2}^q = ( h _{1}+ h _{2})^q   =  h _{1}^q  + h _{2}^q    +   \sum_{i=1}^{q-1}\binom{q}{i} h _{1}^i  h _{2}^{q-i}  \label{eq 2.72}  
\end{align}
  
By integrating  \eqref{eq 2.72}  on the unit interval of $\mathbb{K}^3 $ we get \\
\begin{align}
\displaystyle
& \int_{\mathbb{P}^3} g _{2}^q(t,\eta,\xi )dtd\eta d\xi
 =  2\int_{\{(t,\eta ,\xi )\in \mathbb{P}^3 \mid h _{1}(t,\eta,\xi ) \neq 0 \}}\varphi ^q(t,\eta ,\xi )dtd\eta d\xi 
 + (2^q -2)\int_{\{(t,\eta ,\xi )\in \mathbb{P}^3 \mid h _{1}(t,\eta,\xi )\cdot h _{2}(t,\eta,\xi ) \neq 0 \}}\varphi ^q(t,\eta ,\xi )dtd\eta d\xi  \label{eq 2.73} 
\end{align}
\end{lem} 
 \begin{proof} 
Similar to the proof of Lemma \ref{lem 2.11}, using Lemma \ref{lem 2.1}, \eqref{eq 2.45},\eqref{eq 2.46} and \eqref{eq 2.52}.
\end{proof}

   \begin{lem}
\label{lem 2.16}
We have the following rank formula :\\
\begin{align}
 \label{eq 2.74}
 \displaystyle 
\boxed{ \sigma _{i-1,i,i,i}^{\left[ \beta  \atop{  \alpha  \atop  \gamma   } \right]\times k}
+  \sigma _{i-1,i,i,i+1}^{\left[ \beta   \atop{  \alpha \atop  \gamma  } \right]\times k}=
2\cdot\big[  \sigma _{i-1,i-1,i,i}^{\left[ \alpha  \atop{  \beta  \atop  \gamma  } \right]\times k}
+  \sigma _{i-1,i-1,i,i+1}^{\left[ \alpha   \atop{  \beta  \atop  \gamma  } \right]\times k}  \big] \quad \text{for $0\leq i\leq \inf(k,3s+2m+l-2)$}}
\end{align}
 \end{lem}
 \begin{proof}
We proceed as in Lemma \ref{lem 2.12} by  combining  \eqref{eq 2.51}, \eqref{eq 2.45}, \eqref{eq 2.46}  and \eqref{eq 2.52}.     \\
\begin{align*}
\displaystyle
& \int_{\mathbb{P}^3} g _{2}^q(t,\eta,\xi )dtd\eta d\xi
 =  2\int_{\{(t,\eta ,\xi )\in \mathbb{P}^3 \mid h _{1}(t,\eta,\xi ) \neq 0 \}}\varphi ^q(t,\eta ,\xi )dtd\eta d\xi 
 + (2^q -2)\int_{\{(t,\eta ,\xi )\in \mathbb{P}^3 \mid h _{1}(t,\eta,\xi )\cdot h _{2}(t,\eta,\xi ) \neq 0 \}}\varphi ^q(t,\eta ,\xi )dtd\eta d\xi \\
 & \Longleftrightarrow \\
  &  2^{(k+3s+2m+l-2)q -(3k +3s+2m+l-3)}
 \sum_{i=0}^{\inf(k,3s+2m+l-2)}\big[\sigma _{i,i,i,i}^{\left[ \beta  \atop{ \alpha \atop \gamma    } \right]\times k}
+  \sigma _{i,i,i,i+1}^{\left[ \beta   \atop{ \alpha \atop \gamma   } \right]\times k}
+  \sigma _{i-1,i,i,i}^{\left[ \beta   \atop{ \alpha  \atop \gamma  } \right]\times k}
+  \sigma _{i-1,i,i,i+1}^{\left[ \beta   \atop{ \alpha \atop \gamma  } \right]\times k}\big]\cdot2^{-iq}  \\
  & = 2\cdot 2^{(k+3s+2m+l-3)q -(3k +3s+2m+l-3)}
 \sum_{i=0}^{\inf(k,3s+2m+l-3)}\big[\sigma _{i,i,i,i}^{\left[\alpha  | \alpha  \atop{ \beta  | \gamma \atop \gamma  | \beta   } \right]\times k}
+ \sigma _{i,i,i,i+1}^{\left[\alpha  | \alpha  \atop{ \beta  | \gamma \atop \gamma  | \beta   } \right]\times k}
+ \sigma _{i,i,i+1,i+1}^{\left[\alpha  | \alpha  \atop{ \beta  | \gamma \atop \gamma  | \beta   } \right]\times k}
+ \sigma _{i,i,i+1,i+2}^{\left[\alpha  | \alpha  \atop{ \beta  | \gamma \atop \gamma  | \beta   } \right]\times k} \big]\cdot2^{-iq}  \\
 & + (2^q-2)\cdot 2^{(k+3s+2m+l-3)q -(3k +3s+2m+l-3)}
 \sum_{i=0}^{\inf(k,3s+2m+l-3)}\big[\sigma _{i,i,i,i}^{\left[\alpha  | \beta  \atop{ \beta  | \alpha  \atop \gamma  | \gamma  } \right]\times k}
 + \sigma _{i,i,i,i+1}^{\left[\alpha  | \beta  \atop{ \beta  | \alpha \atop \gamma  | \gamma  } \right]\times k} \big]\cdot2^{-iq} \\
 & \Longleftrightarrow  \\
 &  \sum_{i=0}^{\inf(k,3s+2m+l-2)}\big[\sigma _{i,i,i,i}^{\left[ \beta  \atop{ \alpha \atop \gamma    } \right]\times k}
+  \sigma _{i,i,i,i+1}^{\left[ \beta   \atop{ \alpha \atop \gamma   } \right]\times k}
+  \sigma _{i-1,i,i,i}^{\left[ \beta  \atop{ \alpha  \atop \gamma  } \right]\times k}
+  \sigma _{i-1,i,i,i+1}^{\left[ \beta  \atop{ \alpha \atop \gamma  } \right]\times k}\big]\cdot2^{-iq}  \\
  & = 
 \sum_{i=1}^{\inf(k,3s+2m+l-3)+1}2\cdot\big[\sigma _{i-1,i-1,i-1,i-1}^{\left[\alpha  | \alpha  \atop{ \beta  | \gamma \atop \gamma  | \beta   } \right]\times k}
+ \sigma _{i-1,i-1,i-1,i}^{\left[\alpha  | \alpha  \atop{ \beta  | \gamma \atop \gamma  | \beta   } \right]\times k}
+ \sigma _{i-1,i-1,i,i}^{\left[\alpha  | \alpha  \atop{ \beta  | \gamma \atop \gamma  | \beta   } \right]\times k}
+ \sigma _{i-1,i-1,i,i+1}^{\left[\alpha  | \alpha  \atop{ \beta  | \gamma \atop \gamma  | \beta   } \right]\times k} \big]\cdot2^{-iq}  \\
 & 
+ \sum_{i=0}^{\inf(k,3s+2m+l-3)}\big[\sigma _{i,i,i,i}^{\left[\alpha  | \beta  \atop{ \beta  | \alpha  \atop \gamma   | \gamma  } \right]\times k}
 + \sigma _{i,i,i,i+1}^{\left[\alpha  | \beta   \atop{ \beta  | \alpha \atop \gamma  | \gamma  } \right]\times k} \big]\cdot2^{-iq} 
  -  \sum_{i =1}^{\inf(k,3s+2m+l-3)+1}2\cdot\big[\sigma _{i-1,i-1,i-1,i-1}^{\left[\alpha  | \beta   \atop{ \beta  | \alpha  \atop \gamma   | \gamma  } \right]\times k}
 + \sigma _{i-1,i-1,i-1,i}^{\left[\alpha  | \beta   \atop{ \beta  | \alpha \atop \gamma  | \gamma  } \right]\times k} \big]\cdot2^{-iq} 
\end{align*}

\end{proof}  
   \begin{lem}
\label{lem 2.17}
 We can write : \\
 \begin{equation}
 \label{eq 2.75}
   g_{1}(t,\eta,\xi  ) =  h_{1}(t,\eta,\xi  ) + h_{2}(t,\eta,\xi  ) +  h_{3}(t,\eta,\xi  ) + h_{4}(t,\eta,\xi  ) 
\end{equation}
Let q be a rational integer $\geq 2, $ then we obtain  from the binomial formula   \\
\begin{align}
\displaystyle
& g_{1}^q = ( h_{1}+h_{2}+h_{3}+h_{4})^q = ( h_{1}+h_{2})^q +  ( h_{3}+h_{4})^q + \sum_{i=1}^{q-1}\binom{q}{i}( h_{1}+h_{2})^i ( h_{3}+h_{4})^{q-i} \label{eq 2.76}\\
 & =  h_{1}^q  + h_{2}^q  +  h_{3}^q  +  h_{4}^q   +   \sum_{i=1}^{q-1}\binom{q}{i} h_{1}^i  h_{2}^{q-i}  +  \sum_{i=1}^{q-1}\binom{q}{i} h_{3}^i  h_{4}^{q-i}    
  +  \sum_{i=1}^{q-1}\binom{q}{i} h_{1}^i  h_{3}^{q-i} +  \sum_{i=1}^{q-1}\binom{q}{i} h_{1}^i  h_{4}^{q-i} +\sum_{i=1}^{q-1}\binom{q}{i} h_{2}^i  h_{3}^{q-i} +\sum_{i=1}^{q-1}\binom{q}{i} h_{2}^i  h_{4}^{q-i} \nonumber \\
  & +  \sum_{i=1}^{q-1}\sum_{j =1}^{i-1}\binom{q}{i}\binom{i}{j} h_{1}^{j} h_{2}^{i-j} h_{3}^{q-i}+  \sum_{i=1}^{q-1}\sum_{j =1}^{i-1}\binom{q}{i}\binom{i}{j} h_{1}^{j} h_{2}^{i-j} h_{4}^{q-i}  
  +  \sum_{i=1}^{q-1}\sum_{r =1}^{q-i-1}\binom{q}{i}\binom{q-i}{r} h_{1}^{i} h_{3}^{r} h_{4}^{q-i-r} \nonumber   \\
 &  +   \sum_{i=1}^{q-1}\sum_{r =1}^{q-i-1}\binom{q}{i}\binom{q-i}{r} h_{2}^{i} h_{3}^{r} h_{4}^{q-i-r}
   +  \sum_{i=1}^{q-1}\sum_{j =1}^{i-1}\sum_{r =1}^{q-i-1}\binom{q}{i}\binom{i}{j}\binom{q-i}{r} h_{1}^{j}h_{2}^{i-j}h_{3}^{r} h_{4}^{q-i-r} \nonumber 
\end{align}
 By integrating  \eqref{eq 2.76}  on the unit interval of $\mathbb{K}^3 $ we get \\
\begin{align}
\displaystyle
& \int_{\mathbb{P}^3} g_{1}^q(t,\eta,\xi )dtd\eta d\xi =  4\int_{\{(t,\eta ,\xi )\in \mathbb{P}^3 \mid h_{1}(t,\eta,\xi ) \neq 0 \}}\varphi ^q(t,\eta ,\xi )dtd\eta d\xi 
+4(2^q -2)\int_{\{(t,\eta ,\xi )\in \mathbb{P}^3 \mid h_{1}(t,\eta,\xi )\cdot h_{2}(t,\eta,\xi ) \neq 0 \}}\varphi ^q(t,\eta ,\xi )dtd\eta d\xi \label{eq 2.77} \\
 & + 2(2^q -2)\int_{\{(t,\eta ,\xi )\in \mathbb{P}^3 \mid h_{1}(t,\eta,\xi )\cdot h_{3}(t,\eta,\xi ) \neq 0 \}}\varphi ^q(t,\eta ,\xi )dtd\eta d\xi \nonumber \\
  & +(2^{2q} -6\cdot2^{q} +8)\int_{\{(t,\eta ,\xi )\in \mathbb{P}^3 \mid h_{1}(t,\eta,\xi )\cdot h_{2}(t,\eta,\xi )\cdot h_{3}(t,\eta,\xi ) \neq 0 \}}\varphi ^q(t,\eta ,\xi )dtd\eta d\xi \nonumber 
\end{align}
\begin{align*}
& \Longleftrightarrow \\
 &  2^{(k+3s+2m+l-1)q -(3k +3s+2m+l-3)}
 \sum_{i=0}^{\inf(k,3s+2m+l-1)}\big[\sigma _{i,i,i,i}^{\left[\beta | \gamma  \atop{ \gamma | \beta \atop \alpha | \alpha  } \right]\times k}
+  \sigma _{i-1,i,i,i}^{\left[\beta | \gamma  \atop{ \gamma | \beta \atop \alpha | \alpha  } \right]\times k}
+  \sigma _{i-1,i-1,i,i}^{\left[\beta | \gamma  \atop{ \gamma | \beta \atop \alpha | \alpha  } \right]\times k}
+  \sigma _{i-2,i-1,i,i}^{\left[\beta | \gamma  \atop{ \gamma | \beta \atop \alpha | \alpha  } \right]\times k}\big]\cdot2^{-iq}  \\
& = 4\cdot2^{(k+3s+2m+l-3)q -(3k +3s+2m+l-3)}
 \sum_{i=0}^{\inf(k,3s+2m+l-3)}\big[\sigma _{i,i,i,i}^{\left[\alpha  | \alpha  \atop{ \beta  | \gamma \atop \gamma  | \beta   } \right]\times k}
+ \sigma _{i,i,i,i+1}^{\left[\alpha  | \alpha  \atop{ \beta  | \gamma \atop \gamma  | \beta   } \right]\times k}
+ \sigma _{i,i,i+1,i+1}^{\left[\alpha  | \alpha  \atop{ \beta  | \gamma \atop \gamma  | \beta   } \right]\times k}
+ \sigma _{i,i,i+1,i+2}^{\left[\alpha  | \alpha  \atop{ \beta  | \gamma \atop \gamma  | \beta   } \right]\times k} \big]\cdot2^{-iq}  \\
& + 4\cdot(2^{q} -2)\cdot 2^{(k+3s+2m+l-3)q -(3k +3s+2m+l-3)} \sum_{i=0}^{\inf(k,3s+2m+l-3)}\big[\sigma _{i,i,i,i}^{\left[\alpha  | \beta  \atop{ \beta  | \alpha  \atop \gamma  | \gamma  } \right]\times k}
 + \sigma _{i,i,i,i+1}^{\left[\alpha  | \beta  \atop{ \beta | \alpha \atop \gamma  | \gamma  } \right]\times k} \big]\cdot2^{-iq}   \\
 & + 2\cdot(2^{q} -2)\cdot 2^{(k+3s+2m+l-3)q -(3k +3s+2m+l-3)} \sum_{i=0}^{\inf(k,3s+2m+l-3)}\big[\sigma _{i,i,i,i}^{\left[\alpha  | \gamma  \atop{ \gamma  | \alpha  \atop \beta  | \beta  } \right]\times k}
 + \sigma _{i,i,i,i+1}^{\left[\alpha  | \gamma  \atop{ \gamma  | \alpha \atop \beta  | \beta  } \right]\times k} \big]\cdot2^{-iq}   \\
 & + (2^{2q} -6\cdot2^{q} +8)\cdot2^{(k+3s+2m+l-3)q -(3k +3s+2m+l-3)} \sum_{i=0}^{\inf(k,3s+2m+l-3)} \sigma _{i,i,i,i}^{\left[\alpha \atop{ \beta \atop \gamma } \right]\times k}\cdot2^{-iq} \\
 & \Longleftrightarrow  \nonumber \\
 &  \sum_{i=0}^{\inf(k,3s+2m+l-1)}\big[\sigma _{i,i,i,i}^{\left[\beta | \gamma  \atop{ \gamma | \beta \atop \alpha | \alpha  } \right]\times k}
+  \sigma _{i-1,i,i,i}^{\left[\beta | \gamma  \atop{ \gamma | \beta \atop \alpha | \alpha  } \right]\times k}
+  \sigma _{i-1,i-1,i,i}^{\left[\beta | \gamma  \atop{ \gamma | \beta \atop \alpha | \alpha  } \right]\times k}
+  \sigma _{i-2,i-1,i,i}^{\left[\beta | \gamma  \atop{ \gamma | \beta \atop \alpha | \alpha  } \right]\times k}\big]\cdot2^{-iq}  \\
& =4\cdot
 \sum_{i=2}^{\inf(k,3s+2m+l-3)+2}\big[\sigma _{i-2,i-2,i-2,i-2}^{\left[\alpha  | \alpha  \atop{ \beta  | \gamma \atop \gamma  | \beta   } \right]\times k}
+ \sigma _{i-2,i-2,i-2,i-1}^{\left[\alpha  | \alpha  \atop{ \beta  | \gamma \atop \gamma  | \beta   } \right]\times k}
+ \sigma _{i-2,i-2,i-1,i-1}^{\left[\alpha  | \alpha  \atop{ \beta  | \gamma \atop \gamma  | \beta   } \right]\times k}
+ \sigma _{i-2,i-2,i-1,i}^{\left[\alpha  | \alpha  \atop{ \beta  | \gamma \atop \gamma  | \beta   } \right]\times k} \big]\cdot2^{-iq}  \\
& +  \sum_{i=1}^{\inf(k,3s+2m+l-3)+1}\big[4\cdot\sigma _{i-1,i-1,i-1,i-1}^{\left[\alpha  | \beta  \atop{ \beta  | \alpha  \atop \gamma  | \gamma  } \right]\times k}
 +4\cdot \sigma _{i-1,i-1,i-1,i}^{\left[\alpha  | \beta  \atop{ \beta | \alpha \atop \gamma  | \gamma  } \right]\times k}\big]  
  - \sum_{i=2}^{\inf(k,3s+2m+l-3)+2}\big[ 8\cdot\sigma _{i-2,i-2,i-2,i-2}^{\left[\alpha  | \beta  \atop{ \beta  | \alpha  \atop \gamma  | \gamma  } \right]\times k}
 +8\cdot \sigma _{i-2,i-2,i-2,i-1}^{\left[\alpha  | \beta  \atop{ \beta | \alpha \atop \gamma  | \gamma  } \right]\times k} \big]\cdot2^{-iq}  \\
  & +  \sum_{i=1}^{\inf(k,3s+2m+l-3)+1}\big[ 2\cdot\sigma _{i-1,i-1,i-1,i-1}^{\left[\alpha  | \gamma  \atop{ \gamma  | \alpha  \atop \beta  | \beta  } \right]\times k}
 +2\cdot \sigma _{i-1,i-1,i-1,i}^{\left[\alpha  | \gamma  \atop{ \gamma  | \alpha \atop \beta  | \beta  } \right]\times k}\big]  
  - \sum_{i=2}^{\inf(k,3s+2m+l-3)+2}\big[ 4\cdot\sigma _{i-2,i-2,i-2,i-2}^{\left[\alpha  | \gamma  \atop{ \gamma   | \alpha  \atop \beta   | \beta  } \right]\times k}
 +4\cdot \sigma _{i-2,i-2,i-2,i-1}^{\left[\alpha  | \gamma   \atop{ \gamma  | \alpha \atop \beta  | \beta   } \right]\times k} \big]\cdot2^{-iq}  \\
   & +  \sum_{i=0}^{\inf(k,3s+2m+l-3)} \sigma _{i,i,i,i}^{\left[\alpha \atop{ \beta \atop \gamma } \right]\times k}\cdot2^{-iq} 
-6\cdot\sum_{i=1}^{\inf(k,3s+2m+l-3)+1} \sigma _{i-1,i-1,i-1,i-1}^{\left[\alpha \atop{ \beta \atop \gamma } \right]\times k}\cdot2^{-iq} 
 +8\cdot\sum_{i=2}^{\inf(k,3s+2m+l-3)+2} \sigma _{i-2,i-2,i-2,i-2}^{\left[\alpha \atop{ \beta \atop \gamma } \right]\times k}\cdot2^{-iq}
 \end{align*}
 \begin{align*}
 & \Longrightarrow  \\
 &  \\
 & \text{In the case $k\leq 3s+2m+l-3$}   \\
 &   \\
  &  \sum_{i=2}^{k+1}\big[  \sigma _{i-1,i,i,i}^{\left[\beta | \gamma  \atop{ \gamma | \beta \atop \alpha | \alpha  } \right]\times k}
+  \sigma _{i-1,i-1,i,i}^{\left[\beta | \gamma  \atop{ \gamma | \beta \atop \alpha | \alpha  } \right]\times k}
+  \sigma _{i-2,i-1,i,i}^{\left[\beta | \gamma  \atop{ \gamma | \beta \atop \alpha | \alpha  } \right]\times k}\big]\cdot2^{-iq}  
 = 4\cdot
 \sum_{i=2}^{k+1}\big[ \sigma _{i-2,i-2,i-2,i-1}^{\left[\alpha  | \alpha  \atop{ \beta  | \gamma \atop \gamma  | \beta   } \right]\times k}
+ \sigma _{i-2,i-2,i-1,i-1}^{\left[\alpha  | \alpha  \atop{ \beta  | \gamma \atop \gamma  | \beta   } \right]\times k}
+ \sigma _{i-2,i-2,i-1,i}^{\left[\alpha  | \alpha  \atop{ \beta  | \gamma \atop \gamma  | \beta   } \right]\times k} \big]\cdot2^{-iq}   \\
& +  \sum_{i=2}^{k+1}\big[4\cdot \sigma _{i-1,i-1,i-1,i}^{\left[\alpha  | \beta  \atop{ \beta | \alpha \atop \gamma  | \gamma  } \right]\times k}  
   - 8\cdot \sigma _{i-2,i-2,i-2,i-1}^{\left[\alpha  | \beta  \atop{ \beta | \alpha \atop \gamma  | \gamma  } \right]\times k} \big]\cdot2^{-iq}  
  +  \sum_{i=2}^{k+1}\big[ 2\cdot \sigma _{i-1,i-1,i-1,i}^{\left[\alpha  | \gamma  \atop{ \gamma  | \alpha \atop \beta  | \beta  } \right]\times k}  
   - 4\cdot \sigma _{i-2,i-2,i-2,i-1}^{\left[\alpha  | \gamma   \atop{ \gamma  | \alpha \atop \beta  | \beta   } \right]\times k} \big]\cdot2^{-iq} \quad \text{for all $q\geq 2$} \\
& \\
& \text{In the case $k \geq  3s+2m+l-2$}\\
 & \\
  &  \sum_{i=2}^{3s+2m+l-1}\big[  \sigma _{i-1,i,i,i}^{\left[\beta | \gamma  \atop{ \gamma | \beta \atop \alpha | \alpha  } \right]\times k}
+  \sigma _{i-1,i-1,i,i}^{\left[\beta | \gamma  \atop{ \gamma | \beta \atop \alpha | \alpha  } \right]\times k}
+  \sigma _{i-2,i-1,i,i}^{\left[\beta | \gamma  \atop{ \gamma | \beta \atop \alpha | \alpha  } \right]\times k}\big]\cdot2^{-iq}  
 =4\cdot
 \sum_{i=2}^{3s+2m+l-1}\big[ \sigma _{i-2,i-2,i-2,i-1}^{\left[\alpha  | \alpha  \atop{ \beta  | \gamma \atop \gamma  | \beta   } \right]\times k}
+ \sigma _{i-2,i-2,i-1,i-1}^{\left[\alpha  | \alpha  \atop{ \beta  | \gamma \atop \gamma  | \beta   } \right]\times k}
+ \sigma _{i-2,i-2,i-1,i}^{\left[\alpha  | \alpha  \atop{ \beta  | \gamma \atop \gamma  | \beta   } \right]\times k} \big]\cdot2^{-iq}  \\
& +  \sum_{i=2}^{3s+2m+l-1}\big[4\cdot \sigma _{i-1,i-1,i-1,i}^{\left[\alpha  | \beta  \atop{ \beta | \alpha \atop \gamma  | \gamma  } \right]\times k}  
   - 8\cdot \sigma _{i-2,i-2,i-2,i-1}^{\left[\alpha  | \beta  \atop{ \beta | \alpha \atop \gamma  | \gamma  } \right]\times k} \big]\cdot2^{-iq}  
  +  \sum_{i=2}^{3s+2m+l-1}\big[ 2\cdot \sigma _{i-1,i-1,i-1,i}^{\left[\alpha  | \gamma  \atop{ \gamma  | \alpha \atop \beta  | \beta  } \right]\times k}  
   - 4\cdot \sigma _{i-2,i-2,i-2,i-1}^{\left[\alpha  | \gamma   \atop{ \gamma  | \alpha \atop \beta  | \beta   } \right]\times k} \big]\cdot2^{-iq}\quad \text{for all $q\geq 2$}    
  \end{align*}

\end{lem}

 \begin{proof}
From Lemma \ref{lem 2.1} we get easily \eqref{eq 2.75}, indeed we have \\
\begin{align*}
&  g_{1}(t,\eta,\xi  ) =  \sum_{deg Y\leq k-1}\sum_{deg Z = s-1}E(tYZ)\sum_{deg U \leq s+m-1}E(\eta YU) \sum_{deg V \leq s+m+l-1}E(\xi YV)\\
& =   \sum_{deg Y\leq k-1}E(tYT^{s-1}) \sum_{deg Z\leq  s-2}E(tYZ)\sum_{deg U \leq s+m-1}E(\eta YU) \sum_{deg V \leq s+m+l-1}E(\xi  YV)\\
& = \sum_{deg Y\leq k-1}E(tYT^{s-1}) \sum_{deg Z\leq  s-2}E(tYZ)\sum_{deg U \leq s+m-2}E(\eta YU) \sum_{deg V \leq s+m+l-2}E(\xi  YV) \\
& +   \sum_{deg Y\leq k-1}E(tYT^{s-1})E(\eta YT^{s+m-1}) \sum_{deg Z\leq  s-2}E(tYZ)\sum_{deg U \leq s+m-2}E(\eta YU) \sum_{deg V \leq s+m+l-2}E(\xi YV)   \\
& +   \sum_{deg Y\leq k-1}E(tYT^{s-1})E(\xi YT^{s+m+l-1}) \sum_{deg Z\leq  s-2}E(tYZ)\sum_{deg U \leq s+m-2}E(\eta YU) \sum_{deg V \leq s+m+l-2}E(\xi YV) \\
& +  \sum_{deg Y\leq k-1}E(tYT^{s-1})E(\eta YT^{s+m-1}) E(\xi YT^{s+m+l-1})\sum_{deg Z\leq  s-2}E(tYZ)\sum_{deg U \leq s+m-2}E(\eta YU) \sum_{deg V \leq s+m+l-2}E(\xi  YV ) \\
& =  h_{1}(t,\eta,\xi  ) + h_{2}(t,\eta,\xi  ) +  h_{3}(t,\eta,\xi  ) + h_{4}(t,\eta,\xi  ). 
\end{align*}
To prove \eqref{eq 2.76} we observe that  \\
\begin{align*}
\displaystyle
& \sum_{i=1}^{q-1}\binom{q}{i}( h_{1}+h_{2})^i ( h_{3}+h_{4})^{q-i}
 = \sum_{i=1}^{q-1}\binom{q}{i}\Big [h_{1}^{i} + h_{2}^{i}  +\sum_{j=1}^{i-1}\binom{i}{j}h_{1}^{j}h_{2}^{i-j}\Big ] \Big [h_{3}^{q-i} + h_{4}^{q-i}  +\sum_{r=1}^{q-i-1}\binom{q-i}{r}h_{3}^{r}h_{4}^{q-i-j}\Big ].  
\end{align*}
To prove \eqref{eq 2.77} we proceed as follows  \\

 From \eqref{eq 2.45},  \eqref{eq 2.45},  \eqref{eq 2.45} and  \eqref{eq 2.45} we get \\
  \begin{align}
\label{eq 2.78}
\displaystyle
\int_{\mathbb{P}^{3}}\big( h_{1}^q(t,\eta,\xi  ) + h_{2}^q(t,\eta,\xi  ) +  h_{3}^q(t,\eta,\xi  ) + h_{4}^q(t,\eta,\xi)\big)dtd\eta d\xi  =  4\int_{\{(t,\eta ,\xi )\in \mathbb{P}^3 \mid h_{1}(t,\eta,\xi ) \neq 0 \}}\varphi ^q(t,\eta ,\xi )dtd\eta d\xi \\
& \nonumber
\end{align}

By  \eqref{eq 2.52}, \eqref{eq 2.53}, \eqref{eq 2.54}, \eqref{eq 2.55}, \eqref{eq 2.56} and \eqref{eq 2.57} we obtain  \\

   \begin{align}
\label{eq 2.79}
\displaystyle
& \int_{\mathbb{P}^{3}}\big(\sum_{i=1}^{q-1}\binom{q}{i} h_{1}^i  h_{2}^{q-i}  +  \sum_{i=1}^{q-1}\binom{q}{i} h_{3}^i  h_{4}^{q-i}    
  +  \sum_{i=1}^{q-1}\binom{q}{i} h_{1}^i  h_{3}^{q-i} +  \sum_{i=1}^{q-1}\binom{q}{i} h_{1}^i  h_{4}^{q-i} +\sum_{i=1}^{q-1}\binom{q}{i} h_{2}^i  h_{3}^{q-i} +\sum_{i=1}^{q-1}\binom{q}{i} h_{2}^i  h_{4}^{q-i}\big)dtd\eta d\xi  \\
& = 4(2^q -2)\int_{\{(t,\eta ,\xi )\in \mathbb{P}^3 \mid h_{1}(t,\eta,\xi )\cdot h_{2}(t,\eta,\xi ) \neq 0 \}}\varphi ^q(t,\eta ,\xi )dtd\eta d\xi 
  + 2(2^q -2)\int_{\{(t,\eta ,\xi )\in \mathbb{P}^3 \mid h_{1}(t,\eta,\xi )\cdot h_{3}(t,\eta,\xi ) \neq 0 \}}\varphi ^q(t,\eta ,\xi )dtd\eta d\xi \nonumber\\
  & \nonumber
 \end{align}  
 
 From \eqref{eq 2.59}, \eqref{eq 2.60} we get \\

   \begin{align}
\label{eq 2.80}
\displaystyle
& \int_{\mathbb{P}^{3}}\Bigg[  \sum_{i=1}^{q-1}\sum_{j =1}^{i-1}\binom{q}{i}\binom{i}{j} h_{1}^{j} h_{2}^{i-j} h_{3}^{q-i}+  \sum_{i=1}^{q-1}\sum_{j =1}^{i-1}\binom{q}{i}\binom{i}{j} h_{1}^{j} h_{2}^{i-j} h_{4}^{q-i}  
  +  \sum_{i=1}^{q-1}\sum_{r =1}^{q-i-1}\binom{q}{i}\binom{q-i}{r} h_{1}^{i} h_{3}^{r} h_{4}^{q-i-r} \\
    &  +   \sum_{i=1}^{q-1}\sum_{r =1}^{q-i-1}\binom{q}{i}\binom{q-i}{r} h_{2}^{i} h_{3}^{r} h_{4}^{q-i-r} 
     +  \sum_{i=1}^{q-1}\sum_{j =1}^{i-1}\sum_{r =1}^{q-i-1}\binom{q}{i}\binom{i}{j}\binom{q-i}{r} h_{1}^{j}h_{2}^{i-j}h_{3}^{r} h_{4}^{q-i-r} \Bigg]dtd\eta d\xi  \nonumber \\
 & =(2^{2q} -6\cdot2^{q} +8)\int_{\{(t,\eta ,\xi )\in \mathbb{P}^3 \mid h_{1}(t,\eta,\xi )\cdot h_{2}(t,\eta,\xi )\cdot h_{3}(t,\eta,\xi ) \neq 0 \}}\varphi ^q(t,\eta ,\xi )dtd\eta d\xi \nonumber \\
 & \nonumber
 \end{align}    
 \eqref{eq 2.77} follows now from \eqref{eq 2.78}, \eqref{eq 2.79} and \eqref{eq 2.80} follows by \eqref{eq 2.44},  \eqref{eq 2.45},  \eqref{eq 2.52},   \eqref{eq 2.53} and   \eqref{eq 2.59}.
 \end{proof}
   \begin{lem}
\label{lem 2.18}
We have the following rank formula :\\
\begin{align}
 \label{eq 2.81}
 \displaystyle 
 &  \sigma _{i-1,i,i,i}^{\left[ \gamma  \atop{  \beta \atop \alpha   } \right]\times k}
+  \sigma _{i-1,i-1,i,i}^{\left[ \gamma  \atop{  \beta \atop  \alpha  } \right]\times k}
+  \sigma _{i-2,i-1,i,i}^{\left[ \gamma  \atop{  \beta \atop \alpha  } \right]\times k}  
 = 4\cdot\big[ \sigma _{i-2,i-2,i-2,i-1}^{\left[\alpha   \atop{  \gamma \atop  \beta   } \right]\times k}
+ \sigma _{i-2,i-2,i-1,i-1}^{\left[\alpha   \atop{  \gamma \atop  \beta   } \right]\times k}
+ \sigma _{i-2,i-2,i-1,i}^{\left[\alpha   \atop{  \gamma \atop  \beta   } \right]\times k} \big] \\
 & + \big[4\cdot \sigma _{i-1,i-1,i-1,i}^{\left[\alpha  \atop{ \beta  \atop  \gamma  } \right]\times k}  
   - 8\cdot \sigma _{i-2,i-2,i-2,i-1}^{\left[\alpha  \atop{ \beta  \atop \gamma   } \right]\times k} \big]  
  + \big[ 2\cdot \sigma _{i-1,i-1,i-1,i}^{\left[ \gamma  \atop{  \alpha \atop \beta   } \right]\times k}  
   - 4\cdot \sigma _{i-2,i-2,i-2,i-1}^{\left[ \gamma   \atop{  \alpha \atop \beta    } \right]\times k} \big] \nonumber \quad  \text{for $1 \leq i \leq \inf(k,3s+2m+l-2)+1 $}
\end{align}
\end{lem}
\begin{proof}
The rank formula follows from \eqref{eq 2.77}.
\end{proof}

    \begin{lem}
\label{lem 2.19}
We have the following rank formula \\
 \begin{align}
 \label{eq 2.82}
\boxed{   \sigma _{i-1,i,i,i}^{\left[ \gamma  \atop{  \beta \atop  \alpha  } \right]\times k}
  = 4\cdot \sigma _{i-1,i-1,i-1,i}^{\left[\alpha    \atop{ \beta \atop \gamma  } \right]\times k} \quad \text{for  $1\leq i\leq \inf(k,3s+2m+l-2)$}}
  \end{align}
   \end{lem}
 \begin{proof}
 We obtain by combining \eqref{eq 2.66},\eqref{eq 2.70} and \eqref{eq 2.81}.\\
 
  \begin{align*}
 \displaystyle 
 &  \sigma _{i-1,i,i,i}^{\left[ \gamma  \atop{  \beta \atop \alpha   } \right]\times k}
+  \sigma _{i-1,i-1,i,i}^{\left[ \gamma  \atop{  \beta \atop  \alpha  } \right]\times k}
+  \sigma _{i-2,i-1,i,i}^{\left[ \gamma  \atop{  \beta \atop \alpha  } \right]\times k}  
 = 4\cdot\big[ \sigma _{i-2,i-2,i-2,i-1}^{\left[\alpha   \atop{  \gamma \atop  \beta   } \right]\times k}
+ \sigma _{i-2,i-2,i-1,i-1}^{\left[\alpha   \atop{  \gamma \atop  \beta   } \right]\times k}
+ \sigma _{i-2,i-2,i-1,i}^{\left[\alpha   \atop{  \gamma \atop  \beta   } \right]\times k} \big] \\
 & + \big[4\cdot \sigma _{i-1,i-1,i-1,i}^{\left[\alpha  \atop{ \beta  \atop  \gamma  } \right]\times k}  
   - 8\cdot \sigma _{i-2,i-2,i-2,i-1}^{\left[\alpha  \atop{ \beta  \atop \gamma   } \right]\times k} \big]  
  + \big[ 2\cdot \sigma _{i-1,i-1,i-1,i}^{\left[ \gamma  \atop{  \alpha \atop \beta   } \right]\times k}  
   - 4\cdot \sigma _{i-2,i-2,i-2,i-1}^{\left[ \gamma   \atop{  \alpha \atop \beta    } \right]\times k} \big] \\
   & \Longrightarrow \\
&  \sigma _{i-1,i,i,i}^{\left[ \gamma  \atop{  \beta \atop  \alpha  } \right]\times k} +
2\cdot\big[  \sigma _{i-1,i-1,i-1,i}^{\left[ \gamma  \atop{  \alpha  \atop  \beta  } \right]\times k}
+  \sigma _{i-2,i-1,i-1,i}^{\left[ \gamma  \atop{  \alpha \atop  \beta  } \right]\times k}  \big] 
 = 4\cdot \sigma _{i-2,i-2,i-2,i-1}^{\left[\alpha   \atop{  \gamma \atop  \beta   } \right]\times k}
+2\cdot\big[ \sigma _{i-2,i-1,i-1,i-1}^{\left[ \gamma  \atop{  \alpha  \atop  \beta  } \right]\times k}
+  \sigma _{i-2,i-1,i-1,i}^{\left[ \gamma  \atop{  \alpha  \atop  \beta  } \right]\times k}\big] \\
&  + 4\cdot \sigma _{i-1,i-1,i-1,i}^{\left[\alpha    \atop{ \beta \atop \gamma  } \right]\times k}  
 -8\cdot \sigma _{i-2,i-2,i-2,i-1}^{\left[\alpha   \atop{ \beta \atop \gamma   } \right]\times k}
   +2\cdot \sigma _{i-1,i-1,i-1,i}^{\left[ \gamma  \atop{  \alpha \atop \beta   } \right]\times k}  
 -4\cdot \sigma _{i-2,i-2,i-2,i-1}^{\left[ \gamma   \atop{  \alpha \atop \beta   } \right]\times k} \\ 
  & \Longrightarrow  \\
&  \sigma _{i-1,i,i,i}^{\left[ \gamma  \atop{  \alpha \atop  \beta  } \right]\times k}
  - 4\cdot \sigma _{i-1,i-1,i-1,i}^{\left[\alpha    \atop{ \beta \atop \gamma  } \right]\times k} =
 2\cdot\big[ \sigma _{i-2,i-1,i-1,i-1}^{\left[ \gamma  \atop{  \alpha  \atop  \beta  } \right]\times k}  
  -4\cdot \sigma _{i-2,i-2,i-2,i-1}^{\left[\alpha   \atop{ \beta \atop \gamma   } \right]\times k}\big]  \quad \text{for  $2\leq i\leq \inf(k,3s+2m+l-1)$},\\ 
  & \\
  & \quad \text{which implies }\\
  & \\
   &  \quad  \sigma _{i-1,i,i,i}^{\left[ \gamma  \atop{  \beta \atop  \alpha  } \right]\times k}
  - 4\cdot \sigma _{i-1,i-1,i-1,i}^{\left[\alpha    \atop{ \beta \atop \gamma  } \right]\times k} =
 2^{i-1}\cdot\big[ \sigma _{0,1,1,1}^{\left[ \gamma  \atop{  \alpha  \atop  \beta  } \right]\times k}  
  -4\cdot \sigma _{0,0,0,1}^{\left[\alpha   \atop{ \beta \atop \gamma   } \right]\times k}\big] =0  \quad \text{for  $1\leq i\leq \inf(k,3s+2m+l-2)$}
  \end{align*}
 \end{proof}
  \subsection{A recurrent formula for the number of rank i matrices of the form  $\left[{A\over{B \over C}}\right] ,$ where A, B and C are persymmetric matrices over  $ \mathbb{F}_{2},$ }
  \label{subsec 7}
    \begin{lem}
 \label{lem 2.20} Let $ s\geq 2, \; m\geq 0,\;l\geq 0, \; k\geq 1 $ and  $ 0\leq i\leq \inf{(3s+2m+l,k)}. $ Then we have the following recurrent formula
  for the number $  \Gamma_{i}^{\left[s\atop{ s+m\atop s+m+l} \right]\times k} $ of rank i matrices 
 of the form   $\left[{A\over{B \over C}}\right] $ such that  A is a  $ s\times k $ persymmetric matrix over $ \mathbb{F}_{2},$ B a 
$(s+m)\times k$ persymmetric matrix and C  a $ (s+m+l)\times k $ persymmetric  matrix  \\
  \begin{align}
  \label{eq 2.83}
& \Gamma_{i}^{\left[s\atop{ s+m\atop s+m+l} \right]\times k}  \\
& = \big[2\cdot \Gamma_{i-1}^{\left[s -1 \atop{ s+m\atop s+m+l} \right]\times k}
 +4\cdot \Gamma_{i-1}^{\left[s\atop{ s+m-1\atop s+m+l} \right]\times k}  +8\cdot\Gamma_{i-1}^{\left[s\atop{ s+m\atop s+m+l-1} \right]\times k} \big ] 
 - \big[8\cdot \Gamma_{i-2}^{\left[s -1 \atop{ s+m-1\atop s+m+l} \right]\times k}
 +16\cdot \Gamma_{i-2}^{\left[s -1\atop{ s+m\atop s+m+l-1} \right]\times k}  +32\cdot\Gamma_{i-2}^{\left[s\atop{ s+m-1\atop s+m+l-1} \right]\times k}\big] \nonumber  \\
 & + 64\cdot \Gamma_{i-3}^{\left[s -1 \atop{ s+m-1\atop s+m+l-1} \right]\times k} + \Delta _{i}^{\left[s\atop{ s+m\atop s+m+l} \right]\times k} \nonumber \\
 & \nonumber
 \end{align}
where \\

\begin{align}
  \label{eq 2.84}
  \Delta _{i}^{\left[s\atop{ s+m\atop s+m+l} \right]\times k} 
 = \sigma _{i,i,i,i}^{\left[\alpha \atop{ \beta \atop \gamma } \right]\times k}
-7\cdot\sigma _{i-1,i-1,i-1,i-1}^{\left[\alpha \atop{ \beta \atop \gamma } \right]\times k}
+14\cdot \sigma _{i-2,i-2,i-2,i-2}^{\left[\alpha \atop{ \beta \atop \gamma } \right]\times k}
-8\cdot \sigma _{i-3,i-3,i-3,i-3}^{\left[\alpha \atop{ \beta \atop \gamma } \right]\times k}
 \end{align}
 Recall that 
$$  \sigma _{i,i,i,i}^{\left[\alpha \atop{ \beta \atop \gamma } \right]\times k}$$
 denotes  the cardinality of the following set 
 $$\begin{array}{l}\Big\{ (t,\eta,\xi  ) \in \mathbb{P}/\mathbb{P}_{k+s -1}\times
           \mathbb{P}/\mathbb{P}_{k+s+m-1}\times \mathbb{P}/\mathbb{P}_{k+s+m+l-1}
\mid r(  D^{\left[s-1\atop{ s+m-1\atop s+m+l-1} \right]\times k}(t,\eta,\xi ) ) = i \\
 r(  D^{\left[s\atop{ s+m-1\atop s+m+l-1} \right]\times k}(t,\eta,\xi ) )= i,\quad  
 r(  D^{\left[s\atop{ s+m\atop s+m+l-1} \right]\times k}(t,\eta,\xi ) ) = i,  \quad
  r(  D^{\left[s\atop{ s+m\atop s+m+l} \right]\times k}(t,\eta,\xi ) ) = i \Big\}
    \end{array}$$  
 \end{lem}
\begin{proof} 
 We have obviously\\
   
\begin{itemize}
\item
\begin{align*}
 & \Gamma_{i}^{\left[s\atop{ s+m\atop s+m+l} \right]\times k} \\
 & =\sigma _{i,i,i,i}^{\left[\alpha \atop{ \beta \atop \gamma } \right]\times k}+\sigma _{i-1,i,i,i}^{\left[\alpha \atop{ \beta \atop \gamma } \right]\times k}
+  \sigma _{i-2,i-1,i,i}^{\left[\alpha \atop{ \beta \atop \gamma } \right]\times k} +\sigma _{i-1,i-1,i,i}^{\left[\alpha \atop{ \beta \atop \gamma } \right]\times k}\\
& +\sigma _{i-1,i-1,i-1,i}^{\left[\alpha \atop{ \beta \atop \gamma } \right]\times k} +\sigma _{i-2,i-1,i-1,i}^{\left[\alpha \atop{ \beta \atop \gamma } \right]\times k}
+\sigma _{i-2,i-2,i-1,i}^{\left[\alpha \atop{ \beta \atop \gamma } \right]\times k} +\sigma _{i-3,i-2,i-1,i}^{\left[\alpha \atop{ \beta \atop \gamma } \right]\times k}\\
& =\sigma _{i,i,i,i}^{\left[\alpha \atop{ \gamma \atop \beta } \right]\times k}+\sigma _{i-1,i,i,i}^{\left[\alpha \atop{ \gamma  \atop \beta } \right]\times k}
+  \sigma _{i-2,i-1,i,i}^{\left[\alpha \atop{ \gamma  \atop \beta  } \right]\times k} +\sigma _{i-1,i-1,i,i}^{\left[\alpha \atop{ \gamma  \atop \beta } \right]\times k}\\
& +\sigma _{i-1,i-1,i-1,i}^{\left[\alpha \atop{ \gamma  \atop \beta  } \right]\times k} +\sigma _{i-2,i-1,i-1,i}^{\left[\alpha \atop{ \gamma  \atop \beta  } \right]\times k}
+\sigma _{i-2,i-2,i-1,i}^{\left[\alpha \atop{ \gamma \atop \beta  } \right]\times k} +\sigma _{i-3,i-2,i-1,i}^{\left[\alpha \atop{ \gamma \atop \beta  } \right]\times k}\\
 & =\sigma _{i,i,i,i}^{\left[\beta \atop{ \alpha  \atop \gamma } \right]\times k}+\sigma _{i-1,i,i,i}^{\left[\beta \atop{ \alpha \atop \gamma } \right]\times k}
+  \sigma _{i-2,i-1,i,i}^{\left[\beta  \atop{ \alpha  \atop \gamma } \right]\times k} +\sigma _{i-1,i-1,i,i}^{\left[\beta  \atop{ \alpha  \atop \gamma } \right]\times k}\\
& +\sigma _{i-1,i-1,i-1,i}^{\left[\beta  \atop{ \alpha \atop \gamma } \right]\times k} +\sigma _{i-2,i-1,i-1,i}^{\left[\beta \atop{ \alpha \atop \gamma } \right]\times k}
+\sigma _{i-2,i-2,i-1,i}^{\left[\beta \atop{ \alpha \atop \gamma } \right]\times k} +\sigma _{i-3,i-2,i-1,i}^{\left[\beta \atop{ \alpha  \atop \gamma } \right]\times k}\\
 & =\sigma _{i,i,i,i}^{\left[\beta \atop{ \gamma  \atop \alpha } \right]\times k}+\sigma _{i-1,i,i,i}^{\left[\beta \atop{ \gamma \atop \alpha } \right]\times k}
+  \sigma _{i-2,i-1,i,i}^{\left[\beta  \atop{ \gamma   \atop \alpha  } \right]\times k} +\sigma _{i-1,i-1,i,i}^{\left[\beta  \atop{ \gamma  \atop \alpha } \right]\times k}\\
& +\sigma _{i-1,i-1,i-1,i}^{\left[\beta  \atop{ \gamma  \atop \alpha } \right]\times k} +\sigma _{i-2,i-1,i-1,i}^{\left[\beta \atop{ \gamma  \atop \alpha  } \right]\times k}
+\sigma _{i-2,i-2,i-1,i}^{\left[\beta \atop{ \gamma  \atop \alpha } \right]\times k} +\sigma _{i-3,i-2,i-1,i}^{\left[\beta \atop{ \gamma  \atop \alpha } \right]\times k}\\
 & =\sigma _{i,i,i,i}^{\left[\gamma \atop{ \beta   \atop \alpha } \right]\times k}+\sigma _{i-1,i,i,i}^{\left[\gamma \atop{ \beta \atop \alpha } \right]\times k}
+  \sigma _{i-2,i-1,i,i}^{\left[\gamma  \atop{ \beta    \atop \alpha  } \right]\times k} +\sigma _{i-1,i-1,i,i}^{\left[\gamma  \atop{ \beta  \atop \alpha } \right]\times k}\\
& +\sigma _{i-1,i-1,i-1,i}^{\left[\gamma  \atop{ \beta  \atop \alpha } \right]\times k} +\sigma _{i-2,i-1,i-1,i}^{\left[\gamma  \atop{ \beta   \atop \alpha  } \right]\times k}
+\sigma _{i-2,i-2,i-1,i}^{\left[\gamma \atop{ \beta  \atop \alpha } \right]\times k} +\sigma _{i-3,i-2,i-1,i}^{\left[\gamma  \atop{ \beta  \atop \alpha } \right]\times k}\\
 & =\sigma _{i,i,i,i}^{\left[\gamma \atop{ \alpha   \atop \beta } \right]\times k}+\sigma _{i-1,i,i,i}^{\left[\gamma \atop{ \alpha  \atop \beta  } \right]\times k}
+  \sigma _{i-2,i-1,i,i}^{\left[\gamma  \atop{ \alpha    \atop \beta  } \right]\times k} +\sigma _{i-1,i-1,i,i}^{\left[\gamma  \atop{ \alpha  \atop \beta } \right]\times k}\\
& +\sigma _{i-1,i-1,i-1,i}^{\left[\gamma  \atop{ \alpha   \atop \beta } \right]\times k} +\sigma _{i-2,i-1,i-1,i}^{\left[\gamma  \atop{ \alpha   \atop \beta  } \right]\times k}
+\sigma _{i-2,i-2,i-1,i}^{\left[\gamma \atop{ \alpha  \atop \beta } \right]\times k} +\sigma _{i-3,i-2,i-1,i}^{\left[\gamma  \atop{ \alpha  \atop \beta } \right]\times k}
\end{align*}
\item
\begin{align*}
 & -4\cdot \Gamma_{i-1}^{\left[s\atop{ s+m-1\atop s+m+l} \right]\times k} \\
 & = -2\cdot\left[  \sigma _{i-1,i-1,i-1,i-1}^{\left[\alpha \atop{ \gamma \atop \beta  } \right]\times k}
 +\sigma _{i-2,i-1,i-1,i-1}^{\left[\alpha \atop{ \gamma \atop \beta } \right]\times k}
+  \sigma _{i-2,i-2,i-1,i-1}^{\left[\alpha \atop{ \gamma  \atop \beta } \right]\times k}
 +\sigma _{i-3,i-2,i-1,i-1}^{\left[\alpha \atop{ \gamma  \atop \beta } \right]\times k}\right]\\
 &  -2\cdot\left[   \sigma _{i-1,i-1,i-1,i}^{\left[\alpha \atop{ \gamma  \atop \beta  } \right]\times k}
 +\sigma _{i-2,i-1,i-1,i}^{\left[\alpha \atop{ \gamma \atop \beta  } \right]\times k}
+ \sigma _{i-2,i-2,i-1,i}^{\left[\alpha \atop{ \gamma  \atop \beta } \right]\times k}
 +\sigma _{i-3,i-2,i-1,i}^{\left[\alpha \atop{ \gamma \atop \beta  } \right]\times k}   \right]\\
  & = -2\cdot\left[  \sigma _{i-1,i-1,i-1,i-1}^{\left[\gamma  \atop{ \alpha  \atop \beta  } \right]\times k}
 +\sigma _{i-2,i-1,i-1,i-1}^{\left[\gamma  \atop{ \alpha  \atop \beta } \right]\times k}
+  \sigma _{i-2,i-2,i-1,i-1}^{\left[\gamma  \atop{ \alpha  \atop \beta } \right]\times k}
 +\sigma _{i-3,i-2,i-1,i-1}^{\left[\gamma  \atop{ \alpha  \atop \beta } \right]\times k}\right]\\
 &  -2\cdot\left[   \sigma _{i-1,i-1,i-1,i}^{\left[\gamma \atop{ \alpha  \atop \beta  } \right]\times k}
 +\sigma _{i-2,i-1,i-1,i}^{\left[\gamma  \atop{ \alpha \atop \beta  } \right]\times k}
+ \sigma _{i-2,i-2,i-1,i}^{\left[\gamma  \atop{ \alpha  \atop \beta } \right]\times k}
 +\sigma _{i-3,i-2,i-1,i}^{\left[\gamma  \atop{ \alpha \atop \beta  } \right]\times k}   \right]
\end{align*}
\item
\begin{align*}
 & -8\cdot \Gamma_{i-1}^{\left[s\atop{ s+m\atop s+m+l-1} \right]\times k} \\
 & = -4\cdot\left[  \sigma _{i-1,i-1,i-1,i-1}^{\left[\alpha \atop{ \beta  \atop \gamma  } \right]\times k}
 +\sigma _{i-2,i-1,i-1,i-1}^{\left[\alpha \atop{ \beta  \atop \gamma } \right]\times k}
+  \sigma _{i-2,i-2,i-1,i-1}^{\left[\alpha \atop{ \beta   \atop \gamma  } \right]\times k}
 +\sigma _{i-3,i-2,i-1,i-1}^{\left[\alpha \atop{ \beta  \atop \gamma  } \right]\times k}\right]\\
 &  -4\cdot\left[   \sigma _{i-1,i-1,i-1,i}^{\left[\alpha \atop{\beta   \atop \gamma  } \right]\times k}
 +\sigma _{i-2,i-1,i-1,i}^{\left[\alpha \atop{ \beta  \atop \gamma  } \right]\times k}
+ \sigma _{i-2,i-2,i-1,i}^{\left[\alpha \atop{ \beta   \atop \gamma } \right]\times k}
 +\sigma _{i-3,i-2,i-1,i}^{\left[\alpha \atop{ \beta \atop \gamma  } \right]\times k}   \right]\\
  & = -4\cdot\left[  \sigma _{i-1,i-1,i-1,i-1}^{\left[\beta   \atop{ \alpha  \atop \gamma  } \right]\times k}
 +\sigma _{i-2,i-1,i-1,i-1}^{\left[\beta  \atop{ \alpha  \atop \gamma } \right]\times k}
+  \sigma _{i-2,i-2,i-1,i-1}^{\left[\beta   \atop{ \alpha  \atop \gamma  } \right]\times k}
 +\sigma _{i-3,i-2,i-1,i-1}^{\left[ \beta \atop{ \alpha  \atop \gamma } \right]\times k}\right]\\
 &  -4\cdot\left[   \sigma _{i-1,i-1,i-1,i}^{\left[\beta \atop{ \alpha  \atop \gamma  } \right]\times k}
 +\sigma _{i-2,i-1,i-1,i}^{\left[\beta  \atop{ \alpha \atop \gamma  } \right]\times k}
+ \sigma _{i-2,i-2,i-1,i}^{\left[\beta  \atop{ \alpha  \atop \gamma } \right]\times k}
 +\sigma _{i-3,i-2,i-1,i}^{\left[\beta  \atop{ \alpha \atop \gamma  } \right]\times k}   \right]
\end{align*}
\item
\begin{align*}
 & 32\cdot \Gamma_{i-2}^{\left[s\atop{ s+m-1\atop s+m+l-1} \right]\times k} \\
 & = 8\cdot\left[  \sigma _{i-2,i-2,i-2,i -2}^{\left[\alpha \atop{ \beta  \atop \gamma  } \right]\times k}
 +\sigma _{i-2,i-2,i-2,i-1}^{\left[\alpha \atop{ \beta  \atop \gamma } \right]\times k}
+  \sigma _{i-2,i-2,i-1,i-1}^{\left[\alpha \atop{ \beta   \atop \gamma  } \right]\times k}
 +\sigma _{i-2,i-2,i-1,i}^{\left[\alpha \atop{ \beta  \atop \gamma  } \right]\times k}\right]\\
 & + 8\cdot\left[   \sigma _{i-3,i-2,i-1,i-1}^{\left[\alpha \atop{\beta   \atop \gamma  } \right]\times k}
 +\sigma _{i-3,i-2,i-1,i}^{\left[\alpha \atop{ \beta  \atop \gamma  } \right]\times k}
+ \sigma _{i-3,i-2,i-2,i-2}^{\left[\alpha \atop{ \beta   \atop \gamma } \right]\times k}
 +\sigma _{i-3,i-2,i-2,i-1}^{\left[\alpha \atop{ \beta \atop \gamma  } \right]\times k}   \right]\\
& = 8\cdot\left[  \sigma _{i-2,i-2,i-2,i -2}^{\left[\alpha \atop{ \gamma  \atop \beta  } \right]\times k}
 +\sigma _{i-2,i-2,i-2,i-1}^{\left[\alpha \atop{ \gamma  \atop \beta  } \right]\times k}
+  \sigma _{i-2,i-2,i-1,i-1}^{\left[\alpha \atop{ \gamma   \atop \beta  } \right]\times k}
 +\sigma _{i-2,i-2,i-1,i}^{\left[\alpha \atop{ \gamma   \atop \beta  } \right]\times k}\right]\\
 & + 8\cdot\left[   \sigma _{i-3,i-2,i-1,i-1}^{\left[\alpha \atop{\gamma   \atop \beta  } \right]\times k}
 +\sigma _{i-3,i-2,i-1,i}^{\left[\alpha \atop{ \gamma  \atop \beta  } \right]\times k}
+ \sigma _{i-3,i-2,i-2,i-2}^{\left[\alpha \atop{ \gamma    \atop \beta } \right]\times k}
 +\sigma _{i-3,i-2,i-2,i-1}^{\left[\alpha \atop{ \gamma  \atop \beta  } \right]\times k}   \right]
\end{align*}\vspace{0.1 cm}
\item
\begin{align*}
 & -2\cdot \Gamma_{i-1}^{\left[s-1\atop{ s+m\atop s+m+l} \right]\times k} \\
 & = -\left[  \sigma _{i-1,i-1,i-1,i-1}^{\left[\beta  \atop{ \gamma \atop \alpha  } \right]\times k}
 +\sigma _{i-2,i-1,i-1,i-1}^{\left[\beta  \atop{ \gamma \atop \alpha } \right]\times k}
+  \sigma _{i-2,i-2,i-1,i-1}^{\left[\beta  \atop{ \gamma  \atop \alpha } \right]\times k}
 +\sigma _{i-3,i-2,i-1,i-1}^{\left[\beta  \atop{ \gamma  \atop \alpha } \right]\times k}\right]\\
 &  -\left[   \sigma _{i-1,i-1,i-1,i}^{\left[\beta \atop{ \gamma  \atop \alpha  } \right]\times k}
 +\sigma _{i-2,i-1,i-1,i}^{\left[\beta  \atop{ \gamma \atop \alpha   } \right]\times k}
+ \sigma _{i-2,i-2,i-1,i}^{\left[\beta  \atop{ \gamma  \atop \alpha } \right]\times k}
 +\sigma _{i-3,i-2,i-1,i}^{\left[\beta  \atop{ \gamma \atop \alpha  } \right]\times k}   \right]\\
  & = -\left[  \sigma _{i-1,i-1,i-1,i-1}^{\left[\gamma  \atop{ \beta   \atop \alpha  } \right]\times k}
 +\sigma _{i-2,i-1,i-1,i-1}^{\left[\gamma  \atop{ \beta  \atop \alpha } \right]\times k}
+  \sigma _{i-2,i-2,i-1,i-1}^{\left[\gamma  \atop{ \beta   \atop \alpha } \right]\times k}
 +\sigma _{i-3,i-2,i-1,i-1}^{\left[\gamma  \atop{ \beta   \atop \alpha  } \right]\times k}\right]\\
 &  -\left[   \sigma _{i-1,i-1,i-1,i}^{\left[\gamma \atop{ \beta  \atop \alpha  } \right]\times k}
 +\sigma _{i-2,i-1,i-1,i}^{\left[\gamma  \atop{ \beta  \atop \alpha  } \right]\times k}
+ \sigma _{i-2,i-2,i-1,i}^{\left[\gamma  \atop{ \beta  \atop \alpha } \right]\times k}
 +\sigma _{i-3,i-2,i-1,i}^{\left[\gamma  \atop{ \beta  \atop \alpha  } \right]\times k}   \right]
\end{align*}
\item
\begin{align*}
 & 8\cdot \Gamma_{i-2}^{\left[s -1\atop{ s+m-1\atop s+m+l} \right]\times k} \\
 & = 2\cdot\left[  \sigma _{i-2,i-2,i-2,i -2}^{\left[\gamma  \atop{ \beta  \atop \alpha  } \right]\times k}
 +\sigma _{i-2,i-2,i-2,i-1}^{\left[\gamma  \atop{ \beta  \atop \alpha  } \right]\times k}
+  \sigma _{i-2,i-2,i-1,i-1}^{\left[\gamma  \atop{ \beta   \atop \alpha   } \right]\times k}
 +\sigma _{i-2,i-2,i-1,i}^{\left[\gamma  \atop{ \beta  \atop \alpha  } \right]\times k}\right]\\
 & + 2\cdot\left[   \sigma _{i-3,i-2,i-1,i-1}^{\left[\gamma \atop{\beta   \atop \alpha  } \right]\times k}
 +\sigma _{i-3,i-2,i-1,i}^{\left[\gamma  \atop{ \beta  \atop \alpha   } \right]\times k}
+ \sigma _{i-3,i-2,i-2,i-2}^{\left[\gamma \atop{ \beta   \atop \alpha  } \right]\times k}
 +\sigma _{i-3,i-2,i-2,i-1}^{\left[\gamma  \atop{ \beta \atop \alpha  } \right]\times k}   \right]\\
& = 2\cdot\left[  \sigma _{i-2,i-2,i-2,i -2}^{\left[\gamma  \atop{ \alpha  \atop \beta  } \right]\times k}
 +\sigma _{i-2,i-2,i-2,i-1}^{\left[\gamma \atop{ \alpha  \atop \beta  } \right]\times k}
+  \sigma _{i-2,i-2,i-1,i-1}^{\left[\gamma  \atop{ \alpha   \atop \beta  } \right]\times k}
 +\sigma _{i-2,i-2,i-1,i}^{\left[\gamma  \atop{ \alpha    \atop \beta  } \right]\times k}\right]\\
 & + 2\cdot\left[   \sigma _{i-3,i-2,i-1,i-1}^{\left[\gamma \atop{\alpha   \atop \beta  } \right]\times k}
 +\sigma _{i-3,i-2,i-1,i}^{\left[\gamma  \atop{ \alpha  \atop \beta  } \right]\times k}
+ \sigma _{i-3,i-2,i-2,i-2}^{\left[\gamma  \atop{ \alpha    \atop \beta } \right]\times k}
 +\sigma _{i-3,i-2,i-2,i-1}^{\left[\gamma \atop{ \alpha  \atop \beta  } \right]\times k}   \right]
\end{align*}
\item
\begin{align*}
 & 16\cdot \Gamma_{i-2}^{\left[s -1\atop{ s+m\atop s+m+l-1} \right]\times k} \\
 & = 4\cdot\left[  \sigma _{i-2,i-2,i-2,i -2}^{\left[\beta  \atop{ \alpha   \atop \gamma  } \right]\times k}
 +\sigma _{i-2,i-2,i-2,i-1}^{\left[\beta   \atop{\alpha   \atop \gamma  } \right]\times k}
+  \sigma _{i-2,i-2,i-1,i-1}^{\left[\beta  \atop{ \alpha   \atop \gamma   } \right]\times k}
 +\sigma _{i-2,i-2,i-1,i}^{\left[\beta  \atop{ \alpha  \atop \gamma  } \right]\times k}\right]\\
 & + 4\cdot\left[   \sigma _{i-3,i-2,i-1,i-1}^{\left[\beta  \atop{\alpha    \atop \gamma  } \right]\times k}
 +\sigma _{i-3,i-2,i-1,i}^{\left[\beta  \atop{ \alpha   \atop \gamma   } \right]\times k}
+ \sigma _{i-3,i-2,i-2,i-2}^{\left[\beta \atop{ \alpha   \atop \gamma   } \right]\times k}
 +\sigma _{i-3,i-2,i-2,i-1}^{\left[\beta  \atop{ \alpha \atop \gamma  } \right]\times k}   \right]\\
 & = 4\cdot\left[  \sigma _{i-2,i-2,i-2,i -2}^{\left[\beta  \atop{ \gamma   \atop \alpha   } \right]\times k}
 +\sigma _{i-2,i-2,i-2,i-1}^{\left[\beta   \atop{\gamma   \atop \alpha  } \right]\times k}
+  \sigma _{i-2,i-2,i-1,i-1}^{\left[\beta  \atop{ \gamma   \atop \alpha   } \right]\times k}
 +\sigma _{i-2,i-2,i-1,i}^{\left[\beta  \atop{ \gamma  \atop \alpha  } \right]\times k}\right]\\
 & + 4\cdot\left[   \sigma _{i-3,i-2,i-1,i-1}^{\left[\beta  \atop{\gamma    \atop \alpha  } \right]\times k}
 +\sigma _{i-3,i-2,i-1,i}^{\left[\beta  \atop{ \gamma   \atop \alpha   } \right]\times k}
+ \sigma _{i-3,i-2,i-2,i-2}^{\left[\beta \atop{ \gamma    \atop \alpha   } \right]\times k}
 +\sigma _{i-3,i-2,i-2,i-1}^{\left[\beta  \atop{ \gamma  \atop \alpha  } \right]\times k}   \right]
\end{align*}
\item
\begin{align*}
 & -64\cdot \Gamma_{i-3}^{\left[s-1\atop{ s+m-1\atop s+m+l-1} \right]\times k} \\
 & =-8\cdot\left[\sigma _{i-3,i-3,i-3,i-3}^{\left[\alpha \atop{ \beta \atop \gamma } \right]\times k}
+ \sigma _{i-3,i-3,i-3,i-2}^{\left[\alpha \atop{ \beta \atop \gamma } \right]\times k}
+  \sigma _{i-3,i-3,i-2,i-2}^{\left[\alpha \atop{ \beta \atop \gamma } \right]\times k}
 +\sigma _{i-3,i-3,i-2,i-1}^{\left[\alpha \atop{ \beta \atop \gamma } \right]\times k}\right] \\
&-8\cdot\left[ \sigma _{i-3,i-2,i-2,i-2}^{\left[\alpha \atop{ \beta \atop \gamma } \right]\times k} 
+ \sigma _{i-3,i-2,i-2,i-1}^{\left[\alpha \atop{ \beta \atop \gamma } \right]\times k}
+\sigma _{i-3,i-2,i-1,i-1}^{\left[\alpha \atop{ \beta \atop \gamma } \right]\times k}
 +\sigma _{i-3,i-2,i-1,i}^{\left[\alpha \atop{ \beta \atop \gamma } \right]\times k}\right]\\
& = -8\cdot\left[\sigma _{i-3,i-3,i-3,i-3}^{\left[\alpha \atop{ \gamma \atop \beta } \right]\times k}
+\sigma _{i-3,i-3,i-3,i-2}^{\left[\alpha \atop{ \gamma  \atop \beta } \right]\times k}
+  \sigma _{i-3,i-3,i-2,i-2}^{\left[\alpha \atop{ \gamma  \atop \beta  } \right]\times k}
 +\sigma _{i-3,i-3,i-2,i-1}^{\left[\alpha \atop{ \gamma  \atop \beta } \right]\times k}\right] \\
& -8\cdot\left[\sigma _{i-3,i-2,i-2,i-2}^{\left[\alpha \atop{ \gamma  \atop \beta  } \right]\times k}
 +\sigma _{i-3,i-2,i-2,i-1}^{\left[\alpha \atop{ \gamma  \atop \beta  } \right]\times k}
+\sigma _{i-3,i-2,i-1,i-1}^{\left[\alpha \atop{ \gamma \atop \beta  } \right]\times k}
 +\sigma _{i-3,i-2,i-1,i}^{\left[\alpha \atop{ \gamma \atop \beta  } \right]\times k}\right]\\
  & = -8\cdot\left[\sigma _{i-3,i-3,i-3,i-3}^{\left[\beta \atop{ \alpha  \atop \gamma } \right]\times k}
 +\sigma _{i-3,i-3,i-3,i-2}^{\left[\beta \atop{ \alpha \atop \gamma } \right]\times k}
+  \sigma _{i-3,i-3,i-2,i-2}^{\left[\beta  \atop{ \alpha  \atop \gamma } \right]\times k}
 +\sigma _{i-3,i-3,i-2,i-1}^{\left[\beta  \atop{ \alpha  \atop \gamma } \right]\times k}\right] \\
& -8\cdot\left[ \sigma _{i-3,i-2,i-2,i-2}^{\left[\beta  \atop{ \alpha \atop \gamma } \right]\times k}
 +\sigma _{i-3,i-2,i-2,i-1}^{\left[\beta \atop{ \alpha \atop \gamma } \right]\times k}
+\sigma _{i-3,i-2,i-1,i-1}^{\left[\beta \atop{ \alpha \atop \gamma } \right]\times k}
 +\sigma _{i-3,i-2,i-1,i}^{\left[\beta \atop{ \alpha  \atop \gamma } \right]\times k}\right]\\
 & = -8\cdot\left[\sigma _{i-3,i-3,i-3,i-3}^{\left[\beta \atop{ \gamma  \atop \alpha } \right]\times k}
 +\sigma _{i-3,i-3,i-3,i-2}^{\left[\beta \atop{ \gamma \atop \alpha } \right]\times k}
+  \sigma _{i-3,i-3,i-2,i-2}^{\left[\beta  \atop{ \gamma   \atop \alpha  } \right]\times k}
 +\sigma _{i-3,i-3,i-2,i-1}^{\left[\beta  \atop{ \gamma  \atop \alpha } \right]\times k}\right] \\
 & -8\cdot\left[\sigma _{i-3,i-2,i-2,i-2}^{\left[\beta  \atop{ \gamma  \atop \alpha } \right]\times k}
 +\sigma _{i-3,i-2,i-2,i-1}^{\left[\beta \atop{ \gamma  \atop \alpha  } \right]\times k}
+\sigma _{i-3,i-2,i-1,i-1}^{\left[\beta \atop{ \gamma  \atop \alpha } \right]\times k}
 +\sigma _{i-3,i-2,i-1,i}^{\left[\beta \atop{ \gamma  \atop \alpha } \right]\times k}\right]\\
 & =  -8\cdot\left[\sigma _{i-3,i-3,i-3,i-3}^{\left[\gamma \atop{ \beta   \atop \alpha } \right]\times k}
 +\sigma _{i-3,i-3,i-3,i-2}^{\left[\gamma \atop{ \beta \atop \alpha } \right]\times k}
+  \sigma _{i-3,i-3,i-2,i-2}^{\left[\gamma  \atop{ \beta    \atop \alpha  } \right]\times k}
 +\sigma _{i-3,i-3,i-2,i-1}^{\left[\gamma  \atop{ \beta  \atop \alpha } \right]\times k}\right] \\
 & -8\cdot\left[\sigma _{i-3,i-2,i-2,i-2}^{\left[\gamma  \atop{ \beta  \atop \alpha } \right]\times k}
 +\sigma _{i-3,i-2,i-2,i-1}^{\left[\gamma  \atop{ \beta   \atop \alpha  } \right]\times k}
+\sigma _{i-3,i-2,i-1,i-1}^{\left[\gamma \atop{ \beta  \atop \alpha } \right]\times k}
 +\sigma _{i-3,i-2,i-1,i}^{\left[\gamma  \atop{ \beta  \atop \alpha } \right]\times k}\right]\\
 & = -8\cdot\left[\sigma _{i-3,i-3,i-3,i-3}^{\left[\gamma \atop{ \alpha   \atop \beta } \right]\times k}
 +\sigma _{i-3,i-3,i-3,i-2}^{\left[\gamma \atop{ \alpha  \atop \beta  } \right]\times k}
+  \sigma _{i-3,i-3,i-2,i-2}^{\left[\gamma  \atop{ \alpha    \atop \beta  } \right]\times k}
 +\sigma _{i-3,i-3,i-2,i-1}^{\left[\gamma  \atop{ \alpha  \atop \beta } \right]\times k}\right] \\
 & -8\cdot\left[\sigma _{i-3,i-2,i-2,i-2}^{\left[\gamma  \atop{ \alpha   \atop \beta } \right]\times k}
 +\sigma _{i-3,i-2,i-2,i-1}^{\left[\gamma  \atop{ \alpha   \atop \beta  } \right]\times k}
+\sigma _{i-3,i-2,i-1,i-1}^{\left[\gamma \atop{ \alpha  \atop \beta } \right]\times k}
 +\sigma _{i-3,i-2,i-1,i}^{\left[\gamma  \atop{ \alpha  \atop \beta } \right]\times k}\right]
\end{align*}
\end{itemize}
For every term in the recurrent formula \eqref{eq 2.83} we choose in view of
 \eqref{eq 2.66},\eqref{eq 2.74} and \eqref{eq 2.82} an appropriate permutation of $ \{\alpha, \beta, \gamma \},$
 hence we get\\
 
  \begin{align*}
& \Gamma_{i}^{\left[s\atop{ s+m\atop s+m+l} \right]\times k}  
- \big[2\cdot \Gamma_{i-1}^{\left[s -1 \atop{ s+m\atop s+m+l} \right]\times k}
 +4\cdot \Gamma_{i-1}^{\left[s\atop{ s+m-1\atop s+m+l} \right]\times k}  +8\cdot\Gamma_{i-1}^{\left[s\atop{ s+m\atop s+m+l-1} \right]\times k} \big ]\\
& + \big[8\cdot \Gamma_{i-2}^{\left[s -1 \atop{ s+m-1\atop s+m+l} \right]\times k}
 +16\cdot \Gamma_{i-2}^{\left[s -1\atop{ s+m\atop s+m+l-1} \right]\times k}  +32\cdot\Gamma_{i-2}^{\left[s\atop{ s+m-1\atop s+m+l-1} \right]\times k}\big] 
  - 64\cdot \Gamma_{i-3}^{\left[s -1 \atop{ s+m-1\atop s+m+l-1} \right]\times k}
    -  \Delta _{i}^{\left[s\atop{ s+m\atop s+m+l} \right]\times k} = 0 \\
     & \Longleftrightarrow \\
 & \left[\sigma _{i,i,i,i}^{\left[\gamma \atop{ \beta   \atop \alpha } \right]\times k}+\sigma _{i-1,i,i,i}^{\left[\gamma \atop{ \beta \atop \alpha } \right]\times k}
+  \sigma _{i-2,i-1,i,i}^{\left[\gamma  \atop{ \beta    \atop \alpha  } \right]\times k} +\sigma _{i-1,i-1,i,i}^{\left[\gamma  \atop{ \beta  \atop \alpha } \right]\times k}
 +\sigma _{i-1,i-1,i-1,i}^{\left[\gamma  \atop{ \beta  \atop \alpha } \right]\times k} +\sigma _{i-2,i-1,i-1,i}^{\left[\gamma  \atop{ \beta   \atop \alpha  } \right]\times k}
+\sigma _{i-2,i-2,i-1,i}^{\left[\gamma \atop{ \beta  \atop \alpha } \right]\times k} +\sigma _{i-3,i-2,i-1,i}^{\left[\gamma  \atop{ \beta  \atop \alpha } \right]\times k}\right]\\
  & - \left[\sigma _{i-1,i-1,i-1,i-1}^{\left[\gamma  \atop{ \beta   \atop \alpha  } \right]\times k}
 +\sigma _{i-2,i-1,i-1,i-1}^{\left[\gamma  \atop{ \beta  \atop \alpha } \right]\times k}
+ \sigma _{i-2,i-2,i-1,i-1}^{\left[\gamma  \atop{ \beta   \atop \alpha } \right]\times k}
 +\sigma _{i-3,i-2,i-1,i-1}^{\left[\gamma  \atop{ \beta   \atop \alpha  } \right]\times k}\right]\\
  &   - \left[\sigma _{i-1,i-1,i-1,i}^{\left[\gamma \atop{ \beta  \atop \alpha  } \right]\times k}
 + \sigma _{i-2,i-1,i-1,i}^{\left[\gamma  \atop{ \beta  \atop \alpha  } \right]\times k}
+ \sigma _{i-2,i-2,i-1,i}^{\left[\gamma  \atop{ \beta  \atop \alpha } \right]\times k}
  + \sigma _{i-3,i-2,i-1,i}^{\left[\gamma  \atop{ \beta  \atop \alpha  } \right]\times k}\right] \\
   &  -2\cdot\left[  \sigma _{i-1,i-1,i-1,i-1}^{\left[\gamma  \atop{ \alpha  \atop \beta  } \right]\times k}
 +\sigma _{i-2,i-1,i-1,i-1}^{\left[\gamma  \atop{ \alpha  \atop \beta } \right]\times k}
+  \sigma _{i-2,i-2,i-1,i-1}^{\left[\gamma  \atop{ \alpha  \atop \beta } \right]\times k}
 +\sigma _{i-3,i-2,i-1,i-1}^{\left[\gamma  \atop{ \alpha  \atop \beta } \right]\times k}\right]\\
 &  -2\cdot\left[   \sigma _{i-1,i-1,i-1,i}^{\left[\gamma \atop{ \alpha  \atop \beta  } \right]\times k}
 +\sigma _{i-2,i-1,i-1,i}^{\left[\gamma  \atop{ \alpha \atop \beta  } \right]\times k}
+ \sigma _{i-2,i-2,i-1,i}^{\left[\gamma  \atop{ \alpha  \atop \beta } \right]\times k}
 +\sigma _{i-3,i-2,i-1,i}^{\left[\gamma  \atop{ \alpha \atop \beta  } \right]\times k}   \right]\\
  &  -4\cdot\left[  \sigma _{i-1,i-1,i-1,i-1}^{\left[\beta   \atop{ \alpha  \atop \gamma  } \right]\times k}
 +\sigma _{i-2,i-1,i-1,i-1}^{\left[\beta  \atop{ \alpha  \atop \gamma } \right]\times k}
+  \sigma _{i-2,i-2,i-1,i-1}^{\left[\beta   \atop{ \alpha  \atop \gamma  } \right]\times k}
 +\sigma _{i-3,i-2,i-1,i-1}^{\left[ \beta \atop{ \alpha  \atop \gamma } \right]\times k}\right]\\
 &  -4\cdot\left[   \sigma _{i-1,i-1,i-1,i}^{\left[\beta \atop{ \alpha  \atop \gamma  } \right]\times k}
 +\sigma _{i-2,i-1,i-1,i}^{\left[\beta  \atop{ \alpha \atop \gamma  } \right]\times k}
+ \sigma _{i-2,i-2,i-1,i}^{\left[\beta  \atop{ \alpha  \atop \gamma } \right]\times k}
 +\sigma _{i-3,i-2,i-1,i}^{\left[\beta  \atop{ \alpha \atop \gamma  } \right]\times k}   \right]\\
 & +  2\cdot\left[  \sigma _{i-2,i-2,i-2,i -2}^{\left[\gamma  \atop{ \alpha  \atop \beta  } \right]\times k}
  +\sigma _{i-2,i-2,i-2,i-1}^{\left[\gamma \atop{ \alpha  \atop \beta  } \right]\times k}
+  \sigma _{i-2,i-2,i-1,i-1}^{\left[\gamma  \atop{ \alpha   \atop \beta  } \right]\times k}
 +\sigma _{i-2,i-2,i-1,i}^{\left[\gamma  \atop{ \alpha    \atop \beta  } \right]\times k}\right]\\
 & + 2\cdot\left[   \sigma _{i-3,i-2,i-1,i-1}^{\left[\gamma \atop{\alpha   \atop \beta  } \right]\times k}
 +\sigma _{i-3,i-2,i-1,i}^{\left[\gamma  \atop{ \alpha  \atop \beta  } \right]\times k}
+ \sigma _{i-3,i-2,i-2,i-2}^{\left[\gamma  \atop{ \alpha    \atop \beta } \right]\times k}
 +\sigma _{i-3,i-2,i-2,i-1}^{\left[\gamma \atop{ \alpha  \atop \beta  } \right]\times k}   \right]\\
  & + 4\cdot\left[  \sigma _{i-2,i-2,i-2,i -2}^{\left[\beta  \atop{ \alpha   \atop \gamma  } \right]\times k}
 +\sigma _{i-2,i-2,i-2,i-1}^{\left[\beta   \atop{\alpha   \atop \gamma  } \right]\times k}
+  \sigma _{i-2,i-2,i-1,i-1}^{\left[\beta  \atop{ \alpha   \atop \gamma   } \right]\times k}
 +\sigma _{i-2,i-2,i-1,i}^{\left[\beta  \atop{ \alpha  \atop \gamma  } \right]\times k}\right]\\
 & + 4\cdot\left[   \sigma _{i-3,i-2,i-1,i-1}^{\left[\beta  \atop{\alpha    \atop \gamma  } \right]\times k}
 +\sigma _{i-3,i-2,i-1,i}^{\left[\beta  \atop{ \alpha   \atop \gamma   } \right]\times k}
+ \sigma _{i-3,i-2,i-2,i-2}^{\left[\beta \atop{ \alpha   \atop \gamma   } \right]\times k}
 +\sigma _{i-3,i-2,i-2,i-1}^{\left[\beta  \atop{ \alpha \atop \gamma  } \right]\times k}   \right]\\
  & + 8\cdot\left[  \sigma _{i-2,i-2,i-2,i -2}^{\left[\alpha \atop{ \beta  \atop \gamma  } \right]\times k}
 +\sigma _{i-2,i-2,i-2,i-1}^{\left[\alpha \atop{ \beta  \atop \gamma } \right]\times k}
+  \sigma _{i-2,i-2,i-1,i-1}^{\left[\alpha \atop{ \beta   \atop \gamma  } \right]\times k}
 +\sigma _{i-2,i-2,i-1,i}^{\left[\alpha \atop{ \beta  \atop \gamma  } \right]\times k}\right]\\
 & + 8\cdot\left[   \sigma _{i-3,i-2,i-1,i-1}^{\left[\alpha \atop{\beta   \atop \gamma  } \right]\times k}
 +\sigma _{i-3,i-2,i-1,i}^{\left[\alpha \atop{ \beta  \atop \gamma  } \right]\times k}
+ \sigma _{i-3,i-2,i-2,i-2}^{\left[\alpha \atop{ \beta   \atop \gamma } \right]\times k}
 +\sigma _{i-3,i-2,i-2,i-1}^{\left[\alpha \atop{ \beta \atop \gamma  } \right]\times k}   \right]\\
 & -8\cdot\left[\sigma _{i-3,i-3,i-3,i-3}^{\left[\alpha \atop{ \beta \atop \gamma } \right]\times k}
+ \sigma _{i-3,i-3,i-3,i-2}^{\left[\alpha \atop{ \beta \atop \gamma } \right]\times k}
+  \sigma _{i-3,i-3,i-2,i-2}^{\left[\alpha \atop{ \beta \atop \gamma } \right]\times k}
 +\sigma _{i-3,i-3,i-2,i-1}^{\left[\alpha \atop{ \beta \atop \gamma } \right]\times k}\right] \\
&-8\cdot\left[ \sigma _{i-3,i-2,i-2,i-2}^{\left[\alpha \atop{ \beta \atop \gamma } \right]\times k} 
+ \sigma _{i-3,i-2,i-2,i-1}^{\left[\alpha \atop{ \beta \atop \gamma } \right]\times k}
+\sigma _{i-3,i-2,i-1,i-1}^{\left[\alpha \atop{ \beta \atop \gamma } \right]\times k}
 +\sigma _{i-3,i-2,i-1,i}^{\left[\alpha \atop{ \beta \atop \gamma } \right]\times k}\right]\\ 
 &  - \left[\sigma _{i,i,i,i}^{\left[\alpha \atop{ \beta \atop \gamma } \right]\times k}
- 7\cdot\sigma _{i-1,i-1,i-1,i-1}^{\left[\alpha \atop{ \beta \atop \gamma } \right]\times k}
+14\cdot \sigma _{i-2,i-2,i-2,i-2}^{\left[\alpha \atop{ \beta \atop \gamma } \right]\times k}
- 8\cdot \sigma _{i-3,i-3,i-3,i-3}^{\left[\alpha \atop{ \beta \atop \gamma } \right]\times k}\right] = 0 \\
& \Longleftrightarrow \\
&\left[\sigma _{i-1,i,i,i}^{\left[\gamma \atop{ \beta \atop \alpha } \right]\times k}
+  \sigma _{i-2,i-1,i,i}^{\left[\gamma  \atop{ \beta    \atop \alpha  } \right]\times k} +\sigma _{i-1,i-1,i,i}^{\left[\gamma  \atop{ \beta  \atop \alpha } \right]\times k}\right]\\
&   - \left[ \sigma _{i-2,i-1,i-1,i-1}^{\left[\gamma  \atop{ \beta  \atop \alpha } \right]\times k}
 + \sigma _{i-2,i-2,i-1,i-1}^{\left[\gamma  \atop{ \beta   \atop \alpha } \right]\times k}
 +  \sigma _{i-3,i-2,i-1,i-1}^{\left[\gamma  \atop{ \beta   \atop \alpha  } \right]\times k}\right]\\
  &  -2\cdot\left[ \sigma _{i-2,i-1,i-1,i-1}^{\left[\gamma  \atop{ \alpha  \atop \beta } \right]\times k}
+ \sigma _{i-1,i-1,i-1,i}^{\left[\gamma \atop{ \alpha  \atop \beta  } \right]\times k}
 +\sigma _{i-2,i-1,i-1,i}^{\left[\gamma  \atop{ \alpha \atop \beta  } \right]\times k}\right]\\
 &  -4\cdot\left[ \sigma _{i-2,i-1,i-1,i-1}^{\left[\beta  \atop{ \alpha  \atop \gamma } \right]\times k}
+ \sigma _{i-1,i-1,i-1,i}^{\left[\beta \atop{ \alpha  \atop \gamma  } \right]\times k}
 +\sigma _{i-2,i-1,i-1,i}^{\left[\beta  \atop{ \alpha \atop \gamma  } \right]\times k}\right]\\
 & +  2\cdot\left[ \sigma _{i-2,i-2,i-2,i-1}^{\left[\gamma \atop{ \alpha  \atop \beta  } \right]\times k}
+ \sigma _{i-3,i-2,i-2,i-2}^{\left[\gamma  \atop{ \alpha    \atop \beta } \right]\times k}
 +\sigma _{i-3,i-2,i-2,i-1}^{\left[\gamma \atop{ \alpha  \atop \beta  } \right]\times k}\right]\\
   & + 4\cdot\left[  \sigma _{i-2,i-2,i-2,i-1}^{\left[\beta   \atop{\alpha   \atop \gamma  } \right]\times k}
+ \sigma _{i-3,i-2,i-2,i-2}^{\left[\beta \atop{ \alpha   \atop \gamma   } \right]\times k}
 +\sigma _{i-3,i-2,i-2,i-1}^{\left[\beta  \atop{ \alpha \atop \gamma  } \right]\times k} \right]\\
 & + 8\cdot\left[\sigma _{i-2,i-2,i-2,i-1}^{\left[\alpha \atop{ \beta  \atop \gamma } \right]\times k}
+  \sigma _{i-2,i-2,i-1,i-1}^{\left[\alpha \atop{ \beta   \atop \gamma  } \right]\times k}
 +\sigma _{i-2,i-2,i-1,i}^{\left[\alpha \atop{ \beta  \atop \gamma  } \right]\times k}\right]\\
  & -8\cdot\left[ \sigma _{i-3,i-3,i-3,i-2}^{\left[\alpha \atop{ \beta \atop \gamma } \right]\times k}
+  \sigma _{i-3,i-3,i-2,i-2}^{\left[\alpha \atop{ \beta \atop \gamma } \right]\times k}
 +\sigma _{i-3,i-3,i-2,i-1}^{\left[\alpha \atop{ \beta \atop \gamma } \right]\times k}\right]  = 0 
\end{align*}
Combining \eqref{eq 2.66}, \eqref{eq 2.74} and  \eqref{eq 2.82} we get \\
\begin{align*}
& \eqref{eq 2.83} \Longleftrightarrow 
\left[ \sigma _{i-1,i,i,i}^{\left[\gamma \atop{ \beta \atop \alpha } \right]\times k}
+   2\cdot\big[  \sigma _{i-1,i-1,i-1,i}^{\left[ \gamma  \atop{  \alpha  \atop  \beta  } \right]\times k}
+  \sigma _{i-2,i-1,i-1,i}^{\left[ \gamma  \atop{  \alpha \atop  \beta  } \right]\times k}  \big] \right] \\                                                                                                                                                                                                                        \\
&   - \left[\sigma _{i-2,i-1,i-1,i-1}^{\left[\gamma  \atop{ \beta  \atop \alpha } \right]\times k}
+  2\cdot\big[  \sigma _{i-2,i-2,i-2,i-1}^{\left[ \gamma  \atop{  \alpha  \atop  \beta  } \right]\times k}
+ \sigma _{i-3,i-2,i-2,i-1}^{\left[ \gamma  \atop{  \alpha \atop  \beta  } \right]\times k}  \big]\right] \\               
  &  -2\cdot\left[\sigma _{i-2,i-1,i-1,i-1}^{\left[\gamma  \atop{ \alpha  \atop \beta } \right]\times k}
+ \sigma _{i-1,i-1,i-1,i}^{\left[\gamma \atop{ \alpha  \atop \beta  } \right]\times k}
 +\sigma _{i-2,i-1,i-1,i}^{\left[\gamma  \atop{ \alpha \atop \beta  } \right]\times k}\right]\\
  &  -4\cdot\left[ \sigma _{i-1,i-1,i-1,i}^{\left[\beta \atop{ \alpha  \atop \gamma  } \right]\times k}
 + 2\cdot\big[  \sigma _{i-2,i-2,i-1,i-1}^{\left[ \alpha  \atop{  \beta  \atop  \gamma  } \right]\times k}
+  \sigma _{i-2,i-2,i-1,i}^{\left[ \alpha  \atop{  \beta  \atop  \gamma  } \right]\times k}  \big]   \right]\\
 & +  2\cdot\left[ \sigma _{i-2,i-2,i-2,i-1}^{\left[\gamma \atop{ \alpha  \atop \beta  } \right]\times k}
+ \sigma _{i-3,i-2,i-2,i-2}^{\left[\gamma  \atop{ \alpha    \atop \beta } \right]\times k}
 +\sigma _{i-3,i-2,i-2,i-1}^{\left[\gamma \atop{ \alpha  \atop \beta  } \right]\times k}\right]\\
   & + 4\cdot\left[ \sigma _{i-2,i-2,i-2,i-1}^{\left[\beta   \atop{\alpha   \atop \gamma  } \right]\times k}
+  2\cdot\big[  \sigma _{i-3,i-3,i-2,i-2}^{\left[ \alpha  \atop{  \beta  \atop  \gamma  } \right]\times k}
+  \sigma _{i-3,i-3,i-2,i-1}^{\left[ \alpha  \atop{  \beta  \atop  \gamma  } \right]\times k}  \big]   \right]\\
 & + 8\cdot\left[  \sigma _{i-2,i-2,i-2,i-1}^{\left[\alpha \atop{ \beta  \atop \gamma } \right]\times k}
+  \sigma _{i-2,i-2,i-1,i-1}^{\left[\alpha \atop{ \beta   \atop \gamma  } \right]\times k}
 +\sigma _{i-2,i-2,i-1,i}^{\left[\alpha \atop{ \beta  \atop \gamma  } \right]\times k}\right]\\
  & -8\cdot\left[ \sigma _{i-3,i-3,i-3,i-2}^{\left[\alpha \atop{ \beta \atop \gamma } \right]\times k}
+  \sigma _{i-3,i-3,i-2,i-2}^{\left[\alpha \atop{ \beta \atop \gamma } \right]\times k}
 +\sigma _{i-3,i-3,i-2,i-1}^{\left[\alpha \atop{ \beta \atop \gamma } \right]\times k}\right]  = 0  \\
  & \Longleftrightarrow \\
& \big( \sigma _{i-1,i,i,i}^{\left[\gamma \atop{ \beta \atop \alpha } \right]\times k}
-4\cdot \sigma _{i-1,i-1,i-1,i}^{\left[\beta \atop{ \alpha  \atop \gamma  } \right]\times k}\big)
-\big( \sigma _{i-2,i-1,i-1,i-1}^{\left[\gamma  \atop{ \beta  \atop \alpha } \right]\times k} -
4\cdot\sigma _{i-2,i-2,i-2,i-1}^{\left[\beta   \atop{\alpha   \atop \gamma  } \right]\times k} \big) \\
& -2\cdot\big(  \sigma _{i-2,i-1,i-1,i-1}^{\left[\gamma  \atop{ \alpha  \atop \beta } \right]\times k}
 - 4\cdot\sigma _{i-2,i-2,i-2,i-1}^{\left[\alpha \atop{ \beta  \atop \gamma } \right]\times k}  \big)
+2\cdot\big( \sigma _{i-3,i-2,i-2,i-2}^{\left[\gamma  \atop{ \alpha    \atop \beta } \right]\times k}
 -  4\cdot\sigma _{i-3,i-3,i-3,i-2}^{\left[\alpha \atop{ \beta \atop \gamma } \right]\times k} \big) = 0
 \end{align*}
 \end{proof}
 \section{\textbf{Rank  properties of a  partition  of  triple persymmetric  matrices }  }
 \label{sec 3}
  \subsection{Notation}
  \label{subsec 1}
    \begin{defn}
\label{defn 3.1}We introduce the following definitions in the three - dimensional  $ \mathbb{K} $- vectorspace.\\
\begin{itemize}
\item Let k,s,m and $l$ denote rational integers such that $ k\geq 1,\; s\geq 2 \; and \; m\geq 0, \; l\geq 0 $
\item We denote by $\mathbb{P}/\mathbb{P}_{i}\times \mathbb{P}/\mathbb{P}_{j}\times \mathbb{P}/\mathbb{P}_{r} $
a complete set of coset representatives of $\mathbb{P}_{i}\times\mathbb{P}_{j}\times\mathbb{P}_{r} $
in  $\mathbb{P}\times\mathbb{P}\times\mathbb{P},$ for instance $\mathbb{P}/\mathbb{P}_{s+k-1}\times \mathbb{P}/\mathbb{P}_{s+m +k-1}\times \mathbb{P}/\mathbb{P}_{s+m+l +k-1} $
denotes  a complete set of coset representatives of $\mathbb{P}_{s+k-1}\times\mathbb{P}_{s+m +k-1}\times\mathbb{P}_{s+m+l +k-1} $ in  $\mathbb{P}\times\mathbb{P}\times\mathbb{P}.$
\item $Set \;(t,\eta,\xi  )= (\sum_{i\geq 1}\alpha _{i}T^{-i},\sum_{i\geq 1}\beta  _{i}T^{-i},\sum_{i\geq 1}\gamma  _{i}T^{-i}) \in \mathbb{P}^3$
\item  We denote by 
 $  D^{ \left[s-1\atop{ s+m-1\atop {s+m+l-1 \atop{\overline{\alpha _{s -} \atop{\beta _{s+m -} \atop \gamma_{s+m+l -}} }}}}\right]\times k}(t,\eta ,\xi )  $
  any  $(3s+2m+l)\times k $   matrix, 
such that  after a rearrangement of the rows, if necessary,  we can  obtain a matrix where the first $(3s+2m+l-3)$  rows form  the following triple persymmetric matrix 
$ \left[{D_{s -1 \times k}(t)\over {D_{(s+m-1 )\times k}(\eta )\over D_{(s+m+l-1)\times k}(\xi )}}\right] $ and the last three rows are equal to\\
the following $3\times k$ matrix 
 $\begin{pmatrix}
\alpha  _{s } & \alpha  _{s +1} & \ldots &\alpha  _{s +k-1} \\
\beta _{s+m} & \beta  _{s+m+1} & \ldots & \beta _{s+m+k-1}\\
\gamma  _{s+m+l} & \gamma  _{s+m+l+1}  & \ldots  & \gamma  _{s+m+l+k-1} 
\end{pmatrix},$ that is we get a matrix of the form\\
 $$   \left ( \begin{array} {cccccc}
\alpha _{1} & \alpha _{2}  &  \ldots & \alpha _{k-1}  &  \alpha _{k} \\
\alpha _{2 } & \alpha _{3} &  \ldots  &  \alpha _{k} &  \alpha _{k+1} \\
\vdots & \vdots & \vdots    & \vdots  &  \vdots \\
\alpha _{s-1} & \alpha _{s} & \ldots  &  \alpha _{s+k-3} &  \alpha _{s+k-2}  \\
\hline \\
\beta  _{1} & \beta  _{2}  & \ldots  &  \beta_{k-1} &  \beta _{k}  \\
\beta  _{2} & \beta  _{3}  & \ldots  &  \beta_{k} &  \beta _{k+1}  \\
\vdots & \vdots    &  \vdots & \vdots  &  \vdots \\
\beta  _{m+1} & \beta  _{m+2}  & \ldots  &  \beta_{k+m-1} &  \beta _{k+m}  \\
\vdots & \vdots    &  \vdots & \vdots  &  \vdots \\
\beta  _{s+m-1} & \beta  _{s+m}  & \ldots  &  \beta_{s+m+k-3} &  \beta _{s+m+k-2}  \\
\hline \\
\gamma  _{1} & \gamma   _{2}  & \ldots  & \gamma  _{k-1} &  \gamma  _{k}  \\
\gamma  _{2} & \gamma  _{3}  & \ldots  & \gamma  _{k} &  \gamma  _{k+1}  \\
\vdots & \vdots    &  \vdots & \vdots  &  \vdots \\
 \gamma  _{m+l+1} &  \gamma _{m+l+2}  & \ldots  & \gamma _{k+m+l-1} &  \gamma  _{k+m+l}  \\
\vdots & \vdots   &  \vdots & \vdots  &  \vdots \\
 \gamma  _{s+m+l-1} & \gamma  _{s+m+l}  & \ldots  & \gamma  _{s+m+l+k-3} &  \gamma  _{s+m+l+k-2}  \\
   \hline\\
  \alpha  _{s } & \alpha  _{s +1} & \ldots & \alpha  _{s +k-2}& \alpha  _{s +k-1}\\
   \beta _{s+m} & \beta _{s+m+1} & \ldots & \beta _{s+m+k-2} & \beta _{s+m+k-1}\\
  \gamma  _{s+m+l} & \gamma  _{s+m+l+1}  & \ldots  & \gamma  _{s+m+l+k-2} &  \gamma  _{s+m+l+k-1}  
\end{array}  \right). $$ 
\item
   Let j be a rational integer such that $1\leq j\leq k-1.$\\
  
   $Set \;(t,\eta,\xi  )= (\sum_{i\geq 1}\alpha _{i}T^{-i},\sum_{i\geq 1}\beta  _{i}T^{-i},\sum_{i\geq 1}\gamma  _{i}T^{-i}) \in \mathbb{P}^3$\\
   
   We denote by 
 $  D_{j}^{ \left[s\atop{ s+m\atop s+m+l-1}\right]\times (k-j+1)}(t,\eta ,\xi )  $
  any  $(3s+2m+l)\times (k-j+1) $   matrix, 
such that  after a rearrangement of the rows, if necessary,  we can  obtain  the following triple persymmetric matrix 
$ \left[{D_{s \times (k-j+1)}^{j}(t)\over {D_{(s+m )\times (k-j+1)}^{j}(\eta )\over D_{(s+m+l)\times (k-j+1)}^{j}(\xi )}}\right] $

   $$   \left ( \begin{array} {ccccc|c}
\alpha _{j} & \alpha _{j+1} & \alpha _{j+2} &  \ldots & \alpha _{k-1}  &  \alpha _{k} \\
\alpha _{j+1 } & \alpha _{j+2} & \alpha _{j+3}&  \ldots  &  \alpha _{k} &  \alpha _{k+1} \\
\vdots & \vdots & \vdots    &  \vdots & \vdots  &  \vdots \\
\alpha _{j+s-1} & \alpha _{j+s} & \alpha _{j+s +1} & \ldots  &  \alpha _{k+s-2} &  \alpha _{k+s-1}  \\
\hline 
\beta  _{j} & \beta  _{j+1} & \beta  _{j+2} & \ldots  &  \beta_{k-1} &  \beta _{k}  \\
\beta  _{j+1} & \beta  _{j+2} & \beta  _{j+3} & \ldots  &  \beta_{k} &  \beta _{k+1}  \\
\vdots & \vdots & \vdots    &  \vdots & \vdots  &  \vdots \\
\beta  _{m+j} & \beta  _{m+j+1} & \beta  _{m+j+2} & \ldots  &  \beta_{k+m-1} &  \beta _{k+m}  \\
\vdots & \vdots & \vdots    &  \vdots & \vdots  &  \vdots \\
\beta  _{s+m+j-1} & \beta  _{s+m+j} & \beta  _{s+m+j+1} & \ldots  &  \beta_{s+m+k-2} &  \beta _{s+m+k-1}  \\
\hline
\gamma  _{j} & \gamma   _{j+1} &  \gamma _{j+2} & \ldots  & \gamma  _{k-1} &  \gamma  _{k}  \\
\gamma  _{j+1} & \gamma  _{j+2} & \gamma  _{j+3} & \ldots  & \gamma  _{k} &  \gamma  _{k+1}  \\
\vdots & \vdots & \vdots    &  \vdots & \vdots  &  \vdots \\
 \gamma  _{m+l+j} &  \gamma _{m+l+j+1} &  \gamma _{m+l+j+2} & \ldots  & \gamma _{k+m+l-1} &  \gamma  _{k+m+l}  \\
\vdots & \vdots & \vdots    &  \vdots & \vdots  &  \vdots \\
\gamma  _{s+m+l+j-1} & \gamma  _{s+m+l+j} &  \gamma  _{s+m+l+j+1} & \ldots  & \gamma  _{s+m+l+k-2} &  \gamma  _{s+m+l+k-1}  
 \end{array}  \right). $$ \\
 \item To simplify the notations concerning the exponential sums used in the proofs, we introduce the following definitions.\\
 \item  Let $  \psi _{1}(t,\eta,\xi  ) $ be the quadratic  exponential sum in $\mathbb{P}\times\mathbb{P}\times\mathbb{P}$ defined by
$$ (t,\eta,\xi  ) \in  \mathbb{P}\times \mathbb{P}\times\mathbb{P}\longmapsto  
  \sum_{deg Y = k-1}\sum_{deg Z = s-1}E(tYZ)\sum_{deg U \leq s+m-2}E(\eta YU)
   \sum_{deg V \leq s+m+l-2}E(\xi  YV) \in \mathbb{Z}  $$\\
 \item  Let $  \psi _{2}(t,\eta,\xi  ) $ be the quadratic  exponential sum in $\mathbb{P}\times\mathbb{P}\times\mathbb{P}$ defined by
$$ (t,\eta,\xi  ) \in  \mathbb{P}\times \mathbb{P}\times\mathbb{P}\longmapsto  
  \sum_{deg Y \leq k-2}\sum_{deg Z \leq s-1}E(tYZ)\sum_{deg U = s+m-1}E(\eta YU)
   \sum_{deg V \leq s+m+l-2}E(\xi  YV) \in \mathbb{Z}  $$\\  
  \item  Let $  \psi _{3}(t,\eta,\xi  ) $ be the quadratic  exponential sum in $\mathbb{P}\times\mathbb{P}\times\mathbb{P}$ defined by
$$ (t,\eta,\xi  ) \in  \mathbb{P}\times \mathbb{P}\times\mathbb{P}\longmapsto  
  \sum_{deg Y \leq k-2}\sum_{deg Z \leq s-1}E(tYZ)\sum_{deg U \leq s+m-1}E(\eta YU)
   \sum_{deg V  =  s+m+l-1}E(\xi  YV) \in \mathbb{Z}  $$\\
   \item  Let $ \chi _{1}(t,\eta,\xi  ) $ be the quadratic  exponential sum in $\mathbb{P}\times\mathbb{P}\times\mathbb{P}$ defined by
$$ (t,\eta,\xi  ) \in  \mathbb{P}\times \mathbb{P}\times\mathbb{P}\longmapsto  
  \sum_{deg Y = k-1}\sum_{deg Z \leq s-2}E(tYZ)\sum_{deg U \leq s+m-2}E(\eta YU)
   \sum_{deg V  \leq  s+m+l-2}E(\xi  YV) \in \mathbb{Z}  $$\\
    \item  Let $ \chi _{2}(t,\eta,\xi  ) $ be the quadratic  exponential sum in $\mathbb{P}\times\mathbb{P}\times\mathbb{P}$ defined by
$$ (t,\eta,\xi  ) \in  \mathbb{P}\times \mathbb{P}\times\mathbb{P}\longmapsto  
  \sum_{deg Y\leq  k-2}\sum_{deg Z = s-1}E(tYZ)\sum_{deg U \leq s+m-2}E(\eta YU)
   \sum_{deg V  \leq  s+m+l-2}E(\xi  YV) \in \mathbb{Z}  $$\\ 
     \item We recall the definition of $h_{1}(t,\eta,\xi  )$ (see Definition \ref{defn 2.1}) \\
      Let $ h _{1}(t,\eta,\xi  ) $ be the quadratic  exponential sum in $\mathbb{P}\times\mathbb{P}\times\mathbb{P}$ defined by
$$ (t,\eta,\xi  ) \in  \mathbb{P}\times \mathbb{P}\times\mathbb{P}\longmapsto  
  \sum_{deg Y\leq  k-1}\sum_{deg Z = s-1}E(tYZ)\sum_{deg U \leq s+m-2}E(\eta YU)
   \sum_{deg V  \leq  s+m+l-2}E(\xi  YV) \in \mathbb{Z}  $$\\ 
 \item       Let $ \theta  _{1}(t,\eta,\xi  ) $ be the quadratic  exponential sum in $\mathbb{P}\times\mathbb{P}\times\mathbb{P}$ defined by
$$ (t,\eta,\xi  ) \in  \mathbb{P}\times \mathbb{P}\times\mathbb{P}\longmapsto  
  \sum_{deg Y = k-1}\sum_{deg Z \leq  s-1}E(tYZ)\sum_{deg U = s+m-1}E(\eta YU)
   \sum_{deg V  = s+m+l-2}E(\xi  YV) \in \mathbb{Z}  $$\\  
    \item       Let $ \theta  _{2}(t,\eta,\xi  ) $ be the quadratic  exponential sum in $\mathbb{P}\times\mathbb{P}\times\mathbb{P}$ defined by
$$ (t,\eta,\xi  ) \in  \mathbb{P}\times \mathbb{P}\times\mathbb{P}\longmapsto  
  \sum_{deg Y = k-1}\sum_{deg Z \leq  s-1}E(tYZ)\sum_{deg U \leq  s+m-1}E(\eta YU)
   \sum_{deg V  = s+m+l-1}E(\xi  YV) \in \mathbb{Z}  $$\\  
 \item    Let  $ ( j_{1}, j_{2},   j_{3}, j_{4},j_{5}, j_{6}, j_{7}, j_{8}) \in \mathbb{N}^{8}, $ we define \\
\begin{align*}
  {}^{\#}\left(\begin{array}{c | c}
           j_{1}  &  j_{2} \\
           \hline
           j_{3}  &  j_{4} \\
           \hline
            j_{5}  &  j_{6} \\
            \hline 
           j_{7}   &  j_{8}
           \end{array} \right)_{\mathbb{P}/\mathbb{P}_{k+s -1}\times
           \mathbb{P}/\mathbb{P}_{k+s+m-1}\times  \mathbb{P}/\mathbb{P}_{k+s+m+l-1} }^{{\alpha \over {\beta \over \gamma }}} 
   \end{align*}
    to be the cardinality of the following set
    $$\begin{array}{l}\Big\{ (t,\eta ,\xi ) \in \mathbb{P}/\mathbb{P}_{k+s -1}\times
           \mathbb{P}/\mathbb{P}_{k+s+m-1} \times \mathbb{P}/\mathbb{P}_{k+s+m+l -1} 
\mid r( D^{ \left[s-1\atop{ s+m-1\atop {s+m+l-1}}\right]\times (k-1)}(t,\eta ,\xi  ) )  = j_{1}, \\
r( D^{ \left[s-1\atop{ s+m-1\atop {s+m+l-1}}\right]\times k}(t,\eta ,\xi  ) ) = j_{2},  
 r(  D^{ \left[s\atop{ s+m-1\atop {s+m+l-1}}\right]\times (k-1)}(t,\eta ,\xi  ) ) = j_{3},\\
   r(  D^{ \left[s\atop{ s+m-1\atop {s+m+l-1}}\right]\times k}(t,\eta ,\xi  ) ) = j_{4},\quad
   r(  D^{ \left[s\atop{ s+m\atop {s+m+l-1}}\right]\times (k-1)}(t,\eta ,\xi  ) ) = j_{5},\\
     r(  D^{ \left[s\atop{ s+m\atop {s+m+l-1}}\right]\times k}(t,\eta ,\xi  ) ) = j_{6},\quad
    r(  D^{ \left[s\atop{ s+m\atop {s+m+l}}\right]\times (k-1)}(t,\eta ,\xi  ) ) = j_{7},\quad
     r(  D^{ \left[s\atop{ s+m\atop {s+m+l}}\right]\times k}(t,\eta ,\xi  ) ) = j_{8}. \Big\}  
   \end{array}$$\\  
\item We define    
     \begin{equation*}
\begin{vmatrix}
k-1 & \vline & k  & \vline \\
\hline
\cdot &\vline & \cdot & \vline & D^{\left[s-1\atop{ s+m-1\atop s+m+l-1} \right]\times \cdot}  \\
\hline
j &\vline & j &  \vline  & \alpha _{s}-\\
\hline 
j &\vline & j+1 & \vline  & \beta _{s+m}- \\
\hline 
\cdot &\vline & \cdot & \vline  & \gamma _{s+m+l}-
\end{vmatrix}
\end{equation*}
    to be  the following subset of $\mathbb{P}^{3}$
    $$\begin{array}{l}\Big\{ (t,\eta ,\xi ) = (\sum_{i\geq 1}\alpha _{i}T^{-i},\sum_{i\geq 1}\beta  _{i}T^{-i},\sum_{i\geq 1}\gamma  _{i}T^{-i}) \in \mathbb{P}^3
   \mid r( D^{ \left[s\atop{ s+m-1\atop {s+m+l-1}}\right]\times (k-1)}(t,\eta ,\xi  ) )  = j, \\
r( D^{ \left[s\atop{ s+m-1\atop {s+m+l-1}}\right]\times k}(t,\eta ,\xi  ) )  = j,\quad
r( D^{ \left[s\atop{ s+m\atop {s+m+l-1}}\right]\times (k-1)}(t,\eta ,\xi  ) ) = j, \quad 
 r(  D^{ \left[s\atop{ s+m\atop {s+m+l-1}}\right]\times k}(t,\eta ,\xi  ) ) = j+1 \Big\}  
   \end{array}$$\\  
\item Similar expressions are defined in a similar way.   
    \end{itemize}
\end{defn}

    \subsection{Introduction}
  \label{subsec 2}
We study rank properties of a partition of triple persymmetric matrices by integrating some appropriate exponential sums on 
the unit interval of $\textbf{K}^3$. We adapt the method used in Section 6 and Section 7 of [2] \\[0.01cm]

     \subsection{Computation of exponential sums in $\mathbb{K}^3$ }
  \label{subsec 3} 
    
  \begin{lem}
\label{lem 3.2}
Let $ (t,\eta ,\xi ) \in  \mathbb{P}\times \mathbb{P}\times \mathbb{P}  $ and
\begin{align*}
\psi _{1}(t,\eta,\xi ) & = \sum_{deg Y= k-1}\sum_{deg Z = s-1}E(tYZ)\sum_{deg U \leq s+m-2}E(\eta YU) \sum_{deg V \leq s+m+l-2}E(\xi YV)  \\
 \psi _{2}(t,\eta,\xi) & = \sum_{deg Y\leq k-2} \sum_{deg Z\leq  s-1}E(tYZ)\sum_{deg U = s+m-1}E(\eta YU) \sum_{deg V \leq s+m+l-2}E(\xi YV) \\
  \psi _{3}(t,\eta,\xi) & = \sum_{deg Y\leq k-2} \sum_{deg Z\leq  s-1}E(tYZ)\sum_{deg U \leq s+m-1}E(\eta YU) \sum_{deg V = s+m+l-1}E(\xi YV)   \\
   \chi_{1}(t,\eta,\xi) & = \sum_{deg Y = k-1} \sum_{deg Z\leq  s-2}E(tYZ)\sum_{deg U \leq s+m-2}E(\eta YU) \sum_{deg V \leq s+m+l-2}E(\xi YV) \\
  \chi_{2}(t,\eta,\xi) & = \sum_{deg Y\leq k-2}\sum_{deg Z= s-1}E(tYZ)\sum_{deg U \leq s+m-2}E(\eta YU) \sum_{deg V \leq s+m+l-2}E(\xi  YV ) \\
 h_{1}(t,\eta,\xi) & = \sum_{deg Y\leq k-1}\sum_{deg Z =  s-1}E(tYZ)\sum_{deg U \leq s+m-2}E(\eta YU) \sum_{deg V \leq s+m+l-2}E(\xi  YV ) \\
 \end{align*}
Then \\
  \begin{align}
  \label{eq 3.1}
& \psi _{1}(t,\eta,\xi  ) = 2^{3s+2m+l-3}\cdot
 \sum_{deg Y= k-1\atop {Y\in ker D^{\left[\stackrel{s-1}{\stackrel{s+m-1}{s+m+l-1}}\right] \times k }(t,\eta ,\xi )}}E(tYT^{s-1}) 
 \end{align}
 
 \begin{align}
   \label{eq 3.2}
& \psi _{2}(t,\eta,\xi  ) = 2^{3s+2m+l-2}\cdot
 \sum_{deg Y\leq k-2\atop {Y\in ker D^{\left[\stackrel{s}{\stackrel{s+m-1}{s+m+l-1}}\right] \times (k-1) }(t,\eta ,\xi )}}E(\eta YT^{s+m-1}) \\
& =  \begin{cases}
 2^{k+3s+2m+l-3-  r( D^{\left[\stackrel{s}{\stackrel{s+m-1}{s+m+l-1}}\right] \times k }(t,\eta ,\xi ) )  }  & \text{if }
  r( D^{\left[\stackrel{s}{\stackrel{s+m-1}{s+m+l-1}}\right] \times (k-1) }(t,\eta ,\xi ) ) =  r( D^{\left[\stackrel{s}{\stackrel{s+m}{s+m+l-1}}\right] \times (k-1) }(t,\eta ,\xi ) ) \nonumber \\
     0  & \text{otherwise }.
    \end{cases}
 \end{align}
 
  \begin{align}
    \label{eq 3.3}
& \psi _{3}(t,\eta,\xi  ) = 2^{3s+2m+l-1}\cdot
 \sum_{deg Y\leq k-2\atop {Y\in ker D^{\left[\stackrel{s}{\stackrel{s+m}{s+m+l-1}}\right] \times (k-1) }(t,\eta ,\xi )}}E(\xi  YT^{s+m+l-1}) \\
& =  \begin{cases}
 2^{k+3s+2m+l-2-  r( D^{\left[\stackrel{s}{\stackrel{s+m}{s+m+l-1}}\right] \times (k-1) }(t,\eta ,\xi ) )  }  & \text{if }
  r( D^{\left[\stackrel{s}{\stackrel{s+m}{s+m+l-1}}\right] \times (k-1) }(t,\eta ,\xi ) ) =  r( D^{\left[\stackrel{s}{\stackrel{s+m}{s+m+l}}\right] \times (k-1) }(t,\eta ,\xi ) )  \nonumber \\
     0  & \text{otherwise }.
    \end{cases}
 \end{align}
 
  \begin{align}
    \label{eq 3.4}
& \chi_{1}(t,\eta,\xi  ) 
 =  \begin{cases}
 2^{k+3s+2m+l-4-  r( D^{\left[\stackrel{s-1}{\stackrel{s+m-1}{s+m+l-1}}\right] \times (k-1) }(t,\eta ,\xi ) )  }  & \text{if }
  r( D^{\left[\stackrel{s-1}{\stackrel{s+m-1}{s+m+l-1}}\right] \times (k-1) }(t,\eta ,\xi ) ) =   r( D^{\left[\stackrel{s-1}{\stackrel{s+m-1}{s+m+l-1}}\right] \times k }(t,\eta ,\xi ) )  \\
   0  & \text{otherwise }.
    \end{cases}
 \end{align}
 
  \begin{align}
    \label{eq 3.5}
& \chi_{2}(t,\eta,\xi  ) 
 =  \begin{cases}
 2^{k+3s+2m+l-4-  r( D^{\left[\stackrel{s-1}{\stackrel{s+m-1}{s+m+l-1}}\right] \times (k-1) }(t,\eta ,\xi ) )  }  & \text{if }
  r( D^{\left[\stackrel{s-1}{\stackrel{s+m-1}{s+m+l-1}}\right] \times (k-1) }(t,\eta ,\xi ) ) =   r( D^{\left[\stackrel{s}{\stackrel{s+m-1}{s+m+l-1}}\right] \times (k-1) }(t,\eta ,\xi ) )  \\
   0  & \text{otherwise }.
    \end{cases}
 \end{align}
   \begin{align}
    \label{eq 3.6}
  & \text{We recall that see \eqref{eq 2.2}}\nonumber  \\
    & h_{1}(t,\eta,\xi  ) = 2^{3s+2m+l-3}\cdot
 \sum_{deg Y\leq k-1\atop {Y\in ker D^{\left[\stackrel{s-1}{\stackrel{s+m-1}{s+m+l-1}}\right] \times k }(t,\eta ,\xi )}}E(tYT^{s-1})  \\
& =  \begin{cases}
 2^{k+3s+2m+l-3-  r( D^{\left[\stackrel{s-1}{\stackrel{s+m-1}{s+m+l-1}}\right] \times k }(t,\eta ,\xi ) )  }  & \text{if }
  r( D^{\left[\stackrel{s-1}{\stackrel{s+m-1}{s+m+l-1}}\right] \times k }(t,\eta ,\xi ) ) =   r( D^{\left[\stackrel{s}{\stackrel{s+m-1}{s+m+l-1}}\right] \times k }(t,\eta ,\xi ) )  \nonumber  \\
     0  & \text{otherwise }.
    \end{cases}
 \end{align}
  \end{lem} 
\begin{proof}  
 The proof of \ref{lem 3.2} is somewhat similar to the proofs of the results obtained in [2, see section 4]
  \end{proof}
     \begin{lem}
\label{lem 3.3}
Let $ (t,\eta ,\xi ) \in  \mathbb{P}\times \mathbb{P}\times \mathbb{P}  $  then we have 
\begin{align}
\displaystyle
 &  \psi _{1}^{2}(t,\eta,\xi ) = \chi_{1}(t,\eta,\xi )\cdot \chi_{2}(t,\eta,\xi)  \label{eq 3.7} \\
  & \text{where} \nonumber \\
&   \psi _{1}(t,\eta,\xi )  = \sum_{deg Y= k-1}\sum_{deg Z = s-1}E(tYZ)\sum_{deg U \leq s+m-2}E(\eta YU) \sum_{deg V \leq s+m+l-2}E(\xi YV) \nonumber  \\
&  \chi_{1}(t,\eta,\xi) = \sum_{deg Y = k-1} \sum_{deg Z\leq  s-2}E(tYZ)\sum_{deg U \leq s+m-2}E(\eta YU) \sum_{deg V \leq s+m+l-2}E(\xi YV) \nonumber \\
&    \chi_{2}(t,\eta,\xi)  = \sum_{deg Y\leq k-2}\sum_{deg Z= s-1}E(tYZ)\sum_{deg U \leq s+m-2}E(\eta YU) \sum_{deg V \leq s+m+l-2}E(\xi  YV ) \nonumber 
 \end{align} 
\end{lem} 
 
 \newpage
\begin{proof} 
We have \\

\begin{align*}
 &  \psi_{1}^{2}(t,\eta,\xi  ) 
  =  \big[2^{3s+2m+l-3}\cdot
 \sum_{deg Y_{1}= k-1\atop {Y_{1}\in ker D^{\left[\stackrel{s-1}{\stackrel{s+m-1}{s+m+l-1}}\right] \times k }(t,\eta ,\xi )}}E(tY_{1}T^{s-1}) \big]
\big[2^{3s+2m+l-3}\cdot
 \sum_{deg Y_{2}= k-1\atop {Y_{2}\in ker D^{\left[\stackrel{s-1}{\stackrel{s+m-1}{s+m+l-1}}\right] \times k }(t,\eta ,\xi )}}E(tY_{2}T^{s-1}) \big] \\
&   \text{ We set } \\
  &   \left\{\begin{array}{ccc}
Y_{1} + Y_{2} = Y_{3},  &  deg Y_{1} = k-1, & Y_{1}\in ker D^{\left[\stackrel{s-1}{\stackrel{s+m-1}{s+m+l-1}}\right] \times k }(t,\eta ,\xi )   \\
              Y_{2} = Y_{4},   &   deg Y_{2} = k-1, & Y_{2}\in ker D^{\left[\stackrel{s-1}{\stackrel{s+m-1}{s+m+l-1}}\right] \times k }(t,\eta ,\xi )
\end{array}\right.\\  
& \Longleftrightarrow \\
  &   \left\{\begin{array}{cc}
    deg Y_{3} \leq k-2, & Y_{3}\in ker D^{\left[\stackrel{s-1}{\stackrel{s+m-1}{s+m+l-1}}\right] \times (k-1) }(t,\eta ,\xi )   \\
    deg Y_{4} = k-1, & Y_{4}\in ker D^{\left[\stackrel{s-1}{\stackrel{s+m-1}{s+m+l-1}}\right] \times k }(t,\eta ,\xi )
\end{array}\right.\\ 
& \\
 & \text{Then we obtain}\\
 & \\
  &  \psi_{1}^{2}(t,\eta,\xi  )   = \big[2^{3s+2m+l-3}\cdot
 \sum_{deg Y_{4}= k-1\atop {Y_{4}\in ker D^{\left[\stackrel{s-1}{\stackrel{s+m-1}{s+m+l-1}}\right] \times k }(t,\eta ,\xi )}}1 \big] 
 \big[2^{3s+2m+l-3}\cdot
 \sum_{deg Y_{3}\leq k-2\atop {Y_{3}\in ker D^{\left[\stackrel{s-1}{\stackrel{s+m-1}{s+m+l-1}}\right] \times k }(t,\eta ,\xi )}}E(tY_{3}T^{s-1}) \big] \\
   & = \chi_{1}(t,\eta,\xi)\cdot \chi_{2}(t,\eta,\xi)
 \end{align*}

\end{proof}  
     \begin{lem}
\label{lem 3.4}
We have the following equivalences :\\
\begin{align}
&  \psi_{1}(t,\eta,\xi  ) \neq 0  \Longleftrightarrow \psi_{1}^{2}(t,\eta,\xi  ) \neq 0
\Longleftrightarrow  \chi_{1}(t,\eta,\xi)\cdot \chi_{2}(t,\eta,\xi) \neq 0 \label{eq 3.8}\\
& \nonumber\\
& \Longleftrightarrow  \chi_{1}(t,\eta,\xi) \neq 0 \quad\text{and}\quad \chi_{1}(t,\eta,\xi) \neq 0 \nonumber\\
&\nonumber \\
& \Longleftrightarrow 
  r( D^{\left[\stackrel{s-1}{\stackrel{s+m-1}{s+m+l-1}}\right] \times (k-1) }(t,\eta ,\xi ) ) =
  r( D^{\left[\stackrel{s-1}{\stackrel{s+m-1}{s+m+l-1}}\right] \times k }(t,\eta ,\xi ) ) =
  r( D^{\left[\stackrel{s}{\stackrel{s+m-1}{s+m+l-1}}\right] \times (k-1) }(t,\eta ,\xi ) ) \nonumber
\end{align}
 \end{lem}
 \begin{proof}
 \eqref{eq 3.8} follows from \eqref{eq 3.4}   and   \eqref{eq 3.5}.
 
\end{proof}
  \begin{lem}
\label{lem 3.5}
We have \\
\begin{equation}
 \psi _{1}(t,\eta ,\xi ) = 
\begin{cases}
2^{k+3s+2m+l-4-j}\quad &\text{if }\quad (t,\eta ,\xi ) \in \begin{vmatrix}
k-1 & \vline & k  & \vline \\
\hline
j &\vline & j & \vline & D^{\left[s-1\atop{ s+m-1\atop s+m+l-1} \right]\times \cdot} \label{eq 3.9} \\
\hline
j &\vline & j &  \vline  & \alpha _{s}-\\
\hline 
\cdot &\vline & \cdot & \vline  & \beta _{s+m}- \\
\hline 
\cdot &\vline & \cdot & \vline  & \gamma _{s+m+l}-
\end{vmatrix} \\
& \\
-2^{k+3s+2m+l-4-j} \quad &\text{if }\quad (t,\eta ,\xi ) \in \begin{vmatrix}
k-1 & \vline & k  & \vline \\
\hline
j &\vline & j & \vline & D^{\left[s-1\atop{ s+m-1\atop s+m+l-1} \right]\times \cdot}  \\
\hline
j &\vline & j +1&  \vline  & \alpha _{s}-\\
\hline 
\cdot &\vline & \cdot & \vline  & \beta _{s+m}- \\
\hline 
\cdot &\vline & \cdot & \vline  & \gamma _{s+m+l}-
\end{vmatrix} \\
& \\
0 &\text{otherwise}
\end{cases}
\end{equation}

\begin{equation}
 \psi _{2}(t,\eta ,\xi ) = 
\begin{cases}
2^{k+3s+2m+l-3-j}\quad &\text{if }\quad (t,\eta ,\xi ) \in\begin{vmatrix}
k-1 & \vline & k  & \vline \\
\hline
\cdot &\vline & \cdot & \vline & D^{\left[s-1\atop{ s+m-1\atop s+m+l-1} \right]\times k} \label{eq 3.10} \\
\hline
j &\vline & \cdot &  \vline  & \alpha _{s}-\\
\hline 
j &\vline & \cdot & \vline  & \beta _{s+m}- \\
\hline 
\cdot &\vline & \cdot & \vline  & \gamma _{s+m+l}-
\end{vmatrix} \\
& \\
0 &\text{otherwise}
\end{cases}
\end{equation}

\begin{equation}
 \psi _{3}(t,\eta ,\xi ) = 
\begin{cases}
2^{k+3s+2m+l-2-j}\quad &\text{if }\quad (t,\eta ,\xi ) \in\begin{vmatrix}
k-1 & \vline & k  & \vline \\
\hline
\cdot &\vline & \cdot & \vline & D^{\left[s-1\atop{ s+m-1\atop s+m+l-1} \right]\times  \cdot} \label{eq 3.11} \\
\hline
\cdot &\vline & \cdot &  \vline  & \alpha _{s}-\\
\hline 
j &\vline & \cdot & \vline  & \beta _{s+m}- \\
\hline 
j &\vline & \cdot & \vline  & \gamma _{s+m+l}-
\end{vmatrix} \\
& \\
0 &\text{otherwise}
\end{cases}
\end{equation}

\end{lem} 
\begin{proof}
We deduce respectively from \eqref{eq 3.6} and \eqref{eq 3.6} with $k\rightarrow k-1$ \\
\begin{align}
 &   \psi _{1}(t,\eta,\xi ) 
  = \sum_{deg Y= k-1}\sum_{deg Z = s-1}E(tYZ)\sum_{deg U \leq s+m-2}E(\eta YU) \sum_{deg V \leq s+m+l-2}E(\xi YV) \label{eq 3.12} \\
  &   = \sum_{deg Y\leq  k-1}\sum_{deg Z = s-1}E(tYZ)\sum_{deg U \leq s+m-2}E(\eta YU) \sum_{deg V \leq s+m+l-2}E(\xi YV) \nonumber \\
& -   \sum_{deg Y \leq k-2}\sum_{deg Z = s-1}E(tYZ)\sum_{deg U \leq s+m-2}E(\eta YU) \sum_{deg V \leq s+m+l-2}E(\xi YV) \nonumber \\
& =  \begin{cases}
 2^{k+3s+2m+l-3-  r( D^{\left[\stackrel{s-1}{\stackrel{s+m-1}{s+m+l-1}}\right] \times k }(t,\eta ,\xi ) )  }  & \text{if }
  r( D^{\left[\stackrel{s-1}{\stackrel{s+m-1}{s+m+l-1}}\right] \times k }(t,\eta ,\xi ) ) =   r( D^{\left[\stackrel{s}{\stackrel{s+m-1}{s+m+l-1}}\right] \times k }(t,\eta ,\xi ) ) \nonumber  \\
     0  & \text{otherwise }.
    \end{cases} \nonumber \\
& -    \begin{cases}
 2^{k+3s+2m+l-4-  r( D^{\left[\stackrel{s-1}{\stackrel{s+m-1}{s+m+l-1}}\right] \times (k-1) }(t,\eta ,\xi ) )  }  & \text{if }
  r( D^{\left[\stackrel{s-1}{\stackrel{s+m-1}{s+m+l-1}}\right] \times (k-1) }(t,\eta ,\xi ) ) =   r( D^{\left[\stackrel{s}{\stackrel{s+m-1}{s+m+l-1}}\right] \times (k-1) }(t,\eta ,\xi ) ) \nonumber  \\
     0  & \text{otherwise }.
    \end{cases} \nonumber 
\end{align}

By Lemma \ref{lem 3.4}\\
\begin{align*}
&  \psi_{1}(t,\eta,\xi  ) \neq 0  
 \Longleftrightarrow 
  r( D^{\left[\stackrel{s-1}{\stackrel{s+m-1}{s+m+l-1}}\right] \times (k-1) }(t,\eta ,\xi ) ) =
  r( D^{\left[\stackrel{s-1}{\stackrel{s+m-1}{s+m+l-1}}\right] \times k }(t,\eta ,\xi ) ) =
  r( D^{\left[\stackrel{s}{\stackrel{s+m-1}{s+m+l-1}}\right] \times (k-1) }(t,\eta ,\xi ) ) \\
  & \\
   & \text{We consider now the following two cases in which $ \psi_{1} $ are different from zero :} \\
  & \text{first case :}\\
  & \\
   &  r( D^{\left[\stackrel{s-1}{\stackrel{s+m-1}{s+m+l-1}}\right] \times (k-1) }(t,\eta ,\xi ) ) =
  r( D^{\left[\stackrel{s-1}{\stackrel{s+m-1}{s+m+l-1}}\right] \times k }(t,\eta ,\xi ) ) \\
  & =   r( D^{\left[\stackrel{s}{\stackrel{s+m-1}{s+m+l-1}}\right] \times (k-1) }(t,\eta ,\xi ) ) =
    r( D^{\left[\stackrel{s}{\stackrel{s+m-1}{s+m+l-1}}\right] \times k }(t,\eta ,\xi ) ) = j \\
    & \\
     & \text{second case :}\\
  & \\
   &  r( D^{\left[\stackrel{s-1}{\stackrel{s+m-1}{s+m+l-1}}\right] \times (k-1) }(t,\eta ,\xi ) ) =
  r( D^{\left[\stackrel{s-1}{\stackrel{s+m-1}{s+m+l-1}}\right] \times k }(t,\eta ,\xi ) ) \\
  & =   r( D^{\left[\stackrel{s}{\stackrel{s+m-1}{s+m+l-1}}\right] \times (k-1) }(t,\eta ,\xi ) ) = j, 
 \quad   r( D^{\left[\stackrel{s}{\stackrel{s+m-1}{s+m+l-1}}\right] \times k }(t,\eta ,\xi ) ) = j +1\\   
& \\
& \text{In the first case we obtain from \eqref{eq 3.12}}\\
& \psi_{1}(t,\eta,\xi  ) =  2^{k+3s+2m+l-3-j} -  2^{k+3s+2m+l-4-j} =  2^{k+3s+2m+l-4-j}\\
& \\
& \text{In the second case we obtain from \eqref{eq 3.12}}\\
& \psi_{1}(t,\eta,\xi  ) =  0 -  2^{k+3s+2m+l-4-j} = - 2^{k+3s+2m+l-4-j}
\end{align*}
This completes the proof of \eqref{eq 3.9}.\\

 \eqref{eq 3.10} and \eqref{eq 3.11} follow respectively from  \eqref{eq 3.2} and  \eqref{eq 3.3}.
\end{proof}
  \begin{lem}
\label{lem 3.6}
We have \\
\begin{equation}
 \psi _{1}(t,\eta ,\xi )\cdot \psi _{2}(t,\eta ,\xi )\cdot \psi _{3}^{q-2}(t,\eta ,\xi ) = 
\begin{cases}
2^{(k+3s+2m+l-3-j)q}\cdot2^{q-3} \quad &\text{if }\quad (t,\eta ,\xi ) \in\begin{vmatrix}
k-1 & \vline & k  & \vline \\
\hline
j &\vline & j & \vline & D^{\left[s-1\atop{ s+m-1\atop s+m+l-1} \right]\times  \cdot} \label{eq 3.13} \\
\hline
j &\vline & j &  \vline  & \alpha _{s}-\\
\hline 
j &\vline & \cdot & \vline  & \beta _{s+m}- \\
\hline 
j &\vline & \cdot & \vline  & \gamma _{s+m+l}-
\end{vmatrix} \\
& \\
-2^{(k+3s+2m+l-3-j)q}\cdot2^{q-3} \quad &\text{if }\quad (t,\eta ,\xi ) \in\begin{vmatrix}
k-1 & \vline & k  & \vline \\
\hline
j &\vline & j & \vline & D^{\left[s-1\atop{ s+m-1\atop s+m+l-1} \right]\times  \cdot}  \\
\hline
j &\vline & j +1&  \vline  & \alpha _{s}-\\
\hline 
j &\vline & j+1 & \vline  & \beta _{s+m}- \\
\hline 
j &\vline & j+1 & \vline  & \gamma _{s+m+l}-
\end{vmatrix} \\
& \\
0 &\text{otherwise}
\end{cases}
\end{equation}
\end{lem}
\begin{proof}
Combining \eqref{eq 3.9},\eqref{eq 3.10} and \eqref{eq 3.11}we deduce \eqref{eq 3.13}.

\end{proof}
   \begin{lem}
\label{lem 3.7}

Let $ (t,\eta ,\xi ) \in  \mathbb{P}\times \mathbb{P}\times \mathbb{P}  $ and
\begin{align*}
 \theta _{1}(t,\eta,\xi) & = \sum_{deg Y= k-1}\sum_{deg Z\leq  s-1}E(tYZ)\sum_{deg U = s+m-1}E(\eta YU) \sum_{deg V \leq s+m+l-2}E(\xi  YV ) \\
  \theta _{2}(t,\eta,\xi) & = \sum_{deg Y= k-1}\sum_{deg Z\leq  s-1}E(tYZ)\sum_{deg U \leq s+m-1}E(\eta YU) \sum_{deg V = s+m+l-1}E(\xi  YV ) \\
  \end{align*}
Then \begin{equation}
 \theta  _{1}(t,\eta ,\xi ) = 
\begin{cases}
2^{k+3s+2m+l-3-j}\quad &\text{if }\quad (t,\eta ,\xi ) \in\begin{vmatrix}
k-1 & \vline & k  & \vline \\
\hline
\cdot &\vline &  \cdot  &\vline & D^{\left[s-1\atop{ s+m-1\atop s+m+l-1} \right]\times  \cdot} \label{eq 3.14} \\
\hline
j &\vline & j &  \vline  & \alpha _{s}-\\
\hline 
j &\vline & j & \vline  & \beta _{s+m}- \\
\hline 
\cdot &\vline & \cdot & \vline  & \gamma _{s+m+l}-
\end{vmatrix} \\
& \\
-2^{k+3s+2m+l-3-j} \quad &\text{if }\quad (t,\eta ,\xi ) \in\begin{vmatrix}
k-1 & \vline & k  & \vline \\
\hline
\cdot &\vline & \cdot & \vline & D^{\left[s-1\atop{ s+m-1\atop s+m+l-1} \right]\times  \cdot}  \\
\hline
j &\vline & j &  \vline  & \alpha _{s}-\\
\hline 
j &\vline & j+1 & \vline  & \beta _{s+m}- \\
\hline 
\cdot &\vline & \cdot & \vline  & \gamma _{s+m+l}-
\end{vmatrix} \\
& \\
0 &\text{otherwise}
\end{cases}
\end{equation}

 \begin{equation}
 \theta  _{2}(t,\eta ,\xi ) = 
\begin{cases}
2^{k+3s+2m+l-3-j}\quad &\text{if }\quad (t,\eta ,\xi ) \in\begin{vmatrix}
k-1 & \vline & k  & \vline \\
\hline
\cdot &\vline &  \cdot & \vline & D^{\left[s-1\atop{ s+m-1\atop s+m+l-1} \right]\times  \cdot} \label{eq 3.15} \\
\hline
\cdot &\vline & \cdot &  \vline  & \alpha _{s}-\\
\hline 
j &\vline & j & \vline  & \beta _{s+m}- \\
\hline 
j &\vline & j & \vline  & \gamma _{s+m+l}-
\end{vmatrix} \\
& \\
-2^{k+3s+2m+l-3-j} \quad &\text{if }\quad (t,\eta ,\xi ) \in\begin{vmatrix}
k-1 & \vline & k  & \vline \\
\hline
\cdot &\vline & \cdot & \vline & D^{\left[s-1\atop{ s+m-1\atop s+m+l-1} \right]\times  \cdot}  \\
\hline
\cdot &\vline & \cdot &  \vline  & \alpha _{s}-\\
\hline 
j &\vline & j & \vline  & \beta _{s+m}- \\
\hline 
j &\vline & j+1 & \vline  & \gamma _{s+m+l}-
\end{vmatrix} \\
& \\
0 &\text{otherwise}
\end{cases}
\end{equation}
\end{lem}

\begin{proof}
Similarly to the proof of \eqref{eq 3.9}.
\end{proof}
  \begin{lem}
\label{lem 3.8}
We have \\
\begin{equation}
 h _{1}(t,\eta ,\xi )\cdot \theta  _{1}(t,\eta ,\xi )\cdot \psi _{3}^{q-2}(t,\eta ,\xi ) = 
\begin{cases}
2^{(k+3s+2m+l-3-j)q}\cdot2^{q-2} \quad &\text{if }\quad (t,\eta ,\xi ) \in\begin{vmatrix}
k-1 & \vline & k  & \vline \\
\hline
j &\vline & j & \vline & D^{\left[s-1\atop{ s+m-1\atop s+m+l-1} \right]\times  \cdot} \label{eq 3.16} \\
\hline
j &\vline & j &  \vline  & \alpha _{s}-\\
\hline 
j &\vline & j & \vline  & \beta _{s+m}- \\
\hline 
j &\vline & \cdot & \vline  & \gamma _{s+m+l}-
\end{vmatrix} \\
& \\
-2^{(k+3s+2m+l-3-j)q}\cdot2^{q-2} \quad &\text{if }\quad (t,\eta ,\xi ) \in \begin{vmatrix}
k-1 & \vline & k  & \vline \\
\hline
j &\vline & j & \vline & D^{\left[s-1\atop{ s+m-1\atop s+m+l-1} \right]\times  \cdot}  \\
\hline
j &\vline & j &  \vline  & \alpha _{s}-\\
\hline 
j &\vline & j+1 & \vline  & \beta _{s+m}- \\
\hline 
j &\vline & \cdot & \vline  & \gamma _{s+m+l}-
\end{vmatrix} \\
& \\
0 &\text{otherwise}
\end{cases}
\end{equation}
\end{lem}
\begin{proof}
Combining \eqref{eq 3.6},\eqref{eq 3.14} and \eqref{eq 3.11}we deduce \eqref{eq 3.16}.

\end{proof}
  \begin{lem}
\label{lem 3.9}
We have \\
\begin{equation}
 h _{1}^{q-2}(t,\eta ,\xi )\cdot \psi   _{2}(t,\eta ,\xi )\cdot \theta  _{2}(t,\eta ,\xi ) = 
\begin{cases}
2^{(k+3s+2m+l-3-j)q+1} \quad &\text{if }\quad (t,\eta ,\xi ) \in\begin{vmatrix}
k-1 & \vline & k  & \vline \\
\hline
j &\vline & j & \vline & D^{\left[s-1\atop{ s+m-1\atop s+m+l-1} \right]\times  \cdot} \label{eq 3.17} \\
\hline
j &\vline & j &  \vline  & \alpha _{s}-\\
\hline 
j &\vline & j & \vline  & \beta _{s+m}- \\
\hline 
j &\vline & j & \vline  & \gamma _{s+m+l}-
\end{vmatrix} \\
& \\
-2^{(k+3s+2m+l-3-j)q+1} \quad &\text{if }\quad (t,\eta ,\xi ) \in\begin{vmatrix}
k-1 & \vline & k  & \vline \\
\hline
j &\vline & j & \vline & D^{\left[s-1\atop{ s+m-1\atop s+m+l-1} \right]\times  \cdot}  \\
\hline
j &\vline & j &  \vline  & \alpha _{s}-\\
\hline 
j &\vline & j & \vline  & \beta _{s+m}- \\
\hline 
j &\vline & j+1 & \vline  & \gamma _{s+m+l}-
\end{vmatrix} \\
& \\
0 &\text{otherwise}
\end{cases}
\end{equation}
\end{lem}
\begin{proof}
Combining \eqref{eq 3.6},\eqref{eq 3.10} and \eqref{eq 3.15}we deduce \eqref{eq 3.17}.

\end{proof}

  \subsection{Some rank formulas for partitions of triple persymmetric matrices related to $ \sigma _{i,i,i,i}^{\left[\alpha \atop{ \beta \atop \gamma } \right]\times k}$}
  \label{subsec 4} 

   \begin{lem}
\label{lem 3.10} 
We have for $q\geq 3$
$$\displaystyle 
 \int_{\mathbb{P}\times \mathbb{P}\times \mathbb{P}}  \psi _{1}(t,\eta ,\xi )\cdot \psi _{2}(t,\eta ,\xi )\cdot \psi _{3}^{q-2}(t,\eta ,\xi )  dt d\eta d\xi  =  0. $$
\end{lem}
\begin{proof}
The integral above  is equal to the number of solutions \\
 $ (Y_{1},Z_{1},U_{1},V_{1},Y_{2},Z_{2},U_{2},V_{2},\ldots, Y_{q},Z_{q},U_{q},V_{q}) $
of the polynomial equations
    \[\left\{\begin{array}{cc}
 Y_{1}Z_{1} +  Y_{2}Z_{2}+\ldots + Y_{q}Z_{q}= 0, \\
 Y_{1}U_{1} +  Y_{2}U_{2} +\ldots + Y_{q}U_{q}= 0,\\
  Y_{1}V_{1} +  Y_{2}V_{2} +\ldots + Y_{q}V_{q}= 0.
    \end{array}\right.\]\\ 
 satisfying the degree conditions \\
 \[\left\{\begin{array}{cc}
   \deg Y_{1} = k-1,\quad  \deg Z_{1} = s-1  \quad  \deg U_{1} \leq s+m-2 \quad  \deg V_{1}\leq s+m+l-2  \\
 \deg Y_{2} \leq  k-2,\quad  \deg Z_{2} \leq  s-1,  \quad  \deg U_{2} = s+m-1, \quad  \deg V_{2}\leq s+m+l-2  \\
 \deg Y_{i} \leq  k-2,\quad  \deg Z_{i} \leq  s-1,  \quad  \deg U_{i}\leq  s+m-1, \quad  \deg V_{i}= s+m+l-1\quad \text{for}\quad 3\leq i\leq q 
 \end{array}\right.\] \\
  By degree considerations  we have that
   $ \deg \left[ Y_{1}Z_{1} +  Y_{2}Z_{2}+\ldots + Y_{q}Z_{q}\right] $ is equal to s + k - 2.
 Lemma \ref{lem 3.10} follows.
 \end{proof}
    \begin{lem}
\label{lem 3.11} 
We have for $q\geq 3$
$$\displaystyle 
 \int_{\mathbb{P}\times \mathbb{P}\times \mathbb{P}} h _{1}(t,\eta ,\xi )\cdot \theta  _{1}(t,\eta ,\xi )\cdot \psi _{3}^{q-2}(t,\eta ,\xi )  dt d\eta d\xi  =  0. $$
\end{lem}
\begin{proof}
The integral above  is equal to the number of solutions \\
 $ (Y_{1},Z_{1},U_{1},V_{1},Y_{2},Z_{2},U_{2},V_{2},\ldots, Y_{q},Z_{q},U_{q},V_{q}) $
of the polynomial equations
    \[\left\{\begin{array}{cc}
 Y_{1}Z_{1} +  Y_{2}Z_{2}+\ldots + Y_{q}Z_{q}= 0, \\
 Y_{1}U_{1} +  Y_{2}U_{2} +\ldots + Y_{q}U_{q}= 0,\\
  Y_{1}V_{1} +  Y_{2}V_{2} +\ldots + Y_{q}V_{q}= 0.
    \end{array}\right.\]\\ 
 satisfying the degree conditions \\
 \[\left\{\begin{array}{cc}
   \deg Y_{1} \leq  k-1,\quad  \deg Z_{1} = s-1  \quad  \deg U_{1} \leq s+m-2 \quad  \deg V_{1}\leq s+m+l-2  \\
 \deg Y_{2} = k-1,\quad  \deg Z_{2} \leq  s-1,  \quad  \deg U_{2} = s+m-1, \quad  \deg V_{2}\leq s+m+l-2  \\
 \deg Y_{i} \leq  k-2,\quad  \deg Z_{i} \leq  s-1,  \quad  \deg U_{i}\leq  s+m-1, \quad  \deg V_{i}= s+m+l-1\quad \text{for}\quad 3\leq i\leq q 
 \end{array}\right.\] \\
  By degree considerations  we have that
   $ \deg \left[ Y_{1}U_{1} +  Y_{2}U_{2}+\ldots + Y_{q}U_{q}\right] $ is equal to s+k+m-2.
 Lemma \ref{lem 3.11} follows.
 \end{proof}
     \begin{lem}
\label{lem 3.12} 
We have for $q\geq 3$
$$\displaystyle 
 \int_{\mathbb{P}\times \mathbb{P}\times \mathbb{P}} h _{1}^{q-2}(t,\eta ,\xi )\cdot \psi  _{2}(t,\eta ,\xi )\cdot \theta  _{2}(t,\eta ,\xi )  dt d\eta d\xi  =  0. $$
\end{lem}
\begin{proof}
The integral above  is equal to the number of solutions \\
 $ (Y_{1},Z_{1},U_{1},V_{1},Y_{2},Z_{2},U_{2},V_{2},\ldots, Y_{q},Z_{q},U_{q},V_{q}) $
of the polynomial equations
    \[\left\{\begin{array}{cc}
 Y_{1}Z_{1} +  Y_{2}Z_{2}+\ldots + Y_{q}Z_{q}= 0, \\
 Y_{1}U_{1} +  Y_{2}U_{2} +\ldots + Y_{q}U_{q}= 0,\\
  Y_{1}V_{1} +  Y_{2}V_{2} +\ldots + Y_{q}V_{q}= 0.
    \end{array}\right.\]\\ 
 satisfying the degree conditions \\
 \[\left\{\begin{array}{cc}
   \deg Y_{1} \leq  k-2,\quad  \deg Z_{1} \leq  s-1  \quad  \deg U_{1} = s+m-1 \quad  \deg V_{1}\leq s+m+l-2  \\
 \deg Y_{2} = k-1,\quad  \deg Z_{2} \leq  s-1,  \quad  \deg U_{2} \leq  s+m-1, \quad  \deg V_{2}= s+m+l-1  \\
 \deg Y_{i} \leq  k-1,\quad  \deg Z_{i} = s-1,  \quad  \deg U_{i}\leq  s+m-2, \quad  \deg V_{i}\leq  s+m+l-2\quad \text{for}\quad 3\leq i\leq q 
 \end{array}\right.\] \\
  By degree considerations  we have that
   $ \deg \left[ Y_{1}V_{1} +  Y_{2}V_{2}+\ldots + Y_{q}V_{q}\right] $ is equal to k+s+m+l-2
 Lemma \ref{lem 3.12} follows.
 \end{proof}
 
 \begin{lem}
\label{lem 3.13} We have the following rank identity for  $0\leq j\leq \inf(3s+2m+l-3,k-1) :$
  \begin{align}
  &   {}^{\#}\left(\begin{array}{c | c}
           j & j \\
           \hline
           j & j +1 \\
           \hline
            j & j +1 \\
            \hline 
            j & j+1
           \end{array} \right)_{\mathbb{P}/\mathbb{P}_{k+s -1}\times
           \mathbb{P}/\mathbb{P}_{k+s+m-1}\times  \mathbb{P}/\mathbb{P}_{k+s+m+l-1} }^{{\alpha \over {\beta \over \gamma }}} = 
{}^{\#}\left(\begin{array}{c | c}
         j & j \\
           \hline
           j & j  \\
           \hline
            j & j  \\
            \hline 
            j & j
           \end{array} \right)_{\mathbb{P}/\mathbb{P}_{k+s -1}\times
           \mathbb{P}/\mathbb{P}_{k+s+m-1}\times  \mathbb{P}/\mathbb{P}_{k+s+m+l-1} }^{{\alpha \over {\beta \over \gamma }}} \label{eq 3.18} \\
            &    +   {}^{\#}\left(\begin{array}{c | c}
              j & j \\
           \hline
           j & j  \\
           \hline
            j & j  \\
            \hline 
            j & j+1
           \end{array} \right)_{\mathbb{P}/\mathbb{P}_{k+s -1}\times
           \mathbb{P}/\mathbb{P}_{k+s+m-1}\times  \mathbb{P}/\mathbb{P}_{k+s+m+l-1} }^{{\alpha \over {\beta \over \gamma }}}  
               +  {}^{\#}\left(\begin{array}{c | c}
              j & j \\
           \hline
           j & j  \\
           \hline
            j & j +1 \\
            \hline 
            j & j +1
           \end{array} \right)_{\mathbb{P}/\mathbb{P}_{k+s -1}\times
           \mathbb{P}/\mathbb{P}_{k+s+m-1}\times  \mathbb{P}/\mathbb{P}_{k+s+m+l-1} }^{{\alpha \over {\beta \over \gamma }}} \nonumber
          \end{align} 
\end{lem}
\begin{proof}
   We define  
  $$\begin{array}{l} \mathbb{A} =  \Big\{(t,\eta,\xi ) \in \mathbb{P}^{3}
\mid   r( D^{\left[\stackrel{s-1}{\stackrel{s+m-1}{s+m+l-1}}\right] \times (k-1) }(t,\eta ,\xi ) ) 
  =  r( D^{\left[\stackrel{s-1}{\stackrel{s+m-1}{s+m+l-1}}\right] \times k }(t,\eta ,\xi ) )  
 =      r( D^{\left[\stackrel{s}{\stackrel{s+m-1}{s+m+l-1}}\right] \times (k-1) }(t,\eta ,\xi ) )  \\
    =    r( D^{\left[\stackrel{s}{\stackrel{s+m-1}{s+m+l-1}}\right] \times k }(t,\eta ,\xi ) )  
 =   r( D^{\left[\stackrel{s}{\stackrel{s+m}{s+m+l-1}}\right] \times (k-1) }(t,\eta ,\xi ) ) = r( D^{\left[\stackrel{s}{\stackrel{s+m}{s+m+l}}\right] \times (k-1) }(t,\eta ,\xi ) ) \Big\} \end{array} $$\\
 and 
  $$\begin{array}{l}  \mathbb{B} =   \Big\{(t,\eta,\xi ) \in \mathbb{P}^{3}
\mid   r( D^{\left[\stackrel{s-1}{\stackrel{s+m-1}{s+m+l-1}}\right] \times (k-1) }(t,\eta ,\xi ) ) 
  =  r( D^{\left[\stackrel{s-1}{\stackrel{s+m-1}{s+m+l-1}}\right] \times k }(t,\eta ,\xi ) )  
 =      r( D^{\left[\stackrel{s}{\stackrel{s+m-1}{s+m+l-1}}\right] \times (k-1) }(t,\eta ,\xi ) )  \\
    =    r( D^{\left[\stackrel{s}{\stackrel{s+m}{s+m+l-1}}\right] \times (k-1) }(t,\eta ,\xi ) )  
 =   r( D^{\left[\stackrel{s}{\stackrel{s+m}{s+m+l}}\right] \times (k-1) }(t,\eta ,\xi ) ) \quad \text{and}\\
  r( D^{\left[\stackrel{s}{\stackrel{s+m-1}{s+m+l-1}}\right] \times k }(t,\eta ,\xi ) ) =
  r( D^{\left[\stackrel{s}{\stackrel{s+m}{s+m+l-1}}\right] \times k }(t,\eta ,\xi ) )   =              
  r( D^{\left[\stackrel{s}{\stackrel{s+m}{s+m+l}}\right] \times k }(t,\eta ,\xi ) ) = 
  r( D^{\left[\stackrel{s-1}{\stackrel{s+m-1}{s+m+l-1}}\right] \times (k-1) }(t,\eta ,\xi ) )  +1  \Big\} \end{array} $$\\
  
    By \eqref{eq 3.13} we have,  observing that   $  \psi _{1}(t,\eta ,\xi )\cdot \psi _{2}(t,\eta ,\xi )\cdot \psi _{3}^{q-2}(t,\eta ,\xi )  $
     is constant on cosets of $ \mathbb{P}_{k+s-1}\times \mathbb{P}_{k+s+m-1}\times \mathbb{P}_{k+s+m+l-1} $
     
  \begin{align*}
&  \int_{\mathbb{P}^3} \psi _{1}(t,\eta ,\xi )\cdot \psi _{2}(t,\eta ,\xi )\cdot \psi _{3}^{q-2}(t,\eta ,\xi )dtd\eta d\xi \\
&  = \int_{  \mathbb{A}} \psi _{1}(t,\eta ,\xi )\cdot \psi _{2}(t,\eta ,\xi )\cdot \psi _{3}^{q-2}(t,\eta ,\xi )dtd\eta d\xi 
    +  \int_{ \mathbb{B}}  \psi _{1}(t,\eta ,\xi )\cdot \psi _{2}(t,\eta ,\xi )\cdot \psi _{3}^{q-2}(t,\eta ,\xi )dtd\eta d\xi  \\
&  =\sum_{(t,\eta,\xi  ) \in \mathbb{A}\bigcap\big(\mathbb{P}/\mathbb{P}_{k+s-1}\times\mathbb{P}/\mathbb{P}_{k+s+m-1}\times\mathbb{P}/\mathbb{P}_{k+s+m+l-1}\big) }
2^{q-3}\cdot2^{(k+3s+2m+l-3-  r( D^{\left[\stackrel{s-1}{\stackrel{s+m-1}{s+m+l-1}}\right] \times (k-1) }(t,\eta ,\xi ))q } \\
&  \int_{\mathbb{P}_{s+k-1}}dt \int_{\mathbb{P}_{s+k+m-1}}d \eta \int_{\mathbb{P}_{s+k+m+l-1}}d\xi   \\
&  + \sum_{(t,\eta,\xi  ) \in \mathbb{B}\bigcap\big(\mathbb{P}/\mathbb{P}_{k+s-1}\times\mathbb{P}/\mathbb{P}_{k+s+m-1}\times\mathbb{P}/\mathbb{P}_{k+s+m+l-1}\big) }
-2^{q-3}\cdot2^{(k+3s+2m+l-3-  r( D^{\left[\stackrel{s-1}{\stackrel{s+m-1}{s+m+l-1}}\right] \times (k-1) }(t,\eta ,\xi ))q } \\
&  \int_{\mathbb{P}_{s+k-1}}dt \int_{\mathbb{P}_{s+k+m-1}}d \eta \int_{\mathbb{P}_{s+k+m+l-1}}d\xi \\
& =  \sum_{j=0}^{\inf(3s+2m+l-3,k-1)} 2^{q-3}\cdot 2^{(k+3s+2m+l-3-j) q}\cdot     
     {}^{\#}\left(\begin{array}{c | c}
           j & j \\
           \hline
           j & j  \\
           \hline
            j & j  \\
            \hline 
            j & j
           \end{array} \right)_{\mathbb{P}/\mathbb{P}_{k+s -1}\times
           \mathbb{P}/\mathbb{P}_{k+s+m-1}\times  \mathbb{P}/\mathbb{P}_{k+s+m+l-1} }^{{\alpha \over {\beta \over \gamma }}}
          \cdot2^{-(3k+3s+2m+l-3)}  \\
      &  +  \sum_{j=0}^{\inf(3s+2m+l-3,k-1)} 2^{q-3}\cdot 2^{(k+3s+2m+l-3-j) q}\cdot     
     {}^{\#}\left(\begin{array}{c | c}
           j & j \\
           \hline
           j & j  \\
           \hline
            j & j  \\
            \hline 
            j & j+1
           \end{array} \right)_{\mathbb{P}/\mathbb{P}_{k+s -1}\times
           \mathbb{P}/\mathbb{P}_{k+s+m-1}\times  \mathbb{P}/\mathbb{P}_{k+s+m+l-1} }^{{\alpha \over {\beta \over \gamma }}}
          \cdot2^{-(3k+3s+2m+l-3)}  \\   
          &  + \sum_{j=0}^{\inf(3s+2m+l-3,k-1)} 2^{q-3}\cdot 2^{(k+3s+2m+l-3-j) q}\cdot     
     {}^{\#}\left(\begin{array}{c | c}
           j & j \\
           \hline
           j & j  \\
           \hline
            j & j +1 \\
            \hline 
            j & j+1
           \end{array} \right)_{\mathbb{P}/\mathbb{P}_{k+s -1}\times
           \mathbb{P}/\mathbb{P}_{k+s+m-1}\times  \mathbb{P}/\mathbb{P}_{k+s+m+l-1} }^{{\alpha \over {\beta \over \gamma }}}
          \cdot2^{-(3k+3s+2m+l-3)}  \\ 
          & - \sum_{j=0}^{\inf(3s+2m+l-3,k-1)} 2^{q-3}\cdot 2^{(k+3s+2m+l-3-j) q}\cdot     
     {}^{\#}\left(\begin{array}{c | c}
           j & j \\
           \hline
           j & j +1 \\
           \hline
            j & j +1 \\
            \hline 
            j & j+1
           \end{array} \right)_{\mathbb{P}/\mathbb{P}_{k+s -1}\times
           \mathbb{P}/\mathbb{P}_{k+s+m-1}\times  \mathbb{P}/\mathbb{P}_{k+s+m+l-1} }^{{\alpha \over {\beta \over \gamma }}}
          \cdot2^{-(3k+3s+2m+l-3)}  \\
          \end{align*}
   Now by Lemma \ref{lem 3.10} we get  for all  $ q \geq 3 $ \\
   
   \begin{align*}
&  \sum_{j=0}^{\inf(3s+2m+l-3,k-1)} 2^{-j q}\cdot\Bigg[     
     {}^{\#}\left(\begin{array}{c | c}
           j & j \\
           \hline
           j & j  \\
           \hline
            j & j  \\
            \hline 
            j & j
           \end{array} \right)_{\mathbb{P}/\mathbb{P}_{k+s -1}\times
           \mathbb{P}/\mathbb{P}_{k+s+m-1}\times  \mathbb{P}/\mathbb{P}_{k+s+m+l-1} }^{{\alpha \over {\beta \over \gamma }}} \\
          &      +    {}^{\#}\left(\begin{array}{c | c}
           j & j \\
           \hline
           j & j  \\
           \hline
            j & j  \\
            \hline 
            j & j+1
           \end{array} \right)_{\mathbb{P}/\mathbb{P}_{k+s -1}\times
           \mathbb{P}/\mathbb{P}_{k+s+m-1}\times  \mathbb{P}/\mathbb{P}_{k+s+m+l-1} }^{{\alpha \over {\beta \over \gamma }}}
              +    {}^{\#}\left(\begin{array}{c | c}
           j & j \\
           \hline
           j & j \\
           \hline
            j & j +1 \\
            \hline 
            j & j+1
           \end{array} \right)_{\mathbb{P}/\mathbb{P}_{k+s -1}\times
           \mathbb{P}/\mathbb{P}_{k+s+m-1}\times  \mathbb{P}/\mathbb{P}_{k+s+m+l-1} }^{{\alpha \over {\beta \over \gamma }}}\\
           &     -  {}^{\#}\left(\begin{array}{c | c}
           j & j \\
           \hline
           j & j +1 \\
           \hline
            j & j +1 \\
            \hline 
            j & j+1
           \end{array} \right)_{\mathbb{P}/\mathbb{P}_{k+s -1}\times
           \mathbb{P}/\mathbb{P}_{k+s+m-1}\times  \mathbb{P}/\mathbb{P}_{k+s+m+l-1} }^{{\alpha \over {\beta \over \gamma }}}\Bigg] = 0.
    \end{align*}

\end{proof}

 \begin{lem}
\label{lem 3.14} We have the following rank identities for  $0\leq j\leq \inf(3s+2m+l-3,k-1) :$
  \begin{align}
  &   {}^{\#}\left(\begin{array}{c | c}
           j & j \\
           \hline
           j & j  \\
           \hline
            j & j +1 \\
            \hline 
            j & j+1
           \end{array} \right)_{\mathbb{P}/\mathbb{P}_{k+s -1}\times
           \mathbb{P}/\mathbb{P}_{k+s+m-1}\times  \mathbb{P}/\mathbb{P}_{k+s+m+l-1} }^{{\alpha \over {\beta \over \gamma }}} = 
{}^{\#}\left(\begin{array}{c | c}
         j & j \\
           \hline
           j & j  \\
           \hline
            j & j  \\
            \hline 
            j & j
           \end{array} \right)_{\mathbb{P}/\mathbb{P}_{k+s -1}\times
           \mathbb{P}/\mathbb{P}_{k+s+m-1}\times  \mathbb{P}/\mathbb{P}_{k+s+m+l-1} }^{{\alpha \over {\beta \over \gamma }}} \label{eq 3.19} \\
            &    +   {}^{\#}\left(\begin{array}{c | c}
              j & j \\
           \hline
           j & j  \\
           \hline
            j & j  \\
            \hline 
            j & j+1
           \end{array} \right)_{\mathbb{P}/\mathbb{P}_{k+s -1}\times
           \mathbb{P}/\mathbb{P}_{k+s+m-1}\times  \mathbb{P}/\mathbb{P}_{k+s+m+l-1} }^{{\alpha \over {\beta \over \gamma }}} \nonumber \\
           & \text{and} \nonumber \\
           & \nonumber \\
             &   {}^{\#}\left(\begin{array}{c | c}
           j & j \\
           \hline
           j & j  \\
           \hline
            j & j  \\
            \hline 
            j & j
           \end{array} \right)_{\mathbb{P}/\mathbb{P}_{k+s -1}\times
           \mathbb{P}/\mathbb{P}_{k+s+m-1}\times  \mathbb{P}/\mathbb{P}_{k+s+m+l-1} }^{{\alpha \over {\beta \over \gamma }}} = 
{}^{\#}\left(\begin{array}{c | c}
         j & j \\
           \hline
           j & j  \\
           \hline
            j & j  \\
            \hline 
            j & j+1
           \end{array} \right)_{\mathbb{P}/\mathbb{P}_{k+s -1}\times
           \mathbb{P}/\mathbb{P}_{k+s+m-1}\times  \mathbb{P}/\mathbb{P}_{k+s+m+l-1} }^{{\alpha \over {\beta \over \gamma }}} \label{eq 3.20} 
        \end{align} 
\end{lem}
\begin{proof}
The proofs of \eqref{eq 3.19} and \eqref{eq 3.20} are somewhat similar to the proof of \eqref{eq 3.18}
applying respectively \eqref{eq 3.16}, Lemma \ref{lem 3.11} and  \eqref{eq 3.17}, Lemma  \ref{lem 3.12}

\end{proof}

  \begin{lem}
\label{lem 3.15} We have the following rank identities for all integers j for which  $0\leq j\leq \inf(3s+2m+l-3,k-2) :$
  \begin{align}
  &   {}^{\#}\left(\begin{array}{c | c}
           j & j+1 \\
           \hline
           j & j +1 \\
           \hline
            j & j +1 \\
            \hline 
            j & j+1
           \end{array} \right)_{\mathbb{P}/\mathbb{P}_{k+s -1}\times
           \mathbb{P}/\mathbb{P}_{k+s+m-1}\times  \mathbb{P}/\mathbb{P}_{k+s+m+l-1} }^{{\alpha \over {\beta \over \gamma }}} = 0  \label{eq 3.21}
\end{align}
\end{lem}
\begin{proof}

Proof by contradiction. 
Assume on the contrary that there exists $j_{0} \in [0, \inf(3s+2m+l-3,k-2)]$ such that 
\begin{equation}
\label{eq 3.22}
   {}^{\#}\left(\begin{array}{c | c}
           j_{0} & j_{0}+1 \\
           \hline
           j_{0} & j_{0} +1 \\
           \hline
            j_{0} & j_{0} +1 \\
            \hline 
            j_{0} & j_{0}+1
           \end{array} \right)_{\mathbb{P}/\mathbb{P}_{k+s -1}\times
           \mathbb{P}/\mathbb{P}_{k+s+m-1}\times  \mathbb{P}/\mathbb{P}_{k+s+m+l-1} }^{{\alpha \over {\beta \over \gamma }}} > 0.
\end{equation}
We are going to show that 
\begin{align*}
& {}^{\#}\left(\begin{array}{c | c}
           j_{0} & j_{0}+1 \\
           \hline
           j_{0} & j_{0} +1 \\
           \hline
            j_{0} & j_{0} +1 \\
            \hline 
            j_{0} & j_{0}+1
           \end{array} \right)_{\mathbb{P}/\mathbb{P}_{k+s -1}\times
           \mathbb{P}/\mathbb{P}_{k+s+m-1}\times  \mathbb{P}/\mathbb{P}_{k+s+m+l-1} }^{{\alpha \over {\beta \over \gamma }}} > 0 \\
    &  \Longrightarrow 
 {}^{\#}\left(\begin{array}{c | c}
           j_{0}-1 & j_{0} \\
           \hline
           j_{0}-1 & j_{0}  \\
           \hline
            j_{0}-1 & j_{0}  \\
            \hline 
            j_{0}-1 & j_{0}
           \end{array} \right)_{\mathbb{P}_{1}/\mathbb{P}_{k+s -1}\times
           \mathbb{P}_{1}/\mathbb{P}_{k+s+m-1}\times  \mathbb{P}_{1}/\mathbb{P}_{k+s+m+l-1} }^{{\alpha \over {\beta \over \gamma }}} > 0 \\
      &  \Longrightarrow \ldots \Longrightarrow     {}^{\#}\left(\begin{array}{c | c}
           0 & 1 \\
           \hline
           0 & 1  \\
           \hline
            0 & 1  \\
            \hline 
            0 & 1
           \end{array} \right)_{\mathbb{P}_{j_{0}}/\mathbb{P}_{k+s -1}\times
           \mathbb{P}_{j_{0}}/\mathbb{P}_{k+s+m-1}\times  \mathbb{P}_{j_{0}}/\mathbb{P}_{k+s+m+l-1} }^{{\alpha \over {\beta \over \gamma }}} > 0 \\   
     & \text{which obviously contradicts} \\
     &    {}^{\#}\left(\begin{array}{c | c}
           0 & 1 \\
           \hline
           0 & 1  \\
           \hline
            0 & 1  \\
            \hline 
            0 & 1
           \end{array} \right)_{\mathbb{P}_{j_{0}}/\mathbb{P}_{k+s -1}\times
           \mathbb{P}_{j_{0}}/\mathbb{P}_{k+s+m-1}\times  \mathbb{P}_{j_{0}}/\mathbb{P}_{k+s+m+l-1} }^{{\alpha \over {\beta \over \gamma }}} = 0      
\end{align*}
By \eqref{eq 3.22} there exsists $(t_{0},\eta _{0},\xi _{0}) \in \mathbb{P}/\mathbb{P}_{k+s -1}\times
           \mathbb{P}/\mathbb{P}_{k+s+m-1}\times  \mathbb{P}/\mathbb{P}_{k+s+m+l-1}$ such that \\
 \newpage

 \begin{equation}
\label{eq 3.23}  \left\{  \begin{array}{cc}  
 r( D^{\left[s-1\atop{ s+m-1\atop s+m+l-1} \right]\times (k-1)}(t_{0},\eta_{0} ,\xi_{0} ) )
=  r( D^{\left[s\atop{ s+m-1\atop s+m+l-1} \right]\times (k-1)}(t_{0},\eta_{0} ,\xi_{0} ) ) \\
=   r( D^{\left[s\atop{ s+m\atop s+m+l-1} \right]\times (k-1)}(t_{0},\eta_{0} ,\xi_{0} ) )     
 =  r( D^{\left[s\atop{ s+m\atop s+m+l} \right]\times (k-1)}(t_{0},\eta_{0} ,\xi_{0} ) ) = j_{0},  \\
  r( D^{\left[s-1\atop{ s+m-1\atop s+m+l-1} \right]\times k}(t_{0},\eta_{0} ,\xi_{0} ))
=  r( D^{\left[s\atop{ s+m-1\atop s+m+l-1} \right]\times k}(t_{0},\eta_{0} ,\xi_{0} ) ) \\
=   r( D^{\left[s\atop{ s+m\atop s+m+l-1} \right]\times  k}(t_{0},\eta_{0} ,\xi_{0} ) )     
 =  r( D^{\left[s\atop{ s+m\atop s+m+l} \right]\times k}(t_{0},\eta_{0} ,\xi_{0} ) ) = j_{0}+1  
  \end{array}\right. 
  \end{equation}

 Consider the following partition of the matrix
 $$ \displaystyle  D^{ \left[s-1\atop{ s+m-1\atop {s+m+l-1 \atop{\overline{\alpha _{s -} \atop{\beta _{s+m -} \atop \gamma_{s+m+l -}} }}}}\right]\times k}(t_{0},\eta_{0} ,\xi_{0} ) $$\\
   $$   \left ( \begin{array} {c|cccc|c}
\alpha _{1} & \alpha _{2} & \alpha _{3} &  \ldots & \alpha _{k-1}  &  \alpha _{k} \\
\beta  _{1} & \beta  _{2} & \beta  _{3} & \ldots  &  \beta_{k-1} &  \beta _{k}  \\
\gamma  _{1} & \gamma   _{2} &  \gamma _{3} & \ldots  & \gamma  _{k-1} &  \gamma  _{k}  \\
\hline
\alpha _{2 } & \alpha _{3} & \alpha _{4}&  \ldots  &  \alpha _{k} &  \alpha _{k+1} \\
\beta  _{2} & \beta  _{3} & \beta  _{4} & \ldots  &  \beta_{k} &  \beta _{k+1}  \\
\gamma  _{2} & \gamma  _{3} & \gamma  _{4} & \ldots  & \gamma  _{k} &  \gamma  _{k+1}  \\
\vdots & \vdots & \vdots    &  \vdots & \vdots  &  \vdots \\
\alpha _{s-1} & \alpha _{s} & \alpha _{s +1} & \ldots  &  \alpha _{s+k-3} &  \alpha _{s+k-2}  \\
\beta  _{s-1} & \beta  _{s} & \beta  _{s+1} & \ldots  &  \beta_{s+k-3} &  \beta _{s+k-2}  \\
\gamma  _{s-1} & \gamma  _{s} &  \gamma  _{s+1} & \ldots  & \gamma  _{s+k-3} &  \gamma  _{s+k-2}  \\
\beta  _{s} & \beta  _{s+1} & \beta  _{s+2} & \ldots  &  \beta_{s+k-2} &  \beta _{s+k-1}  \\
\beta  _{s+1} & \beta  _{s+2} & \beta  _{s+3} & \ldots  &  \beta_{k+s-1} &  \beta _{k+s}  \\
\vdots & \vdots & \vdots    &  \vdots & \vdots  &  \vdots \\
\beta  _{s+m-1} & \beta  _{s+m} & \beta  _{s+m+1} & \ldots  &  \beta_{s+m+k-3} &  \beta _{s+m+k-2}  \\
 \gamma  _{s} &  \gamma _{s+1} &  \gamma _{s+2} & \ldots  & \gamma _{s+k-2} &  \gamma  _{s+k-1}  \\
 \gamma  _{s+1} &  \gamma _{s+2} &  \gamma _{s+3} & \ldots  & \gamma _{s+k-1} &  \gamma  _{s+k}  \\
\vdots & \vdots & \vdots    &  \vdots & \vdots  &  \vdots \\
\gamma  _{s+m+l-1} & \gamma  _{s+m+l} &  \gamma  _{s+m+l+1} & \ldots  & \gamma  _{s+m+l+k-3} &  \gamma  _{s+m+l+k-2}  \\
\hline 
  \alpha _{s} & \alpha _{s+1} & \alpha _{s +2} & \ldots  &  \alpha _{s+k-2} &  \alpha _{s+k-1}  \\
   \beta  _{s+m} & \beta  _{s+m+1} & \beta  _{s+m+2} & \ldots  &  \beta_{s+m+k-2} &  \beta _{s+m+k-1} \\
   \gamma  _{s+m+l} & \gamma   _{s+m+l+1} &  \gamma  _{s+m+l+2} & \ldots  &  \gamma _{s+m+l+k-2} &  \gamma  _{s+m+l+k-1} 
 \end{array}  \right). $$

   Let j be a rational integer such that $1\leq j\leq k-1.$\\
  
   $Set \;(t,\eta,\xi  )= (\sum_{i\geq 1}\alpha _{i}T^{-i},\sum_{i\geq 1}\beta  _{i}T^{-i},\sum_{i\geq 1}\gamma  _{i}T^{-i}) \in \mathbb{P}^3$\\
   
  Recall that we  denote by 
 $  D_{j}^{ \left[s\atop{ s+m\atop s+m+l}\right]\times (k-j+1)}(t,\eta ,\xi )  $
  any  $(3s+2m+l)\times (k-j+1) $   matrix, 
such that  after a rearrangement of the rows, if necessary,  we can  obtain  the following triple persymmetric matrix 
$ \left[{D_{s \times (k-j+1)}^{j}(t)\over {D_{(s+m )\times (k-j+1)}^{j}(\eta )\over D_{(s+m+l)\times (k-j+1)}^{j}(\xi )}}\right] $

   $$   \left ( \begin{array} {ccccc|c}
\alpha _{j} & \alpha _{j+1} & \alpha _{j+2} &  \ldots & \alpha _{k-1}  &  \alpha _{k} \\
\alpha _{j+1 } & \alpha _{j+2} & \alpha _{j+3}&  \ldots  &  \alpha _{k} &  \alpha _{k+1} \\
\vdots & \vdots & \vdots    &  \vdots & \vdots  &  \vdots \\
\alpha _{j+s-1} & \alpha _{j+s} & \alpha _{j+s +1} & \ldots  &  \alpha _{k+s-2} &  \alpha _{k+s-1}  \\
\hline 
\beta  _{j} & \beta  _{j+1} & \beta  _{j+2} & \ldots  &  \beta_{k-1} &  \beta _{k}  \\
\beta  _{j+1} & \beta  _{j+2} & \beta  _{j+3} & \ldots  &  \beta_{k} &  \beta _{k+1}  \\
\vdots & \vdots & \vdots    &  \vdots & \vdots  &  \vdots \\
\beta  _{m+j} & \beta  _{m+j+1} & \beta  _{m+j+2} & \ldots  &  \beta_{k+m-1} &  \beta _{k+m}  \\
\vdots & \vdots & \vdots    &  \vdots & \vdots  &  \vdots \\
\beta  _{s+m+j-1} & \beta  _{s+m+j} & \beta  _{s+m+j+1} & \ldots  &  \beta_{s+m+k-2} &  \beta _{s+m+k-1}  \\
\hline
\gamma  _{j} & \gamma   _{j+1} &  \gamma _{j+2} & \ldots  & \gamma  _{k-1} &  \gamma  _{k}  \\
\gamma  _{j+1} & \gamma  _{j+2} & \gamma  _{j+3} & \ldots  & \gamma  _{k} &  \gamma  _{k+1}  \\
\vdots & \vdots & \vdots    &  \vdots & \vdots  &  \vdots \\
 \gamma  _{m+l+j} &  \gamma _{m+l+j+1} &  \gamma _{m+l+j+2} & \ldots  & \gamma _{k+m+l-1} &  \gamma  _{k+m+l}  \\
\vdots & \vdots & \vdots    &  \vdots & \vdots  &  \vdots \\
\gamma  _{s+m+l+j-1} & \gamma  _{s+m+l+j} &  \gamma  _{s+m+l+j+1} & \ldots  & \gamma  _{s+m+l+k-2} &  \gamma  _{s+m+l+k-1}  
 \end{array}  \right). $$ \\
  By \eqref{eq 3.23} we have  $ r( D_{2}^{ \left[s-1\atop{ s+m-1\atop s+m+l-1}\right]\times (k-1)}(t_{0},\eta_{0} ,\xi_{0} ) ) \leq 
 j_{0}, \quad   r( D^{ \left[s-1\atop{ s+m-1\atop s+m+l-1}\right]\times k}(t_{0},\eta_{0} ,\xi_{0} ) ) = j_{0} +1.$\\
 It follows that  $ r( D_{2}^{ \left[s-1\atop{ s+m-1\atop s+m+l-1}\right]\times (k-1)}(t_{0},\eta_{0} ,\xi_{0} ) ) = j_{0}$
   Let $ a_{1},a_{2},\ldots,a_{k}$ denote the  columns of \quad $  D^{ \left[s-1\atop{ s+m-1\atop s+m+l-1}\right]\times k}(t_{0},\eta_{0} ,\xi_{0} )  $  that is,\\
     
\begin{displaymath}
a_{i} = \left(\begin{array}{c}
\alpha _{i} \\
 \beta  _{i}  \\
 \gamma _{i} \\
 \alpha _{i+1 } \\
\beta  _{i+1}\\
 \gamma _{i+1} \\
\vdots \\
\alpha _{i+s-2}  \\
\beta  _{i+s-2} \\
 \gamma _{i+s-2} \\
\beta  _{i+s-1}  \\
\beta  _{i+s}  \\
\vdots \\
\beta  _{i+s+m-3}  \\
\beta  _{i+s+m-2} \\
\gamma  _{i+s-1}\\
\gamma  _{i+s}\\
\vdots \\
\gamma   _{i+s+m+l-3}  \\
\gamma  _{i+s+m+l-2} 
\end{array}\right).
\end{displaymath} \\
Since $r( D_{2}^{ \left[s-1\atop{ s+m-1\atop s+m+l-1}\right]\times (k-1)}(t_{0},\eta_{0} ,\xi_{0} ) ) = j_{0}\quad
and \quad  r( D^{ \left[s-1\atop{ s+m-1\atop s+m+l-1}\right]\times k}(t_{0},\eta_{0} ,\xi_{0} ) ) = j_{0} +1, $\\
 we  have 
$ a_{1}\notin span\left\{a_{2},a_{3},\ldots,a_{k}\right\},
    \text{therefore  }\;r( D_{2}^{ \left[s-1\atop{ s+m-1\atop s+m+l-1}\right]\times (k-2)}(t_{0},\eta_{0} ,\xi_{0} ) ) = j_{0}-1  $\\
We have then 
  \begin{equation}
\label{eq 3.24}  \left\{  \begin{array}{cc}  
 r( D_{2}^{\left[s-1\atop{ s+m-1\atop s+m+l-1} \right]\times (k-2)}(t_{0},\eta_{0} ,\xi_{0} ) ) = j_{0}-1,\quad
   r( D_{2}^{\left[s-1\atop{ s+m-1\atop s+m+l-1} \right]\times (k-1)}(t_{0},\eta_{0} ,\xi_{0} )) = j_{0}\\
    r( D^{\left[s-1\atop{ s+m-1\atop s+m+l-1} \right]\times (k-1)}(t_{0},\eta_{0} ,\xi_{0} ) ) 
=  r( D^{\left[s\atop{ s+m-1\atop s+m+l-1} \right]\times (k-1)}(t_{0},\eta_{0} ,\xi_{0} ) ) \\
   = r( D^{\left[s\atop{ s+m\atop s+m+l-1} \right]\times (k-1)}(t_{0},\eta_{0} ,\xi_{0} ) )     
 =  r( D^{\left[s\atop{ s+m\atop s+m+l} \right]\times (k-1)}(t_{0},\eta_{0} ,\xi_{0} ) ) = j_{0}.
  \end{array}\right. 
  \end{equation}         
Now we consider the matrix obtained from the matrix 
$$  \displaystyle  D^{ \left[s-1\atop{ s+m-1\atop {s+m+l-1 \atop{\overline{\alpha _{s -} \atop{\beta _{s+m -} \atop \gamma_{s+m+l -}} }}}}\right]\times k}(t_{0},\eta_{0} ,\xi_{0} ) $$\\
by delating the first column and replacing the last column by the first one :  
    $$   \left ( \begin{array} {cccc|c}
 \alpha _{2} & \alpha _{3} &  \ldots & \alpha _{k-1}  & \alpha _{1} \\
 \beta  _{2} & \beta  _{3} & \ldots  &  \beta_{k-1} &  \beta  _{1}  \\
 \gamma   _{2} &  \gamma _{3} & \ldots  & \gamma  _{k-1} & \gamma  _{1}    \\
 \alpha _{3} & \alpha _{4}&  \ldots  &  \alpha _{k} & \alpha _{2 }  \\
  \beta  _{3} & \beta  _{4} & \ldots  &  \beta_{k} &  \beta  _{2} \\
 \gamma  _{3} & \gamma  _{4} & \ldots  & \gamma  _{k} & \gamma  _{2}  \\
 \vdots & \vdots    &  \vdots & \vdots  &  \vdots \\
  \alpha _{s} & \alpha _{s +1} & \ldots  &  \alpha _{s+k-3} &  \alpha _{s-1}  \\
 \beta  _{s} & \beta  _{s+1} & \ldots  &  \beta_{s+k-3} &  \beta  _{s-1}  \\
 \gamma  _{s} &  \gamma  _{s+1} & \ldots  & \gamma  _{s+k-3} &  \gamma  _{s-1} \\
 \beta  _{s+1} & \beta  _{s+2} & \ldots  &  \beta_{s+k-2} &   \beta  _{s} \\
 \beta  _{s+2} & \beta  _{s+3} & \ldots  &  \beta_{k+s-1} & \beta  _{s+1}  \\
 \vdots & \vdots    &  \vdots & \vdots  &  \vdots \\
  \beta  _{s+m} & \beta  _{s+m+1} & \ldots  &  \beta_{s+m+k-3} &  \beta  _{s+m-1}  \\
   \gamma _{s+1} &  \gamma _{s+2} & \ldots  & \gamma _{s+k-2} &  \gamma  _{s} \\
   \gamma _{s+2} &  \gamma _{s+3} & \ldots  & \gamma _{s+k-1} &  \gamma  _{s+1} \\
  \vdots & \vdots    &  \vdots & \vdots  &  \vdots \\
 \gamma  _{s+m+l} &  \gamma  _{s+m+l+1} & \ldots  & \gamma  _{s+m+l+k-3} & \gamma  _{s+m+l-1}  \\
\hline 
 \alpha _{s+1} & \alpha _{s +2} & \ldots  &  \alpha _{s+k-2} &  \alpha _{s}  \\
    \beta  _{s+m+1} & \beta  _{s+m+2} & \ldots  &  \beta_{s+m+k-2} &     \beta  _{s+m} \\
 \gamma   _{s+m+l+1} &  \gamma  _{s+m+l+2} & \ldots  &  \gamma _{s+m+l+k-2} &   \gamma  _{s+m+l}  
 \end{array}  \right). $$ 
From \eqref{eq 3.24} we obtain by elementary rank considerations
  \begin{equation}
\label{eq 3.25}  \left\{  \begin{array}{cc}  
 r( D_{2}^{\left[s-1\atop{ s+m-1\atop s+m+l-1} \right]\times (k-2)}(t_{0},\eta_{0} ,\xi_{0} ) ) = 
   r( D_{2}^{\left[s\atop{ s+m-1\atop s+m+l-1} \right]\times (k-2)}(t_{0},\eta_{0} ,\xi_{0} )) = \\
    r( D^{\left[s\atop{ s+m\atop s+m+l-1} \right]\times (k-2)}(t_{0},\eta_{0} ,\xi_{0} ) ) 
=  r( D^{\left[s\atop{ s+m\atop s+m+l} \right]\times (k-2)}(t_{0},\eta_{0} ,\xi_{0} ) ) = j_{0}-1 \\
 r( D_{2}^{\left[s-1\atop{ s+m-1\atop s+m+l-1} \right]\times (k-1)}(t_{0},\eta_{0} ,\xi_{0} ) ) = j_{0}
   \end{array}\right. 
  \end{equation}         
 Consider now the matrix 
 $$  \displaystyle  D_{2}^{ \left[s-1\atop{ s+m-1\atop {s+m+l-1 \atop{\overline{\alpha _{s -} \atop{\beta _{s+m -} \atop \gamma_{s+m+l -}} }}}}\right]\times (k-1)}(t_{0},\eta_{0} ,\xi_{0} ) $$\\
 obtained by the matrix
  $$  \displaystyle  D^{ \left[s-1\atop{ s+m-1\atop {s+m+l-1 \atop{\overline{\alpha _{s -} \atop{\beta _{s+m -} \atop \gamma_{s+m+l -}} }}}}\right]\times k}(t_{0},\eta_{0} ,\xi_{0} ) $$\\
  by delating the first column, that is 
     $$   \left ( \begin{array} {cccc|c}
 \alpha _{2} & \alpha _{3} &  \ldots & \alpha _{k-1}  &  \alpha _{k} \\
 \beta  _{2} & \beta  _{3} & \ldots  &  \beta_{k-1} &  \beta _{k}  \\
 \gamma   _{2} &  \gamma _{3} & \ldots  & \gamma  _{k-1} &  \gamma  _{k}  \\
 \alpha _{3} & \alpha _{4}&  \ldots  &  \alpha _{k} &  \alpha _{k+1} \\
 \beta  _{3} & \beta  _{4} & \ldots  &  \beta_{k} &  \beta _{k+1}  \\
 \gamma  _{3} & \gamma  _{4} & \ldots  & \gamma  _{k} &  \gamma  _{k+1}  \\
 \vdots & \vdots    &  \vdots & \vdots  &  \vdots \\
  \alpha _{s} & \alpha _{s +1} & \ldots  &  \alpha _{s+k-3} &  \alpha _{s+k-2}  \\
 \beta  _{s} & \beta  _{s+1} & \ldots  &  \beta_{s+k-3} &  \beta _{s+k-2}  \\
 \gamma  _{s} &  \gamma  _{s+1} & \ldots  & \gamma  _{s+k-3} &  \gamma  _{s+k-2}  \\
 \beta  _{s+1} & \beta  _{s+2} & \ldots  &  \beta_{s+k-2} &  \beta _{s+k-1}  \\
 \beta  _{s+2} & \beta  _{s+3} & \ldots  &  \beta_{k+s-1} &  \beta _{k+s}  \\
 \vdots & \vdots    &  \vdots & \vdots  &  \vdots \\
 \beta  _{s+m} & \beta  _{s+m+1} & \ldots  &  \beta_{s+m+k-3} &  \beta _{s+m+k-2}  \\
   \gamma _{s+1} &  \gamma _{s+2} & \ldots  & \gamma _{s+k-2} &  \gamma  _{s+k-1}  \\
   \gamma _{s+2} &  \gamma _{s+3} & \ldots  & \gamma _{s+k-1} &  \gamma  _{s+k}  \\
\vdots & \vdots    &  \vdots & \vdots  &  \vdots \\
 \gamma  _{s+m+l} &  \gamma  _{s+m+l+1} & \ldots  & \gamma  _{s+m+l+k-3} &  \gamma  _{s+m+l+k-2}  \\
\hline 
   \alpha _{s+1} & \alpha _{s +2} & \ldots  &  \alpha _{s+k-2} &  \alpha _{s+k-1}  \\
   \hline
    \beta  _{s+m+1} & \beta  _{s+m+2} & \ldots  &  \beta_{s+m+k-2} &  \beta _{s+m+k-1} \\
    \hline
    \gamma   _{s+m+l+1} &  \gamma  _{s+m+l+2} & \ldots  &  \gamma _{s+m+l+k-2} &  \gamma  _{s+m+l+k-1} 
 \end{array}  \right). $$ 
 From \eqref{eq 3.25} we then deduce 
   \begin{equation*}
  \left\{  \begin{array}{cc}  
 r( D_{2}^{\left[s-1\atop{ s+m-1\atop s+m+l-1} \right]\times (k-1)}(t_{0},\eta_{0} ,\xi_{0} ) ) = 
   r( D_{2}^{\left[s\atop{ s+m-1\atop s+m+l-1} \right]\times (k-1)}(t_{0},\eta_{0} ,\xi_{0} )) = \\
    r( D_{2}^{\left[s\atop{ s+m\atop s+m+l-1} \right]\times (k-1)}(t_{0},\eta_{0} ,\xi_{0} ) ) 
=  r( D_{2}^{\left[s\atop{ s+m\atop s+m+l} \right]\times (k-1)}(t_{0},\eta_{0} ,\xi_{0} ) ) = j_{0} 
  \end{array}\right. 
  \end{equation*} 
  We get 
  \begin{align}
  \label{eq 3.26}
  {}^{\#}\left(\begin{array}{c | c}
           j_{0}-1 & j_{0} \\
           \hline
           j_{0}-1 & j_{0}  \\
           \hline
            j_{0}-1 & j_{0}  \\
            \hline 
            j_{0}-1 & j_{0}
           \end{array} \right)_{\mathbb{P}_{1}/\mathbb{P}_{k+s -1}\times
           \mathbb{P}_{1}/\mathbb{P}_{k+s+m-1}\times  \mathbb{P}_{1}/\mathbb{P}_{k+s+m+l-1} }^{{\alpha \over {\beta \over \gamma }}} > 0 
     \end{align}

  We have now proved that 
  \begin{align*}
& {}^{\#}\left(\begin{array}{c | c}
           j_{0} & j_{0}+1 \\
           \hline
           j_{0} & j_{0} +1 \\
           \hline
            j_{0} & j_{0} +1 \\
            \hline 
            j_{0} & j_{0}+1
           \end{array} \right)_{\mathbb{P}/\mathbb{P}_{k+s -1}\times
           \mathbb{P}/\mathbb{P}_{k+s+m-1}\times  \mathbb{P}/\mathbb{P}_{k+s+m+l-1} }^{{\alpha \over {\beta \over \gamma }}} > 0 \\
    &  \Longrightarrow 
 {}^{\#}\left(\begin{array}{c | c}
           j_{0}-1 & j_{0} \\
           \hline
           j_{0}-1 & j_{0}  \\
           \hline
            j_{0}-1 & j_{0}  \\
            \hline 
            j_{0}-1 & j_{0}
           \end{array} \right)_{\mathbb{P}_{1}/\mathbb{P}_{k+s -1}\times
           \mathbb{P}_{1}/\mathbb{P}_{k+s+m-1}\times  \mathbb{P}_{1}/\mathbb{P}_{k+s+m+l-1} }^{{\alpha \over {\beta \over \gamma }}} > 0 \\
    \end{align*}
  We repeat this procedure and obtain after finitely many steps 
  \begin{align*}
          {}^{\#}\left(\begin{array}{c | c}
           0 & 1 \\
           \hline
           0 & 1  \\
           \hline
            0 & 1  \\
            \hline 
            0 & 1
           \end{array} \right)_{\mathbb{P}_{j_{0}}/\mathbb{P}_{k+s -1}\times
           \mathbb{P}_{j_{0}}/\mathbb{P}_{k+s+m-1}\times  \mathbb{P}_{j_{0}}/\mathbb{P}_{k+s+m+l-1} }^{{\alpha \over {\beta \over \gamma }}} > 0     
\end{align*}
  From 
 $$  D_{j_{0}+1}^{ \left[s-1\atop{ s+m-1\atop {s+m+l-1 \atop{\overline{\alpha _{s -} \atop{\beta _{s+m -} \atop \gamma_{s+m+l -}} }}}}\right]\times (k-j_{0})}(t_{0},\eta_{0} ,\xi_{0} ) =
      \left ( \begin{array} {cccc|c}
 \alpha _{j_{0}+1} & \alpha _{j_{0}+2} &  \ldots & \alpha _{k-1}  &  \alpha _{k} \\
  \beta  _{j_{0}+1} & \beta  _{j_{0}+2} & \ldots  &  \beta_{k-1} &  \beta _{k}  \\
  \gamma   _{j_{0}+1} &  \gamma _{j_{0}+2} & \ldots  & \gamma  _{k-1} &  \gamma  _{k}  \\
  \alpha _{j_{0}+2} & \alpha _{j_{0}+3}&  \ldots  &  \alpha _{k} &  \alpha _{k+1} \\
  \beta  _{j_{0}+2} & \beta  _{j_{0}+3} & \ldots  &  \beta_{k} &  \beta _{k+1}  \\
  \gamma  _{j_{0}+2} & \gamma  _{j_{0}+3} & \ldots  & \gamma  _{k} &  \gamma  _{k+1}  \\
  \vdots & \vdots    &  \vdots & \vdots  &  \vdots \\
   \alpha _{s+j_{0}-1} & \alpha _{s +j_{0}} & \ldots  &  \alpha _{s+k-3} &  \alpha _{s+k-2}  \\
   \beta  _{s+j_{0}-1} & \beta  _{s+j_{0}} & \ldots  &  \beta_{s+k-3} &  \beta _{s+k-2}  \\
  \gamma  _{s+j_{0}-1} &  \gamma  _{s+j_{0}} & \ldots  & \gamma  _{s+k-3} &  \gamma  _{s+k-2}  \\
  \beta  _{s+j_{0}} & \beta  _{s+j_{0}+1} & \ldots  &  \beta_{s+k-2} &  \beta _{s+k-1}  \\
  \beta  _{s+j_{0}} & \beta  _{s+j_{0}+1} & \ldots  &  \beta_{k+s-1} &  \beta _{k+s}  \\
  \vdots & \vdots    &  \vdots & \vdots  &  \vdots \\
  \beta  _{s+j_{0}+m-1} & \beta  _{s+j_{0}+m} & \ldots  &  \beta_{s+m+k-3} &  \beta _{s+m+k-2}  \\
  \gamma _{s+j_{0}} &  \gamma _{s+j_{0}+1} & \ldots  & \gamma _{s+k-2} &  \gamma  _{s+k-1}  \\
     \gamma _{s+j_{0}+1} &  \gamma _{s+j_{0}+2} & \ldots  & \gamma _{s+k-1} &  \gamma  _{s+k}  \\
   \vdots & \vdots    &  \vdots & \vdots  &  \vdots \\
 \gamma  _{s+j_{0}+m+l-1} &  \gamma  _{s+j_{0}+m+l} & \ldots  & \gamma  _{s+m+l+k-3} &  \gamma  _{s+m+l+k-2}  \\
\hline 
   \alpha _{s+j_{0}} & \alpha _{s +j_{0}+1} & \ldots  &  \alpha _{s+k-2} &  \alpha _{s+k-1}  \\
   \hline
    \beta  _{s+j_{0}+m} & \beta  _{s+j_{0}+m+1} & \ldots  &  \beta_{s+m+k-2} &  \beta _{s+m+k-1} \\
    \hline
    \gamma   _{s+j_{0}+m+l} &  \gamma  _{s+j_{0}+m+l+1} & \ldots  &  \gamma _{s+m+l+k-2} &  \gamma  _{s+m+l+k-1} 
 \end{array}  \right), $$ 
 we obviously get 
  \begin{align*}
          {}^{\#}\left(\begin{array}{c | c}
           0 & 1 \\
           \hline
           0 & 1  \\
           \hline
            0 & 1  \\
            \hline 
            0 & 1
           \end{array} \right)_{\mathbb{P}_{j_{0}}/\mathbb{P}_{k+s -1}\times
           \mathbb{P}_{j_{0}}/\mathbb{P}_{k+s+m-1}\times  \mathbb{P}_{j_{0}}/\mathbb{P}_{k+s+m+l-1} }^{{\alpha \over {\beta \over \gamma }}} = 0,     
\end{align*}
which obviously contradicts \eqref{eq 3.26}
   \end{proof}
  \section{\textbf{Study of the remainder $  \Delta _{i}^{\left[s\atop{ s+m\atop s+m+l} \right]\times k} $  in the recurrent formula  }}
 \label{sec 4}
  \subsection{Notation}
  \label{subsec 1}
We define $  \Delta _{i}^{\left[s\atop{ s+m\atop s+m+l} \right]\times k} $ to be the sum\\[0.01 cm]
\begin{align*}
 \sigma _{i,i,i,i}^{\left[\alpha \atop{ \beta \atop \gamma } \right]\times k}
-7\cdot\sigma _{i-1,i-1,i-1,i-1}^{\left[\alpha \atop{ \beta \atop \gamma } \right]\times k}
+14\cdot \sigma _{i-2,i-2,i-2,i-2}^{\left[\alpha \atop{ \beta \atop \gamma } \right]\times k}
-8\cdot \sigma _{i-3,i-3,i-3,i-3}^{\left[\alpha \atop{ \beta \atop \gamma } \right]\times k}
 \end{align*}
 Recall that 
$$  \sigma _{i,i,i,i}^{\left[\alpha \atop{ \beta \atop \gamma } \right]\times k}$$
 denotes  the cardinality of the following set 
 $$\begin{array}{l}\Big\{ (t,\eta,\xi  ) \in \mathbb{P}/\mathbb{P}_{k+s -1}\times
           \mathbb{P}/\mathbb{P}_{k+s+m-1}\times \mathbb{P}/\mathbb{P}_{k+s+m+l-1}
\mid r(  D^{\left[s-1\atop{ s+m-1\atop s+m+l-1} \right]\times k}(t,\eta,\xi ) ) = i \\
 r(  D^{\left[s\atop{ s+m-1\atop s+m+l-1} \right]\times k}(t,\eta,\xi ) )= i,\quad  
 r(  D^{\left[s\atop{ s+m\atop s+m+l-1} \right]\times k}(t,\eta,\xi ) ) = i,  \quad
  r(  D^{\left[s\atop{ s+m\atop s+m+l} \right]\times k}(t,\eta,\xi ) ) = i \Big\}
    \end{array}$$

   \subsection{Introduction}
  \label{subsec 2}

Adapting the methods used in Section 8 of [2] and using the rank formulas established in Section 3, we deduce by elementary rank considerations the following formula
 for $3\leqslant i \leqslant 3s+2m+l ,\quad k\geqslant i+1$\\
$$   \Delta _{i}^{\left[s\atop{ s+m\atop s+m+l} \right]\times k} =  \sum_{j=i-3}^{i+1} a_{j}\cdot  \Gamma_{j}^{\left[s-1\atop{ s+m -1\atop s+m+l-1} \right]\times j} $$
where the $a_{j} \in \textbf{Z} $ are explicitly  determined.\\

 \subsection{\textbf{A rank formula  for $ \sigma _{i,i,i,i}^{\left[s-1\atop{ s+m-1\atop {s+m+l-1 \atop{\overline{\alpha_{s -} \atop{\beta_{s+m -} \atop \gamma_{s+m+l -}} }}}}\right]\times k} $}}
  \label{subsec 3}
 \begin{lem}
   \label{lem 4.1}
   \label{eq 4.1}
     For $ 1\leq i\leq 2s+m+l-4,\quad k \geq  i+1  $ we have
    \begin{equation}
     \sigma _{i,i,i,i}^{\left[s-1\atop{ s+m-1\atop {s+m+l-1 \atop{\overline{\alpha_{s -} \atop{\beta_{s+m -} \atop \gamma_{s+m+l -}} }}}}\right]\times k} 
    =  \sigma _{i,i,i,i}^{\left[s-1\atop{ s+m-1\atop {s+m+l-1 \atop{\overline{\alpha_{s -} \atop{\beta_{s+m -} \atop \gamma_{s+m+l -}} }}}}\right]\times (i+1)} 
    \end{equation}
\end{lem} 
\begin{proof}
The formula \eqref{eq 4.1} is obvious for k=i+1.\\
We consider the following partition of the matrix
$$D^{\left[s-1\atop{ s+m-1\atop {s+m+l-1 \atop{\overline{\alpha_{s -} \atop{\beta_{s+m -} \atop \gamma_{s+m+l -}} }}}}\right]\times k}(t,\eta ,\xi )$$

   $$   \left ( \begin{array} {cccc|c}
\alpha _{1} & \alpha _{2}  &  \ldots & \alpha _{k-1}  &  \alpha _{k} \\
\alpha _{2 } & \alpha _{3} &  \ldots  &  \alpha _{k} &  \alpha _{k+1} \\
\vdots & \vdots & \vdots    & \vdots  &  \vdots \\
\alpha _{s-1} & \alpha _{s} & \ldots  &  \alpha _{s+k-3} &  \alpha _{s+k-2}  \\
\beta  _{1} & \beta  _{2}  & \ldots  &  \beta_{k-1} &  \beta _{k}  \\
\beta  _{2} & \beta  _{3}  & \ldots  &  \beta_{k} &  \beta _{k+1}  \\
\vdots & \vdots    &  \vdots & \vdots  &  \vdots \\
\beta  _{m+1} & \beta  _{m+2}  & \ldots  &  \beta_{k+m-1} &  \beta _{k+m}  \\
\vdots & \vdots    &  \vdots & \vdots  &  \vdots \\
\beta  _{s+m-1} & \beta  _{s+m}  & \ldots  &  \beta_{s+m+k-3} &  \beta _{s+m+k-2}  \\
\gamma  _{1} & \gamma   _{2}  & \ldots  & \gamma  _{k-1} &  \gamma  _{k}  \\
\gamma  _{2} & \gamma  _{3}  & \ldots  & \gamma  _{k} &  \gamma  _{k+1}  \\
\vdots & \vdots    &  \vdots & \vdots  &  \vdots \\
 \gamma  _{m+l+1} &  \gamma _{m+l+2}  & \ldots  & \gamma _{k+m+l-1} &  \gamma  _{k+m+l}  \\
\vdots & \vdots   &  \vdots & \vdots  &  \vdots \\
 \gamma  _{s+m+l-1} & \gamma  _{s+m+l}  & \ldots  & \gamma  _{s+m+l+k-3} &  \gamma  _{s+m+l+k-2}  \\
   \hline
  \alpha  _{s } & \alpha  _{s +1} & \ldots & \alpha  _{s +k-2}& \alpha  _{s +k-1}\\
  \hline
   \beta _{s+m} & \beta _{s+m+1} & \ldots & \beta _{s+m+k-2} & \beta _{s+m+k-1}\\
   \hline
  \gamma  _{s+m+l} & \gamma  _{s+m+l+1}  & \ldots  & \gamma  _{s+m+l+k-2} &  \gamma  _{s+m+l+k-1}  
\end{array}  \right). $$ 
Let $k\geq i+2$. Using Lemma \ref{lem 3.15} we obtain by elementary rank considerations \\
\begin{align}
&   \sigma _{i,i,i,i}^{\left[s-1\atop{ s+m-1\atop {s+m+l-1 \atop{\overline{\alpha_{s -} \atop{\beta_{s+m -} \atop \gamma_{s+m+l -}} }}}}\right]\times k}  \label{eq 4.2} \\
&   =      {}^{\#}\left(\begin{array}{c | c}
           i & i \\
           \hline
           i & i \\
           \hline
            i & i \\
            \hline 
            i & i
           \end{array} \right)_{\mathbb{P}/\mathbb{P}_{k+s -1}\times
           \mathbb{P}/\mathbb{P}_{k+s+m-1}\times  \mathbb{P}/\mathbb{P}_{k+s+m+l-1} }^{{\alpha \over {\beta \over \gamma }}}   
     +     {}^{\#}\left(\begin{array}{c | c}
           i-1 & i \\
           \hline
           i & i \\
           \hline
            i & i \\
            \hline 
            i & i
           \end{array} \right)_{\mathbb{P}/\mathbb{P}_{k+s -1}\times
           \mathbb{P}/\mathbb{P}_{k+s+m-1}\times  \mathbb{P}/\mathbb{P}_{k+s+m+l-1} }^{{\alpha \over {\beta \over \gamma }}}  \nonumber   \\
           &   +     {}^{\#}\left(\begin{array}{c | c}
           i-1 & i \\
           \hline
           i-1 & i \\
           \hline
            i & i \\
            \hline 
            i & i
           \end{array} \right)_{\mathbb{P}/\mathbb{P}_{k+s -1}\times
           \mathbb{P}/\mathbb{P}_{k+s+m-1}\times  \mathbb{P}/\mathbb{P}_{k+s+m+l-1} }^{{\alpha \over {\beta \over \gamma }}}   
     +     {}^{\#}\left(\begin{array}{c | c}
           i-1 & i \\
           \hline
           i-1 & i \\
           \hline
            i-1 & i \\
            \hline 
            i & i
           \end{array} \right)_{\mathbb{P}/\mathbb{P}_{k+s -1}\times
           \mathbb{P}/\mathbb{P}_{k+s+m-1}\times  \mathbb{P}/\mathbb{P}_{k+s+m+l-1} }^{{\alpha \over {\beta \over \gamma }}}  \nonumber  \\
         &   +   {}^{\#}\left(\begin{array}{c | c}
           i-1 & i \\
           \hline
           i-1 & i \\
           \hline
            i-1 & i \\
            \hline 
            i-1 & i
           \end{array} \right)_{\mathbb{P}/\mathbb{P}_{k+s -1}\times
           \mathbb{P}/\mathbb{P}_{k+s+m-1}\times  \mathbb{P}/\mathbb{P}_{k+s+m+l-1} }^{{\alpha \over {\beta \over \gamma }}}  \nonumber  \\
           &   =      {}^{\#}\left(\begin{array}{c | c}
           i & i \\
           \hline
           i & i \\
           \hline
            i & i \\
            \hline 
            i & i
           \end{array} \right)_{\mathbb{P}/\mathbb{P}_{k+s -1}\times
           \mathbb{P}/\mathbb{P}_{k+s+m-1}\times  \mathbb{P}/\mathbb{P}_{k+s+m+l-1} }^{{\alpha \over {\beta \over \gamma }}}  \nonumber 
           \end{align}
Further, by Lemmas \ref{lem 3.13}, \ref{lem 3.14} and \ref{lem 3.15} we have \\

\begin{align}
& 8\cdot \sigma _{i,i,i,i}^{\left[s-1\atop{ s+m-1\atop {s+m+l-1 \atop{\overline{\alpha_{s -} \atop{\beta_{s+m -} \atop \gamma_{s+m+l -}} }}}}\right]\times (k-1)}  \label{eq 4.3} \\
&   =      {}^{\#}\left(\begin{array}{c | c}
           i & i \\
           \hline
           i & i \\
           \hline
            i & i \\
            \hline 
            i & i
           \end{array} \right)_{\mathbb{P}/\mathbb{P}_{k+s -1}\times
           \mathbb{P}/\mathbb{P}_{k+s+m-1}\times  \mathbb{P}/\mathbb{P}_{k+s+m+l-1} }^{{\alpha \over {\beta \over \gamma }}}   
     +     {}^{\#}\left(\begin{array}{c | c}
           i & i \\
           \hline
           i & i +1\\
           \hline
            i & i+1 \\
            \hline 
            i & i+1
           \end{array} \right)_{\mathbb{P}/\mathbb{P}_{k+s -1}\times
           \mathbb{P}/\mathbb{P}_{k+s+m-1}\times  \mathbb{P}/\mathbb{P}_{k+s+m+l-1} }^{{\alpha \over {\beta \over \gamma }}}  \nonumber   \\
           &   +     {}^{\#}\left(\begin{array}{c | c}
           i & i \\
           \hline
           i & i \\
           \hline
            i & i+1 \\
            \hline 
            i & i+1
           \end{array} \right)_{\mathbb{P}/\mathbb{P}_{k+s -1}\times
           \mathbb{P}/\mathbb{P}_{k+s+m-1}\times  \mathbb{P}/\mathbb{P}_{k+s+m+l-1} }^{{\alpha \over {\beta \over \gamma }}}   
     +     {}^{\#}\left(\begin{array}{c | c}
           i & i \\
           \hline
           i & i \\
           \hline
            i & i \\
            \hline 
            i & i+1
           \end{array} \right)_{\mathbb{P}/\mathbb{P}_{k+s -1}\times
           \mathbb{P}/\mathbb{P}_{k+s+m-1}\times  \mathbb{P}/\mathbb{P}_{k+s+m+l-1} }^{{\alpha \over {\beta \over \gamma }}}  \nonumber  \\
         &   +   {}^{\#}\left(\begin{array}{c | c}
           i & i +1\\
           \hline
           i & i +1\\
           \hline
            i & i +1\\
            \hline 
            i & i+1
           \end{array} \right)_{\mathbb{P}/\mathbb{P}_{k+s -1}\times
           \mathbb{P}/\mathbb{P}_{k+s+m-1}\times  \mathbb{P}/\mathbb{P}_{k+s+m+l-1} }^{{\alpha \over {\beta \over \gamma }}}  \nonumber  \\
           &   =   8\cdot{}^{\#}\left(\begin{array}{c | c}
           i & i \\
           \hline
           i & i \\
           \hline
            i & i \\
            \hline 
            i & i
           \end{array} \right)_{\mathbb{P}/\mathbb{P}_{k+s -1}\times
           \mathbb{P}/\mathbb{P}_{k+s+m-1}\times  \mathbb{P}/\mathbb{P}_{k+s+m+l-1} }^{{\alpha \over {\beta \over \gamma }}}  \nonumber 
           \end{align}
By \eqref{eq 4.2},\eqref{eq 4.3} we deduce \\
\begin{equation}
\label{eq 4.4}
\sigma _{i,i,i,i}^{\left[s-1\atop{ s+m-1\atop {s+m+l-1 \atop{\overline{\alpha_{s -} \atop{\beta_{s+m -} \atop \gamma_{s+m+l -}} }}}}\right]\times (k-1)}
 = \sigma _{i,i,i,i}^{\left[s-1\atop{ s+m-1\atop {s+m+l-1 \atop{\overline{\alpha_{s -} \atop{\beta_{s+m -} \atop \gamma_{s+m+l -}} }}}}\right]\times k} \quad \text{for all $k\geq i+2$}
\end{equation}
Hence by \eqref{eq 4.4}we obtain successively\\
\begin{equation*}
\sigma _{i,i,i,i}^{\left[s-1\atop{ s+m-1\atop {s+m+l-1 \atop{\overline{\alpha_{s -} \atop{\beta_{s+m -} \atop \gamma_{s+m+l -}} }}}}\right]\times k}
= \sigma _{i,i,i,i}^{\left[s-1\atop{ s+m-1\atop {s+m+l-1 \atop{\overline{\alpha_{s -} \atop{\beta_{s+m -} \atop \gamma_{s+m+l -}} }}}}\right]\times (k-1)}
=\ldots = \sigma _{i,i,i,i}^{\left[s-1\atop{ s+m-1\atop {s+m+l-1 \atop{\overline{\alpha_{s -} \atop{\beta_{s+m -} \atop \gamma_{s+m+l -}} }}}}\right]\times (i+1)}.
\end{equation*}
\end{proof}
 \begin{lem}
   \label{lem 4.2}
   For all i such that $1\leq i\leq 3s+2m+l-4$ we have\\
     \begin{equation}
  \label{eq 4.5}
\sigma _{i,i,i,i}^{\left[s-1\atop{ s+m-1\atop {s+m+l-1 \atop{\overline{\alpha_{s -} \atop{\beta_{s+m -} \atop \gamma_{s+m+l -}} }}}}\right]\times (i+1)}  =
8\cdot \Gamma_{i}^{\left[s-1\atop{ s-1+m\atop s-1+m+l} \right]\times i}- \Gamma_{i+1}^{\left[s-1\atop{ s-1+m\atop s-1+m+l} \right]\times (i+1)}
   \end{equation}
\end{lem}
\begin{proof}
We consider the following partition of the matrix
$$D^{\left[s-1\atop{ s+m-1\atop {s+m+l-1 \atop{\overline{\alpha_{s -} \atop{\beta_{s+m -} \atop \gamma_{s+m+l -}} }}}}\right]\times (i+1)}(t,\eta ,\xi )$$

   $$   \left ( \begin{array} {cccc|c}
\alpha _{1} & \alpha _{2}  &  \ldots & \alpha _{i}  &  \alpha _{i+1} \\
\alpha _{2 } & \alpha _{3} &  \ldots  &  \alpha _{i+1} &  \alpha _{i+2} \\
\vdots & \vdots & \vdots    & \vdots  &  \vdots \\
\alpha _{s-1} & \alpha _{s} & \ldots  &  \alpha _{s+i-2} &  \alpha _{s+i-1}  \\
\beta  _{1} & \beta  _{2}  & \ldots  &  \beta_{i} &  \beta _{i+1}  \\
\beta  _{2} & \beta  _{3}  & \ldots  &  \beta_{i+1} &  \beta _{i+2}  \\
\vdots & \vdots    &  \vdots & \vdots  &  \vdots \\
\beta  _{m+1} & \beta  _{m+2}  & \ldots  &  \beta_{i+m} &  \beta _{i+m+1}  \\
\vdots & \vdots    &  \vdots & \vdots  &  \vdots \\
\beta  _{s+m-1} & \beta  _{s+m}  & \ldots  &  \beta_{s+m+i-2} &  \beta _{s+m+i-1}  \\
\gamma  _{1} & \gamma   _{2}  & \ldots  & \gamma  _{i} &  \gamma  _{i+1}  \\
\gamma  _{2} & \gamma  _{3}  & \ldots  & \gamma  _{i+1} &  \gamma  _{i+2}  \\
\vdots & \vdots    &  \vdots & \vdots  &  \vdots \\
 \gamma  _{m+l+1} &  \gamma _{m+l+2}  & \ldots  & \gamma _{i+m+l} &  \gamma  _{i+m+l+1}  \\
\vdots & \vdots   &  \vdots & \vdots  &  \vdots \\
 \gamma  _{s+m+l-1} & \gamma  _{s+m+l}  & \ldots  & \gamma  _{s+m+l+i-2} &  \gamma  _{s+m+l+i-1}  \\
   \hline
  \alpha  _{s } & \alpha  _{s +1} & \ldots & \alpha  _{s +i-1}& \alpha  _{s +i}\\
  \hline
   \beta _{s+m} & \beta _{s+m+1} & \ldots & \beta _{s+m+i-1} & \beta _{s+m+i}\\
   \hline
  \gamma  _{s+m+l} & \gamma  _{s+m+l+1}  & \ldots  & \gamma  _{s+m+l+i-1} &  \gamma  _{s+m+l+i}  
\end{array}  \right). $$ 
By elementary rank considerations and using Lemma \ref{lem 3.15} we obtain 
\begin{align}
&   \sigma _{i,i,i,i}^{\left[s-1\atop{ s+m-1\atop {s+m+l-1 \atop{\overline{\alpha_{s -} \atop{\beta_{s+m -} \atop \gamma_{s+m+l -}} }}}}\right]\times (i+1)}  \label{eq 4.6} \\
&   =      {}^{\#}\left(\begin{array}{c | c}
           i & i \\
           \hline
           i & i \\
           \hline
            i & i \\
            \hline 
            i & i
           \end{array} \right)_{\mathbb{P}/\mathbb{P}_{i+s}\times
           \mathbb{P}/\mathbb{P}_{i+s+m}\times  \mathbb{P}/\mathbb{P}_{i+s+m+l} }^{{\alpha \over {\beta \over \gamma }}}   
     +     {}^{\#}\left(\begin{array}{c | c}
           i-1 & i \\
           \hline
           i & i \\
           \hline
            i & i \\
            \hline 
            i & i
           \end{array} \right)_{\mathbb{P}/\mathbb{P}_{i+s }\times
           \mathbb{P}/\mathbb{P}_{i+s+m}\times  \mathbb{P}/\mathbb{P}_{i+s+m+l} }^{{\alpha \over {\beta \over \gamma }}}  \nonumber   \\
           &   +     {}^{\#}\left(\begin{array}{c | c}
           i-1 & i \\
           \hline
           i-1 & i \\
           \hline
            i & i \\
            \hline 
            i & i
           \end{array} \right)_{\mathbb{P}/\mathbb{P}_{i+s }\times
           \mathbb{P}/\mathbb{P}_{i+s+m}\times  \mathbb{P}/\mathbb{P}_{i+s+m+l} }^{{\alpha \over {\beta \over \gamma }}}   
     +     {}^{\#}\left(\begin{array}{c | c}
           i-1 & i \\
           \hline
           i-1 & i \\
           \hline
            i-1 & i \\
            \hline 
            i & i
           \end{array} \right)_{\mathbb{P}/\mathbb{P}_{i+s }\times
           \mathbb{P}/\mathbb{P}_{i+s+m}\times  \mathbb{P}/\mathbb{P}_{i+s+m+l} }^{{\alpha \over {\beta \over \gamma }}}  \nonumber  \\
         &   +   {}^{\#}\left(\begin{array}{c | c}
           i-1 & i \\
           \hline
           i-1 & i \\
           \hline
            i-1 & i \\
            \hline 
            i-1 & i
           \end{array} \right)_{\mathbb{P}/\mathbb{P}_{i+s }\times
           \mathbb{P}/\mathbb{P}_{i+s+m}\times  \mathbb{P}/\mathbb{P}_{i+s+m+l} }^{{\alpha \over {\beta \over \gamma }}}  \nonumber  \\
           &   =      {}^{\#}\left(\begin{array}{c | c}
           i & i \\
           \hline
           i & i \\
           \hline
            i & i \\
            \hline 
            i & i
           \end{array} \right)_{\mathbb{P}/\mathbb{P}_{i+s }\times
           \mathbb{P}/\mathbb{P}_{i+s+m}\times  \mathbb{P}/\mathbb{P}_{i+s+m+l} }^{{\alpha \over {\beta \over \gamma }}}  \nonumber 
           \end{align}
Further by  by Lemmas \ref{lem 3.13}, \ref{lem 3.14} and \ref{lem 3.15} we have 
\begin{align}
&  8\cdot \sigma _{i,i,i,i}^{\left[s-1\atop{ s+m-1\atop {s+m+l-1 \atop{\overline{\alpha_{s -} \atop{\beta_{s+m -} \atop \gamma_{s+m+l -}} }}}}\right]\times i} 
 =   8^2\cdot \Gamma_{i}^{\left[s-1\atop{ s-1+m\atop s-1+m+l} \right]\times i} \label{eq 4.7} \\
&   =      {}^{\#}\left(\begin{array}{c | c}
           i & i \\
           \hline
           i & i \\
           \hline
            i & i \\
            \hline 
            i & i
           \end{array} \right)_{\mathbb{P}/\mathbb{P}_{i+s}\times
           \mathbb{P}/\mathbb{P}_{i+s+m}\times  \mathbb{P}/\mathbb{P}_{i+s+m+l} }^{{\alpha \over {\beta \over \gamma }}}   
     +     {}^{\#}\left(\begin{array}{c | c}
           i & i \\
           \hline
           i & i +1\\
           \hline
            i & i +1\\
            \hline 
            i & i+1
           \end{array} \right)_{\mathbb{P}/\mathbb{P}_{i+s }\times
           \mathbb{P}/\mathbb{P}_{i+s+m}\times  \mathbb{P}/\mathbb{P}_{i+s+m+l} }^{{\alpha \over {\beta \over \gamma }}}  \nonumber   \\
           &   +     {}^{\#}\left(\begin{array}{c | c}
           i & i \\
           \hline
           i & i \\
           \hline
            i & i +1\\
            \hline 
            i & i+1
           \end{array} \right)_{\mathbb{P}/\mathbb{P}_{i+s }\times
           \mathbb{P}/\mathbb{P}_{i+s+m}\times  \mathbb{P}/\mathbb{P}_{i+s+m+l} }^{{\alpha \over {\beta \over \gamma }}}   
     +     {}^{\#}\left(\begin{array}{c | c}
           i & i \\
           \hline
           i & i \\
           \hline
            i & i \\
            \hline 
            i & i+1
           \end{array} \right)_{\mathbb{P}/\mathbb{P}_{i+s }\times
           \mathbb{P}/\mathbb{P}_{i+s+m}\times  \mathbb{P}/\mathbb{P}_{i+s+m+l} }^{{\alpha \over {\beta \over \gamma }}}  \nonumber  \\
         &   +   {}^{\#}\left(\begin{array}{c | c}
           i & i+1 \\
           \hline
           i & i +1\\
           \hline
            i & i+1 \\
            \hline 
            i & i+1
           \end{array} \right)_{\mathbb{P}/\mathbb{P}_{i+s }\times
           \mathbb{P}/\mathbb{P}_{i+s+m}\times  \mathbb{P}/\mathbb{P}_{i+s+m+l} }^{{\alpha \over {\beta \over \gamma }}}  \nonumber  \\
           &   =     8\cdot {}^{\#}\left(\begin{array}{c | c}
           i & i \\
           \hline
           i & i \\
           \hline
            i & i \\
            \hline 
            i & i
           \end{array} \right)_{\mathbb{P}/\mathbb{P}_{i+s }\times
           \mathbb{P}/\mathbb{P}_{i+s+m}\times  \mathbb{P}/\mathbb{P}_{i+s+m+l} }^{{\alpha \over {\beta \over \gamma }}} 
           +   8\cdot \Gamma_{i+1}^{\left[s-1\atop{ s-1+m\atop s-1+m+l} \right]\times (i+1)} \nonumber 
           \end{align}
Combining \eqref{eq 4.6}, \eqref{eq 4.7}we deduce \\
\begin{equation*}
 8^2\cdot \Gamma_{i}^{\left[s-1\atop{ s-1+m\atop s-1+m+l} \right]\times i} = 
8\cdot \sigma _{i,i,i,i}^{\left[s-1\atop{ s+m-1\atop {s+m+l-1 \atop{\overline{\alpha_{s -} \atop{\beta_{s+m -} \atop \gamma_{s+m+l -}} }}}}\right]\times (i+1)} 
+  8\cdot \Gamma_{i+1}^{\left[s-1\atop{ s-1+m\atop s-1+m+l} \right]\times (i+1)}
\end{equation*}
and \eqref{eq 4.5} is established.
\end{proof}
\begin{lem}
\label{lem 4.3}Let $ s\geq 2 $ and $ m\geq 0,\;l\geq 0 $ we have in the following two cases : \\

\underline {The case $1 \leq  k \leq 3s+2m+l-3 $}\\
\begin{equation}
\label{eq 4.8}
 \sigma _{i,i,i,i}^{\left[s-1\atop{ s+m-1\atop {s+m+l-1 \atop{\overline{\alpha_{s -} \atop{\beta_{s+m -} \atop \gamma_{s+m+l -}} }}}}\right]\times k} =
  \begin{cases}
1 & \text{if  } i = 0,\quad k\geq 1, \\
 8\cdot \Gamma_{i}^{\left[s-1\atop{ s-1+m\atop s-1+m+l} \right]\times i}- \Gamma_{i+1}^{\left[s-1\atop{ s-1+m\atop s-1+m+l} \right]\times (i+1)}
      & \text{if   } 1 \leq i \leq k-1, \\
   8 \cdot \Gamma_{k}^{\left[s-1\atop{ s-1+m\atop s-1+m+l} \right]\times k}     & \text{if   } i = k.
\end{cases}
\end{equation}

\underline {The case $ k\geq  3s+2m+l-3 $ }\\
\begin{equation}
\label{eq 4.9}
  \sigma _{i,i,i,i}^{\left[s-1\atop{ s+m-1\atop {s+m+l-1 \atop{\overline{\alpha_{s -} \atop{\beta_{s+m -} \atop \gamma_{s+m+l -}} }}}}\right]\times k} =
 \begin{cases}
1 & \text{if  } i = 0,          \\
  8\cdot \Gamma_{i}^{\left[s-1\atop{ s-1+m\atop s-1+m+l} \right]\times i}- \Gamma_{i+1}^{\left[s-1\atop{ s-1+m\atop s-1+m+l} \right]\times (i+1)}  
      & \text{if   } 1 \leq i \leq 3s+2m+l-4,  \\
  8 \cdot \Gamma_{3s+2m+l-3}^{\left[s-1\atop{ s-1+m\atop s-1+m+l} \right]\times (3s+2m+l-3)}       & \text{if   } i = 3s+2m+l-3. 
\end{cases}
\end{equation}
\end{lem}
\begin{proof}
Combining \eqref{eq 4.1} and \eqref{eq 4.5} and recalling the definition of 
$ \sigma _{i,i,i,i}^{\left[s-1\atop{ s+m-1\atop {s+m+l-1 \atop{\overline{\alpha_{s -} \atop{\beta_{s+m -} \atop \gamma_{s+m+l -}} }}}}\right]\times k} $
we deduce easily Lemma \ref{lem 4.3}
\end{proof}
 \subsection{\textbf{A rank formula  for $ \Delta _{i}^{\left[s\atop{ s+m\atop s+m+l} \right]\times k} $}}
  \label{subsec 4}
\begin{lem}
\label{lem 4.4}
 The remainder $\Delta _{i}^{\left[s\atop{ s+m\atop s+m+l} \right]\times k}$ in the recurrent formula is equal to 
   \begin{equation}
\label{eq 4.10} 
  \begin{cases}
1 & \text{if  } i = 0,\; k \geq 1,         \\
   8 \cdot\Gamma_{1}^{\left[s-1\atop{ s-1+m\atop s-1+m+l} \right]\times 1} -\Gamma_{2}^{\left[s-1\atop{ s-1+m\atop s-1+m+l} \right]\times 2} -7
      & \text{if   } i = 1,\; k\geq 2, \\
   8 \cdot\Gamma_{1}^{\left[s-1\atop{ s-1+m\atop s-1+m+l} \right]\times 1}  -7  
    & \text{if   } i = 1,\; k =1, \\
  15\cdot\Gamma_{2}^{\left[s-1\atop{ s-1+m\atop s-1+m+l} \right]\times 2}-56\cdot \Gamma_{1}^{\left[s-1\atop{ s-1+m\atop s-1+m+l} \right]\times 1}
  -  \Gamma_{3}^{\left[s-1\atop{ s-1+m\atop s-1+m+l} \right]\times 3}+14
    & \text{if   } i = 2,\; k\geq 3, \\
      15\cdot\Gamma_{2}^{\left[s-1\atop{ s-1+m\atop s-1+m+l} \right]\times 2}-56\cdot \Gamma_{1}^{\left[s-1\atop{ s-1+m\atop s-1+m+l} \right]\times 1}
 +14 & \text{if   } i = 2,\; k = 2, \\
     15\cdot\Gamma_{3}^{\left[s-1\atop{ s-1+m\atop s-1+m+l} \right]\times 3}-80\cdot \Gamma_{2}^{\left[s-1\atop{ s-1+m\atop s-1+m+l} \right]\times 2}
 +  112 \cdot\Gamma_{1}^{\left[s-1\atop{ s-1+m\atop s-1+m+l} \right]\times 1}\\
   - \Gamma_{4}^{\left[s-1\atop{ s-1+m\atop s-1+m+l} \right]\times 4} -8
  & \text{if   } i = 3,\; k \geq 4, \\
      15\cdot\Gamma_{3}^{\left[s-1\atop{ s-1+m\atop s-1+m+l} \right]\times 3}-80\cdot \Gamma_{2}^{\left[s-1\atop{ s-1+m\atop s-1+m+l} \right]\times 2}
 +  112 \cdot\Gamma_{1}^{\left[s-1\atop{ s-1+m\atop s-1+m+l} \right]\times 1} -8
    & \text{if   } i = 3,\; k = 3, \\
 15\cdot\Gamma_{i}^{\left[s-1\atop{ s-1+m\atop s-1+m+l} \right]\times i}-80\cdot \Gamma_{i-1}^{\left[s-1\atop{ s-1+m\atop s-1+m+l} \right]\times (i-1)}
 +  120 \cdot\Gamma_{i-2}^{\left[s-1\atop{ s-1+m\atop s-1+m+l} \right]\times (i-2)} \\
  - 64\cdot\Gamma_{i-3}^{\left[s-1\atop{ s-1+m\atop s-1+m+l} \right]\times (i-3)}  -   \Gamma_{i+1}^{\left[s-1\atop{ s-1+m\atop s-1+m+l} \right]\times (i+1)}
     & \text{if } 4\leq i\leq 3s+2m+l-4,  k\geq i+1, \\
     15\cdot\Gamma_{i}^{\left[s-1\atop{ s-1+m\atop s-1+m+l} \right]\times i}-80\cdot \Gamma_{i-1}^{\left[s-1\atop{ s-1+m\atop s-1+m+l} \right]\times (i-1)}
 +  120 \cdot\Gamma_{i-2}^{\left[s-1\atop{ s-1+m\atop s-1+m+l} \right]\times (i-2)} \\
  - 64\cdot\Gamma_{i-3}^{\left[s-1\atop{ s-1+m\atop s-1+m+l} \right]\times (i-3)}    
      & \text{if } 4\leq i\leq 3s+2m+l-4,  k= i, \\ 
         15\cdot\Gamma_{3s+2m+l-3}^{\left[s-1\atop{ s-1+m\atop s-1+m+l} \right]\times (3s+2m+l-3)}-80\cdot \Gamma_{3s+2m+l-4}^{\left[s-1\atop{ s-1+m\atop s-1+m+l} \right]\times (3s+2m+l-4)}\\
          +  120 \cdot\Gamma_{3s+2m+l-5}^{\left[s-1\atop{ s-1+m\atop s-1+m+l} \right]\times (3s+2m+l-5)} 
  - 64\cdot\Gamma_{3s+2m+l-6}^{\left[s-1\atop{ s-1+m\atop s-1+m+l} \right]\times (3s+2m+l-6)}      
       & \text{if } i = 3s+2m+l-3,  k\geq i, \\  
   -80\cdot \Gamma_{3s+2m+l-3}^{\left[s-1\atop{ s-1+m\atop s-1+m+l} \right]\times (3s+2m+l-3)}\\
          +  120 \cdot\Gamma_{3s+2m+l-4}^{\left[s-1\atop{ s-1+m\atop s-1+m+l} \right]\times (3s+2m+l-4)} 
  - 64\cdot\Gamma_{3s+2m+l-5}^{\left[s-1\atop{ s-1+m\atop s-1+m+l} \right]\times (3s+2m+l-5)}          
        & \text{if } i = 3s+2m+l-2,  k\geq i, \\ 
     120 \cdot\Gamma_{3s+2m+l-3}^{\left[s-1\atop{ s-1+m\atop s-1+m+l} \right]\times (3s+2m+l-3)} 
  - 64\cdot\Gamma_{3s+2m+l-4}^{\left[s-1\atop{ s-1+m\atop s-1+m+l} \right]\times (3s+2m+l-4)}             
         & \text{if } i = 3s+2m+l-1,  k\geq i, \\ 
      - 64\cdot\Gamma_{3s+2m+l-3}^{\left[s-1\atop{ s-1+m\atop s-1+m+l} \right]\times (3s+2m+l-3)}             
      & \text{if } i = 3s+2m+l,  k\geq i .
     \end{cases}
\end{equation}
  \end{lem} 
\begin{proof}  
  From \eqref{eq 4.8}and \eqref{eq 4.9}and recalling that $\Delta _{i}^{\left[s\atop{ s+m\atop s+m+l} \right]\times k} $
  is equal to \begin{align*}
 &  \sigma _{i,i,i,i}^{\left[s-1\atop{ s+m-1\atop {s+m+l-1 \atop{\overline{\alpha_{s -} \atop{\beta_{s+m -} \atop \gamma_{s+m+l -}} }}}}\right]\times k}
-7\cdot \sigma _{i-1,i-1,i-1,i-1}^{\left[s-1\atop{ s+m-1\atop {s+m+l-1 \atop{\overline{\alpha_{s -} \atop{\beta_{s+m -} \atop \gamma_{s+m+l -}} }}}}\right]\times k}
+14\cdot \sigma _{i-2,i-2,i-2,i-2}^{\left[s-1\atop{ s+m-1\atop {s+m+l-1 \atop{\overline{\alpha_{s -} \atop{\beta_{s+m -} \atop \gamma_{s+m+l -}} }}}}\right]\times k}
-8\cdot \sigma _{i-3,i-3,i-3,i-3}^{\left[s-1\atop{ s+m-1\atop {s+m+l-1 \atop{\overline{\alpha_{s -} \atop{\beta_{s+m -} \atop \gamma_{s+m+l -}} }}}}\right]\times k}
 \end{align*}
 we deduce easily Lemma \ref{lem 4.4}
  \end{proof}

 \begin{lem}
\label{lem  4.5} We have 
\begin{align}
&\Delta _{i}^{\left[s\atop{ s+m\atop s+m+l} \right]\times k}
 =\Delta _{i}^{\left[s\atop{ s+m\atop s+m+l} \right]\times (i+1)} & \text{for $0\leq i\leq \inf(3s+2m+l-4,k-1)$}\label{eq 4.11}\\
 &\Delta _{i}^{\left[s\atop{ s+m\atop s+m+l} \right]\times k}
 =\Delta _{i}^{\left[s\atop{ s+m\atop s+m+l} \right]\times i} 
 & \text{for $i\in \{3s+2m+l-3,3s+2m+l-2,3s+2m+l-1,3s+2m+l\}, k\geq i $}\label{eq 4.12}
\end{align}
\end{lem}
\begin{proof}
Follows immediately from Lemma \ref{lem 4.4}.
\end{proof}
 \section{\textbf{A recurrent formula for $\omega _{i}(s,m,l,k)$ }}
 \label{sec 5}
  \subsection{Notation}
  \label{subsec 1}
  \begin{defn}
  \label{defn 5.1}
   Set for $ 0\leq i\leq \inf(3s+2m+l,k) $\\
  \begin{align*}
 \omega _{i}(s,m,l,k) & = \Gamma_{i}^{\left[s\atop{ s+m\atop s+m+l} \right]\times k} -
 4\cdot\Gamma_{i-1}^{\left[s\atop{ s+m-1\atop s+m+l} \right]\times k} -
8\cdot \Gamma_{i-1}^{\left[s\atop{ s+m\atop s+m+l-1} \right]\times k} +
32\cdot\Gamma_{i-2}^{\left[s\atop{ s+m-1\atop s+m+l-1} \right]\times k}
  \end{align*}
  \end{defn}
  
   \subsection{Introduction}
  \label{subsec 2}
In this section we prove that $\omega_{i}(s,m,l,k)$ can be partly obtained from \\
$\omega_{i}(s,m,l,k)$ with $i \rightarrow  i-(s-1),\;   s \rightarrow 1,\;m \rightarrow m+s-1,\;l \rightarrow l,\; k \rightarrow k ,$ that is from the sum 
$$ \Gamma_{i-s+1}^{\left[1\atop{ s+m \atop s+m+l} \right]\times k}  -4\cdot \Gamma_{i-s}^{\left[1\atop{ s+m-1 \atop s+m+l} \right]\times k} -8\cdot \Gamma_{i-s}^{\left[1\atop{ s+m \atop s+m+l-1} \right]\times k} +32\cdot \Gamma_{i-s-1}^{\left[1\atop{ s+m -1\atop s+m+l-1} \right]\times k} . $$
These sums were studied in the first section.

 \subsection{\textbf{The recursive formula for $\omega _{i}(s,m,l,k)$} }
  \label{subsec 3}
  \begin{lem}
\label{lem  5.2} $ Let \; s \geq 4,  m\geq 0, \; l\geq 0\;  and \; k \geq 1, $ we have in the following two cases: \\
\begin{equation}
\label{eq 5.1}
 \omega _{i}(s,m,l,k) = \begin{cases} \displaystyle
 2^{i-3}\cdot\omega _{3}(s-(i-3),m+(i-3),l,k)  +   \sum_{j = 0}^{i-4}2^{j}\Delta _{i-j}^{\left[s-j \atop{ s+m\atop s+m+l} \right]\times k}  & \text{if  }  4\leq i\leq \inf(s-1,k)   \\
\displaystyle  2^{s-1}\cdot \omega _{i-(s-1)}(1,m+s-1,l,k) +  \sum_{j = 0}^{s-2}2^{j}\Delta _{i-j}^{\left[s-j \atop{ s+m\atop s+m+l} \right]\times k} & \text{if  }  s\leq i\leq \inf(3s+2m+l,k)
  \end{cases}
\end{equation}
\end{lem}
\begin{proof}From the recurrent formula \eqref{eq 2.83} we obtain the following equivalence
\begin{equation}
\label{eq 5.2}
 \omega _{i}(s,m,l,k)   = 2\cdot\omega _{i-1}(s-1,m+1,l,k)  +  \Delta _{i}^{\left[s\atop{ s+m\atop s+m+l} \right]\times k}   
\end{equation}
\underline{The first case $ 4 \leq i\leq \inf(s-1,k) $} \\

Using successively \eqref{eq 5.2}we get
 \begin{align*}
 \omega _{i}(s,m,l,k) &  = 2\cdot\omega _{i-1}(s-1,m+1,l,k)  +  \Delta _{i}^{\left[s\atop{ s+m\atop s+m+l} \right]\times k}   \\
    2\cdot \omega _{i-1}(s-1,m+1,l,k) & =2\cdot\left( 2\cdot \omega _{i-2}(s-2,m+2,l,k) + \Delta _{i -1}^{\left[s -1\atop{ s+m\atop s+m+l} \right]\times k}\right) \\
                                                       &  \vdots       \\
      2^{j}\cdot \omega _{i-j}(s-j,m+j,l,k) & = 2^{j}\cdot\left( 2\cdot \omega _{i-(j+1)}(s-(j+1),m+(j+1),l,k) +  \Delta _{i-j}^{\left[s -j\atop{ s+m\atop s+m+l} \right]\times k}\right) \\
                                                         &  \vdots       \\
       2^{i-4}\cdot \omega _{i-(i-4)}(s-(i-4),m+(i-4),l,k) & =2^{i-4}\cdot\left( 2\cdot \omega _{i-(i-3)}(s-(i-3),m+(i-3),l,k) +    \Delta _{i-(i-4)}^{\left[s -(i-4)\atop{ s+m\atop s+m+l} \right]\times k}\right)                                                              
  \end{align*}
  By summing the above equations we obtain:\\
  \begin{align*}
 & \sum_{j = 0}^{i-4}2^{j} \omega _{i-j}(s-j,m+j,l,k)  = \sum_{j = 0}^{i-4} 2^{j}\cdot\left( 2\cdot \omega _{i-(j+1)}(s-(j+1),m+(j+1),l,k) +  \Delta _{i-j}^{\left[s -j\atop{ s+m\atop s+m+l} \right]\times k}\right) \\
&\Longleftrightarrow  \omega _{i}(s,m,l,k)  = 2^{i-3}\cdot \omega _{3}(s-(i-3),m+(i-3),l,k) + \sum_{j = 0}^{i -4}2^{j} \Delta _{i-j}^{\left[s -j\atop{ s+m\atop s+m+l} \right]\times k}\\
\end{align*}
 \underline{The second  case $ s \leq i\leq \inf(3s+2m+l,k) $} \\
 \begin{align*}
 \omega _{i}(s,m,l,k) &  = 2\cdot\omega _{i-1}(s-1,m+1,l,k)  +  \Delta _{i}^{\left[s\atop{ s+m\atop s+m+l} \right]\times k}   \\
    2\cdot \omega _{i-1}(s-1,m+1,l,k) & =2\cdot\left( 2\cdot \omega _{i-2}(s-2,m+2,l,k) + \Delta _{i -1}^{\left[s -1\atop{ s+m\atop s+m+l} \right]\times k}\right) \\
      2^{2}\cdot \omega _{i-2}(s-2,m+2,l,k) & =2^{2}\cdot\left( 2\cdot \omega _{i-3}(s-3,m+3,l,k) + \Delta _{i -2}^{\left[s -2\atop{ s+m\atop s+m+l} \right]\times k}\right) \\
                                                       &  \vdots       \\
      2^{s-3}\cdot \omega _{i-(s-3)}(s-(s-3),m+(s-3),l,k) & = 2^{s-3}\cdot\left( 2\cdot \omega _{i-(s- 2)}(s-(s-2),m+(s-2),l,k) +  \Delta _{i-(s-3)}^{\left[s - (s-3) \atop{ s+m\atop s+m+l} \right]\times k}\right) \\
                                                         &  \vdots       \\
       2^{ s- 2}\cdot \omega _{i-(s-2)}(s-(s-2),m+(s-2),l,k) & =2^{s-2}\cdot\left( 2\cdot \omega _{i-(s-1)}(s-(s-1),m+(s-1),l,k) +    \Delta _{i-(s-2)}^{\left[s -(s-2)\atop{ s+m\atop s+m+l} \right]\times k}\right)                                                              
  \end{align*}
  By summing the above equations we obtain:\\
    \begin{align*}
 & \sum_{j = 0}^{s-2}2^{j} \omega _{i-j}(s-j,m+j,l,k)  = \sum_{j = 0}^{s-2} 2^{j}\cdot\left( 2\cdot \omega _{i-(j+1)}(s-(j+1),m+(j+1),l,k) +  \Delta _{i-j}^{\left[s -j\atop{ s+m\atop s+m+l} \right]\times k}\right) \\
&\Longleftrightarrow  \omega _{i}(s,m,1,k)  = 2^{s-1}\cdot \omega _{i-(s-1)}(1,m+ s-1,1,k) + \sum_{j = 0}^{s-2}2^{j} \Delta _{i-j}^{\left[s -j\atop{ s+m\atop s+m+l} \right]\times k}\\
\end{align*} 
     \end{proof}
  
 \section{\textbf{ COMPUTATION OF  $  \Gamma  _{i}^{\left[s\atop{ s+m\atop s+m+l} \right]\times k}\;  for\; i\leq s-1,\;i\leq k,\;m\geq 0,\;l\geq 0 $}}   
\label{sec 6}
  \subsection{Notation}
  \label{subsec 1}
We define $\gamma_{j}$ to be equal to the number $ \Gamma_{j}^{\left[j\atop{ j\atop j} \right]\times ( j+1)} $ of rank j matrices of the form   $\left[{A\over{B \over C}}\right], $ 
where A,B and C are $ j\times(j+1) $ persymmetric matrices.

 \subsection{Introduction}
  \label{subsec 2}
In this section we adapt the method used in Section 11 of [2]
to compute by induction on i   the number $ \Gamma_{i}^{\left[s\atop{ s+m\atop s+m+l} \right]\times k} $ of rank i( $ \leqslant $s-1)
matrices of the form  $\left[{A\over{B \over C}}\right],$ where A,B and C are persymmetric.

 \subsection{Computation of   $  \Gamma  _{i}^{\left[s\atop{ s+m\atop s+m+l} \right]\times k}\;  for\; i\leq s-1,\;i\leq k,\;m\geq 0,\;l\geq 0 $ }
  \label{subsec 3}
\begin{lem}
\label{lem 6.1}Consider the matrix\\
  $$ D^{\left[s\atop{ s+m\atop {s+m+l }}\right]\times k}(t,\eta ,\xi ) = \left ( \begin{array} {ccccc}
\alpha _{1} & \alpha _{2}  &  \ldots & \alpha _{k-1}  &  \alpha _{k} \\
\alpha _{2 } & \alpha _{3} &  \ldots  &  \alpha _{k} &  \alpha _{k+1} \\
\vdots & \vdots & \vdots    & \vdots  &  \vdots \\
\alpha _{s-1} & \alpha _{s} & \ldots  &  \alpha _{s+k-3} &  \alpha _{s+k-2}  \\
  \alpha  _{s } & \alpha  _{s +1} & \ldots & \alpha  _{s +k-2}& \alpha  _{s +k-1}\\
\beta  _{1} & \beta  _{2}  & \ldots  &  \beta_{k-1} &  \beta _{k}  \\
\beta  _{2} & \beta  _{3}  & \ldots  &  \beta_{k} &  \beta _{k+1}  \\
\vdots & \vdots    &  \vdots & \vdots  &  \vdots \\
\beta  _{m+1} & \beta  _{m+2}  & \ldots  &  \beta_{k+m-1} &  \beta _{k+m}  \\
\vdots & \vdots    &  \vdots & \vdots  &  \vdots \\
\beta  _{s+m-1} & \beta  _{s+m}  & \ldots  &  \beta_{s+m+k-3} &  \beta _{s+m+k-2}  \\
 \beta _{s+m} & \beta _{s+m+1} & \ldots & \beta _{s+m+k-2} & \beta _{s+m+k-1}\\
\gamma  _{1} & \gamma   _{2}  & \ldots  & \gamma  _{k-1} &  \gamma  _{k}  \\
\gamma  _{2} & \gamma  _{3}  & \ldots  & \gamma  _{k} &  \gamma  _{k+1}  \\
\vdots & \vdots    &  \vdots & \vdots  &  \vdots \\
 \gamma  _{m+l+1} &  \gamma _{m+l+2}  & \ldots  & \gamma _{k+m+l-1} &  \gamma  _{k+m+l}  \\
\vdots & \vdots   &  \vdots & \vdots  &  \vdots \\
  \gamma  _{s+m+l} & \gamma  _{s+m+l+1}  & \ldots  & \gamma  _{s+m+l+k-2} &  \gamma  _{s+m+l+k-1}  
\end{array}  \right). $$

Let $ g_{k,s,m,l}(t,\eta,\xi  )  $ be the quadratic exponential sum in  $ \mathbb{P}^{3} $ defined by 
  $$ (t,\eta,\xi  ) \in  \mathbb{P}^{3}\longmapsto  
  \sum_{deg Y\leq k-1}\sum_{deg Z \leq s-1}E(tYZ)\sum_{deg U \leq  s+m-1}E(\eta YU)\sum_{deg V \leq  s+m+l-1}E(\xi  YV) \in \mathbb{Z}.  $$
   We have
   \begin{align}
  \sum_{ i = 0}^{\inf(3s+2m+l,k)}  \Gamma  _{i}^{\left[s\atop{ s+m\atop s+m+l} \right]\times k} & = 2^{3k+3s+2m+l-3} \label{eq 6.1}\\
 &  \text{and}\nonumber \\
  \sum_{i = 0}^{\inf(3s+2m+l,k)}  \Gamma  _{i}^{\left[s\atop{ s+m\atop s+m+l} \right]\times k}\cdot2^{-i} & =  2^{2k+3s+2m+l-3} + 2^{3k-3} - 2^{2k-3}. \label{eq 6.2}
\end{align}
\end{lem}

\begin{proof}
The proof of \eqref{eq 6.1} is obvious.\\

By \eqref{eq 2.6} with  $ s-2\rightarrow s-1,\; s+m-2\rightarrow s+m-1,\; s+m+l-2\rightarrow s+m+l-1 $ we obtain,
 observing that  $ g_{k,s,m,l}(t,\eta,\xi)  $ is constant on cosets of $ \mathbb{P}_{k+s-1}\times \mathbb{P}_{k+s+m-1}\times \mathbb{P}_{k+s+m+l-1} $ \\
 
 \begin{align*}
&\int_{\mathbb{P}^{3}} g_{k,s,m,l}(t,\eta,\xi  )dtd\eta d\xi \\
&  =\sum_{(t,\eta,\xi  )\in \mathbb{P}/\mathbb{P}_{k+s-1}\times \mathbb{P}/\mathbb{P}_{k+s +m-1}\times \mathbb{P}/\mathbb{P}_{k+s +m+l-1}}
 2^{3s+2m+l+k- r( D^{\left[s\atop{ s+m\atop {s+m+l }}\right]\times k}(t,\eta ,\xi ))}\int_{\mathbb{P}_{k+s-1}}dt \int_{\mathbb{P}_{k+s+m-1}}d\eta \int_{\mathbb{P}_{k+s+m+l-1}}d\xi  \\
 & = \sum_{i=0}^{\inf(3s+2m+l,k)} \Gamma  _{i}^{\left[s\atop{ s+m\atop s+m+l} \right]\times k}\cdot 2^{3s+2m+l+k-i}\int_{\mathbb{P}_{k+s-1}}dt \int_{\mathbb{P}_{k+s+m-1}}d\eta \int_{\mathbb{P}_{k+s+m+l-1}}d\xi \\
 & = 2^{-2k+3}\cdot \sum_{i = 0}^{\inf(3s+2m+l,k)} \Gamma  _{i}^{\left[s\atop{ s+m\atop s+m+l} \right]\times k} \cdot2^{-i}. 
\end{align*}
On the other hand \\
 \begin{align*}
& \int_{\mathbb{P}^{3}} g_{k,s,m,l}(t,\eta,\xi  ) dt d \eta d\xi \\
 & =  Card \left\{(Y,Z,U), deg Y \leq k-1, deg Z \leq s-1,  deg U \leq s+m-1,  deg V \leq s+m+l-1 \mid Y \cdot Z = Y \cdot U = Y \cdot V = 0 \right\}\\
 & = 2^{3s+2m+l} + 2^{k} -1. 
 \end{align*}
 The above equations imply  \eqref{eq 6.2}.
\end{proof}
\begin{lem}
\label{lem 6.2}  We have  for all $ j\geq 1 $\vspace{0.1 cm}\\
 \begin{equation}
 \label{eq 6.3}
 (H_{j}) \hspace{3 cm}\Gamma_{j}^{\left[s\atop{ s+m\atop s+m+l} \right]\times k}  = \Gamma_{j}^{\left[j\atop{ j\atop j} \right]\times (j+1)}
  \overset{ \text{def}}{=} \gamma _{j}
    \quad   \text{ for   \quad   $ s\geq j+1,\;k\geq j+1,\; m\geq 0.\;l\geq 0  $}\\
 \end{equation}
\end{lem}
\begin{proof}The proof is by strong induction, that is:\\
If
\begin{itemize}
\item \text{$ (H_{1}) $ is true, and  }
\item \text{for all $ j\geq 1, \quad (H_{1})\wedge (H_{2})\wedge\ldots\wedge (H_{j})\quad implies \; (H_{j+1}) $}
\end{itemize}
then $ (H_{j}) $ is true for all $ j \geq 1. $\\
\underline {$ (H_{1})\, is \;true\; $}\\

Indeed from  \eqref{eq 2.83}  with   i = 1 and \eqref{eq 4.11}   \\
$\Gamma_{1}^{\left[s\atop{ s+m\atop s+m+l} \right]\times (k+1)} -  \Gamma_{1}^{\left[s\atop{ s+m\atop s+m+l} \right]\times k}
 = 14 +  \Delta _{1}^{\left[s\atop{ s+m\atop s+m+l} \right]\times (k+1)} -\big( 14 +  \Delta _{1}^{\left[s\atop{ s+m\atop s+m+l} \right]\times k}\big ) = 0 $
  for all $ s \geq 2,\; k\geq 2,\; m\geq  0,\;l\geq 0 $ \vspace{0.1 cm}\\
    which implies that 
  \begin{equation}
  \label{eq 6.4}
\Gamma_{1}^{\left[s\atop{ s+m\atop s+m+l} \right]\times k}   = \Gamma_{1}^{\left[s\atop{ s+m\atop s+m+l} \right]\times 2}
 \quad  \text{for all}\; s\geq 2, \; k \geq 2, \; m\geq  0,\;l\geq 0.
\end{equation}
Consider the matrix\\
   $$   \left ( \begin{array} {cc}
\alpha _{1} & \alpha _{2}   \\
\alpha _{2 } & \alpha _{3}  \\
\vdots & \vdots     \\
\alpha _{s} & \alpha _{s+1}   \\
\beta  _{1} & \beta  _{2}    \\
\beta  _{2} & \beta  _{3}   \\
\vdots & \vdots   \\
\beta  _{m+1} & \beta  _{m+2}   \\
\vdots & \vdots  \\
\beta  _{s+m} & \beta  _{s+m+1} \\
\gamma _{1}& \gamma _{2}\\
\gamma _{2}& \gamma _{3}\\
\vdots & \vdots    \\
\gamma  _{m+l+1} & \gamma  _{m+l+2}   \\
\vdots & \vdots  \\
\gamma  _{s+m+l} & \gamma  _{s+m+l+1} 
   \end{array}  \right). $$ \vspace{0.1 cm}\\
 We have  respectively  by \eqref{eq 6.1},  \eqref{eq 6.2} with k = 2 \\
 \begin{align}
  \sum_{i=0}^{2}\Gamma_{i}^{\left[s\atop{ s+m\atop s+m+l} \right]\times 2}  &  = 2^{3s+2m+l+3}, \label{eq 6.5} \\
  & \nonumber \\
    \sum_{i=0}^{2}\Gamma_{i}^{\left[s\atop{ s+m\atop s+m+l} \right]\times 2} \cdot 2^{-i} &  =  2^{3s+2m+l+1} + 2^{3} - 2.  \label{eq 6.6}
\end{align}
From \eqref{eq 6.5}  and  \eqref{eq 6.6} using   \eqref{eq 6.4}  we deduce \vspace{0.1 cm}\\
\begin{align}
& \Gamma_{1}^{\left[s\atop{ s+m\atop s+m+l} \right]\times k} = \Gamma_{1}^{\left[1\atop{ 1\atop 1} \right]\times 2} = 21
& \text{for all $ s\geq 2,\; k \geq 2,\; m \geq 0,\;l \geq 0 $} \label{eq 6.7}\\
&  \Gamma_{2}^{\left[s\atop{ s+m\atop s+m+l} \right]\times 2} =  2^{3s+2m+l+3}-22  & \text{for all $ s\geq 2,\; k = 2 ,\; m \geq 0,\;l \geq 0 $}\label{eq 6.8}
\end{align}
\underline {$ (H_{1}) \Rightarrow  (H_{2}) $}\\
 From  \eqref{eq 2.83}  with   i = 2  we obtain   \\
\begin{align*}
 \Gamma_{2}^{\left[s\atop{ s+m\atop s+m+l} \right]\times k} = 2 \cdot \Gamma_{1}^{\left[s-1\atop{ s+m\atop s+m+l} \right]\times k} 
+ 4 \cdot\Gamma_{1}^{\left[s\atop{ s+m-1\atop s+m+l} \right]\times k}+ 8 \cdot \Gamma_{1}^{\left[s\atop{ s+m\atop s+m+l-1} \right]\times k}
-56 + \Delta _{2}^{\left[s\atop{ s+m\atop s+m+l} \right]\times k}
\end{align*}

 Applying  \eqref{eq 4.11} with i=2 and $k\geq 3$ we get using \eqref{eq 6.7} with $s\geq 3,\;m\geq 0,\;l\geq 0$

\begin{align}
& \Gamma_{2}^{\left[s\atop{ s+m\atop s+m+l} \right]\times (k+1)} -  \Gamma_{2}^{\left[s\atop{ s+m\atop s+m+l} \right]\times k}
=2\cdot\big[\Gamma_{1}^{\left[s-1\atop{ s+m\atop s+m+l} \right]\times (k+1)} -  \Gamma_{1}^{\left[s-1\atop{ s+m\atop s+m+l} \right]\times k}\big] \nonumber\\
& +4\cdot\big[\Gamma_{1}^{\left[s\atop{ s+m-1\atop s+m+l} \right]\times (k+1)} -  \Gamma_{1}^{\left[s\atop{ s+m-1\atop s+m+l} \right]\times k}\big]
+8\cdot\big[\Gamma_{1}^{\left[s\atop{ s+m\atop s+m+l-1} \right]\times (k+1)} -  \Gamma_{1}^{\left[s\atop{ s+m\atop s+m+l-1} \right]\times k}\big] = 0 \nonumber\\
& \text{which implies} \nonumber\\
& \Gamma_{2}^{\left[s\atop{ s+m\atop s+m+l} \right]\times k} = \Gamma_{2}^{\left[s\atop{ s+m\atop s+m+l} \right]\times 3} \quad \text{for all  $ k\geq 3,\;s\geq 3,\;m\geq 0,\;l\geq 0$} \label{eq 6.9}
\end{align}
Consider the matrix\\
   $$   \left ( \begin{array} {ccc}
\alpha _{1} & \alpha _{2} & \alpha _{3}  \\
\alpha _{2 } & \alpha _{3} & \alpha _{4}  \\
\vdots & \vdots   & \vdots    \\
\alpha _{s} & \alpha _{s+1} &  \alpha _{s+2}   \\
\beta  _{1} & \beta  _{2}  & \beta  _{3}     \\
\beta  _{2} & \beta  _{3}  & \beta  _{4}    \\
\vdots & \vdots  &\vdots \\
\beta  _{m+1} & \beta  _{m+2}  & \beta  _{m+3}    \\
\vdots & \vdots &\vdots \\
\beta  _{s+m} & \beta  _{s+m+1} & \beta  _{s+m+2}   \\
\gamma _{1}& \gamma _{2} & \gamma _{3}\\
\gamma _{2}& \gamma _{3}& \gamma _{4}\\
\vdots & \vdots  &\vdots  \\
\gamma  _{m+l+1} & \gamma  _{m+l+2}  & \gamma  _{m+l+3}  \\
\vdots & \vdots &\vdots \\
\gamma  _{s+m+l} & \gamma  _{s+m+l+1} & \gamma  _{s+m+l+2}
   \end{array}  \right). $$ \vspace{0.1 cm}\\

 We have  respectively  by \eqref{eq 6.1},  \eqref{eq 6.2} with k = 3 \\
 \begin{align}
  \sum_{i=0}^{3}\Gamma_{i}^{\left[s\atop{ s+m\atop s+m+l} \right]\times 3}  &  = 2^{3s+2m+l+6}, \label{eq 6.10} \\
  & \nonumber \\
    \sum_{i=0}^{3}\Gamma_{i}^{\left[s\atop{ s+m\atop s+m+l} \right]\times 3} \cdot 2^{-i} &  =  2^{3s+2m+l+3} + 2^{6} - 2^{3}.  \label{eq 6.11}
\end{align}
From \eqref{eq 6.10}  and  \eqref{eq 6.11} using  \eqref{eq 6.7} and \eqref{eq 6.9}  we deduce \vspace{0.1 cm}\\
\begin{align}
& \Gamma_{2}^{\left[s\atop{ s+m\atop s+m+l} \right]\times k} = \Gamma_{2}^{\left[2\atop{ 2\atop 2} \right]\times 3} = 378 
& \text{for all $ s\geq 3,\; k \geq 3,\; m \geq 0,\;l \geq 0 $} \label{eq 6.12} \\
&  \Gamma_{3}^{\left[s\atop{ s+m\atop s+m+l} \right]\times 3} =  2^{3s+2m+l+6} - 400  & \text{for all $ s\geq 3,\; k = 3 ,\; m \geq 0,\;l \geq 0 $} \label{eq 6.13}\\
& \nonumber
\end{align}

 \underline {$ \quad (H_{1})\wedge (H_{2})\; implies \; (H_{3}) $}\\
 From  \eqref{eq 2.83}  with   i = 3  we obtain   \\
   \begin{align*}
  & \Gamma_{3}^{\left[s\atop{ s+m\atop s+m+l} \right]\times k}  \\
& = \big[2\cdot \Gamma_{2}^{\left[s -1 \atop{ s+m\atop s+m+l} \right]\times k}
 +4\cdot \Gamma_{2}^{\left[s\atop{ s+m-1\atop s+m+l} \right]\times k}  +8\cdot\Gamma_{2}^{\left[s\atop{ s+m\atop s+m+l-1} \right]\times k} \big ] 
 - \big[8\cdot \Gamma_{1}^{\left[s -1 \atop{ s+m-1\atop s+m+l} \right]\times k}
 +16\cdot \Gamma_{1}^{\left[s -1\atop{ s+m\atop s+m+l-1} \right]\times k}  +32\cdot\Gamma_{1}^{\left[s\atop{ s+m-1\atop s+m+l-1} \right]\times k}\big] \nonumber  \\
 & + 64\cdot \Gamma_{0}^{\left[s -1 \atop{ s+m-1\atop s+m+l-1} \right]\times k} + \Delta _{3}^{\left[s\atop{ s+m\atop s+m+l} \right]\times k} \nonumber \\
 & \nonumber
 \end{align*}
 
 Applying  \eqref{eq 4.11} with i=3 and $k\geq 4$ we get using \eqref{eq 6.7} and \eqref{eq 6.12} with $s\geq 4,\;m\geq 0,\;l\geq 0$

\begin{align*}
& \Gamma_{3}^{\left[s\atop{ s+m\atop s+m+l} \right]\times (k+1)} -  \Gamma_{3}^{\left[s\atop{ s+m\atop s+m+l} \right]\times k}
=2\cdot\big[\Gamma_{2}^{\left[s-1\atop{ s+m\atop s+m+l} \right]\times (k+1)} -  \Gamma_{2}^{\left[s-1\atop{ s+m\atop s+m+l} \right]\times k}\big]\\
& +4\cdot\big[\Gamma_{2}^{\left[s\atop{ s+m-1\atop s+m+l} \right]\times (k+1)} -  \Gamma_{2}^{\left[s\atop{ s+m-1\atop s+m+l} \right]\times k}\big]
+8\cdot\big[\Gamma_{2}^{\left[s\atop{ s+m\atop s+m+l-1} \right]\times (k+1)} -  \Gamma_{2}^{\left[s\atop{ s+m\atop s+m+l-1} \right]\times k}\big]  \\
& -8 \cdot\big[ \Gamma_{1}^{\left[s -1 \atop{ s+m-1\atop s+m+l} \right]\times (k+1)} - \Gamma_{1}^{\left[s -1 \atop{ s+m-1\atop s+m+l} \right]\times k}\big]
 -16 \cdot\big[ \Gamma_{1}^{\left[s -1 \atop{ s+m\atop s+m+l-1} \right]\times (k+1)} - \Gamma_{1}^{\left[s -1 \atop{ s+m\atop s+m+l-1} \right]\times k}\big]\\
& -32 \cdot\big[ \Gamma_{1}^{\left[s  \atop{ s+m-1\atop s+m+l-1} \right]\times (k+1)} - \Gamma_{1}^{\left[s \atop{ s+m-1\atop s+m+l-1} \right]\times k}\big] = 0
& \text{which implies}
\end{align*}

\begin{equation}
\label{eq 6.14}
 \Gamma_{3}^{\left[s\atop{ s+m\atop s+m+l} \right]\times k} = \Gamma_{3}^{\left[s\atop{ s+m\atop s+m+l} \right]\times 4} \quad \text{for all  $ k\geq 4,\;s\geq 4,\;m\geq 0,\;l\geq 0$}
\end{equation}
Consider the matrix\\
   $$   \left ( \begin{array} {cccc}
\alpha _{1} & \alpha _{2} & \alpha _{3}  & \alpha _{4} \\
\alpha _{2 } & \alpha _{3} & \alpha _{4}  & \alpha _{5} \\
\vdots & \vdots   & \vdots  & \vdots  \\
\alpha _{s} & \alpha _{s+1} &  \alpha _{s+2} & \alpha _{s+3}   \\
\beta  _{1} & \beta  _{2}  & \beta  _{3}   & \beta  _{4}  \\
\beta  _{2} & \beta  _{3}  & \beta  _{4}  & \beta  _{5}  \\
\vdots & \vdots  &\vdots &\vdots \\
\beta  _{m+1} & \beta  _{m+2}  & \beta  _{m+3} & \beta  _{m+4}   \\
\vdots & \vdots &\vdots &\vdots \\
\beta  _{s+m} & \beta  _{s+m+1} & \beta  _{s+m+2} & \beta  _{s+m+3}  \\
\gamma _{1}& \gamma _{2} & \gamma _{3} & \gamma _{4}\\
\gamma _{2}& \gamma _{3}& \gamma _{4} & \gamma _{5}\\
\vdots & \vdots  &\vdots &\vdots \\
\gamma  _{m+l+1} & \gamma  _{m+l+2}  & \gamma  _{m+l+3} & \gamma  _{m+l+4} \\
\vdots & \vdots &\vdots & \vdots \\
\gamma  _{s+m+l} & \gamma  _{s+m+l+1} & \gamma  _{s+m+l+2}& \gamma  _{s+m+l+3}
   \end{array}  \right). $$ \vspace{0.1 cm}\\

 We have  respectively  by \eqref{eq 6.1},  \eqref{eq 6.2} with k = 4 \\
 \begin{align}
  \sum_{i=0}^{4}\Gamma_{i}^{\left[s\atop{ s+m\atop s+m+l} \right]\times 4}  &  = 2^{3s+2m+l+9}, \label{eq 6.15} \\
  & \nonumber \\
    \sum_{i=0}^{4}\Gamma_{i}^{\left[s\atop{ s+m\atop s+m+l} \right]\times 4} \cdot 2^{-i} &  =  2^{3s+2m+l+5} + 2^{9} - 2^{5}.  \label{eq 6.16}
\end{align}
From \eqref{eq 6.15}  and  \eqref{eq 6.16} using  \eqref{eq 6.14},  \eqref{eq 6.7} \eqref{eq 6.12} we deduce \vspace{0.1 cm}\\
\begin{align}
& \Gamma_{3}^{\left[s\atop{ s+m\atop s+m+l} \right]\times k} = \Gamma_{3}^{\left[3\atop{ 3\atop 3} \right]\times 4} = 6384
& \text{for all $ s\geq 4,\; k \geq 4,\; m \geq 0,\;l \geq 0 $} \label{eq 6.17}\\
&  \Gamma_{4}^{\left[s\atop{ s+m\atop s+m+l} \right]\times 4} =  2^{3s+2m+l+9} - 6784 & \text{for all $ s\geq 4,\; k = 4 ,\; m \geq 0,\;l \geq 0 $} \label{eq 6.18}
\end{align}
 \underline{Computation of $ \omega _{3}(s-(i-3),m+(i-3),l,k) $ in the formula \eqref{eq 5.1} for  $4\leq i\leq \inf(s-1,k)$}

    \begin{align}
 & \text{ By  Definition \ref{defn 5.1} }\nonumber \\
   & \omega _{3}(s-(i-3),m+(i-3),l,k)  = \Gamma_{3}^{\left[s-(i-3)\atop{ s+m\atop s+m+l} \right]\times k} -
 4\cdot\Gamma_{2}^{\left[s-(i-3)\atop{ s+m-1\atop s+m+l} \right]\times k} -
8\cdot \Gamma_{2}^{\left[s-(i-3)\atop{ s+m\atop s+m+l-1} \right]\times k} +
32\cdot\Gamma_{1}^{\left[s-(i-3)\atop{ s+m-1\atop s+m+l-1} \right]\times k}\nonumber \\
& \nonumber  \\
& \text{ from \eqref{eq 6.7},\eqref{eq 6.12} and \eqref{eq 6.17} we obtain}\nonumber \\
& \nonumber  \\
&  \omega _{3}(s-(i-3),m+(i-3),l,k) = 6384 -4\cdot378 -8\cdot378 +32\cdot21 = 2520 \quad  \text{for $4\leq i\leq \inf(s-1,k)$}\label{eq 6.19}\\
& \nonumber
  \end{align}
 \underline{Computation of $ \omega _{i}(s,m,l,k+1) -  \omega _{i}(s,m,l,k)$  for  $4\leq i\leq \inf(s-1,k-1)$}
 \begin{align}
 & \text{By \eqref{eq 5.1}and \eqref{eq 6.19} we get for $ 4\leq i\leq \inf(s-1,k)$} \nonumber\\
  & \omega _{i}(s,m,l,k) = \displaystyle
 2^{i-3}\cdot\omega _{3}(s-(i-3),m+(i-3),l,k)  +   \sum_{j = 0}^{i-4}2^{j}\Delta _{i-j}^{\left[s-j \atop{ s+m\atop s+m+l} \right]\times k} \nonumber   \\
 & = 315\cdot2^{i} +  \sum_{j = 0}^{i-4}2^{j}\Delta _{i-j}^{\left[s-j \atop{ s+m\atop s+m+l} \right]\times k}\nonumber \\
 & \text{From \eqref{eq 4.11} we then deduce}\nonumber\\
 & \nonumber\\
  & \omega _{i}(s,m,l,k+1) -  \omega _{i}(s,m,l,k) = 0 \quad \text{for $4\leq i\leq \inf(s-1,k-1)$},\nonumber  \\
  & \text{that is }\nonumber \\
 &  \Gamma_{i}^{\left[s\atop{ s+m\atop s+m+l} \right]\times k+1} - \Gamma_{i}^{\left[s\atop{ s+m\atop s+m+l} \right]\times k}\label{eq 6.20}\\
    &  = 4\cdot\big(\Gamma_{i-1}^{\left[s\atop{ s+m-1\atop s+m+l} \right]\times k+1} - \Gamma_{i-1}^{\left[s\atop{ s+m-1\atop s+m+l} \right]\times k}\big)
   +  8\cdot\big(\Gamma_{i-1}^{\left[s\atop{ s+m\atop s+m+l-1} \right]\times k+1} - \Gamma_{i-1}^{\left[s\atop{ s+m\atop s+m+l-1} \right]\times k}\big) \nonumber \\
     &  - 32\cdot\big(\Gamma_{i-2}^{\left[s\atop{ s+m-1\atop s+m+l-1} \right]\times k+1} - \Gamma_{i-2}^{\left[s\atop{ s+m-1\atop s+m+l-1} \right]\times k}\big)  \quad \text{for $4\leq i\leq \inf(s-1,k-1)$} \nonumber
 \end{align}
  \underline { $  (H_{1})\wedge (H_{2})\wedge\ldots\wedge (H_{j})\; implies \; (H_{j+1})  $}\\
 
 From  \eqref{eq 6.20} with  i $\rightarrow $ j+1 we get  \\
 \begin{align}
  &  \Gamma_{j+1}^{\left[s\atop{ s+m\atop s+m+l} \right]\times k+1} - \Gamma_{j+1}^{\left[s\atop{ s+m\atop s+m+l} \right]\times k}\label{eq 6.21} \\
 &  = 4\cdot\big(\Gamma_{j}^{\left[s\atop{ s+m-1\atop s+m+l} \right]\times k+1} - \Gamma_{j}^{\left[s\atop{ s+m-1\atop s+m+l} \right]\times k}\big)
 +  8\cdot\big(\Gamma_{j}^{\left[s\atop{ s+m\atop s+m+l-1} \right]\times k+1} - \Gamma_{j}^{\left[s\atop{ s+m\atop s+m+l-1} \right]\times k}\big) \nonumber \\
 &  - 32\cdot\big(\Gamma_{j-1}^{\left[s\atop{ s+m-1\atop s+m+l-1} \right]\times k+1} - \Gamma_{j-1}^{\left[s\atop{ s+m-1\atop s+m+l-1} \right]\times k}\big)  \quad \text{for $s\geq j+2,\;k\geq j+2,\;m\geq 0,\;l\geq 0 $} \nonumber\\ 
  & \text{ By $ (H_{j}),\;  (H_{j-1}) $ we obtain respectively}\nonumber \\
  & \Gamma_{j}^{\left[s\atop{ s+m-1\atop s+m+l} \right]\times k} = \Gamma_{j}^{\left[s\atop{ s+m\atop s+m+l-1} \right]\times k} = \Gamma_{j}^{\left[j\atop{ j\atop j} \right]\times (j+1)} = \gamma _{j} \quad \text{ for  $s\geq j+2,\;k\geq j+2,\;m \geq 0,\;l \geq 0 $}\label{eq 6.22}\\
 & \Gamma_{j-1}^{\left[s\atop{ s+m-1\atop s+m+l-1} \right]\times k}  = \Gamma_{j-1}^{\left[j-1\atop{ j-1\atop j-1} \right]\times j} = \gamma _{j-1} \quad \text{ for  $s\geq j+1,\;k\geq j+1,\;m \geq 0,\;l \geq 0 $}\label{eq 6.23}\\   
 & \text{ we deduce from \eqref{eq 6.21}, \eqref{eq 6.22} and \eqref{eq 6.23}} \nonumber \\
 &   \Gamma_{j+1}^{\left[s\atop{ s+m\atop s+m+l} \right]\times k} =   \Gamma_{j+1}^{\left[s\atop{ s+m\atop s+m+l} \right]\times (j+2)} \quad \text{for $k\geq j+2,\;s\geq j+2,\;m\geq 0,\;l\geq 0$}\label{eq 6.24}
       \end{align}
Consider now the matrix\\
  $$ D^{\left[s\atop{ s+m\atop {s+m+l }}\right]\times (j+2)}(t,\eta ,\xi ) = \left ( \begin{array} {ccccc}
\alpha _{1} & \alpha _{2}  &  \ldots & \alpha _{j+1}  &  \alpha _{j+2} \\
\alpha _{2 } & \alpha _{3} &  \ldots  &  \alpha _{j+2} &  \alpha _{j+3} \\
\vdots & \vdots & \vdots    & \vdots  &  \vdots \\
\alpha _{s-1} & \alpha _{s} & \ldots  &  \alpha _{s+j-1} &  \alpha _{s+j}  \\
  \alpha  _{s } & \alpha  _{s +1} & \ldots & \alpha  _{s +j}& \alpha  _{s +j+1}\\
\beta  _{1} & \beta  _{2}  & \ldots  &  \beta_{j+1} &  \beta _{j+2}  \\
\beta  _{2} & \beta  _{3}  & \ldots  &  \beta_{j+2} &  \beta _{j+3}  \\
\vdots & \vdots    &  \vdots & \vdots  &  \vdots \\
\beta  _{m+1} & \beta  _{m+2}  & \ldots  &  \beta_{j+m+1} &  \beta _{j+m+2}  \\
\vdots & \vdots    &  \vdots & \vdots  &  \vdots \\
\beta  _{s+m-1} & \beta  _{s+m}  & \ldots  &  \beta_{s+m+j-1} &  \beta _{s+m+j}  \\
 \beta _{s+m} & \beta _{s+m+1} & \ldots & \beta _{s+m+j} & \beta _{s+m+j+1}\\
\gamma  _{1} & \gamma   _{2}  & \ldots  & \gamma  _{j+1} &  \gamma  _{j+2}  \\
\gamma  _{2} & \gamma  _{3}  & \ldots  & \gamma  _{j+2} &  \gamma  _{j+3}  \\
\vdots & \vdots    &  \vdots & \vdots  &  \vdots \\
 \gamma  _{m+l+1} &  \gamma _{m+l+2}  & \ldots  & \gamma _{j+m+l+1} &  \gamma  _{j+m+l+2}  \\
\vdots & \vdots   &  \vdots & \vdots  &  \vdots \\
  \gamma  _{s+m+l} & \gamma  _{s+m+l+1}  & \ldots  & \gamma  _{s+m+l+j} &  \gamma  _{s+m+l+j+1}  
\end{array}  \right). $$ 
 We have  respectively  by \eqref{eq 6.1},  \eqref{eq 6.2} with k = j+2 \\
 \begin{align}
  \sum_{i=0}^{j+2}\Gamma_{i}^{\left[s\atop{ s+m\atop s+m+l} \right]\times (j+2)}  &  = 2^{3j+3s+2m+l+3}, \label{eq 6.25} \\
  & \nonumber \\
    \sum_{i=0}^{j+2}\Gamma_{i}^{\left[s\atop{ s+m\atop s+m+l} \right]\times (j+2)} \cdot 2^{-i} &  =  2^{2j+3s+2m+l+1} + 2^{3j+3} - 2^{2j+1}.  \label{eq 6.26}
\end{align}
Setting $ s = j+1,\; m = 0,\; l = 0 $ in \eqref{eq 6.25}, \eqref{eq 6.26} we have \\
\begin{align}
  \sum_{i=0}^{j+2}\Gamma_{i}^{\left[j+1\atop{ j+1\atop j+1} \right]\times (j+2)}  &  = 2^{6j+6}, \label{eq 6.27} \\
  & \nonumber \\
    \sum_{i=0}^{j+2}\Gamma_{i}^{\left[j+1\atop{ j+1\atop j+1} \right]\times (j+2)} \cdot 2^{-i} &  =  2^{5j+4} + 2^{3j+3} - 2^{2j+1}.  \label{eq 6.28}
\end{align}
From   $  (H_{1})\wedge (H_{2})\wedge\ldots\wedge (H_{j})  $ it follows  that \\
 \begin{equation}
 \label{eq 6.29}
\Gamma_{i}^{\left[s\atop{ s+m\atop s+m+l} \right]\times (j+2)}= \Gamma_{i}^{\left[j+1\atop{ j+1 \atop j+1} \right]\times (j+2)}= \Gamma_{i}^{\left[i\atop{ i\atop i} \right]\times (i+1)} \quad \text{for $ 1\leq i\leq j $}
\end{equation}
We then  deduce from  \eqref{eq 6.25},  \eqref{eq 6.26} and \eqref{eq 6.29}  \\
\begin{equation}
\label{eq 6.30}
\Gamma_{j+1}^{\left[s\atop{ s+m\atop s+m+l} \right]\times (j+2)} =  2^{4j+5} - 2^{3j+3} + \sum_{i = 0}^{j}\gamma _{i}(1-2^{j+2-i})
\end{equation}
and from  \eqref{eq 6.27},  \eqref{eq 6.28} and \eqref{eq 6.29}   \\
\begin{equation}
\label{eq 6.31}
\Gamma_{j+1}^{\left[j+1\atop{ j+1\atop j+1} \right]\times (j+2)} =  2^{4j+5} - 2^{3j+3} + \sum_{i = 0}^{j}\gamma _{i}(1-2^{j+2-i})
\end{equation}
 
By  \eqref{eq 6.24},  \eqref{eq 6.30} and   \eqref{eq 6.31}  it follows that    \\ 
\begin{equation}
 \label{eq 6.32}
(H_{j+1})\quad \Gamma_{j+1}^{\left[s\atop{ s+m\atop s+m+l} \right]\times k} = \Gamma_{j+1}^{\left[j+1\atop{ j+1 \atop j+1} \right]\times (j+2)}= \gamma _{j+1}\quad \text{for all $ s\geq j+2,\;k\geq j+2,\;m\geq 0,\;l\geq 0 $}
\end{equation}
\end{proof}
\begin{lem}
\label{lem 6.3}We have  \\
\begin{equation}
\label{eq 6.33}
\Gamma_{i}^{\left[s\atop{ s+m\atop s+m+l} \right]\times k} = \Gamma_{i}^{\left[i\atop{ i \atop i} \right]\times (i+1)} = \begin{cases}
 1  & \text{if  } i = 0, \\
 105\cdot2^{4i-6} - 21\cdot2^{3i -5}  &  \text{if  }  1\leq i\leq s-1,\;k\geq i+1,\;m\geq 0,\;l\geq 0.
\end{cases}
\end{equation}
\end{lem}
\begin{proof}
   We have respectively by \eqref{eq 6.1}, \eqref{eq 6.2} with k = i+1, s = i and m=0,$\; l=0$ using \eqref{eq 6.3}
   \begin{align}
&  \sum_{ j = 0}^{i+1}  \Gamma  _{j}^{\left[i\atop{ i\atop i} \right]\times (i+1)}  = 
   \sum_{ j = 0}^{i}\gamma _{j} + \Gamma  _{i+1}^{\left[i\atop{ i\atop i} \right]\times (i+1)} = 2^{6i} \label{eq 6.34}\\
 &  \text{and}\nonumber \\
&  \sum_{j = 0}^{i+1}  \Gamma  _{j}^{\left[i\atop{ i\atop i} \right]\times (i+1)}\cdot2^{-j} 
   =  \sum_{ j = 0}^{i}\gamma _{j}\cdot2^{-j}  +  \Gamma  _{i+1}^{\left[i\atop{ i\atop i} \right]\times (i+1)}\cdot2^{-(i+1)}=
    2^{5i-1} + 2^{3i} - 2^{2i-1}. \label{eq 6.35} \\
    & \text{Subtracting \eqref{eq 6.34} from \eqref{eq 6.35} multiplied with $2^{i+1}$ we deduce }\nonumber \\
    &  \gamma _{i}= 2^{4i+1} - 2^{3i} +\sum_{ j = 0}^{i-1}\gamma _{j} - \sum_{j = 0}^{i-1}\gamma _{j}\cdot2^{i+1-j} \label{eq 6.36} \\
    & \text{Substituting i for i-1 in \eqref{eq 6.36} we get}\nonumber \\
    &  \gamma _{i-1}= 2^{4i-3} - 2^{3i-3} +\sum_{ j = 0}^{i-2}\gamma _{j} - \sum_{j = 0}^{i-2}\gamma _{j}\cdot2^{i-j} \label{eq 6.37} \\
    & \text{From \eqref{eq 6.36}, \eqref{eq 6.37} we obtain}\nonumber \\
    & \gamma _{i} -2\cdot\gamma _{i-1}=  2^{4i+1} - 2^{3i} +\sum_{ j = 0}^{i-1}\gamma _{j} - \sum_{j = 0}^{i-1}\gamma _{j}\cdot2^{i+1-j}
    -2\cdot\big( 2^{4i-3} - 2^{3i-3} +\sum_{ j = 0}^{i-2}\gamma _{j} - \sum_{j = 0}^{i-2}\gamma _{j}\cdot2^{i-j}  \big) \nonumber \\
    & \Longleftrightarrow \nonumber \\
     & \gamma _{i} -2\cdot\gamma _{i-1} =  7\cdot2^{4i-2} - 3\cdot2^{3i-2} - \sum_{ j = 0}^{i-2}\gamma _{j} -3\cdot\gamma _{i-1}\nonumber \\
      & \Longleftrightarrow \nonumber \\
     & \sum_{ j = 0}^{i}\gamma _{j}  =  7\cdot2^{4i-2} - 3\cdot2^{3i-2} \label{eq 6.38} \\
     & \text{Finally from \eqref{eq 6.38} we deduce}\nonumber \\
     & \gamma _{i} =  \sum_{ j = 0}^{i}\gamma _{j} -  \sum_{ j = 0}^{i-1}\gamma _{j}
      =  7\cdot2^{4i-2} - 3\cdot2^{3i-2} - \big( 7\cdot2^{4i-6} - 3\cdot2^{3i-5} \big) =  105\cdot2^{4i-6} - 21\cdot2^{3i -5}\nonumber 
\end{align}

\end{proof}

\begin{lem}
\label{lem 6.4}We have  \\
\begin{align}
&  \sum_{ j = 0}^{s-1}\gamma _{j} = 7\cdot2^{4s-6} - 3\cdot2^{3s-5}\label{eq 6.39} \\
&  \sum_{ j = 0}^{s-1}\gamma _{j}\cdot2^{-j} = 15\cdot2^{3s-6} - 7\cdot2^{2s-5}\label{eq 6.40} \\
& \Gamma_{i}^{\left[s\atop{ s+m\atop s+m+l} \right]\times i} = 2^{3s+2m+l+3i-3} -7\cdot2^{4i-6} +3\cdot2^{3i-5} \quad \text{for $ 1\leq i\leq s+1,\;m\geq 0,\;l\geq 0 $} \label{eq 6.41}
\end{align}
\end{lem}
\begin{proof}
\begin{align*}
&  \text{Proof of \eqref{eq 6.39} }\\
& \\
& \text{From \eqref{eq 6.33} we obtain}\\
&  \sum_{ j = 0}^{s-1}\gamma _{j} =  \sum_{ j = 0}^{s-1}\big( 105\cdot2^{4j-6} - 21\cdot2^{3j -5}\big) =  7\cdot2^{4s-6} - 3\cdot2^{3s-5} \\
&  \text{Proof of \eqref{eq 6.40} }\\
& \\
&  \sum_{ j = 0}^{s-1}\gamma _{j}\cdot2^{-j} =  \sum_{ j = 0}^{s-1}\big( 105\cdot2^{4j-6} - 21\cdot2^{3j -5}\big)\cdot2^{-j}
 = \sum_{ j = 0}^{s-1}\big( 105\cdot2^{3j-6} - 21\cdot2^{2j -5}\big) = 15\cdot2^{3s-6} - 7\cdot2^{2s-5} \\
&  \text{Proof of \eqref{eq 6.41} in the case $1 \leq i \leq s .$}\\
& \\
 &\text{ From \eqref{eq 6.1} with k = i  using \eqref{eq 6.33} we obtain} \\
 & \Gamma_{i}^{\left[s\atop{ s+m\atop s+m+l} \right]\times i} = 2^{3s+2m+l+3i-3} - \sum_{ j = 0}^{i-1}\gamma _{j} 
=  2^{3s+2m+l+3i-3} - \sum_{ j = 0}^{i-1}(105 \cdot2^{4j-6} - 21\cdot 2^{3j -5}) \\
& =  2^{3s+2m+l+3i-3} -7\cdot2^{4i-6} +3\cdot2^{3i-5}\\
&  \text{Proof of \eqref{eq 6.41} in the case  i = s+1}\\
& \\
 & \text{ From \eqref{eq 6.1} with k = s+1  using \eqref{eq 6.39} we obtain }\\
  & \Gamma_{s}^{\left[s\atop{ s+m\atop s+m+l} \right]\times (s+1)} + \Gamma_{s+1}^{\left[s\atop{ s+m\atop s+m+l} \right]\times (s+1)} = 2^{6s+2m+l} - \sum_{ j = 0}^{s-1}\gamma _{j}
  =  2^{6s+2m+l} -7\cdot2^{4s-6} +3\cdot2^{3s-5} \\
  & \text{ From \eqref{eq 6.2} with k = s+1  using \eqref{eq 6.40} we obtain }\\
 & 2\cdot \Gamma_{s}^{\left[s\atop{ s+m\atop s+m+l} \right]\times (s+1)} + \Gamma_{s+1}^{\left[s\atop{ s+m\atop s+m+l} \right]\times (s+1)}
  = 2^{6s+2m+l} +2^{4s+1} -2^{3s}   - 2^{s+1}\cdot\sum_{ j = 0}^{s-1}\gamma _{j}\cdot2^{-j}\\
  &  =  2^{6s+2m+l}  +2^{4s+1} -2^{3s} - 2^{s+1}\cdot\big(15\cdot2^{3s-6} -7\cdot2^{2s-5}\big)   \\
  & =  2^{6s+2m+l}  +2^{4s+1} -2^{3s} -15\cdot2^{4s-5} + 7\cdot2^{3s-4} =   2^{6s+2m+l} + 49\cdot2^{4s-5}-9\cdot2^{3s-4}\\
  & \text{We then deduce from the above equations}
   \end{align*}
   
   \begin{align}
  & \Gamma_{s}^{\left[s\atop{ s+m\atop s+m+l} \right]\times (s+1)} =  105\cdot2^{4s-6} - 21\cdot2^{3s -5},\quad \text{remark that \eqref{eq 6.33} holds for i = s, \; k = s+1} \label{eq 6.42}\\
  & \Gamma_{s+1}^{\left[s\atop{ s+m\atop s+m+l} \right]\times (s+1)} =  2^{6s+2m+l}  -7\cdot2^{4s-2} +3\cdot2^{3s-2}\label{eq 6.43}
 \end{align}

 \end{proof}
  \section{\textbf{Computation of  $\omega _{i}(s,m,l,k+1) -\omega _{i}(s,m,l,k) $  in the case $ s\leq i\leq \inf(3s+2m+l,k),\; k>i$}}
 \label{sec 7}
  \subsection{Notation}
  \label{subsec 1}
  \begin{defn}
  \label{defn 7.1}
 \begin{align*}
&\text{We define $ A_{i}^{[s,m,l,k]} $ to be equal to}
 &    \Gamma_{i-s}^{\left[s+m \atop s+m+l \right]\times k} -
 4\cdot\Gamma_{i-s-1}^{\left[ s+m-1\atop s+m+l \right]\times k} -
8\cdot \Gamma_{i-s-1}^{\left[ s+m\atop s+m+l-1 \right]\times k} +
32\cdot\Gamma_{i-s-2}^{\left[ s+m-1\atop s+m+l-1 \right]\times k} \\
& \text{for} \quad s\leqslant i \leqslant \inf(k, 3s+2m+l)
 \end{align*}
We recall that ( see \ref{defn 5.1} )
  \begin{align*}
&   \omega _{i}(s,m,l,k+1) -  \omega _{i}(s,m,l,k) \\
  & =  \big(\Gamma_{i}^{\left[s\atop{ s+m\atop s+m+l} \right]\times( k+1)} - \Gamma_{i}^{\left[s\atop{ s+m\atop s+m+l} \right]\times k} \big) -4\cdot \big(\Gamma_{i-1}^{\left[s\atop{ s+m-1\atop s+m+l} \right]\times( k+1)} -\Gamma_{i-1}^{\left[s\atop{ s+m-1\atop s+m+l} \right]\times k} \big)\\
  & -8\cdot\big( \Gamma_{i-1}^{\left[s\atop{ s+m\atop s+m+l-1} \right]\times( k+1)} - \Gamma_{i-1}^{\left[s\atop{ s+m\atop s+m+l-1} \right]\times k} \big) + 32\cdot\big(\Gamma_{i-2}^{\left[s\atop{ s+m-1\atop s+m+l-1} \right]\times (k+1)} -\Gamma_{i-2}^{\left[s\atop{ s+m-1\atop s+m+l-1} \right]\times k}\big)
    \end{align*}

  \end{defn}
  
   \subsection{Introduction}
  \label{subsec 2}
   We adapt the method used in Section 10 of [2] to compute explicitly \\[0.02 cm]
 $\omega _{i}(s,m,l,k+1) -  \omega _{i}(s,m,l,k) $ for $ s\leqslant i \leqslant \inf(k, 3s+2m+l),\; k>i $ \\[0.01 cm]

 \subsection{\textbf{An equivalent to the formula for  $\omega _{i}(s,m,l,k)$ in the case  $ s\leq i\leq \inf(3s+2m+l,k) $}}
  \label{subsec 3}
   \begin{lem}
\label{lem  7.2} $ Let \; s \geq 2,  m\geq 0, \; l\geq 0\;  and \; k \geq 1, $then  we have  \\
\begin{align}
 \displaystyle
&  \omega _{i}(s,m,l,k) = 
  2^{k+s-1} A_{i}^{[s,m,l,k]} + 2^{i} \Delta _{i-s+1}^{\left[s+m \atop s+m+l \right]\times k}
 - 2^{i-1} \Delta _{i-s}^{\left[s+m \atop s+m+l \right]\times k} +  \sum_{j = 0}^{s-2}2^{j}\Delta _{i-j}^{\left[s-j \atop{ s+m\atop s+m+l} \right]\times k}\quad \text{if  }  s\leq i\leq \inf(3s+2m+l,k) \label{eq 7.1}\\
 &  \omega _{i}(s,m,l,k+1) - \omega _{i}(s,m,l,k) =  2^{k+s} A_{i}^{[s,m,l,k+1]} - 2^{k+s-1} A_{i}^{[s,m,l,k]} \quad \text{if }  s\leq i\leq \inf(3s+2m+l,k), \;k >i  \label{eq 7.2}  \\
 & \text{where } \nonumber \\
 &  A_{i}^{[s,m,l,k]} =   \Gamma_{i-s}^{\left[s+m \atop s+m+l \right]\times k} -
 4\cdot\Gamma_{i-s-1}^{\left[ s+m-1\atop s+m+l \right]\times k} -
8\cdot \Gamma_{i-s-1}^{\left[ s+m\atop s+m+l-1 \right]\times k} +
32\cdot\Gamma_{i-s-2}^{\left[ s+m-1\atop s+m+l-1 \right]\times k} \nonumber
 \end{align}
\end{lem}
\begin{proof}

\begin{align*}
  &  \underline{proof \; of  \; \eqref{eq 7.1}}\\
  & \text{From Definition \ref{defn 5.1} and  Theorem \ref{thm 1.4} we obtain }\\
&  \omega _{i-(s-1)}(1,m+s-1,l,k)  = \Gamma_{i-(s-1)}^{\left[1\atop{ s+m\atop s+m+l} \right]\times k} -
 4\cdot\Gamma_{i-s}^{\left[1\atop{ s+m-1\atop s+m+l} \right]\times k} -
8\cdot \Gamma_{i-s}^{\left[1\atop{ s+m\atop s+m+l-1} \right]\times k} +
32\cdot\Gamma_{i-s-1}^{\left[1\atop{ s+m-1\atop s+m+l-1} \right]\times k}\\
& = (2^{k} -2^{i-s}) \Gamma_{i-s}^{\left[ s+m\atop s+m+l \right]\times k} +2^{i-s+1}\Gamma_{i-s+1}^{\left[ s+m\atop s+m+l \right]\times k}\\
& -4\cdot\big[ (2^{k} -2^{i-s-1})\Gamma_{i-s-1}^{\left[ s+m-1\atop s+m+l \right]\times k} +2^{i-s}\Gamma_{i-s}^{\left[ s+m-1\atop s+m+l \right]\times k}\big]\\
& -8\cdot\big[ (2^{k} -2^{i-s-1})\Gamma_{i-s-1}^{\left[ s+m\atop s+m+l-1 \right]\times k} +2^{i-s}\Gamma_{i-s}^{\left[ s+m\atop s+m+l-1 \right]\times k}\big]\\
& +32\cdot\big[ (2^{k} -2^{i-s-2})\Gamma_{i-s-2}^{\left[ s+m-1\atop s+m+l-1 \right]\times k} +2^{i-s-11}\Gamma_{i-s-1}^{\left[ s+m-1\atop s+m+l \right]\times k}\big]\\
& = 2^{k}\cdot\big[ \Gamma_{i-s}^{\left[s+m \atop s+m+l \right]\times k} -
 4\cdot\Gamma_{i-s-1}^{\left[ s+m-1\atop s+m+l \right]\times k} -
8\cdot \Gamma_{i-s-1}^{\left[ s+m\atop s+m+l-1 \right]\times k} +
32 \cdot\Gamma_{i-s-2}^{\left[ s+m-1\atop s+m+l-1 \right]\times k}\big]\\
& +  2^{i-s+1}\cdot\big[ \Gamma_{i-s+1}^{\left[s+m \atop s+m+l \right]\times k} -
 2 \cdot\Gamma_{i-s}^{\left[ s+m-1\atop s+m+l \right]\times k} -
4\cdot \Gamma_{i-s}^{\left[ s+m\atop s+m+l-1 \right]\times k} +
8 \cdot\Gamma_{i-s-1}^{\left[ s+m-1\atop s+m+l-1 \right]\times k}\big]\\
& - 2^{i-s}\cdot\big[ \Gamma_{i-s}^{\left[s+m \atop s+m+l \right]\times k} -
 2 \cdot\Gamma_{i-s-1}^{\left[ s+m-1\atop s+m+l \right]\times k} -
4 \cdot\Gamma_{i-s-1}^{\left[ s+m\atop s+m+l-1 \right]\times k} +
8\cdot \Gamma_{i-s-2}^{\left[ s+m-1\atop s+m+l-1 \right]\times k}\big]\\
& =  2^{k} A_{i}^{[s,m,l,k]} + 2^{i-s+1} \Delta _{i-s+1}^{\left[s+m \atop s+m+l \right]\times k}
 - 2^{i-s} \Delta _{i-s}^{\left[s+m \atop s+m+l \right]\times k} \\
 & \text{ Hence from \eqref{eq 5.1} we  get easily \eqref{eq 7.1}}\\
 &  \underline{proof \; of  \; \eqref{eq 7.2}}\\
 & \text{From Lemma \ref{lem 4.5} we obtain }\\
 &  \sum_{j = 0}^{s-2}2^{j}\Delta _{i-j}^{\left[s-j \atop{ s+m\atop s+m+l} \right]\times (k+1)}- \sum_{j = 0}^{s-2}2^{j}\Delta _{i-j}^{\left[s-j \atop{ s+m\atop s+m+l} \right]\times k}= 0 \quad  \text{for}\quad  k>i, \\
 & \text{and from Lemma 8.5 (see [2] section 8) we get}\\
 & \Delta _{i-s+1}^{\left[s+m \atop s+m+l \right]\times (k+1)}-\Delta _{i-s+1}^{\left[s+m \atop s+m+l \right]\times k} = 0,\quad
  \Delta _{i-s}^{\left[s+m \atop s+m+l \right]\times (k+1)}-\Delta _{i-s}^{\left[s+m \atop s+m+l \right]\times k} = 0 \quad \text{for}\quad k>i \\
  & \text{Hence from the above equations we obtain easily \eqref{eq 7.2}}
 \end{align*}
\end{proof}

 \subsection{\textbf{Computation of $A_{i}^{[s,m,l,k]}$}}
  \label{subsec 4}
   \begin{lem}
\label{lem  7.3}
Set for $ s\leq i\leq \inf(k,3s+2m+l) $\\
$$A_{i}^{[s,m,l,k]} = \Gamma _{i-s}^{\Big[\substack{m+s \\ m+s+l }\Big] \times k} - 4\cdot \Gamma _{i-s-1}^{\Big[\substack{m+s-1 \\ m+s+l }\Big] \times k}
   -8\cdot\Gamma _{i-s-1}^{\Big[\substack{m+s \\ m+s+l-1 }\Big] \times k} +32\cdot\Gamma _{i-s-2}^{\Big[\substack{m+s-1 \\ m+s+l-1}\Big] \times k}$$

   We have in the case \underline{ $l = 0, m =  0, s = 2 $}\\
 $$A_{i}^{[2,0,0,k]} = \Gamma _{i-2}^{\Big[\substack{2 \\ 2}\Big] \times k} - 12\cdot \Gamma _{i-3}^{\Big[\substack{1 \\ 2}\Big] \times k}
 +32\cdot\Gamma _{i-4}^{\Big[\substack{1 \\ 1}\Big] \times k}$$
 \begin{equation}
 \label{eq 7.3}
   A_{i}^{[2,0,0,k]} =  
\begin{cases}
1  &\text{if } i = 2 \\
-3  &\text{if  } i = 3\\
-3\cdot2^{k+1} +2 &\text{if  } i = 4 \\
3\cdot2^{k+ 1}  &\text{if  } i = 5 \\
3\cdot2^{2k +2}  &\text{if  } i = 6
\end{cases}        
\end{equation}
  We have in the case \underline{ $l = 0, m =  1, s = 2 $}\\
 $$A_{i}^{[2,1,0,k]} = \Gamma _{i-2}^{\Big[\substack{3 \\ 3}\Big] \times k} - 12\cdot \Gamma _{i-3}^{\Big[\substack{2 \\ 3}\Big] \times k}
 +32\cdot\Gamma _{i-4}^{\Big[\substack{2 \\ 2}\Big] \times k}$$
 \begin{equation}
\label{eq 7.4}
   A_{i}^{[2,1,0,k]} =  
\begin{cases}
1  &\text{if } i = 2 \\
-3  &\text{if  } i = 3\\
2  &\text{if  } i = 4 \\
-3\cdot2^{k+2}  &\text{if  } i = 5 \\
3\cdot2^{k+ 2}  &\text{if  } i = 6 \\
0 &\text{if  } i = 7 \\
3\cdot2^{2k +4}  &\text{if  } i = 8
\end{cases}            
\end{equation}

We have in the case \underline{ $l = 0, m\geq 2, s = 2 $}\\
$$A_{i}^{[2,m,0,k]} = \Gamma _{i-2}^{\Big[\substack{m+2 \\ m+2 }\Big] \times k} - 12\cdot \Gamma _{i-3}^{\Big[\substack{m+ 1 \\ m+2 }\Big] \times k}
   +32\cdot\Gamma _{i- 4}^{\Big[\substack{m+1 \\ m+1}\Big] \times k}$$
 \begin{equation}
\label{eq 7.5}
 A_{i}^{[2,m,0,k]} =  
\begin{cases}
1  &\text{if } i = 2 \\
-3  &\text{if  } i = 3 \\
2 &\text{if } i = 4 \\
0  &\text{if  } 5 \leq i\leq m+3 \\
-3\cdot2^{k+m+1}  &\text{if  } i = m+ 4 \\
3\cdot2^{k+m +1}  &\text{if  } i = m+5 \\
0  &\text{if  } m +7\leq i\leq 2m+5 \\
3\cdot2^{2k+2m +2}  &\text{if  } i = 2m+6 
\end{cases}
\end{equation}
    We have in the case \underline{ $l = 1, m =  0, s = 2 $}\\
 $$A_{i}^{[2,0,1,k]} = \Gamma _{i-2}^{\Big[\substack{2 \\ 3}\Big] \times k} - 4\cdot \Gamma _{i-3}^{\Big[\substack{1 \\ 3}\Big] \times k} -8\cdot\Gamma _{i-3}^{\Big[\substack{2 \\ 2}\Big] \times k}
 +32\cdot\Gamma _{i-4}^{\Big[\substack{2 \\ 2}\Big] \times k}$$
 \begin{equation}
\label{eq 7.6}
   A_{i}^{[2,0,1,k]} =  
\begin{cases}
1  &\text{if } i = 2 \\
-3  &\text{if  } i = 3\\
-2^{k+1} + 2  &\text{if  } i = 4 \\
-3\cdot2^{k+1}  &\text{if  } i = 5 \\
2^{k+ 3}  &\text{if  } i = 6 \\
3\cdot2^{2k +3}  &\text{if  } i = 7
\end{cases}
\end{equation}

 We have in the case \underline{ $l = 1, m =  1, s = 2 $}\\  
 $$A_{i}^{[2,1,1,k]} = \Gamma _{i-2}^{\Big[\substack{3 \\ 4}\Big] \times k} - 4\cdot \Gamma _{i-3}^{\Big[\substack{2 \\ 4}\Big] \times k} -8\cdot\Gamma _{i-3}^{\Big[\substack{3 \\ 3}\Big] \times k}
 +32\cdot\Gamma _{i-4}^{\Big[\substack{2 \\ 3}\Big] \times k}$$
  \begin{equation}
  \label{eq 7.7}
   A_{i}^{[2,1,1,k]} =  
\begin{cases}
1  &\text{if } i = 2 \\
-3  &\text{if  } i = 3\\
2  &\text{if  } i = 4\\
-2^{k+2}  &\text{if  } i = 5 \\
-3\cdot2^{k+2}  &\text{if  } i = 6 \\
2^{k+ 4}  &\text{if  } i = 7 \\
0     &\text{if  } i = 8 \\
3\cdot2^{2k +5}  &\text{if  } i = 9
\end{cases}  
\end{equation}

  We have in the case \underline{ $l = 1, m\geq 2, s =2 $}\\
 $$A_{i}^{[2,m,1,k]} = \Gamma _{i-2}^{\Big[\substack{m+2 \\ m+3 }\Big] \times k} - 4\cdot \Gamma _{i-3}^{\Big[\substack{m+1 \\ m+3 }\Big] \times k}
   -8\cdot\Gamma _{i -3}^{\Big[\substack{m+2 \\ m+2 }\Big] \times k} +32\cdot\Gamma _{i-4}^{\Big[\substack{m+1 \\ m+2}\Big] \times k}$$
 \begin{equation}
\label{eq 7.8}
  A_{i}^{[2,m,1,k]} =  
\begin{cases}
1  &\text{if } i = 2 \\
-3  &\text{if  } i = 3\\
2 &\text{if } i = 4 \\
0  &\text{if  } 5\leq i\leq m +3 \\
-2^{k+m+1}  &\text{if  } i = m+4\\
-3\cdot2^{k+m+1}  &\text{if  } i = m+ 5\\
2^{k+m+3}  &\text{if  } i = m+6 \\
0  &\text{if  } m+7\leq i\leq 2m+6\\
3\cdot2^{2k+2m +3}  &\text{if  } i = 2m+7 
\end{cases}
\end{equation}

We have in the case \underline{ $l = 2, m = 0 , s = 2 $}\\
$$A_{i}^{[2,0,2,k]} = \Gamma _{i-2}^{\Big[\substack{2 \\ 4 }\Big] \times k} - 4\cdot \Gamma _{i-3}^{\Big[\substack{1 \\ 4 }\Big] \times k}
   -8\cdot\Gamma _{i- 3}^{\Big[\substack{2 \\ 3 }\Big] \times k} +32\cdot\Gamma _{i- 4}^{\Big[\substack{1 \\ 3}\Big] \times k}$$
\begin{equation}
\label{eq 7.9}
 A_{i}^{[2,0,2,k]} =  
\begin{cases}
1  &\text{if } i = 2 \\
-3  &\text{if  } i = 3 \\
-2^{k+1} + 2 &\text{if } i = 4 \\
5\cdot2^{k+1}  &\text{if  } i =  5\\
-5\cdot2^{k+3}  &\text{if  }   i = 6\\
2^{k+5}  &\text{if  } i = 7 \\
3\cdot2^{2k+4}  &\text{if  } i = 8 
\end{cases}
\end{equation}
   
    We have in the case \underline{ $l = 2, m = 1, s = 2 $}\\
$$A_{i}^{[2,1,2,k]} = \Gamma _{i-2}^{\Big[\substack{3 \\ 5 }\Big] \times k} - 4\cdot \Gamma _{i-3}^{\Big[\substack{2 \\ 5 }\Big] \times k}
   -8\cdot\Gamma _{i- 3}^{\Big[\substack{3 \\ 4 }\Big] \times k} +32\cdot\Gamma _{i- 4}^{\Big[\substack{2 \\ 4}\Big] \times k}$$
  \begin{equation}
\label{eq 7.10}  
  A_{i}^{[2,1,2,k]} =  
\begin{cases}
1  &\text{if } i = 2 \\
-3  &\text{if  } i = 3 \\
2 &\text{if } i = 4 \\
-2^{k+2}  &\text{if  } i = 5 \\
5\cdot2^{k+2}  &\text{if  } i = 6 \\
-5\cdot2^{k+4}  &\text{if  } i = 7 \\
2^{k+6}  &\text{if  } i =  8 \\
0  &\text{if  } i = 9 \\
3\cdot2^{2k+6}  &\text{if  } i = 10
\end{cases}
\end{equation}
  We have in the case \underline{ $l = 2, m\geq 2, s = 2 $}\\
$$A_{i}^{[2,m,2,k]} = \Gamma _{i-2}^{\Big[\substack{m+2 \\ m+4 }\Big] \times k} - 4\cdot \Gamma _{i-3}^{\Big[\substack{m+ 1 \\ m+4 }\Big] \times k}
   -8\cdot\Gamma _{i- 3}^{\Big[\substack{m+2 \\ m+3 }\Big] \times k} +32\cdot\Gamma _{i- 4}^{\Big[\substack{m+1 \\ m+3}\Big] \times k}$$
\begin{equation}
\label{eq 7.11}
 A_{i}^{[2,m,2,k]} =  
\begin{cases}
1  &\text{if } i = 2 \\
-3  &\text{if  } i = 3 \\
2 &\text{if } i = 4 \\
0  &\text{if  } 5 \leq i\leq m+3 \\
-2^{k+m+1}  &\text{if  } i = m+ 4 \\
5\cdot2^{k+m+1}  &\text{if  } i = m+ 5\\
-5\cdot2^{k+m+3}  &\text{if  } i = m+6 \\
2^{k+m+5}  &\text{if  } i = m+ 7 \\
0  &\text{if  } m+ 8\leq i\leq 2m+7 \\
3\cdot2^{2k+2m+4}  &\text{if  } i = 2m+8
\end{cases}
\end{equation}

   We have in the case \underline{ $l\geq 3, m = 0 , s = 2 $}\\
$$A_{i}^{[2,0,l,k]} = \Gamma _{i-2}^{\Big[\substack{2 \\ 2+l }\Big] \times k} - 4\cdot \Gamma _{i-3}^{\Big[\substack{1 \\ 2+l }\Big] \times k}
   -8\cdot\Gamma _{i- 3}^{\Big[\substack{2 \\ 1+l }\Big] \times k} +32\cdot\Gamma _{i- 4}^{\Big[\substack{1 \\ 1+l}\Big] \times k}$$
\begin{equation}
\label{eq 7.12}
 A_{i}^{[2,0,l,k]} =  
\begin{cases}
1  &\text{if } i = 2 \\
-3  &\text{if  } i = 3 \\
-2^{k+1} + 2 &\text{if } i = 4 \\
5\cdot2^{k+1}  &\text{if  } i =  5\\
3\cdot2^{k+2i-9}  &\text{if  }  6\leq i\leq  l+3 \\
-5\cdot2^{k+2l-1}  &\text{if  } i = l+4 \\
2^{k+2l +1}  &\text{if  } i =  l+5 \\
3\cdot2^{2k+l +2}  &\text{if  } i = l+6 
\end{cases}
\end{equation}
    We have in the case \underline{ $l\geq 3, m = 1, s = 2 $}\\
$$A_{i}^{[2,1,l,k]} = \Gamma _{i-2}^{\Big[\substack{3 \\ 3+l }\Big] \times k} - 4\cdot \Gamma _{i-3}^{\Big[\substack{2 \\ 3+l }\Big] \times k}
   -8\cdot\Gamma _{i- 3}^{\Big[\substack{3 \\ 2+l }\Big] \times k} +32\cdot\Gamma _{i- 4}^{\Big[\substack{2 \\ 2+l}\Big] \times k}$$
 \begin{equation}
\label{eq 7.13}
 A_{i}^{[2,1,l,k]} =  
\begin{cases}
1  &\text{if } i = 2 \\
-3  &\text{if  } i = 3 \\
2 &\text{if } i = 4 \\
-2^{k+2}  &\text{if  } i = 5 \\
5\cdot2^{k+2}  &\text{if  } i = 6 \\
3\cdot2^{k +2i- 10}  &\text{if  } 7 \leq i\leq  l+4 \\
-5\cdot2^{k+2l}  &\text{if  } i = l+5 \\
2^{k+2l +2}  &\text{if  } i =  l+6 \\
0  &\text{if  } i = l+7 \\
3\cdot2^{2k+l +4}  &\text{if  } i = l+8
\end{cases} 
\end{equation}

 We have in the case \underline{ $l\geq 3, m\geq 2, s = 2 $}\\
$$A_{i}^{[2,m,l,k]} = \Gamma _{i-2}^{\Big[\substack{m+2 \\ m+2+l }\Big] \times k} - 4\cdot \Gamma _{i-3}^{\Big[\substack{m+ 1 \\ m+2+l }\Big] \times k}
   -8\cdot\Gamma _{i- 3}^{\Big[\substack{m+2 \\ m+1+l }\Big] \times k} +32\cdot\Gamma _{i- 4}^{\Big[\substack{m+1 \\ m+1+l}\Big] \times k}$$
  \begin{equation}
  \label{eq 7.14}
 A_{i}^{[2,m,l,k]} =  
\begin{cases}
1  &\text{if } i = 2 \\
-3  &\text{if  } i = 3 \\
2 &\text{if } i = 4 \\
0  &\text{if  } 5 \leq i\leq m+3 \\
-2^{k+m+1}  &\text{if  } i = m+ 4 \\
5\cdot2^{k+m+1}  &\text{if  } i = m+ 5\\
3\cdot2^{k-m +2i-9}  &\text{if  } m+ 6\leq i\leq m+ l+3 \\
-5\cdot2^{k+m+2l-1}  &\text{if  } i = m+l+4 \\
2^{k+m+2l +1}  &\text{if  } i = m+ l+5 \\
0  &\text{if  } m+ l+6\leq i\leq 2m+l+5 \\
3\cdot2^{2k+2m+l +2}  &\text{if  } i = 2m+l+6 
\end{cases}
\end{equation}
 
     We have in the case \underline{ $l = 0, m =  0, s = 3 $}\\
 $$A_{i}^{[3,0,0,k]} = \Gamma _{i-3}^{\Big[\substack{3 \\ 3}\Big] \times k} - 12\cdot \Gamma _{i-4}^{\Big[\substack{2 \\ 3 }\Big] \times k}
 +32\cdot\Gamma _{i-5}^{\Big[\substack{2 \\ 2}\Big] \times k}$$
 \begin{equation}
\label{eq 7.15}
    A_{i}^{[3,0,0,k]} =  
\begin{cases}
1  &\text{if } i = 3 \\
-3  &\text{if  } i = 4\\
2 &\text{if } i = 5 \\
-3\cdot2^{k+ 2}  &\text{if  } i = 6 \\
3\cdot2^{k+ 2}  &\text{if  } i = 7 \\
0  &\text{if  }             i = 8 \\
3\cdot2^{2k+ 4}  &\text{if  } i = 9
\end{cases}     
\end{equation}

  We have in the case \underline{ $l = 0, m\geq 1, s = 3 $}\\
 $$A_{i}^{[s,m,0,k]} = \Gamma _{i-3}^{\Big[\substack{m+3 \\ m+3}\Big] \times k} - 12\cdot \Gamma _{i- 4}^{\Big[\substack{m+2 \\ m+3}\Big] \times k}
 +32\cdot\Gamma _{i- 5}^{\Big[\substack{m+2 \\ m+2}\Big] \times k}$$
\begin{equation}
\label{eq 7.16}
   A_{i}^{[3,m,0,k]} =  
\begin{cases}
1  &\text{if } i = 3 \\
-3  &\text{if  } i = 4\\
2 &\text{if } i = 5 \\
0  &\text{if  } 6  \leq i\leq m+ 5 \\
-3\cdot2^{k+m+2}  &\text{if  } i = m+6 \\
3\cdot2^{k+m+ 2}  &\text{if  } i = m+7\\
0  &\text{if  } m+ 8\leq i\leq 2m+ 8\\
3\cdot2^{2k+2m +4}  &\text{if  } i = 2m+ 9 
\end{cases}
\end{equation}
   We have in the case \underline{ $l = 1, m = 0, s =3 $}\\
 $$A_{i}^{[3,0,1,k]} = \Gamma _{i-3}^{\Big[\substack{3 \\ 4 }\Big] \times k} - 4\cdot \Gamma _{i-4}^{\Big[\substack{2 \\ 4 }\Big] \times k}
   -8\cdot\Gamma _{i-4}^{\Big[\substack{3 \\ 3 }\Big] \times k} +32\cdot\Gamma _{i- 5}^{\Big[\substack{2 \\ 3}\Big] \times k}$$
  \begin{equation}
 \label{eq 7.17} 
  A_{i}^{[3,0,1,k]} =  
\begin{cases}
1  &\text{if } i = 3 \\
-3  &\text{if  } i = 4\\
2 &\text{if } i = 5  \\
-2^{k+2}  &\text{if  } i = 6 \\
-3\cdot2^{k+2}  &\text{if  } i = 7\\
2^{k+4}  &\text{if  } i = 8 \\
0  &\text{if  } i = 9 \\
3\cdot2^{2k +5}  &\text{if  } i = 10 
\end{cases}
 \end{equation}
 
  We have in the case \underline{ $l = 1, m\geq 1, s =3 $}\\
 $$A_{i}^{[3,m,1,k]} = \Gamma _{i-3}^{\Big[\substack{m+3 \\ m+4 }\Big] \times k} - 4\cdot \Gamma _{i-4}^{\Big[\substack{m+2 \\ m+4 }\Big] \times k}
   -8\cdot\Gamma _{i-4}^{\Big[\substack{m+3 \\ m+3 }\Big] \times k} +32\cdot\Gamma _{i- 5}^{\Big[\substack{m+2 \\ m+3}\Big] \times k}$$
   \begin{equation}
  \label{eq 7.18} 
    A_{i}^{[3,m,1,k]} =  
\begin{cases}
1  &\text{if } i = 3 \\
-3  &\text{if  } i = 4\\
2 &\text{if } i = 5  \\
0  &\text{if  } 6 \leq i\leq m+5  \\
-2^{k+m+2}  &\text{if  } i = m+6 \\
-3\cdot2^{k+m+2}  &\text{if  } i = m+ 7\\
2^{k+m+4}  &\text{if  } i = m+8 \\
0  &\text{if  } m+9\leq i\leq 2m+9 \\
3\cdot2^{2k+2m +5}  &\text{if  } i = 2m+10 
\end{cases}
   \end{equation}

  We have in the case \underline{ $l = 2, m = 0, s = 3 $}\\
 $$A_{i}^{[3,0,2,k]} = \Gamma _{i-3}^{\Big[\substack{3 \\ 5 }\Big] \times k} - 4\cdot \Gamma _{i-4}^{\Big[\substack{2 \\ 5 }\Big] \times k}
   -8\cdot\Gamma _{i-4}^{\Big[\substack{3 \\ 4 }\Big] \times k} +32\cdot\Gamma _{i-5}^{\Big[\substack{2 \\ 4}\Big] \times k}$$
  \begin{equation}
 \label{eq 7.19} 
  A_{i}^{[3,0,2,k]} =  
\begin{cases}
1  &\text{if } i = 3 \\
-3  &\text{if  } i = 4\\
2 &\text{if } i = 5 \\
-2^{k+2}  &\text{if  } i = 6 \\
5\cdot2^{k+2}  &\text{if  } i =  7\\
-5\cdot2^{k+4}  &\text{if  } i =  8 \\
2^{k+ 6}  &\text{if  } i =  9 \\
0  &\text{if  } i = 10  \\
3\cdot2^{2k +6}  &\text{if  } i =  11 
\end{cases} 
  \end{equation}
 
  We have in the case \underline{ $l = 2, m\geq 1, s = 3 $}\\
 $$A_{i}^{[3,m,2,k]} = \Gamma _{i-3}^{\Big[\substack{m+3 \\ m+5 }\Big] \times k} - 4\cdot \Gamma _{i-4}^{\Big[\substack{m+2 \\ m+5 }\Big] \times k}
   -8\cdot\Gamma _{i-4}^{\Big[\substack{m+3 \\ m+4 }\Big] \times k} +32\cdot\Gamma _{i-5}^{\Big[\substack{m+2 \\ m+4}\Big] \times k}$$
\begin{equation}
\label{eq 7.20}
  A_{i}^{[3,m,2,k]} =  
\begin{cases}
1  &\text{if } i = 3 \\
-3  &\text{if  } i = 4\\
2 &\text{if } i = 5 \\
0  &\text{if  } 6 \leq i\leq m+5 \\
-2^{k+m+2}  &\text{if  } i = m+6 \\
5\cdot2^{k+m+2}  &\text{if  } i = m+ 7\\
-5\cdot2^{k+m+4}  &\text{if  } i = m+ 8 \\
2^{k+m+ 6}  &\text{if  } i = m+ 9 \\
0  &\text{if  } m+ 10 \leq i\leq 2m+10 \\
3\cdot2^{2k+2m +6}  &\text{if  } i = 2m+ 11 
\end{cases}
\end{equation}

   We have in the case \underline{ $l\geq 3, m = 0 , s = 3 $}\\
$$A_{i}^{[3,0,l,k]} = \Gamma _{i-3}^{\Big[\substack{3 \\ 3+l }\Big] \times k} - 4\cdot \Gamma _{i- 4}^{\Big[\substack{ 2 \\ 3+l }\Big] \times k}
   -8\cdot\Gamma _{i- 4}^{\Big[\substack{3 \\ 2+l }\Big] \times k} +32\cdot\Gamma _{i- 5}^{\Big[\substack{2 \\ 2+l}\Big] \times k}$$
  \begin{equation}
  \label{eq 7.21}
  A_{i}^{[3,0,l,k]} =  
\begin{cases}
1  &\text{if } i = 3 \\
-3  &\text{if  } i = 4\\
2 &\text{if } i = 5 \\
-2^{k+2}  &\text{if  } i = 6 \\
5\cdot2^{k+2}  &\text{if  } i = 7\\
3\cdot2^{k+2i-9}  &\text{if  }  8\leq i\leq l +5 \\
-5\cdot2^{k+2l}  &\text{if  } i = 6 +l \\
2^{k+2l +2}  &\text{if  } i = 7 +l \\
0  &\text{if  }  8+l  \\
3\cdot2^{2k+l+4}  &\text{if  } i = 9 +l 
\end{cases}
  \end{equation}
 We have in the case \underline{ $l\geq 3, m\geq 1, s = 3 $}\\
$$A_{i}^{[3,m,l,k]} = \Gamma _{i-3}^{\Big[\substack{m+3 \\ m+3+l }\Big] \times k} - 4\cdot \Gamma _{i- 4}^{\Big[\substack{m+ 2 \\ m+3+l }\Big] \times k}
   -8\cdot\Gamma _{i- 4}^{\Big[\substack{m+3 \\ m+2+l }\Big] \times k} +32\cdot\Gamma _{i- 5}^{\Big[\substack{m+2 \\ m+2+l}\Big] \times k}$$
 \begin{equation}
 \label{eq 7.22}
 A_{i}^{[3,m,l,k]} =  
\begin{cases}
1  &\text{if } i = 3 \\
-3  &\text{if  } i = 4\\
2 &\text{if } i = 5 \\
0  &\text{if  } 6 \leq i\leq m+ 5 \\
-2^{k+m+2}  &\text{if  } i = m+6 \\
5\cdot2^{k+m+2}  &\text{if  } i = m+7\\
3\cdot2^{k-m+2i-9}  &\text{if  } m+ 8\leq i\leq m+l +5 \\
-5\cdot2^{k+m+2l}  &\text{if  } i = m+6 +l \\
2^{k+m+2l +2}  &\text{if  } i = m+7 +l \\
0  &\text{if  } m+8+l\leq i\leq 2m+l +8 \\
3\cdot2^{2k+2m+l+4}  &\text{if  } i = 2m+9 +l 
\end{cases}
\end{equation}
  
     We have in the case \underline{ $l = 0, m\geq 0, s\geq 4 $}\\
 $$A_{i}^{[s,m,0,k]} = \Gamma _{i-s}^{\Big[\substack{m+s \\ m+s }\Big] \times k} - 12\cdot \Gamma _{i-s-1}^{\Big[\substack{m+s-1 \\ m+s }\Big] \times k}
 +32\cdot\Gamma _{i-s-2}^{\Big[\substack{m+s-1 \\ m+s-1}\Big] \times k}$$
 \begin{equation}
 \label{eq 7.23}
    A_{i}^{[s,m,0,k]} =  
\begin{cases}
1  &\text{if } i = s \\
-3  &\text{if  } i = s +1\\
2 &\text{if } i = s +2 \\
0  &\text{if  } s+3\leq i\leq m+2s-1 \\
-3\cdot2^{k+m+s-1}  &\text{if  } i = m+2s \\
3\cdot2^{k+m+s-1}  &\text{if  } i = m+2s +1\\
0  &\text{if  } m+2s+2\leq i\leq 2m+3s -1\\
3\cdot2^{2k+2s+2m -2}  &\text{if  } i = 2m+3s 
\end{cases}       
 \end{equation}

 We have in the case \underline{ $l = 1, m\geq 0, s\geq 4 $}\\
 $$A_{i}^{[s,m,1,k]} = \Gamma _{i-s}^{\Big[\substack{m+s \\ m+s+1 }\Big] \times k} - 4\cdot \Gamma _{i-s-1}^{\Big[\substack{m+s-1 \\ m+s+1 }\Big] \times k}
   -8\cdot\Gamma _{i-s-1}^{\Big[\substack{m+s \\ m+s }\Big] \times k} +32\cdot\Gamma _{i-s-2}^{\Big[\substack{m+s-1 \\ m+s}\Big] \times k}$$
 \begin{equation}
\label{eq 7.24} 
  A_{i}^{[s,m,1,k]} =  
\begin{cases}
1  &\text{if } i = s \\
-3  &\text{if  } i = s +1\\
2 &\text{if } i = s +2 \\
0  &\text{if  } s+3\leq i\leq m+2s-1 \\
-2^{k+m+s-1}  &\text{if  } i = m+2s \\
-3\cdot2^{k+m+s-1}  &\text{if  } i = m+2s +1\\
2^{k+m+s+ 1}  &\text{if  } i = m+2s +2 \\
0  &\text{if  } m+2s+3\leq i\leq 2m+3s \\
3\cdot2^{2k+2s+2m -1}  &\text{if  } i = 2m+3s+1 
\end{cases}     
\end{equation}

 We have in the case \underline{ $l = 2, m\geq 0, s\geq 4 $}\\
 $$A_{i}^{[s,m,2,k]} = \Gamma _{i-s}^{\Big[\substack{m+s \\ m+s+2 }\Big] \times k} - 4\cdot \Gamma _{i-s-1}^{\Big[\substack{m+s-1 \\ m+s+2 }\Big] \times k}
   -8\cdot\Gamma _{i-s-1}^{\Big[\substack{m+s \\ m+s+1 }\Big] \times k} +32\cdot\Gamma _{i-s-2}^{\Big[\substack{m+s-1 \\ m+s+1}\Big] \times k}$$
 \begin{equation}
 \label{eq 7.25}
  A_{i}^{[s,m,2,k]} =  
\begin{cases}
1  &\text{if } i = s \\
-3  &\text{if  } i = s +1\\
2 &\text{if } i = s +2 \\
0  &\text{if  } s+3\leq i\leq m+2s-1 \\
-2^{k+m+s-1}  &\text{if  } i = m+2s \\
5\cdot2^{k+m+s-1}  &\text{if  } i = m+2s +1\\
-5\cdot2^{k+m+s+ 1}  &\text{if  } i = m+2s +2 \\
2^{k+m+s+ 3}  &\text{if  } i = m+2s +3 \\
0  &\text{if  } m+2s+4\leq i\leq 2m+3s+1 \\
3\cdot2^{2k+2s+2m}  &\text{if  } i = 2m+3s+2 
\end{cases} 
   \end{equation}
We have in the case \underline{ $l\geq 3, m\geq 0, s\geq 4 $}\\
$$A_{i}^{[s,m,l,k]} = \Gamma _{i-s}^{\Big[\substack{m+s \\ m+s+l }\Big] \times k} - 4\cdot \Gamma _{i-s-1}^{\Big[\substack{m+s-1 \\ m+s+l }\Big] \times k}
   -8\cdot\Gamma _{i-s-1}^{\Big[\substack{m+s \\ m+s+l-1 }\Big] \times k} +32\cdot\Gamma _{i-s-2}^{\Big[\substack{m+s-1 \\ m+s+l-1}\Big] \times k}$$
  \begin{equation}
\label{eq 7.26}
 A_{i}^{[s,m,l,k]} =  
\begin{cases}
1  &\text{if } i = s \\
-3  &\text{if  } i = s +1\\
2 &\text{if } i = s +2 \\
0  &\text{if  } s+3\leq i\leq m+2s-1 \\
-2^{k+m+s-1}  &\text{if  } i = m+2s \\
5\cdot2^{k+m+s-1}  &\text{if  } i = m+2s +1\\
3\cdot2^{k-m-3s+2i-3}  &\text{if  } m+2s+2\leq i\leq m+2s+l-1 \\
-5\cdot2^{k+m+s+2l-3}  &\text{if  } i = m+2s +l \\
2^{k+m+s+2l-1}  &\text{if  } i = m+2s +l +1 \\
0  &\text{if  } m+2s+l+2\leq i\leq 2m+3s+l-1 \\
3\cdot2^{2k+2s+2m+l-2}  &\text{if  } i = 2m+3s+l 
\end{cases}
\end{equation}
 \end{lem}
\begin{proof}

 By way of example we prove \eqref{eq 7.26}\\
 Let $s\geq 4,\;m\geq 0,\; l\geq 4,\; i = s+j,\;\text{for}\;0\leq j\leq 2m+2s+l $\\
  We can write\\
    \begin{equation}
 \label{eq 7.27}   
   A_{s+j}^{[s,m,l,k]} = \Gamma _{j}^{\Big[\substack{m+s \\ m+s+l }\Big] \times k} - 4\cdot \Gamma _{j-1}^{\Big[\substack{m+s-1 \\ m+s+l }\Big] \times k}
   -8\cdot\Gamma _{j-1}^{\Big[\substack{m+s \\ m+s+l-1 }\Big] \times k} +32\cdot\Gamma _{j-2}^{\Big[\substack{m+s-1 \\ m+s+l-1}\Big] \times k} \;\text{for}\;0\leq j\leq 2m+2s+l 
 \end{equation}
 From Theorem 3.13  (see section 3 in [2]) we obtain \\
 
\begin{align*}
& \underline{j = 0}\quad  A_{s}^{[s,m,l,k]} = \Gamma _{0}^{\Big[\substack{m+s \\ m+s+l }\Big] \times k} = 1 \\
& \underline{j = 1}\quad  A_{s+1}^{[s,m,l,k]}  = \Gamma _{1}^{\Big[\substack{m+s \\ m+s+l }\Big] \times k} - 4\cdot \Gamma _{0}^{\Big[\substack{m+s-1 \\ m+s+l }\Big] \times k}
   -8\cdot\Gamma _{0}^{\Big[\substack{m+s \\ m+s+l-1 }\Big] \times k}= 9-12=-3 \\
& \underline{j = 2} \quad  A_{s+2}^{[s,m,l,k]} =    \Gamma _{2}^{\Big[\substack{m+s \\ m+s+l }\Big] \times k} - 4\cdot \Gamma _{1}^{\Big[\substack{m+s-1 \\ m+s+l }\Big] \times k}
   -8\cdot\Gamma _{1}^{\Big[\substack{m+s \\ m+s+l-1 }\Big] \times k} +32\cdot\Gamma _{0}^{\Big[\substack{m+s-1 \\ m+s+l-1}\Big] \times k} 
   = 78-12\cdot9+32=2 \\
   & \\
  & \underline{3 \leq j \leq m+s-1 }\\
  &   A_{s+j}^{[s,m,l,k]} = \Gamma _{j}^{\Big[\substack{m+s \\ m+s+l }\Big] \times k} - 4\cdot \Gamma _{j-1}^{\Big[\substack{m+s-1 \\ m+s+l }\Big] \times k}
   -8\cdot\Gamma _{j-1}^{\Big[\substack{m+s \\ m+s+l-1 }\Big] \times k} +32\cdot\Gamma _{j-2}^{\Big[\substack{m+s-1 \\ m+s+l-1}\Big] \times k} \\
   & \\
    & = \big(21\cdot2^{3j-4} -3\cdot2^{2j-3}\big) -12\cdot\big(21\cdot2^{3j-7} -3\cdot2^{2j-5}\big) +32\cdot\big(21\cdot2^{3j-10} -3\cdot2^{2j-7}\big) = 0 \\
     & \underline{j = m+s }\\
  &   A_{2s+m}^{[s,m,l,k]} = \Gamma _{m+s}^{\Big[\substack{m+s \\ m+s+l }\Big] \times k} - 4\cdot \Gamma _{m+s-1}^{\Big[\substack{m+s-1 \\ m+s+l }\Big] \times k}
   -8\cdot\Gamma _{m+s-1}^{\Big[\substack{m+s \\ m+s+l-1 }\Big] \times k} +32\cdot\Gamma _{m+s-2}^{\Big[\substack{m+s-1 \\ m+s+l-1}\Big] \times k} \\
   & \\
  & = \big(2^{k+m+s-1} +21\cdot2^{3(m+s)-4} -11\cdot2^{2(m+s)-3}\big)
   -4\cdot\big(2^{k+m+s-2} +21\cdot2^{3(m+s-1)-4} -11\cdot2^{2(m+s-1)-3}\big) \\
    &   -8\cdot\big(21\cdot2^{3(m+s-1)-4} -3\cdot2^{2(m+s-1)-3}\big) 
       +32\cdot\big(21\cdot2^{3(m+s-2)-4} -3\cdot2^{2(m+s-2)-3}\big) = - 2^{k+m+s-1}\\
   & \underline{j = m+s+1}\\ 
      &   A_{2s+m+1}^{[s,m,l,k]} = \Gamma _{m+s+1}^{\Big[\substack{m+s \\ m+s+l }\Big] \times k} - 4\cdot \Gamma _{m+s}^{\Big[\substack{m+s-1 \\ m+s+l }\Big] \times k}
   -8\cdot\Gamma _{m+s}^{\Big[\substack{m+s \\ m+s+l-1 }\Big] \times k} +32\cdot\Gamma _{m+s-1}^{\Big[\substack{m+s-1 \\ m+s+l-1}\Big] \times k} \\  
   & \\
    & = \big(3\cdot2^{k+m+s-1} +21\cdot(2^{3m+3s-1} -2^{2m+2s-1})\big)
    -\big(3\cdot2^{k+m+s} +21\cdot(2^{3m+3s-2} -2^{2m+2s-1})\big)\\
  &  -\big(2^{k+m+s+2} +21\cdot2^{3m+3s-1} -11\cdot2^{2m+2s}\big)
    +\big(2^{k+m+s+3} +21\cdot2^{3m+3s-2} -11\cdot2^{2m+2s}\big) = 5\cdot2^{k+m+s-1} \\
    & \underline{j = m+s+2}\\ 
      &   A_{2s+m+2}^{[s,m,l,k]} = \Gamma _{m+s+2}^{\Big[\substack{m+s \\ m+s+l }\Big] \times k} - 4\cdot \Gamma _{m+s+1}^{\Big[\substack{m+s-1 \\ m+s+l }\Big] \times k}
   -8\cdot\Gamma _{m+s+1}^{\Big[\substack{m+s \\ m+s+l-1 }\Big] \times k} +32\cdot\Gamma _{m+s}^{\Big[\substack{m+s-1 \\ m+s+l-1}\Big] \times k} \\  
   & \\
    & = \big(3\cdot2^{k+m+s+1} +21\cdot(2^{3m+3s+2} -2^{2m+2s+2})\big)
    -\big(3\cdot2^{k+m+s+2} +21\cdot(2^{3m+3s+1} -2^{2m+2s+2})\big)\\
  &   -\big(3\cdot2^{k+m+s+2} +21\cdot(2^{3m+3s+2} -2^{2m+2s+1})\big)
  +   \big(3\cdot2^{k+m+s+3} +21\cdot(2^{3m+3s+1} -2^{2m+2s+1})\big) = 3\cdot2^{k+m+s+1}\\
  & \underline{m+s+2\leq j\leq m+s+l-1}\\
 &  A_{s+j}^{[s,m,l,k]} = \Gamma _{j}^{\Big[\substack{m+s \\ m+s+l }\Big] \times k} - 4\cdot \Gamma _{j-1}^{\Big[\substack{m+s-1 \\ m+s+l }\Big] \times k}
   -8\cdot\Gamma _{j-1}^{\Big[\substack{m+s \\ m+s+l-1 }\Big] \times k} +32\cdot\Gamma _{j-2}^{\Big[\substack{m+s-1 \\ m+s+l-1}\Big] \times k} \\
    & \\
    & = \big(3\cdot2^{k-m-s+2j-3} +21\cdot(2^{3j-4} -2^{-m-s+3j-4})\big)
      -\big(3\cdot2^{k-m-s+2j-2} +21\cdot(2^{3j-5} -2^{-m-s+3j-4})\big)\\
     &   -\big(3\cdot2^{k-m-s+2j-2} +21\cdot(2^{3j-4} -2^{-m-s+3j-4})\big)
    +   \big(3\cdot2^{k-m-s+2j-1} +21\cdot(2^{3j-5} -2^{-m-s+3j-4})\big) \\
    & = 3\cdot2^{k-m-s+2j-3}\\
     & \underline{j = m+s+l}\\
 &  A_{m+2s+l}^{[s,m,l,k]} = \Gamma _{m+s+l}^{\Big[\substack{m+s \\ m+s+l }\Big] \times k} - 4\cdot \Gamma _{m+s+l-1}^{\Big[\substack{m+s-1 \\ m+s+l }\Big] \times k}
   -8\cdot\Gamma _{m+s+l-1}^{\Big[\substack{m+s \\ m+s+l-1 }\Big] \times k} +32\cdot\Gamma _{m+s+l-2}^{\Big[\substack{m+s-1 \\ m+s+l-1}\Big] \times k} \\
    & \\
    & = \big(11\cdot2^{k+m+s+2l-3} +21\cdot2^{3m+3s+3l-4} -53\cdot2^{2m+2s+3l-4}\big)
      -\big(3\cdot2^{k+m+s+2l-2} +21\cdot(2^{3m+3s+3l-5} -2^{2m+2s+3l-4})\big)\\
       &   -\big(11\cdot2^{k+m+s+2l-2} +21\cdot2^{3m+3s+3l-4} -53\cdot2^{2m+2s+3l-4}\big)
      +  \big(3\cdot2^{k+m+s+2l-1} +21\cdot(2^{3m+3s+3l-5} -2^{2m+2s+3l-4})\big) \\
       & = -5\cdot2^{k+m+s+2l-3}\\
        & \underline{j = m+s+l+1}\\
 &  A_{m+2s+l+1}^{[s,m,l,k]} = \Gamma _{m+s+l+1}^{\Big[\substack{m+s \\ m+s+l }\Big] \times k} - 4\cdot \Gamma _{m+s+l}^{\Big[\substack{m+s-1 \\ m+s+l }\Big] \times k}
   -8\cdot\Gamma _{m+s+l}^{\Big[\substack{m+s \\ m+s+l-1 }\Big] \times k} +32\cdot\Gamma _{m+s+l-1}^{\Big[\substack{m+s-1 \\ m+s+l-1}\Big] \times k} \\
    & \\
    & = \big(21\cdot(2^{k+m+s+2l-1} +2^{3m+3s+3l-1} -5\cdot2^{2m+2s+3l-1})\big)
      -\big(11\cdot2^{k+m+s+2l} +21\cdot2^{3m+3s+3l-2} -53\cdot2^{2m+2s+3l-1}\big)\\
       &  - \big(21\cdot(2^{k+m+s+2l} +2^{3m+3s+3l-1} -5\cdot2^{2m+2s+3l-1})\big)
        + \big(11\cdot2^{k+m+s+2l+1} +21\cdot2^{3m+3s+3l-2} -53\cdot2^{2m+2s+3l-1}\big) \\
       & = 2^{k+m+s+2l-1}\\  
         & \underline{m+s+l+2\leq j\leq 2m+2s+l-2}\\
 &  A_{s+j}^{[s,m,l,k]} = \Gamma _{j}^{\Big[\substack{m+s \\ m+s+l }\Big] \times k} - 4\cdot \Gamma _{j-1}^{\Big[\substack{m+s-1 \\ m+s+l }\Big] \times k}
   -8\cdot\Gamma _{j-1}^{\Big[\substack{m+s \\ m+s+l-1 }\Big] \times k} +32\cdot\Gamma _{j-2}^{\Big[\substack{m+s-1 \\ m+s+l-1}\Big] \times k} \\
    & \\
    & = \big(21\cdot(2^{k-2m-2s-l+3j-4} +2^{3j-4} -5\cdot2^{4j-2s-2m-l-5})\big)
     - \big(21\cdot(2^{k-2m-2s-l+3j-4} +2^{3j-5} -5\cdot2^{4j-2s-2m-l-6})\big)\\
       &  - \big(21\cdot(2^{k-2m-2s-l+3j-3} +2^{3j-4} -5\cdot2^{4j-2s-2m-l-5})\big)
         +  \big(21\cdot(2^{k-2m-2s-l+3j-3} +2^{3j-5} -5\cdot2^{4j-2s-2m-l-6})\big) \\
       & = 0 \\ 
     & \underline{j = 2m+2s+l}\\
 &  A_{2m+3s+l}^{[s,m,l,k]} = \Gamma _{2m+2s+l}^{\Big[\substack{m+s \\ m+s+l }\Big] \times k} - 4\cdot \Gamma _{2m+2s+l-1}^{\Big[\substack{m+s-1 \\ m+s+l }\Big] \times k}
   -8\cdot\Gamma _{2m+2s+l-1}^{\Big[\substack{m+s \\ m+s+l-1 }\Big] \times k} +32\cdot\Gamma _{2m+2s+l-2}^{\Big[\substack{m+s-1 \\ m+s+l-1}\Big] \times k} \\
    & \\
    & = \big(2^{2k+2s+2m+l-2} -3\cdot2^{k+4s+4m+2l-4} +2^{6s+6m+3l-5})\big)
     - \big(2^{2k+2s+2m+l-1} -3\cdot2^{k+4s+4m+2l-4} +2^{6s+6m+3l-6})\big)\\
      &  - \big(2^{2k+2s+2m+l} -3\cdot2^{k+4s+4m+2l-3} +2^{6s+6m+3l-5})\big)
        +  \big(2^{2k+2s+2m+l+1} -3\cdot2^{k+4s+4m+2l-3} +2^{6s+6m+3l-6})\big) \\
         & =3\cdot2^{2k+2s+2m+l-2}    
         \end{align*}
 
\end{proof} 

 \subsection{\textbf{Computation of $ \omega _{i}(s,m,l,k+1) - \omega _{i}(s,m,l,k) $}}
  \label{subsec 5}
   \begin{lem}
\label{lem  7.4} 
 We have \\
 
 \huge The case s = 2 \vspace{0.1 cm}\\
\normalsize
$\underline{The\; case \; l = 0, m =  0,  s = 2, k>i}$ \\
\begin{equation}
\label{eq 7.28}
\omega _{i}(2,0,0,k+1) - \omega _{i}(2,0,0,k)  = 
\begin{cases}
2^{k+1} &\text{if } i = 2 \\
-3\cdot2^{k+1}  &\text{if  } i = 3\\
-9\cdot2^{2k+ 2}+ 2^{k+2}  &\text{if  } i = 4 \\
9\cdot2^{2k+2}  &\text{if  } i = 5 \\
21\cdot2^{3k+3}  &\text{if  } i = 6
\end{cases}
\end{equation}

$\underline{The\; case \; l = 0, m =  1,  s = 2, k>i}$\\
\begin{equation}
\label{eq 7.29}
\omega _{i}(2,1,0,k+1) - \omega _{i}(2,1,0,k)  = 
\begin{cases}
2^{k+1} &\text{if } i = 2 \\
-3\cdot2^{k+1}  &\text{if  } i = 3\\
 2^{k+2}  &\text{if  } i = 4 \\
-9\cdot2^{2k+ 3}  &\text{if  } i = 5 \\
9\cdot2^{2k+3}  &\text{if  } i = 6 \\
0   &\text{if  } i = 7 \\
21\cdot2^{3k+5}  &\text{if  } i = 8
\end{cases}
\end{equation}

$\underline{The\; case \; l = 0, m \geq 2,  s = 2, k>i}$\\
\begin{equation}
\label{eq 7.30}
\omega _{i}(2,m,0,k+1) - \omega _{i}(2,m,0,k)  = 
\begin{cases}
2^{k+1} &\text{if } i = 2 \\
-3\cdot2^{k+1}  &\text{if  } i = 3\\
 2^{k+2}  &\text{if  } i = 4 \\
 0  &\text{if  } 5\leq i\leq m+3 \\
-9\cdot2^{2k+ m+2}  &\text{if  } i = m+4 \\
9\cdot2^{2k+m+2}  &\text{if  } i = m+5 \\
0   &\text{if  } m+7\leq i\leq 2m+5 \\
21\cdot2^{3k+2m+3}  &\text{if  } i = 2m+6
\end{cases}
\end{equation}

$\underline{The\; case \; l = 1, m =  0,  s = 2, k>i}$\\
\begin{equation}
\label{eq 7.31}
\omega _{i}(2,0,1,k+1) - \omega _{i}(2,0,1,k)  = 
\begin{cases}
2^{k+1} &\text{if } i = 2 \\
-3\cdot2^{k+1}  &\text{if  } i = 3\\
-3\cdot 2^{2k+2} +2^{k+2}  &\text{if  } i = 4 \\
-9\cdot2^{2k+ 2}  &\text{if  } i = 5 \\
3\cdot2^{2k+ 4}  &\text{if  } i = 6 \\
21\cdot2^{3k+4}  &\text{if  } i = 7
\end{cases}
\end{equation}
$\underline{The\; case \; l = 1, m =  1,  s = 2, k>i}$\\
\begin{equation}
\label{eq 7.32}
\omega _{i}(2,1,1,k+1) - \omega _{i}(2,1,1,k)  = 
\begin{cases}
2^{k+1} &\text{if } i = 2 \\
-3\cdot2^{k+1}  &\text{if  } i = 3\\
2^{k+2}  &\text{if  } i = 4\\
-3\cdot 2^{2k+3}   &\text{if  } i = 5 \\
-9\cdot2^{2k+ 3}  &\text{if  } i = 6 \\
3\cdot2^{2k+ 5}  &\text{if  } i = 7 \\
0 &\text{if  } i = 8 \\
21\cdot2^{3k+6}  &\text{if  } i = 9
\end{cases}
\end{equation}

$\underline{The\; case \; l = 1, m \geq 2,  s = 2, k>i}$\\
\begin{equation}
\label{eq 7.33}
\omega _{i}(2,m,1,k+1) - \omega _{i}(2,m,1,k)  = 
\begin{cases}
2^{k+1} &\text{if } i = 2 \\
-3\cdot2^{k+1}  &\text{if  } i = 3\\
 2^{k+2}  &\text{if  } i = 4 \\
 0  &\text{if  } 5\leq i\leq m+3 \\
3\cdot2^{2k+ m+2}  &\text{if  } i = m+4 \\
- 9\cdot2^{2k+m+2}  &\text{if  } i = m+5 \\
3\cdot2^{2k+ m+4}  &\text{if  } i = m+6 \\
0   &\text{if  } m+7\leq i\leq 2m+6 \\
21\cdot2^{3k+2m+4}  &\text{if  } i = 2m+7
\end{cases}
\end{equation}

$\underline{The\; case \; l = 2, m =  0,  s = 2, k>i}$\\
\begin{equation}
\label{eq 7.34}
\omega _{i}(2,0,2,k+1) - \omega _{i}(2,0,2,k)  = 
\begin{cases}
2^{k+1} &\text{if } i = 2 \\
-3\cdot2^{k+1}  &\text{if  } i = 3\\
  -3\cdot 2^{2k+2} + 2^{k+2}  &\text{if  } i = 4\\
15\cdot 2^{2k+2}   &\text{if  } i = 5 \\
-15\cdot2^{2k+ 4}  &\text{if  } i = 6 \\
3\cdot2^{2k+ 6}  &\text{if  } i = 7 \\
21\cdot2^{3k+5}  &\text{if  } i = 8
\end{cases}
\end{equation}
$\underline{The\; case \; l = 2, m =  1,  s = 2, k>i}$\\
\begin{equation}
\label{eq 7.35}
\omega _{i}(2,1,2,k+1) - \omega _{i}(2,1,2,k)  = 
\begin{cases}
2^{k+1} &\text{if } i = 2 \\
-3\cdot2^{k+1}  &\text{if  } i = 3\\
 2^{k+2}  &\text{if  } i = 4\\
  -3\cdot 2^{2k+3}    &\text{if  } i = 5\\
15\cdot 2^{2k+3}   &\text{if  } i = 6 \\
-15\cdot2^{2k+ 5}  &\text{if  } i = 7 \\
3\cdot2^{2k+ 7}  &\text{if  } i = 8 \\
0  &\text{if  } i = 9 \\
21\cdot2^{3k+7}  &\text{if  } i = 10
\end{cases}
\end{equation}

$\underline{The\; case \; l = 2, m \geq 2,  s = 2, k>i}$\\
\begin{equation}
\label{eq 7.36}
\omega _{i}(2,m,2,k+1) - \omega _{i}(2,m,2,k)  = 
\begin{cases}
2^{k+1} &\text{if } i = 2 \\
-3\cdot2^{k+1}  &\text{if  } i = 3\\
 2^{k+2}  &\text{if  } i = 4\\
 0   &\text{if  } 5\leq i\leq m+3  \\
  -3\cdot 2^{2k+m+2}    &\text{if  } i = m+4 \\
15\cdot 2^{2k+m+2}   &\text{if  } i = m+5  \\
-15\cdot2^{2k+ m+4}  &\text{if  } i = m+6  \\
3\cdot2^{2k+ m+6}  &\text{if  } i =  m+7\\
0   &\text{if  } m+8\leq i\leq 2m+7 \\
21\cdot2^{3k+2m +5}  &\text{if  } i = 2m+8
\end{cases}
\end{equation}

$\underline{The\; case \; l \geq  3, m = 0,  s = 2, k>i}$\\
\begin{equation}
\label{eq 7.37}
\omega _{i}(2,0,l,k+1) - \omega _{i}(2,0,l,k)  = 
\begin{cases}
2^{k+1} &\text{if } i = 2 \\
-3\cdot2^{k+1}  &\text{if  } i = 3\\
-3\cdot2^{2k +2} + 2^{k+2}  &\text{if  } i = 4\\
   15\cdot 2^{2k+2}    &\text{if  } i = 5 \\
9\cdot 2^{2k+ 2i-8}   &\text{if  } 6\leq i\leq l+3  \\
-15\cdot2^{2k+ 2l}  &\text{if  } i = l+4 \\
3\cdot2^{2k+ 2l +2}  &\text{if  } i =  l+5 \\
21\cdot2^{3k+l+3}  &\text{if  } i = l+6
\end{cases}
\end{equation}

$\underline{The\; case \; l \geq  3, m = 1,  s = 2, k>i}$\\
\begin{equation}
\label{eq 7.38}
\omega _{i}(2,1,l,k+1) - \omega _{i}(2,1,l,k)  = 
\begin{cases}
2^{k+1} &\text{if } i = 2 \\
-3\cdot2^{k+1}  &\text{if  } i = 3\\
2^{k+2}  &\text{if  } i = 4\\
-3\cdot2^{2k +3}   &\text{if  } i = 5\\
   15\cdot 2^{2k+3}    &\text{if  } i = 6 \\
9\cdot 2^{2k+ 2i-9}   &\text{if  } 7\leq i\leq l+4  \\
-15\cdot2^{2k+ 2l +1}  &\text{if  } i = l+5 \\
3\cdot2^{2k+ 2l +3}  &\text{if  } i =  l+6 \\
0 &\text{if  } i =  l+7 \\
21\cdot2^{3k+l+5}  &\text{if  } i = l+8
\end{cases}
\end{equation}

$\underline{The\; case \; l \geq  3, m \geq 2,  s = 2, k>i}$\\
\begin{equation}
\label{eq 7.39}
\omega _{i}(2,m,l,k+1) - \omega _{i}(2,m,l,k)  = 
\begin{cases}
2^{k+1} &\text{if } i = 2 \\
-3\cdot2^{k+1}  &\text{if  } i = 3\\
2^{k+2}  &\text{if  } i = 4\\
0  &\text{if  } 5\leq i\leq m+3 \\
-3\cdot2^{2k + m+2}   &\text{if  } i = m+4 \\
   15\cdot 2^{2k+m+2}    &\text{if  } i = m+5 \\
9\cdot 2^{2k-m + 2i-8}   &\text{if  } m+6  \leq i\leq m+l+3  \\
-15\cdot2^{2k+ m+2l}  &\text{if  } i = m+l+4 \\
3\cdot2^{2k+m + 2l + 2}  &\text{if  } i =  m+l+5 \\
0 &\text{if  } m+l+6\leq i\leq 2m+l+5 \\
21\cdot2^{3k +2m +l+3}  &\text{if  } i = 2m+l+6
\end{cases}
\end{equation}

\huge The case s = 3 \vspace{0.1 cm}\\
\normalsize
$\underline{The\; case \; l = 0, m =  0,  s = 3, k>i}$ \\
\begin{equation}
\label{eq 7.40}
\omega _{i}(3,0,0,k+1) - \omega _{i}(3,0,0,k)  = 
\begin{cases}
2^{k+2} &\text{if } i = 3 \\
-3\cdot2^{k+2}  &\text{if  } i = 4\\
2^{k+3}  &\text{if  } i = 5 \\
-9\cdot2^{2k+ 4}  &\text{if  } i = 6 \\
9\cdot2^{2k+4}  &\text{if  } i = 7 \\
0  &\text{if  } i = 8 \\
21\cdot2^{3k+6}  &\text{if  } i = 9
\end{cases}
\end{equation}

$\underline{The\; case \; l = 0, m \geq 1,  s = 3, k>i}$ \\
\begin{equation}
\label{eq 7.41}
\omega _{i}(3,m,0,k+1) - \omega _{i}(3,m,0,k)  = 
\begin{cases}
2^{k+2} &\text{if } i = 3 \\
-3\cdot2^{k+2}  &\text{if  } i = 4\\
2^{k+3}  &\text{if  } i = 5 \\
0  &\text{if  } 6\leq i\leq m+5 \\
-9\cdot2^{2k+m+ 4}  &\text{if  } i = m+ 6 \\
9\cdot2^{2k+m+4}  &\text{if  } i = m+ 7 \\
0  &\text{if  } m+8\leq i\leq 2m+8 \\
21\cdot2^{3k+2m+6}  &\text{if  } i = 2m+9
\end{cases}
\end{equation}

$\underline{The\; case \; l = 1, m = 0,  s = 3, k>i}$ \\
\begin{equation}
\label{eq 7.42} 
\omega _{i}(3,0,1,k+1) - \omega _{i}(3,0,1,k)  = 
\begin{cases}
2^{k+2} &\text{if } i = 3 \\
-3\cdot2^{k+2}  &\text{if  } i = 4\\
2^{k+3}  &\text{if  } i = 5 \\
-3\cdot2^{2k+4} &\text{if  } i = 6 \\
-9\cdot2^{2k+ 4}  &\text{if  } i = 7 \\
3\cdot2^{2k+6}  &\text{if  } i = 8 \\
0  &\text{if  } i = 9 \\
21\cdot2^{3k+7}  &\text{if  } i = 10
\end{cases}
\end{equation}

$\underline{The\; case \; l = 1, m \geq 1,  s = 3, k>i}$ \\
\begin{equation}
\label{eq 7.43} 
\omega _{i}(3,m,1,k+1) - \omega _{i}(3,m,1,k)  = 
\begin{cases}
2^{k+2} &\text{if } i = 3 \\
-3\cdot2^{k+2}  &\text{if  } i = 4\\
2^{k+3}  &\text{if  } i = 5 \\
0  &\text{if  } 6\leq i\leq m+5 \\
-3\cdot2^{2k+m+4} &\text{if  } i = m+ 6 \\
-9\cdot2^{2k+m+ 4}  &\text{if  } i = m+ 7 \\
3\cdot2^{2k+m+6}  &\text{if  } i = m+8 \\
0  &\text{if  } m+9\leq i\leq 2m +9\\
21\cdot2^{3k+2m+7}  &\text{if  } i = 2m+ 10
\end{cases}
\end{equation}

$\underline{The\; case \; l = 2, m = 0,  s = 3, k>i}$ \\
\begin{equation}
\label{eq 7.44}
\omega _{i}(3,0,2,k+1) - \omega _{i}(3,0,2,k)  = 
\begin{cases}
2^{k+2} &\text{if } i = 3 \\
-3\cdot2^{k+2}  &\text{if  } i = 4\\
2^{k+3}  &\text{if  } i = 5 \\
-3\cdot2^{2k+4} &\text{if  } i =  6 \\
15\cdot2^{2k+ 4}  &\text{if  } i =  7 \\
-15\cdot2^{2k+6}  &\text{if  } i = 8 \\
3\cdot2^{2k+8}  &\text{if  } i = 9 \\
0  &\text{if  } i = 10 \\
21\cdot2^{3k+ 8}  &\text{if  } i = 11
\end{cases}
\end{equation}

$\underline{The\; case \; l = 2, m \geq 1,  s = 3, k>i}$ \\
\begin{equation}
\label{eq 7.45}
\omega _{i}(3,m,2,k+1) - \omega _{i}(3,m,2,k)  = 
\begin{cases}
2^{k+2} &\text{if } i = 3 \\
-3\cdot2^{k+2}  &\text{if  } i = 4\\
2^{k+3}  &\text{if  } i = 5 \\
0  &\text{if  } 6\leq i\leq m+5 \\
-3\cdot2^{2k+m+4} &\text{if  } i = m+ 6 \\
15\cdot2^{2k+m+ 4}  &\text{if  } i = m+ 7 \\
-15\cdot2^{2k+m+6}  &\text{if  } i =m+ 8 \\
3\cdot2^{2k+m+8}  &\text{if  } i =m+ 9 \\
0  &\text{if  } m+10\leq i\leq 2m+10 \\
21\cdot2^{3k+2m+ 8}  &\text{if  } i = 2m+11
\end{cases}
\end{equation}

$\underline{The\; case \; l \geq 3, m = 0,  s = 3, k>i}$ \\
\begin{equation}
\label{eq 7.46}
\omega _{i}(3,0,l,k+1) - \omega _{i}(3,0,l,k)  = 
\begin{cases}
2^{k+2} &\text{if } i = 3 \\
-3\cdot2^{k+2}  &\text{if  } i = 4\\
2^{k+3}  &\text{if  } i = 5 \\
- 3\cdot2^{2k+4}  &\text{if  } i = 6 \\
15\cdot2^{2k+ 4}  &\text{if  } i =  7 \\
9\cdot2^{2k+2i- 10} &\text{if  } 8\leq i\leq l+5 \\
-15\cdot2^{2k+2l+2}  &\text{if  } i = l+6 \\
3\cdot2^{2k+2l+4}  &\text{if  } i = l+7 \\
0  &\text{if  } i = l+8 \\
21\cdot2^{3k+ l+6}  &\text{if  } i = l+9
\end{cases}
\end{equation}

$\underline{The\; case \; l \geq 3, m \geq 1,  s = 3, k>i}$ \\
\begin{equation}
\label{eq 7.47}
\omega _{i}(3,m,l,k+1) - \omega _{i}(3,m,l,k)  = 
\begin{cases}
2^{k+2} &\text{if } i = 3 \\
-3\cdot2^{k+2}  &\text{if  } i = 4\\
2^{k+3}  &\text{if  } i = 5 \\
0 &\text{if  }6\leq i\leq m+5 \\
-3\cdot2^{2k+ m+4}  &\text{if  } i = m+6 \\
15\cdot2^{2k+m+ 4}  &\text{if  } i = m+ 7 \\
9\cdot2^{2k-m+2i-10} &\text{if  }m+ 8\leq i\leq m+ l+5 \\
-15\cdot2^{2k+m+2l+2}  &\text{if  } i = m+ l+6 \\
3\cdot2^{2k+m+2l+4}  &\text{if  } i = m+ l+7 \\
0  &\text{if  } m+l+8\leq i\leq 2m+l+8 \\
21\cdot2^{3k+2m+ l+6}  &\text{if  } i =2m+ l+9
\end{cases}
\end{equation}

\huge The case $ s \geq 4  $ \vspace{0.1 cm}\\
\normalsize

$\underline{The\; case \; l = 0, m\geq 0,  s\geq 4, k>i}$\\
\begin{equation}
\label{eq 7.48}
\omega _{i}(s,m,0,k+1) - \omega _{i}(s,m,0,k)  = 
\begin{cases}
2^{k+s-1} &\text{if } i = s \\
-3\cdot2^{k+s-1}  &\text{if  } i = s +1\\
2^{k+s} &\text{if } i = s +2 \\
0  &\text{if  } s+3\leq i\leq m+2s-1 \\
-9\cdot2^{2k+2s +m-2}  &\text{if  } i = m+2s \\
9\cdot2^{2k+2s +m-2}  &\text{if  } i = m+2s +1\\
0  &\text{if  } m+2s+2\leq i\leq 2m+3s-1 \\
21\cdot2^{3k+3s+2m-3}  &\text{if  } i = 2m+3s
\end{cases}
\end{equation}

$\underline{The\; case \; l = 1, m\geq 0,  s\geq 4, k>i}$\\
\begin{equation}
\label{eq 7.49}
\omega _{i}(s,m,1,k+1) - \omega _{i}(s,m,1,k)  = 
\begin{cases}
2^{k+s-1} &\text{if } i = s \\
-3\cdot2^{k+s-1}  &\text{if  } i = s +1\\
2^{k+s} &\text{if } i = s +2 \\
0  &\text{if  } s+3\leq i\leq m+2s-1 \\
-3\cdot2^{2k+2s +m-2}  &\text{if  } i = m+2s \\
- 9\cdot2^{2k+2s +m-2}  &\text{if  } i = m+2s +1\\
3\cdot2^{2k+2s +m}  &\text{if  } i = m+2s +2\\
0  &\text{if  } m+2s+3\leq i\leq 2m+3s \\
21\cdot2^{3k+3s+2m-2}  &\text{if  } i = 2m+3s +1
\end{cases}
\end{equation}

$\underline{The\; case \; l = 2, m\geq 0,  s\geq 4, k>i}$\\
\begin{equation}
\label{eq 7.50}
\omega _{i}(s,m,2,k+1) - \omega _{i}(s,m,2,k)  = 
\begin{cases}
2^{k+s-1} &\text{if } i = s \\
-3\cdot2^{k+s-1}  &\text{if  } i = s +1\\
2^{k+s} &\text{if } i = s +2 \\
0  &\text{if  } s+3\leq i\leq m+2s-1 \\
-3\cdot2^{2k+2s +m-2}  &\text{if  } i = m+2s \\
15\cdot2^{2k+2s +m-2}  &\text{if  } i = m+2s +1\\
-15\cdot2^{2k+2s +m}  &\text{if  } i = m+2s +2\\
3\cdot2^{2k+2s +m +2}  &\text{if  } i = m+2s +3 \\
0  &\text{if  } m+2s+4\leq i\leq 2m+3s +1\\
21\cdot2^{3k+3s+2m-1}  &\text{if  } i = 2m+3s +2
\end{cases}
\end{equation}

$\underline{The\; case \; l\geq 3, m\geq 0,  s\geq 4, k>i}$\\
\begin{equation}
\label{eq 7.51}
\omega _{i}(s,m,l,k+1) - \omega _{i}(s,m,l,k)  = 
\begin{cases}
2^{k+s-1} &\text{if } i = s \\
-3\cdot2^{k+s-1}  &\text{if  } i = s +1\\
2^{k+s} &\text{if } i = s +2 \\
0  &\text{if  } s+3\leq i\leq m+2s-1 \\
-3\cdot2^{2k+2s +m-2}  &\text{if  } i = m+2s \\
15\cdot2^{2k+2s +m-2}  &\text{if  } i = m+2s +1\\
9\cdot2^{2k -2s-m+2i-4}  &\text{if  } m+2s+2\leq i\leq m+2s+l-1 \\
-15\cdot2^{2k+2s +m +2l-4}  &\text{if  } i = m+2s +l \\
3\cdot2^{2k+2s +m +2l-2}  &\text{if  } i = m+2s +l +1 \\
0  &\text{if  } m+2s+l+2\leq i\leq 2m+3s+l-1 \\
21\cdot2^{3k+3s+2m+l-3}  &\text{if  } i = 2m+3s+l 
\end{cases}
\end{equation}

 \end{lem}
  \begin{proof}
 By way of example we prove \eqref{eq 7.51}
 \begin{align*}
 & \text{From  \eqref{eq 7.2} we have}\\
&  \omega _{i}(s,m,l,k+1) - \omega _{i}(s,m,l,k) =  2^{k+s} A_{i}^{[s,m,l,k+1]} - 2^{k+s-1} A_{i}^{[s,m,l,k]} \\
& \text{We now deduce from \eqref{eq 7.26}}\\
& \underline{i = s}\\
&  \omega _{s}(s,m,l,k+1) - \omega _{s}(s,m,l,k) = 2^{k+s}\cdot1- 2^{k+s-1}\cdot1=2^{k+s-1}\\
& \underline{i = s+1}\\
&  \omega _{s+1}(s,m,l,k+1) - \omega _{s+1}(s,m,l,k) = 2^{k+s}\cdot(-3)- 2^{k+s-1}\cdot(-3)=-3\cdot2^{k+s-1}\\
& \underline{i = s+2}\\
&  \omega _{s+2}(s,m,l,k+1) - \omega _{s+2}(s,m,l,k) = 2^{k+s}\cdot2- 2^{k+s-1}\cdot2=2^{k+s}\\
& \underline{s+3\leq i\leq m+2s-1}\\
&  \omega _{i}(s,m,l,k+1) - \omega _{i}(s,m,l,k) = 2^{k+s}\cdot 0- 2^{k+s-1}\cdot 0=0\\
& \underline{i= m+2s}\\
&  \omega _{m+2s}(s,m,l,k+1) - \omega _{m+2s}(s,m,l,k) = 2^{k+s}\cdot (-2^{k+m+s})- 2^{k+s-1}\cdot (-2^{k+m+s-1})= -3\cdot2^{2k+m+2s-2} \\
& \underline{i= m+2s+1}\\
&  \omega _{m+2s+1}(s,m,l,k+1) - \omega _{m+2s+1}(s,m,l,k) = 2^{k+s}\cdot (5\cdot2^{k+m+s})- 2^{k+s-1}\cdot (5\cdot2^{k+m+s-1})= 15\cdot2^{2k+m+2s-2} \\
& \underline{ m+2s+2\leq i\leq m+2s+l-1}\\
&  \omega _{i}(s,m,l,k+1) - \omega _{i}(s,m,l,k) = 2^{k+s}\cdot (3\cdot2^{k-m-3s+2i-2})- 2^{k+s-1}\cdot (3\cdot2^{k-m-3s+2i-3})= 9\cdot2^{2k-2s-m+2i-4} \\
& \underline{ i= m+2s+l}\\
&  \omega _{m+2s+l}(s,m,l,k+1) - \omega _{m+2s+l}(s,m,l,k) = 2^{k+s}\cdot (-5\cdot2^{k+m+s+2l-2})- 2^{k+s-1}\cdot (-5\cdot2^{k+m+s+2l-3})= -15\cdot2^{2k+2s+m+2l-4} \\
& \underline{ i= m+2s+l+1}\\
&  \omega _{m+2s+l+1}(s,m,l,k+1) - \omega _{m+2s+l+1}(s,m,l,k) = 2^{k+s}\cdot 2^{k+m+s+2l}- 2^{k+s-1}\cdot 2^{k+m+s+2l-1})=3\cdot2^{2k+2s+m+2l-2} \\
& \underline{ m+2s+l+2\leq i\leq 2m+3s+l-1}\\
&  \omega _{i}(s,m,l,k+1) - \omega _{i}(s,m,l,k) = 2^{k+s}\cdot 0 - 2^{k+s-1}\cdot 0 = 0 \\
& \underline{  i = 2m+3s+l}\\
&  \omega _{2m+3s+l}(s,m,l,k+1) - \omega _{2m+3s+l}(s,m,l,k) = 2^{k+s}\cdot(3\cdot2^{2k+2s+2m+l}) - 2^{k+s-1}\cdot (3\cdot2^{2k+2s+2m+l-2})\\
&  = 21\cdot2^{3k+3s+2m+l-3} 
\end{align*}
 \end{proof} 
   \section{\textbf{Computation of $ \Delta _{k}\Gamma_{s+j}^{\left[s\atop{ s\atop s} \right]} - 8\cdot \Delta _{k}\Gamma_{s+j-1}^{\left[s-1\atop{ s\atop s}\right]}\quad \text{for\quad $ 0\leq j\leq 2s $}  $}}
  \label{sec 8}
  \subsection{Notation}
  \label{subsec 1}
  \begin{defn}
  \label{defn 8.1}
   We define \\
   \begin{itemize}
\item  
   $$ \Delta _{k}\Gamma_{i}^{\left[s\atop{ s +m\atop s +m+l} \right]} =  
     \Gamma_{i}^{\left[s\atop{ s +m\atop s+m+l} \right]\times (k+1)} - \Gamma_{i}^{\left[s\atop{ s +m\atop s +m+l} \right]\times k} $$\\
\item     
  $$\Delta _{k}\omega_{i} (s,m,l) = \omega _{i}(s,m,l,k+1) - \omega_{i} (s,m,l,k) $$
  \end{itemize}
\end{defn}

   \subsection{Introduction}
  \label{subsec 2}

   We adapt the method used in Section 12 of [2] to compute explicitly \\[0.02 cm]
 $   \Delta_{k} \Gamma_{s+j}^{\left[s\atop{ s\atop s} \right]} -8\cdot \Delta_{k} \Gamma_{s+j-1}^{\left[s-1\atop{ s\atop s} \right]} \quad
 \text{for}\quad 0\leqslant j\leqslant 2s$

   \subsection{\textbf{A reduction formula for  $ \Delta _{k}\Gamma_{s+j}^{\left[s\atop{ s\atop s} \right]} - 8\cdot \Delta _{k}\Gamma_{s+j-1}^{\left[s-1\atop{ s\atop s}\right]}\quad \text{for\quad $ 0\leq j\leq 2s $}$} }
   \label{subsec 3}
   \begin{lem}
   \label{lem 8.2} We have the following reduction formulas : \\
   \begin{align}
      &  \Delta _{k}\Gamma_{s+j}^{\left[s\atop{ s \atop s } \right]}
 -  8\cdot \Delta _{k}\Gamma_{s+j -1}^{\left[s -1\atop{ s \atop s } \right]}
      =  4^{2s-3}\cdot\left[ \Delta _{k}\Gamma_{j- s  + 3}^{\left[1 \atop{ 2 \atop s }\right]} 
 - 8\cdot \Delta _{k}\Gamma_{j -s +2}^{\left[1 \atop{ 2  \atop s -1} \right]}\right] \label{eq 8.1}\\
    &    + \Delta _{k}\omega_{s+j} (s,0,0) + \sum_{q =1}^{s-2}  4^{2q-1}\cdot\Delta _{k}\omega_{s - q + (j-q+1)} (s-q,1, q-1) 
    +  \sum_{q =1}^{s-2}  4^{2q}\cdot\Delta _{k}\omega_{s - q +(j-q)} (s-q,0, q) \quad \text{for $0\leq j\leq s-1$} \nonumber  \\
     &  \Delta _{k}\Gamma_{2s+j}^{\left[s\atop{ s \atop s } \right]}
 -  8\cdot \Delta _{k}\Gamma_{2s+j -1}^{\left[s -1\atop{ s \atop s } \right]}
       =  4^{2s-3}\cdot\left[ \Delta _{k}\Gamma_{j + 3}^{\left[1 \atop{ 2 \atop s }\right]} 
 - 8\cdot \Delta _{k}\Gamma_{j+2}^{\left[1 \atop{ 2  \atop s -1} \right]}\right] \label{eq 8.2}\\
    &    + \Delta _{k}\omega_{2s+j} (s,0,0) + \sum_{q =1}^{s-2}  4^{2q-1}\cdot\Delta _{k}\omega_{2(s-q) + (j+1)} (s-q,1, q-1) 
    +  \sum_{q =1}^{s-2}  4^{2q}\cdot\Delta _{k}\omega_{2(s-q) + j} (s-q,0, q)   \quad \text{for $0\leq j\leq s $}  \nonumber
   \end{align}
 \end{lem}
 \begin{proof}
 From Definition \ref{defn 5.1} and Definition \ref{defn 8.1} we obtain the following formula \\
 \begin{align}
  \Delta_{k}\omega _{i}(s,m,l)  = \Delta _{k}\Gamma_{i}^{\left[s\atop{ s+m\atop s+m+l} \right]}
- 4\cdot\Delta_{k} \Gamma_{i-1}^{\left[s\atop{ s+m-1\atop s+m+l} \right]} -
8\cdot\Delta _{k} \Gamma_{i-1}^{\left[s\atop{ s+m\atop s+m+l-1} \right]} +
32\cdot\Delta _{k}\Gamma_{i-2}^{\left[s\atop{ s+m-1\atop s+m+l-1} \right]} \label{eq 8.3}
\end{align}
\newpage
  \underline{The case $ i = s+j,\quad 0\leq j\leq s-1 $} \vspace{0.1 cm}\\
   From \eqref{eq 8.3} we obtain successively the following equations \\
   
   \small   
  \begin{align*}
  \Delta _{k}\Gamma_{s+j}^{\left[s\atop{ s \atop s } \right]}
   - 4\cdot \Delta _{k}\Gamma_{s+j -1}^{\left[s-1\atop{ s \atop s } \right]} & = 
   8\cdot\left[  \Delta _{k}\Gamma_{s+j -1}^{\left[s-1\atop{ s \atop s } \right]}
   - 4\cdot \Delta _{k}\Gamma_{s+j -2}^{\left[s-1\atop{ s -1\atop s } \right]}\right] + \Delta _{k}\omega_{s+j} (s,0,0) \\
  4\cdot\left[  \Delta _{k}\Gamma_{s+j-1}^{\left[s -1\atop{ s \atop s } \right]}
   - 4\cdot \Delta _{k}\Gamma_{s+j -2}^{\left[s-1\atop{ s-1 \atop s } \right]}\right] & = 
4\cdot\left[   8\cdot\left[  \Delta _{k}\Gamma_{s+j -2}^{\left[s -1\atop{ s-1 \atop s } \right]}
   - 4\cdot \Delta _{k}\Gamma_{s+j -3}^{\left[s-1\atop{ s -1\atop s-1 } \right]}\right] + \Delta _{k}\omega_{s+j-1} (s-1,1,0) \right]\\
    4^{2}\cdot\left[   \Delta _{k}\Gamma_{s+j-2}^{\left[s -1\atop{ s -1\atop s } \right]}
   - 4\cdot \Delta _{k}\Gamma_{s+j -3}^{\left[s- 2\atop{ s-1 \atop s } \right]}\right] & = 
 4^{2}\cdot\left[  8\cdot\left[  \Delta _{k}\Gamma_{s+j -3}^{\left[s -1\atop{ s-1 \atop s -1} \right]}
   - 4\cdot \Delta _{k}\Gamma_{s+j - 4}^{\left[s-2 \atop{ s -1\atop s-1 } \right]}\right] + \Delta _{k}\omega_{s+j-2} (s-1,0,1)\right] \\
    4^{3}\cdot\left[   \Delta _{k}\Gamma_{s+j-3}^{\left[s -2\atop{ s -1\atop s } \right]}
   - 4\cdot \Delta _{k}\Gamma_{s+j -4}^{\left[s- 2\atop{ s-2 \atop s } \right]}\right] & = 
 4^{3}\cdot\left[  8\cdot\left[  \Delta _{k}\Gamma_{s+j -4}^{\left[s -2\atop{ s-1 \atop s -1} \right]}
   - 4\cdot \Delta _{k}\Gamma_{s+j - 5}^{\left[s-2 \atop{ s -2\atop s- 1 } \right]}\right] + \Delta _{k}\omega_{s+j-3} (s-2,1,1)\right] \\
     4^{4}\cdot\left[   \Delta _{k}\Gamma_{s+j-4}^{\left[s -2\atop{ s -2\atop s } \right]}
   - 4\cdot \Delta _{k}\Gamma_{s+j -5}^{\left[s- 3\atop{ s-2 \atop s } \right]}\right] & = 
 4^{4}\cdot\left[  8\cdot\left[  \Delta _{k}\Gamma_{s+j - 5}^{\left[s -2\atop{ s- 2 \atop s -1} \right]}
   - 4\cdot \Delta _{k}\Gamma_{s+j - 6}^{\left[s-3 \atop{ s -2\atop s-1 } \right]}\right] + \Delta _{k}\omega_{s+j-4} (s-2,0,2)\right] \\
     &  \vdots       \\
      4^{2p-1}\cdot\left[   \Delta _{k}\Gamma_{s+j-(2p-1)}^{\left[s -p\atop{ s - (p-1)\atop s } \right]}
   - 4\cdot \Delta _{k}\Gamma_{s+j - 2p}^{\left[s- p\atop{ s- p \atop s } \right]}\right] & = 
 4^{2p-1}\cdot\left[  8\cdot\left[  \Delta _{k}\Gamma_{s+j -2p}^{\left[s -p\atop{ s-(p-1) \atop s -1} \right]}
   - 4\cdot \Delta _{k}\Gamma_{s+j - 2p-1}^{\left[s- p \atop{ s -p\atop s- 1 } \right]}\right] + \Delta _{k}\omega_{s+j- (2p-1)} (s-p,1, p-1)\right] \\ 
     4^{2p}\cdot\left[   \Delta _{k}\Gamma_{s+j- 2p}^{\left[s -p\atop{ s - p \atop s } \right]}
   - 4\cdot \Delta _{k}\Gamma_{s+j - (2p+1)}^{\left[s- p -1 \atop{ s- p \atop s } \right]}\right] & = 
 4^{2p}\cdot\left[  8\cdot\left[  \Delta _{k}\Gamma_{s+j - (2p+1)}^{\left[s -p\atop{ s- p \atop s -1} \right]}
   - 4\cdot \Delta _{k}\Gamma_{s+j - (2p+2)}^{\left[s- p-1 \atop{ s -p\atop s- 1 } \right]}\right] + \Delta _{k}\omega_{s+j- 2p} (s-p,0, p)\right] 
    \end{align*}
      By summing the left- hand side of  the above equations we obtain: \\ 
     \begin{align*}
   &   \Delta _{k}\Gamma_{s+j}^{\left[s\atop{ s \atop s } \right]}
   - 4\cdot \Delta _{k}\Gamma_{s+j -1}^{\left[s-1\atop{ s \atop s } \right]}
    + \sum_{q=1}^{p}   4^{2q-1}\cdot\left[   \Delta _{k}\Gamma_{s+j-(2q-1)}^{\left[s -q\atop{ s - (q-1)\atop s } \right]}
   - 4\cdot \Delta _{k}\Gamma_{s+j - 2q}^{\left[s- q\atop{ s- q \atop s } \right]}\right] 
       +   \sum_{q=1}^{p}   4^{2q}\cdot\left[   \Delta _{k}\Gamma_{s+j- 2q}^{\left[s -q \atop{ s -  q \atop s } \right]}
   - 4\cdot \Delta _{k}\Gamma_{s+j - (2q+1)}^{\left[s- q -1\atop{ s- q \atop s } \right]}\right] \\
   & =  \Delta _{k}\Gamma_{s+j}^{\left[s\atop{ s \atop s } \right]} +
     \sum_{q=1}^{p}   4^{2q-1}\cdot   \Delta _{k}\Gamma_{s+j- (2q-1)}^{\left[s -q \atop{ s -  (q-1) \atop s }\right]}+ 
     \sum_{q=1}^{p}   4^{2q}\cdot   \Delta _{k}\Gamma_{s+j- 2q}^{\left[s -q \atop{ s -  q \atop s }\right]} 
      -   \sum_{q=1}^{p}   4^{2q}\cdot   \Delta _{k}\Gamma_{s+j- 2q}^{\left[s -q \atop{ s -  q \atop s }\right]} - 
     \sum_{q=1}^{p}   4^{2q+1}\cdot \Delta _{k}\Gamma_{s+j- (2q+1)}^{\left[s -q-1 \atop{ s - q \atop s }\right]} 
       - 4\cdot \Delta _{k}\Gamma_{s+j -1}^{\left[s-1\atop{ s \atop s } \right]}\\
      & =   \Delta _{k}\Gamma_{s+j}^{\left[s\atop{ s \atop s } \right]} +
     \sum_{q=1}^{p}   4^{2q-1}\cdot   \Delta _{k}\Gamma_{s+j- (2q-1)}^{\left[s -q \atop{ s -  (q-1) \atop s }\right]}+ 
     \sum_{q=1}^{p}   4^{2q}\cdot   \Delta _{k}\Gamma_{s+j- 2q}^{\left[s -q \atop{ s -  q \atop s }\right]} 
      -   \sum_{q=1}^{p}   4^{2q}\cdot   \Delta _{k}\Gamma_{s+j- 2q}^{\left[s -q \atop{ s -  q \atop s }\right]} - 
     \sum_{q=0}^{p}   4^{2q+1}\cdot \Delta _{k}\Gamma_{s+j- (2q+1)}^{\left[s -q-1 \atop{ s - q \atop s }\right]}\\
         & =   \Delta _{k}\Gamma_{s+j}^{\left[s\atop{ s \atop s } \right]} +
     \sum_{q=1}^{p}   4^{2q-1}\cdot   \Delta _{k}\Gamma_{s+j- (2q-1)}^{\left[s -q \atop{ s -  (q-1) \atop s }\right]}+ 
     \sum_{q=1}^{p}   4^{2q}\cdot   \Delta _{k}\Gamma_{s+j- 2q}^{\left[s -q \atop{ s -  q \atop s }\right]} 
      -   \sum_{q=1}^{p}   4^{2q}\cdot   \Delta _{k}\Gamma_{s+j- 2q}^{\left[s -q \atop{ s -  q \atop s }\right]} - 
     \sum_{q= 1}^{p+1}   4^{2q-1}\cdot \Delta _{k}\Gamma_{s+j- (2q-1)}^{\left[s -q \atop{ s - (q-1) \atop s }\right]}\\
     & =  \Delta _{k}\Gamma_{s+j}^{\left[s\atop{ s \atop s } \right]} - 4^{2p +1}\cdot \Delta _{k}\Gamma_{s+j- 2p-1}^{\left[s - p-1 \atop{ s - p \atop s }\right]}
       \end{align*}
       
     By summing the right-hand side of  the above equations we obtain: \\ 
    \begin{align*}
     &  8\cdot\left[  \Delta _{k}\Gamma_{s+j -1}^{\left[s-1\atop{ s \atop s } \right]}
   - 4\cdot \Delta _{k}\Gamma_{s+j -2}^{\left[s-1\atop{ s -1\atop s } \right]}\right] + \Delta _{k}\omega_{s+j} (s,0,0) \\
  &+ \sum_{q =1}^{p}  4^{2q-1}\cdot\left[  8\cdot\left[  \Delta _{k}\Gamma_{s+j -2q}^{\left[s -q\atop{ s-(q-1) \atop s -1} \right]}
   - 4\cdot \Delta _{k}\Gamma_{s+j - 2q-1}^{\left[s- q \atop{ s -q\atop s- 1 } \right]}\right] + \Delta _{k}\omega_{s+j- (2q-1)} (s-q,1, q-1)\right] \\ 
   & +\sum_{q =1}^{p}  4^{2q}\cdot\left[  8\cdot\left[  \Delta _{k}\Gamma_{s+j -(2q+1)}^{\left[s -q\atop{ s- q \atop s -1} \right]}
   - 4\cdot \Delta _{k}\Gamma_{s+j - (2q+2)}^{\left[s- q-1 \atop{ s -q\atop s- 1 } \right]}\right] + \Delta _{k}\omega_{s+j- 2q} (s-q,0, q)\right]  \\
  & =   8\cdot\left[  \Delta _{k}\Gamma_{s+j -1}^{\left[s -1\atop{ s \atop s } \right]}
   - 4\cdot \Delta _{k}\Gamma_{s+j -2}^{\left[s-1\atop{ s -1\atop s } \right]}\right] + \Delta _{k}\omega_{s+j} (s,0,0) 
      + \sum_{q =1}^{p} 2\cdot 4^{2q}\cdot \Delta _{k}\Gamma_{s+j - 2q}^{\left[s -q\atop{ s- (q-1) \atop s -1} \right]} +
    \sum_{q =1}^{p} 2\cdot 4^{2q+1}\cdot \Delta _{k}\Gamma_{s+j - (2q+1)}^{\left[s -q\atop{ s- q \atop s -1} \right]}\\
    & -  \sum_{q =1}^{p} 2\cdot 4^{2q+1}\cdot \Delta _{k}\Gamma_{s+j - (2q+1)}^{\left[s -q\atop{ s- q \atop s -1} \right]}
    -  \sum_{q =1}^{p} 2\cdot 4^{2q+2}\cdot \Delta _{k}\Gamma_{s+j - (2q+2)}^{\left[s -q-1\atop{ s- q \atop s -1} \right]}\\
    & +  \sum_{q =1}^{p}  4^{2q-1}\cdot\Delta _{k}\omega_{s+j- (2q-1)} (s-q,1, q-1) 
    +  \sum_{q =1}^{p}  4^{2q}\cdot\Delta _{k}\omega_{s+j- 2q} (s-q,0, q) \\
     & =   8\cdot \Delta _{k}\Gamma_{s+j -1}^{\left[s -1\atop{ s \atop s } \right]} + \Delta _{k}\omega_{s+j} (s,0,0)
      + \sum_{q =1}^{p} 2\cdot 4^{2q}\cdot \Delta _{k}\Gamma_{s+j - 2q}^{\left[s -q\atop{ s- (q-1) \atop s -1} \right]} 
        -  \sum_{q =0}^{p} 2\cdot 4^{2q+2}\cdot \Delta _{k}\Gamma_{s+j - (2q+2)}^{\left[s -q-1\atop{ s- q \atop s -1} \right]}\\
    & +  \sum_{q =1}^{p}  4^{2q-1}\cdot\Delta _{k}\omega_{s+j- (2q-1)} (s-q,1, q-1) 
    +  \sum_{q =1}^{p}  4^{2q}\cdot\Delta _{k}\omega_{s+j- 2q} (s-q,0, q) \\
     & =   8\cdot \Delta _{k}\Gamma_{s+j -1}^{\left[s -1\atop{ s \atop s } \right]} + \Delta _{k}\omega_{s+j} (s,0,0)
      + \sum_{q =1}^{p} 2\cdot 4^{2q}\cdot \Delta _{k}\Gamma_{s+j - 2q}^{\left[s -q\atop{ s- (q-1) \atop s -1} \right]} 
        -  \sum_{q = 1}^{p+1} 2\cdot 4^{2q}\cdot \Delta _{k}\Gamma_{s+j - 2q}^{\left[s -q\atop{ s- (q-1) \atop s -1} \right]}\\
    & +  \sum_{q =1}^{p}  4^{2q-1}\cdot\Delta _{k}\omega_{s+j- (2q-1)} (s-q,1, q-1) 
    +  \sum_{q =1}^{p}  4^{2q}\cdot\Delta _{k}\omega_{s+j- 2q} (s-q,0, q) \\
    & =  8\cdot \Delta _{k}\Gamma_{s+j -1}^{\left[s -1\atop{ s \atop s } \right]} + \Delta _{k}\omega_{s+j} (s,0,0)
     - 2\cdot 4^{2(p+1)}\cdot \Delta _{k}\Gamma_{s+j - 2(p+1)}^{\left[s - (p+1)\atop{ s- p  \atop s -1} \right]}\\
    & +   \sum_{q =1}^{p}  4^{2q-1}\cdot\Delta _{k}\omega_{s+j- (2q-1)} (s-q,1, q-1) 
    +  \sum_{q =1}^{p}  4^{2q}\cdot\Delta _{k}\omega_{s+j- 2q} (s-q,0, q) 
     \end{align*}
  Therefore \\
   \begin{align}
&  \Delta _{k}\Gamma_{s+j}^{\left[s\atop{ s \atop s } \right]} - 4^{2p +1}\cdot \Delta _{k}\Gamma_{s+j- 2p-1}^{\left[s - p-1 \atop{ s - p \atop s }\right]} \\
& =  8\cdot \Delta _{k}\Gamma_{s+j -1}^{\left[s -1\atop{ s \atop s } \right]} + \Delta _{k}\omega_{s+j} (s,0,0)
     - 2\cdot 4^{2(p+1)}\cdot \Delta _{k}\Gamma_{s+j - 2(p+1)}^{\left[s - (p+1)\atop{ s- p  \atop s -1} \right]} \nonumber \\
    & +   \sum_{q =1}^{p}  4^{2q-1}\cdot\Delta _{k}\omega_{s+j- (2q-1)} (s-q,1, q-1) 
    +  \sum_{q =1}^{p}  4^{2q}\cdot\Delta _{k}\omega_{s+j- 2q} (s-q,0, q) \nonumber  \\
& \Longleftrightarrow  \nonumber \\
&  \Delta _{k}\Gamma_{s+j}^{\left[s\atop{ s \atop s } \right]}
 -  8\cdot \Delta _{k}\Gamma_{s+j -1}^{\left[s -1\atop{ s \atop s } \right]}\label{eq 8.5}\\
 & =  4^{2p +1}\cdot \Delta _{k}\Gamma_{s+j- 2p-1}^{\left[s - p-1 \atop{ s - p \atop s }\right]} 
 - 2\cdot 4^{2(p+1)}\cdot \Delta _{k}\Gamma_{s+j - 2(p+1)}^{\left[s - (p+1)\atop{ s- p  \atop s -1} \right]} \nonumber \\
    &    + \Delta _{k}\omega_{s+j} (s,0,0) + \sum_{q =1}^{p}  4^{2q-1}\cdot\Delta _{k}\omega_{s - q + (j-q+1)} (s-q,1, q-1) 
    +  \sum_{q =1}^{p}  4^{2q}\cdot\Delta _{k}\omega_{s - q +(j-q)} (s-q,0, q) \nonumber
    \end{align}
   Setting  p = s -2  in  \eqref{eq 8.5} we obtain the formula \eqref{eq 8.1}, that is \\
     \begin{align*}
    &  \Delta _{k}\Gamma_{s+j}^{\left[s\atop{ s \atop s } \right]}
 -  8\cdot \Delta _{k}\Gamma_{s+j -1}^{\left[s -1\atop{ s \atop s } \right]} \\
     & =  4^{2s-3}\cdot\left[ \Delta _{k}\Gamma_{j- s  + 3}^{\left[1 \atop{ 2 \atop s }\right]} 
 - 8\cdot \Delta _{k}\Gamma_{j -s +2}^{\left[1 \atop{ 2  \atop s -1} \right]}\right]  \\
    &    + \Delta _{k}\omega_{s+j} (s,0,0) + \sum_{q =1}^{s-2}  4^{2q-1}\cdot\Delta _{k}\omega_{s - q + (j-q+1)} (s-q,1, q-1) 
    +  \sum_{q =1}^{s-2}  4^{2q}\cdot\Delta _{k}\omega_{s - q +(j-q)} (s-q,0, q) 
  \end{align*}
  \underline{The case $ i = 2s+j,\quad 0\leq j\leq s $} \vspace{0.1 cm}\\
  Proceeding as in the first case we obtain \\
   \begin{align}
&  \Delta _{k}\Gamma_{2s+j}^{\left[s\atop{ s \atop s } \right]}
 -  8\cdot \Delta _{k}\Gamma_{2s+j -1}^{\left[s -1\atop{ s \atop s } \right]} \label{eq 8.6}\\
 & =  4^{2p +1}\cdot \Delta _{k}\Gamma_{2s+j- 2p-1}^{\left[s - p-1 \atop{ s - p \atop s }\right]} 
 - 2\cdot 4^{2(p+1)}\cdot \Delta _{k}\Gamma_{2s+j - 2(p+1)}^{\left[s - (p+1)\atop{ s- p  \atop s -1} \right]} \nonumber \\
    &    + \Delta _{k}\omega_{2s+j} (s,0,0) + \sum_{q =1}^{p}  4^{2q-1}\cdot\Delta _{k}\omega_{2(s-q) + (j+1)} (s-q,1, q-1) 
    +  \sum_{q =1}^{p}  4^{2q}\cdot\Delta _{k}\omega_{2(s-q) + j} (s-q,0, q) \nonumber 
  \end{align} 
   Setting p = s -2 in \eqref{eq 8.6} we get the formula \eqref{eq 8.2}, that is :\\
   \begin{align*}
   &  \Delta _{k}\Gamma_{2s+j}^{\left[s\atop{ s \atop s } \right]}
 -  8\cdot \Delta _{k}\Gamma_{2s+j -1}^{\left[s -1\atop{ s \atop s } \right]} \\
     & =  4^{2s-3}\cdot\left[ \Delta _{k}\Gamma_{j + 3}^{\left[1 \atop{ 2 \atop s }\right]} 
 - 8\cdot \Delta _{k}\Gamma_{j+2}^{\left[1 \atop{ 2  \atop s -1} \right]}\right]   \\
    &    + \Delta _{k}\omega_{2s+j} (s,0,0) + \sum_{q =1}^{s-2}  4^{2q-1}\cdot\Delta _{k}\omega_{2(s-q) + (j+1)} (s-q,1, q-1) 
    +  \sum_{q =1}^{s-2}  4^{2q}\cdot\Delta _{k}\omega_{2(s-q) + j} (s-q,0, q)  
\end{align*}

 \end{proof}
 
   \subsection{\textbf{An explicit formula for the difference \\
    $ \Delta _{k}\Gamma_{s+j}^{\left[s\atop{ s\atop s} \right]} - 8\cdot \Delta _{k}\Gamma_{s+j-1}^{\left[s-1\atop{ s\atop s}\right]}\quad \text{for\quad $ 0\leq j\leq 2s $}$} }
   \label{subsec 4}
       \begin{lem}
   \label{lem  8.3}
   We have for $s\geq 5 $ \\
   
   \begin{equation}
   \label{eq 8.7}
    \Delta _{k}\Gamma_{j-s+3}^{\left[1\atop{ 2\atop s} \right]} - 8\cdot \Delta _{k}\Gamma_{j-s+2}^{\left[1 \atop{ 2 \atop s-1} \right]}=
 \begin{cases}
0  & \text{if   } 0 \leq j \leq s-3\\
2^{k}  & \text{if   } j = s -2\\
9\cdot2^{k} & \text{if   } j = s -1 
\end{cases}
\end{equation}
\begin{equation}
\label{eq 8.8}
 \Delta _{k}\Gamma_{j+3}^{\left[1\atop{ 2 \atop s} \right]} - 8\cdot \Delta _{k}\Gamma_{j+2}^{\left[1 \atop{ 2 \atop s-1} \right]}=
 \begin{cases}
3\cdot2^{2k+1} - 17\cdot2^{k+1}  & \text{if   } j = 0 \\
 -15\cdot2^{2k+1} +3 \cdot2^{k+3}  & \text{if   } j = 1 \\
 -9\cdot2^{2k+2j-1} & \text{if   } 2\leq j\leq s-3 \\
  -33\cdot2^{2k+2s-5} & \text{if   } j = s-2 \\
-159\cdot2^{2k+2s-3} & \text{if   } j = s-1\\
-21\cdot2^{3k+s} + 21\cdot2^{2k+2s}  & \text{if   } j = s \\
\end{cases}
 \end{equation}
 \end{lem}
\begin{proof} 
 From Lemma \ref{lem 1.22} we deduce easily Lemma \ref{lem 8.3}. 
 \end{proof}

   \begin{lem}
   \label{lem  8.4}
   We have :\\
   
   \begin{equation}
   \label{eq 8.9}
    \Delta _{k}\Gamma_{s+j}^{\left[s\atop{ s\atop s} \right]} - 8\cdot \Delta _{k}\Gamma_{s+j-1}^{\left[s-1\atop{ s\atop s} \right]}=
 \begin{cases}
3\cdot2^{k+s-1} &\text{if  }  j = 0 \\
15\cdot2^{k+s-1}& \text{if   } j = 1 \\
63\cdot2^{k+s+3j-6} & \text{if   } 2\leq j\leq s-1\\
\end{cases}
\end{equation}
\begin{equation}
\label{eq 8.10}
 \Delta _{k}\Gamma_{2s+j}^{\left[s\atop{ s\atop s} \right]} - 8\cdot \Delta _{k}\Gamma_{2s+j-1}^{\left[s-1\atop{ s\atop s} \right]}=
 \begin{cases}
 -15\cdot2^{2k+2s-2} - 33\cdot2^{k+4s-6}  & \text{if   } j = 0 \\
-321\cdot2^{2k+2s-2}  + 3\cdot2^{k+4s-3}  & \text{if   } j = 1 \\
 -315\cdot2^{2k+2s+4j-6} & \text{if   } 2\leq j\leq s-1\\
-21\cdot2^{3k+3s-3} + 21\cdot2^{2k+6s-6}  & \text{if   } j = s \\
\end{cases}
 \end{equation}
 \end{lem}
 \begin{proof}
 \underline{proof of \eqref{eq 8.9}}\\
 From the reduction formula \eqref{eq 8.1} using  Lemma \ref{lem 7.4} and the formula \eqref{eq 8.7} we deduce \\
 
  \underline{The case  j = 0 }\\
  \begin{align*}
 &  \Delta _{k}\Gamma_{s}^{\left[s\atop{ s \atop s } \right]}
 -  8\cdot \Delta _{k}\Gamma_{s -1}^{\left[s -1\atop{ s \atop s } \right]} \\
     & =  2^{4s-6}\cdot\left[ \Delta _{k}\Gamma_{- s  + 3}^{\left[1 \atop{ 2 \atop s }\right]} 
 - 8\cdot \Delta _{k}\Gamma_{ -s +2}^{\left[1 \atop{ 2  \atop s -1} \right]}\right] \\
    &    + \Delta _{k}\omega_{s} (s,0,0) + \sum_{q =1}^{s-2}  2^{4q-2}\cdot\Delta _{k}\omega_{s - q + (-q+1)} (s-q,1, q-1) 
    +  \sum_{q =1}^{s-2}  2^{4q}\cdot\Delta _{k}\omega_{s - q + (-q)} (s-q,0, q) \\
       &    =  \Delta _{k}\omega_{s} (s,0,0) + 2^{2}\cdot\Delta _{k}\omega_{s - 1} (s-1,1,0) 
        = 2^{k+s-1} +2^{2}\cdot2^{k+s-2} = 3\cdot 2^{k+s-1}
    \end{align*}
  \underline{The case  j = 1 }\\
  \begin{align*}
 &  \Delta _{k}\Gamma_{s+1}^{\left[s\atop{ s \atop s } \right]}
 -  8\cdot \Delta _{k}\Gamma_{s }^{\left[s -1\atop{ s \atop s } \right]} \\
     & =  2^{4s-6}\cdot\left[ \Delta _{k}\Gamma_{- s  + 4}^{\left[1 \atop{ 2 \atop s }\right]} 
 - 8\cdot \Delta _{k}\Gamma_{ -s + 3}^{\left[1 \atop{ 2  \atop s -1} \right]}\right] \\
    &    + \Delta _{k}\omega_{s +1} (s,0,0) + \sum_{q =1}^{s-2}  2^{4q-2}\cdot\Delta _{k}\omega_{s - q + (-q+ 2)} (s-q,1, q-1) 
    +  \sum_{q =1}^{s-2}  2^{4q}\cdot\Delta _{k}\omega_{s - q + (-q +1)} (s-q,0, q) \\
     &    =  \Delta _{k}\omega_{s +1} (s,0,0) + \sum_{q =1}^{2}  2^{4q-2}\cdot\Delta _{k}\omega_{s - q + (-q+ 2)} (s-q,1, q-1) 
    +  \sum_{q =1}^{1}  2^{4q}\cdot\Delta _{k}\omega_{s - q + (-q +1)} (s-q,0, q) \\
    & =  \Delta _{k}\omega_{s +1} (s,0,0) + 2^{2}\cdot\Delta _{k}\omega_{s - 1 + 1} (s-1,1, 0) +  2^{6}\cdot\Delta _{k}\omega_{s - 2 } (s-2,1, 1) + 2^{4}\cdot\Delta _{k}\omega_{s - 1 } (s-1,0, 1) \\
     & = - 3\cdot 2^{k+s-1}+ 2^{2}\cdot(- 3\cdot 2^{k+s-2})+ 2^{6}\cdot 2^{k+s-3} +2^{4}\cdot 2^{k+s-2}=15\cdot2^{k+s-1}
    \end{align*}
     \underline{The case  j = 2}\\
     \begin{align*}
     &  \Delta _{k}\Gamma_{s+2}^{\left[s\atop{ s \atop s } \right]}
 -  8\cdot \Delta _{k}\Gamma_{s+ 1}^{\left[s -1\atop{ s \atop s } \right]}\\
     & =  2^{4s-6}\cdot\left[ \Delta _{k}\Gamma_{- s  + 5}^{\left[1 \atop{ 2 \atop s }\right]} 
 - 8\cdot \Delta _{k}\Gamma_{ -s + 4}^{\left[1 \atop{ 2  \atop s -1} \right]}\right]\\
    &    + \Delta _{k}\omega_{s+2} (s,0,0) + \sum_{q =1}^{s-2}  2^{4q-2}\cdot\Delta _{k}\omega_{s - q + (-q+3)} (s-q,1, q-1) 
    +  \sum_{q =1}^{s-2}  2^{4q}\cdot\Delta _{k}\omega_{s - q +( -q+2)} (s-q,0, q) \\
    & =   \Delta _{k}\omega_{s+2} (s,0,0) + \sum_{q =1}^{3}  2^{4q-2}\cdot\Delta _{k}\omega_{s - q + (-q+3)} (s-q,1, q-1) 
    +  \sum_{q =1}^{2}  2^{4q}\cdot\Delta _{k}\omega_{s - q +( -q+2)} (s-q,0, q) \\
    & =  \Delta _{k}\omega_{s + 2} (s,0, 0) +  2^{2}\cdot\Delta _{k}\omega_{s - 1 + 2} (s-1,1, 0) 
   +  2^{6}\cdot\Delta _{k}\omega_{s - 2 + 1} (s-2,1, 1) +  2^{10}\cdot\Delta _{k}\omega_{s - 3} (s-3,1, 2) \\
    &  +  2^{4}\cdot\Delta _{k}\omega_{s - 1 + 1} (s-1,0, 1) +  2^{8}\cdot\Delta _{k}\omega_{s - 2} (s-2,0, 2) \\
    & = 2^{k+s}+  2^{2}\cdot2^{k+s-1} + 2^{6}\cdot(- 3\cdot2^{k+s-3})   + 2^{10}\cdot2^{k+s-4} 
  + 2^{4}\cdot(- 3\cdot2^{k+s-2} ) + 2^{8} \cdot2^{k+s -3}= 63\cdot2^{k+s}
  \end{align*} 
    \underline{The case $ 3\leq j\leq s-3  $}\vspace{0.1 cm}\\
  \begin{align*}
 & \Delta _{k}\omega_{s - q + (j-q+1)} (s-q,1, q-1) = 0 \\
 & \Longleftrightarrow  s-q+3 \leq s - q + (j-q+1)\leq 1+2\cdot(s-q) -1 \Longleftrightarrow q\leq j-2 \\
&  \Delta _{k}\omega_{s - q +(j-q)} (s-q,0, q) = 0 \\
 & \Longleftrightarrow  s-q+3 \leq s - q + ( j-q ) \leq 2\cdot(s-q) -1 \Longleftrightarrow q\leq j-3 
\end{align*}
 
     \begin{align*}
      &  \Delta _{k}\Gamma_{s+j}^{\left[s\atop{ s \atop s } \right]}
 -  8\cdot \Delta _{k}\Gamma_{s+j -1}^{\left[s -1\atop{ s \atop s } \right]}\\
     & =  4^{2s-3}\cdot\left[ \Delta _{k}\Gamma_{j- s  + 3}^{\left[1 \atop{ 2 \atop s }\right]} 
 - 8\cdot \Delta _{k}\Gamma_{j -s +2}^{\left[1 \atop{ 2  \atop s -1} \right]}\right]\\
    &    + \Delta _{k}\omega_{s+j} (s,0,0) + \sum_{q =1}^{s-2}  4^{2q-1}\cdot\Delta _{k}\omega_{s - q + (j-q+1)} (s-q,1, q-1) 
    +  \sum_{q =1}^{s-2}  4^{2q}\cdot\Delta _{k}\omega_{s - q +(j-q)} (s-q,0, q) \\
   &  =    \Delta _{k}\omega_{s+j} (s,0,0) + \sum_{q =1}^{j+1}  4^{2q-1}\cdot\Delta _{k}\omega_{s - q + (j-q+1)} (s-q,1, q-1) 
    +  \sum_{q =1}^{j}  4^{2q}\cdot\Delta _{k}\omega_{s - q +(j-q)} (s-q,0, q) \\  
      &  =    \Delta _{k}\omega_{s+j} (s,0,0) + \sum_{q = j-1}^{j+1}  4^{2q-1}\cdot\Delta _{k}\omega_{s - q + (j-q+1)} (s-q,1, q-1) 
    +  \sum_{q = j-2}^{j}  4^{2q}\cdot\Delta _{k}\omega_{s - q +(j-q)} (s-q,0, q) \\  
    & = [2^{4(j-1) -2} \cdot2^{k+s -(j-1)} + 2^{4j -2} \cdot(-3\cdot2^{k+ (s-j) -1})  + 2^{4j + 2} \cdot2^{k+s -(j+1)-1} ] \\
    & +  [2^{4(j-2)} \cdot2^{k+s -(j-2)} + 2^{4(j-1)} \cdot(-3\cdot2^{k+ s -(j-1)-1})  + 2^{4j } \cdot2^{k+ (s- j) -1} ] = 63\cdot2^{k+s+3j-6}
  \end{align*}
  
    \underline{The case  j = s -2  }\vspace{0.1 cm}\\
     \begin{align*}
 & \Delta _{k}\omega_{s - q + (s -q -1)} (s-q,1, q-1) = 0 \\
 & \Longleftrightarrow  s-q+3 \leq s - q + (s-q-1)\leq 1+2\cdot(s-q) -1 \Longleftrightarrow q\leq s-4 \\
&  \Delta _{k}\omega_{s - q +(s-q -2)} (s-q,0, q) = 0 \\
 & \Longleftrightarrow  s-q+3 \leq s - q + ( s-q -2) \leq 2\cdot(s-q) -1 \Longleftrightarrow q\leq s-5 
\end{align*}
    
      \begin{align*}
      &  \Delta _{k}\Gamma_{2s -2}^{\left[s\atop{ s \atop s } \right]}
 -  8\cdot \Delta _{k}\Gamma_{2s -3}^{\left[s -1\atop{ s \atop s } \right]}\\
     & =  2^{4s-6}\cdot\left[ \Delta _{k}\Gamma_{1}^{\left[1 \atop{ 2 \atop s }\right]} 
 - 8\cdot \Delta _{k}\Gamma_{0}^{\left[1 \atop{ 2  \atop s -1} \right]}\right]\\
    &    + \Delta _{k}\omega_{2s -2} (s,0,0) + \sum_{q =1}^{s-2}  2^{4q-2}\cdot\Delta _{k}\omega_{s - q + (s - q -1)} (s-q,1, q-1) 
    +  \sum_{q =1}^{s-2}  2^{4q}\cdot\Delta _{k}\omega_{s - q +(s - q - 2)} (s-q,0, q) \\
    & =  2^{4s-6}\cdot2^{k}     + \Delta _{k}\omega_{2s -2} (s,0,0) + 
    \sum_{q = s-3}^{s-2}  2^{4q-2}\cdot\Delta _{k}\omega_{s - q + (s - q -1)} (s-q,1, q-1) 
    +  \sum_{q = s-4}^{s-2}  2^{4q}\cdot\Delta _{k}\omega_{s - q +(s - q - 2)} (s-q,0, q) \\  
    & = 2^{k+4s -6} + [2^{4(s-3)-2}\cdot2^{k +s -(s-3)} +  2^{4(s-2)-2}\cdot(-3\cdot2^{k +s -(s-2)-1})]\\
    & +  [2^{4(s-4)}\cdot2^{k +s -(s-4)} +  2^{4(s-3)}\cdot(-3\cdot2^{k +s -(s-3)-1}) +2^{4(s-2)}\cdot2^{k +s -(s-2)-1} ]\\
    & = 63\cdot2^{k+4s-12}
      \end{align*}
      
    \underline{The case  j = s -1  }\\
     \begin{align*}
 & \Delta _{k}\omega_{s - q + (s -q )} (s-q,1, q-1) = 0 \\
 & \Longleftrightarrow  s-q+3 \leq s - q + (s-q)\leq 1+2\cdot(s-q) -1 \Longleftrightarrow q\leq s-3 \\
&  \Delta _{k}\omega_{s - q +(s-q -1)} (s-q,0, q) = 0 \\
 & \Longleftrightarrow  s-q+3 \leq s - q + ( s-q -1) \leq 2\cdot(s-q) -1 \Longleftrightarrow q\leq s-4 
\end{align*}
      \begin{align*}
      &  \Delta _{k}\Gamma_{2s -1}^{\left[s\atop{ s \atop s } \right]}
 -  8\cdot \Delta _{k}\Gamma_{2s -2}^{\left[s -1\atop{ s \atop s } \right]}\\
     & =  2^{4s-6}\cdot\left[ \Delta _{k}\Gamma_{2}^{\left[1 \atop{ 2 \atop s }\right]} 
 - 8\cdot \Delta _{k}\Gamma_{1}^{\left[1 \atop{ 2  \atop s -1} \right]}\right]\\
    &    + \Delta _{k}\omega_{2s -1} (s,0,0) + \sum_{q =1}^{s-2}  2^{4q-2}\cdot\Delta _{k}\omega_{s - q + (s - q )} (s-q,1, q-1) 
    +  \sum_{q =1}^{s-2}  2^{4q}\cdot\Delta _{k}\omega_{s - q +(s - q - 1)} (s-q,0, q) \\
       & =   2^{4s-6}\cdot\left[ \Delta _{k}\Gamma_{2}^{\left[1 \atop{ 2 \atop s }\right]} 
 - 8\cdot \Delta _{k}\Gamma_{1}^{\left[1 \atop{ 2  \atop s -1} \right]}\right]     + \Delta _{k}\omega_{2s -1} (s,0,0) + 
    \sum_{q = s-2}^{s-2}  2^{4q-2}\cdot\Delta _{k}\omega_{s - q + (s - q )} (s-q,1, q-1) \\
     &  +  \sum_{q = s-3}^{s-2}  2^{4q}\cdot\Delta _{k}\omega_{s - q +(s - q - 1)} (s-q,0, q) \\  
     & = 2^{k+4s -6}\cdot9\cdot2^{k}  + 2^{4(s-2)-2}\cdot2^{k +s -(s-2)} 
    +  2^{4(s-3)}\cdot2^{k + 3}  +  2^{4(s-2)}\cdot(-3\cdot2^{k + 1})  = 63\cdot2^{k+4s-9} 
   \end{align*}
    \underline{proof of \eqref{eq 8.10}}\\
   From the reduction formula \eqref{eq 8.2} using Lemma \ref{lem 7.4} and the formula \eqref{eq 8.8} we deduce \\
\underline{The case  j = 0}\\
\begin{align*}
&  \Delta _{k}\Gamma_{2s}^{\left[s\atop{ s \atop s } \right]}
 -  8\cdot \Delta _{k}\Gamma_{2s -1}^{\left[s -1\atop{ s \atop s } \right]}\\
  & =  4^{2s-3}\cdot\left[ \Delta _{k}\Gamma_{ 3}^{\left[1 \atop{ 2 \atop s }\right]} 
 - 8\cdot \Delta _{k}\Gamma_{2}^{\left[1 \atop{ 2  \atop s -1} \right]}\right]\\
    &    + \Delta _{k}\omega_{2s} (s,0,0) + \sum_{q =1}^{s-2}  4^{2q-1}\cdot\Delta _{k}\omega_{2(s-q) + 1} (s-q,1, q-1) 
    +  \sum_{q =1}^{s-2}  4^{2q}\cdot\Delta _{k}\omega_{2(s-q) } (s-q,0, q)  \\
    & = 2^{4s-6}\cdot[3\cdot2^{2k+1} - 17\cdot2^{k+1}]  -9\cdot2^{2k+2s-2} +2^2\cdot[-9\cdot2^{2k+2s-3}] \\
    &   +  \sum_{q =2}^{s-2}  2^{4q-2}\cdot[-3\cdot2^{2k+2(s-q)+1-2}]
    + \sum_{q =1}^{s-3}  2^{4q}\cdot[-3\cdot2^{2k+2(s-q)-2}] + 2^{4s-8}\cdot[-3\cdot2^{2k+2} +2^{k+2}]\\
    & = -15\cdot2^{2k+2s-2} -33\cdot2^{k+4s-6}
\end{align*}
\underline{The case j = 1} \\

We have 
\begin{align*}
  \Delta _{k}\omega_{2s+1} (s,0,0) & = 9\cdot2^{2k+2s-2}\\
  & \\
   \Delta _{k}\omega_{2(s-1)+2} (s-1,1,0) & = 9\cdot2^{2k+2s-3}\\
 & \\
  \Delta _{k}\omega_{2(s-2)+2} (s-2,1,1) &  = - 9\cdot2^{2k+2s-5}\\
 & \\
 \Delta _{k}\omega_{2(s-q)+2} (s-q,1,q-1) & = 15\cdot2^{2k+2(s-q) -1} & \text{if } 3\leq q\leq s-2 \\
& \\
 \Delta _{k}\omega_{2(s-1)+1} (s-1,0,1) & = - 9\cdot2^{2k+2s-4}\\
 & \\
  \Delta _{k}\omega_{2(s-q)+1} (s-q,0,q) & = 15\cdot2^{2k+2(s-q) -2} & \text{if } 2\leq q\leq s-2 
\end{align*}
 \begin{align*}
&  \Delta _{k}\Gamma_{2s+1}^{\left[s\atop{ s \atop s } \right]}
 -  8\cdot \Delta _{k}\Gamma_{2s }^{\left[s -1\atop{ s \atop s } \right]}
  =  2^{4s-6}\cdot\left[ \Delta _{k}\Gamma_{4}^{\left[1 \atop{ 2 \atop s }\right]} 
 - 8\cdot \Delta _{k}\Gamma_{3}^{\left[1 \atop{ 2  \atop s -1} \right]}\right]\\
    &    + \Delta _{k}\omega_{2s+1} (s,0,0) + \sum_{q =1}^{s-2}  2^{4q-2}\cdot\Delta _{k}\omega_{2(s-q) + 2} (s-q,1, q-1) 
    +  \sum_{q =1}^{s-2}  2^{4q}\cdot\Delta _{k}\omega_{2(s-q) + 1} (s-q,0, q) \\
    & =    2^{4s-6}\cdot( 3\cdot2^{2k+1} - 17 \cdot2^{k+1}) 
     +  9\cdot2^{2k+2s-2} + 2^{2}\cdot 9\cdot2^{2k+2s - 3} + 2^{6}\cdot( - 9\cdot2^{2k+2s-5}) \\
    & +  \sum_{q =3}^{s-2}  2^{4q-2}\cdot 15\cdot2^{2k+2(s-q) -1} + 2^{4}\cdot( - 9\cdot2^{2k+2s-4}) 
     +  \sum_{q =2}^{s-2}  2^{4q}\cdot15\cdot2^{2k+2(s-q) -2}\\
     & = -15\cdot2^{2k+4s-5} +3\cdot2^{k+4s-3} - 81\cdot2^{2k+2s- 2} + 15\cdot2^{2k+4s-5}
      - 5\cdot2^{2k+2s+3}  - 5\cdot2^{2k+2s+2}\\
      &  = -321\cdot2^{2k+2s-2} + 3\cdot2^{k+4s-3}
\end{align*}

\underline{The case  $ 2\leq j\leq s-4 $}\\

     \begin{align*}
&  \Delta _{k}\Gamma_{2s+j}^{\left[s\atop{ s \atop s } \right]}
 -  8\cdot \Delta _{k}\Gamma_{2s+j -1}^{\left[s -1\atop{ s \atop s } \right]}\\
  & =  2^{4s-6}\cdot\left[ \Delta _{k}\Gamma_{j + 3}^{\left[1 \atop{ 2 \atop s }\right]} 
 - 8\cdot \Delta _{k}\Gamma_{j+2}^{\left[1 \atop{ 2  \atop s -1} \right]}\right]\\
    &    + \Delta _{k}\omega_{2s+j} (s,0,0) + \sum_{q =1}^{s-2}  2^{4q-2}\cdot\Delta _{k}\omega_{2(s-q) + (j+1)} (s-q,1, q-1) 
    +  \sum_{q =1}^{s-2}  2^{4q}\cdot\Delta _{k}\omega_{2(s-q) + j} (s-q,0, q)  
\end{align*}
We have \\
\begin{align*}
\Delta _{k}\omega_{2s+j} (s,0,0) = 0 \quad  & \text{if    }  2\leq j\leq s - 4,
\end{align*}

   \[ \Delta _{k} \omega_{2(s-q) + (j+1)} (s-q,1, q-1) =
\begin{cases}
0  &\text{if } 1\leq q\leq j-1 \\
3\cdot2^{2k+2s-3}  &\text{if  }  q = j \\
- 15\cdot2^{2k+ 2s-5}  &\text{if  } q = j+1 \\
9\cdot2^{2k+ 2s -2q +2j-3}  &\text{if  } q\geq j+2, 
\end{cases}\]\\         
 
  \[ \Delta _{k} \omega_{2(s-q) + j } (s-q,0, q) =
\begin{cases}
0  &\text{if } 1\leq q\leq j-2 \\
3\cdot2^{2k+2s-2}  &\text{if  }  q = j-1 \\
- 15\cdot2^{2k+ 2s-4}  &\text{if  } q = j \\
9\cdot2^{2k+ 2s -2q +2j-4}  &\text{if  } q\geq j+1. 
\end{cases}\]\\            
   \begin{align*}
& \sum_{q =1}^{s-2}  2^{4q-2}\cdot\Delta _{k}\omega_{2(s-q) + (j+1)} (s-q,1, q-1)  =   \sum_{q = j}^{s-2}  2^{4q-2}\cdot\Delta _{k}\omega_{2(s-q) + (j+1)} (s-q,1, q-1)\\
& = 2^{4j-2}\cdot(3\cdot2^{2k+2s-3}) +  2^{4(j+1)-2}\cdot(-15\cdot2^{2k+2s-5}) +  \sum_{q = j+2}^{s-2}  2^{4q-2}\cdot9\cdot2^{2k+2s-2q+2j-3}\\
& = 3\cdot2^{2k+2s +4j -5}  - 15\cdot2^{2k+2s +4j-3} + 9\cdot2^{2k+ 2s +2j-5}\sum_{q = j+2}^{s-2}  2^{2q}
 =  3\cdot2^{2k+4s +2j -7}  - 105\cdot2^{2k+2s +4j-5}
\end{align*}
and \\
  \begin{align*}
& \sum_{q =1}^{s-2}  2^{4q}\cdot\Delta _{k}\omega_{2(s-q) + j} (s-q,0, q)  =   \sum_{q = j-1}^{s-2}  2^{4q}\cdot\Delta _{k}\omega_{2(s-q) + j} (s-q,0, q)\\
& = 2^{4(j-1)}\cdot(3\cdot2^{2k+2s-2}) +  2^{4j}\cdot(-15\cdot2^{2k+2s-4}) +  \sum_{q = j+1}^{s-2}  2^{4q}\cdot9\cdot2^{2k+2s-2q+2j-4}\\
& = 3\cdot2^{2k+2s +4j -6}  - 15\cdot2^{2k+2s +4j-4} + 9\cdot2^{2k+ 2s +2j-4}\sum_{q = j+1}^{s-2}  2^{2q}
 =  3\cdot2^{2k+4s +2j -6}  - 105\cdot2^{2k+2s +4j- 6}
\end{align*}

Combining the above results we obtain for  $ 2\leq j\leq s-4 $ \\
\begin{align*}
&  \Delta _{k}\Gamma_{2s+j}^{\left[s\atop{ s \atop s } \right]}
 -  8\cdot \Delta _{k}\Gamma_{2s+j -1}^{\left[s -1\atop{ s \atop s } \right]}\\
  & =  2^{4s-6}\cdot\left[ \Delta _{k}\Gamma_{j + 3}^{\left[1 \atop{ 2 \atop s }\right]} 
 - 8\cdot \Delta _{k}\Gamma_{j+2}^{\left[1 \atop{ 2  \atop s -1} \right]}\right]\\
    &    + \Delta _{k}\omega_{2s+j} (s,0,0) + \sum_{q =1}^{s-2}  2^{4q-2}\cdot\Delta _{k}\omega_{2(s-q) + (j+1)} (s-q,1, q-1) 
    +  \sum_{q =1}^{s-2}  2^{4q}\cdot\Delta _{k}\omega_{2(s-q) + j} (s-q,0, q) \\
    & =   -9\cdot2^{2k+4s +2j -7} + ( 3\cdot2^{2k+4s +2j -7}  - 105\cdot2^{2k+2s +4j-5} )+ ( 3\cdot2^{2k+4s +2j -6}  - 105\cdot2^{2k+2s +4j- 6}) \\
    & = -315\cdot2^{2k+2s +4j-6}
\end{align*}
  \underline{The case j = s -3}
     \begin{align*}
&  \Delta _{k}\Gamma_{2s+ (s-3)}^{\left[s\atop{ s \atop s } \right]}
 -  8\cdot \Delta _{k}\Gamma_{2s+ (s-4)}^{\left[s -1\atop{ s \atop s } \right]}\\
  & =  2^{4s-6}\cdot\left[ \Delta _{k}\Gamma_{s}^{\left[1 \atop{ 2 \atop s }\right]} 
 - 8\cdot \Delta _{k}\Gamma_{s-1}^{\left[1 \atop{ 2  \atop s -1} \right]}\right]\\
    &    + \Delta _{k}\omega_{2s+ (s-3)} (s,0,0) + \sum_{q =1}^{s-2}  2^{4q-2}\cdot\Delta _{k}\omega_{2(s-q) + (s-2)} (s-q,1, q-1) 
    +  \sum_{q =1}^{s-2}  2^{4q}\cdot\Delta _{k}\omega_{2(s-q) + (s-3)} (s-q,0, q)  
\end{align*}
We have \\
\begin{align*}
\Delta _{k}\omega_{2s+(s-3)} (s,0,0) = 0 \quad  & \text{if    }\; j = s -3,
\end{align*}
  \[ \Delta _{k} \omega_{2(s-q) + (s-2)} (s-q,1, q-1) =
\begin{cases}
0  &\text{if } 1\leq q\leq s-4 \\
3\cdot2^{2k+2s-3}  &\text{if  }  q = s-3 \\
- 15\cdot2^{2k+ 2s-5}  &\text{if  } q = s-2, 
\end{cases}\]\\         
 
  \[ \Delta _{k} \omega_{2(s-q) + (s-3) } (s-q,0, q) =
\begin{cases}
0  &\text{if } 1\leq q\leq  s -5 \\
3\cdot2^{2k+2s-2}  &\text{if  }  q = s - 4 \\
- 15\cdot2^{2k+ 2s-4}  &\text{if  } q = s -3 \\
9\cdot2^{2k+ 2s - 6}  &\text{if  } q = s -2,
\end{cases}\]\\            
   \begin{align*}
& \sum_{q =1}^{s-2}  2^{4q-2}\cdot\Delta _{k}\omega_{2(s-q) + (s-2)} (s-q,1, q-1)  =   \sum_{q = s-3}^{s-2}  2^{4q-2}\cdot\Delta _{k}\omega_{2(s-q) + (s-2)} (s-q,1, q-1)\\
& = 2^{4(s-3)-2}\cdot(3\cdot2^{2k+2s-3}) +  2^{4(s-2)-2}\cdot(-15\cdot2^{2k+2s-5}) 
 =  -57\cdot2^{2k+6s -17},  
\end{align*}

  \begin{align*}
& \sum_{q =1}^{s-2}  2^{4q}\cdot\Delta _{k}\omega_{2(s-q) + (s-3)} (s-q,0, q)  =   \sum_{q = s-4}^{s-2}  2^{4q}\cdot\Delta _{k}\omega_{2(s-q) + (s-3)} (s-q,0, q)\\
& = 2^{4(s-4)}\cdot(3\cdot2^{2k+2s-2}) +  2^{4(s-3)}\cdot(-15\cdot2^{2k+2s-4}) +  2^{4(s-2)}\cdot(9\cdot2^{2k+2s-6}) 
 =  87 \cdot2^{2k+6s  - 18}. \\
 & 
\end{align*}
As before we deduce from the above results \\
    \begin{align*}
&  \Delta _{k}\Gamma_{2s+ (s-3)}^{\left[s\atop{ s \atop s } \right]}
 -  8\cdot \Delta _{k}\Gamma_{2s+ (s-4)}^{\left[s -1\atop{ s \atop s } \right]}
   =  2^{4s-6}\cdot\left[ \Delta _{k}\Gamma_{s}^{\left[1 \atop{ 2 \atop s }\right]} 
 - 8\cdot \Delta _{k}\Gamma_{s-1}^{\left[1 \atop{ 2  \atop s -1} \right]}\right]\\
    &    + \Delta _{k}\omega_{2s+ (s-3)} (s,0,0) + \sum_{q =1}^{s-2}  2^{4q-2}\cdot\Delta _{k}\omega_{2(s-q) + (s-2)} (s-q,1, q-1) 
    +  \sum_{q =1}^{s-2}  2^{4q}\cdot\Delta _{k}\omega_{2(s-q) + (s-3)} (s-q,0, q)  \\
    & =  -9\cdot2^{2k+ 6s - 13} -57\cdot2^{2k+6s -17}  + 87 \cdot2^{2k+6s - 18} 
     = -315\cdot2^{2k+6s-18}
\end{align*}
  \underline{The case j = s -2}
    \begin{align*}
&  \Delta _{k}\Gamma_{2s+ (s-2)}^{\left[s\atop{ s \atop s } \right]}
 -  8\cdot \Delta _{k}\Gamma_{2s+ (s-3)}^{\left[s -1\atop{ s \atop s } \right]}\\
  & =  2^{4s-6}\cdot\left[ \Delta _{k}\Gamma_{s+1}^{\left[1 \atop{ 2 \atop s }\right]} 
 - 8\cdot \Delta _{k}\Gamma_{s}^{\left[1 \atop{ 2  \atop s -1} \right]}\right]\\
    &    + \Delta _{k}\omega_{2s+ (s-2)} (s,0,0) + \sum_{q =1}^{s-2}  2^{4q-2}\cdot\Delta _{k}\omega_{2(s-q) + (s-1)} (s-q,1, q-1) 
    +  \sum_{q =1}^{s-2}  2^{4q}\cdot\Delta _{k}\omega_{2(s-q) + (s-2)} (s-q,0, q)  
\end{align*}
We have \\
\begin{align*}
\Delta _{k}\omega_{2s+(s-2)} (s,0,0) = 0 \quad  & \text{if    }\; j = s-2,
\end{align*}
   \[ \Delta _{k} \omega_{2(s-q) + (s-1)} (s-q,1, q-1) =
\begin{cases}
0  &\text{if } 1\leq q\leq s-3 \\
3\cdot2^{2k+2s-3}  &\text{if  }  q = s-2, 
\end{cases}\]\\         
   \[ \Delta _{k}  \omega_{2(s-q) + (s-2) } (s-q,0, q) =
\begin{cases}
0  &\text{if } 1\leq q\leq  s -4 \\
3\cdot2^{2k+2s-2}  &\text{if  }  q = s - 3 \\
- 15\cdot2^{2k+ 2s-4}  &\text{if  } q = s - 2, 
\end{cases}\]\\            
  \begin{align*}
& \sum_{q =1}^{s-2}  2^{4q-2}\cdot\Delta _{k}\omega_{2(s-q) + (s-1)} (s-q,1, q-1)  =   \sum_{q = s-2}^{s-2}  2^{4q-2}\cdot\Delta _{k}\omega_{2(s-q) + (s-1)} (s-q,1, q-1)\\
& = 2^{4(s-2)-2}\cdot(3\cdot2^{2k+2s-3}) =  3\cdot2^{2k+6s -13},
\end{align*}

  \begin{align*}
& \sum_{q =1}^{s-2}  2^{4q}\cdot\Delta _{k}\omega_{2(s-q) + (s-2)} (s-q,0, q)  =   \sum_{q = s-3}^{s-2}  2^{4q}\cdot\Delta _{k}\omega_{2(s-q) + (s-2)} (s-q,0, q)\\
& = 2^{4(s-3)}\cdot(3\cdot2^{2k+2s-2}) +  2^{4(s-2)}\cdot(-15\cdot2^{2k+2s-4}) = 3\cdot2^{2k+6s-14} -15\cdot2^{2k+6s - 12} 
 =  -57 \cdot2^{2k+6s  - 14}. 
\end{align*}

As before we deduce from the above results \\
    \begin{align*}
&  \Delta _{k}\Gamma_{2s+ (s-2)}^{\left[s\atop{ s \atop s } \right]}
 -  8\cdot \Delta _{k}\Gamma_{2s+ (s-3)}^{\left[s -1\atop{ s \atop s } \right]}\\
  & =  2^{4s-6}\cdot\left[ \Delta _{k}\Gamma_{s+1}^{\left[1 \atop{ 2 \atop s }\right]} 
 - 8\cdot \Delta _{k}\Gamma_{s}^{\left[1 \atop{ 2  \atop s -1} \right]}\right]\\
    &    + \Delta _{k}\omega_{2s+ (s-2)} (s,0,0) + \sum_{q =1}^{s-2}  2^{4q-2}\cdot\Delta _{k}\omega_{2(s-q) + (s-1)} (s-q,1, q-1) 
    +  \sum_{q =1}^{s-2}  2^{4q}\cdot\Delta _{k}\omega_{2(s-q) + (s-2)} (s-q,0, q)  \\
 & =     -33\cdot2^{2k+ 6s - 11} + 3\cdot2^{2k+6s -13}  -57 \cdot2^{2k+6s  - 14} \\
 & = -315\cdot2^{2k+6s -14}
 \end{align*} 
    
 \underline{The case j = s }

    \begin{align*}
&  \Delta _{k}\Gamma_{2s+ s}^{\left[s\atop{ s \atop s } \right]}
 -  8\cdot \Delta _{k}\Gamma_{2s+ (s-1)}^{\left[s -1\atop{ s \atop s } \right]}\\
  & =  2^{4s-6}\cdot\left[ \Delta _{k}\Gamma_{s+3}^{\left[1 \atop{ 2 \atop s }\right]} 
 - 8\cdot \Delta _{k}\Gamma_{s+2}^{\left[1 \atop{ 2  \atop s -1} \right]}\right]\\
    &    + \Delta _{k}\omega_{3s} (s,0,0) + \sum_{q =1}^{s-2}  2^{4q-2}\cdot\Delta _{k}\omega_{2(s-q) + s +1} (s-q,1, q-1) 
    +  \sum_{q =1}^{s-2}  2^{4q}\cdot\Delta _{k}\omega_{2(s-q) + s} (s-q,0, q)  
\end{align*}
We have \\
\begin{align*}
    \Delta _{k}\omega_{3s} (s,0,0) & = 21\cdot2^{3k+3s-3},\\
    \Delta _{k} \omega_{2(s-q) + s+1} (s-q,1, q-1) &  = 21\cdot2^{3k+3(s-q) +q-2} \quad  &\text{if } 1\leq q\leq s-2, \\
   \Delta _{k} \omega_{2(s-q) + s} (s-q,0, q) &  = 21\cdot2^{3k+3(s-q) +q-3}\quad   &\text{if } 1\leq q\leq s-2, 
\end{align*}
 \begin{align*}
& \sum_{q =1}^{s-2}  2^{4q-2}\cdot\Delta _{k}\omega_{2(s-q) + s +1} (s-q,1, q-1)  =   \sum_{q =1}^{s-2}  2^{4q-2}\cdot21\cdot2^{3k+3(s-q) +q-2}\\
& = 21\cdot2^{3k+3s-4}\cdot \sum_{q =1}^{s-2}  2^{2q}=  21\cdot2^{3k+3s-2}\cdot \sum_{q =0}^{s-3}  2^{2q}
 = 7\cdot2^{3k+5s-6} -7\cdot2^{3k+3s-2},
 \end{align*}

 \begin{align*}
& \sum_{q =1}^{s-2}  2^{4q-2}\cdot\Delta _{k}\omega_{2(s-q) + s} (s-q,0, q)  =   \sum_{q =1}^{s-2}  2^{4q}\cdot21\cdot2^{3k+3(s-q) +q- 3}\\
& = 21\cdot2^{3k+3s-3}\cdot \sum_{q =1}^{s-2}  2^{2q}=  21\cdot2^{3k+3s-1}\cdot \sum_{q =0}^{s-3}  2^{2q}
 = 7\cdot2^{3k+5s-5} -7\cdot2^{3k+3s-1}.
 \end{align*}

As before we deduce from the above results \\
    \begin{align*}
&  \Delta _{k}\Gamma_{2s+ s}^{\left[s\atop{ s \atop s } \right]}
 -  8\cdot \Delta _{k}\Gamma_{2s+ (s-1)}^{\left[s -1\atop{ s \atop s } \right]}\\
  & =  2^{4s-6}\cdot\left[ \Delta _{k}\Gamma_{s+3}^{\left[1 \atop{ 2 \atop s }\right]} 
 - 8\cdot \Delta _{k}\Gamma_{s+2}^{\left[1 \atop{ 2  \atop s -1} \right]}\right]\\
    &    + \Delta _{k}\omega_{3s} (s,0,0) + \sum_{q =1}^{s-2}  2^{4q-2}\cdot\Delta _{k}\omega_{2(s-q) + s+1} (s-q,1, q-1) 
     +  \sum_{q =1}^{s-2}  2^{4q}\cdot\Delta _{k}\omega_{2(s-q) + s} (s-q,0, q)  \\
& =   - 21\cdot2^{3k+ 5s - 6} + 21\cdot2^{2k+ 6s - 6} +  21\cdot2^{3k+3s-3} + 7\cdot2^{3k+5s-6} -7\cdot2^{3k+3s-2} + 7\cdot2^{3k+5s-5} -7\cdot2^{3k+3s-1}\\
   & = 21\cdot2^{2k+6s-6} -21\cdot2^{3k+3s-3}
    \end{align*}
    \end{proof}
 
   \section{\textbf{An explicit formula for $ \Delta _{k}\Gamma_{i}^{\left[s\atop{ s\atop s} \right]} , \quad for\; 0\leq i \leq 3s $}}
\label{sec 9} 
    \subsection{Notation}
  \label{subsec 1}
  \begin{defn}
  \label{defn 9.1}

  We recall the formula obtained in Section \ref{sec 8} (see \eqref{eq 8.3}):\\
  \begin{align*}
  \Delta_{k}\omega _{i}(s,m,l)  = \Delta _{k}\Gamma_{i}^{\left[s\atop{ s+m\atop s+m+l} \right]}
- 4\cdot\Delta_{k} \Gamma_{i-1}^{\left[s\atop{ s+m-1\atop s+m+l} \right]} -
8\cdot\Delta _{k} \Gamma_{i-1}^{\left[s\atop{ s+m\atop s+m+l-1} \right]} +
32\cdot\Delta _{k}\Gamma_{i-2}^{\left[s\atop{ s+m-1\atop s+m+l-1} \right]} 
\end{align*}

  \end{defn}
  
   \subsection{Introduction}
  \label{subsec 2}
 Using the results concerning  $ \Delta_{k}\omega _{i}(s,m,l) $ obtained in Lemma \ref{lem 7.4} we adapt the method 
 of Section 12 in [2] to compute  $ \Delta _{k}\Gamma_{i}^{\left[s\atop{ s\atop s} \right]} \; \text{for}\; 0\leqslant i \leqslant 3s$
 applying successively the formula \ref{eq 8.3}.

  \subsection{\textbf{Computation of $    \Delta _{k}\Gamma_{s+j}^{\left[s\atop{ s+m\atop s+m+l} \right]}  =  \Gamma_{s+j}^{\left[s\atop{ s+m\atop s+m+l} \right]\times (k+1)} - \Gamma_{s+j}^{\left[s\atop{ s+m\atop s+m+l} \right]\times k}$\quad \text{for j =0,1,2}}}
\label{subsec 3}
   \begin{lem} We have \\
  \label{lem 9.2}
  \begin{equation}
  \label{eq 9.1}
 \Delta _{k}\Gamma_{s}^{\left[s\atop{ s+m\atop s+m+l} \right]} =  
\begin{cases}
2^{k+s-1}  &\text{if }  m\geq 1, l\geq 0\\
3\cdot2^{k+s-1}  &\text{if }  m = 0, l\geq 1\\
7\cdot2^{k+s-1}  &\text{if }  m = 0, l = 0,
\end{cases}
   \end{equation}
\begin{equation}
\label{eq 9.2}
 \Delta _{k}\Gamma_{s+1}^{\left[s\atop{ s+m\atop s+m+l} \right]}=
\begin{cases}
9\cdot2^{k+s-1}  &\text{if }  m\geq 2, l\geq 0\\
17\cdot2^{k+s-1}  &\text{if }  m = 1, l\geq 1\\
33\cdot2^{k+s-1}  &\text{if }  m = 1, l = 0 \\
147\cdot2^{k+s-1}  &\text{if }  m = 0, l = 0 \\
71\cdot2^{k+s-1}  &\text{if }  m = 0, l = 1 \\
39\cdot2^{k+s-1}  &\text{if }  m = 0, l \geq 2, 
\end{cases} 
\end{equation}
\begin{equation}
\label{eq 9.3}
 \Delta _{k}\Gamma_{s+2}^{\left[s\atop{ s+m\atop s+m+l} \right]}=
\begin{cases}
78\cdot2^{k+s-1}  &\text{if }  m\geq 3, l\geq 0\\
110\cdot2^{k+s-1}  &\text{if }  m = 2, l\geq 1\\
174\cdot2^{k+s-1}  &\text{if }  m = 2, l = 0 \\
198\cdot2^{k+s-1}  &\text{if }  m = 1, l \geq 2 \\
326\cdot2^{k+s-1}  &\text{if }  m = 1, l = 1 \\
630\cdot2^{k+s-1}  &\text{if }  m = 1, l = 0 \\
438\cdot2^{k+s-1}  &\text{if }  m = 0, l \geq 3 \\
694\cdot2^{k+s-1}  &\text{if }  m = 0, l = 2 \\
1302\cdot2^{k+s-1}  &\text{if }  m = 0, l = 1 \\
2646\cdot2^{k+s-1}  &\text{if }  m = 0, l = 0. 
\end{cases}
  \end{equation}
  \end{lem}
  \begin{proof}
  The proof is based on  Lemma \ref{lem 7.4}  and  \eqref{eq 6.33}. \\
  We remark that to obtain the result in one case we are using some of the results in the precedent cases.\\
   \underline{proof of \eqref{eq 9.1}}
    \begin{align*}
 & \text{The case  $ m\geq 1, l\geq 0 $}   \\
&   \Delta _{k}\Gamma_{s}^{\left[s\atop{ s+m\atop s+m+l} \right]} -4\cdot \Delta _{k}\Gamma_{s-1}^{\left[s\atop{ s+m-1\atop s+m+l} \right]}
-8\cdot\Delta _{k}\Gamma_{s-1}^{\left[s\atop{ s+m\atop s+m+l-1} \right]}+32\cdot\Delta _{k}\Gamma_{s-2}^{\left[s\atop{ s+m-1\atop s+m+l-1} \right]}= \Delta _{k} \omega _{s}(s,m,l) = 2^{k+s-1}\\
& \text{since} \quad \Delta _{k}\Gamma_{s-1}^{\left[s\atop{ s+m-1\atop s+m+l}\right]}=  \Delta _{k}\Gamma_{s-1}^{\left[s\atop{ s+m\atop s+m+l-1} \right]} = \Delta _{k}\Gamma_{s-2}^{\left[s\atop{ s+m-1\atop s+m+l-1} \right]} = 0 
\quad  \text{we get}\quad \Delta _{k}\Gamma_{s}^{\left[s\atop{ s+m\atop s+m+l} \right]} =  2^{k+s-1}\\
& \\
 & \text{The case  $ m = 0, l\geq 1 $}   \\
 &   \Delta _{k}\Gamma_{s}^{\left[s\atop{ s \atop s+l} \right]} -4\cdot \Delta _{k}\Gamma_{s-1}^{\left[s -1\atop{ s\atop s+l} \right]}
-8\cdot\Delta _{k}\Gamma_{s-1}^{\left[s\atop{ s\atop s+l-1} \right]}+32\cdot\Delta _{k}\Gamma_{s-2}^{\left[s-1\atop{ s\atop s+l-1} \right]}= \Delta _{k} \omega _{s}(s,0,l) = 2^{k+s-1}\\
& \text{since} \quad \Delta _{k}\Gamma_{s-1}^{\left[s\atop{ s\atop s+(l-1)}\right]}=  \Delta _{k}\Gamma_{s-2}^{\left[s-1\atop{ s\atop s+(l-1)} \right]} = 0 \; \text{and} \; \Delta _{k}\Gamma_{s-1}^{\left[s-1\atop{ s\atop s+l} \right]} =  2^{k+s-2}
\quad  \text{we get}\quad \Delta _{k}\Gamma_{s}^{\left[s\atop{ s\atop s+l} \right]} = 3\cdot 2^{k+s-1}\\
&  \\
 & \text{The case  $ m = 0, l = 0 $}   \\
 &   \Delta _{k}\Gamma_{s}^{\left[s\atop{ s \atop s} \right]} -4\cdot \Delta _{k}\Gamma_{s-1}^{\left[s -1\atop{ s\atop s} \right]}
-8\cdot\Delta _{k}\Gamma_{s-1}^{\left[s-1\atop{ s\atop s} \right]}+32\cdot\Delta _{k}\Gamma_{s-2}^{\left[s-1\atop{ s-1\atop s} \right]}= \Delta _{k} \omega _{s}(s,0,0) = 2^{k+s-1}\\
& \text{since} \quad \Delta _{k}\Gamma_{s-1}^{\left[s-1\atop{ s\atop s}\right]}= 2^{k+s-2}, \quad  \Delta _{k}\Gamma_{s-2}^{\left[s-1\atop{ s-1\atop s} \right]} = 0
\quad  \text{we get}\quad \Delta _{k}\Gamma_{s}^{\left[s\atop{ s\atop s} \right]} = 2^{k+s-1} +12\cdot 2^{k+s-2} = 7\cdot2^{k+s-1}
 \end{align*}
  \underline{proof of \eqref{eq 9.2}} 
 \begin{align*}
  & \text{The case  $ m \geq  2, l \geq  0 $}   \\
 &   \Delta _{k}\Gamma_{s+1}^{\left[s\atop{ s+m\atop s+m+l} \right]} -4\cdot \Delta _{k}\Gamma_{s}^{\left[s\atop{ s+m-1\atop s+m+l} \right]}
-8\cdot\Delta _{k}\Gamma_{s}^{\left[s\atop{ s+m\atop s+m+l-1} \right]}+32\cdot\Delta _{k}\Gamma_{s-1}^{\left[s\atop{ s+m-1\atop s+m+l-1} \right]}= \Delta _{k} \omega _{s+1}(s,m,l) =-3\cdot2^{k+s-1}\\
& \text{since} \quad  \Delta _{k}\Gamma_{s}^{\left[s\atop{ s+m-1\atop s+m+l} \right]} = \Delta _{k}\Gamma_{s}^{\left[s\atop{ s+m\atop s+m+l-1} \right]}= 2^{k+s-1}\; \text{and} \; \Delta _{k}\Gamma_{s-1}^{\left[s\atop{ s+m-1\atop s+m+l-1} \right]}= 0\\
& \quad  \text{we get}\quad  \Delta _{k}\Gamma_{s+1}^{\left[s\atop{ s+m\atop s+m+l} \right]}= 12\cdot2^{k+s-1}-3\cdot2^{k+s-1} = 9\cdot2^{k+s-1} \\
& \\
  & \text{The case  $ m=1, l \geq  1 $}   \\
 &   \Delta _{k}\Gamma_{s+1}^{\left[s\atop{ s+1\atop s+1+l} \right]} -4\cdot \Delta _{k}\Gamma_{s}^{\left[s\atop{ s\atop s+1+l} \right]}
-8\cdot\Delta _{k}\Gamma_{s}^{\left[s\atop{ s+1\atop s+l} \right]}+32\cdot\Delta _{k}\Gamma_{s-1}^{\left[s\atop{ s\atop s+l} \right]}= \Delta _{k} \omega _{s+1}(s,1,l) =-3\cdot2^{k+s-1}\\
& \text{since} \quad  \Delta _{k}\Gamma_{s}^{\left[s\atop{ s\atop s+1+l} \right]}=3\cdot2^{k+s-1},\quad \Delta _{k}\Gamma_{s}^{\left[s\atop{ s+1\atop s+l} \right]}= 2^{k+s-1}\; \text{and} \; \Delta _{k}\Gamma_{s-1}^{\left[s\atop{ s\atop s+l} \right]}= 0\\
& \quad  \text{we get}\quad  \Delta _{k}\Gamma_{s+1}^{\left[s\atop{ s+m\atop s+m+l} \right]}= 12\cdot2^{k+s-1}+8\cdot2^{k+s-1}-3\cdot2^{k+s-1} = 17\cdot2^{k+s-1} \\
& \\
  & \text{The case  $ m = 1, l =  0 $}   \\
 &   \Delta _{k}\Gamma_{s+1}^{\left[s\atop{ s+1\atop s+1} \right]} -4\cdot \Delta _{k}\Gamma_{s}^{\left[s\atop{ s\atop s+1} \right]}
-8\cdot\Delta _{k}\Gamma_{s}^{\left[s\atop{ s\atop s+1} \right]}+32\cdot\Delta _{k}\Gamma_{s-1}^{\left[s\atop{ s\atop s} \right]}= \Delta _{k} \omega _{s+1}(s,1,0) =-3\cdot2^{k+s-1}\\
& \text{since} \quad  \Delta _{k}\Gamma_{s}^{\left[s\atop{ s\atop s+1} \right]} = 3\cdot2^{k+s-1}\; \text{and} \; \Delta _{k}\Gamma_{s-1}^{\left[s\atop{ s\atop s} \right]}= 0
 \quad  \text{we get}\quad  \Delta _{k}\Gamma_{s+1}^{\left[s\atop{ s+1\atop s+1} \right]}= 36\cdot2^{k+s-1}-3\cdot2^{k+s-1} = 33\cdot2^{k+s-1} \\
& \\
  & \text{The case  $ m = 0, l =  0 $}   \\
 &   \Delta _{k}\Gamma_{s+1}^{\left[s\atop{ s\atop s} \right]} -4\cdot \Delta _{k}\Gamma_{s}^{\left[s-1\atop{ s\atop s} \right]}
-8\cdot\Delta _{k}\Gamma_{s}^{\left[s-1\atop{ s\atop s} \right]}+32\cdot\Delta _{k}\Gamma_{s-1}^{\left[s-1\atop{ s-1\atop s} \right]}= \Delta _{k} \omega _{s+1}(s,0,0) =-3\cdot2^{k+s-1}\\
& \text{since} \quad  \Delta _{k}\Gamma_{s}^{\left[s-1\atop{ s\atop s} \right]} = 33\cdot2^{k+s-2}\; \text{and} \; \Delta _{k}\Gamma_{s-1}^{\left[s-1\atop{ s-1\atop s} \right]}= 3\cdot2^{k+s-2}
 \quad  \text{we get}\quad  \Delta _{k}\Gamma_{s+1}^{\left[s\atop{ s\atop s} \right]}\\
  & = 12\cdot33\cdot2^{k+s-2}-32\cdot3\cdot2^{k+s-2}-3\cdot2^{k+s-1} = 147\cdot2^{k+s-1} \\
& \\
  & \text{The case  $ m = 0, l =  1 $}   \\
 &   \Delta _{k}\Gamma_{s+1}^{\left[s\atop{ s\atop s+1} \right]} -4\cdot \Delta _{k}\Gamma_{s}^{\left[s-1\atop{ s\atop s+1} \right]}
-8\cdot\Delta _{k}\Gamma_{s}^{\left[s\atop{ s\atop s} \right]}+32\cdot\Delta _{k}\Gamma_{s-1}^{\left[s-1\atop{ s\atop s} \right]}= \Delta _{k} \omega _{s+1}(s,0,1) =-3\cdot2^{k+s-1}\\
& \text{since} \quad  \Delta _{k}\Gamma_{s}^{\left[s-1\atop{ s\atop s+1} \right]} = 17\cdot2^{k+s-2},\quad \Delta _{k}\Gamma_{s}^{\left[s\atop{ s\atop s} \right]}= 7\cdot2^{k+s-1}
 \quad \text{and}\quad \Delta _{k}\Gamma_{s-1}^{\left[s-1\atop{ s\atop s} \right]} =2^{k+s-2}\quad \text{we get}\\
&  \Delta _{k}\Gamma_{s+1}^{\left[s\atop{ s\atop s+1} \right]}
   = 4\cdot17\cdot2^{k+s-2}  + 8\cdot7\cdot2^{k+s-1}  -32\cdot2^{k+s-2} -3\cdot2^{k+s-1} = 71 \cdot2^{k+s-1} \\
  & \\
  & \text{The case  $ m = 0, l \geq  2 $}   \\
 &   \Delta _{k}\Gamma_{s+1}^{\left[s\atop{ s\atop s+l} \right]} -4\cdot \Delta _{k}\Gamma_{s}^{\left[s-1\atop{ s\atop s+l} \right]}
-8\cdot\Delta _{k}\Gamma_{s}^{\left[s\atop{ s\atop s+l-1} \right]}+32\cdot\Delta _{k}\Gamma_{s-1}^{\left[s-1\atop{ s\atop s+l-1} \right]}= \Delta _{k} \omega _{s+1}(s,0,l) =-3\cdot2^{k+s-1}\\
& \text{since} \quad  \Delta _{k}\Gamma_{s}^{\left[s-1\atop{ s\atop s+l} \right]} = 17\cdot2^{k+s-2},\quad \Delta _{k}\Gamma_{s}^{\left[s\atop{ s\atop s+l-1} \right]}= 3\cdot2^{k+s-1}
 \quad \text{and}\quad \Delta _{k}\Gamma_{s-1}^{\left[s-1\atop{ s\atop s+l-1} \right]} = 2^{k+s-2}\quad \text{we get}\\
&  \Delta _{k}\Gamma_{s+1}^{\left[s\atop{ s\atop s+l} \right]}
   = 4\cdot17\cdot2^{k+s-2}  + 8\cdot3\cdot2^{k+s-1}  -32\cdot2^{k+s-2} -3\cdot2^{k+s-1} = 39 \cdot2^{k+s-1} 
 \end{align*}
  \underline{proof of \eqref{eq 9.3}} 
  \begin{align*}
  & \text{The case  $ m \geq  3, l \geq  0 $}   \\
 &   \Delta _{k}\Gamma_{s+2}^{\left[s\atop{ s+m\atop s+m+l} \right]} -4\cdot \Delta _{k}\Gamma_{s+1}^{\left[s\atop{ s+m-1\atop s+m+l} \right]}
-8\cdot\Delta _{k}\Gamma_{s+1}^{\left[s\atop{ s+m\atop s+m+l-1} \right]}+32\cdot\Delta _{k}\Gamma_{s}^{\left[s\atop{ s+m-1\atop s+m+l-1} \right]}= \Delta _{k} \omega _{s+2}(s,m,l) = 2^{k+s}\\
& \text{since} \quad  \Delta _{k}\Gamma_{s+1}^{\left[s\atop{ s+m-1\atop s+m+l} \right]}
 = \Delta _{k}\Gamma_{s+1}^{\left[s\atop{ s+m\atop s+m+l-1} \right]}=9\cdot2^{k+s-1}\; \text{and} \; \Delta _{k}\Gamma_{s}^{\left[s\atop{ s+m-1\atop s+m+l-1} \right]}= 2^{k+s-1}\\
& \quad  \text{we get}\quad  \Delta _{k}\Gamma_{s+2}^{\left[s\atop{ s+m\atop s+m+l} \right]}
= 12\cdot9\cdot2^{k+s-1}-32\cdot2^{k+s-1}+2^{k+s} = 78\cdot2^{k+s-1} \\
& \\ 
  & \text{The case  $ m =2, l \geq  1 $}   \\
 &   \Delta _{k}\Gamma_{s+2}^{\left[s\atop{ s+2\atop s+2+l} \right]} -4\cdot \Delta _{k}\Gamma_{s+1}^{\left[s\atop{ s+1\atop s+2+l} \right]}
-8\cdot\Delta _{k}\Gamma_{s+1}^{\left[s\atop{ s+2\atop s+2+l-1} \right]}+32\cdot\Delta _{k}\Gamma_{s}^{\left[s\atop{ s+1\atop s+l+1} \right]}= \Delta _{k} \omega _{s+2}(s,2,l) = 2^{k+s}\\
& \text{since} \quad  \Delta _{k}\Gamma_{s+1}^{\left[s\atop{ s+1\atop s+2+l} \right]}= 17\cdot2^{k+s-1}, \quad
  \Delta _{k}\Gamma_{s+1}^{\left[s\atop{ s+2\atop s+l+1} \right]}=9\cdot2^{k+s-1}\; \text{and} \; \Delta _{k}\Gamma_{s}^{\left[s\atop{ s+1\atop s+l+1} \right]}= 2^{k+s-1}\\
& \quad  \text{we get}\quad  \Delta _{k}\Gamma_{s+2}^{\left[s\atop{ s+2\atop s+2+l} \right]}
= 4\cdot17\cdot2^{k+s-1}+ 8\cdot9\cdot2^{k+s-1} -32\cdot2^{k+s-1}+2^{k+s} = 110\cdot2^{k+s-1} \\
& \\ 
  & \text{The case  $ m = 2, l = 0 $}   \\
 &   \Delta _{k}\Gamma_{s+2}^{\left[s\atop{ s+2\atop s+2} \right]} -4\cdot \Delta _{k}\Gamma_{s+1}^{\left[s\atop{ s+1\atop s+2} \right]}
-8\cdot\Delta _{k}\Gamma_{s+1}^{\left[s\atop{ s+1\atop s+2} \right]}+32\cdot\Delta _{k}\Gamma_{s}^{\left[s\atop{ s+1\atop s+1} \right]}= \Delta _{k} \omega _{s+2}(s,2,0) = 2^{k+s}\\
& \text{since} \quad  \Delta _{k}\Gamma_{s+1}^{\left[s\atop{ s+1\atop s+2} \right]}=17\cdot2^{k+s-1}
\; \text{and} \; \Delta _{k}\Gamma_{s}^{\left[s\atop{ s+1\atop s+1} \right]}= 2^{k+s-1}\\
& \quad  \text{we get}\quad  \Delta _{k}\Gamma_{s+2}^{\left[s\atop{ s+2\atop s+2} \right]}
= 12\cdot17\cdot2^{k+s-1}-32\cdot2^{k+s-1}+2^{k+s} = 174\cdot2^{k+s-1} \\
& \\
  & \text{The case  $ m = 1, l \geq  2 $}   \\
 &   \Delta _{k}\Gamma_{s+2}^{\left[s\atop{ s+1\atop s+1+l} \right]} -4\cdot \Delta _{k}\Gamma_{s+1}^{\left[s\atop{ s\atop s+1+l} \right]}
-8\cdot\Delta _{k}\Gamma_{s+1}^{\left[s\atop{ s+1\atop s+l} \right]}+32\cdot\Delta _{k}\Gamma_{s}^{\left[s\atop{ s\atop s+l} \right]}= \Delta _{k} \omega _{s+2}(s,1,l) = 2^{k+s}\\
& \text{since} \quad  \Delta _{k}\Gamma_{s+1}^{\left[s\atop{ s\atop s+1+l} \right]}= 39\cdot2^{k+s-1},\quad
 \Delta _{k}\Gamma_{s+1}^{\left[s\atop{ s+1\atop s+l} \right]}=17\cdot2^{k+s-1}\; \text{and} \; \Delta _{k}\Gamma_{s}^{\left[s\atop{ s\atop s+l} \right]}=3\cdot2^{k+s-1}\\
& \quad  \text{we get}\quad  \Delta _{k}\Gamma_{s+2}^{\left[s\atop{ s+1\atop s+1+l} \right]}
= 4\cdot39\cdot2^{k+s-1} +8\cdot17\cdot2^{k+s-1} -32\cdot3 \cdot2^{k+s-1} +2^{k+s} = 198\cdot2^{k+s-1} \\
& \\ 
  & \text{The case  $ m = 1, l =1 $}   \\
 &   \Delta _{k}\Gamma_{s+2}^{\left[s\atop{ s+1\atop s+2} \right]} -4\cdot \Delta _{k}\Gamma_{s+1}^{\left[s\atop{ s\atop s+2} \right]}
-8\cdot\Delta _{k}\Gamma_{s+1}^{\left[s\atop{ s+1\atop s+1} \right]}+32\cdot\Delta _{k}\Gamma_{s}^{\left[s\atop{ s\atop s+1} \right]}= \Delta _{k} \omega _{s+2}(s,1,1) = 2^{k+s}\\
& \text{since} \quad  \Delta _{k}\Gamma_{s+1}^{\left[s\atop{ s\atop s+2} \right]}= 39\cdot2^{k+s-1},\quad
 \Delta _{k}\Gamma_{s+1}^{\left[s\atop{ s+1\atop s+1} \right]}=33\cdot2^{k+s-1}\; \text{and} \; \Delta _{k}\Gamma_{s}^{\left[s\atop{ s\atop s+1} \right]}=3\cdot2^{k+s-1}\\
& \quad  \text{we get}\quad  \Delta _{k}\Gamma_{s+2}^{\left[s\atop{ s+1\atop s+2} \right]}
= 4\cdot39\cdot2^{k+s-1} +8\cdot33\cdot2^{k+s-1} -32\cdot3 \cdot2^{k+s-1} +2^{k+s} = 326\cdot2^{k+s-1} \\
& \\
 & \text{The case  $ m = 1, l = 0 $}   \\
 &   \Delta _{k}\Gamma_{s+2}^{\left[s\atop{ s+1\atop s+1} \right]} -4\cdot \Delta _{k}\Gamma_{s+1}^{\left[s\atop{ s\atop s+1} \right]}
-8\cdot\Delta _{k}\Gamma_{s+1}^{\left[s\atop{ s\atop s+1} \right]}+32\cdot\Delta _{k}\Gamma_{s}^{\left[s\atop{ s\atop s} \right]}= \Delta _{k} \omega _{s+2}(s,1,0) = 2^{k+s}\\
& \text{since} \quad  \Delta _{k}\Gamma_{s+1}^{\left[s\atop{ s\atop s+1} \right]}= 71\cdot2^{k+s-1}, \quad
  \text{and} \; \Delta _{k}\Gamma_{s}^{\left[s\atop{ s\atop s} \right]}=7\cdot 2^{k+s-1}\\
& \quad  \text{we get}\quad  \Delta _{k}\Gamma_{s+2}^{\left[s\atop{ s+1\atop s+1} \right]}
= 12\cdot71\cdot2^{k+s-1}-32\cdot7\cdot2^{k+s-1}+2^{k+s} = 630\cdot2^{k+s-1} \\
& \\
  & \text{The case  $ m = 0, l \geq  3 $}   \\
 &   \Delta _{k}\Gamma_{s+2}^{\left[s\atop{ s\atop s+l} \right]} -4\cdot \Delta _{k}\Gamma_{s+1}^{\left[s-1\atop{ s\atop s+l} \right]}
-8\cdot\Delta _{k}\Gamma_{s+1}^{\left[s\atop{ s\atop s+l-1} \right]}+32\cdot\Delta _{k}\Gamma_{s}^{\left[s-1\atop{ s\atop s+l-1} \right]}= \Delta _{k} \omega _{s+2}(s,0,l) = 2^{k+s}\\
& \text{since} \quad  \Delta _{k}\Gamma_{s+1}^{\left[s-1\atop{ s\atop s+l} \right]}= 198\cdot2^{k+s-2}, \quad
  \Delta _{k}\Gamma_{s+1}^{\left[s\atop{ s\atop s+l-1} \right]}=39\cdot2^{k+s-1}\; \text{and} \; \Delta _{k}\Gamma_{s}^{\left[s-1\atop{ s\atop s+l-1} \right]}=17\cdot2^{k+s-2}\\
& \quad  \text{we get}\quad  \Delta _{k}\Gamma_{s+2}^{\left[s\atop{ s\atop s+l} \right]}
= 4\cdot198\cdot2^{k+s-2}+ 8\cdot39\cdot2^{k+s-1} -32\cdot17\cdot2^{k+s-2}+2^{k+s} = 438\cdot2^{k+s-1} \\
& \\
  & \text{The case  $ m = 0, l = 2 $}   \\
 &   \Delta _{k}\Gamma_{s+2}^{\left[s\atop{ s\atop s+2} \right]} -4\cdot \Delta _{k}\Gamma_{s+1}^{\left[s-1\atop{ s\atop s+2} \right]}
-8\cdot\Delta _{k}\Gamma_{s+1}^{\left[s\atop{ s\atop s+1} \right]}+32\cdot\Delta _{k}\Gamma_{s}^{\left[s-1\atop{ s\atop s+1} \right]}= \Delta _{k} \omega _{s+2}(s,0,2) = 2^{k+s}\\
& \text{since} \quad  \Delta _{k}\Gamma_{s+1}^{\left[s-1\atop{ s\atop s+2} \right]}= 198\cdot2^{k+s-2}, \quad
  \Delta _{k}\Gamma_{s+1}^{\left[s\atop{ s\atop s+1} \right]}=71\cdot2^{k+s-1}\; \text{and} \; \Delta _{k}\Gamma_{s}^{\left[s-1\atop{ s\atop s+1} \right]}=17\cdot2^{k+s-2}\\
& \quad  \text{we get}\quad  \Delta _{k}\Gamma_{s+2}^{\left[s\atop{ s\atop s+2} \right]}
= 4\cdot198\cdot2^{k+s-2}+ 8\cdot71\cdot2^{k+s-1} -32\cdot17\cdot2^{k+s-2}+2^{k+s} = 694\cdot2^{k+s-1} \\
& \\
  & \text{The case  $ m = 0, l=1 $}   \\
 &   \Delta _{k}\Gamma_{s+2}^{\left[s\atop{ s\atop s+1} \right]} -4\cdot \Delta _{k}\Gamma_{s+1}^{\left[s-1\atop{ s\atop s+1} \right]}
-8\cdot\Delta _{k}\Gamma_{s+1}^{\left[s\atop{ s\atop s} \right]}+32\cdot\Delta _{k}\Gamma_{s}^{\left[s-1\atop{ s\atop s} \right]}= \Delta _{k} \omega _{s+2}(s,0,1) = 2^{k+s}\\
& \text{since} \quad  \Delta _{k}\Gamma_{s+1}^{\left[s-1\atop{ s\atop s+1} \right]}= 326\cdot2^{k+s-2}, \quad
  \Delta _{k}\Gamma_{s+1}^{\left[s\atop{ s\atop s} \right]}=147\cdot2^{k+s-1}\; \text{and} \; \Delta _{k}\Gamma_{s}^{\left[s-1\atop{ s\atop s} \right]}=33\cdot2^{k+s-2}\\
& \quad  \text{we get}\quad  \Delta _{k}\Gamma_{s+2}^{\left[s\atop{ s\atop s+1} \right]}
= 4\cdot326\cdot2^{k+s-2}+ 8\cdot147\cdot2^{k+s-1} -32\cdot33\cdot2^{k+s-2}+2^{k+s} = 1302\cdot2^{k+s-1} \\
& \\
  & \text{The case  $ m = 0, l = 0 $}   \\
 &   \Delta _{k}\Gamma_{s+2}^{\left[s\atop{ s\atop s} \right]} -4\cdot \Delta _{k}\Gamma_{s+1}^{\left[s-1\atop{ s\atop s} \right]}
-8\cdot\Delta _{k}\Gamma_{s+1}^{\left[s-1\atop{ s\atop s} \right]}+32\cdot\Delta _{k}\Gamma_{s}^{\left[s-1\atop{ s-1\atop s} \right]}= \Delta _{k} \omega _{s+2}(s,0,0) = 2^{k+s}\\
& \text{since} \quad  \Delta _{k}\Gamma_{s+1}^{\left[s-1\atop{ s\atop s} \right]}= 630\cdot2^{k+s-2}, \quad
 \text{and} \; \Delta _{k}\Gamma_{s}^{\left[s-1\atop{ s-1\atop s} \right]}=71\cdot2^{k+s-2}\\
& \quad  \text{we get}\quad  \Delta _{k}\Gamma_{s+2}^{\left[s\atop{ s\atop s} \right]}
= 12\cdot630\cdot2^{k+s-2} -32\cdot71\cdot2^{k+s-2}+2^{k+s} = 2646\cdot2^{k+s-1} 
 \end{align*}
 \end{proof} 
  \subsection{\textbf{Computation of $ \Delta _{k}\Gamma_{s+j}^{\left[s\atop{ s\atop s} \right]} - 2^{9}\cdot \Delta _{k}\Gamma_{s+j-3}^{\left[s-1\atop{ s-1\atop s-1}\right]}, \quad for\; 0\leq j \leq s-1 $}}
\label{subsec 4} 
  \begin{lem} We have \\
  \label{lem 9.3}
  \begin{equation}
\label{eq 9.4}
 \Delta _{k}\Gamma_{s+j}^{\left[s\atop{ s\atop s} \right]} - 2^{9}\cdot \Delta _{k}\Gamma_{s+j-3}^{\left[s-1\atop{ s-1\atop s-1} \right]}=
 \begin{cases}
7\cdot2^{k+s-1} &\text{if  }  j = 0 \\
147\cdot2^{k+s-1}& \text{if   } j = 1 \\
427\cdot2^{k+s} & \text{if   }  j = 2 \\
441\cdot2^{k+s+3j -6} & \text{if   }  3\leq j\leq s-1 
\end{cases}
\end{equation}

\end{lem}
 \begin{proof}
We use Lemma \ref{lem 7.4},  Lemma \ref{lem 9.2},  \eqref{eq 6.33} and \eqref{eq 8.9}\\

\underline{The case j = 0}\\
\begin{equation*}
\Delta _{k}\Gamma_{s}^{\left[s\atop{ s\atop s} \right]} - 2^{9}\cdot \Delta _{k}\Gamma_{s-3}^{\left[s-1\atop{ s-1\atop s-1} \right]} =
\Delta _{k}\Gamma_{s}^{\left[s\atop{ s\atop s} \right]} = 7\cdot2^{k+s-1}
\end{equation*}
\underline{The case j = 1}\\
\begin{equation*}
\Delta _{k}\Gamma_{s +1}^{\left[s\atop{ s\atop s} \right]} - 2^{9}\cdot \Delta _{k}\Gamma_{s-2}^{\left[s-1\atop{ s-1\atop s-1} \right]} =
\Delta _{k}\Gamma_{s +1}^{\left[s\atop{ s\atop s} \right]} = 147\cdot2^{k+s-1}
\end{equation*}
\underline{The case j = 2}\\
\begin{align*}
& \Delta _{k}\Gamma_{s +2}^{\left[s\atop{ s\atop s} \right]}
 - 8\cdot \Delta _{k}\Gamma_{s +1}^{\left[s-1\atop{ s\atop s} \right]} = 63\cdot2^{k+s}
 \Longleftrightarrow 4\cdot\Delta _{k}\Gamma_{s +1}^{\left[s-1\atop{ s\atop s} \right]}
 - 32\cdot \Delta _{k}\Gamma_{s }^{\left[s-1\atop{ s-1\atop s} \right]} +  \Delta _{k}\omega_{s+2} (s,0,0) = 63\cdot2^{k+s}\\
 & \Longleftrightarrow \Delta _{k}\Gamma_{s +1}^{\left[s-1\atop{ s\atop s} \right]}
 - 8\cdot \Delta _{k}\Gamma_{s }^{\left[s-1\atop{ s-1\atop s} \right]} = 2^{-2}\cdot[ 63\cdot2^{k+s} - \Delta _{k}\omega_{s+2} (s,0,0)]
  = 2^{-2}\cdot[ 63\cdot2^{k+s} - 2^{k+s}] =  31\cdot2^{k+s-1}\\
  & \Longleftrightarrow 4\cdot\Delta _{k}\Gamma_{s }^{\left[s-1\atop{ s-1\atop s} \right]}
 - 32\cdot \Delta _{k}\Gamma_{s -1}^{\left[s-1\atop{ s-1\atop s-1} \right]}+  \Delta _{k}\omega_{s+1} (s-1,1,0)  = 31\cdot2^{k+s-1}\\
 & \Longleftrightarrow \Delta _{k}\Gamma_{s }^{\left[s-1\atop{ s-1\atop s} \right]}
 - 8\cdot \Delta _{k}\Gamma_{s -1}^{\left[s-1\atop{ s-1\atop s-1} \right]} = 2^{-2}\cdot[ 31\cdot2^{k+s-1} - \Delta _{k}\omega_{s+1} (s-1,1,0)]
  = 2^{-2}\cdot[ 31\cdot2^{k+s-1} - 2^{k+s-1}] =  15\cdot2^{k+s-2}
 \end{align*}
We then obtain the following equalities, \\
\begin{align*}
 \Delta _{k}\Gamma_{s +2}^{\left[s\atop{ s\atop s} \right]}
 - 8\cdot \Delta _{k}\Gamma_{s +1}^{\left[s-1\atop{ s\atop s} \right]}& = 63\cdot2^{k+s}\\
  2^{3}\cdot\big(\Delta _{k}\Gamma_{s +1}^{\left[s-1\atop{ s\atop s} \right]}
 - 8\cdot \Delta _{k}\Gamma_{s }^{\left[s-1\atop{ s -1\atop s} \right]}\big) &  =  2^{3}\cdot 31\cdot2^{k+s-1}\\
   2^{6}\cdot\big(\Delta _{k}\Gamma_{s }^{\left[s-1\atop{ s-1\atop s} \right]}
 - 8\cdot \Delta _{k}\Gamma_{s-1 }^{\left[s-1\atop{ s -1\atop s-1} \right]}\big) & =  2^{6}\cdot 15\cdot2^{k+s-2}.
\end{align*}
Summing the above equalities we get \\
\begin{equation*}
\Delta _{k}\Gamma_{s +2}^{\left[s\atop{ s\atop s} \right]}
 - 2^{9}\cdot \Delta _{k}\Gamma_{s -1}^{\left[s-1\atop{ s-1\atop s-1} \right]} = 2^{k+s}\cdot[63 +31\cdot2^{2} +15\cdot2^{4}] = 427\cdot2^{k+s}
\end{equation*}
\vspace{0.3 cm}\\
\underline{The case $3\leq j\leq s-1 $}\\

\begin{align*}
& \Delta _{k}\Gamma_{s +j}^{\left[s\atop{ s\atop s} \right]}
 - 8\cdot \Delta _{k}\Gamma_{s +j-1}^{\left[s-1\atop{ s\atop s} \right]} = 63\cdot2^{k+s +3j-6}\\
 & \Updownarrow \\
 & 4\cdot\Delta _{k}\Gamma_{s +j-1}^{\left[s-1\atop{ s\atop s} \right]}
 - 32\cdot \Delta _{k}\Gamma_{s +j-2}^{\left[s-1\atop{ s-1\atop s} \right]} +  \Delta _{k}\omega_{s+j} (s,0,0) = 63\cdot2^{k+s+3j-6}\\
 & \Updownarrow \\
 &  \Delta _{k}\Gamma_{s +j-1}^{\left[s-1\atop{ s\atop s} \right]}
 - 8\cdot \Delta _{k}\Gamma_{s +j-2}^{\left[s-1\atop{ s-1\atop s} \right]} = 2^{-2}\cdot[ 63\cdot2^{k+s +3j -6} - \Delta _{k}\omega_{s+j} (s,0,0)]
  = 2^{-2}\cdot[ 63\cdot2^{k+s +3j-6} -  0] =  63\cdot2^{k+s +3j -8}\\
  & \Updownarrow \\
   &  4\cdot\Delta _{k}\Gamma_{s +j-2}^{\left[s-1\atop{ s-1\atop s} \right]}
 - 32\cdot \Delta _{k}\Gamma_{s +j-3}^{\left[s-1\atop{ s-1\atop s-1} \right]}+  \Delta _{k}\omega_{s+j-1} (s-1,1,0)  = 63\cdot2^{k+s +3j-8}\\
 & \Updownarrow \\
  & \Delta _{k}\Gamma_{s +j-2}^{\left[s-1\atop{ s-1\atop s} \right]}
 - 8\cdot \Delta _{k}\Gamma_{s -1}^{\left[s-1\atop{ s-1\atop s-1} \right]} = 2^{-2}\cdot[ 63\cdot2^{k+s +3j-8} - \Delta _{k}\omega_{s+j-1} (s-1,1,0)]
  = [ 63\cdot2^{k+s +3j-10}-0] = 63\cdot2^{k+s +3j -10}
 \end{align*}
 \vspace{0.3 cm}\\
 We then obtain the following equalities \\
\begin{align*}
 \Delta _{k}\Gamma_{s +j}^{\left[s\atop{ s\atop s} \right]}
 - 8\cdot \Delta _{k}\Gamma_{s +j-1}^{\left[s-1\atop{ s\atop s} \right]}& = 63\cdot2^{k+s +3j-6}\\
  2^{3}\cdot\big(\Delta _{k}\Gamma_{s +j-1}^{\left[s-1\atop{ s\atop s} \right]}
 - 8\cdot \Delta _{k}\Gamma_{s +j-2}^{\left[s-1\atop{ s -1\atop s} \right]}\big) & =  2^{3}\cdot 63\cdot2^{k+s +3j-8}\\
   2^{6}\cdot\big(\Delta _{k}\Gamma_{s +j-2}^{\left[s-1\atop{ s-1\atop s} \right]}
 - 8\cdot \Delta _{k}\Gamma_{s +j-3 }^{\left[s-1\atop{ s -1\atop s-1} \right]}\big) & =  2^{6}\cdot 63\cdot2^{k+s +3j -10}
\end{align*}
\vspace{0.3 cm}\\
Summing the above equalities we get \\
\begin{equation*}
\Delta _{k}\Gamma_{s + j}^{\left[s\atop{ s\atop s} \right]}
 - 2^{9}\cdot \Delta _{k}\Gamma_{s -1}^{\left[s-1\atop{ s-1\atop s-1} \right]} = 63\cdot2^{k+s +3j-6}\cdot[1+2+2^{2}] = 63\cdot 7 \cdot2^{k+s+3j-6}=441\cdot2^{k+s+3j-6}
\end{equation*}
\end{proof}
\subsection{\textbf{Computation of $ \Delta _{k}\Gamma_{s+j}^{\left[s\atop{ s\atop s} \right]} \quad for\; 0\leq j \leq s-1 $}}
\label{subsec 5} 
  \begin{lem} We have \\
  \label{lem 9.4}
  \begin{equation}
\label{eq 9.5}
 \Delta _{k}\Gamma_{s+j}^{\left[s\atop{ s\atop s} \right]}=
 \begin{cases}
7\cdot2^{k+s-1}& \text{if    }\quad j = 0 \\
147\cdot[5\cdot2^{k+s+4j-7} - 2^{k+s+3j-6}]= 147\cdot(5\cdot2^{j-1} - 1)\cdot2^{k+s+3j-6}  & \text{if    }\quad 1\leq j\leq s-1\\
\end{cases}
\end{equation}
\end{lem}
\begin{proof}
We obtain from \eqref{eq 9.4}\\

\underline{The case j = 0}\vspace{0.1 cm}\\
\begin{equation*}
\Delta _{k}\Gamma_{s}^{\left[s\atop{ s\atop s} \right]}= 7\cdot2^{k+s-1}
\end{equation*}
\underline{The case j = 1}\vspace{0.1 cm}\\
\begin{equation*}
\Delta _{k}\Gamma_{s +1}^{\left[s\atop{ s\atop s} \right]}= 147\cdot2^{k+s-1}
\end{equation*}
\underline{The case j = 2}\vspace{0.1 cm}\\
\begin{equation*}
\Delta _{k}\Gamma_{s +2}^{\left[s\atop{ s\atop s} \right]}= 427 \cdot2^{k+s}
 +2^{9}\cdot\Delta _{k}\Gamma_{s -1}^{\left[s-1\atop{ s-1\atop s-1} \right]} =  427 \cdot2^{k+s} + 2^{9}\cdot\cdot7\cdot2^{k+s-2}= 1323\cdot2^{k+s}
\end{equation*}\vspace{0.1 cm}\\

\underline{The case $3\leq j\leq s-1,\quad j \equiv 0 \pmod{2}$}\vspace{0.1 cm}\\
\begin{align*}
\Delta _{k}\Gamma_{s + j}^{\left[s\atop{ s\atop s} \right]}
 - 2^{9}\cdot \Delta _{k}\Gamma_{s +j -3}^{\left[s-1\atop{ s-1\atop s-1} \right]} & =  441\cdot2^{k+s+3j-6} \\
 2^{9}\cdot\big(\Delta _{k}\Gamma_{s + j -3}^{\left[s-1\atop{ s-1\atop s-1} \right]}
 - 2^{9}\cdot \Delta _{k}\Gamma_{s +j - 6}^{\left[s-2\atop{ s-2\atop s-2} \right]}\big) & =  441\cdot2^{k+(s-1)+3(j-2) - 6}\cdot2^{9}\\
 (2^{9})^{2}\cdot\big(\Delta _{k}\Gamma_{s + j -6}^{\left[s-2\atop{ s-2\atop s-2} \right]}
 - 2^{9}\cdot \Delta _{k}\Gamma_{s +j - 9}^{\left[s-3\atop{ s-3\atop s-3} \right]}\big) & =  441\cdot2^{k+(s-2)+3(j-4) - 6}\cdot(2^{9})^{2}\\
 & \vdots \\
   (2^{9})^{\frac{j}{2}-2}\cdot\big(\Delta _{k}\Gamma_{s +j -3(\frac{j}{2}-2) }^{\left[s- (\frac{j}{2}-2)\atop{ s- (\frac{j}{2}-2)\atop s- (\frac{j}{2}-2)} \right]}
 - 2^{9}\cdot \Delta _{k}\Gamma_{s +j -3(\frac{j}{2}-1) }^{\left[s- (\frac{j}{2}-1)\atop{ s- (\frac{j}{2}-1)\atop s-(\frac{j}{2}-1)} \right]}\big) & =  441\cdot2^{k+(s-(\frac{j}{2}-2)+3(j- (j-4)) - 6}\cdot(2^{9})^{\frac{j}{2}-2}\\
   (2^{9})^{\frac{j}{2}-1}\cdot\big(\Delta _{k}\Gamma_{s +j -3(\frac{j}{2}-1) }^{\left[s- (\frac{j}{2}-1)\atop{ s- (\frac{j}{2}-1)\atop s- (\frac{j}{2}-1)} \right]}
 - 2^{9}\cdot \Delta _{k}\Gamma_{s +j -3(\frac{j}{2}) }^{\left[s- (\frac{j}{2})\atop{ s- (\frac{j}{2})\atop s-(\frac{j}{2})} \right]}\big) & =  427\cdot2^{k+s -\frac{j}{2} +1}\cdot(2^{9})^{\frac{j}{2}-1}\\
    (2^{9})^{\frac{j}{2}}\cdot\big(\Delta _{k}\Gamma_{s +j -3(\frac{j}{2}) }^{\left[s- (\frac{j}{2})\atop{ s- (\frac{j}{2})\atop s- (\frac{j}{2})} \right]}
 - 2^{9}\cdot \Delta _{k}\Gamma_{s +j -3(\frac{j}{2}+1) }^{\left[s- (\frac{j}{2}+1)\atop{ s- (\frac{j}{2}+1)\atop s-(\frac{j}{2}+1)} \right]}\big) & =  7\cdot2^{k+s -\frac{j}{2} -1}\cdot(2^{9})^{\frac{j}{2}}
  \end{align*}
We obtain by summing the above equations \vspace{0.1 cm}\\
\begin{align*}
&\Delta _{k}\Gamma_{s + j}^{\left[s\atop{ s\atop s} \right]} =
 \sum_{d = 0}^{\frac{j}{2} -2}441\cdot2^{k + (s-d) +3(j -2d)-6}\cdot2^{qd} 
 + 427\cdot2^{k+s-\frac{j}{2}+1}\cdot(2^{9})^{\frac{j}{2}-1}
  +  7\cdot2^{k+s-\frac{j}{2}-1}\cdot(2^{9})^{\frac{j}{2}}\\
  & = 441\cdot2^{k + s +3j -6} \sum_{d = 0}^{\frac{j}{2} -2}2^{2d}
   + 427\cdot2^{k+s +4j-8} +  7\cdot2^{k+s +4j-1}\\
   & =  441\cdot2^{k + s +3j -6}\cdot\frac{2^{j-2}-1}{2^{2}-1} + 427\cdot2^{k+s +4j-8} +  7\cdot2^{k+s +4j-1}\\
   & = 147\cdot2^{k+ s +4j-8} - 147\cdot2^{k+ s +3j-6}  + 427\cdot2^{k+s +4j-8} +  7\cdot2^{k+s +4j-1}\\
   & =  - 147\cdot2^{k+ s +3j-6}  + 735\cdot2^{k+s +4j-7}\\
   & = 147\cdot[5\cdot2^{k+s+4j-7} - 2^{k+s+3j-6}]
\end{align*}
\underline{The case $3\leq j\leq s-1,\quad j \equiv 1 \pmod{2}$}\vspace{0.1 cm}\\
 Similarly to the proof in the case $ j \equiv 0 \pmod{2} $    \\
 \end{proof}
  \subsection{\textbf{Computation of $ \Delta _{k}\Gamma_{2s+j}^{\left[s\atop{ s\atop s} \right]} - 2^{9}\cdot \Delta _{k}\Gamma_{2s+j-3}^{\left[s-1\atop{ s-1\atop s-1}\right]}, \quad for\; 0\leq j \leq s $}}
\label{subsec 6} 
  \begin{lem} We have \\
  \label{lem 9.5}
  \begin{equation}
\label{eq 9.6}
 \Delta _{k}\Gamma_{2s+j}^{\left[s\atop{ s\atop s} \right]} - 2^{9}\cdot \Delta _{k}\Gamma_{2s+j-3}^{\left[s-1\atop{ s-1\atop s-1} \right]}=
 \begin{cases}
21\cdot2^{2k+2s-2}-231\cdot2^{k+4s-6} &\text{if  }  j = 0 \\
-2373\cdot2^{2k+2s-2}+ 21\cdot2^{k+4s-3} & \text{if   } j = 1 \\
-7\cdot315\cdot2^{2k+2s+4j-6}& \text{if   } 2\leq j\leq s-1 \\
-441\cdot2^{3k+3s-3}+147\cdot2^{2k+6s-6} & \text{if   }  j=s 
\end{cases}
\end{equation}

\end{lem}

  \begin{proof}
We use Lemma \ref{lem 7.4},  \eqref{eq 8.10} and \eqref{eq 8.3} \\
\underline{The case j = 0}
\begin{align*}
& \Delta _{k}\Gamma_{2s}^{\left[s\atop{ s\atop s} \right]}
 - 8\cdot \Delta _{k}\Gamma_{2s -1}^{\left[s-1\atop{ s\atop s} \right]} = -15\cdot2^{2k+2s-2} -  33\cdot2^{k+4s -6}\\
& \Updownarrow  \\
&   4\cdot\Delta _{k}\Gamma_{2s -1}^{\left[s-1\atop{ s\atop s} \right]}
 - 32\cdot \Delta _{k}\Gamma_{2s -2}^{\left[s-1\atop{ s-1\atop s} \right]} +  \Delta _{k}\omega_{2s} (s,0,0) =  -15\cdot2^{2k+2s-2} -  33\cdot2^{k+4s -6} \\
 & \Updownarrow  \\
&  \Delta _{k}\Gamma_{2s -1}^{\left[s-1\atop{ s\atop s} \right]}
 - 8\cdot \Delta _{k}\Gamma_{2s -2}^{\left[s-1\atop{ s-1\atop s} \right]} = 2^{-2}\cdot[ -15\cdot2^{2k+2s-2} -  33\cdot2^{k+4s -6} - \Delta _{k}\omega_{2s} (s,0,0)]\\
 & = 2^{-2}\cdot[  -15\cdot2^{2k+2s-2} -  33\cdot2^{k+4s -6} + 9\cdot2^{2k+2s-2}] =   -3\cdot2^{2k+2s-3} -  33\cdot2^{k+4s - 8}                                              \\
  & \Updownarrow \\
  &   4\cdot\Delta _{k}\Gamma_{2s -2}^{\left[s-1\atop{ s-1\atop s} \right]}
 - 32\cdot \Delta _{k}\Gamma_{2s -3}^{\left[s-1\atop{ s-1\atop s-1} \right]}+  \Delta _{k}\omega_{2s-1} (s-1,1,0)  =  -3\cdot2^{2k+2s-3} -  33\cdot2^{k+4s - 8}  \\
 & \Updownarrow \\
 &  \Delta _{k}\Gamma_{2s -2}^{\left[s-1\atop{ s-1\atop s} \right]}
 - 8\cdot \Delta _{k}\Gamma_{2s-3}^{\left[s-1\atop{ s-1\atop s-1} \right]} = 2^{-2}\cdot[ -3\cdot2^{2k+2s-3} -  33\cdot2^{k+4s - 8}- \Delta _{k}\omega_{2s-1} (s-1,1,0)] \\
 & = 2^{-2}\cdot[ -3\cdot2^{2k+2s-3} -  33\cdot2^{k+4s - 8} +  9\cdot2^{2k+2s-3}] =  -3\cdot2^{2k+2s-4} -  33\cdot2^{k+4s - 10}
 \end{align*}
 
We then obtain the following equalities \\
\begin{align*}
& \Delta _{k}\Gamma_{2s}^{\left[s\atop{ s\atop s} \right]}
 - 8\cdot \Delta _{k}\Gamma_{2s -1}^{\left[s-1\atop{ s\atop s} \right]} = -15\cdot2^{2k+2s-2} -  33\cdot2^{k+4s -6} \\
 & 2^{3}\cdot\big(\Delta _{k}\Gamma_{2s -1}^{\left[s-1\atop{ s\atop s} \right]}
 - 8\cdot \Delta _{k}\Gamma_{2s -2}^{\left[s-1\atop{ s -1\atop s} \right]}\big) =  2^{3}\cdot [ -3\cdot2^{2k+2s-3} -  33\cdot2^{k+4s - 8} ] \\
  & 2^{6}\cdot\big(\Delta _{k}\Gamma_{2s -2}^{\left[s-1\atop{ s-1\atop s} \right]}
 - 8\cdot \Delta _{k}\Gamma_{2s -3 }^{\left[s-1\atop{ s -1\atop s-1} \right]}\big) =  2^{6}\cdot [ 3\cdot2^{2k+2s- 4} -  33\cdot2^{k+4s - 10} ] 
\end{align*}
Summing the above equalities we get \\
\begin{equation*}
\Delta _{k}\Gamma_{2s }^{\left[s\atop{ s\atop s} \right]}
 - 2^{9}\cdot \Delta _{k}\Gamma_{2s -3}^{\left[s-1\atop{ s-1\atop s-1} \right]} = 2^{2k+2s -2}\cdot[-15 -12 +48]  +  2^{k+4s -6}\cdot[ -33 -33\cdot2 -33\cdot2^{2}]  
 =  21\cdot2^{2k+2s-2} -  231\cdot2^{k+4s - 6}
\end{equation*}
\vspace{0.2 cm}\\
\underline{The case j = 1}

\begin{align*}
& \Delta _{k}\Gamma_{2s+1}^{\left[s\atop{ s\atop s} \right]}
 - 8\cdot \Delta _{k}\Gamma_{2s }^{\left[s-1\atop{ s\atop s} \right]} = -321\cdot2^{2k+2s-2}  +3\cdot2^{k+4s - 3}\\
& \Updownarrow \\
 & 4\cdot\Delta _{k}\Gamma_{2s}^{\left[s-1\atop{ s\atop s} \right]}
 - 32\cdot \Delta _{k}\Gamma_{2s -1 }^{\left[s-1\atop{ s-1\atop s} \right]} +  \Delta _{k}\omega_{2s +1} (s,0,0) =   -321\cdot2^{2k+2s-2}  +3\cdot2^{k+4s - 3} \\
 & \Updownarrow \\
&  \Delta _{k}\Gamma_{2s }^{\left[s-1\atop{ s\atop s} \right]}
 - 8\cdot \Delta _{k}\Gamma_{2s -1}^{\left[s-1\atop{ s-1\atop s} \right]} = 2^{-2}\cdot[ -321\cdot2^{2k+2s-2}  +3\cdot2^{k+4s - 3} - \Delta _{k}\omega_{2s +1} (s,0,0)] \\
 & = 2^{-2}\cdot[ -321\cdot2{2k+2s-2}  +3\cdot2^{k+4s - 3}- 9\cdot2^{2k+2s-2}] =   -165\cdot2^{2k+2s-3}  + 3\cdot2^{k+4s - 5}                                              \\
  & \Updownarrow \\
  &  4\cdot\Delta _{k}\Gamma_{2s -1}^{\left[s-1\atop{ s-1\atop s} \right]}
 - 32\cdot \Delta _{k}\Gamma_{2s -2}^{\left[s-1\atop{ s-1\atop s-1} \right]}+  \Delta _{k}\omega_{2s} (s-1,1,0)  =    -165\cdot2^{2k+2s-3}  + 3\cdot2^{k+4s - 5}  \\
 & \Updownarrow \\
&  \Delta _{k}\Gamma_{2s -1}^{\left[s-1\atop{ s-1\atop s} \right]}
 - 8\cdot \Delta _{k}\Gamma_{2s-2}^{\left[s-1\atop{ s-1\atop s-1} \right]} = 2^{-2}\cdot[  -165\cdot2^{2k+2s-3}  + 3\cdot2^{k+4s - 5}  - \Delta _{k}\omega_{2s} (s-1,1,0)] \\
 & = 2^{-2}\cdot[  -165\cdot2^{2k+2s-3} + 3\cdot2^{k+4s - 5}  -  9\cdot2^{2k+2s-3}] =   -87\cdot2^{2k+2s-4} + 3\cdot2^{k+4s - 7}
 \end{align*}
We then obtain the following equalities \vspace{0.1 cm}\\
\begin{align*}
 \Delta _{k}\Gamma_{2s+1}^{\left[s\atop{ s\atop s} \right]}
 - 8\cdot \Delta _{k}\Gamma_{2s }^{\left[s-1\atop{ s\atop s} \right]}& = - 321\cdot2^{2k+2s-2}  +3\cdot2^{k+4s -3} \\
  2^{3}\cdot\big(\Delta _{k}\Gamma_{2s }^{\left[s-1\atop{ s\atop s} \right]}
 - 8\cdot \Delta _{k}\Gamma_{2s -1}^{\left[s-1\atop{ s -1\atop s} \right]}\big) &  =  2^{3}\cdot [ -165\cdot2^{2k+2s-3} +3 \cdot2^{k+4s - 5} ] \\
   2^{6}\cdot\big(\Delta _{k}\Gamma_{2s -1}^{\left[s-1\atop{ s-1\atop s} \right]}
 - 8\cdot \Delta _{k}\Gamma_{2s -2 }^{\left[s-1\atop{ s -1\atop s-1} \right]}\big) & =  2^{6}\cdot [ -87\cdot2^{2k+2s- 4}  +3\cdot2^{k+4s - 7} ] 
\end{align*}
Summing the above equalities we get \\
\begin{align*}
& \Delta _{k}\Gamma_{2s+1 }^{\left[s\atop{ s\atop s} \right]}
 - 2^{9}\cdot \Delta _{k}\Gamma_{2s -2}^{\left[s-1\atop{ s-1\atop s-1} \right]} = 2^{2k+2s -2}\cdot[-321 +(-165)\cdot2^{2} +(-87)\cdot2^{4}]  +  2^{k+4s -3}\cdot[ 3 + 3\cdot2 + 3\cdot2^{2}] \\
&   =  -2373\cdot2^{2k+2s-2}  +21\cdot2^{k+4s - 3}
\end{align*} 
\underline{The case $2\leq j\leq s-1$}

\begin{align*}
& \Delta _{k}\Gamma_{2s+j}^{\left[s\atop{ s\atop s} \right]}
 - 8\cdot \Delta _{k}\Gamma_{2s +j-1}^{\left[s-1\atop{ s\atop s} \right]} = -315\cdot2^{2k+2s +4j-6}  \\
& \Updownarrow  \\
&  4\cdot\Delta _{k}\Gamma_{2s +j-1}^{\left[s-1\atop{ s\atop s} \right]}
 - 32\cdot \Delta _{k}\Gamma_{2s + j-2 }^{\left[s-1\atop{ s-1\atop s} \right]} +  \Delta _{k}\omega_{2s +j} (s,0,0) =  -315\cdot2^{2k+2s +4j-6}    \\
 & \Updownarrow  \\
 &  \Delta _{k}\Gamma_{2s +j-1}^{\left[s-1\atop{ s\atop s} \right]}
 - 8\cdot \Delta _{k}\Gamma_{2s +j-2}^{\left[s-1\atop{ s-1\atop s} \right]} = 2^{-2}\cdot[ -315\cdot2^{2k+2s +4j-6}  - \Delta _{k}\omega_{2s +j} (s,0,0)] \\
 & = 2^{-2}\cdot[ -315\cdot2{2k+2s +4j -6}  - 0] =   -315\cdot2^{2k+2s +4j-8}                                             \\
  & \Updownarrow \\
  & 4\cdot\Delta _{k}\Gamma_{2s j-2}^{\left[s-1\atop{ s-1\atop s} \right]}
 - 32\cdot \Delta _{k}\Gamma_{2s +j-3}^{\left[s-1\atop{ s-1\atop s-1} \right]}+  \Delta _{k}\omega_{2s+j-1} (s-1,1,0)  =    -315\cdot2^{2k+2s +4j -8}    \\
 & \Updownarrow  \\
 &  \Delta _{k}\Gamma_{2s +j-2}^{\left[s-1\atop{ s-1\atop s} \right]}
 - 8\cdot \Delta _{k}\Gamma_{2s +j-3}^{\left[s-1\atop{ s-1\atop s-1} \right]} = 2^{-2}\cdot[  -315\cdot2^{2k+2s +4j -8}    - \Delta _{k}\omega_{2s+j-1} (s-1,1,0)] \\
 & = 2^{-2}\cdot[  -315\cdot2^{2k+2s +4j-8}  - 0] =   -315\cdot2^{2k+2s +4j-10} 
 \end{align*}
 \vspace{0.2 cm}\\
 We then obtain the following equalities \\
\begin{align*}
 \Delta _{k}\Gamma_{2s+j}^{\left[s\atop{ s\atop s} \right]}
 - 8\cdot \Delta _{k}\Gamma_{2s +j-1}^{\left[s-1\atop{ s\atop s} \right]} &  = - 315\cdot2^{2k+2s +4j-6} \\
  2^{3}\cdot\big(\Delta _{k}\Gamma_{2s+j-1 }^{\left[s-1\atop{ s\atop s} \right]}
 - 8\cdot \Delta _{k}\Gamma_{2s +j-2}^{\left[s-1\atop{ s -1\atop s} \right]}\big)  & =  2^{3}\cdot [ -315\cdot2^{2k+2s +4j-8} ] \\
   2^{6}\cdot\big(\Delta _{k}\Gamma_{2s  +j-2}^{\left[s-1\atop{ s-1\atop s} \right]}
 - 8\cdot \Delta _{k}\Gamma_{2s + j-3 }^{\left[s-1\atop{ s -1\atop s-1} \right]}\big) &  =  2^{6}\cdot [ -315\cdot2^{2k+2s +4j-10}  ] 
\end{align*}
 \vspace{0.2 cm}\\
Summing the above equalities we get \\
\begin{align*}
& \Delta _{k}\Gamma_{2s+j }^{\left[s\atop{ s\atop s} \right]}
 - 2^{9}\cdot \Delta _{k}\Gamma_{2s +j -3}^{\left[s-1\atop{ s-1\atop s-1} \right]} =-315\cdot 2^{2k+2s +4j -6}\cdot[1 + 2 + 2^{2}]  
  =  7\cdot(-315)\cdot2^{2k+2s +4j-6}  
\end{align*} 
\vspace{0.2 cm}\\
\underline{The case j = s} 
\begin{align*}
& \Delta _{k}\Gamma_{3s}^{\left[s\atop{ s\atop s} \right]}
 - 8\cdot \Delta _{k}\Gamma_{3s -1}^{\left[s-1\atop{ s\atop s} \right]} = -21\cdot2^{3k+3s-3}  +21\cdot2^{2k+ 6s-6}  \\
& \Updownarrow \\
&  4\cdot\Delta _{k}\Gamma_{3s -1}^{\left[s-1\atop{ s\atop s} \right]}
 - 32\cdot \Delta _{k}\Gamma_{3s-2}^{\left[s-1\atop{ s-1\atop s} \right]} +  \Delta _{k}\omega_{3s} (s,0,0) =   -21\cdot2^{3k+3s-3}  +21\cdot2^{2k+ 6s-6}  \\
 & \Updownarrow  \\
 & \Delta _{k}\Gamma_{3s-1}^{\left[s-1\atop{ s\atop s} \right]}
 - 8\cdot \Delta _{k}\Gamma_{3s -2}^{\left[s-1\atop{ s-1\atop s} \right]} = 2^{-2}\cdot[  -21\cdot2^{3k+3s-3}  +21\cdot2^{2k+ 6s-6}   - \Delta _{k}\omega_{3s } (s,0,0)] \\
 & =   2^{-2}\cdot[  -21\cdot2^{3k+3s-3}  +21\cdot2^{2k+ 6s-6}  -21\cdot2^{3k+3s-3}] =    -21\cdot2^{3k+3s-4}  +21\cdot2^{2k+ 6s-8}    \\
  & \Updownarrow \\
  & 4\cdot\Delta _{k}\Gamma_{3s-2}^{\left[s-1\atop{ s-1\atop s} \right]}
 - 32\cdot \Delta _{k}\Gamma_{3s -3}^{\left[s-1\atop{ s-1\atop s-1} \right]}+  \Delta _{k}\omega_{3s-1} (s-1,1,0)  =    -21\cdot2^{3k+3s-4}  +21\cdot2^{2k+ 6s-8}      \\
 & \Updownarrow  \\
 & \Delta _{k}\Gamma_{3s -2}^{\left[s-1\atop{ s-1\atop s} \right]}
 - 8\cdot \Delta _{k}\Gamma_{3s -3}^{\left[s-1\atop{ s-1\atop s-1} \right]} = 2^{-2}\cdot[  -21\cdot2^{3k+3s-4}  +21\cdot2^{2k+ 6s-8}   - \Delta _{k}\omega_{3s -1} (s-1,1,0)] \\
 & = 2^{-2}\cdot[ -21\cdot2^{3k+3s-4}  +21\cdot2^{2k+ 6s-8}   -21\cdot2^{3k+3s-4} ]  =    -21\cdot2^{3k+3s-5}  +21\cdot2^{2k+ 6s- 10}  
 \end{align*}
 \vspace{0.2 cm}\\
 We then obtain the following equalities \vspace{0.1 cm}\\
\begin{align*}
 \Delta _{k}\Gamma_{3s}^{\left[s\atop{ s\atop s} \right]}
 - 8\cdot \Delta _{k}\Gamma_{3s-1}^{\left[s-1\atop{ s\atop s} \right]} & =   -21\cdot2^{3k+3s-3}  +21\cdot2^{2k+ 6s- 6}  \\
  2^{3}\cdot\big(\Delta _{k}\Gamma_{3s-1 }^{\left[s-1\atop{ s\atop s} \right]}
 - 8\cdot \Delta _{k}\Gamma_{3s-2}^{\left[s-1\atop{ s -1\atop s} \right]}\big) & =  2^{3}\cdot [   -21\cdot2^{3k+3s-4}  +21\cdot2^{2k+ 6s- 8}   ] \\
   2^{6}\cdot\big(\Delta _{k}\Gamma_{3s  -2}^{\left[s-1\atop{ s-1\atop s} \right]}
 - 8\cdot \Delta _{k}\Gamma_{3s -3 }^{\left[s-1\atop{ s -1\atop s-1} \right]}\big) & =  2^{6}\cdot [   -21\cdot2^{3k+3s-5}  +21\cdot2^{2k+ 6s- 10}   ] 
\end{align*}
 \vspace{0.2 cm}\\
Summing the above equalities we get \vspace{0.1 cm}\\
\begin{align*}
& \Delta _{k}\Gamma_{3s }^{\left[s\atop{ s\atop s} \right]}
 - 2^{9}\cdot \Delta _{k}\Gamma_{3(s-1)}^{\left[s-1\atop{ s-1\atop s-1} \right]} =  -21\cdot2^{3k+3s-3}\cdot[1+2^{2} +2^{4}]   
 + 21\cdot2^{2k+ 6s- 6}\cdot[1+2 +2^{2}]\\
 &  = - 441\cdot2^{3k+3s-3} + 7\cdot21\cdot2^{2k+ 6s- 6} = -(21)^{2}\cdot2^{3k+3s-3} + 7\cdot21\cdot2^{2k+ 6s- 6}
\end{align*}
\end{proof}
\subsection{\textbf{Computation of $ \Delta _{k}\Gamma_{2s+j}^{\left[s\atop{ s\atop s} \right]} \quad for\; 0\leq j \leq s $}}
\label{subsec 7} 
 \begin{lem} We have \\
  \label{lem 9.6}
  \begin{equation}
\label{eq 9.7}
 \Delta _{k}\Gamma_{2s+j}^{\left[s\atop{ s\atop s} \right]}=
 \begin{cases}
21\cdot[2^{2k+2s-2} + 35\cdot2^{k+5s-7} - 39\cdot2^{k+4s-6}] & \text{if   } j = 0 \\
105\cdot[3\cdot2^{2k+2s +4j-6} + 7\cdot2^{k+ 5s +4j -7} - 31\cdot2^{k+4s +5j -8}] & \text{if   }  1\leq j\leq s-1\\
7\cdot2^{3k+3s-3} -21\cdot2^{2k+6s-6} + 7\cdot2^{k+9s-8} & \text{if  } j = s
\end{cases}
\end{equation}
\end{lem}
\begin{proof}
\underline{The case j = 0}
\vspace{0.1 cm}

We have from \eqref{eq 9.5} with $s\rightarrow s-1,\;j = (s-1) -1$\\
\begin{align*}
\Delta _{k}\Gamma_{2s -3}^{\left[s-1\atop{ s-1\atop s-1} \right]}=\Delta _{k}\Gamma_{(s-1)+ (s-2)}^{\left[s-1\atop{ s-1\atop s-1} \right]}
= 147\cdot(5\cdot2^{s-3} -1)\cdot2^{k+4s-13} = 147\cdot5\cdot2^{k+5s- 16} - 147\cdot2^{k+4s-13}\\
&
\end{align*}

We then deduce from \eqref{eq 9.6} with j = 0 \vspace{0.1 cm}\\
\begin{align}
& \Delta _{k}\Gamma_{2s }^{\left[s\atop{ s\atop s} \right]}
=  2^{9}\cdot \Delta _{k}\Gamma_{2s -3}^{\left[s-1\atop{ s-1\atop s-1} \right]} +  21\cdot2^{2k+2s-2} -  231\cdot2^{k+4s - 6} \label{eq 9.8}\\
& =   147\cdot5\cdot2^{k+5s- 7} - 147\cdot2^{k+4s-4} +  21\cdot2^{2k+2s-2} -  231\cdot2^{k+4s - 6} \nonumber \\
& = 21\cdot[2^{2k+2s-2} +35\cdot2^{k+5s-7} -39\cdot2^{k+4s-6}] \nonumber
\end{align}
\underline{The case j = 1}
\vspace{0.1 cm}

From \eqref{eq 9.8} with $s\rightarrow s-1$ we obtain \\

\begin{align*}
\Delta _{k}\Gamma_{2s -2}^{\left[s-1\atop{ s-1\atop s-1} \right]}=\Delta _{k}\Gamma_{2(s-1)}^{\left[s-1\atop{ s-1\atop s-1} \right]}
=  21\cdot[2^{2k+2s-4} +35\cdot2^{k+5s- 12} -39\cdot2^{k+4s-10}]\\
& 
\end{align*}
We then deduce from \eqref{eq 9.6} with j = 1 \vspace{0.1 cm}\\
\begin{align}
& \Delta _{k}\Gamma_{2s +1}^{\left[s\atop{ s\atop s} \right]}
=  2^{9}\cdot \Delta _{k}\Gamma_{2s -2}^{\left[s-1\atop{ s-1\atop s-1} \right]}  -2373\cdot2^{2k+2s-2}  +21\cdot2^{k+4s - 3} \label{eq 9.9}\\
& =  21\cdot[2^{2k+2s +5} +35\cdot2^{k+5s- 3} -39\cdot2^{k+4s-1}]  -21\cdot113\cdot2^{2k+2s-2}  +21\cdot2^{k+4s - 3}\nonumber  \\
& = 105\cdot[3\cdot2^{2k+2s-2} + 7\cdot2^{k+5s- 3} -31\cdot2^{k+4s- 3}]\nonumber
\end{align}
\underline{The case j = 2}
\vspace{0.1 cm}

From \eqref{eq 9.9} with $s\rightarrow s-1$ we obtain \\
\begin{align*}
 & \Delta _{k}\Gamma_{2(s-1)+1}^{\left[s-1\atop{ s-1\atop s-1} \right]}
 = 105\cdot[3\cdot2^{2k+2(s-1)-2} + 7\cdot2^{k+5(s-1)- 3} -31\cdot2^{k+4(s-1)- 3}] \\
& =  105\cdot[3\cdot2^{2k+2s-4} + 7\cdot2^{k+5s-8} -31\cdot2^{k+4s-7}] 
\end{align*}
We then deduce from \eqref{eq 9.6} with j = 2 \vspace{0.1 cm}\\
\begin{align}
& \Delta _{k}\Gamma_{2s+2 }^{\left[s\atop{ s\atop s} \right]} =
   2^{9}\cdot \Delta _{k}\Gamma_{2(s-1)+1}^{\left[s-1\atop{ s-1\atop s-1} \right]} +  7\cdot(-315)\cdot2^{2k+2s +2}\label{eq 9.10} \\
   & =  2^{9}\cdot105\cdot[3\cdot2^{2k+2s-4} + 7\cdot2^{k+5s-8} -31\cdot2^{k+4s-7}] + 105\cdot(-21)\cdot2^{2k+2s +2} \nonumber \\
   & = 105\cdot[3\cdot2^{2k+2s +2} + 7\cdot2^{k+5s+1} -31\cdot2^{k+4s+2}] \nonumber 
\end{align} \vspace{0.1 cm}\\

We observe that \\
 \begin{align*}
&16^{3}\cdot \Delta _{k-3}\Gamma_{2(s-1)+1 }^{\left[s -1\atop{ s-1\atop s-1} \right]} =
   2^{12}\cdot 105\cdot[3\cdot2^{2(k-3)+2(s-1)  -2} + 7\cdot2^{k-3+5(s-1) -3} -31\cdot2^{k-3+4(s-1) -3}]  \\
 & =   105\cdot[3\cdot2^{2k+2s+2} + 7\cdot2^{k +5s +1} -31\cdot2^{k +4s +2}]
\end{align*} \vspace{0.1 cm}\\
From \eqref{eq 9.10} and the the above equation we deduce \\
\begin{equation}
\label{eq 9.11}
\Delta _{k}\Gamma_{2s+2 }^{\left[s\atop{ s\atop s} \right]}= 16^{3}\cdot \Delta _{k-3}\Gamma_{2(s-1)+1 }^{\left[s -1\atop{ s-1\atop s-1} \right]}
\end{equation}
\underline{The case $2\leq j\leq s-1$}
\vspace{0.1 cm}

We have for $ 1\leq j\leq s -2 $ \vspace{0.1 cm}\\
\begin{equation}
\label{eq 9.12}
 (H_{j}) \quad 
\Delta _{k}\Gamma_{2s+1 +j }^{\left[s\atop{ s\atop s} \right]}= 16^{3j}\cdot \Delta _{k-3j}\Gamma_{2(s-j)+1 }^{\left[s -j\atop{ s-j\atop s-j} \right]}
\end{equation}
 We show \eqref{eq 9.12} by induction on j.\\
By \eqref{eq 9.11} $ (H_{j}) $ holds for j = 1. \\
We shall show that $ (H_{j}) $ holds implies that $ (H_{j+1 }) $ holds. \\
We have \\
\begin{align*}
& \Delta _{k}\Gamma_{2s+ 1+ (j+1) }^{\left[s\atop{ s\atop s} \right]}
 = 2^{9}\cdot \Delta _{k}\Gamma_{2s + 1 + (j+1) -3}^{\left[s-1\atop{ s-1\atop s-1} \right]} + 7\cdot(-315)\cdot2^{2k+2s +4(j+2)-6}    \\
& =   2^{9}\cdot \Delta _{k}\Gamma_{2(s-1) + 1 +j}^{\left[s-1\atop{ s-1\atop s-1} \right]} + 7\cdot(-315)\cdot2^{2k+2s +4j +2}  \\
& =    2^{9}\cdot2^{12j}\cdot \Delta _{k-3j}\Gamma_{2(s-(j+1)) + 1 }^{\left[s-1-j\atop{ s-1-j\atop s-1-j} \right]} + 7\cdot(-315)\cdot2^{2k+2s +4j +2}  \\
& = 2^{12j+9}\cdot105\cdot\big[3\cdot2^{2(k-3j)+2(s-1-j) - 2} + 7\cdot2^{k-3j +5(s-1-j)  -3} -31\cdot2^{k-3j +4(s-1-j) -3}  \big] + 21\cdot(- 105)\cdot2^{2k+2s +4j +2}\\
& = 105\cdot\big[3\cdot2^{2k+2s+4j+2} + 7\cdot2^{k +5s +4j +1} -31\cdot2^{k +4s+5j+2}  \big]
 \end{align*} \vspace{0.01 cm}\\
 On the other hand \\
 
  \begin{align*}
& 16^{3(j+1)}\cdot\Delta _{k-3(j+1)}\Gamma_{2(s-(j+1)) + 1 }^{\left[s- (j+1)\atop{ s- (j+1)\atop s- (j+1)} \right]}\\
& = 2^{12j+12}\cdot105\cdot\big[3\cdot2^{2(k-3(j+1))+2(s-j-1) - 2} + 7\cdot2^{k-3(j+1) +5(s -j-1)  -3} -31\cdot2^{k-3(j+1) +4(s- j-1) -3}  \big] \\
& = 105\cdot\big[3\cdot2^{2k+2s+4j+2} + 7\cdot2^{k+5s+4j+1} -31\cdot2^{k +4s +5j+2}  \big] 
\end{align*}\vspace{0.01 cm}\\

We have now proved that for $0\leq j\leq s-2$ we have \\
\begin{equation*}
\Delta _{k}\Gamma_{2s+1 +j }^{\left[s\atop{ s\atop s} \right]}=  105\cdot\big[3\cdot2^{2k+2s+4j-2} + 7\cdot2^{k+5s+4j -3} -31\cdot2^{k +4s +5j -3}  \big] 
\end{equation*}
\underline{The case  j = s}
\vspace{0.1 cm}

We have successively \\
\begin{align*}
&\Delta _{k}\Gamma_{3s}^{\left[s\atop{ s\atop s} \right]}
 - 2^{9}\cdot \Delta _{k}\Gamma_{3(s-1)}^{\left[s-1\atop{ s-1\atop s-1} \right]}  = - 441\cdot2^{3k+3s-3}+ 147\cdot2^{2k+6s-6} \\
 &2^{9}\cdot\big(\Delta _{k}\Gamma_{3(s-1)}^{\left[s-1\atop{ s-1\atop s-1} \right]}
 - 2^{9}\cdot \Delta _{k}\Gamma_{3(s-2)}^{\left[s-2\atop{ s-2\atop s-2} \right]}\big)  =\big(  - 441\cdot2^{3k+3(s-1)-3}+ 147\cdot2^{2k+6(s-1)-6}\big)\cdot2^{9}\\
 &(2^{9})^{2}\cdot\big(\Delta _{k}\Gamma_{3(s-2)}^{\left[s-2\atop{ s-2\atop s-2} \right]}
 - 2^{9}\cdot \Delta _{k}\Gamma_{3(s-3)}^{\left[s-3\atop{ s-3\atop s-3} \right]}\big)  = \big(  - 441\cdot2^{3k+3(s-2)-3}+ 147\cdot2^{2k+6(s-2)-6}\big) \cdot(2^{9})^{2}\\
 & \vdots \\
  &(2^{9})^{s-2}\cdot\big(\Delta _{k}\Gamma_{3(s-(s-2))}^{\left[s- (s-2)\atop{ s- (s-2)\atop s- (s-2)} \right]}
 - 2^{9}\cdot \Delta _{k}\Gamma_{3(s-(s-1))}^{\left[s-(s-1)\atop{ s- (s-1)\atop s- (s-1)} \right]}\big)  = \big(  - 441\cdot2^{3k+3(s-(s-2))-3}+ 147\cdot2^{2k+6(s-(s-2))-6}\big) \cdot(2^{9})^{s-2}\\
  \end{align*}
We obtain by summing the above equations \vspace{0.1 cm}\\
\begin{align*}
&\Delta _{k}\Gamma_{3s}^{\left[s\atop{ s\atop s} \right]} =
 \sum_{d = 0}^{s-2}\big(  - 441\cdot2^{3k+3(s- d)-3}+ 147\cdot2^{2k+6(s-d)-6}\big) \cdot(2^{9})^{d}
  +( 2^{9})^{s-1}\cdot \Delta _{k}\Gamma_{3}^{\left[ 1\atop{ 1 \atop 1} \right]} \\
  & = ( 2^{9})^{s-1}\cdot \Delta _{k}\Gamma_{3}^{\left[ 1\atop{ 1 \atop 1} \right]}+ \sum_{d = 0}^{s-2} - 441\cdot2^{3k+3s+6d -3}+ 147\cdot2^{2k +3d+6s-6}\\
  & = 2^{9s-9}\cdot\big(7\cdot2^{3k}  -21\cdot2^{2k} +7\cdot2^{k+1} \big)  - 441\cdot2^{3k+3s -3}\sum_{d = 0}^{s-2}2^{6d} +147\cdot2^{2k+6s-6}\sum_{d = 0}^{s-2}2^{3d}\\
  & = 7\cdot2^{3k+9s-9}  -21\cdot2^{2k+9s-9} +7\cdot2^{k+9s-8}  - 441\cdot2^{3k+3s -3}\cdot{2^{6s-6}-1\over2^{6} -1}
  + 147\cdot2^{2k+6s-6}\cdot{2^{3s-3} -1 \over 2^{3}-1}\\
  & =  7\cdot2^{3k+9s-9}  -21\cdot2^{2k+9s-9} +7\cdot2^{k+9s-8} - 7\cdot2^{3k+9s-9} + 7\cdot2^{3k+3s-3} + 21\cdot2^{2k+9s-9} -  21\cdot2^{2k+6s-9}\\
  & =  7\cdot2^{3k+3s-3}  -21\cdot2^{2k+6s - 6} + 7\cdot2^{k+9s-8}
  \end{align*}
\end{proof}
 \begin{lem} We have \\
  \label{lem 9.7}
 \begin{equation}
\label{eq 9.13}
 \Delta _{k}\Gamma_{s+j-1}^{\left[s-1\atop{ s\atop s} \right]}=
 \begin{cases}
2^{k+s- 2}  &\text{if  }  j = 0 \\
33\cdot2^{k+s-2}& \text{if   } j = 1\\
7\cdot(105\cdot2^{j-2} - 15)\cdot2^{k+s+3j-8}  & \text{if   } 2\leq j\leq s-1\\
9\cdot2^{2k+2s-3} + 21\cdot35\cdot2^{k+5s-10} - 393\cdot2^{k+4s-8} & \text{if   } j = s \\
159\cdot2^{2k+2s-3} + 735\cdot2^{k+5s- 6} - 1629\cdot2^{k+4s-5} & \text{if   } j = s +1
\end{cases}
\end{equation}

  \begin{equation}
\label{eq 9.14}
 \Delta _{k}\Gamma_{2s+j-1}^{\left[s-1\atop{ s\atop s} \right]}=
 \begin{cases}
105\cdot(3\cdot2^{2k+2s +4j-8} + 7\cdot2^{k+5s+4j-10} - 31\cdot2^{k+4s +5j-11}) & \text{if   } 2\leq j\leq s-1\\
7\cdot2^{3k+3s-4} -21\cdot2^{2k+6s- 8} + 7\cdot2^{k+9s-11} & \text{if   } j = s 
\end{cases}
\end{equation}
\begin{equation}
\label{eq 9.15}
 \Delta _{k}\Gamma_{s+j-2}^{\left[s-1\atop{ s-1\atop s} \right]}=
 \begin{cases}
3\cdot2^{k+s- 2}  &\text{if  }  j = 1 \\
71\cdot2^{k+s-2}& \text{if   } j = 2\\
21\cdot(35\cdot2^{j-2} - 8)\cdot2^{k+s+3j-11}  & \text{if   } 3\leq j\leq s-1\\
3\cdot2^{2k+2s-4} + 735\cdot2^{k+5s-13} - 45\cdot2^{k+4s-8} & \text{if   } j = s \\
81\cdot2^{2k+2s-4} + 735\cdot2^{k+5s- 9} - 51\cdot2^{k+4s-3} & \text{if   } j = s +1
\end{cases}
\end{equation}

\begin{equation}
\label{eq 9.16}
 \Delta _{k}\Gamma_{2s+j-2}^{\left[s-1\atop{ s-1\atop s} \right]}=
 \begin{cases}
315\cdot2^{2k+2s +4j-10} + 735\cdot2^{k+5s+4j-13} - 3255\cdot2^{k+4s +5j-14} & \text{if   } 2\leq j\leq s-1\\
7\cdot2^{3k+3s-5} -21\cdot2^{2k+6s- 10} + 7\cdot2^{k+9s-14} & \text{if   } j = s 
\end{cases}
\end{equation}

\end{lem}
\begin{proof}
\textbf{proof of \eqref{eq 9.13}, \eqref{eq 9.14}}\\[0.2 cm]
Combining \eqref{eq 8.9}, \eqref{eq 8.10} and \eqref{eq 9.5}, \eqref{eq 9.7} give the desired formulas.\\[0.2 cm]
\textbf{proof of \eqref{eq 9.15}, \eqref{eq 9.16}}\\[0.2 cm]

We obtain from \eqref{eq 8.3} with $m=l=0$\\
\begin{align}
 \Delta _{k}\Gamma_{s+j}^{\left[s\atop{ s \atop s } \right]}
   - 4\cdot \Delta _{k}\Gamma_{s+j -1}^{\left[s-1\atop{ s \atop s } \right]} & = 
   8\cdot\left[  \Delta _{k}\Gamma_{s+j -1}^{\left[s-1\atop{ s \atop s } \right]}
   - 4\cdot \Delta _{k}\Gamma_{s+j -2}^{\left[s-1\atop{ s -1\atop s } \right]}\right] + \Delta _{k}\omega_{s+j} (s,0,0) \nonumber\\
   \Longleftrightarrow  32\cdot \Delta _{k}\Gamma_{s+j -2}^{\left[s-1\atop{ s -1\atop s } \right]} & = 
   12\cdot \Delta _{k}\Gamma_{s+j -1}^{\left[s-1\atop{ s \atop s } \right]} - \Delta _{k}\Gamma_{s+j}^{\left[s\atop{ s \atop s } \right]} 
    + \Delta _{k}\omega_{s+j} (s,0,0) \label{eq 9.17}\\
    \Delta _{k}\Gamma_{2s+j}^{\left[s\atop{ s \atop s } \right]}
   - 4\cdot \Delta _{k}\Gamma_{2s+j -1}^{\left[s-1\atop{ s \atop s } \right]} & = 
   8\cdot\left[  \Delta _{k}\Gamma_{2s+j -1}^{\left[s-1\atop{ s \atop s } \right]}
   - 4\cdot \Delta _{k}\Gamma_{2s+j -2}^{\left[s-1\atop{ s -1\atop s } \right]}\right] + \Delta _{k}\omega_{2s+j} (s,0,0) \nonumber\\
   \Longleftrightarrow  32\cdot \Delta _{k}\Gamma_{2s+j -2}^{\left[s-1\atop{ s -1\atop s } \right]} & = 
   12\cdot \Delta _{k}\Gamma_{2s+j -1}^{\left[s-1\atop{ s \atop s } \right]} - \Delta _{k}\Gamma_{2s+j}^{\left[s\atop{ s \atop s } \right]} 
    + \Delta _{k}\omega_{2s+j} (s,0,0) \label{eq 9.18}\\
    & \nonumber
   \end{align}
\eqref{eq 9.15}, \eqref{eq 9.16} follows now from \eqref{eq 9.17}, \eqref{eq 9.18} and \eqref{eq 7.48}
\end{proof}

 \section{\textbf{An explicit formula for $ \Gamma_{i}^{\left[s\atop{ s\atop s} \right]\times k} , \quad for\; s \leq i \leq 3s,\; k\geq i $}}
\label{sec 10} 
    \subsection{Notation}
  \label{subsec 1}
  \begin{defn}
  \label{defn 10.1}
We recall that  $ \Gamma_{i}^{\left[s\atop{ s\atop s} \right]\times k} $ denote the number of rank i matrices of the form\\
 $$   \left ( \begin{array} {cccccc}
\alpha _{1} & \alpha _{2} & \alpha _{3} &  \ldots & \alpha _{k-1}  &  \alpha _{k} \\
\alpha _{2 } & \alpha _{3} & \alpha _{4}&  \ldots  &  \alpha _{k} &  \alpha _{k+1} \\
\vdots & \vdots & \vdots    &  \vdots & \vdots  &  \vdots \\
\alpha _{s-1} & \alpha _{s} & \alpha _{s +1} & \ldots  &  \alpha _{s+k-3} &  \alpha _{s+k-2}  \\
\alpha _{s} & \alpha _{s+1} & \alpha _{s +2} & \ldots  &  \alpha _{s+k-2} &  \alpha _{s+k-1}  \\
\hline \\
\beta  _{1} & \beta  _{2} & \beta  _{3} & \ldots  &  \beta_{k-1} &  \beta _{k}  \\
\beta  _{2} & \beta  _{3} & \beta  _{4} & \ldots  &  \beta_{k} &  \beta _{k+1}  \\
\vdots & \vdots & \vdots    &  \vdots & \vdots  &  \vdots \\
\beta  _{s-1} & \beta  _{s} & \beta  _{s+1} & \ldots  &  \beta_{s+k-3} &  \beta _{s+k-2}  \\
\beta  _{s} & \beta  _{s+1} & \beta  _{s+2} & \ldots  &  \beta_{s+k-2} &  \beta _{s+k-1} \\
\hline \\
\gamma  _{1} & \gamma   _{2} & \gamma  _{3} & \ldots  &  \gamma _{k-1} &  \gamma  _{k}  \\
\gamma   _{2} & \gamma  _{3} & \gamma   _{4} & \ldots  &  \gamma _{k} &  \gamma  _{k+1}  \\
\vdots & \vdots & \vdots    &  \vdots & \vdots  &  \vdots \\
\gamma  _{s -1} & \gamma  _{s} & \gamma  _{s+1} & \ldots  &  \gamma _{s+k-3} &  \gamma  _{s+k-2}  \\
\gamma  _{s} & \gamma  _{s+1} & \gamma  _{s+2} & \ldots  &  \gamma _{s+k-2} &  \gamma  _{s+k-1} 
\end{array}  \right) $$

  \end{defn}

   \subsection{Introduction}
  \label{subsec 2}
We adapt the method of Sections 13, 14, 15, 16 and 17 in [2]  to compute\\
 $  \Gamma_{i}^{\left[s\atop{ s\atop s} \right]\times k} \; \text{for} \; s \leqslant i \leqslant \inf(3s,k) $
using the results obtained in Section \ref{sec 6} and Section \ref{sec 9}.

  \subsection{Computation of  $ \Gamma_{s+j}^{\left[s\atop{ s\atop s} \right] \times k }\quad for\quad j=0,1,2,\quad j \leq s-1, \quad k \geq s+j $}
  \label{subsec 3}
  \begin{align*}
&
\end{align*}
 We have respectively  from \eqref{eq 6.1} and \eqref{eq 6.2} with $ m=l=0 $\\
       \begin{align}
  \sum_{ i = 0}^{\inf(3s,k)}  \Gamma  _{i}^{\left[s\atop{ s\atop s} \right]\times k} & = 2^{3k+3s-3} \label{eq 10.1}\\
 &  \text{and}\nonumber \\
  \sum_{i = 0}^{\inf(3s,k)}  \Gamma  _{i}^{\left[s\atop{ s\atop s} \right]\times k}\cdot2^{-i} & =  2^{2k+3s-3} + 2^{3k-3} - 2^{2k-3}. \label{eq 10.2}
\end{align}
We get   respectively from \eqref{eq 6.39} and \eqref{eq 6.40} \\
    \begin{align}
  \sum_{ i = 0}^{s-1}  \Gamma  _{i}^{\left[s\atop{ s\atop s} \right]\times k}& = 7\cdot2^{4s-6}-3\cdot2^{3s-5} \quad \text{if }\quad k\geq s \label{eq 10.3}\\
 &  \text{and}\nonumber \\
  \sum_{i = 0}^{s-1}  \Gamma  _{i}^{\left[s\atop{ s\atop s} \right]\times k}\cdot2^{-i} & = 15\cdot2^{3s-6}-7\cdot2^{2s-5} \quad \text{if } \quad k\geq s  \label{eq 10.4}
\end{align}

\begin{lem}
\label{lem 10.2}
We have \\

 \begin{align}
  \Gamma  _{s}^{\left[s\atop{ s\atop s} \right]\times k} & = 
  7\cdot2^{k+s-1} -7\cdot2^{2s} +105\cdot2^{4s-6} -21\cdot2^{3s-5}\quad \text{if}\quad k\geq s+1,\quad s\geq 2   \label{eq 10.5}\\
  &  \nonumber \\
   \Gamma  _{s}^{\left[s\atop{ s\atop s} \right]\times s} & = 2^{6s-3}-7\cdot2^{4s-6}+3\cdot2^{3s-5}\quad \text{if} \quad s\geq 2  \label{eq 10.6}\\
     &  \nonumber \\
   \Gamma  _{s+1}^{\left[s\atop{ s\atop s} \right]\times k}
 &  =  147\cdot2^{k+s-1} +105\cdot2^{4s-2} -21\cdot2^{3s-2} -315\cdot2^{2s} \quad \text{if}\quad k\geq s+2,\quad s\geq 2  \label{eq 10.7}\\
   &  \nonumber \\
     \Gamma  _{s+1}^{\left[s\atop{ s\atop s} \right]\times (s+1)}& =  2^{6s}-7\cdot2^{4s-2}+3\cdot2^{3s-2}\quad \text{if} \quad s\geq 2  \label{eq 10.8} \\
       &  \nonumber \\
   \Gamma  _{s+2}^{\left[s\atop{ s\atop s} \right]\times k} & = 147 \cdot 9 \cdot2^{k+s} +21\cdot\big[5\cdot2^{4s+2} -2^{3s+1} -275\cdot2^{2s+1}\big] \quad \text{if}\quad k\geq s+3,\quad s\geq 3  \label{eq 10.9}\\
     &  \nonumber \\
      \Gamma  _{s+2}^{\left[s\atop{ s\atop s} \right]\times (s+2)} &  =  2^{6s+3}-7\cdot2^{4s+2}+3\cdot2^{3s+1}+7\cdot2^{2s+1} \quad \text{if} \quad s\geq 2 \label{eq 10.10} \\
 &  \nonumber \\
      \Gamma  _{s+3}^{\left[s\atop{ s\atop s} \right]\times (s+3)} &  =  2^{6s+6}-7\cdot2^{4s+6}+3\cdot2^{3s+4}+21\cdot2^{2s+5} \quad \text{if}\quad s\geq 3  \label{eq 10.11} 
\end{align}
\end{lem}
\begin{proof}
\underline{Proof of \eqref{eq 10.5}}\\
By \eqref{eq 6.42}, we get\\
\begin{equation}
\label{eq 10.12}
 \Gamma  _{s}^{\left[s\atop{ s\atop s} \right]\times (s+1)}  = 
   105\cdot2^{4s-6} -21\cdot2^{3s-5}
\end{equation}
We have from \eqref{eq 9.5} with j=0, $k>s$\\
\begin{equation}
\label{eq 10.13}
 \Delta _{k}\Gamma  _{s}^{\left[s\atop{ s\atop s} \right]}
  =   \Gamma  _{s}^{\left[s\atop{ s\atop s} \right]\times (k+1)} -  \Gamma  _{s}^{\left[s\atop{ s\atop s} \right]\times k} = 7\cdot2^{k+s-1}
  \end{equation}
 We obtain from \eqref{eq 10.12} and  \eqref{eq 10.13} for $k\geq s+1$
  \begin{align}
& \sum_{i=s+1}^{k}\big( \Gamma  _{s}^{\left[s\atop{ s\atop s} \right]\times (i+1)} -  \Gamma  _{s}^{\left[s\atop{ s\atop s} \right]\times i})
= \sum_{i=s+1}^{k}7\cdot2^{k+s-1} \nonumber \\
& \Updownarrow  \nonumber \\
& \Gamma  _{s}^{\left[s\atop{ s\atop s} \right]\times (k+1)} -  \Gamma  _{s}^{\left[s\atop{ s\atop s} \right]\times (s+1)} = 7\cdot2^{k+s} - 7\cdot2^{2s}
\nonumber \\
& \Updownarrow  \nonumber \\
& \Gamma  _{s}^{\left[s\atop{ s\atop s} \right]\times k} = 7\cdot2^{k+s-1} - 7\cdot2^{2s} +  105\cdot2^{4s-6} -21\cdot2^{3s-5}\quad \text{for $k\geq s+2$}\label{eq 10.14}
\end{align}
 It follows from \eqref{eq 10.12} that \eqref{eq 10.14} is valid for $k\geq s+1$\\
 
\underline{proof of \eqref{eq 10.6}}
By \eqref{eq 6.41} with  $ m = l = 0, i = s.$ \\

\underline{proof of \eqref{eq 10.7} and \eqref{eq 10.10}}
 We get respectively  from \eqref{eq 10.1} and \eqref{eq 10.2} with $ k=s+2 $\\
 
       \begin{align}
  \sum_{ i = 0}^{s+2}  \Gamma  _{i}^{\left[s\atop{ s\atop s} \right]\times (s+2)} & = 2^{6s+3} \label{eq 10.15}\\
 &  \text{and}\nonumber \\
  \sum_{i = 0}^{s+2}  \Gamma  _{i}^{\left[s\atop{ s\atop s} \right]\times (s+2)}\cdot2^{-i} & =  2^{5s+1} + 2^{3s+3} - 2^{2s+1}. \label{eq 10.16}
\end{align}
Applying \eqref{eq 10.3} and \eqref{eq 10.5} with $k = s+2$ we get from \eqref{eq 10.15}\\
\begin{align*}
&   \sum_{ i = 0}^{s+2}  \Gamma  _{i}^{\left[s\atop{ s\atop s} \right]\times (s+2)}=
  \sum_{ i = 0}^{s -1}  \Gamma  _{i}^{\left[s\atop{ s\atop s} \right]\times (s+2)}  
  + \Gamma  _{s}^{\left[s\atop{ s\atop s} \right]\times (s+2)} + \Gamma  _{s+1}^{\left[s\atop{ s\atop s} \right]\times (s+2)} + \Gamma  _{s+2}^{\left[s\atop{ s\atop s} \right]\times (s+2)}\\
  & = \big( 7\cdot2^{4s-6}-3\cdot2^{3s-5}\big) +  \big(7\cdot2^{2s+1} - 7\cdot2^{2s} +  105\cdot2^{4s-6} -21\cdot2^{3s-5}\big)
  +  \Gamma  _{s+1}^{\left[s\atop{ s\atop s} \right]\times (s+2)} + \Gamma  _{s+2}^{\left[s\atop{ s\atop s} \right]\times (s+2)} = 2^{6s+3} 
\end{align*}
  Hence we obtain from the above equation \\
 \begin{equation}
\label{eq 10.17}
 \Gamma  _{s+1}^{\left[s\atop{ s\atop s} \right]\times (s+2)} + \Gamma  _{s+2}^{\left[s\atop{ s\atop s} \right]\times (s+2)}
  = 2^{6s+3} - 7\cdot2^{4s-2}+ 3\cdot2^{3s-2} - 7\cdot2^{2s}
\end{equation}
Applying \eqref{eq 10.4} and \eqref{eq 10.5} with $k = s+2$ we get from \eqref{eq 10.16}\\
\begin{align*}
&   \sum_{ i = 0}^{s+2}  \Gamma  _{i}^{\left[s\atop{ s\atop s} \right]\times (s+2)}\cdot2^{-i}=
  \sum_{ i = 0}^{s -1}  \Gamma  _{i}^{\left[s\atop{ s\atop s} \right]\times (s+2)}\cdot2^{-i}  
  + \Gamma  _{s}^{\left[s\atop{ s\atop s} \right]\times (s+2)}\cdot2^{-s} + \Gamma  _{s+1}^{\left[s\atop{ s\atop s} \right]\times (s+2)}\cdot2^{-(s+1)} + \Gamma  _{s+2}^{\left[s\atop{ s\atop s} \right]\times (s+2)}\cdot2^{-(s+2}\\
  & = \big( 15\cdot2^{3s-6}-7\cdot2^{2s-5}\big) +  \big(7\cdot2^{s} +  105\cdot2^{3s-6} -21\cdot2^{2s-5}\big)\\
  &  +  \Gamma  _{s+1}^{\left[s\atop{ s\atop s} \right]\times (s+2)}\cdot2^{-(s+1)} + \Gamma  _{s+2}^{\left[s\atop{ s\atop s} \right]\times (s+2)}\cdot2^{-(s+2)} = 2^{5s+1} + 2^{3s+3} - 2^{2s+1}
\end{align*}
We then get from the above equation \\
\begin{align}
& \Gamma  _{s+1}^{\left[s\atop{ s\atop s} \right]\times (s+2)}\cdot2^{-(s+1)} + \Gamma  _{s+2}^{\left[s\atop{ s\atop s} \right]\times (s+2)}\cdot2^{-(s+2)}
  = 2^{5s+1} + 49\cdot2^{3s-3} - 9\cdot2^{2s-3}-7\cdot2^{s} \nonumber \\
&  \Updownarrow  \nonumber \\
& 2\cdot \Gamma  _{s+1}^{\left[s\atop{ s\atop s} \right]\times (s+2)} + \Gamma  _{s+2}^{\left[s\atop{ s\atop s} \right]\times (s+2)}
  = 2^{6s+3} + 49\cdot2^{4s-1} - 9\cdot2^{3s-1}-7\cdot2^{2s+2} \label{eq 10.18} 
\end{align}
Now we deduce from \eqref{eq 10.17} and \eqref{eq 10.18}\\
\begin{align}
& \Gamma  _{s+1}^{\left[s\atop{ s\atop s} \right]\times (s+2)} 
  =105\cdot 2^{4s-2} -21\cdot2^{3s-2} - 21\cdot2^{2s} \label{eq 10.19} \\
&  \Gamma  _{s+2}^{\left[s\atop{ s\atop s} \right]\times (s+2)}
  = 2^{6s+3} -7\cdot2^{4s+2} +3\cdot2^{3s+1}+7\cdot2^{2s+1} \label{eq 10.20} 
\end{align}
We have from \eqref{eq 9.5} with j = 1, $k\geq s+2$ \\
  \begin{equation}
\label{eq 10.21}
 \Delta _{k}\Gamma  _{s+1}^{\left[s\atop{ s\atop s} \right]}
  =   \Gamma  _{s+1}^{\left[s\atop{ s\atop s} \right]\times (k+1)} -  \Gamma  _{s+1}^{\left[s\atop{ s\atop s} \right]\times k} =147\cdot2^{k+s-3}
  \end{equation}
 We obtain from \eqref{eq 10.21} and  \eqref{eq 10.19} for $k\geq s+2$
  \begin{align}
& \sum_{i=s+2}^{k}\big( \Gamma  _{s+1}^{\left[s\atop{ s\atop s} \right]\times (i+1)} -  \Gamma  _{s+1}^{\left[s\atop{ s\atop s} \right]\times i})
= \sum_{i=s+1}^{k}147\cdot2^{k+s-3} \nonumber \\
& \Updownarrow  \nonumber \\
& \Gamma  _{s+1}^{\left[s\atop{ s\atop s} \right]\times (k+1)} -  \Gamma  _{s+1}^{\left[s\atop{ s\atop s} \right]\times (s+2)} =147\cdot2^{k+s} -147\cdot2^{2s+1}
\nonumber \\
& \Updownarrow  \nonumber \\
& \Gamma  _{s+1}^{\left[s\atop{ s\atop s} \right]\times k}
 =147\cdot2^{k+s-1} - 147\cdot2^{2s+1} +  \big[105\cdot2^{4s-2} -21\cdot2^{3s-2} -21\cdot2^{2s}\big] \quad \text{for $k\geq s+3$}\label{eq 10.22}
\end{align}
 It follows from \eqref{eq 10.19} that \eqref{eq 10.22} is valid for $k\geq s+2$\\
 
 \underline{proof of \eqref{eq 10.6}} 
 
By \eqref{eq 6.41} with  $ m = l = 0, i = s.$ \\
 
  \underline{Proof of \eqref{eq 10.9} and \eqref{eq 10.11}}
  .
  We get respectively  from \eqref{eq 10.1} and \eqref{eq 10.2} with $ k=s+3 $\\
 
       \begin{align}
  \sum_{ i = 0}^{s+3}  \Gamma  _{i}^{\left[s\atop{ s\atop s} \right]\times (s+3)} & = 2^{6s+6} \label{eq 10.23}\\
 &  \text{and}\nonumber \\
  \sum_{i = 0}^{s+3}  \Gamma  _{i}^{\left[s\atop{ s\atop s} \right]\times (s+3)}\cdot2^{-i} & =  2^{5s+3} + 2^{3s+6} - 2^{2s+3}. \label{eq 10.24}
\end{align}
Applying \eqref{eq 10.3},  \eqref{eq 10.5} and  \eqref{eq 10.7} with $k = s+3$ we get from \eqref{eq 10.23}\\
\begin{align*}
&   \sum_{ i = 0}^{s+3}  \Gamma  _{i}^{\left[s\atop{ s\atop s} \right]\times (s+3)}=
  \sum_{ i = 0}^{s -1}  \Gamma  _{i}^{\left[s\atop{ s\atop s} \right]\times (s+3)}  
  + \Gamma  _{s}^{\left[s\atop{ s\atop s} \right]\times (s+3)} + \Gamma  _{s+1}^{\left[s\atop{ s\atop s} \right]\times (s+3)} + \Gamma _{s+2}^{\left[s\atop{ s\atop s} \right]\times (s+3)}
  + \Gamma  _{s+3}^{\left[s\atop{ s\atop s} \right]\times (s+3)} \\
   & = \big( 7\cdot2^{4s-6}-3\cdot2^{3s-5}\big) +  \big( 21\cdot2^{2s} +  105\cdot2^{4s-6} -21\cdot2^{3s-5}\big) \\
   &  +  \big(147\cdot2^{2s+2} + 105\cdot2^{4s-2} -21\cdot2^{3s-2} -315\cdot2^{2s}\big)
     +  \Gamma  _{s+2}^{\left[s\atop{ s\atop s} \right]\times (s+3)} + \Gamma  _{s+3}^{\left[s\atop{ s\atop s} \right]\times (s+3)} = 2^{6s+6} 
\end{align*}
  Hence we obtain from the above equation \\
 \begin{equation}
\label{eq 10.25}
 \Gamma  _{s+2}^{\left[s\atop{ s\atop s} \right]\times (s+3)} + \Gamma  _{s+3}^{\left[s\atop{ s\atop s} \right]\times (s+3)}
  = 2^{6s+6} - 7\cdot2^{4s+2}+ 3\cdot2^{3s+1} - 147\cdot2^{2s+1}
\end{equation}
Applying \eqref{eq 10.4},  \eqref{eq 10.5} and \eqref{eq 10.7} with $k = s+3$ we get from \eqref{eq 10.24}\\
\begin{align*}
&   \sum_{ i = 0}^{s+3}  \Gamma  _{i}^{\left[s\atop{ s\atop s} \right]\times (s+3)}\cdot2^{-i}=
  \sum_{ i = 0}^{s -1}  \Gamma  _{i}^{\left[s\atop{ s\atop s} \right]\times (s+3)}\cdot2^{-i}  
  + \Gamma  _{s}^{\left[s\atop{ s\atop s} \right]\times (s+3)}\cdot2^{-s} + \Gamma  _{s+1}^{\left[s\atop{ s\atop s} \right]\times (s+3)}\cdot2^{-(s+1)}\\
  &   + \Gamma  _{s+2}^{\left[s\atop{ s\atop s} \right]\times (s+3)}\cdot2^{-(s+2)}
  + \Gamma  _{s+3}^{\left[s\atop{ s\atop s} \right]\times (s+3)}\cdot2^{-(s+3)} \\
   & = \big( 15\cdot2^{3s-6}-7\cdot2^{2s-5}\big) +  \big(21\cdot2^{s} +  105\cdot2^{3s-6} -21\cdot2^{2s-5}\big)
    +  \big(147\cdot2^{s+1} + 105\cdot2^{3s-3} -21\cdot2^{2s-3} -315\cdot2^{s-1}\big) \\
      &  +  \Gamma  _{s+2}^{\left[s\atop{ s\atop s} \right]\times (s+3)}\cdot2^{-(s+2)} + \Gamma  _{s+3}^{\left[s\atop{ s\atop s} \right]\times (s+3)}\cdot2^{-(s+3)} = 2^{5s+3} + 2^{3s+6} - 2^{2s+3}
\end{align*}
We then get from the above equation \\
\begin{align}
& \Gamma  _{s+2}^{\left[s\atop{ s\atop s} \right]\times (s+3)}\cdot2^{-(s+2)} + \Gamma  _{s+3}^{\left[s\atop{ s\atop s} \right]\times (s+3)}\cdot2^{-(s+3)}
  = 2^{5s+3} + 49\cdot2^{3s} - 9\cdot2^{2s-1}-315\cdot2^{s-1} \nonumber \\
&  \Updownarrow  \nonumber \\
& 2\cdot \Gamma  _{s+2}^{\left[s\atop{ s\atop s} \right]\times (s+3)} + \Gamma  _{s+3}^{\left[s\atop{ s\atop s} \right]\times (s+3)}
  = 2^{6s+6} + 49\cdot2^{4s+3} - 9\cdot2^{3s+2}-315\cdot2^{2s+2} \label{eq 10.26} 
\end{align}
Now we deduce from \eqref{eq 10.25} and \eqref{eq 10.26}\\
\begin{align}
& \Gamma  _{s+2}^{\left[s\atop{ s\atop s} \right]\times (s+3)} 
  = 21\cdot\big[5\cdot 2^{4s+2} -2^{3s+1} - 23\cdot2^{2s+1}\big] \label{eq 10.27} \\
&  \Gamma  _{s+3}^{\left[s\atop{ s\atop s} \right]\times (s+3)}
  = 2^{6s+6} -7\cdot2^{4s+6} +3\cdot2^{3s+4}+21\cdot2^{2s+5} \label{eq 10.28} 
\end{align}
We have from \eqref{eq 9.5} with j = 2, $k\geq s+3$ \\
  \begin{equation}
\label{eq 10.29}
 \Delta _{k}\Gamma  _{s+2}^{\left[s\atop{ s\atop s} \right]}
  =   \Gamma  _{s+2}^{\left[s\atop{ s\atop s} \right]\times (k+1)} -  \Gamma  _{s+2}^{\left[s\atop{ s\atop s} \right]\times k} =147\cdot9\cdot2^{k+s}
  \end{equation}
 We obtain from \eqref{eq 10.29} and  \eqref{eq 10.27} for $k\geq s+3$
  \begin{align}
& \sum_{i=s+3}^{k}\big( \Gamma  _{s+2}^{\left[s\atop{ s\atop s} \right]\times (i+1)} -  \Gamma  _{s+2}^{\left[s\atop{ s\atop s} \right]\times i})
= \sum_{i=s+3}^{k}147\cdot9\cdot2^{k+s} \nonumber \\
& \Updownarrow  \nonumber \\
& \Gamma  _{s+2}^{\left[s\atop{ s\atop s} \right]\times (k+1)} -  \Gamma  _{s+2}^{\left[s\atop{ s\atop s} \right]\times (s+3)} =147\cdot9\cdot2^{k+s+1} -147\cdot9\cdot2^{2s+3}
\nonumber \\
& \Updownarrow  \nonumber \\
& \Gamma  _{s+2}^{\left[s\atop{ s\atop s} \right]\times k}
 =147\cdot9\cdot2^{k+s} +21\cdot\big[5\cdot2^{4s+2} -2^{3s+1} -275\cdot2^{2s+1}\big] \quad \text{for $k\geq s+4$}\label{eq 10.30}
\end{align}
 It follows from \eqref{eq 10.27} that \eqref{eq 10.30} is valid for $k\geq s+3$\\
  \end{proof}
   \subsection{An explicit formula  for  $ \Gamma_{s+j}^{\left[s\atop{ s\atop s} \right] \times k }\quad for\quad 1\leq j\leq s-1,  \quad k \geq s+j $}
  \label{subsec 4}
\begin{lem}
\label{lem 10.3}
We have \\
 \begin{align}
(H_{j}) \quad \Gamma  _{s+j}^{\left[s\atop{ s\atop s} \right]\times k} & = 
  147\cdot(5\cdot2^{j-1} - 1)\cdot2^{k+s+3j-6} \label{eq 10.31} \\
   & + 21\cdot\big[5\cdot2^{4s+4j-6} - 2^{3s+3j-5} -(155\cdot2^{j-1} - 35)\cdot2^{2s+4j-7} \big] \quad \text{for}\quad 1\leq j\leq s-1,\; k\geq s+j+1.\nonumber \\
    \Gamma  _{s+j}^{\left[s\atop{ s\atop s} \right]\times (s+j)} & =  
    2^{6s+3j-3} +7\cdot2^{2s+5j-8} -7\cdot2^{2s+4j-7} -7\cdot2^{4s+4j-6} + 3\cdot2^{3s+3j-5}\quad \text{for}\quad 1\leq j\leq s. \label{eq 10.32} 
\end{align}
\end{lem}
\begin{proof}
The proof is by strong induction, that is:\\
If
\begin{itemize}
\item \text{$ (H_{1}) $ is true, and  }
\item \text{for all $ j\geq 1, \quad (H_{1})\wedge (H_{2})\wedge\ldots\wedge (H_{j})\quad implies \; (H_{j+1}) $}
\end{itemize}
then $ (H_{j}) $ is true for all $ j \geq 1. $\\

\underline {$ (H_{1})\, is \;true\; $}\\
Indeed from  \eqref{eq 10.7}  we have  $$  \Gamma  _{s+1}^{\left[s\atop{ s\atop s} \right]\times k}
  =  147\cdot2^{k+s-1} +105\cdot2^{4s-2} -21\cdot2^{3s-2} -315\cdot2^{2s} \quad \text{if}\quad k\geq s+2,\;s\geq 2 $$\\
  and $ (H_{j})$ with j=1 gives \\
 $$  147\cdot(5 - 1)\cdot2^{k+s-3} + 21\cdot\big[5\cdot2^{4s-2} - 2^{3s-2} -(155 - 35)\cdot2^{2s-3} \big] =  \Gamma  _{s+1}^{\left[s\atop{ s\atop s} \right]\times k} $$
 From \eqref{eq 10.8} we see that $\Gamma  _{s+1}^{\left[s\atop{ s\atop s} \right]\times (s+1)}$ is equal to  \eqref{eq 10.32} with j=1. \\
   \underline {$ (H_{1})\; implies \; (H_{2}) $}
   
  Recall the proof of \eqref{eq 10.9}.
 
  \begin{align*}
 \eqref{eq 10.9}\quad  \Gamma  _{s+2}^{\left[s\atop{ s\atop s} \right]\times k} 
  = 147 \cdot 9 \cdot2^{k+s} +21\cdot\big[5\cdot2^{4s+2} -2^{3s+1} -275\cdot2^{2s+1}\big] \quad \text{if}\quad k\geq s+3, \quad s\geq 3. 
\end{align*}

 $  (H_{j})  $ with j = 2 gives \\
 \begin{align*}
    147\cdot(5\cdot2 - 1)\cdot2^{k+s}  + 21\cdot\big[5\cdot2^{4s+2} - 2^{3s+1} -(155\cdot2 - 35)\cdot2^{2s+1} \big] 
=  \Gamma  _{s+2}^{\left[s\atop{ s\atop s} \right]\times k} 
\end{align*}
 From \eqref{eq 10.10} we see that $\Gamma  _{s+2}^{\left[s\atop{ s\atop s} \right]\times (s+2)}$ is equal to  \eqref{eq 10.32} with j=2. \\

  \underline {$ (H_{1})\wedge (H_{2})\wedge\ldots\wedge (H_{j})\quad implies \; (H_{j+1}) $}

 We get respectively from \eqref{eq 10.1} and \eqref{eq 10.2} with  k = s+j+2 \\
    \begin{align}
  \sum_{ i = 0}^{s+j+2}  \Gamma  _{i}^{\left[s\atop{ s\atop s} \right]\times (s+j+2)} & = 2^{6s+3j+3} \label{eq 10.33}\\
 &  \text{and}\nonumber \\
  \sum_{i = 0}^{s+j+2}  \Gamma  _{i}^{\left[s\atop{ s\atop s} \right]\times (s+j+2)}\cdot2^{-i} & =  2^{5s+2j+1} + 2^{3s+3j+3} - 2^{2s+2j+1}. \label{eq 10.34}
\end{align}

 We have respectively  from \eqref{eq 10.3} and \eqref{eq 10.4} with  k=s+j+2 \\
 \begin{align}
   \sum_{ i = 0}^{s-1}  \Gamma  _{i}^{\left[s\atop{ s\atop s} \right]\times (s+j+2)}& = 7\cdot2^{4s-6}-3\cdot2^{3s-5} \label{eq 10.35}\\
    \sum_{i = 0}^{s-1}  \Gamma  _{i}^{\left[s\atop{ s\atop s} \right]\times (s+j+2)}\cdot2^{-i} & = 15\cdot2^{3s-6}-7\cdot2^{2s-5}\label{eq 10.36} 
\end{align}

 By \eqref{eq 10.5} with k=s+j+2 \\
\begin{align}
  \Gamma  _{s}^{\left[s\atop{ s\atop s} \right]\times (s+j+2)}= 7\cdot2^{2s+j+1} -7\cdot2^{2s} +105\cdot2^{4s-6} -21\cdot2^{3s-5}\label{eq 10.37} \\
  & \nonumber \\
    \Gamma  _{s}^{\left[s\atop{ s\atop s} \right]\times (s+j+2)}\cdot2^{-s}= 7\cdot2^{s+j+1} -7\cdot2^{s} +105\cdot2^{3s-6} -21\cdot2^{2s-5}\label{eq 10.38} 
\end{align}

We get from  $ (H_{q}) $ with k = s+j+2 and $1 \leq q \leq  j $ \\
\begin{align}
 & \Gamma  _{s+q}^{\left[s\atop{ s\atop s} \right]\times (s+j+2)}=  147\cdot(5\cdot2^{q-1} - 1)\cdot2^{(s+j+2)+s+3q-6}
   + 21\cdot\big[5\cdot2^{4s+4q-6} - 2^{3s+3q-5} -(155\cdot2^{q-1} - 35)\cdot2^{2s+4q-7} \big] \nonumber \\
   & \Updownarrow  \nonumber \\
  & \Gamma  _{s+q}^{\left[s\atop{ s\atop s} \right]\times (s+j+2)} = 2^{5q}\cdot\big[-3255\cdot2^{2s-8}\big] +
   2^{4q}\cdot\big[735\cdot2^{2s+j-5} +105\cdot2^{4s-6} + 735\cdot2^{2s-7}\big]
   +  2^{3q}\cdot\big[-147\cdot2^{2s+j-4} -21\cdot2^{3s-5} \big] \label{eq 10.39}\\
   & \nonumber \\
   & \Gamma  _{s+q}^{\left[s\atop{ s\atop s} \right]\times (s+j+2)}\cdot2^{-(s+q)} = 
   2^{4q}\cdot\big[-3255\cdot2^{s-8}\big] + 2^{3q}\cdot\big[735\cdot2^{s+j-5} +105\cdot2^{3s-6} + 735\cdot2^{s-7}\big]
   +  2^{2q}\cdot\big[-147\cdot2^{s+j-4} -21\cdot2^{2s-5} \big] \label{eq 10.40}         
    \end{align}
     We obtain from  \eqref{eq 10.33}, \eqref{eq 10.35},   \eqref{eq 10.37}  and  \eqref{eq 10.39}.
\begin{align}
&   \sum_{ i = 0}^{s+j+2}  \Gamma  _{i}^{\left[s\atop{ s\atop s} \right]\times (s+j+2)}=
  \sum_{ i = 0}^{s -1}  \Gamma  _{i}^{\left[s\atop{ s\atop s} \right]\times (s+j+2)}  
  + \Gamma  _{s}^{\left[s\atop{ s\atop s} \right]\times (s+j+2)} +
   \sum_{q=1}^{j}\Gamma  _{s+q}^{\left[s\atop{ s\atop s} \right]\times (s+j+2)} + \Gamma _{s+j+1}^{\left[s\atop{ s\atop s} \right]\times (s+j+2)}
  + \Gamma  _{s+j+2}^{\left[s\atop{ s\atop s} \right]\times (s+j+2)}\nonumber \\
   & = \big( 7\cdot2^{4s-6}-3\cdot2^{3s-5}\big) +  \big(7\cdot2^{2s+j+1} -7\cdot2^{2s} +105\cdot2^{4s-6} -21\cdot2^{3s-5}\big) \nonumber \\
& +   \big[-3255\cdot2^{2s-8}\big]\cdot \sum_{q=1}^{j} 2^{5q} + \big[735\cdot2^{2s+j-5}
 +105\cdot2^{4s-6} + 735\cdot2^{2s-7}\big]\cdot \sum_{q=1}^{j} 2^{4q}
   + \big[-147\cdot2^{2s+j-4} -21\cdot2^{3s-5} \big]\cdot\sum_{q=1}^{j}2^{3q} \nonumber \\
    & +  \Gamma  _{s+j+1}^{\left[s\atop{ s\atop s} \right]\times (s+j+2)} + \Gamma  _{s+j+2}^{\left[s\atop{ s\atop s} \right]\times (s+j+2)}
     = 2^{6s+3j+3}\nonumber \\
    & \Updownarrow  \nonumber\\
   &  \Gamma  _{s+j+1}^{\left[s\atop{ s\atop s} \right]\times (s+j+2)} + \Gamma  _{s+j+2}^{\left[s\atop{ s\atop s} \right]\times (s+j+2)} = 2^{6s+3j+3}
     -91\cdot2^{2s+5j-3}+35\cdot2^{2s+4j-3}-7\cdot2^{4s+4j-2} +3\cdot2^{3s+3j-2} \label{eq 10.41}
\end{align} 
  We get in the same way  from  \eqref{eq 10.34}, \eqref{eq 10.36},   \eqref{eq 10.38}  and  \eqref{eq 10.40}.
  \begin{align*}
&   \sum_{ i = 0}^{s+j+2}  \Gamma  _{i}^{\left[s\atop{ s\atop s} \right]\times (s+j+2)}\cdot2^{-i}=
  \sum_{ i = 0}^{s -1}  \Gamma  _{i}^{\left[s\atop{ s\atop s} \right]\times (s+j+2)}\cdot2^{-i}  
  + \Gamma  _{s}^{\left[s\atop{ s\atop s} \right]\times (s+j+2)}\cdot2^{-s} +
   \sum_{q=1}^{j}\Gamma  _{s+q}^{\left[s\atop{ s\atop s} \right]\times (s+j+2)}\cdot2^{-(s+q)}\\
  & + \Gamma _{s+j+1}^{\left[s\atop{ s\atop s} \right]\times (s+j+2)}\cdot2^{-(s+j+1)}
   + \Gamma  _{s+j+2}^{\left[s\atop{ s\atop s} \right]\times (s+j+2)}\cdot2^{-(s+j+2)} \\
    & = \big(15\cdot2^{3s-6}-7\cdot2^{2s-5} \big) +  \big( 7\cdot2^{s+j+1} -7\cdot2^{s} +105\cdot2^{3s-6} -21\cdot2^{2s-5}  \big) \\
& +   \big[-3255\cdot2^{s-8}\big]\cdot \sum_{q=1}^{j} 2^{4q} + \big[735\cdot2^{s+j-5} +105\cdot2^{3s-6}
 + 735\cdot2^{s-7}\big]\cdot \sum_{q=1}^{j} 2^{3q}
   + \big[-147\cdot2^{s+j-4} -21\cdot2^{2s-5} \big]\cdot\sum_{q=1}^{j}2^{2q}\\
    & +  \Gamma  _{s+j+1}^{\left[s\atop{ s\atop s} \right]\times (s+j+2)}\cdot2^{-(s+j+1)}
     + \Gamma  _{s+j+2}^{\left[s\atop{ s\atop s} \right]\times (s+j+2)}\cdot2^{-(s+j+2)} =  2^{5s+2j+1} + 2^{3s+3j+3} - 2^{2s+2j+1} \\
    & \Updownarrow \\
   &  \Gamma  _{s+j+1}^{\left[s\atop{ s\atop s} \right]\times (s+j+2)}\cdot2^{-(s+j+1)}
    + \Gamma  _{s+j+2}^{\left[s\atop{ s\atop s} \right]\times (s+j+2)}\cdot2^{-(s+j+2)} \\
    & =   2^{5s+2j+1} +49 \cdot2^{3s+3j-3}-203\cdot2^{s+4j-4} +91\cdot2^{s+3j-4} -9 \cdot2^{2s+2j-3}
    \end{align*}
    \begin{align}
     & \Updownarrow  \nonumber \\
    & 2\cdot\Gamma  _{s+j+1}^{\left[s\atop{ s\atop s} \right]\times (s+j+2)} + \Gamma  _{s+j+2}^{\left[s\atop{ s\atop s} \right]\times (s+j+2)} = 
    2^{6s+3j+3} +49 \cdot2^{4s+4j-1}-203\cdot2^{2s+5j-2} +91\cdot2^{2s+4j-2} -9 \cdot2^{3s+3j-1}\label{eq 10.42}
     \end{align}
Hence by \eqref{eq 10.41}, \eqref{eq 10.42} we deduce \\
\begin{align}
& \Gamma  _{s+j+1}^{\left[s\atop{ s\atop s} \right]\times (s+j+2)} = 
 105 \cdot2^{4s+4j-2}-315\cdot2^{2s+5j-3} +147\cdot2^{2s+4j-3} -21 \cdot2^{3s+3j-2}\label{eq 10.43} \\
 & \nonumber \\
 & \Gamma  _{s+j+2}^{\left[s\atop{ s\atop s} \right]\times (s+j+2)} = 
 2^{6s+3j+3} +7\cdot2^{2s+5j+2} -7\cdot2^{2s+4j+1} -7\cdot2^{4s+4j+2} +3 \cdot2^{3s+3j+1}\label{eq 10.44}
\end{align}
We have from \eqref{eq 9.5} with $ j\rightarrow j+1,\quad j+1\leq s-1,\quad k\geq s+j+2$ \\
  \begin{equation}
\label{eq 10.45}
 \Delta _{k}\Gamma  _{s+j+1}^{\left[s\atop{ s\atop s} \right]}
  =   \Gamma  _{s+j+1}^{\left[s\atop{ s\atop s} \right]\times (k+1)}
   -  \Gamma  _{s+j+1}^{\left[s\atop{ s\atop s} \right]\times k} =147\cdot(5\cdot2^{j}-1)\cdot2^{k+s+3j-3}
  \end{equation}
 We obtain from \eqref{eq 10.45} and  \eqref{eq 10.43} for $k\geq s+j+2$
  \begin{align}
& \sum_{i=s+j+2}^{k}\big( \Gamma  _{s+j+1}^{\left[s\atop{ s\atop s} \right]\times (i+1)} -  \Gamma  _{s+j+1}^{\left[s\atop{ s\atop s} \right]\times i})
= \sum_{i=s+j+2}^{k}147\cdot(5\cdot2^{j}-1)\cdot2^{k+s +3j-3} \nonumber \\
& \Updownarrow  \nonumber \\
& \Gamma  _{s+j+1}^{\left[s\atop{ s\atop s} \right]\times (k+1)} - 
 \Gamma  _{s+j+1}^{\left[s\atop{ s\atop s} \right]\times (s+j+2)} =147\cdot(5\cdot2^{j}-1)\cdot2^{k+s+3j-2} -147\cdot(5\cdot2^{j}-1)\cdot2^{2s+4j-1}
\nonumber \\
& \Updownarrow  \nonumber \\
& \Gamma  _{s+j+1}^{\left[s\atop{ s\atop s} \right]\times k}
 = 147\cdot(5\cdot2^{j}-1)\cdot2^{k+s+3j-3} -147\cdot(5\cdot2^{j}-1)\cdot2^{2s+4j-1}\label{eq 10.46} \\
 & +  105 \cdot2^{4s+4j-2}-315\cdot2^{2s+5j-3} +147\cdot2^{2s+4j-3} -21 \cdot2^{3s+3j-2}\quad \text{for}\quad k\geq s+j+3 \nonumber
\end{align}
 It follows from \eqref{eq 10.43} that \eqref{eq 10.46} is valid for $k\geq s+j+2$\\
   By \eqref{eq 10.31} we have  with $j\rightarrow j+1$\\
  \begin{align}
(H_{j+1}) \quad \Gamma  _{s+j+1}^{\left[s\atop{ s\atop s} \right]\times k} & = 
  147\cdot(5\cdot2^{j} - 1)\cdot2^{k+s+3j-3}\label{eq 10.47} \\
   & + 21\cdot\big[5\cdot2^{4s+4j-2} - 2^{3s+3j-2} -(155\cdot2^{j}
    - 35)\cdot2^{2s+4j-3} \big] \quad \text{for}\quad 1\leq j+1\leq s-1,\; k\geq s+j+2.\nonumber 
\end{align}
By comparing \eqref{eq 10.46} with \eqref{eq 10.47} we see that $ (H_{j+1}) $ holds.

\end{proof}

   \subsection{ Explicit formulas  for 
    $ \Gamma_{2s}^{\left[s\atop{ s\atop s} \right] \times k }\quad 
     where \quad k \geq 2s+1\quad and \quad  for \quad \Gamma_{2s+1}^{\left[s\atop{ s\atop s} \right] \times (2s+1) } $}
  \label{subsec 5}
  \begin{lem}
  \label{lem 10.4}
  We have \\
 \begin{align}
 \Gamma  _{2s}^{\left[s\atop{ s\atop s} \right]\times k} & = 
  7\cdot2^{2k+2s-2}+ 21\cdot\big[35\cdot2^{k+5s-7} - 39\cdot2^{k+4s-6} \big] \label{eq 10.48} \\
  &  + 7\cdot\big[15\cdot2^{8s-6} - 465\cdot2^{7s-8} +349\cdot2^{6s-7} \big] \quad \text{for}\quad  k\geq 2s+1.\nonumber \\
    \Gamma  _{2s+1}^{\left[s\atop{ s\atop s} \right]\times (2s+1)} & =  
    2^{9s} -7\cdot2^{8s-2} +7\cdot2^{7s-3} - 2^{6s-3}  \label{eq 10.49} 
\end{align}
 \end{lem}
 \begin{proof}
  We get respectively from \eqref{eq 10.43} and \eqref{eq 10.44} with j = s-1 \\
   \begin{align}
& \Gamma  _{2s}^{\left[s\atop{ s\atop s} \right]\times (2s+1)} = 
 105 \cdot2^{8s-6}-315\cdot2^{7s-8} + 63\cdot2^{6s-7}  \label{eq 10.50} \\
& \Gamma  _{2s+1}^{\left[s\atop{ s\atop s} \right]\times (2s+1)} = 
 2^{9s} +7\cdot2^{7s-3} -7\cdot2^{8s-2}  - 2^{6s-3}\label{eq 10.51}
\end{align}
 We have from \eqref{eq 9.7} with j = 0 and $k \geq 2s+1$  \\
  \begin{equation}
\label{eq 10.52}
 \Delta _{k}\Gamma  _{2s}^{\left[s\atop{ s\atop s} \right]}
  =   \Gamma  _{2s}^{\left[s\atop{ s\atop s} \right]\times (k+1)} -  \Gamma  _{2s}^{\left[s\atop{ s\atop s} \right]\times k} =21\cdot\big(2^{2k+2s-2} +35\cdot2^{k+5s-7} -39\cdot2^{k+4s-6}\big) 
  \end{equation}
   We obtain from \eqref{eq 10.52} for $k\geq 2s+1$ and  \eqref{eq 10.50} 
  \begin{align*}
& \sum_{i=2s+1}^{k}\big( \Gamma  _{2s}^{\left[s\atop{ s\atop s} \right]\times (i+1)} -  \Gamma  _{2s}^{\left[s\atop{ s\atop s} \right]\times i})
= \sum_{i= 2s+1}^{k} 21\cdot\big(2^{2i+2s-2} +35\cdot2^{i+5s-7} -39\cdot2^{i+4s-6}\big)   \\
& \Updownarrow   \\
& \Gamma  _{2s}^{\left[s\atop{ s\atop s} \right]\times (k+1)}
= 21\cdot2^{2s-2}\cdot \sum_{i= 2s+1}^{k}2^{2i} + \big(735\cdot2^{5s-7} -819\cdot2^{4s-6} \big)\cdot \sum_{i= 2s+1}^{k}2^{i}
 + 105 \cdot2^{8s-6}-315\cdot2^{7s-8} + 63\cdot2^{6s-7}\\
& =  7\cdot2^{2k+2s}+ 735\cdot2^{k+5s-6} - 819\cdot2^{k+4s-5}     + 105\cdot2^{8s-6} - 3255\cdot2^{7s-8} +2443\cdot2^{6s-7} 
\end{align*}
\begin{align}
& \Updownarrow  \nonumber  \\
&  \Gamma  _{2s}^{\left[s\atop{ s\atop s} \right]\times k} = 
  7\cdot2^{2k+2s-2}+ 21\cdot\big[35\cdot2^{k+5s-7} - 39\cdot2^{k+4s-6} \big]  
    + 7\cdot\big[15\cdot2^{8s-6} - 465\cdot2^{7s-8} +349\cdot2^{6s-7} \big]\quad \text{for}\quad k\geq 2s+2 \label{eq 10.53}
\end{align}
From \eqref{eq 10.51} we see that \eqref{eq 10.53} is valid for k = 2s+1.

\end{proof}
  
    \subsection{ Explicit formulas  for 
    $ \Gamma_{2s+1}^{\left[s\atop{ s\atop s} \right] \times k }\quad  where \quad k \geq 2s+2\quad and \quad  for \quad \Gamma_{2s+2}^{\left[s\atop{ s\atop s} \right] \times (2s+2) } $}
  \label{subsec 6}
  \begin{lem}
  \label{lem 10.5}
  We have \\
 \begin{align}
 \Gamma  _{2s+1}^{\left[s\atop{ s\atop s} \right]\times k} & = 
  105\cdot\big(2^{2k+2s-2}+ 7\cdot2^{k+5s-3} - 31\cdot2^{k+4s-3} \big) \label{eq 10.54} \\
  &  + 105\cdot\big(2^{8s-2} - 31\cdot2^{7s-3} +93\cdot2^{6s-3} \big) \quad \text{for}\quad  k\geq 2s+2.\nonumber \\
    \Gamma  _{2s+2}^{\left[s\atop{ s\atop s} \right]\times (2s+2)} & =  
    2^{9s+3} -7\cdot2^{8s+2} +7\cdot2^{7s+2} - 2^{6s+3}  \label{eq 10.55} 
\end{align}
 \end{lem}
  \begin{proof}
   We get respectively from \eqref{eq 10.1} and \eqref{eq 10.2} for  k = 2s+2 \\
    \begin{align}
  \sum_{ i = 0}^{2s+2}  \Gamma  _{i}^{\left[s\atop{ s\atop s} \right]\times (2s+2)} & = 2^{9s+3} \label{eq 10.56}\\
 &  \text{and}\nonumber \\
  \sum_{i = 0}^{2s+2}  \Gamma  _{i}^{\left[s\atop{ s\atop s} \right]\times (2s+2)}\cdot2^{-i} & =  2^{7s+1} + 2^{6s+3} - 2^{4s+1}. \label{eq 10.57}
\end{align}
  We have  from \eqref{eq 10.3},  \eqref{eq 10.4} with  k= 2s+2 \\
 \begin{align}
   \sum_{ i = 0}^{s-1}  \Gamma  _{i}^{\left[s\atop{ s\atop s} \right]\times (2s+2)}& = 7\cdot2^{4s-6}-3\cdot2^{3s-5} \label{eq 10.58}\\
    \sum_{i = 0}^{s-1}  \Gamma  _{i}^{\left[s\atop{ s\atop s} \right]\times (2s+2)}\cdot2^{-i} & = 15\cdot2^{3s-6}-7\cdot2^{2s-5} \label{eq 10.59}
 \end{align}
  By \eqref{eq 10.5} with k=2s+2 \\
\begin{align}
  \Gamma  _{s}^{\left[s\atop{ s\atop s} \right]\times (2s+2)}= 427\cdot2^{3s-5} -7\cdot2^{2s} +105\cdot2^{4s-6}  \label{eq 10.60} \\
  & \nonumber \\
    \Gamma  _{s}^{\left[s\atop{ s\atop s} \right]\times (2s+2)}\cdot2^{-s}= 427\cdot2^{2s-5} -7\cdot2^{s} +105\cdot2^{3s-6} \label{eq 10.61} 
\end{align}
We obtain from \eqref{eq 10.31} with $j\rightarrow q$ and $k\geq 2s$\\
\begin{align}
& \sum_{q=1}^{s-1} \Gamma  _{s+q}^{\left[s\atop{ s\atop s} \right]\times k} \label{eq 10.62} \\
& = \sum_{q=1}^{s-1}2^{4q}\big[735\cdot2^{k+s-7} +105\cdot2^{4s-6} + 735\cdot2^{2s-7}\big] \nonumber  \\
&  + \sum_{q=1}^{s-1}2^{3q}\big[-147\cdot2^{k+s-6} - 21\cdot2^{3s-5}\big]    +  \sum_{q=1}^{s-1}2^{5q}\big[-3255\cdot2^{2s-8}\big] \nonumber  \\
& = 49\cdot2^{k+5s-7} - 7\cdot2^{k+s-1} -21\cdot2^{k+4s-6}+7\cdot2^{8s-6} +37\cdot2^{6s-7}-7\cdot2^{4s-2}
 -105\cdot2^{7s-8} +3\cdot2^{3s-2} +7\cdot2^{2s} \nonumber
 \end{align}
\begin{align}
& \sum_{q=1}^{s-1} \Gamma  _{s+q}^{\left[s\atop{ s\atop s} \right]\times k}\cdot2^{-(s+q)} \label{eq 10.63} \\
& = \sum_{q=1}^{s-1}2^{3q}\big[735\cdot2^{k-7} +105\cdot2^{3s-6} + 735\cdot2^{s-7}\big] \nonumber  \\
&  + \sum_{q=1}^{s-1}2^{2q}\big[-147\cdot2^{k-6} - 21\cdot2^{2s-5}\big]    +  \sum_{q=1}^{s-1}2^{4q}\big[-3255\cdot2^{s-8}\big] \nonumber  \\
& = 105\cdot2^{k+3s-7} - 49\cdot2^{k+2s-6} -7\cdot2^{k-1} -217\cdot2^{5s-8} +15\cdot2^{6s-6} +77\cdot2^{4s-7}
 -15\cdot2^{3s-3} +7\cdot2^{2s-3} +7\cdot2^{s} \nonumber
 \end{align}
 From \eqref{eq 10.62}  and  \eqref{eq 10.63} with k = 2k+2 we obtain \\
 \begin{align}
  & \sum_{q=1}^{s-1} \Gamma  _{s+q}^{\left[s\atop{ s\atop s} \right]\times (2s+2)} \label{eq 10.64} \\
 & =  287\cdot2^{7s-8} - 53\cdot2^{3s-2} -131\cdot2^{6s-7}  +7\cdot2^{8s-6} -7\cdot2^{4s-2} +7\cdot2^{2s} \nonumber \\
  & \sum_{q=1}^{s-1} \Gamma  _{s+q}^{\left[s\atop{ s\atop s} \right]\times (2s+2)}\cdot2^{-(s+q)} \label{eq 10.65} \\
 & =  623\cdot2^{5s-8} +15\cdot2^{6s-6} - 315\cdot2^{4s-7}  -105\cdot2^{2s-3} -15\cdot2^{3s-3} +7\cdot2^{s} \nonumber 
  \end{align}
We get from \eqref{eq 10.48} with k = 2s+2 \\
\begin{align}
& \Gamma  _{2s}^{\left[s\atop{ s\atop s} \right]\times (2s+2)} = 
 105 \cdot2^{8s-6} +2625\cdot2^{7s-8} -525\cdot2^{6s-7}  \label{eq 10.66} \\
& \Gamma  _{2s}^{\left[s\atop{ s\atop s} \right]\times (2s+2)}\cdot2^{-2s} = 
 105 \cdot2^{6s-6} +2625\cdot2^{5s-8} -525\cdot2^{4s-7}  \label{eq 10.67} 
 \end{align}
  We obtain from  \eqref{eq 10.56}, \eqref{eq 10.58}, \eqref{eq 10.60},  \eqref{eq 10.64}  and  \eqref{eq 10.66}.
\begin{align*}
&   \sum_{ i = 0}^{2s+2}  \Gamma  _{i}^{\left[s\atop{ s\atop s} \right]\times (2s+2)}=
  \sum_{ i = 0}^{s -1}  \Gamma  _{i}^{\left[s\atop{ s\atop s} \right]\times (2s+2)}  
  + \Gamma  _{s}^{\left[s\atop{ s\atop s} \right]\times (2s+2)} + \sum_{q=1}^{s-1}\Gamma  _{s+q}^{\left[s\atop{ s\atop s} \right]\times (2s+2)} + \Gamma _{2s}^{\left[s\atop{ s\atop s} \right]\times (2s+2)}
  + \Gamma  _{2s+1}^{\left[s\atop{ s\atop s} \right]\times (2s+2)} + \Gamma  _{2s+2}^{\left[s\atop{ s\atop s} \right]\times (2s+2)} \\
   & = \big( 7\cdot2^{4s-6}-3\cdot2^{3s-5}\big) +  \big( 427\cdot2^{3s-5} -7\cdot2^{2s} +105\cdot2^{4s-6} \big) \\
 & +  \big( 287\cdot2^{7s-8} - 53\cdot2^{3s-2} -131\cdot2^{6s-7}  +7\cdot2^{8s-6} -7\cdot2^{4s-2} +7\cdot2^{2s} \big) \\
 & + \big( 105 \cdot2^{8s-6} +2625\cdot2^{7s-8} -525\cdot2^{6s-7}\big) 
   + \Gamma  _{2s+1}^{\left[s\atop{ s\atop s} \right]\times (2s+2)} + \Gamma  _{2s+2}^{\left[s\atop{ s\atop s} \right]\times (2s+2)} = 2^{9s+3} \\
    & \Updownarrow 
\end{align*}   
 \begin{align}
     \Gamma  _{2s+1}^{\left[s\atop{ s\atop s} \right]\times (2s+2)} + \Gamma  _{2s+2}^{\left[s\atop{ s\atop s} \right]\times (2s+2)} 
  =  2^{9s+3} -7\cdot2^{8s-2}-91\cdot2^{7s-3} +41\cdot2^{6s-3}\label{eq 10.68}
\end{align}
We get in the same way from  \eqref{eq 10.57}, \eqref{eq 10.59},   \eqref{eq 10.61},   \eqref{eq 10.65} and \eqref{eq 10.67}
\begin{align*}
&   \sum_{ i = 0}^{2s+2}  \Gamma  _{i}^{\left[s\atop{ s\atop s} \right]\times (2s+2)}\cdot2^{-i}=
  \sum_{ i = 0}^{s -1}  \Gamma  _{i}^{\left[s\atop{ s\atop s} \right]\times (2s+2)}\cdot2^{-i}  
  + \Gamma  _{s}^{\left[s\atop{ s\atop s} \right]\times (2s+2)}\cdot2^{-s}   + \sum_{q=1}^{s-1}\Gamma  _{s+q}^{\left[s\atop{ s\atop s} \right]\times (2s+2)}\cdot2^{-(s+q)} \\
  &  + \Gamma _{2s}^{\left[s\atop{ s\atop s} \right]\times (2s+2)}\cdot2^{-2s}  
  + \Gamma  _{2s+1}^{\left[s\atop{ s\atop s} \right]\times (2s+2)}\cdot2^{-(2s+1)}   + \Gamma  _{2s+2}^{\left[s\atop{ s\atop s} \right]\times (2s+2)}\cdot2^{-(2s+2)}   \\
   & = \big(  15\cdot2^{3s-6}-7\cdot2^{2s-5}  \big) +  \big(427\cdot2^{2s-5} -7\cdot2^{s} +105\cdot2^{3s-6}  \big) \\
    & +  \big( 623\cdot2^{5s-8} +15\cdot2^{6s-6} - 315\cdot2^{4s-7}  -105\cdot2^{2s-3} -15\cdot2^{3s-3} +7\cdot2^{s} \big) \\
  & + \big( 105 \cdot2^{6s-6} +2625\cdot2^{5s-8} -525\cdot2^{4s-7}  \big) 
   + \Gamma  _{2s+1}^{\left[s\atop{ s\atop s} \right]\times (2s+2)}\cdot2^{-(2s+1)}  + \Gamma  _{2s+2}^{\left[s\atop{ s\atop s} \right]\times (2s+2)}\cdot2^{-(2s+2)}  =  2^{7s+1} + 2^{6s+3} - 2^{4s+1} \\
    & \Updownarrow 
\end{align*}   
\begin{align}
 2\cdot \Gamma  _{2s+1}^{\left[s\atop{ s\atop s} \right]\times (2s+2)} + \Gamma  _{2s+2}^{\left[s\atop{ s\atop s} \right]\times (2s+2)} 
  =  2^{9s+3}  +49\cdot2^{8s-1} - 203\cdot2^{7s-2} +73\cdot2^{6s-2}\label{eq 10.69}
\end{align}
Hence by \eqref{eq 10.68}, \eqref{eq 10.69}
\begin{align}
& \Gamma  _{2s+1}^{\left[s\atop{ s\atop s} \right]\times (2s+2)}= 105\cdot\big(2^{8s-2} - 3\cdot2^{7s-3} + 2^{6s-3}  \big) \label{eq 10.70}\\
&  \Gamma  _{2s+2}^{\left[s\atop{ s\atop s} \right]\times (2s+2)} = 2^{9s+3} -7\cdot2^{8s+2} +7\cdot2^{7s+2} -2^{6s+3}\label{eq 10.71}
\end{align}
We have from \eqref{eq 9.7} with j = 1 and $k \geq 2s+2$  \\
  \begin{equation}
\label{eq 10.72}
 \Delta _{k}\Gamma  _{2s+1}^{\left[s\atop{ s\atop s} \right]}
  =   \Gamma  _{2s+1}^{\left[s\atop{ s\atop s} \right]\times (k+1)} -  \Gamma  _{2s+1}^{\left[s\atop{ s\atop s} \right]\times k} =105\cdot\big(3\cdot2^{2k+2s-2} +7\cdot2^{k+5s-3} -31\cdot2^{k+4s-3}\big) 
  \end{equation}
   We obtain from \eqref{eq 10.72} for $k\geq 2s+2$ and  \eqref{eq 10.70} 
  \begin{align*}
& \sum_{i=2s+2}^{k}\big( \Gamma  _{2s+1}^{\left[s\atop{ s\atop s} \right]\times (i+1)} -  \Gamma  _{2s+1}^{\left[s\atop{ s\atop s} \right]\times i}\big)
= \sum_{i= 2s+2}^{k} 105\cdot\big(3\cdot2^{2i+2s-2} +7\cdot2^{i+5s-3} -31\cdot2^{i+4s-3}\big)   \\
& \Updownarrow   \\
& \Gamma  _{2s+1}^{\left[s\atop{ s\atop s} \right]\times (k+1)} -\Gamma  _{2s+1}^{\left[s\atop{ s\atop s} \right]\times (2s+2)}
= 315\cdot2^{2s-2}\cdot \sum_{i= 2s+2}^{k}2^{2i} + \big(735\cdot2^{5s-3} -3255\cdot2^{4s-3} \big)\cdot \sum_{i= 2s+2}^{k}2^{i} \\
 & =  105\cdot2^{2k+2s}+ 735\cdot2^{k+5s-2} - 3255\cdot2^{k+4s-2}   - 735\cdot2^{7s-1} +2415\cdot2^{6s-1} 
\end{align*}
\begin{align}
& \Updownarrow  \nonumber  \\
&  \Gamma  _{2s+1}^{\left[s\atop{ s\atop s} \right]\times k} = 
  105\cdot2^{2k+2s-2}+ 735\cdot2^{k+5s-3} - 3255\cdot2^{k+4s-3}   - 735\cdot2^{7s-1} +2415\cdot2^{6s-1} 
   +  105\cdot\big(2^{8s-2} - 3\cdot2^{7s-3} + 2^{6s-3}  \big) \nonumber \\
  & \Updownarrow  \nonumber  \\
  & \Gamma  _{2s+1}^{\left[s\atop{ s\atop s} \right]\times k} =  105\cdot\big(2^{2k+2s-2}+ 7\cdot2^{k+5s-3} - 31\cdot2^{k+4s-3} \big)
    + 105\cdot\big(2^{8s-2} - 31\cdot2^{7s-3} +93\cdot2^{6s-3} \big)  \quad \text{for}\quad k\geq 2s+3 \label{eq 10.73}
\end{align}
From \eqref{eq 10.70} we see that \eqref{eq 10.73} is valid for k = 2s+2.

 \end{proof}
  
     \subsection{ Explicit formulas  for 
    $ \Gamma_{2s+2}^{\left[s\atop{ s\atop s} \right] \times k }\quad  where \quad k \geq 2s+3\quad and \quad  for \quad \Gamma_{2s+3}^{\left[s\atop{ s\atop s} \right] \times (2s+3) } $}
  \label{subsec 7}
  \begin{lem}
  \label{lem 10.6}
  We have \\
 \begin{align}
 \Gamma  _{2s+2}^{\left[s\atop{ s\atop s} \right]\times k} & = 
  105\cdot\big(2^{2k+2s +2}+ 7\cdot2^{k+5s+1} - 31\cdot2^{k+4s+2} \big) \label{eq 10.74} \\
  &  + 105\cdot\big(2^{8s+2} - 31\cdot2^{7s+2} +93\cdot2^{6s+3} \big) \quad \text{for}\quad  k\geq 2s+3.\nonumber \\
    \Gamma  _{2s+3}^{\left[s\atop{ s\atop s} \right]\times (2s+3)} & =  
    2^{9s+6} -7\cdot2^{8s+6} +7\cdot2^{7s+7} - 2^{6s+9}  \label{eq 10.75} \\
    &  \Gamma  _{2s+2}^{\left[s\atop{ s\atop s} \right]\times k}= 16^{3}\cdot \Gamma  _{2(s-1)+1}^{\left[s-1\atop{ s-1\atop s-1} \right]\times (k-3)}\quad \text{for}\quad  k\geq 2s+2.\label{eq 10.76} 
\end{align}
 \end{lem}
  \begin{proof}
  We get respectively from \eqref{eq 10.1} and \eqref{eq 10.2} for  k = 2s+3 \\
    \begin{align}
  \sum_{ i = 0}^{2s+3}  \Gamma  _{i}^{\left[s\atop{ s\atop s} \right]\times (2s+3)} & = 2^{9s+6} \label{eq 10.77}\\
 &  \text{and}\nonumber \\
  \sum_{i = 0}^{2s+3}  \Gamma  _{i}^{\left[s\atop{ s\atop s} \right]\times (2s+3)}\cdot2^{-i} & =  2^{7s+3} + 2^{6s+6} - 2^{4s+3}. \label{eq 10.78}
\end{align}
  We have  from \eqref{eq 10.3},  \eqref{eq 10.4} with  k= 2s+3 \\
 \begin{align}
   \sum_{ i = 0}^{s-1}  \Gamma  _{i}^{\left[s\atop{ s\atop s} \right]\times (2s+3)}& = 7\cdot2^{4s-6}-3\cdot2^{3s-5} \label{eq 10.79}\\
    \sum_{i = 0}^{s-1}  \Gamma  _{i}^{\left[s\atop{ s\atop s} \right]\times (2s+3)}\cdot2^{-i} & = 15\cdot2^{3s-6}-7\cdot2^{2s-5} \label{eq 10.80}
 \end{align}
   By \eqref{eq 10.5} with k=2s+3 \\
\begin{align}
  \Gamma  _{s}^{\left[s\atop{ s\atop s} \right]\times (2s+3)}= 875\cdot2^{3s-5} -7\cdot2^{2s} +105\cdot2^{4s-6}  \label{eq 10.81} \\
  & \nonumber \\
    \Gamma  _{s}^{\left[s\atop{ s\atop s} \right]\times (2s+3)}\cdot2^{-s}= 875\cdot2^{2s-5} -7\cdot2^{s} +105\cdot2^{3s-6} \label{eq 10.82} 
\end{align}
   From \eqref{eq 10.62}  and  \eqref{eq 10.63} with k = 2k+3 we obtain \\
 \begin{align}
  & \sum_{q=1}^{s-1} \Gamma  _{s+q}^{\left[s\atop{ s\atop s} \right]\times (2s+3)} \label{eq 10.83} \\
& =   7\cdot2^{8s-6} + 679\cdot2^{7s-8} -299\cdot2^{6s-7}-7\cdot2^{4s-2}- 109\cdot2^{3s-2}+7\cdot2^{2s} \nonumber \\
   & \sum_{q=1}^{s-1} \Gamma  _{s+q}^{\left[s\atop{ s\atop s} \right]\times (2s+3)}\cdot2^{-(s+q)} \label{eq 10.84} \\
 & =  1463\cdot2^{5s-8} +15\cdot2^{6s-6} - 707\cdot2^{4s-7}  -15\cdot2^{3s-3} -217\cdot2^{2s-3}  +7\cdot2^{s} \nonumber 
  \end{align}
 We get from \eqref{eq 10.48} with k = 2s+3 \\
\begin{align}
& \Gamma  _{2s}^{\left[s\atop{ s\atop s} \right]\times (2s+3)} = 
 105 \cdot2^{8s-6} +8505\cdot2^{7s-8} +3675\cdot2^{6s-7}  \label{eq 10.85} \\
& \Gamma  _{2s}^{\left[s\atop{ s\atop s} \right]\times (2s+3)}\cdot2^{-2s} = 
 105 \cdot2^{6s-6} +8505\cdot2^{5s-8} +3675\cdot2^{4s-7}  \label{eq 10.86} 
 \end{align}
  We get from \eqref{eq 10.54} with k = 2s+3 \\
\begin{align}
& \Gamma  _{2s+1}^{\left[s\atop{ s\atop s} \right]\times (2s+3)} = 
 105 \cdot2^{8s-2} +2625\cdot2^{7s-3} -2835\cdot2^{6s-3}  \label{eq 10.87} \\
& \Gamma  _{2s+1}^{\left[s\atop{ s\atop s} \right]\times (2s+3)}\cdot2^{-(2s+1)} = 
 105 \cdot2^{6s-3} +2625\cdot2^{5s-4} -2835\cdot2^{4s-4}  \label{eq 10.88} 
 \end{align}
    We obtain from  \eqref{eq 10.77}, \eqref{eq 10.79}, \eqref{eq 10.81},  \eqref{eq 10.83}, \eqref{eq 10.85}  and  \eqref{eq 10.87}.
\begin{align*}
&   \sum_{ i = 0}^{2s+3}  \Gamma  _{i}^{\left[s\atop{ s\atop s} \right]\times (2s+3)}=
  \sum_{ i = 0}^{s -1}  \Gamma  _{i}^{\left[s\atop{ s\atop s} \right]\times (2s+3)}  
  + \Gamma  _{s}^{\left[s\atop{ s\atop s} \right]\times (2s+3)} + \sum_{q=1}^{s-1}\Gamma  _{s+q}^{\left[s\atop{ s\atop s} \right]\times (2s+3)}\\
  &  + \Gamma _{2s}^{\left[s\atop{ s\atop s} \right]\times (2s+3)} + \Gamma  _{2s+1}^{\left[s\atop{ s\atop s} \right]\times (2s+3)}
   + \Gamma  _{2s+2}^{\left[s\atop{ s\atop s} \right]\times (2s+3)} +  \Gamma  _{2s+3}^{\left[s\atop{ s\atop s} \right]\times (2s+3)} \\
  & = \big( 7\cdot2^{4s-6}-3\cdot2^{3s-5}\big) +  \big( 875\cdot2^{3s-5} -7\cdot2^{2s} +105\cdot2^{4s-6} \big) \\
  & +  \big( 7\cdot2^{8s-6} + 679\cdot2^{7s-8} -299\cdot2^{6s-7}-7\cdot2^{4s-2}- 109\cdot2^{3s-2}+7\cdot2^{2s} \big) \\
   & + \big(105 \cdot2^{8s-6} +8505\cdot2^{7s-8} +3675\cdot2^{6s-7}\big) + \big( 105 \cdot2^{8s-2} +2625\cdot2^{7s-3} -2835\cdot2^{6s-3} \big) \\
  &  + \Gamma  _{2s+2}^{\left[s\atop{ s\atop s} \right]\times (2s+3)} + \Gamma  _{2s+3}^{\left[s\atop{ s\atop s} \right]\times (2s+3)} = 2^{9s+6} \\
  & \Updownarrow 
\end{align*}   
 \begin{align}
     \Gamma  _{2s+2}^{\left[s\atop{ s\atop s} \right]\times (2s+3)} + \Gamma  _{2s+3}^{\left[s\atop{ s\atop s} \right]\times (2s+3)} 
  =  2^{9s+6} -7\cdot2^{8s+2}-91\cdot2^{7s+2} +41\cdot2^{6s+3}\label{eq 10.89}
\end{align}

We get in the same way from  \eqref{eq 10.78}, \eqref{eq 10.80},   \eqref{eq 10.82},  \eqref{eq 10.84}, \eqref{eq 10.86} and \eqref{eq 10.88}.
\begin{align*}
&   \sum_{ i = 0}^{2s+3}  \Gamma  _{i}^{\left[s\atop{ s\atop s} \right]\times (2s+3)}\cdot2^{-i}=
  \sum_{ i = 0}^{s -1}  \Gamma  _{i}^{\left[s\atop{ s\atop s} \right]\times (2s+3)}\cdot2^{-i}  
  + \Gamma  _{s}^{\left[s\atop{ s\atop s} \right]\times (2s+3)}\cdot2^{-s}   + \sum_{q=1}^{s-1}\Gamma  _{s+q}^{\left[s\atop{ s\atop s} \right]\times (2s+3)}\cdot2^{-(s+q)} \\
  &  + \Gamma _{2s}^{\left[s\atop{ s\atop s} \right]\times (2s+3)}\cdot2^{-2s}  
  + \Gamma  _{2s+1}^{\left[s\atop{ s\atop s} \right]\times (2s+3)}\cdot2^{-(2s+1)} 
    + \Gamma  _{2s+2}^{\left[s\atop{ s\atop s} \right]\times (2s+3)}\cdot2^{-(2s+2)}  + \Gamma  _{2s+3}^{\left[s\atop{ s\atop s} \right]\times (2s+3)}\cdot2^{-(2s+3)}  \\
  & = \big(  15\cdot2^{3s-6}-7\cdot2^{2s-5}  \big) +  \big(  875\cdot2^{2s-5} -7\cdot2^{s} +105\cdot2^{3s-6}  \big) \\  & +  \big(  1463\cdot2^{5s-8} +15\cdot2^{6s-6} - 707\cdot2^{4s-7}  -15\cdot2^{3s-3} -217\cdot2^{2s-3}  +7\cdot2^{s} \big) \\
  & + \big( 105 \cdot2^{6s-6} +8505\cdot2^{5s-8} +3675\cdot2^{4s-7}  \big) + \big( 105 \cdot2^{6s-3} +2625\cdot2^{5s-4} -2835\cdot2^{4s-4}  \big) 
  + \Gamma  _{2s+2}^{\left[s\atop{ s\atop s} \right]\times (2s+3)}\cdot2^{-(2s+2)}\\
 &   + \Gamma  _{2s+3}^{\left[s\atop{ s\atop s} \right]\times (2s+3)}\cdot2^{-(2s+3)}  =  2^{7s+3} + 2^{6s+6} - 2^{4s+3} \\
  & \Updownarrow 
\end{align*}  
\begin{align}
 2\cdot \Gamma  _{2s+2}^{\left[s\atop{ s\atop s} \right]\times (2s+3)} + \Gamma  _{2s+3}^{\left[s\atop{ s\atop s} \right]\times (2s+3)} 
  =  2^{9s+6}  +49\cdot2^{8s+3} - 203\cdot2^{7s+3} +73\cdot2^{6s+4}\label{eq 10.90}
\end{align}
Hence by \eqref{eq 10.89}, \eqref{eq 10.90}
\begin{align}
& \Gamma  _{2s+2}^{\left[s\atop{ s\atop s} \right]\times (2s+3)}= 105\cdot\big(2^{8s+2} - 3\cdot2^{7s+2} + 2^{6s+3}  \big) \label{eq 10.91}\\
&  \Gamma  _{2s+3}^{\left[s\atop{ s\atop s} \right]\times (2s+3)} = 2^{9s+6} -7\cdot2^{8s+6} +7\cdot2^{7s+7} -2^{6s+9}\label{eq 10.92}
\end{align}
 We have from \eqref{eq 9.7} with j = 2 and $k \geq 2s+3$  \\
  \begin{equation}
\label{eq 10.93}
 \Delta _{k}\Gamma  _{2s+2}^{\left[s\atop{ s\atop s} \right]}
  =   \Gamma  _{2s+2}^{\left[s\atop{ s\atop s} \right]\times (k+1)} -  \Gamma  _{2s+2}^{\left[s\atop{ s\atop s} \right]\times k} =105\cdot\big(3\cdot2^{2k+2s+2} +7\cdot2^{k+5s+1} -31\cdot2^{k+4s+2}\big) 
  \end{equation}
   We obtain from \eqref{eq 10.93} for $k\geq 2s+3$ and  \eqref{eq 10.91} 
  \begin{align*}
& \sum_{i=2s+3}^{k}\big( \Gamma  _{2s+2}^{\left[s\atop{ s\atop s} \right]\times (i+1)} -  \Gamma  _{2s+2}^{\left[s\atop{ s\atop s} \right]\times i}\big)
= \sum_{i= 2s+3}^{k} 105\cdot\big(3\cdot2^{2i+2s+2} +7\cdot2^{i+5s+1} -31\cdot2^{i+4s+2}\big)   \\
& \Updownarrow   \\
& \Gamma  _{2s+2}^{\left[s\atop{ s\atop s} \right]\times (k+1)} -\Gamma  _{2s+2}^{\left[s\atop{ s\atop s} \right]\times (2s+3)}
= 315\cdot2^{2s+2}\cdot \sum_{i= 2s+3}^{k}2^{2i} + \big(735\cdot2^{5s+1} -3255\cdot2^{4s+2} \big)\cdot \sum_{i= 2s+3}^{k}2^{i} \\
 & =  105\cdot2^{2k+2s+4}+ 735\cdot2^{k+5s+2} - 3255\cdot2^{k+4s+3}   - 735\cdot2^{7s+4} +2415\cdot2^{6s+5} 
\end{align*}

 \begin{align}
& \Updownarrow  \nonumber  \\
&  \Gamma  _{2s+2}^{\left[s\atop{ s\atop s} \right]\times (k+1)} = 
  105\cdot2^{2k+2s+4}+ 735\cdot2^{k+5s+2} - 3255\cdot2^{k+4s+3}   - 735\cdot2^{7s+4} +2415\cdot2^{6s+5} 
   +  105\cdot\big(2^{8s+2} - 3\cdot2^{7s+2} + 2^{6s+3}  \big) \nonumber \\
  & \Updownarrow  \nonumber  \\
  & \Gamma  _{2s+2}^{\left[s\atop{ s\atop s} \right]\times k} =  105\cdot\big(2^{2k+2s+2}+ 7\cdot2^{k+5s+1} - 31\cdot2^{k+4s+2} \big)
    + 105\cdot\big(2^{8s+2} - 31\cdot2^{7s+2} +93\cdot2^{6s+3} \big)  \quad \text{for}\quad k\geq 2s+4 \label{eq 10.94}
\end{align}

From \eqref{eq 10.91} we see that \eqref{eq 10.94} is valid for k = 2s+3.
  It remains to prove \eqref{eq 10.76}.\\
We have respectively  from \eqref{eq 10.49} with $s\rightarrow s-1$ and \eqref{eq 10.55} \\
\begin{align*}
 & 16^{3}\cdot \Gamma  _{2(s-1)+1}^{\left[s-1\atop{ s-1\atop s-1} \right]\times (2(s-1)+1)}
 = 2^{12}\cdot\big( 2^{9(s-1)} -7\cdot2^{8(s-1)-2} +7\cdot2^{7(s-1)-3} - 2^{6(s-1)-3} \big)\\
 &  = 2^{9s+3} -7\cdot2^{8s+2} +7\cdot2^{7s+2} - 2^{6s+3}=  \Gamma  _{2s+3}^{\left[s\atop{ s\atop s} \right]\times (2s+3)} 
\end{align*}

We get  respectively  from \eqref{eq 10.54} with $s\rightarrow s-1,\quad k\rightarrow k-3,\quad k\geq 2s+3 $ and \eqref{eq 10.74} \\
\begin{align*}
 & 16^{3}\cdot \Gamma  _{2(s-1)+1}^{\left[s-1\atop{ s-1\atop s-1} \right]\times (k-3) }
 = 2^{12}\cdot105\cdot\big(2^{2(k-3)+2(s-1)-2}+ 7\cdot2^{k-3+5(s-1)-3} - 31\cdot2^{k-3+4(s-1)-3} \big)\\
    &  + 2^{12}\cdot105\cdot\big(2^{8(s-1)-2} - 31\cdot2^{7(s-1)-3} +93\cdot2^{6(s-1)-3} \big)  \\
  &  =   105\cdot\big(2^{2k+2s+2}+ 7\cdot2^{k+5s+1} - 31\cdot2^{k+4s+2} \big) 
     + 105\cdot\big(2^{8s+2} - 31\cdot2^{7s+2} +93\cdot2^{6s+3} \big)
      =  \Gamma  _{2s+2}^{\left[s\atop{ s\atop s} \right]\times k}   
\end{align*}
\end{proof}

    \subsection{Un explicit reduction  formula  for 
    $ \Gamma_{2s+1+j}^{\left[s\atop{ s\atop s} \right] \times k }\quad for \quad 0\leq j\leq s-1,\quad k\geq 2s+1+j $}
  \label{subsec 8} 
   \begin{lem}
  \label{lem 10.7}
  We have for $0\leq j\leq s-2,\quad k\geq 2s+2+j$\\
 \begin{align}
(H_{j}) \quad \Gamma  _{2s+1+j}^{\left[s\atop{ s\atop s} \right]\times k} & = 
 16^{3j}\cdot\Gamma  _{2(s-j)+1}^{\left[s-j\atop{ s-j\atop s-j} \right]\times (k-3j)} \label{eq 10.95}\\
 & = 2^{12j}\cdot105\cdot\big(2^{2(k-3j)+2(s-j) -2}+ 7\cdot2^{k-3j+5(s-j)-3} - 31\cdot2^{k-3j+4(s-j)-3} \big) \nonumber \\
  &  + 2^{12j}\cdot 105\cdot\big(2^{8(s-j)-2} - 31\cdot2^{7(s-j)-3} +93\cdot2^{6(s-j)-3} \big) \nonumber \\
  & =  105\cdot\big(2^{2k+2s+4j -2}+ 7\cdot2^{k+5s+4j-3} - 31\cdot2^{k+4s+5j-3} \big) \nonumber \\
  &  + 105\cdot\big(2^{8s+4j-2} - 31\cdot2^{7s+5j-3} +93\cdot2^{6s+6j-3} \big) \nonumber
 \end{align}
 \end{lem}
  \begin{lem}
  \label{lem 10.8}
  We have for $0\leq j\leq s-1.$\\
 \begin{align}
 \Gamma  _{2s+1+j}^{\left[s\atop{ s\atop s} \right]\times (2s+1+j)} & = 
 16^{3j}\cdot\Gamma  _{2(s-j)+1}^{\left[s-j\atop{ s-j\atop s-j} \right]\times (2(s-j)+1)} \label{eq 10.96}\\
 &  =   2^{12j}\cdot\big(2^{9(s-j)} -7\cdot2^{8(s-j)-2} +7\cdot2^{7(s-j)-3} - 2^{6(s-j)-3}\big)\nonumber \\
 & =  2^{9s+3j} -7\cdot2^{8s+4j-2} +7\cdot2^{7s+5j-3} - 2^{6s+6j-3}\nonumber
 \end{align}
 We have for $k\geq 3s. $\\
  \begin{align}
 \Gamma  _{3s}^{\left[s\atop{ s\atop s} \right]\times k} & = 
 16^{3(s-1)}\cdot\Gamma  _{3}^{\left[1\atop{ 1\atop 1} \right]\times (k-3(s-1))} \label{eq 10.97}\\
 &  =   2^{12s-12}\cdot\big(2^{3k-9s+9} -7\cdot2^{2k-6s+6} +7\cdot2^{k-3s+4} - 2^{3}\big) \nonumber \\
 & =  2^{3k+3s-3} -7\cdot2^{2k+6s-6} +7\cdot2^{k+9s-8} - 2^{12s-9} \nonumber 
 \end{align}
  \end{lem}
   \begin{proof}
   \underline{Proof of \eqref{eq 10.95} and  \eqref{eq 10.96}}
   
 The proof is by strong induction, that is:\\
If
\begin{itemize}
\item \text{$ (H_{0}) $ is true, and  }
\item \text{for all $ j\geq 1, \quad (H_{0})\wedge (H_{1})\wedge\ldots\wedge (H_{j})\quad implies \; (H_{j+1}) $}
\end{itemize}
then $ (H_{j}) $ is true for all $ j \geq 1. $\\

\underline{$(H_{0})$ is true}

Obviously \\

\underline{$(H_{0})  \Rightarrow  (H_{1})$ }
Immediately by the proof of \eqref{eq 10.76}\\

 \underline{$ (H_{0})\wedge (H_{1})\wedge\ldots\wedge (H_{j})\quad implies \; (H_{j+1}) $}\\

We are going to show that $(H_{j+1})$ is true assuming that $(H_{q})$ is true for all $0\leq q\leq j$
   \begin{align*}
 & (H_{j+1}) \quad \Gamma  _{2s+1+(j+1)}^{\left[s\atop{ s\atop s} \right]\times k}  = 
 16^{3(j+1)}\cdot\Gamma  _{2(s-(j+1))+1}^{\left[s-(j+1)\atop{ s-(j+1)\atop s-(j+1)} \right]\times (k-3(j+1))} \\
 & = 2^{12j+12}\cdot105\cdot\big(2^{2(k-3(j+1))+2(s-(j+1)) -2}+ 7\cdot2^{k-3(j+1)+5(s-(j+1))-3} - 31\cdot2^{k-3(j+1)+4(s-(j+1))-3} \big)  \\
  &  + 2^{12j+12}\cdot 105\cdot\big(2^{8(s-(j+1))-2} - 31\cdot2^{7(s-(j+1))-3} +93\cdot2^{6(s-(j+1))-3} \big)  \\
   & =  105\cdot\big(2^{2k+2s+4j +2}+ 7\cdot2^{k+5s+4j+1} - 31\cdot2^{k+4s+5j+2} \big)  \\
  &  + 105\cdot\big(2^{8s+4j+2} - 31\cdot2^{7s+5j+2} +93\cdot2^{6s+6j+3} \big) 
 \end{align*}
 We get respectively from \eqref{eq 10.1} and \eqref{eq 10.2} for  k = 2s+j+3 \\
    \begin{align}
  \sum_{ i = 0}^{2s+j+3}  \Gamma  _{i}^{\left[s\atop{ s\atop s} \right]\times (2s+j+3)} & = 2^{9s+3j+6} \label{eq 10.98}\\
 &  \text{and}\nonumber \\
  \sum_{i = 0}^{2s +j+3}  \Gamma  _{i}^{\left[s\atop{ s\atop s} \right]\times (2s+j+3)}\cdot2^{-i}
   & =  2^{7s+2j+3} + 2^{6s+3j+6} - 2^{4s+2j+3}. \label{eq 10.99}
\end{align}
  We have  from \eqref{eq 10.3},  \eqref{eq 10.4} with  k= 2s+j+3 \\
 \begin{align}
   \sum_{ i = 0}^{s-1}  \Gamma  _{i}^{\left[s\atop{ s\atop s} \right]\times (2s+j+3)}& = 7\cdot2^{4s-6}-3\cdot2^{3s-5} \label{eq 10.100}\\
    \sum_{i = 0}^{s-1}  \Gamma  _{i}^{\left[s\atop{ s\atop s} \right]\times (2s+j+3)}\cdot2^{-i} & = 15\cdot2^{3s-6}-7\cdot2^{2s-5} \label{eq 10.101}
 \end{align}
 By \eqref{eq 10.5} with k=2s+j+3 \\
\begin{align}
  \Gamma  _{s}^{\left[s\atop{ s\atop s} \right]\times (2s+j+3)}= 7\cdot2^{3s+j+2} -7\cdot2^{2s} +105\cdot2^{4s-6}-21\cdot2^{3s-5}  \label{eq 10.102} \\
  & \nonumber \\
    \Gamma  _{s}^{\left[s\atop{ s\atop s} \right]\times (2s+j+3)}\cdot2^{-s}= 7\cdot2^{2s+j+2} -7\cdot2^{s} +105\cdot2^{3s-6}-21\cdot2^{2s-5} \label{eq 10.103} 
\end{align}
  From \eqref{eq 10.62}  and  \eqref{eq 10.63} with k = 2k+j+3 we obtain \\
 \begin{align}
  & \sum_{q=1}^{s-1} \Gamma  _{s+q}^{\left[s\atop{ s\atop s} \right]\times (2s+j+3)} \label{eq 10.104} \\
  & = 49\cdot2^{7s+j-4} -7\cdot2^{3s+j+2} -21\cdot2^{6s+j-3}+7\cdot2^{8s-6} - 105\cdot2^{7s-8} +37\cdot2^{6s-7}
    -7\cdot2^{4s-2} + 3\cdot2^{3s-2} + 7\cdot2^{2s} \nonumber \\
  & \sum_{q=1}^{s-1} \Gamma  _{s+q}^{\left[s\atop{ s\atop s} \right]\times (2s+j+3)}\cdot2^{-(s+q)} \label{eq 10.105} \\
   & = 105\cdot2^{5s+j-4} -49\cdot2^{4s+j-3} -7\cdot2^{2s+j+2}  +15\cdot2^{6s-6} -217\cdot2^{5s-8}
    +77\cdot2^{4s-7} -15\cdot2^{3s-3} + 7\cdot2^{2s-3}+ 7\cdot2^{s} \nonumber 
    \end{align}
We get from \eqref{eq 10.48} with k = 2s+j+3 \\
\begin{align}
& \Gamma  _{2s}^{\left[s\atop{ s\atop s} \right]\times (2s+j+3)} = 
7\cdot2^{6s+2j+4} +735\cdot2^{7s+j-4} -819\cdot2^{6s+j-3}
+ 105 \cdot2^{8s-6} -3255\cdot2^{7s-8} +2443\cdot2^{6s-7}  \label{eq 10.106} \\
& \Gamma  _{2s}^{\left[s\atop{ s\atop s} \right]\times (2s+j+3)}\cdot2^{-2s} = 
  7\cdot2^{4s+2j+4} +735\cdot2^{5s+j-4} -819\cdot2^{4s+j-3}
+ 105 \cdot2^{6s-6} -3255\cdot2^{5s-8} +2443\cdot2^{4s-7} \label{eq 10.107} 
 \end{align}
From the induction hypothesis $(H_{q})$ for $0\leq q\leq j$ we get by \eqref{eq 10.95}
with  k = 2s+j+3.
\begin{align}
& \sum_{q=0}^{j}\Gamma  _{2s+1+q}^{\left[s\atop{ s\atop s} \right]\times (2s+j+3)} \label{eq 10.108}\\
 & =  \sum_{q=0}^{j} 105\cdot\big(2^{2(2s+j+3)+2s+4q -2}+ 7\cdot2^{(2s+j+3)+5s+4q-3} - 31\cdot2^{(2s+j+3)+4s+5q-3} \big) \nonumber \\
  &  +  \sum_{q=0}^{j}105\cdot\big(2^{8s+4q-2} - 31\cdot2^{7s+5q-3} +93\cdot2^{6s+6q-3} \big) \nonumber \\
  & = 91\cdot2^{7s+5j+2} - 41\cdot2^{6s+6j+3} + 7\cdot2^{8s+4j+2} - 7\cdot2^{6s+2j+4} \nonumber\\
   &   - 49\cdot2^{7s+j} +  105\cdot2^{6s+j}  -7\cdot2^{8s-2} +105\cdot2^{7s-3}- 155\cdot2^{6s-3} \nonumber
\end{align}
\begin{align}
& \sum_{q=0}^{j}\Gamma  _{2s+1+q}^{\left[s\atop{ s\atop s} \right]\times (2s+j+3)}\cdot2^{-(2s+1+q)} \label{eq 10.109}\\
 & =  \sum_{q=0}^{j} 105\cdot\big(2^{2(2s+j+3)+2s+4q -2-(2s+1+q)}+ 7\cdot2^{(2s+j+3)+5s+4q-3-(2s+1+q)} - 31\cdot2^{(2s+j+3)+4s+5q-3-(2s+1+q)} \big) \nonumber \\
  &  +  \sum_{q=0}^{j}105\cdot\big(2^{8s+4q-2-(2s+1+q)} - 31\cdot2^{7s+5q-3-(2s+1+q)} +93\cdot2^{6s+6q-3-(2s+1+q)} \big) \nonumber \\
  & =  \sum_{q=0}^{j} 105\cdot\big(2^{4s+2j+3+3q}+ 7\cdot2^{5s+j+3q-1} - 31\cdot2^{4s+j+4q-1} \big) 
    +  \sum_{q=0}^{j}105\cdot\big(2^{6s+3q-3} - 31\cdot2^{5s+4q-4} +93\cdot2^{4s+5q-4} \big) \nonumber \\
    & = 203\cdot2^{5s+4j} -73\cdot2^{4s+5j+1} -15\cdot2^{4s+2j+3}+ 15\cdot2^{6s+3j} -15\cdot2^{6s-3} \nonumber\\
      &  -105\cdot2^{5s+j-1} +217\cdot2^{4s+j-1} +217\cdot2^{5s-4} -315\cdot2^{4s-4}\nonumber
  \end{align}
 We obtain from  \eqref{eq 10.98}, \eqref{eq 10.100}, \eqref{eq 10.102},  \eqref{eq 10.104},  \eqref{eq 10.106}  and  \eqref{eq 10.108}.
\begin{align*}
&   \sum_{ i = 0}^{2s+j+3}  \Gamma  _{i}^{\left[s\atop{ s\atop s} \right]\times (2s+j+3)}=
  \sum_{ i = 0}^{s -1}  \Gamma  _{i}^{\left[s\atop{ s\atop s} \right]\times (2s+j+3)}  
  + \Gamma  _{s}^{\left[s\atop{ s\atop s} \right]\times (2s+j+3)} + \sum_{q=1}^{s-1}\Gamma  _{s+q}^{\left[s\atop{ s\atop s} \right]\times (2s+j+3)}\\
  &  + \Gamma _{2s}^{\left[s\atop{ s\atop s} \right]\times (2s+j+3)}
  +  \sum_{ q = 0}^{j}\Gamma  _{2s+1+q}^{\left[s\atop{ s\atop s} \right]\times (2s+j+3)} + \Gamma  _{2s+j+2}^{\left[s\atop{ s\atop s} \right]\times (2s+j+3)} + \Gamma  _{2s+j+3}^{\left[s\atop{ s\atop s} \right]\times (2s+j+3)} \\
    & = \big( 7\cdot2^{4s-6}-3\cdot2^{3s-5}\big) +  \big( 7\cdot2^{3s+j+2} -7\cdot2^{2s} +105\cdot2^{4s-6}-21\cdot2^{3s-5} \big) \\
    & +  \big( 287\cdot2^{7s-8} - 53\cdot2^{3s-2} -131\cdot2^{6s-7}  +7\cdot2^{8s-6} -7\cdot2^{4s-2} +7\cdot2^{2s} \big) \\
  & + \big(  49\cdot2^{7s+j-4} -7\cdot2^{3s+j+2} -21\cdot2^{6s+j-3}+7\cdot2^{8s-6} - 105\cdot2^{7s-8} +37\cdot2^{6s-7}
    -7\cdot2^{4s-2} + 3\cdot2^{3s-2} + 7\cdot2^{2s} \big) \\
     & +  \big( 7\cdot2^{6s+2j+4} +735\cdot2^{7s+j-4} -819\cdot2^{6s+j-3}
+ 105 \cdot2^{8s-6} -3255\cdot2^{7s-8} +2443\cdot2^{6s-7}  \big) \\
&    + \big(  91\cdot2^{7s+5j+2} - 41\cdot2^{6s+6j+3} + 7\cdot2^{8s+4j+2} - 7\cdot2^{6s+2j+4} \\
&  - 49\cdot2^{7s+j} +  105\cdot2^{6s+j}  -7\cdot2^{8s-2} +105\cdot2^{7s-3}- 155\cdot2^{6s-3} \big)
 + \Gamma  _{2s+j+2}^{\left[s\atop{ s\atop s} \right]\times (2s+j+3)} + \Gamma  _{2s+j+3}^{\left[s\atop{ s\atop s} \right]\times (2s+j+3)} = 2^{9s+3j+6} \\
    & \Updownarrow 
\end{align*}   
 \begin{align}
     \Gamma  _{2s+j+2}^{\left[s\atop{ s\atop s} \right]\times (2s+j+3)} + \Gamma  _{2s+j+3}^{\left[s\atop{ s\atop s} \right]\times (2s+j+3)} 
  =  2^{9s+3j+6} -7\cdot2^{8s+4j+2}-91\cdot2^{7s+5j+2} +41\cdot2^{6s+6j+3}\label{eq 10.110}
\end{align}
 We get in the same way  from  \eqref{eq 10.99}, \eqref{eq 10.101}, \eqref{eq 10.103},  \eqref{eq 10.105},  \eqref{eq 10.107}  and  \eqref{eq 10.109}.
\begin{align*}
&   \sum_{ i = 0}^{2s+j+3}  \Gamma  _{i}^{\left[s\atop{ s\atop s} \right]\times (2s+j+3)}\cdot2^{-i} =
  \sum_{ i = 0}^{s -1}  \Gamma  _{i}^{\left[s\atop{ s\atop s} \right]\times (2s+j+3)}\cdot2^{-i}  
  + \Gamma  _{s}^{\left[s\atop{ s\atop s} \right]\times (2s+j+3)}\cdot2^{-s} + \sum_{q=1}^{s-1}\Gamma  _{s+q}^{\left[s\atop{ s\atop s} \right]\times (2s+j+3)}\cdot2^{-(s+q)}\\
  &  + \Gamma _{2s}^{\left[s\atop{ s\atop s} \right]\times (2s+j+3)}\cdot2^{-2s}
  +  \sum_{ q = 0}^{j}\Gamma  _{2s+1+q}^{\left[s\atop{ s\atop s} \right]\times (2s+j+3)}\cdot2^{-(2s+1+q)} + \Gamma  _{2s+j+2}^{\left[s\atop{ s\atop s} \right]\times (2s+j+3)}\cdot2^{-(2s+j+2)} + \Gamma  _{2s+j+3}^{\left[s\atop{ s\atop s} \right]\times (2s+j+3)}\cdot2^{-(2s+j+3)} \\
  & = \big( 15\cdot2^{3s-6}-7\cdot2^{2s-5} \big) +  \big( 7\cdot2^{2s+j+2} -7\cdot2^{s} +105\cdot2^{3s-6}-21\cdot2^{2s-5}  \big) \\
 & +  \big( 105\cdot2^{5s+j-4} -49\cdot2^{4s+j-3} -7\cdot2^{2s+j+2}  +15\cdot2^{6s-6} -217\cdot2^{5s-8}
    +77\cdot2^{4s-7} -15\cdot2^{3s-3} + 7\cdot2^{2s-3}+ 7\cdot2^{s}  \big) \\
 & + \big(7\cdot2^{4s+2j+4} +735\cdot2^{5s+j-4} -819\cdot2^{4s+j-3}
+ 105 \cdot2^{6s-6} -3255\cdot2^{5s-8} +2443\cdot2^{4s-7}   \big) \\
 &    + \big( 203\cdot2^{5s+4j} -73\cdot2^{4s+5j+1} -15\cdot2^{4s+2j+3}+ 15\cdot2^{6s+3j} -15\cdot2^{6s-3} \\
  &  -105\cdot2^{5s+j-1} +217\cdot2^{4s+j-1} +217\cdot2^{5s-4} -315\cdot2^{4s-4} \big) \\
& + \Gamma  _{2s+j+2}^{\left[s\atop{ s\atop s} \right]\times (2s+j+3)}\cdot2^{-(2s++j+2)} + \Gamma  _{2s+j+3}^{\left[s\atop{ s\atop s} \right]\times (2s+j+3)}\cdot2^{-(2s+j+3)} = 2^{7s+2j+3}+ 2^{6s+3j+6} - 2^{4s+2j+3} \\
    & \Updownarrow 
\end{align*}   
 \begin{align}
  &   \Gamma  _{2s+j+2}^{\left[s\atop{ s\atop s} \right]\times (2s+j+3)}\cdot2^{-(2s+j+2)} + \Gamma  _{2s+j+3}^{\left[s\atop{ s\atop s} \right]\times (2s+j+3)}\cdot2^{-(2s+j+3)} 
  =  2^{7s+2j+3} +73\cdot2^{4s+5j+1}- 203\cdot2^{5s+4j} +49\cdot2^{6s+3j}\nonumber \\
  & \Updownarrow  \nonumber \\
  &   2\cdot\Gamma  _{2s+j+2}^{\left[s\atop{ s\atop s} \right]\times (2s+j+3)} + \Gamma  _{2s+j+3}^{\left[s\atop{ s\atop s} \right]\times (2s+j+3)} 
  =  2^{9s+3j+6} +73\cdot2^{6s+6j+4}- 203\cdot2^{7s+5j+3} +49\cdot2^{8s+4j+3}\label{eq 10.111} 
\end{align}
Hence by \eqref{eq 10.110}, \eqref{eq 10.111}
\begin{align}
& \Gamma  _{2s+j+2}^{\left[s\atop{ s\atop s} \right]\times (2s+j+3)}= 105\cdot\big(2^{6s+6j+3} - 3\cdot2^{7s+5j+2} + 2^{8s+4j+2}  \big) \label{eq 10.112}\\
&  \Gamma  _{2s+j+3}^{\left[s\atop{ s\atop s} \right]\times (2s+j+3)} = 2^{9s+3j+6} -7\cdot2^{8s+4j+6} +7\cdot2^{7s+5j+7} -2^{6s+6j+9}\label{eq 10.113}
\end{align}
 We have from \eqref{eq 9.7} with $ j \rightarrow  j+2,\quad j+2\leq s-1 $ and $k \geq 2s+j+3$  \\
  \begin{equation}
\label{eq 10.114}
 \Delta _{k}\Gamma  _{2s+j+2}^{\left[s\atop{ s\atop s} \right]}
  =   \Gamma  _{2s+j+2}^{\left[s\atop{ s\atop s} \right]\times (k+1)} -  \Gamma  _{2s+j+2}^{\left[s\atop{ s\atop s} \right]\times k} =105\cdot\big(3\cdot2^{2k+2s+4j+2} +7\cdot2^{k+5s+4j+1} -31\cdot2^{k+4s+5j+2}\big) 
  \end{equation}
    We obtain from \eqref{eq 10.114} for $k\geq 2s+j+3$ and  \eqref{eq 10.112} 
  \begin{align*}
& \sum_{i=2s+j+3}^{k}\big( \Gamma  _{2s+j+2}^{\left[s\atop{ s\atop s} \right]\times (i+1)} -  \Gamma  _{2s+j+2}^{\left[s\atop{ s\atop s} \right]\times i}\big)
= \sum_{i= 2s+j+3}^{k} 105\cdot\big(3\cdot2^{2i+2s+4j+2} +7\cdot2^{i+5s+4j+1} -31\cdot2^{i+4s+5j+2}\big)   \\
& \Updownarrow   \\
& \Gamma  _{2s+j+2}^{\left[s\atop{ s\atop s} \right]\times (k+1)} -\Gamma  _{2s+j+2}^{\left[s\atop{ s\atop s} \right]\times (2s+j+3)}
= 315\cdot2^{2s+4j+2}\cdot \sum_{i= 2s+j++3}^{k}2^{2i} + \big(735\cdot2^{5s+4j+1} -3255\cdot2^{4s+5j+2} \big)\cdot \sum_{i= 2s+j+3}^{k}2^{i} \\
 & =  105\cdot2^{2k+2s+4j+4}+ 735\cdot2^{k+5s+4j+2} - 3255\cdot2^{k+4s+5j+3}   - 735\cdot2^{7s+5j+4} +2415\cdot2^{6s+6j+5} 
\end{align*}
\begin{align}
& \Updownarrow  \nonumber  \\
&  \Gamma  _{2s+j+2}^{\left[s\atop{ s\atop s} \right]\times k} = 
  105\cdot2^{2k+2s+4j+2}+ 735\cdot2^{k+5s+4j+1} - 3255\cdot2^{k+4s+5j+2}   - 735\cdot2^{7s+5j+4} +2415\cdot2^{6s+6j+5} \nonumber \\
  &   +    105\cdot\big(2^{6s+6j+3} - 3\cdot2^{7s+5j+2} + 2^{8s+4j+2}  \big)    \quad \text{for}\quad k\geq 2s+j+4    \nonumber \\
  & \Updownarrow  \nonumber  \\
  & \Gamma  _{2s+j+2}^{\left[s\atop{ s\atop s} \right]\times k} =  105\cdot\big(2^{2k+2s+4j+2}+ 7\cdot2^{k+5s+4j+1} - 31\cdot2^{k+4s+5j+2} \big) \label{eq 10.115} \\
   &   + 105\cdot\big(2^{8s+4j+2} - 31\cdot2^{7s+5j+2} +93\cdot2^{6s+6j+3} \big)  \quad \text{for}\quad k\geq 2s+j+4 \nonumber
\end{align}
From \eqref{eq 10.112} we see that \eqref{eq 10.115} is valid for k = 2s+j+3.\\

\underline{Proof of \eqref{eq 10.97}}

We then obtain from \eqref{eq 10.110} with j = s-2 \\

 \begin{align}
  &   \Gamma  _{3s}^{\left[s\atop{ s\atop s} \right]\times (3s+1)} 
  =  2^{12s} -7\cdot2^{8s+4(s-2)+2}-91\cdot2^{7s+5(s-2)+2} +41\cdot2^{6s+6(s-2)+3}  = 315\cdot2^{12s-9}    \label{eq 10.116}
\end{align}
 We have from \eqref{eq 9.7} with $ j = s $ and $k \geq 3s+1$  \\
  \begin{equation}
\label{eq 10.117}
 \Delta _{k}\Gamma  _{3s}^{\left[s\atop{ s\atop s} \right]}
  =   \Gamma  _{3s}^{\left[s\atop{ s\atop s} \right]\times (k+1)} -  \Gamma  _{3s}^{\left[s\atop{ s\atop s} \right]\times k} =7\cdot2^{3k+3s-3} -21\cdot2^{2k+6s-6}+7\cdot2^{k+9s-8} 
  \end{equation}
    We obtain from \eqref{eq 10.117} for $k\geq 3s+1$ and  \eqref{eq 10.116} 
  \begin{align*}
& \sum_{i=3s+1}^{k}\big( \Gamma  _{3s}^{\left[s\atop{ s\atop s} \right]\times (i+1)} -  \Gamma  _{3s}^{\left[s\atop{ s\atop s} \right]\times i}\big)
= \sum_{i= 3s+1}^{k} 7\cdot2^{3i+3s-3} -21\cdot2^{2i+6s-6}+7\cdot2^{i+9s-8}   \\
& \Updownarrow   \\
& \Gamma  _{3s}^{\left[s\atop{ s\atop s} \right]\times (k+1)} -\Gamma  _{3s}^{\left[s\atop{ s\atop s} \right]\times (3s+1)}
= 7\cdot2^{3s-3}\cdot \sum_{i= 3s+1}^{k}2^{3i}  -21\cdot2^{6s-6}\cdot \sum_{i= 3s+1}^{k}2^{2i}+ 7\cdot2^{9s-8}\cdot \sum_{i= 3s+1}^{k}2^{i}\\
 & =  2^{3k+3s} -7\cdot2^{2k+6s-4} +7\cdot2^{k+9s-7} -79\cdot2^{12s-7}
\end{align*}
\begin{align}
& \Updownarrow  \nonumber  \\
&  \Gamma  _{3s}^{\left[s\atop{ s\atop s} \right]\times k} = 
     2^{3k+3s-3} -7\cdot2^{2k+6s-6} +7\cdot2^{k+9s-8} -79\cdot2^{12s-7}  +  315\cdot2^{12s-9}   \quad \text{for}\quad k\geq 3s+2   \nonumber \\
    & \Updownarrow  \nonumber  \\
  & \Gamma  _{3s}^{\left[s\atop{ s\atop s} \right]\times k} =    2^{3k+3s-3} -7\cdot2^{2k+6s-6} +7\cdot2^{k+9s-8} - 2^{12s-9}    \label{eq 10.118} 
  \end{align}
From \eqref{eq 10.116} we see that \eqref{eq 10.118} is valid for k = 3s+1.
  \end{proof}
 
 \begin{thm}
 \label{thm 10.9}
 We have \\
 \begin{equation}
\label{eq 10.119}
\Gamma_{i}^{\left[s\atop{ s\atop s} \right]\times k} 
 = \begin{cases}
 1  & \text{if  } i = 0, \\
 105\cdot2^{4i-6} - 21\cdot2^{3i -5}  &  \text{if  }  1\leq i\leq s-1,\;k\geq i+1, \\
  7\cdot2^{k+s-1} -7\cdot2^{2s} +105\cdot2^{4s-6} -21\cdot2^{3s-5}   &    \text{if }  i = s, \; k\geq s+1,\; s\geq 1,  \\
   147\cdot(5\cdot2^{j-1} - 1)\cdot2^{k+s+3j-6} \\
    + 21\cdot\big[5\cdot2^{4s+4j-6} - 2^{3s+3j-5} -(155\cdot2^{j-1} - 35)\cdot2^{2s+4j-7} \big] &  \text{if  }  i = s+j,  \; 1\leq j\leq s-1,\; k\geq s+j+1,\\
     7\cdot2^{2k+2s-2}+ 21\cdot\big[35\cdot2^{k+5s-7} - 39\cdot2^{k+4s-6} \big]  \\
    + 7\cdot\big[15\cdot2^{8s-6} - 465\cdot2^{7s-8} +349\cdot2^{6s-7} \big]  &  \text{if  } i=2s, \; k\geq 2s+1, \\
   105\cdot\big(2^{2k+2s+4j -2}+ 7\cdot2^{k+5s+4j-3} - 31\cdot2^{k+4s+5j-3} \big)  \\
     + 105\cdot\big(2^{8s+4j-2} - 31\cdot2^{7s+5j-3} +93\cdot2^{6s+6j-3} \big) 
     &  \text{if   }\quad  i=2s +1+j, \; k\geq 2s+2+j, \\
     &  0 \leq j\leq s-2, \\
    2^{3k+3s-3} -7\cdot2^{2k+6s-6} +7\cdot2^{k+9s-8} - 2^{12s-9}  &  \text{if  } i=3s, \; k\geq 3s  
     \end{cases}
\end{equation}  
  \begin{equation}
\label{eq 10.120}
\Gamma_{i}^{\left[s\atop{ s\atop s} \right]\times i} 
 = \begin{cases}
  2^{3s+3i-3} -7\cdot2^{4i-6} +3\cdot2^{3i-5}  &  \text{if  }\; 1\leq i\leq s+1, \\
   2^{6s+3j-3} +7\cdot2^{2s+5j-8} -7\cdot2^{2s+4j-7} -7\cdot2^{4s+4j-6} + 3\cdot2^{3s+3j-5} &  \text{if  }\;i = s+j, \;1\leq j \leq s+1,\\
     2^{9s+3j} -7\cdot2^{8s+4j-2} +7\cdot2^{7s+5j-3} - 2^{6s+6j-3}  &  \text{if  }\;i = 2s+1+j, \;0 \leq j \leq s-1,\\
  \end{cases}
\end{equation}  
  We have for $0\leq j\leq s-2,\quad k\geq 2s+2+j$\\
 \begin{align}
 \Gamma  _{2s+1+j}^{\left[s\atop{ s\atop s} \right]\times k} & = 
 16^{3j}\cdot\Gamma  _{2(s-j)+1}^{\left[s-j\atop{ s-j\atop s-j} \right]\times (k-3j)} \label{eq 10.121}\\
 & = 2^{12j}\cdot105\cdot\big(2^{2(k-3j)+2(s-j) -2}+ 7\cdot2^{k-3j+5(s-j)-3} - 31\cdot2^{k-3j+4(s-j)-3} \big) \nonumber \\
  &  + 2^{12j}\cdot 105\cdot\big(2^{8(s-j)-2} - 31\cdot2^{7(s-j)-3} +93\cdot2^{6(s-j)-3} \big) \nonumber \\
  & =  105\cdot\big(2^{2k+2s+4j -2}+ 7\cdot2^{k+5s+4j-3} - 31\cdot2^{k+4s+5j-3} \big) \nonumber \\
  &  + 105\cdot\big(2^{8s+4j-2} - 31\cdot2^{7s+5j-3} +93\cdot2^{6s+6j-3} \big) \nonumber
 \end{align}
 
  We have for $0\leq j\leq s-1.$\\
 \begin{align}
 \Gamma  _{2s+1+j}^{\left[s\atop{ s\atop s} \right]\times (2s+1+j)} & = 
 16^{3j}\cdot\Gamma  _{2(s-j)+1}^{\left[s-j\atop{ s-j\atop s-j} \right]\times (2(s-j)+1)} \label{eq 10.122}\\
 &  =   2^{12j}\cdot\big(2^{9(s-j)} -7\cdot2^{8(s-j)-2} +7\cdot2^{7(s-j)-3} - 2^{6(s-j)-3}\big)\nonumber \\
 & =  2^{9s+3j} -7\cdot2^{8s+4j-2} +7\cdot2^{7s+5j-3} - 2^{6s+6j-3}\nonumber
 \end{align}
 We have for $k\geq 3s. $\\
  \begin{align}
 \Gamma  _{3s}^{\left[s\atop{ s\atop s} \right]\times k} & = 
 16^{3(s-1)}\cdot\Gamma  _{3}^{\left[1\atop{ 1\atop 1} \right]\times (k-3(s-1))} \label{eq 10.123}\\
 &  =   2^{12s-12}\cdot\big(2^{3k-9s+9} -7\cdot2^{2k-6s+6} +7\cdot2^{k-3s+4} - 2^{3}\big) \nonumber \\
 & =  2^{3k+3s-3} -7\cdot2^{2k+6s-6} +7\cdot2^{k+9s-8} - 2^{12s-9} \nonumber 
 \end{align}

\end{thm} 
 \begin{thm}
\label{thm 10.10} Let q be a rational integer $ \geq 1,$ then \\
 \begin{align}
    g_{k,s}(t,\eta,\xi ) =  g(t,\eta,\xi ) & = \sum_{deg Y\leq k-1}\sum_{deg Z \leq  s-1}E(tYZ)\sum_{deg U \leq s-1}E(\eta YU)\sum_{deg V \leq s-1}E(\eta YV)   =
   2^{3s+k- r(D^{\left[s\atop{ s\atop s} \right]\times k}(t,\eta,\xi)) }, \label{eq 10.124} \\
   & \nonumber \\
    \int_{\mathbb{P}\times \mathbb{P}\times \mathbb{P}} g^{q}(t,\eta,\xi  )dtd\eta d\xi  
 & =2^{(3s+k)q}\cdot 2^{-3k-3s+3}\cdot \sum_{i = 0}^{\inf(3s,k)} \Gamma_{i}^{\left[s\atop{ s\atop s} \right]\times k}  \cdot2^{- qi}. \label{eq 10.125} \\
  & \nonumber 
 \end{align}
 \end{thm}
   \begin{thm}
\label{thm 10.11}
 We denote by  $ R_{q}(k,s) $ the number of solutions \\
 $(Y_1,Z_1,U_{1},V_{1}, \ldots,Y_q,Z_q,U_{q},V_{q}) $  of the polynomial equations
   \[\left\{\begin{array}{c}
 Y_{1}Z_{1} +Y_{2}Z_{2}+ \ldots + Y_{q}Z_{q} = 0,  \\
   Y_{1}U_{1} + Y_{2}U_{2} + \ldots  + Y_{q}U_{q} = 0,\\
    Y_{1}V_{1} + Y_{2}V_{2} + \ldots  + Y_{q}V_{q} = 0,\\  
 \end{array}\right.\]
  satisfying the degree conditions \\
                   $$  degY_i \leq k-1 , \quad degZ_i \leq s-1 ,\quad degU_{i}\leq s-1,\quad degV_{i}\leq s-1 \quad for \quad 1\leq i \leq q. $$ \\                           
Then \\
\begin{align}
  R_{q}(k,s) & =  \int_{\mathbb{P}\times \mathbb{P}\times \mathbb{P}} g_{k,s}^{q}(t,\eta,\xi  )dtd\eta d\xi  
 = 2^{(3s+k)q}\cdot 2^{-3k-3s+3}\cdot \sum_{i = 0}^{\inf(3s,k)} \Gamma_{i}^{\left[s\atop{ s\atop s} \right]\times k}  \cdot2^{- qi}. \label{eq 10.126}
\end{align}
  \end{thm}

\begin{example}  s = 1, $k \geq i+1$ for $0\leq i\leq 2$
  \[ \Gamma_{i}^{\left[1\atop{ 1\atop 1} \right]\times k}=
 \begin{cases}
1  &\text{if  } i = 0 \\
7\cdot(2^{k} - 1)   & \text{if   } i = 1 \\
7\cdot(2^{k}-1)\cdot(2^{k} - 2)  & \text{if   } i = 2 \\
2^{3k} - 7\cdot2^{2k} + 7\cdot2^{k+ 1}- 2^{3} & \text{if  } i = 3,\;k\geq 3
\end{cases}
\]
\end{example}

\begin{example}  s = 2, $k \geq i+1$ for $0\leq i\leq 5$
\[ \Gamma_{i}^{\left[2\atop{ 2\atop 2} \right]\times k}=
 \begin{cases}
1  &\text{if  }  i = 0 \\
21   &\text{if  }  i = 1 \\
7\cdot2^{k+ 1} + 266  & \text{if   } i = 2 \\
147\cdot2^{k+ 1} + 1344   & \text{if   } i = 3 \\
7\cdot2^{2k+2} + 651\cdot2^{k+ 2}  - 22624 & \text{if   } i = 4 \\
105\cdot2^{2k+2} - 315\cdot2^{k+ 5} + 53760   & \text{if   } i = 5 \\
2^{3k+3} - 7\cdot2^{2k+6} + 7\cdot2^{k+ 10} - 32768  & \text{if  } i = 6,\;k\geq 6
\end{cases}
\]

\end{example}
\begin{example} s = 2, k = 6.\\

 The number  $\Gamma_{i}^{\left[2\atop{ 2\atop 2} \right]\times 6}$ of rank i matrices of the form \\
  $$  \left ( \begin{array} {cccccc}
\alpha _{1} & \alpha _{2} & \alpha _{3} & \alpha _{4} & \alpha _{5} & \alpha _{6}  \\
\alpha _{2 } & \alpha _{3} & \alpha _{4} & \alpha _{5} & \alpha _{6} & \alpha _{7} \\
\beta  _{1} & \beta  _{2} & \beta  _{3}  & \beta  _{4} & \beta  _{5} & \beta  _{6}\\
\beta  _{2} & \beta  _{3} & \beta  _{4}  & \beta  _{5} & \beta  _{6} & \beta  _{7}\\
\gamma  _{1} & \gamma  _{2} & \gamma  _{3}  & \gamma  _{4} & \gamma  _{5} & \gamma  _{6} \\
\gamma  _{2} & \gamma  _{3} & \gamma  _{4}  & \gamma  _{5} & \gamma  _{6} & \gamma  _{7}
 \end{array}  \right) $$
is equal to 
\[ \begin{cases}
1  &\text{if  }  i = 0 \\
21   &\text{if  }  i = 1 \\
1162     & \text{if   } i = 2 \\
20160          & \text{if   } i = 3 \\
258720           & \text{if   } i = 4 \\
1128960           & \text{if   } i = 5 \\
688128  & \text{if  } i = 6
\end{cases}
\]
  The number of solutions \\
 $(Y_1,Z_1,U_{1},V_{1}, \ldots,Y_q,Z_q,U_{q},V_{q}) $  of the polynomial equations
   \[\left\{\begin{array}{c}
 Y_{1}Z_{1} +Y_{2}Z_{2}+ \ldots + Y_{q}Z_{q} = 0,  \\
   Y_{1}U_{1} + Y_{2}U_{2} + \ldots  + Y_{q}U_{q} = 0,\\
    Y_{1}V_{1} + Y_{2}V_{2} + \ldots  + Y_{q}V_{q} = 0,\\  
 \end{array}\right.\]
  satisfying the degree conditions \\
                   $$  degY_i \leq 5 , \quad degZ_i \leq 1 ,\quad degU_{i}\leq 1,\quad degV_{i}\leq 1 \quad for \quad 1\leq i \leq q. $$ \\                           
is equal to 
\begin{align*}
&  R_{q}(6,2)  =  \int_{\mathbb{P}\times \mathbb{P}\times \mathbb{P}} g_{6,2}^{q}(t,\eta,\xi  )dtd\eta d\xi  
 = 2^{12q-21}\cdot \sum_{i = 0}^{6} \Gamma_{i}^{\left[2\atop{ 2\atop 2} \right]\times 6}  \cdot2^{- qi}\\
& = 2^{12q-21}\cdot\big(1 + 21\cdot2^{-q} + 1162\cdot2^{-2q} + 20160\cdot2^{-3q} + 258720\cdot2^{-4q} + 1128960 \cdot2^{-5q}
 + 688128 \cdot2^{-6q}\big ) \\
 & =   2^{6q-21}\cdot\big(2^{6q} + 21\cdot2^{5q} + 1162\cdot2^{4q} + 20160\cdot2^{3q} + 258720\cdot2^{2q} + 1128960 \cdot2^{q}
 + 688128 \big ) 
\end{align*}

\end{example}

\begin{example}  s = 3, $k \geq i+1$ for $0\leq i\leq 8$
\[ \Gamma_{i}^{\left[3\atop{ 3\atop 3} \right]\times k}=
 \begin{cases}
1 &\text{if  }  i = 0 \\
21 &\text{if  }  i = 1 \\
378  &\text{if  }  i = 2 \\
7\cdot2^{k+ 2} +  5936  & \text{if   } i = 3 \\
147\cdot2^{k+ 2} + 84672  & \text{if   } i = 4 \\
147\cdot9\cdot2^{k+ 3} + 959616  & \text{if   } i = 5 \\
7\cdot2^{2k+4} + 2121\cdot2^{k+ 6} + 5863424   & \text{if   } i = 6\\
105\cdot2^{2k+4} + 2625\cdot2^{k+ 9} - 92897280  & \text{if   } i = 7 \\
105\cdot2^{2k+8} - 315\cdot2^{k+ 14} + 220200960  & \text{if   } i = 8 \\
2^{3k+6} - 7\cdot2^{2k+12} + 7\cdot2^{k+ 19} - 134217728  & \text{if  } i = 9,\;k\geq 9
\end{cases}
\]
\end{example}

 \begin{example}  s = 3, k = 5, q = 3. \\
 
 The number  $\Gamma_{i}^{\left[3\atop{ 3\atop 3} \right]\times 5}$ of rank i matrices of the form \\
  $$  \left ( \begin{array} {ccccc}
\alpha _{1} & \alpha _{2} & \alpha _{3} & \alpha _{4} & \alpha _{5}  \\
\alpha _{2 } & \alpha _{3} & \alpha _{4} & \alpha _{5} & \alpha _{6} \\
\alpha _{3} & \alpha _{4} & \alpha _{5} & \alpha _{6} & \alpha _{7} \\
\beta  _{1} & \beta  _{2} & \beta  _{3}  & \beta  _{4} & \beta  _{5}\\
\beta  _{2} & \beta  _{3} & \beta  _{4}  & \beta  _{5} & \beta  _{6}\\
\beta  _{3} & \beta  _{4} & \beta  _{5}  & \beta  _{6} & \beta  _{7}\\
\gamma  _{1} & \gamma  _{2} & \gamma  _{3}  & \gamma  _{4} & \gamma  _{5}\\
\gamma  _{2} & \gamma  _{3} & \gamma  _{4}  & \gamma  _{5} & \gamma  _{6}\\
\gamma  _{3} & \gamma  _{4} & \gamma  _{5}  & \gamma  _{6} & \gamma  _{7}
 \end{array}  \right) $$
is equal to 
 \[  \begin{cases}
1 &\text{if  }  i = 0 \\
21 &\text{if  }  i = 1 \\
378  &\text{if  }  i = 2 \\
6832  & \text{if   } i = 3 \\
103488 & \text{if   } i = 4 \\
1986432 & \text{if   } i = 5 
\end{cases}
\]

 The number of solutions \\
 $(Y_1,Z_1,U_{1},V_{1},Y_2,Z_2,U_{2},V_{2} ,Y_3,Z_3,U_{3},V_{3}) $  of the polynomial equations
   \[\left\{\begin{array}{c}
 Y_{1}Z_{1} +Y_{2}Z_{2} + Y_{3}Z_{3} = 0,  \\
   Y_{1}U_{1} + Y_{2}U_{2}  + Y_{3}U_{3} = 0,\\
    Y_{1}V_{1} + Y_{2}V_{2} + Y_{3}V_{3} = 0,\\  
 \end{array}\right.\]
  satisfying the degree conditions \\
                   $$  degY_i \leq 4 , \quad degZ_i \leq 2 ,\quad degU_{i}\leq 2,\quad degV_{i}\leq 2 \quad for \quad 1\leq i \leq 3. $$ \\                           

is equal to 
\begin{align*}
&  R_{3}(5,3)  =  \int_{\mathbb{P}\times \mathbb{P}\times \mathbb{P}} g_{5,3}^{3}(t,\eta,\xi  )dtd\eta d\xi  
 = 2^{33}\cdot \sum_{i = 0}^{5} \Gamma_{i}^{\left[3\atop{ 3\atop 3} \right]\times 5}  \cdot2^{- 3i}\\
& = 2^{33}\cdot\big(1 + 21\cdot2^{-3} + 378\cdot2^{-6} + 6832\cdot2^{-9} + 103488\cdot2^{-12} + 1986432 \cdot2^{-15}\big )  = 3563904\times 2^{18}
\end{align*}

\end{example}
\newpage
 \section{\textbf{Computation of $ \Delta _{k}\Gamma_{i}^{\left[s\atop{ s+m\atop s+m} \right]} \quad for \; 0 \leq i \leq  3s+2m $}}
\label{sec 11} 
     \subsection{Notation}
  \label{subsec 1}
  \begin{defn}
  \label{defn 11.1}
In this section we define $ \gamma _{j} $ to be equal to $(21\cdot2^{3j-4} - 3\cdot2^{2j-3})$ for $1\leq j\leq m-1$. 
\end{defn}
   \subsection{Introduction}
We adapt the method in Section \ref{sec 9} to compute explicitly 
$ \Delta _{k}\Gamma_{i}^{\left[s\atop{ s+m\atop s+m} \right]} \quad for \; s \leq i \leq  3s+2m $

  \label{subsec 2}
   \subsection{\textbf{Computation of $ \Delta _{k}\Gamma_{s+j}^{\left[s\atop{ s+m\atop s+m+l} \right]} \quad for \; 0 \leq j \leq  m-1,\; l\geq 0 $}}
  \label{subsec 3}
   \begin{lem}
\label{lem 11.2} We have the following reduction formula for $0\leq j\leq m,\; l\geq 0 $\\

\begin{align}
  &   \Delta _{k}\Gamma_{s+j}^{\left[s\atop{ s +m \atop s +m+l } \right]}
    - 8\cdot \Delta _{k}\Gamma_{s+j-1}^{\left[s\atop{ s +m \atop s +m+l-1 } \right]}
    =  4^{j}\cdot\Delta _{k}\Gamma_{s}^{\left[s\atop{ s +(m-j)\atop s +m+l} \right]} 
      +\sum_{q = 0}^{j-1}4^{q}\cdot \Delta _{k}\omega_{s+j- q} (s,m - q, l+ q) \label{eq 11.1}
\end{align}
From \eqref{eq 11.1} we deduce successively for  $ j = 0,1,\ldots,m-1,m $

\begin{equation}
\label{eq 11.2}
  \Delta _{k}\Gamma_{s+j}^{\left[s\atop{ s+m\atop s+m+l} \right]} =
 \begin{cases}
 2^{k+s-1} & \text{if   } j = 0 \\
  \gamma _{j}\cdot2^{k+s-1} = (21\cdot2^{3j-4} - 3\cdot2^{2j-3})\cdot2^{k+s-1} & \text{if } 1\leq j\leq m-1,\; l\geq 0\\
    21\cdot2^{k+s+3m-5} + 13\cdot2^{k+s+2m-4} & \text{if }j = m\geq 1,\;l \geq 1\\
    21\cdot2^{k+s+3m-5} + 45\cdot2^{k+s+2m-4} & \text{if }j =  m\geq 1,\;l = 0
  \end{cases}
\end{equation}

\end{lem}
\begin{proof} 
From \eqref{eq 8.3} we obtain successively the following equations \\
  \small   
  \begin{align*}
 &  \Delta _{k}\Gamma_{s+j}^{\left[s\atop{ s +m\atop s+m+l } \right]}
   - 4\cdot \Delta _{k}\Gamma_{s+j -1}^{\left[s\atop{ s +m-1\atop s+m+l } \right]}  = 
   8\cdot\Big[  \Delta _{k}\Gamma_{s+j -1}^{\left[s\atop{ s+m \atop s +m+l-1} \right]}
   - 4\cdot \Delta _{k}\Gamma_{s+j -2}^{\left[s\atop{ s +m-1\atop s +m+l-1} \right]}\Big] + \Delta _{k}\omega_{s+j} (s,m,l) \\
    &  4\cdot\bigg[  \Delta _{k}\Gamma_{s+j-1}^{\left[s \atop{ s +m-1\atop s+m+l } \right]}
   - 4\cdot \Delta _{k}\Gamma_{s+j -2}^{\left[s\atop{ s+m-2 \atop s +m+l } \right]}\bigg]  = 
4\cdot\bigg[   8\cdot\Big[  \Delta _{k}\Gamma_{s+j -2}^{\left[s\atop{ s+m-1 \atop s +m+l-1 } \right]}
   - 4\cdot \Delta _{k}\Gamma_{s+j -3}^{\left[s\atop{ s +m-2\atop s +m+l-1 } \right]}\Big] + \Delta _{k}\omega_{s+j-1} (s,m-1, l+1) \bigg]\\
    &  4^{2}\cdot\bigg[   \Delta _{k}\Gamma_{s+j-2}^{\left[s \atop{ s +m-2 \atop s+m+l } \right]}
   - 4\cdot \Delta _{k}\Gamma_{s+j -3}^{\left[s\atop{ s +m-3 \atop s +m+l} \right]}\bigg]  = 
 4^{2}\cdot\bigg[  8\cdot\Big[  \Delta _{k}\Gamma_{s+j -3}^{\left[s \atop{ s +m-2 \atop s +m+l-1} \right]}
   - 4\cdot \Delta _{k}\Gamma_{s+j - 4}^{\left[s \atop{ s +m-3 \atop s +m+l-1 } \right]}\Big] + \Delta _{k}\omega_{s+j-2} (s,m-2,l+2)\bigg] \\
     &  \vdots       \\
   &  4^{j-1}\cdot\bigg[   \Delta _{k}\Gamma_{s+j- (j-1)}^{\left[s \atop{ s +m -(j-1) \atop s+m+l } \right]}
   - 4\cdot \Delta _{k}\Gamma_{s+j -j}^{\left[s\atop{ s +m-j \atop s +m+l} \right]}\bigg]\\
    & =  4^{j-1}\cdot\bigg[  8\cdot\Big[  \Delta _{k}\Gamma_{s+j -j}^{\left[s \atop{ s +m-(j-1) \atop s +m+l-1} \right]}
   - 4\cdot \Delta _{k}\Gamma_{s+j - (j+1)}^{\left[s \atop{ s +m- j \atop s +m+l-1 } \right]}\Big]
    + \Delta _{k}\omega_{s+j- (j-1)} (s,m-(j-1),l+(j-1))\bigg]  
 \end{align*}
 By summing  the left-hand side and the right-hand side of the above equations  we obtain the following equality 
     \begin{align*}
   &   \Delta _{k}\Gamma_{s+j}^{\left[s\atop{ s +m \atop s +m+l } \right]}
   - 4^{j}\cdot \Delta _{k}\Gamma_{s}^{\left[s\atop{ s +(m-j)\atop s +m+l} \right]}\\
    &  = 8\cdot\bigg[ \Delta _{k}\Gamma_{s+j-1}^{\left[s\atop{ s +m \atop s +m+l-1 } \right]} - 4^{j}\cdot \Delta _{k}\Gamma_{s -1}^{\left[s\atop{ s +(m-j)\atop s +m+l-1} \right]}\bigg]
   +\sum_{q = 0}^{j-1}4^{q}\cdot \Delta _{k}\omega_{s+j- q} (s,m - q, l+ q) \\
   & \Updownarrow \\
   &  \Delta _{k}\Gamma_{s+j}^{\left[s\atop{ s +m \atop s +m+l } \right]}
    - 8\cdot \Delta _{k}\Gamma_{s+j-1}^{\left[s\atop{ s +m \atop s +m+l-1 } \right]}
     =  4^{j}\cdot\bigg[ \Delta _{k}\Gamma_{s}^{\left[s\atop{ s +(m-j)\atop s +m+l} \right]}
      - 8\cdot \Delta _{k}\Gamma_{s -1}^{\left[s\atop{ s +(m-j)\atop s +m+l-1} \right]}\bigg]
       +\sum_{q = 0}^{j-1}4^{q}\cdot \Delta _{k}\omega_{s+j- q} (s,m - q, l+ q) \\
       & \Updownarrow \\
       &   \Delta _{k}\Gamma_{s+j}^{\left[s\atop{ s +m \atop s +m+l } \right]}
    - 8\cdot \Delta _{k}\Gamma_{s+j-1}^{\left[s\atop{ s +m \atop s +m+l-1 } \right]}
    =  4^{j}\cdot\Delta _{k}\Gamma_{s}^{\left[s\atop{ s +(m-j)\atop s +m+l} \right]} 
      +\sum_{q = 0}^{j-1}4^{q}\cdot \Delta _{k}\omega_{s+j- q} (s,m - q, l+ q) 
       \end{align*}
       To prove \eqref{eq 11.2} we shall apply Lemma \ref{lem 7.4}.\\
       
      \textbf{The case j = 1} 
       \begin{align*}   
          & \Delta _{k}\Gamma_{s+1}^{\left[s\atop{ s +m \atop s +m+l } \right]}
    =  8\cdot \Delta _{k}\Gamma_{s}^{\left[s\atop{ s +m \atop s +m+l-1 } \right]}
    +  4\cdot\Delta _{k}\Gamma_{s}^{\left[s\atop{ s +(m-1)\atop s +m+l} \right]} 
      + \Delta _{k}\omega_{s+1} (s,m , l) \\
      & = 8\cdot2^{k+s-1} +4\cdot2^{k+s-1} -3\cdot2^{k+s-1}= 9\cdot2^{k+s-1}= \gamma _{1}\cdot2^{k+s-1} \quad \text{if  }\quad m\geq 2,\quad l\geq 0 \\
      \end{align*}
    \textbf{The case j = 2}    
    \begin{align*} 
     & \\
       & \Delta _{k}\Gamma_{s+2}^{\left[s\atop{ s +m \atop s +m+l } \right]}
    =  8\cdot \Delta _{k}\Gamma_{s +1}^{\left[s\atop{ s +m \atop s +m+l-1 } \right]}
    +  4^{2}\cdot\Delta _{k}\Gamma_{s}^{\left[s\atop{ s +(m-2)\atop s +m+l} \right]} 
      + \Delta _{k}\omega_{s+2} (s,m , l)  + 4\cdot\Delta _{k}\omega_{s+1} (s,m-1 , l+1)\\
      & = 8\cdot9\cdot2^{k+s-1}+ 4^{2}\cdot2^{k+s-1} +  2^{k+s} -  4\cdot3\cdot2^{k+s-1} = 78\cdot2^{k+s-1}=\gamma _{2}\cdot2^{k+s-1} \quad \text{if  }\quad m\geq 3,\quad l\geq 0 
        \end{align*}
     \textbf{The case j = 3} 
         \begin{align*}  
         & \\
         & \Delta _{k}\Gamma_{s+3}^{\left[s\atop{ s +m \atop s +m+l } \right]}
    =  8\cdot \Delta _{k}\Gamma_{s +2}^{\left[s\atop{ s +m \atop s +m+l-1 } \right]}
    +  4^{3}\cdot\Delta _{k}\Gamma_{s}^{\left[s\atop{ s +(m- 3)\atop s +m+l} \right]} \\
     & + \Delta _{k}\omega_{s+3} (s,m , l)  + 4\cdot\Delta _{k}\omega_{s+2} (s,m-1 , l+1) + 4^{2}\cdot\Delta _{k}\omega_{s+1} (s,m-2 , l+2)\\
      & = 8\cdot78\cdot2^{k+s-1}+ 4^{3}\cdot2^{k+s-1} + 4\cdot 2^{k+s} -  4^{2}\cdot3\cdot2^{k+s-1} = 648\cdot2^{k+s-1} =\gamma _{3}\cdot2^{k+s-1} \quad \text{if  }\quad m\geq 4,\quad l\geq 0 
        \end{align*}
      \textbf{The case  general $1\leq j\leq m-1$}  \\
        
      Proof by induction on $ j \leq m-1 $.
      Assume that $ \Delta _{k}\Gamma_{s +j-1}^{\left[s\atop{ s +m \atop s +m+l-1 } \right]} =\gamma _{j-1}\cdot2^{k+s-1}$
       for $1\leq j-1\leq m-2 $.\\
       We shall prove that  $ \Delta _{k}\Gamma_{s +j}^{\left[s\atop{ s +m \atop s +m+l } \right]} =\gamma _{j}\cdot2^{k+s-1}$ for 
          $1\leq j\leq m-1 $.\\
      \begin{align*}     
      &  \\
      & \Delta _{k}\Gamma_{s+j}^{\left[s\atop{ s +m \atop s +m+l } \right]}
    =  8\cdot \Delta _{k}\Gamma_{s +j-1}^{\left[s\atop{ s +m \atop s +m+l-1 } \right]}
    +  4^{j}\cdot\Delta _{k}\Gamma_{s}^{\left[s\atop{ s +(m-j)\atop s +m+l} \right]} 
    + \sum_{q = 0}^{j -1} 4^{q}\cdot\Delta _{k}\omega_{s+ j -q} (s,m-q , l+q) \\
      & = 8\cdot\gamma _{j-1}\cdot2^{k+s-1} + 2^{2j}\cdot2^{k+s-1} 
    +  \sum_{q = j -2}^{j -1}4^{q}\cdot\Delta _{k}\omega_{s+ j -q} (s,m-q , l+q) \\
     & = 8\cdot(21\cdot2^{3(j-1)}  - 3\cdot2^{2(j-1)-3})\cdot2^{k+s-1} +  2^{2j}\cdot2^{k+s-1}\\
      &  + 4^{j-2}\cdot\Delta _{k}\omega_{s+ j - (j-2)} (s,m- (j-2) , l+ (j-2)) +  4^{j-1}\cdot\Delta _{k}\omega_{s+ j - (j-1)} (s,m-(j-1) , l+(j-1)) \\
        & =  8\cdot[21\cdot2^{3(j-1)}  - 3\cdot2^{2(j-1)-3}]\cdot2^{k+s-1} +  2^{2j}\cdot2^{k+s-1} + 2^{2j-4}\cdot2^{k+s}+ 2^{2j-2}\cdot(-3\cdot2^{k+s-1})\\
        & = [21\cdot2^{3j - 4}  - 3\cdot2^{2j -3}]\cdot2^{k+s-1} = \gamma _{j}\cdot2^{k+s-1}\quad \text{if  }\quad m\geq j+1,\quad l\geq 0 
       \end{align*} 
    \textbf{The case $j=m,\;l\geq 1 $}  \\ 
        \begin{align*}   
        &  \\
       &  \Delta _{k}\Gamma_{s+m}^{\left[s\atop{ s +m \atop s +m+l } \right]}
    - 8\cdot \Delta _{k}\Gamma_{s+m-1}^{\left[s\atop{ s +m \atop s +m+l-1 } \right]}
    =  4^{m}\cdot\Delta _{k}\Gamma_{s}^{\left[s\atop{ s \atop s +m+l} \right]} 
      +\sum_{q = 0}^{m -1}4^{q}\cdot \Delta _{k}\omega_{s+m- q} (s,m - q, l+ q) \\
      &  \Delta _{k}\Gamma_{s+m}^{\left[s\atop{ s +m \atop s +m+l } \right]}
       = 2^{2m}\cdot3\cdot2^{k+s-1}+ 8\cdot\gamma _{m-1}\cdot2^{k+s-1}+\sum_{q = 0}^{m -1}4^{q}\cdot \Delta _{k}\omega_{s+m- q} (s,m - q, l+ q) \\
       & = 2^{2m}\cdot3\cdot2^{k+s-1}+ 8\cdot\gamma _{m-1}\cdot2^{k+s-1}+\sum_{q = m-2}^{m -1}4^{q}\cdot \Delta _{k}\omega_{s+m- q} (s,m - q, l+ q) \\
       & = 3\cdot2^{k+s+2m-1} + 8\cdot[21\cdot2^{3m-7} - 3\cdot2^{2m -5}]\cdot2^{k+s-1}\\
       &  + 2^{2m-4}\cdot \Delta _{k}\omega_{s+ 2} (s, 2, l+ (m-2)) +  2^{2m-2}\cdot \Delta _{k}\omega_{s+ 1} (s, 1, l+ (m-1)) \\
       & = 3\cdot2^{k+s+2m-1} + [21\cdot2^{3m-7} - 3\cdot2^{2m -5}]\cdot2^{k+s +2}
        +  2^{2m-4}\cdot2^{k+s} +  2^{2m-2}\cdot (-3\cdot2^{k+s-1}) \\
        & = 21\cdot2^{k+s+3m-5} +13\cdot2^{k+s+2m-4}\\
       & \textbf{The case $j=m,\;l = 0 $}\\
          &  \Delta _{k}\Gamma_{s+m}^{\left[s\atop{ s +m \atop s +m } \right]}
       = 2^{2m}\cdot3\cdot2^{k+s-1}+ 8\cdot \Delta _{k}\Gamma_{s+m-1}^{\left[s\atop{ s +m -1\atop s +m } \right]}
       +\sum_{q = 0}^{m -1}4^{q}\cdot \Delta _{k}\omega_{s+m- q} (s,m - q,  q) \\
        & = 2^{2m}\cdot3\cdot2^{k+s-1}+ 8\cdot[ 21\cdot2^{k+s+3m-8} +13\cdot2^{k+s+2m-6}]
        +\sum_{q = m-2}^{m -1}4^{q}\cdot \Delta _{k}\omega_{s+m- q} (s,m - q,  q) \\
           & = 2^{2m}\cdot3\cdot2^{k+s-1}+ 8\cdot[ 21\cdot2^{k+s+3m-8} +13\cdot2^{k+s+2m-6}]
          + 2^{2m-4}\cdot \Delta _{k}\omega_{s+ 2} (s, 2,  m-2) +  2^{2m-2}\cdot \Delta _{k}\omega_{s+ 1} (s, 1,  m-1) \\
           & = 2^{2m}\cdot3\cdot2^{k+s-1}+ 8\cdot[ 21\cdot2^{k+s+3m-8} +13\cdot2^{k+s+2m-6}] 
            +  2^{2m-4}\cdot2^{k+s} +  2^{2m-2}\cdot (-3\cdot2^{k+s-1}) \\
       & = 21\cdot2^{k+s+3m-5} + 45\cdot2^{k+s+2m-4}
      \end{align*}      
       \end{proof}

  \subsection{\textbf{Computation of $ \Delta _{k}\Gamma_{s+j}^{\left[s\atop{ s\atop s+l} \right]}-8\cdot\Delta _{k}\Gamma_{s+j-1}^{\left[s\atop{ s\atop s+l-1} \right]} \quad for \; 0 \leq j \leq  s-1,\; l\geq 0 $}}
  \label{subsec 4}
    \begin{lem}
\label{lem 11.3} 
    Let $j$ and $l$ be integers such that $2\leq j\leq s-1,\; l\geq 0. $\\
 We have the following reduction formula for $2 \leq q\leq s-2 $\\
 \begin{align}
   &   \Delta _{k}\Gamma_{s+j}^{\left[s\atop{ s  \atop s +l } \right]}
    -8\cdot\Delta _{k}\Gamma_{s+j-1}^{\left[s\atop{ s  \atop s +l -1 } \right]}
     =   4^{2q+2}\cdot \Delta _{k}\Gamma_{s +j-2q-2}^{\left[s -q-1\atop{ s -q-1\atop s + l} \right]}
   - 8\cdot 4^{2q+2}\cdot \Delta _{k}\Gamma_{s +j-2q-3}^{\left[s -q-1\atop{ s -q-1\atop s + l -1} \right]} \label{eq 11.3}\\
   &   +\sum_{i = 0}^{q}4^{2i}\cdot \Delta _{k}\omega_{s -i +(j-i)} (s-i,0, l+ i)  + 
   \sum_{i = 0}^{q}4^{2i+1}\cdot \Delta _{k}\omega_{s - i- 1+(j-i)} (s-i-1,1, l+ i) \nonumber
\end{align}
\begin{equation*}
\label{eq 11.4}
\text{Setting}\quad  \begin{cases}
 q=j & \text{if   } 2\leq j\leq s-3, \\
 q= s -3   & \text{if   } j = s-2, \\
  q= s -2     &   \text{if  } j = s -1.
 \end{cases}
 \end{equation*}
in \eqref{eq 11.3} we obtain 
\begin{equation}
\label{eq 11.4}
 \Delta _{k}\Gamma_{s+j}^{\left[s\atop{ s \atop s+l} \right]} -   8\cdot\Delta _{k}\Gamma_{s+j -1}^{\left[s\atop{ s \atop s+l-1} \right]} 
= \begin{cases}
3\cdot2^{k+s-1} & \text{if   } j = 0,\; l\geq 0. \\
15\cdot2^{k+s-1} & \text{if   } j = 1,\;  l\geq 0.\\
 63\cdot2^{k+s+3j-6} &   \text{if  }\; 2\leq j\leq s-1, \: l\geq 0.
 \end{cases}
 \end{equation}
 \end{lem}
 
 \begin{proof}
  From \eqref{eq 8.3} we obtain successively the followings equations
     
 \begin{align*}
 &  \Delta _{k}\Gamma_{s+j}^{\left[s\atop{ s \atop s+l } \right]}
   - 4\cdot \Delta _{k}\Gamma_{s+j -1}^{\left[s -1\atop{ s \atop s+l } \right]}  = 
   8\cdot\left[  \Delta _{k}\Gamma_{s+j -1}^{\left[s\atop{ s \atop s + l-1} \right]}
   - 4\cdot \Delta _{k}\Gamma_{s+j -2}^{\left[s\atop{ s -1\atop s + l-1} \right]}\right] + \Delta _{k}\omega_{s+j} (s,0,l) \\
     &  4\cdot\left[  \Delta _{k}\Gamma_{s+j-1}^{\left[s -1\atop{ s \atop s+l } \right]}
   - 4\cdot \Delta _{k}\Gamma_{s+j -2}^{\left[s -1\atop{ s -1 \atop s +l } \right]}\right]  = 
4\cdot\left[   8\cdot\Big[  \Delta _{k}\Gamma_{s+j -2}^{\left[s-1\atop{ s  \atop s +l-1 } \right]}
   - 4\cdot \Delta _{k}\Gamma_{s+j -3}^{\left[s -1\atop{ s -1\atop s +l-1 } \right]}\Big] + \Delta _{k}\omega_{s+j-1} (s-1,1, l) \right]\\
     &  4^{2}\cdot\left[   \Delta _{k}\Gamma_{s+j-2}^{\left[s -1\atop{ s -1 \atop s+l } \right]}
   - 4\cdot \Delta _{k}\Gamma_{s+j -3}^{\left[s\atop{ s -2 \atop s +l} \right]}\right]  = 
 4^{2}\cdot\left[  8\cdot\Big[  \Delta _{k}\Gamma_{s+j -3}^{\left[s -1\atop{ s -1 \atop s +l-1} \right]}
   - 4\cdot \Delta _{k}\Gamma_{s+j - 4}^{\left[s -2\atop{ s -1 \atop s +l-1 } \right]}\Big] + \Delta _{k}\omega_{s+j-2} (s-1,0,l+1)\right] \\
     &  \vdots       \\
      &  4^{2q}\cdot\left[   \Delta _{k}\Gamma_{s+j -2q}^{\left[s -q\atop{ s -q \atop s+l } \right]}
   - 4\cdot \Delta _{k}\Gamma_{s+j -2q-1}^{\left[s -q-1\atop{ s -q \atop s +l} \right]}\right]
     =  4^{2q}\cdot\left[  8\cdot\Big[  \Delta _{k}\Gamma_{s+j -2q-1}^{\left[s -q\atop{ s -q \atop s +l-1} \right]}
   - 4\cdot \Delta _{k}\Gamma_{s+j -2q-2}^{\left[s -q-1\atop{ s -q \atop s +l -1 } \right]}\Big] + \Delta _{k}\omega_{s+j- 2q} (s -q,0,l+ q)\right] \\
    &  4^{2q +1}\cdot\left[   \Delta _{k}\Gamma_{s+j -2q-1}^{\left[s -q-1\atop{ s -q \atop s+l } \right]}
   - 4\cdot \Delta _{k}\Gamma_{s+j -2q-2}^{\left[s -q-1\atop{ s -q-1 \atop s +l} \right]}\right] \\
     &  =  4^{2q+1}\cdot\left[  8\cdot\Big[  \Delta _{k}\Gamma_{s+j -2q-2}^{\left[s -q-1\atop{ s -q \atop s +l-1} \right]}
   - 4\cdot \Delta _{k}\Gamma_{s+j -2q-3}^{\left[s -q-1\atop{ s -q-1 \atop s +l -1 } \right]}\Big] + \Delta _{k}\omega_{s+j- 2q-1} (s -q-1,1,l+ q)\right] 
   \end{align*}
By summing  the above equations we obtain after reduction \eqref{eq 11.3}. 
     \begin{align*}
   &   \Delta _{k}\Gamma_{s+j}^{\left[s\atop{ s  \atop s +l } \right]}
    -8\cdot\Delta _{k}\Gamma_{s+j-1}^{\left[s\atop{ s  \atop s +l -1 } \right]}
     =   4^{2q+2}\cdot \Delta _{k}\Gamma_{s +j-2q-2}^{\left[s -q-1\atop{ s -q-1\atop s + l} \right]}
   - 8\cdot 4^{2q+2}\cdot \Delta _{k}\Gamma_{s +j-2q-3}^{\left[s -q-1\atop{ s -q-1\atop s + l -1} \right]}\\
   &   +\sum_{i = 0}^{q}4^{2i}\cdot \Delta _{k}\omega_{s -i +(j-i)} (s-i,0, l+ i)  + 
   \sum_{i = 0}^{q}4^{2i+1}\cdot \Delta _{k}\omega_{s - i- 1+(j-i)} (s-i-1,1, l+ i) 
\end{align*}
\textbf{The two first cases $j=0,\;l\geq 0$ and $j=1,\;l\geq 0$  follows from \eqref{eq 9.1} and \eqref{eq 9.2}} \\
\textbf{The case  $ 2\leq j\leq s-4 $}\\
\textbf{Setting $q=j$ in the reduction formula \eqref{eq 11.3} we obtain using Lemma \ref{lem 6.3} and results in Lemma \ref{lem 7.4}}
   \begin{align*}
   &   \Delta _{k}\Gamma_{s+j}^{\left[s\atop{ s  \atop s +l } \right]}  -8\cdot\Delta _{k}\Gamma_{s+j-1}^{\left[s\atop{ s  \atop s +l -1 } \right]}
   - 4^{2j+2}\cdot \Delta _{k}\Gamma_{s -j -2}^{\left[s -j-1\atop{ s -j-1\atop s + l} \right]}\\
    &  = - 8\cdot 4^{2j+2}\cdot \Delta _{k}\Gamma_{s  -j -3}^{\left[s -j-1\atop{ s -j-1\atop s + l -1} \right]}
   +\sum_{i = 0}^{j}4^{2i}\cdot \Delta _{k}\omega_{s -i +(j-i)} (s-i,0, l+ i)  
   + \sum_{i = 0}^{j}4^{2i+1}\cdot \Delta _{k}\omega_{s - i - 1+(j-i)} (s-i-1,1, l+ i) \\
   & \text{we have }\\
 &   \Delta _{k}\omega_{s -i +(j-i)} (s-i,0, l+ i) = 0 \Leftrightarrow  s-i+3\leq s-i +(j-i)\leq 2\cdot(s-i) -1\Leftrightarrow i\leq j-3 \\
 & \Delta _{k}\omega_{s -i -1+(j-i)} (s-i-1,1, l+ i) = 0 \Leftrightarrow  s-i-1+3\leq s-i-1 +(j-i)\leq 1+ 2\cdot(s-i-1) -1\Leftrightarrow i\leq j-3 \\
 &\text{We then deduce}\\
  &   \Delta _{k}\Gamma_{s+j}^{\left[s\atop{ s  \atop s +l } \right]}  -8\cdot\Delta _{k}\Gamma_{s+j-1}^{\left[s\atop{ s  \atop s +l -1 } \right]}\\
  & =  \sum_{i = j-2}^{j}4^{2i}\cdot \Delta _{k}\omega_{s -i +(j-i)} (s-i,0, l+ i)  
   + \sum_{i = j-2}^{j}4^{2i+1}\cdot \Delta _{k}\omega_{s - i - 1+(j-i)} (s-i-1,1, l+ i) \\
   & = 2^{4(j-2)}\cdot \Delta _{k}\omega_{s -j +2 +2} (s- j+2,0, l+ j-2)   + 2^{4(j-1)}\cdot \Delta _{k}\omega_{s -j + 1+1} (s- j+1,0, l+ j-1) \\
    & +  2^{4j}\cdot \Delta _{k}\omega_{s -j } (s- j,0, l+ j) + 2^{4(j-2)+2}\cdot \Delta _{k}\omega_{s -j +1 +2} (s- j+1,1, l+ j-2) \\
 &   2^{4(j-1)+2}\cdot \Delta _{k}\omega_{s -j + 1} (s- j,1, l+ j-1) +  2^{4j +2}\cdot \Delta _{k}\omega_{s -j-1 } (s- j-1,1, l+ j)\\
    & = 2^{4j-8}\cdot2^{k + s-j+2} + (-3\cdot2^{k+(s-j+1)-1})\cdot2^{4j-4} + 2^{4j}\cdot2^{k + s -j -1}
     + 2^{4j-6}\cdot2^{k + s-j+1} + (-3\cdot2^{k+(s-j)-1})\cdot2^{4j-2}\\
     & + 2^{4j+2}\cdot2^{k +( s -j -1) -1} = 21\cdot2^{k+s+3j-6}  + 21\cdot2^{k+s+3j-5} = 63\cdot2^{k+s+3j-6}
  \end{align*}
  
  \textbf{The case  $ j = s-3 $}\\
  \textbf{Setting $q = s-3$ in the reduction formula \eqref{eq 11.3} we obtain similarly}
   \begin{align*}
    &   \Delta _{k}\Gamma_{2s-3}^{\left[s\atop{ s  \atop s +l } \right]}  -8\cdot\Delta _{k}\Gamma_{2s-4}^{\left[s\atop{ s  \atop s +l -1 } \right]}
   - 4^{2s-4}\cdot \Delta _{k}\Gamma_{1}^{\left[2\atop{ 2 \atop s + l} \right]}\\
    &  = - 8\cdot 4^{2s-4}\cdot \Delta _{k}\Gamma_{0}^{\left[2\atop{ 2 \atop s + l -1} \right]}
   +\sum_{i = 0}^{s - 3}4^{2i}\cdot \Delta _{k}\omega_{2(s-i) -3} (s-i,0, l+ i)  
   + \sum_{i = 0}^{s-3}4^{2i+1}\cdot \Delta _{k}\omega_{2(s-i-1)-2} (s-i-1,1, l+ i) \\
   & \text{we have }\\
 &   \Delta _{k}\omega_{2(s-i) -3} (s-i,0, l+ i) = 0 \Leftrightarrow  s-i+3\leq 2(s-i)-3\leq 2\cdot(s-i) -1\Leftrightarrow i\leq s-6 \\
 & \Delta _{k}\omega_{2(s-i-1)-2} (s-i-1,1, l+ i) = 0 \Leftrightarrow  s-i-1+3\leq 2(s-i-1)-2\leq 1+ 2\cdot(s-i-1) -1\Leftrightarrow i\leq s-6 \\
 &\text{We then deduce}\\
  &   \Delta _{k}\Gamma_{2s-3}^{\left[s\atop{ s  \atop s +l } \right]}  -8\cdot\Delta _{k}\Gamma_{2s-4}^{\left[s\atop{ s  \atop s +l -1 } \right]}\\
  & =  \sum_{i = s-5}^{s-3}4^{2i}\cdot \Delta _{k}\omega_{2(s-i) -3} (s-i,0, l+ i)  
   + \sum_{i = s-5}^{s-3}4^{2i+1}\cdot \Delta _{k}\omega_{2(s-i-1) -2} (s-i-1,1, l+ i) \\
   & = 2^{4(s-5)}\cdot \Delta _{k}\omega_{7} (5,0, l+ s-5) + 
    2^{4(s-4)}\cdot \Delta _{k}\omega_{5} (4,0, l+ s-4) +  2^{4(s-3)}\cdot \Delta _{k}\omega_{3 } (3,0, l+ s-3) \\
      & + 2^{4(s-5)+2}\cdot \Delta _{k}\omega_{6} (4,1, l+s-5) + 
    2^{4(s-4)+2}\cdot \Delta _{k}\omega_{4} (3,1, l+ s-4)\\
    & +  2^{4(s-3) +2}\cdot \Delta _{k}\omega_{2 } (2,1, l+ s-3)\\
     & = 2^{4s -20}\cdot2^{k + 5} + (-3\cdot2^{k+ 3})\cdot2^{4s -16} + 2^{4s-12}\cdot2^{k + 2}
      + 2^{4s - 18}\cdot2^{k + 4} + (-3\cdot2^{k+ 2})\cdot2^{4s -14} + 2^{4s-10}\cdot2^{k + 1}\\
    & = 21\cdot2^{k+4s -15}  + 21\cdot2^{k+4s -14} = 63\cdot2^{k+ 4s-15}
  \end{align*}
\textbf{The case  $ j = s-2 $}\\
  \textbf{Setting $q=s-3$ in the reduction formula \eqref{eq 11.3} we obtain similarly}
 \begin{align*}
  &   \Delta _{k}\Gamma_{2s-2}^{\left[s\atop{ s  \atop s +l } \right]}  -8\cdot\Delta _{k}\Gamma_{2s-3}^{\left[s\atop{ s  \atop s +l -1 } \right]}
   - 4^{2s-4}\cdot \Delta _{k}\Gamma_{2}^{\left[2\atop{ 2 \atop s + l} \right]}\\
    &  = - 8\cdot 4^{2s-4}\cdot \Delta _{k}\Gamma_{1}^{\left[2\atop{ 2 \atop s + l -1} \right]}
   +\sum_{i = 0}^{s - 3}4^{2i}\cdot \Delta _{k}\omega_{2(s-i) -2} (s-i,0, l+ i)  
   + \sum_{i = 0}^{s-3}4^{2i+1}\cdot \Delta _{k}\omega_{2(s-i-1)-1} (s-i-1,1, l+ i) \\
   & \text{we have }\\
 &   \Delta _{k}\omega_{2(s-i) -2} (s-i,0, l+ i) = 0 \Leftrightarrow  s-i+3\leq 2(s-i)-2\leq 2\cdot(s-i) -1\Leftrightarrow i\leq s-5 \\
 & \Delta _{k}\omega_{2(s-i-1)-1} (s-i-1,1, l+ i) = 0 \Leftrightarrow  s-i-1+3\leq 2(s-i-1)-1\leq 1+ 2\cdot(s-i-1) -1\Leftrightarrow i\leq s-5 \\
&   \Delta _{k}\Gamma_{2}^{\left[2\atop{ 2 \atop s + l} \right]} = 3\cdot2^{k+1}\\
  &\text{We then deduce}\\
  &   \Delta _{k}\Gamma_{2s-2}^{\left[s\atop{ s  \atop s +l } \right]}  -8\cdot\Delta _{k}\Gamma_{2s-3}^{\left[s\atop{ s  \atop s +l -1 } \right]}\\
  & = 4^{2s-4}\cdot \Delta _{k}\Gamma_{2}^{\left[2\atop{ 2 \atop s + l} \right]}  + \sum_{i = s-4}^{s-3}4^{2i}\cdot \Delta _{k}\omega_{2(s-i) -2} (s-i,0, l+ i)  
   + \sum_{i = s-4}^{s-3}4^{2i+1}\cdot \Delta _{k}\omega_{2(s-i-1) -1} (s-i-1,1, l+ i) \\
   & =  2^{4s-8}\cdot3\cdot2^{k+1} +  2^{4(s-4)}\cdot \Delta _{k}\omega_{6} (4,0, l+ s-4) +  2^{4(s-3)}\cdot \Delta _{k}\omega_{4 } (3,0, l+ s-3) \\
    & + 2^{4(s-4)+2}\cdot \Delta _{k}\omega_{5} (3,1, l+s-4) +  2^{4(s-3)+2}\cdot \Delta _{k}\omega_{3} (2,1, l+ s-3)\\
    & =  2^{4s-8}\cdot3\cdot2^{k+1}  +  2^{4s -16}\cdot2^{k + 4} + (-3\cdot2^{k+ 2})\cdot2^{4s -12} + 2^{4s-14}\cdot2^{k + 3}
    +  (-3\cdot2^{k+ 1})\cdot2^{4s -10} \\
    &  = 63\cdot2^{k+ 4s-12}
  \end{align*}
  \textbf{The case  $ j = s-1 $}\\
  \textbf{Setting $q=s-2$ in the reduction formula \eqref{eq 11.3} we obtain similarly}
   \begin{align*}
  &   \Delta _{k}\Gamma_{2s-1}^{\left[s\atop{ s  \atop s +l } \right]}  -8\cdot\Delta _{k}\Gamma_{2s-2}^{\left[s\atop{ s  \atop s +l -1 } \right]}
   - 4^{2s-2}\cdot \Delta _{k}\Gamma_{1}^{\left[1\atop{ 1 \atop s + l} \right]}\\
    &  = - 8\cdot 4^{2s-2}\cdot \Delta _{k}\Gamma_{0}^{\left[1\atop{ 1 \atop s + l -1} \right]}
   +\sum_{i = 0}^{s - 2}4^{2i}\cdot \Delta _{k}\omega_{2(s-i) -1} (s-i,0, l+ i)  
   + \sum_{i = 0}^{s- 2}4^{2i+1}\cdot \Delta _{k}\omega_{2(s-i-1)} (s-i-1,1, l+ i) \\
   & \text{we have }\\
 &   \Delta _{k}\omega_{2(s-i) -1} (s-i,0, l+ i) = 0 \Leftrightarrow  s-i+3\leq 2(s-i)-1\leq 2\cdot(s-i) -1\Leftrightarrow i\leq s-4 \\
 & \Delta _{k}\omega_{2(s-i-1)} (s-i-1,1, l+ i) = 0 \Leftrightarrow  s-i-1+3\leq 2(s-i-1)\leq 1+ 2\cdot(s-i-1) -1\Leftrightarrow i\leq s-4,\; s-i-1\geq 4 \\
&   \Delta _{k}\Gamma_{1}^{\left[1\atop{ 1 \atop s + l} \right]} = 3\cdot2^{k}\\
  &\text{We then deduce}\\
  &   \Delta _{k}\Gamma_{2s-1}^{\left[s\atop{ s  \atop s +l } \right]}  -8\cdot\Delta _{k}\Gamma_{2s-2}^{\left[s\atop{ s  \atop s +l -1 } \right]}\\
  & = 4^{2s-2}\cdot \Delta _{k}\Gamma_{1}^{\left[1\atop{ 1 \atop s + l} \right]}  + \sum_{i = s-3}^{s-2}4^{2i}\cdot \Delta _{k}\omega_{2(s-i) -1} (s-i,0, l+ i)  
   + \sum_{i = s-4}^{s-2}4^{2i+1}\cdot \Delta _{k}\omega_{2(s-i-1) } (s-i-1,1, l+ i) \\
    & =  2^{4s-4}\cdot3\cdot2^{k} +  2^{4(s-3)}\cdot \Delta _{k}\omega_{5} (3,0, l+ s-3) +  2^{4(s-2)}\cdot \Delta _{k}\omega_{3 } (2,0, l+ s-2) \\
    & + 2^{4(s-4)+2}\cdot \Delta _{k}\omega_{6} (3,1, l+s-4) +  2^{4(s-3)+2}\cdot \Delta _{k}\omega_{4} (2,1, l+ s-3) +  2^{4(s-2)+2}\cdot \Delta _{k}\omega_{2} (1,1, l+s-2)\\
     & =  2^{4s-4}\cdot3\cdot2^{k}  +  2^{4s -12}\cdot2^{k + 3} + (-3\cdot2^{k+ 1})\cdot2^{4s - 8}+  2^{4s-14}\cdot 0  +
    2^{4s -10}\cdot2^{k + 2} +  (-3\cdot2^{k})\cdot2^{4s - 6} \\
    &  = 2^{k+4s-9}\cdot[1 + 3\cdot2^{5}  -3\cdot2^{2} + 2 - 3\cdot2^{3}]  = 63\cdot2^{k+ 4s- 9}
  \end{align*}
\end{proof}
      \subsection{\textbf{Computation of $ \Delta _{k}\Gamma_{s+ m+ j}^{\left[s\atop{ s+m \atop s+m+l} \right]} -   8\cdot\Delta _{k}\Gamma_{s+m + j -1}^{\left[s\atop{ s +m \atop s+m+ l-1} \right]}  \quad for\; 0\leq j \leq s-1 ,\; m\geq 1,\;l\geq 0 $}}
\label{subsec 5} 
\begin{lem}
\label{lem 11.4}
We have for $ m\geq 1,\; l\geq 0 $ the following reduction formula  \\
\begin{align}
& \Delta _{k}\Gamma_{s+m+j}^{\left[s\atop{ s+m \atop s+m+l} \right]} -   8\cdot\Delta _{k}\Gamma_{s+m+ j -1}^{\left[s\atop{ s+m \atop s+m+l-1} \right]} \label{eq 11.5} \\
& =  2^{2m}\cdot\left[\Delta _{k}\Gamma_{s+j}^{\left[s\atop{ s \atop s+m+l} \right]} -   8\cdot\Delta _{k}\Gamma_{s+ j -1}^{\left[s\atop{ s \atop s+m+l-1} \right]}\right ]
  + \sum_{i=0}^{m-1}2^{2i}\cdot\Delta _{k}\omega_{s+m+j-i } (s,m-i, l+ i) \nonumber
\end{align}
where 
\begin{equation}
\label{eq 11.6}
\sum_{i=0}^{m-1}2^{2i}\cdot\Delta _{k}\omega_{s+m+j-i } (s,m-i, l+ i) =
 \begin{cases}
(-3)\cdot2^{k+s-1} & \text{if   } j = 0,\; m = 1,\; l\geq 0 \\
(-5)\cdot2^{k+s +2m-4} & \text{if   } j = 0,\; m \geq 2,\; l\geq 0 \\
2^{k+s +2m -2} & \text{if   } j = 1, \; m\geq 1,\; l\geq 0 \\
 0  &   \text{if  }\;  2\leq j\leq s-1,\; m\geq 1, \; l\geq 0.
 \end{cases}
\end{equation}

We get from \eqref{eq 11.4} applying the  reduction formula \eqref{eq 11.5}:
\begin{equation}
\label{eq 11.7}
\Delta _{k}\Gamma_{s+m+j}^{\left[s\atop{ s+m \atop s+m+l} \right]} -   8\cdot\Delta _{k}\Gamma_{s+m+ j -1}^{\left[s\atop{ s+m \atop s+m+l-1} \right]}    =
 \begin{cases}
9\cdot2^{k+s-1} & \text{if   } j = 0,\; m = 1,\; l\geq 0 \\
19\cdot2^{k+s +2m-4} & \text{if   } j = 0,\; m \geq 2,\; l\geq 0 \\
31\cdot2^{k+s +2m -2} & \text{if   } j = 1, \; m\geq 1,\; l\geq 0 \\
 63\cdot2^{k+s+2m +3j-6} &   \text{if  }\;  2\leq j\leq s-1,\; m\geq 0, \; l\geq 0.
 \end{cases}
\end{equation}
\end{lem}
\begin{proof}
Let $q $ be a rational integer between one and $m-1.$
 From \eqref{eq 8.3} we obtain  successively the following equations

  \begin{align*}
 &  \Delta _{k}\Gamma_{s+m+j}^{\left[s\atop{ s +m\atop s+m+l } \right]}
   - 4\cdot \Delta _{k}\Gamma_{s+m+j -1}^{\left[s\atop{ s +m-1\atop s+m+l } \right]}  = 
   8\cdot\left[  \Delta _{k}\Gamma_{s+m+j -1}^{\left[s\atop{ s+m \atop s +m+l-1} \right]}
   - 4\cdot \Delta _{k}\Gamma_{s+m+j -2}^{\left[s\atop{ s +m-1\atop s +m+l-1} \right]}\right] + \Delta _{k}\omega_{s+m+j} (s,m,l) \\
  &  4\cdot\left[  \Delta _{k}\Gamma_{s+m+j-1}^{\left[s \atop{ s +m-1\atop s+m+l } \right]}
   - 4\cdot \Delta _{k}\Gamma_{s+m+j -2}^{\left[s\atop{ s+m-2 \atop s +m+l } \right]}\right]  = 
4\cdot\left[   8\cdot\Big[  \Delta _{k}\Gamma_{s+m+j -2}^{\left[s\atop{ s+m-1 \atop s +m+l-1 } \right]}
   - 4\cdot \Delta _{k}\Gamma_{s+m+j -3}^{\left[s\atop{ s +m-2\atop s +m+l-1 } \right]}\Big] + \Delta _{k}\omega_{s+m+j-1} (s,m-1, l+1) \right]\\
    &  4^{2}\cdot\left[   \Delta _{k}\Gamma_{s+m+j-2}^{\left[s \atop{ s +m-2 \atop s+m+l } \right]}
   - 4\cdot \Delta _{k}\Gamma_{s+m+j -3}^{\left[s\atop{ s +m-3 \atop s +m+l} \right]}\right]  = 
 4^{2}\cdot\left[  8\cdot\Big[  \Delta _{k}\Gamma_{s+m+j -3}^{\left[s \atop{ s +m-2 \atop s +m+l-1} \right]}
   - 4\cdot \Delta _{k}\Gamma_{s+m+j - 4}^{\left[s \atop{ s +m-3 \atop s +m+l-1 } \right]}\Big] + \Delta _{k}\omega_{s+m+j-2} (s,m-2,l+2)\right] \\
     &  \vdots       \\
   &  4^{q}\cdot\left[   \Delta _{k}\Gamma_{s+m+j - q}^{\left[s \atop{ s +m - q) \atop s+m+l } \right]}
   - 4\cdot \Delta _{k}\Gamma_{s+m+j - (q+1)}^{\left[s\atop{ s +m - (q+1) \atop s +m+l} \right]}\right]\\
    & =  4^{q}\cdot\left[  8\cdot\Big[  \Delta _{k}\Gamma_{s+m+j -(q+1)}^{\left[s \atop{ s +m- q \atop s +m+l-1} \right]}
   - 4\cdot \Delta _{k}\Gamma_{s+m+j - (q+2)}^{\left[s \atop{ s +m- (q+1) \atop s +m+l-1 } \right]}\Big] + \Delta _{k}\omega_{s+m +j-q} (s,m- q,l+ q)\right]  
 \end{align*}
 By summing  the above equations we obtain after reduction the following formula  \\ 
     \begin{align*}
    &  \Delta _{k}\Gamma_{s+m+j}^{\left[s\atop{ s+m  \atop s+m +l } \right]}
    -8\cdot\Delta _{k}\Gamma_{s+m+j-1}^{\left[s\atop{ s +m \atop s +m+ l -1 } \right]}\\
     & =   4^{q+1}\cdot \Delta _{k}\Gamma_{s +m+j-(q+1)}^{\left[s \atop{ s +m-(q+1)\atop s +m+ l} \right]}
   - 8\cdot 4^{q+1}\cdot \Delta _{k}\Gamma_{s + m+j -(q+2)}^{\left[s \atop{ s  +m - (q+1)\atop s + m+l-1} \right]}
      +\sum_{i = 0}^{q}4^{i}\cdot \Delta _{k}\omega_{s +m+j-i} (s,m-i, l+ i) \\
      & \textbf{Setting $ q = m-1$ in the above formula we get the reduction formula \eqref{eq 11.5}} 
     \end{align*} 
     \textbf{Proof of \eqref{eq 11.6} } 
 The proof is based on  results in Lemma \ref{lem 7.4} \\
     \begin{align*}
    & \underline{The\; case\; j=0,\; m=1,\; l\geq 0} \\
& \quad \sum_{i=0}^{m-1}2^{2i}\cdot\Delta _{k}\omega_{s+m+j-i } (s,m-i, l+ i) 
=\Delta _{k}\omega_{s+1 } (s,1, l) = (-3)\cdot2^{k+s-1} \\
& \\
&  \underline{The\; case\; j=0,\; m\geq 2,\; l\geq 0} \\
& \quad \sum_{i=0}^{m-1}2^{2i}\cdot\Delta _{k}\omega_{s+m+j-i } (s,m-i, l+ i)
 = \quad \sum_{i=m-2}^{m-1}2^{2i}\cdot\Delta _{k}\omega_{s+m-i } (s,m-i, l+ i) \\
 & = 2^{2(m-2)}\cdot\Delta _{k}\omega_{s+2 } (s,2, l+m-2) + 2^{2(m-1)}\cdot\Delta _{k}\omega_{s+1 } (s,1, l+ m-1) \\
 & = 2^{2(m-2)}\cdot2^{k+s} +  2^{2(m-1)}\cdot(-3\cdot2^{k+s-1}) = (-5)\cdot2^{k+s+2m-4}\\
 & \text{Since \quad $\Delta _{k}\omega_{s+m-i } (s,m-i, l+ i) = 0 \Longleftrightarrow  i\leq m-3 $}\\
 & \\
   &  \underline{The\; case\; j=1,\; m\geq 1,\; l\geq 0} \\
   & \quad \sum_{i=0}^{m-1}2^{2i}\cdot\Delta _{k}\omega_{s+m+j-i } (s,m-i, l+ i) = 
   \quad \sum_{i=m-1}^{m-1}2^{2i}\cdot\Delta _{k}\omega_{s+m+1-i } (s,m-i, l+ i)\\
    &  = 2^{2(m-1)}\cdot\Delta _{k}\omega_{s+2} (s,1, l+ m-1) =  2^{2(m-1)}\cdot2^{k+s}=2^{k+s+2m-2}\\
    & \text{Since \quad $\Delta _{k}\omega_{s+m+1-i } (s,m-i, l+ i) = 0 \Longleftrightarrow  i\leq m-2 $}\\ 
 & \\
  &  \underline{The\; case\; 2\leq j\leq s-1,\; m\geq 1,\; l\geq 0} \\
   & \quad \sum_{i=0}^{m-1}2^{2i}\cdot\Delta _{k}\omega_{s+m+j-i } (s,m-i, l+ i) = 0 \\
    & \text{Since \quad $\Delta _{k}\omega_{s+m+j-i } (s,m-i, l+ i) = 0 \Longleftrightarrow 0\leq  i\leq m-1 $} 
     \end{align*}
\textbf{Proof of \eqref{eq 11.7} } 
 
 We deduce easily \eqref{eq 11.7} from the reduction formula \eqref{eq 11.5} using \eqref{eq 11.6}
 and \eqref{eq 11.4} with $l \longrightarrow m+l$
  \end{proof}
 \subsection{\textbf{Computation of $ \Delta _{k}\Gamma_{s+m+1}^{\left[s\atop{ s+m\atop s+m} \right]} \quad for\;  m\geq 0 $}}
\label{subsec 6} 
 \begin{lem}
 \label{lem 11.5} We have \\
 \begin{equation}
\label{eq 11.8}
  \Delta _{k}\Gamma_{s+m+1}^{\left[s\atop{ s+m\atop s+m} \right]} =
 \begin{cases}
147\cdot 2^{k+s-1} & \text{if   } m = 0 \\
315\cdot 2^{k+s} & \text{if   } m = 1 \\
  21\cdot2^{k+s+3m -2} + 273\cdot2^{k+s+2m- 2} & \text{if }  m\geq 2\;(holds \; for \; m\geq 0 )\\
  \end{cases}
  \end{equation}
 \end{lem}
 \begin{proof}
 From \eqref{eq 8.3} and, applying  \eqref{eq 11.7} in the case $j=1,\;m\geq 1,\;l=0:$ 

\begin{align*}
  &  \Delta _{k}\Gamma_{s+m+1}^{\left[s\atop{ s+m  \atop s+m  } \right]}
    -8\cdot\Delta _{k}\Gamma_{s+m}^{\left[s\atop{ s +m -1\atop s +m } \right]} \\
     & =  4\cdot\Big[ \Delta _{k}\Gamma_{s+m}^{\left[s\atop{ s+m -1 \atop s+m  } \right]}
    -8\cdot\Delta _{k}\Gamma_{s+m -1}^{\left[s\atop{ s +m -1\atop s +m -1} \right]}\Big]
     +  \Delta _{k}\omega_{s + m+1 } (s, m, 0) =  31\cdot2^{k+ s +2m -2}\\
   &   \Delta _{k}\Gamma_{s+m}^{\left[s\atop{ s+m -1 \atop s+m  } \right]}
    = 8\cdot\Delta _{k}\Gamma_{s+m -1}^{\left[s\atop{ s +m -1\atop s +m -1} \right]} 
    + {1\over 4}\cdot\big( 31\cdot2^{k+ s +2m -2} - \Delta _{k}\omega_{s + m+1 } (s, m, 0)\big) \\
    & \text{Combining the above equations we have by \eqref{eq 11.2} with $m\longrightarrow m-1$
     and \eqref{eq 7.48} in the case $ m\geq 2 $ }\\
    &  \Delta _{k}\Gamma_{s+m+1}^{\left[s\atop{ s+m  \atop s+m  } \right]} = 
    8\cdot\Delta _{k}\Gamma_{s+m}^{\left[s\atop{ s +m -1\atop s +m } \right]} + 31\cdot2^{k+ s +2m -2}\\
     & =    8\cdot\Big( 8\cdot\Delta _{k}\Gamma_{s+m -1}^{\left[s\atop{ s +m -1\atop s +m -1} \right]} 
    + {1\over 4}\cdot\big( 31\cdot2^{k+ s +2m -2} - \Delta _{k}\omega_{s + m+1 } (s, m, 0)\big)\Big ) + 31\cdot2^{k+ s +2m -2} \\
    & = 2^{6}\cdot\big[ 21\cdot2^{k+s+3(m-1) -5} + 45\cdot2^{k+s+2(m-1) - 4}\big] +  93\cdot2^{k+s+2m-2}\quad \text{since $\Delta _{k}\omega_{s + m+1 } (s, m, 0)=0$} \\
    & =  21\cdot2^{k+s+3m -2} + 273\cdot2^{k+s+2m- 2}\\
    & \text{similarly in the case $m=1$}\\
      &  \Delta _{k}\Gamma_{s+2}^{\left[s\atop{ s+1  \atop s+1 } \right]} = 
    8\cdot\Delta _{k}\Gamma_{s+1}^{\left[s\atop{ s \atop s +1 } \right]} + 31\cdot2^{k+ s}\\
     & =    8\cdot\Big( 8\cdot\Delta _{k}\Gamma_{s}^{\left[s\atop{ s \atop s } \right]} 
    + {1\over 4}\cdot\big( 31\cdot2^{k+ s } - \Delta _{k}\omega_{s + 2 } (s, 1, 0)\big)\Big ) + 31\cdot2^{k+ s } \\
    & = 2^{6}\cdot7\cdot2^{k+s-1} +  91\cdot2^{k+s} =    315\cdot2^{k+s}\quad \text{since $\Delta _{k}\omega_{s +2 } (s, 1, 0)= 2^{k+s}$} 
       \end{align*}

 \end{proof}
   \subsection{\textbf{Computation of $ \Delta _{k}\Gamma_{2s+j}^{\left[s\atop{ s \atop s+l} \right]} -   8\cdot\Delta _{k}\Gamma_{2s+j -1}^{\left[s\atop{ s \atop s+l-1} \right]}  \quad for\; 0\leq j \leq l $}}
\label{subsec 7} 
 \begin{lem}
\label{lem 11.6}
 We have for $ l\geq 3 $ the following reduction formula holding for $0\leq j\leq l.$\\
\begin{align}
 &   \Delta _{k}\Gamma_{2s+j}^{\left[s\atop{ s  \atop s +l } \right]}
    -8\cdot\Delta _{k}\Gamma_{2s+j-1}^{\left[s\atop{ s  \atop s +l -1 } \right]} 
     =   4^{2s-2}\cdot \Delta _{k}\Gamma_{2+j}^{\left[1\atop{ 1 \atop s + l} \right]}
   - 8\cdot 4^{2s-2}\cdot \Delta _{k}\Gamma_{1+j}^{\left[ 1\atop{ 1 \atop s + l -1} \right]}\label{eq 11.9} \\
   &   +\sum_{i = 0}^{s - 2}4^{2i}\cdot \Delta _{k}\omega_{2(s-i)+j} (s-i,0, l+ i)  + 
   \sum_{i = 0}^{s - 2}4^{2i+1}\cdot \Delta _{k}\omega_{2(s-i-1)+ (j+1)} (s-i-1,1, l+ i) \nonumber
\end{align}
 where 
 \begin{align}
 \label{eq 11.10}
   \Delta _{k}\Gamma_{2+j}^{\left[1\atop{ 1 \atop s + l} \right]}
   - 8\cdot \Delta _{k}\Gamma_{1+j}^{\left[ 1\atop{ 1 \atop s + l -1} \right]} = 
  \begin{cases}
 3\cdot2^{2k} - 9\cdot2^{k} & \text{if  } j = 0, \\
  -15\cdot2^{2k}  +3\cdot2^{k+1}   & \text{if  } j = 1, \\
   -9\cdot2^{2k +2j-2}      & \text{if  } 2\leq j\leq l, 
   \end{cases}
 \end{align}
\begin{align}
\label{eq 11.11}
\sum_{i = 0}^{s - 2}4^{2i}\cdot \Delta _{k}\omega_{2(s-i)+j} (s-i,0, l+ i) =
 \begin{cases}
     - 2^{2k +4s -4} + 2^{2k +2s -2} + 2^{k +4s -6}  & \text{if  } j = 0, \\
   5\cdot 2^{2k +4s -4} -5\cdot 2^{2k +2s -2} & \text{if  } j = 1, \\
     3\cdot 2^{2k +4s +2j-6} -3\cdot 2^{2k +2s +2j-4}   & \text{if  } 2\leq j\leq l-1,  \\
   3\cdot 2^{2k +4s +2l-6} - 27\cdot 2^{2k +2s +2l-4}  & \text{if  } j = l,  
    \end{cases}
\end{align}

\begin{align}
\label{eq 11.12}
\sum_{i = 0}^{s - 2}4^{2i+1}\cdot \Delta _{k}\omega_{2(s-i-1)+j+1} (s-i-1,1, l+ i) =
 \begin{cases}
     - 2^{2k +4s -3} + 2^{2k +2s -1} + 2^{k +4s -5}  & \text{if  } j = 0, \\
   5\cdot 2^{2k +4s -3} -5\cdot 2^{2k +2s -1} & \text{if  } j = 1, \\
     3\cdot 2^{2k +4s +2j-5} -3\cdot 2^{2k +2s +2j-3}   & \text{if  } 2\leq j\leq l-1,  \\
   3\cdot 2^{2k +4s +2l-5} - 27\cdot 2^{2k +2s +2l-3}  & \text{if  } j = l.  
    \end{cases}
\end{align}
Combining \eqref{eq 11.9}, \eqref{eq 11.10}, \eqref{eq 11.11} and  \eqref{eq 11.12} we obtain \\

\begin{equation}
\label{eq 11.13}
 \Delta _{k}\Gamma_{2s+j}^{\left[s\atop{ s \atop s+l} \right]} -   8\cdot\Delta _{k}\Gamma_{2s+j -1}^{\left[s\atop{ s \atop s+l-1} \right]}    =
 \begin{cases}
 3\cdot2^{2k+2s-2} - 33\cdot2^{k+4s-6} & \text{if  } j = 0 \\
  -15\cdot2^{2k+2s-2}  +3\cdot2^{k+4s-3}   & \text{if  } j = 1 \\
   -9\cdot2^{2k+2s +2j-4}      & \text{if  } 2\leq j\leq l-1  \\
     -81\cdot2^{2k+2s +2j-4}      & \text{if  } j = l   
    \end{cases}
\end{equation}
 \end{lem}
 \begin{proof}
  \underline{Proof of \eqref{eq 11.9}} \\
 Let q be a rational integer between zero and s-2.
 From \eqref{eq 8.3} we obtain successively the following equations. \\
 
  \small   
  \begin{align*}
 &  \Delta _{k}\Gamma_{2s+j}^{\left[s\atop{ s \atop s+l } \right]}
   - 4\cdot \Delta _{k}\Gamma_{2s+j -1}^{\left[s -1\atop{ s \atop s+l } \right]}  = 
   8\cdot\left[  \Delta _{k}\Gamma_{2s+j -1}^{\left[s\atop{ s \atop s + l-1} \right]}
   - 4\cdot \Delta _{k}\Gamma_{2s+j -2}^{\left[s\atop{ s -1\atop s + l-1} \right]}\right] + \Delta _{k}\omega_{2s+j} (s,0,l) \\
     &  4\cdot\left[  \Delta _{k}\Gamma_{2s+j-1}^{\left[s -1\atop{ s \atop s+l } \right]}
   - 4\cdot \Delta _{k}\Gamma_{2s+j -2}^{\left[s -1\atop{ s -1 \atop s +l } \right]}\right]  = 
4\cdot\left[   8\cdot\Big[  \Delta _{k}\Gamma_{2s+j -2}^{\left[s-1\atop{ s  \atop s +l-1 } \right]}
   - 4\cdot \Delta _{k}\Gamma_{2s+j -3}^{\left[s -1\atop{ s -1\atop s +l-1 } \right]}\Big] + \Delta _{k}\omega_{2s+j-1} (s-1,1, l) \right]\\
     &  4^{2}\cdot\left[   \Delta _{k}\Gamma_{2s+j-2}^{\left[s -1\atop{ s -1 \atop s+l } \right]}
   - 4\cdot \Delta _{k}\Gamma_{2s+j -3}^{\left[s\atop{ s -2 \atop s +l} \right]}\right]  = 
 4^{2}\cdot\left[  8\cdot\Big[  \Delta _{k}\Gamma_{2s+j -3}^{\left[s -1\atop{ s -1 \atop s +l-1} \right]}
   - 4\cdot \Delta _{k}\Gamma_{2s+j - 4}^{\left[s -2\atop{ s -1 \atop s +l-1 } \right]}\Big] + \Delta _{k}\omega_{2s+j-2} (s-1,0,l+1)\right] \\
     &  \vdots       \\
      &  4^{2q}\cdot\left[   \Delta _{k}\Gamma_{2s+j -2q}^{\left[s -q\atop{ s -q \atop s+l } \right]}
   - 4\cdot \Delta _{k}\Gamma_{2s+j -2q-1}^{\left[s -q-1\atop{ s -q \atop s +l} \right]}\right]\\
    & =  4^{2q}\cdot\left[  8\cdot\Big[  \Delta _{k}\Gamma_{2s+j -2q-1}^{\left[s -q\atop{ s -q \atop s +l-1} \right]}
   - 4\cdot \Delta _{k}\Gamma_{2s+j -2q-2}^{\left[s -q-1\atop{ s -q \atop s +l -1 } \right]}\Big] + \Delta _{k}\omega_{2s+j- 2q} (s -q,0,l+ q)\right] \\
    &  4^{2q +1}\cdot\left[   \Delta _{k}\Gamma_{2s+j -2q-1}^{\left[s -q-1\atop{ s -q \atop s+l } \right]}
   - 4\cdot \Delta _{k}\Gamma_{2s+j -2q-2}^{\left[s -q-1\atop{ s -q-1 \atop s +l} \right]}\right]\\
    & =  4^{2q+1}\cdot\left[  8\cdot\Big[  \Delta _{k}\Gamma_{2s+j -2q-2}^{\left[s -q-1\atop{ s -q \atop s +l-1} \right]}
   - 4\cdot \Delta _{k}\Gamma_{2s+j -2q-3}^{\left[s -q-1\atop{ s -q-1 \atop s +l -1 } \right]}\Big] + \Delta _{k}\omega_{2s+j- 2q-1} (s -q-1,1,l+ q)\right] 
   \end{align*}
 By summing  the above equations we obtain after reduction the following formula \\ 
  \begin{align*}
   &   \Delta _{k}\Gamma_{2s+j}^{\left[s\atop{ s  \atop s +l } \right]}
    -8\cdot\Delta _{k}\Gamma_{2s+j-1}^{\left[s\atop{ s  \atop s +l -1 } \right]}\\
    & =   4^{2q+2}\cdot \Delta _{k}\Gamma_{2(s-q-1) +j}^{\left[s -q-1\atop{ s -q-1\atop s + l} \right]}
   - 8\cdot 4^{2q+2}\cdot \Delta _{k}\Gamma_{2(s-q-1) + (j-1)}^{\left[s -q-1\atop{ s -q-1\atop s + l -1} \right]}\\
   &   +\sum_{i = 0}^{q}4^{2i}\cdot \Delta _{k}\omega_{2(s-i)+j} (s-i,0, l+ i)  + 
   \sum_{i = 0}^{q}4^{2i+1}\cdot \Delta _{k}\omega_{2(s-i-1)+ (j+1)} (s-i-1,1, l+ i) \\
   &  \textbf{Setting q = s - 2 in the above formula we obtain the reduction formula \eqref{eq 11.9}}\\
   & 
   \end{align*}
   \underline{Proof of \eqref{eq 11.10}} 
     Let $2+j < \inf(k, 2+s+l)$.  
  We obtain from Corollary 3.10 for n = 2 applying Theorem 3.1 (see [2] section 3) \\
  \begin{equation}
  \label{eq 11.14}
 \Gamma _{2+j}^{\left[\stackrel{1}{\stackrel{1}{s+l}}\right]\times k}= \begin{cases}
2^{2k}+15\cdot2^k +158  &  \text{if  }  j = 0  \\
3\cdot2^{2k+ 2j - 2 } + 63\cdot[2^{k +3j -2} +5\cdot2^{4j-1}] & \text{if   } 1\leq j\leq l \\
\end{cases}
\end{equation}
\begin{equation}
\label{eq 11.15}
 \Gamma _{1+j}^{\left[\stackrel{1}{\stackrel{1}{s+l-1}}\right]\times k}= \begin{cases}
3\cdot2^k + 9   &  \text{if  }    j = 0 \\
2^{2k}+15\cdot2^k +158  &  \text{if  }   j = 1  \\
3\cdot2^{2k+ 2j -4 } + 63\cdot[2^{k +3j -5} +5\cdot2^{4j-5}] & \text{if   } 2\leq j\leq l \\
\end{cases}
\end{equation}
\textbf{We deduce \eqref{eq 11.10} from  \eqref{eq 11.14} and  \eqref{eq 11.15}.}
\begin{align*}
&
\end{align*}
 \underline{Proof of \eqref{eq 11.11}}
The proof is based on results in  Lemma \ref{lem 7.4} \\
The case j = 0, \\
\begin{align*}
 &  \Delta _{k}\omega_{2(s-i)} (s-i,0, l+ i) = -3\cdot2^{2k +2(s-i)-2}\quad  \text{if $ 0\leq i\leq s -3 $},\\
 &  \Delta _{k}\omega_{4} (2,0, l+ (s-2)) = -3\cdot2^{2k +2} +2^{k+2}   \quad  \text{if  i = s -2 }.\\
 & \text{It follows that}\\
 &  \sum_{i = 0}^{s - 2}4^{2i}\cdot \Delta _{k}\omega_{2(s-i)} (s-i,0, l+ i) =  \sum_{i = 0}^{s - 2}2^{4i}\cdot( -3\cdot2^{2k +2(s-i)-2})  +  2^{k+2}\cdot2^{4(s-2)}\\
 & = - 2^{2k +4s -4} + 2^{2k +2s -2} + 2^{k +4s -6}. 
\end{align*}

The case j = 1, \\
\begin{align*}
 &  \Delta _{k}\omega_{2(s-i)+1} (s-i,0, l+ i) = 15\cdot2^{2k +2(s-i)-2}\quad  \text{if $ 0\leq i\leq s -2 $},\\
 & \text{It follows that}\\
 &  \sum_{i = 0}^{s - 2}4^{2i}\cdot \Delta _{k}\omega_{2(s-i)+1} (s-i,0, l+ i) =  \sum_{i = 0}^{s - 2}2^{4i}\cdot( 15\cdot2^{2k +2(s-i)-2}) \\
 & = 5\cdot2^{2k +4s -4} -5\cdot2^{2k +2s -2}.
\end{align*}

The case $2\leq j\leq l-1$, \\
\begin{align*}
 &  \Delta _{k}\omega_{2(s-i)+j} (s-i,0, l+ i) = 9\cdot2^{2k +2(s-i)+2j-4}\quad  \text{if $ 0\leq i\leq s -2 $},\\
 & \text{It follows that}\\
 &  \sum_{i = 0}^{s - 2}4^{2i}\cdot \Delta _{k}\omega_{2(s-i)+j} (s-i,0, l+ i) =  \sum_{i = 0}^{s - 2}2^{4i}\cdot( 9\cdot2^{2k +2(s-i) +2j-4}) \\
 & = 3\cdot2^{2k +4s +2j-6} -3\cdot2^{2k +2s +2j-4}.
\end{align*}

The case $j = l$, \\
\begin{align*}
 &  \Delta _{k}\omega_{2(s-i)+l} (s-i,0, l+ i) = 9\cdot2^{2k +2(s-i)+2l-4}\quad  \text{if $ 1\leq i\leq s -2 $},\\
  &  \Delta _{k}\omega_{2s+l} (s,0, l) = -15\cdot2^{2k +2s+2l-4}\quad  \text{if $ i = 0 $},\\
  & \text{It follows that}\\
 &  \sum_{i = 0}^{s - 2}4^{2i}\cdot \Delta _{k}\omega_{2(s-i)+j} (s-i,0, l+ i) =  \sum_{i = 1}^{s - 2}2^{4i}\cdot( 9\cdot2^{2k +2(s-i) +2l-4}) +( -15\cdot2^{2k +2s+2l-4}) \\
 & = 3\cdot2^{2k +4s +2l-6} -27\cdot2^{2k +2s +2l-4}.
\end{align*}

\underline{Proof of \eqref{eq 11.12}}

Similarly to the proof of \eqref{eq 11.11} applying \eqref{eq 11.14}, \eqref{eq 11.15} and Lemma \ref{lem 1.22}.

The case j = 0, \\
\begin{align*}
 &  \Delta _{k}\omega_{2(s-i-1)+1} (s-i-1,1, l+ i) = -3\cdot2^{2k +2(s-i-1)-1}\quad  \text{if $ 0\leq i\leq s -3 $},\\
 &  \Delta _{k}\omega_{3} (1,1, l+ (s-2)) = \Delta _{k}\Gamma_{3}^{\left[ 1\atop{ 2 \atop s + l } \right]}
  -4\cdot \Delta _{k}\Gamma_{2}^{\left[ 1\atop{ 1 \atop s + l } \right]} - 8\cdot \Delta _{k}\Gamma_{2}^{\left[ 1\atop{ 2 \atop s + l -1} \right]}
  +32\cdot \Delta _{k}\Gamma_{1}^{\left[ 1\atop{ 1 \atop s + l -1} \right]}\\
 & =  -3\cdot2^{2k +1} +2^{k+1}   \quad  \text{if  i = s -2 }.\\
 & \text{It follows that}\\
 &  \sum_{i = 0}^{s - 2}4^{2i+1}\cdot \Delta _{k}\omega_{2(s-i-1)+1} (s-i-1,1, l+ i) =  \sum_{i = 0}^{s - 2}2^{4i+2}\cdot( -3\cdot2^{2k +2(s-i-1)-1})  +  2^{k+1}\cdot2^{4(s-2)+2}\\
 & = - 2^{2k +4s -3} + 2^{2k +2s -1} + 2^{k +4s -5}. 
\end{align*}

The case j = 1, \\
\begin{align*}
 &  \Delta _{k}\omega_{2(s-i-1)+2} (s-i-1,1, l+ i) = 15\cdot2^{2k +2(s-i-1)-1}\quad  \text{if $ 0\leq i\leq s -2 $},\\
  & \Delta _{k}\omega_{4} (1,1, l + s-2) =  \Delta _{k}\Gamma_{4}^{\left[ 1\atop{ 2 \atop s + l } \right]}
  -4\cdot \Delta _{k}\Gamma_{3}^{\left[ 1\atop{ 1 \atop s + l } \right]} - 8\cdot \Delta _{k}\Gamma_{3}^{\left[ 1\atop{ 2 \atop s + l -1} \right]}
  +32\cdot \Delta _{k}\Gamma_{2}^{\left[ 1\atop{ 1 \atop s + l -1} \right]}\\
  & = [9\cdot2^{2k+1} +105\cdot2^{k+3}] -4\cdot[9\cdot2^{2k} +63\cdot2^{k+1}] -8\cdot[3\cdot2^{2k+1} + 102\cdot2^{k}] +32\cdot[3\cdot2^{2k} +15\cdot2^{k}] = 15\cdot2^{2k+1}\\
 & \text{It follows that}\\
 &  \sum_{i = 0}^{s - 2}4^{2i+1}\cdot \Delta _{k}\omega_{2(s-i-1)+2} (s-i-1,1, l+ i) =  \sum_{i = 0}^{s - 2}2^{4i+2}\cdot( 15\cdot2^{2k +2(s-i-1)-1}) \\
 & = 5\cdot2^{2k +4s -3} -5\cdot2^{2k +2s -1}.
\end{align*}

The case $2\leq j\leq l-1$, \\
\begin{align*}
 &  \Delta _{k}\omega_{2(s-i-1)+j+1} (s-i-1,1, l+ i) = 9\cdot2^{2k +2(s-i)+2j-5}\quad  \text{if $ 0\leq i\leq s -2 $},\\
  & \Delta _{k}\omega_{j+3} (1,1, l + s-2) =  \Delta _{k}\Gamma_{j+3}^{\left[ 1\atop{ 2 \atop s + l } \right]}
  -4\cdot \Delta _{k}\Gamma_{j+2}^{\left[ 1\atop{ 1 \atop s + l } \right]} - 8\cdot \Delta _{k}\Gamma_{j+2}^{\left[ 1\atop{ 2 \atop s + l -1} \right]}
  +32\cdot \Delta _{k}\Gamma_{j+1}^{\left[ 1\atop{ 1 \atop s + l -1} \right]}\\
  & = [9\cdot2^{2k+2j-1} +105\cdot2^{k+3j}] -4\cdot[9\cdot2^{2k+2j-2} +63\cdot2^{k+3j-2}]\\
   &  -8\cdot[9\cdot2^{2k+ 2j-3} + 105\cdot2^{k+3j-3}] +32\cdot[9\cdot2^{2k+2j-4} +63\cdot2^{k+3j -5}] = 9\cdot2^{2k+2j-1}\\
  & \text{It follows that}\\
 &  \sum_{i = 0}^{s - 2}4^{2i+1}\cdot \Delta _{k}\omega_{2(s-i-1)+j+1} (s-i-1,1, l+ i) =  \sum_{i = 0}^{s - 2}2^{4i+2}\cdot( 9\cdot2^{2k +2(s-i) +2j-5}) \\
 & = 3\cdot2^{2k +4s +2j-5} -3\cdot2^{2k +2s +2j-3}.
\end{align*}

The case $j = l$, \\
\begin{align*}
 &  \Delta _{k}\omega_{2(s-i-1)+l+1} (s-i-1,1, l+ i) = 9\cdot2^{2k +2(s-i)+2l-5}\quad  \text{if $ 1\leq i\leq s -2 $},\\
  &  \Delta _{k}\omega_{2(s-1)+l+1} (s-1,1, l) = -15\cdot2^{2k +2s+2l-5}\quad  \text{if $ i = 0 $},\\
  & \Delta _{k}\omega_{l+3} (1,1, l + s-2) =  \Delta _{k}\Gamma_{l+3}^{\left[ 1\atop{ 2 \atop s + l } \right]}
  -4\cdot \Delta _{k}\Gamma_{l+2}^{\left[ 1\atop{ 1 \atop s + l } \right]} - 8\cdot \Delta _{k}\Gamma_{l+2}^{\left[ 1\atop{ 2 \atop s + l -1} \right]}
  +32\cdot \Delta _{k}\Gamma_{l+1}^{\left[ 1\atop{ 1 \atop s + l -1} \right]}\\
   & = [9\cdot2^{2k+2l-1} +105\cdot2^{k+3l}] -4\cdot[9\cdot2^{2k+2l-2} +63\cdot2^{k+3l-2}] \\
    &  -8\cdot[9\cdot2^{2k+ 2l-3} + 105\cdot2^{k+3l-3}] +32\cdot[9\cdot2^{2k+2l-4} +63\cdot2^{k+3l -5}] = 9\cdot2^{2k+2l-1}\\
   & \text{It follows that}\\
 &  \sum_{i = 0}^{s - 2}4^{2i+1}\cdot \Delta _{k}\omega_{2(s-i-1)+l+1} (s-i-1,1, l+ i) =  \sum_{i = 1}^{s - 2}2^{4i+2}\cdot( 9\cdot2^{2k +2(s-i) +2l-5}) +( -15\cdot2^{2k +2s+2l-5}) \\
 & = 3\cdot2^{2k +4s +2l-5} -27\cdot2^{2k +2s +2l-3}.
\end{align*}
\underline{\eqref{eq 11.13}} is obtained as indicated in Lemma \ref{lem 11.6}
  \end{proof}
 
    \subsection{\textbf{Computation of $ \Delta _{k}\Gamma_{2s+l+1+j}^{\left[s\atop{ s \atop s+l} \right]} -   8\cdot\Delta _{k}\Gamma_{2s+l+j}^{\left[s\atop{ s \atop s+l-1} \right]}  \quad for\; 0\leq j \leq s-1 $}}
\label{subsec 8} 
\begin{lem}
\label{lem 11.7}
We have for $ l\geq 3 $ the following reduction formula holding for $0\leq j\leq s -3.$\\
\begin{align}
    &   \Delta _{k}\Gamma_{2s+l+1+j}^{\left[s\atop{ s  \atop s +l } \right]}
    -8\cdot\Delta _{k}\Gamma_{2s+ l+j}^{\left[s\atop{ s  \atop s +l -1 } \right]} 
     =   2^{4j+8}\cdot \Delta _{k}\Gamma_{2(s-j-2)+l+j+1}^{\left[s -j-2\atop{ s - j-2 \atop s-j-2 +(j+2+l)} \right]}
   - 8\cdot 2^{4j+8}\cdot \Delta _{k}\Gamma_{2(s-j-2)+l+j}^{\left[s -j-2\atop{ s - j-2\atop s -j+2 + (j+1+l)} \right]} \label{eq 11.16}\\
   &   +\sum_{i = 0}^{j+1}2^{4i}\cdot \Delta _{k}\omega_{2(s-i)+ l+ j+1} (s-i,0, l+ i)  + 
   \sum_{i = 0}^{j+1}2^{4i+2}\cdot \Delta _{k}\omega_{2(s-i-1)+ l+j+2} (s-i-1,1, l+ i) \nonumber
\end{align}
 where 
 \begin{equation}
 \label{eq 11.17}
   \Delta _{k}\Gamma_{2(s-j-2)+l+j+1}^{\left[s -j-2\atop{ s - j-2 \atop s-j-2 +(j+2+l)} \right]}
   - 8\cdot  \Delta _{k}\Gamma_{2(s-j-2)+l+j}^{\left[s -j-2\atop{ s - j-2\atop s -j+2 + (j+1+l)} \right]} =
   -9\cdot2^{2k +2s+2l-6},  
 \end{equation}
  \begin{align}
& \sum_{i = 0}^{j+1}2^{4i}\cdot \Delta _{k}\omega_{2(s-i)+ l+ j+1} (s-i,0, l+ i) = 
  3\cdot2^{2k+2s+2l+4j -2} - 15\cdot2^{2k+2s+2l+4j} \label{eq 11.18},\\
&   \sum_{i = 0}^{j+1}2^{4i+2}\cdot \Delta _{k}\omega_{2(s-i-1)+ l+j+2} (s-i-1,1, l+ i)
 =3\cdot2^{2k+2s+2l+4j -1} - 15\cdot2^{2k+2s+2l+4j +1},\label{eq 11.19}
\end{align}
Combining \eqref{eq 11.16}, \eqref{eq 11.17}, \eqref{eq 11.18} and  \eqref{eq 11.19} we obtain \\
\begin{equation}
\label{eq 11.20}
 \Delta _{k}\Gamma_{2s+l+1+j}^{\left[s\atop{ s  \atop s +l } \right]}
    -8\cdot\Delta _{k}\Gamma_{2s+ l+j}^{\left[s\atop{ s  \atop s +l -1 } \right]} =   -315\cdot2^{2k+2s+2l+4j-2} \quad for \quad 0\leq j\leq s-3
\end{equation}
 \end{lem}
 
 \begin{proof}
\underline{proof of \eqref{eq 11.16}}
\begin{align*}
&
\end{align*}
 Let q be a rational integer between zero and s-2.
 From \eqref{eq 8.3} we obtain successively the following equations. \\
\small   
  \begin{align*}
 &  \Delta _{k}\Gamma_{2s+l+1+j}^{\left[s\atop{ s \atop s+l } \right]}
   - 4\cdot \Delta _{k}\Gamma_{2s+ l+j}^{\left[s -1\atop{ s \atop s+l } \right]}  = 
   8\cdot\left[  \Delta _{k}\Gamma_{2s+l+j}^{\left[s\atop{ s \atop s + l-1} \right]}
   - 4\cdot \Delta _{k}\Gamma_{2s+l+j-1}^{\left[s\atop{ s -1\atop s + l-1} \right]}\right] + \Delta _{k}\omega_{2s+l+1+j} (s,0,l) \\
     &  4\cdot\left[  \Delta _{k}\Gamma_{2s+l+j}^{\left[s -1\atop{ s \atop s+l } \right]}
   - 4\cdot \Delta _{k}\Gamma_{2s+l+j-1}^{\left[s -1\atop{ s -1 \atop s +l } \right]}\right]  = 
   4\cdot\left[   8\cdot\Big[  \Delta _{k}\Gamma_{2s+l+j-1}^{\left[s-1\atop{ s  \atop s +l-1 } \right]}
   - 4\cdot \Delta _{k}\Gamma_{2s+l+j-2}^{\left[s -1\atop{ s -1\atop s +l-1 } \right]}\Big] + \Delta _{k}\omega_{2s+l+j} (s-1,1, l) \right]\\
     &  4^{2}\cdot\left[   \Delta _{k}\Gamma_{2s+l+j-1}^{\left[s -1\atop{ s -1 \atop s+l } \right]}
   - 4\cdot \Delta _{k}\Gamma_{2s+l+j-2}^{\left[s\atop{ s -2 \atop s +l} \right]}\right]  =
    4^{2}\cdot\left[  8\cdot\Big[  \Delta _{k}\Gamma_{2s+l+j-2}^{\left[s -1\atop{ s -1 \atop s +l-1} \right]}
   - 4\cdot \Delta _{k}\Gamma_{2s+ l+j-3}^{\left[s -2\atop{ s -1 \atop s +l-1 } \right]}\Big] + \Delta _{k}\omega_{2s+ l+j-1} (s-1,0,l+1)\right] \\
     &  \vdots       \\
      &  4^{2q}\cdot\left[   \Delta _{k}\Gamma_{2s+ l +1+j-2q}^{\left[s -q\atop{ s -q \atop s+l } \right]}
       - 4\cdot \Delta _{k}\Gamma_{2s+ l+j -2q}^{\left[s -q-1\atop{ s -q \atop s +l} \right]}\right]\\
         & =  4^{2q}\cdot\left[  8\cdot\Big[  \Delta _{k}\Gamma_{2s+ l+j -2q}^{\left[s -q\atop{ s -q \atop s +l-1} \right]}
          - 4\cdot \Delta _{k}\Gamma_{2s+l+j-1-2q}^{\left[s -q-1\atop{ s -q \atop s +l -1 } \right]}\Big] + \Delta _{k}\omega_{2s+ l+1+j-2q} (s -q,0,l+ q)\right] \\
    &  4^{2q +1}\cdot\left[   \Delta _{k}\Gamma_{2s+ l+j -2q}^{\left[s -q-1\atop{ s -q \atop s+l } \right]}
   - 4\cdot \Delta _{k}\Gamma_{2s+l +j-1-2q}^{\left[s -q-1\atop{ s -q-1 \atop s +l} \right]}\right]\\
    & =  4^{2q+1}\cdot\left[  8\cdot\Big[  \Delta _{k}\Gamma_{2s+ l +j-1-2q}^{\left[s -q-1\atop{ s -q \atop s +l-1} \right]}
   - 4\cdot \Delta _{k}\Gamma_{2s+ l+j-2-2q}^{\left[s -q-1\atop{ s -q-1 \atop s +l -1 } \right]}\Big] + \Delta _{k}\omega_{2s+ l+j-2q} (s -q-1,1,l+ q)\right] 
   \end{align*}
 By summing  the above equations we obtain the following formula holding for $0\leq q\leq s-2$ \\ 
  \begin{align}
   &   \Delta _{k}\Gamma_{2s+l+1+j}^{\left[s\atop{ s  \atop s +l } \right]}
    -8\cdot\Delta _{k}\Gamma_{2s+ l+j}^{\left[s\atop{ s  \atop s +l -1 } \right]} 
     =   4^{2q+2}\cdot \Delta _{k}\Gamma_{2s+l+j-1-2q}^{\left[s -q-1\atop{ s -q-1\atop s + l} \right]}
   - 8\cdot 4^{2q+2}\cdot \Delta _{k}\Gamma_{2s +l+j-2-2q}^{\left[s -q-1\atop{ s -q-1\atop s + l -1} \right]}\label{eq 11.21} \\
   &   +\sum_{i = 0}^{q}4^{2i}\cdot \Delta _{k}\omega_{2(s-i)+ l+ j+1} (s-i,0, l+ i)  + 
   \sum_{i = 0}^{q}4^{2i+1}\cdot \Delta _{k}\omega_{2(s-i-1)+ l+j+2} (s-i-1,1, l+ i)\nonumber \\
   & \text{Setting q = j+1 for $ 0\leq j\leq s -3 $ in \eqref{eq 11.21} we get the reduction formula \eqref{eq 11.16}} \nonumber 
    \end{align}
   \underline{Proof of \eqref{eq 11.17}} 
 \begin{align*}
  &
\end{align*}
 We get from  \eqref{eq 11.13} with $s\rightarrow s -j-2,\;l\rightarrow j+2+l,\;j\rightarrow l+j+1:$\\
  \begin{align*}
 &  \Delta _{k}\Gamma_{2(s-j-2)+l+j+1}^{\left[s -j-2\atop{ s - j-2 \atop s-j-2 +(j+2+l)} \right]}
   - 8\cdot  \Delta _{k}\Gamma_{2(s-j-2)+l+j}^{\left[s -j-2\atop{ s - j-2\atop s -j-2 + (j+1+l)} \right]} =
   -9\cdot2^{2k +2(s-j-2)+2(j+2+l)-4} =  -9\cdot2^{2k +2s+2l-6}, \\
   & \\
     & \text{observing that $2\leq l+j+1\leq j+2+l-1$ and $s-j-2\geq s-(s-3)-2=1 $}  \\
     & \text{Remark that in the case $ j = s-3  $ we have from Lemma \ref{lem 1.23}}\\
      &  \Delta _{k}\Gamma_{1+(l+s-1)}^{\left[1\atop{ 1 \atop 1+(s+l-1) } \right]}
   - 8\cdot  \Delta _{k}\Gamma_{1+(l+s-2)}^{\left[1\atop{ 1\atop 1+(s+l-2)} \right]} =  -9\cdot2^{2k +2s+2l-6}\\
   & 
   \end{align*}
   \underline{Proof of \eqref{eq 11.18}} 
 \begin{align*}
  &
\end{align*}
We obtain from Lemma \ref{lem 7.4} (see \eqref{eq 7.51})
\begin{align*}
 &   \Delta _{k}\omega_{2(s-i) +l+j+1} (s-i,0, l+ i) = 0 \\
 &\Updownarrow \\
 &  2(s-i) +l+i+2\leq 2(s-i) +l+j+1\leq 3\cdot(s-i) +l+i-1\\
   & \Updownarrow \\
  &   i\leq j-1\; and\;j\leq s-2 
  \end{align*}
We then deduce :\\
\begin{align*}
 & \sum_{i = 0}^{j+1}2^{4i}\cdot \Delta _{k}\omega_{2(s-i)+ l+ j+1} (s-i,0, l+ i) = 
 \sum_{i = j}^{j+1}2^{4i}\cdot \Delta _{k}\omega_{2(s-i) +(l+j+1)} (s-i,0, l+ i)  = \\
 & =  2^{4j}\cdot \Delta _{k}\omega_{2(s-j)+l+j+1} (s-j,0, l+ j) 
    +  2^{4(j+1)}\cdot \Delta _{k}\omega_{2(s-j-1)+ l+j+1 } (s-j-1,0, l+ j+1) \\
&  =   2^{4j}\cdot3\cdot2^{2k+2(s-j)+2(l+j) -2} + 2^{4(j+1)}\cdot(-15)\cdot2^{2k+2(s-j-1)+2(l+j+1)-4} \\
& =  3\cdot2^{2k+2s+2l+4j -2} - 15\cdot2^{2k+2s+2l+4j} \\
&
\end{align*}
 \underline{Proof of \eqref{eq 11.19}} 
 \begin{align*}
  &
\end{align*}
Similar to the proof of \eqref{eq 11.18}.\\
\begin{align*}
 &   \Delta _{k}\omega_{2(s-i-1) +l+j+2} (s-i-1,1, l+ i) = 0 \\
 &\Updownarrow \\
 &1+  2(s-i-1) +l+i+2\leq 2(s-i-1) +l+j+2\leq 2+3\cdot(s-i-1) +l+i-1\\
   & \Updownarrow \\
  &   i\leq j-1\; and\;j\leq s-2 
  \end{align*}
We then deduce :\\
\begin{align*}
 & \sum_{i = 0}^{j+1}2^{4i+2}\cdot \Delta _{k}\omega_{2(s-i-1)+ l+ j+2} (s-i-1,1, l+ i) = 
 \sum_{i = j}^{j+1}2^{4i+2}\cdot \Delta _{k}\omega_{2(s-i-1) +(l+j+2)} (s-i-1,1, l+ i)  = \\
 & =  2^{4j+2}\cdot \Delta _{k}\omega_{2(s-j-1)+l+j+2} (s-j-1,1, l+ j) 
    +  2^{4(j+1)+2}\cdot \Delta _{k}\omega_{2(s-j-1)+ l+j+2 } (s-j-2,1, l+ j+1) \\
&  =   2^{4j+2}\cdot3\cdot2^{2k+2(s-j-1)+1+2(l+j) -2} + 2^{4j+6}\cdot(-15)\cdot2^{2k+2(s-j-2)+1+2(l+j+1)-4} \\
& =  3\cdot2^{2k+2s+2l+4j -1} - 15\cdot2^{2k+2s+2l+4j+1} \\
&
\end{align*}
\underline{\eqref{eq 11.19}} is obtained as indicated in Lemma \ref{lem 11.7}
    \end{proof}
 
 \begin{lem}
\label{lem 11.8}
We have for $ l\geq 3 $ the following reduction formula \\
\begin{align}
    &   \Delta _{k}\Gamma_{3s+l-1}^{\left[s\atop{ s  \atop s +l } \right]}
    -8\cdot\Delta _{k}\Gamma_{3s+l-2}^{\left[s\atop{ s  \atop s +l -1 } \right]} 
     =   2^{4s-4}\cdot \Delta _{k}\Gamma_{s+l+1}^{\left[1\atop{ 1 \atop s+l} \right]}
   - 8\cdot 2^{4s-4}\cdot \Delta _{k}\Gamma_{s+l}^{\left[1\atop{ 1\atop s+l-1} \right]} \label{eq 11.22}\\
   &   +\sum_{i = 0}^{s-2}2^{4i}\cdot \Delta _{k}\omega_{2(s-i)+ l+s-1} (s-i,0, l+ i)  + 
   \sum_{i = 0}^{s-2}2^{4i+2}\cdot \Delta _{k}\omega_{2(s-i-1)+ l+s} (s-i-1,1, l+ i) \nonumber
\end{align}
 where 
 \begin{equation}
 \label{eq 11.23}
   \Delta _{k}\Gamma_{s+l+1}^{\left[1\atop{ 1 \atop s+l} \right]}
   - 8\cdot  \Delta _{k}\Gamma_{s+l}^{\left[1\atop{ 1\atop s+l-1} \right]} =
   -81\cdot2^{2k +2s+2l-4},  
 \end{equation}
  \begin{align}
& \sum_{i = 0}^{s-2}2^{4i}\cdot \Delta _{k}\omega_{2(s-i)+ l+ s-1} (s-i,0, l+ i) = 
  3\cdot2^{2k+2l+6s-10}  \label{eq 11.24},\\
&   \sum_{i = 0}^{s-2}2^{4i+2}\cdot \Delta _{k}\omega_{2(s-i-1)+ l+s} (s-i-1,1, l+ i)
 =3\cdot2^{2k+2s+2l+6s-9} ,\label{eq 11.25}
\end{align}
Combining \eqref{eq 11.22}, \eqref{eq 11.23}, \eqref{eq 11.24} and  \eqref{eq 11.25} we obtain \\
\begin{equation}
\label{eq 11.26}
 \Delta _{k}\Gamma_{3s+l-1}^{\left[s\atop{ s  \atop s +l } \right]}
    -8\cdot\Delta _{k}\Gamma_{3s+l-2}^{\left[s\atop{ s  \atop s +l -1 } \right]} =   -315\cdot2^{2k+2l+6s-10} 
\end{equation}
\end{lem}
\begin{proof}
\begin{align*}
 &
\end{align*}
\underline{Proof of \eqref{eq 11.22}}

Setting $q=j=s-2$ in \eqref{eq 11.21} we get the reduction formula  \eqref{eq 11.22} \\
\begin{align*}
 &
\end{align*}
\underline{Proof of \eqref{eq 11.23}}

 From  Lemma \ref{lem 1.23} we obtain :\\
 
  \begin{align*}
   &  \Delta _{k}\Gamma_{s+l+1}^{\left[1\atop{ 1 \atop 1+s+l-1} \right]}
   - 8\cdot  \Delta _{k}\Gamma_{s+l}^{\left[1\atop{ 1\atop 1+s+l-2} \right]} \\
    &  = 81\cdot(2^{2k+2s+2l-4} -2^{k+3s+3l- 5}) - 8\cdot(81\cdot(2^{2k+2s+2l-6} -2^{k+3s+3l- 8})) 
     =  -81\cdot2^{2k +2s+2l-4}.\\
     &
  \end{align*}
\underline{Proof of \eqref{eq 11.24}}
 We obtain from Lemma \ref{lem 7.4} (see \eqref{eq 7.37})
\begin{align*}
 &   \Delta _{k}\omega_{2(s-i) +l+s-1} (s-i,0, l+ i) = 0 \\
 &\Updownarrow \\
 &  2(s-i) +l+i+2\leq 2(s-i) +l+s-1\leq 3\cdot(s-i) +l+i-1\\
   & \Updownarrow \\
  &   i\leq s-3
  \end{align*}
We then deduce :\\
\begin{align*}
 & \sum_{i = 0}^{s-2}2^{4i}\cdot \Delta _{k}\omega_{2(s-i)+ l+s-1} (s-i,0, l+ i) = 
 \sum_{i = s-2}^{s-2}2^{4i}\cdot \Delta _{k}\omega_{2(s-i) +(l+s-1)} (s-i,0, l+ i)  = \\
 & =  2^{4(s-2)}\cdot \Delta _{k}\omega_{l+s+3} (2,0, l+s-2) = 2^{4(s-2)}\cdot3\cdot2^{2k+2(l+s-2)+2}
  =  3\cdot2^{2k+2l+6s-10}  \\
&
\end{align*}
\underline{Proof of \eqref{eq 11.25}}
 We obtain from Lemmas \ref{lem 1.22}, \ref{lem 1.23} 
\begin{align*}
 &   \Delta _{k}\omega_{2(s-i-1) +l+s} (s-i-1,1, l+ i) = 0 \\
 &\Updownarrow \\
 & 1+ 2(s-i-1) +l+i+2\leq 2(s-i-1) +l+s \leq  2+3\cdot(s-i-1) +l+i-1\\
   & \Updownarrow \\
  &   i\leq s-3
  \end{align*}
We then deduce :\\
\begin{align*}
 & \sum_{i = 0}^{s-2}2^{4i+2}\cdot \Delta _{k}\omega_{2(s-i-1)+ l+s} (s-i-1,1, l+ i) = 
 \sum_{i = s-2}^{s-2}2^{4i+2}\cdot \Delta _{k}\omega_{2(s-i-1) +l+s} (s-i-1,1, l+ i)  = \\
 & =  2^{4s-6}\cdot \Delta _{k}\omega_{l+s+2} (1,1, l+s-2) \\
 & = 2^{4s-6}\cdot\big(  \Delta _{k}\Gamma_{l+s+2}^{\left[ 1\atop{ 2 \atop s + l } \right]}
  -4\cdot \Delta _{k}\Gamma_{l+s+1}^{\left[ 1\atop{ 1 \atop s + l } \right]} - 8\cdot \Delta _{k}\Gamma_{l+s+1}^{\left[ 1\atop{ 2 \atop s + l -1} \right]}
  +32\cdot \Delta _{k}\Gamma_{l+s}^{\left[ 1\atop{ 1 \atop s + l -1} \right]}\big)\\
  & =  2^{4s-6}\cdot\Big(\big ( \Delta _{k}\Gamma_{l+s+2}^{\left[ 1\atop{ 1+1 \atop 1+1+ s + l-2 } \right]} - 
   8\cdot \Delta _{k}\Gamma_{l+s+1}^{\left[ 1\atop{ 1+1 \atop 1+1+ s + l -3} \right]} \big)
    -4 \cdot\big( \Delta _{k}\Gamma_{l+s+1}^{\left[ 1\atop{ 1 \atop 1+ s + l -1} \right]}  - 
     8 \cdot \Delta _{k}\Gamma_{l+s}^{\left[ 1\atop{ 1 \atop 1+ s + l -2} \right]} \big)  \Big) \\
     & =  2^{4s-6}\cdot\Big( -159\cdot2^{2k+2s+2l-3} -4\cdot(-81)\cdot2^{2k+2l+2s-4}  \Big) =  2^{4s-6}\cdot3\cdot2^{2k+2l+2s-3}\\
     & = 3\cdot2^{2k+2l+6s-9},\\
     & \text{or otherwise using \eqref{eq 5.2} and \eqref{eq 7.37}}\\
     &  \Delta _{k}\omega_{l+s+2} (1,1, l+s-2) = {1\over 2}\cdot \Delta _{k}\omega_{l+s+3} (2,0, l+s-2) = 3\cdot2^{2k+2l+2s-3}
\end{align*}
\end{proof}
\begin{lem}
\label{lem 11.9}
We have for $ l\geq 3 $ the following reduction formula 
\begin{align}
    &   \Delta _{k}\Gamma_{3s+l}^{\left[s\atop{ s  \atop s +l } \right]}
    -8\cdot\Delta _{k}\Gamma_{3s+l-1}^{\left[s\atop{ s  \atop s +l -1 } \right]} 
     =   2^{4s-4}\cdot \Delta _{k}\Gamma_{s+l+2}^{\left[1\atop{ 1 \atop s+l} \right]}
   - 8\cdot 2^{4s-4}\cdot \Delta _{k}\Gamma_{s+l+1}^{\left[1\atop{ 1\atop s+l-1} \right]} \label{eq 11.27}\\
   &   +\sum_{i = 0}^{s-2}2^{4i}\cdot \Delta _{k}\omega_{2(s-i)+ l+s} (s-i,0, l+ i)  + 
   \sum_{i = 0}^{s-2}2^{4i+2}\cdot \Delta _{k}\omega_{2(s-i-1)+ l+s+1} (s-i-1,1, l+ i) \nonumber
\end{align}
 where 
 \begin{equation}
 \label{eq 11.28}
   \Delta _{k}\Gamma_{s+l+2}^{\left[1\atop{ 1 \atop 1+s+l-1} \right]}
   - 8\cdot  \Delta _{k}\Gamma_{s+l+1}^{\left[1\atop{ 1\atop 1+s+l-2} \right]} =
   -21\cdot2^{3k+s+l-1} + 21\cdot2^{2k+2s+2l-2} ,  
 \end{equation}
  \begin{align}
& \sum_{i = 0}^{s-2}2^{4i}\cdot \Delta _{k}\omega_{2(s-i)+ l+ s} (s-i,0, l+ i) = 
 7\cdot2^{3k+5s+l-5} - 7\cdot2^{3k+3s+l-3}  \label{eq 11.29},\\
&   \sum_{i = 0}^{s-2}2^{4i+2}\cdot \Delta _{k}\omega_{2(s-i-1)+ l+s+1} (s-i-1,1, l+ i)
 = 7\cdot2^{3k+5s+l-6} - 7\cdot2^{3k+3s+l-4} ,\label{eq 11.30}
\end{align}
Combining \eqref{eq 11.27}, \eqref{eq 11.28}, \eqref{eq 11.29} and  \eqref{eq 11.30} we obtain 
\begin{equation}
\label{eq 11.31}
 \Delta _{k}\Gamma_{3s+l}^{\left[s\atop{ s  \atop s +l } \right]}
    -8\cdot\Delta _{k}\Gamma_{3s+l-1}^{\left[s\atop{ s  \atop s +l -1 } \right]} =  -21\cdot2^{3k+3s+l-3} + 21\cdot2^{2k+6s+2l-6}
\end{equation}
\end{lem}
\begin{proof}

\underline{Proof of \eqref{eq 11.27}}

Setting $q = s-2,\;j=s-1$ in \eqref{eq 11.21} we get the reduction formula  \eqref{eq 11.27} \\

\medskip
\underline{Proof of \eqref{eq 11.28}}
 From  Lemma \ref{lem 1.23} we obtain :\\
   \begin{align*}
   &  \Delta _{k}\Gamma_{s+l+2}^{\left[1\atop{ 1 \atop 1+s+l-1} \right]}
   - 8\cdot  \Delta _{k}\Gamma_{s+l+1}^{\left[1\atop{ 1\atop 1+s+l-2} \right]} \\
  &  = (7\cdot2^{3k+ s+l-1} - 21\cdot2^{2k+ 2s+2l -2} +7\cdot2^{k+3s+3l-2})\\
   & - 8\cdot(7\cdot2^{3k+ s+l-2} - 21\cdot2^{2k+ 2s+2l -4} +7\cdot2^{k+3s+3l-5})\\
& = - 21\cdot2^{3k+ s+ l-1} + 21\cdot2^{2k+ 2s+2l -2}\\
     &
  \end{align*}
\underline{Proof of \eqref{eq 11.29}}
 We obtain from Lemma \ref{lem 7.4} (see \eqref{eq 7.37}, \eqref{eq 7.46} and  \eqref{eq 7.51})
\begin{align*}
    \Delta _{k}\omega_{2(s-i) +l+s} (s-i,0, l+ i) =  21\cdot2^{3k+ 3s-2i+ l-3}\quad for\quad 0\leq i\leq s-2
 \end{align*}
We then deduce :\\
\begin{align*}
 & \sum_{i = 0}^{s-2}2^{4i}\cdot \Delta _{k}\omega_{2(s-i)+ l+s} (s-i,0, l+ i) = 
 \sum_{i = 0}^{s-2}2^{4i}\cdot21\cdot2^{3k+ 3s-2i+ l-3}
  = 7\cdot2^{3k+ 5s+l-5} - 7\cdot2^{3k+ 3s+l-3}
 & \\
 \end{align*}
 
\underline{Proof of \eqref{eq 11.30}}
 We obtain from Lemma \ref{lem 7.4} (see \eqref{eq 7.38}, \eqref{eq 7.47} and  \eqref{eq 7.51})
\begin{align*}
  &  \Delta _{k}\omega_{2(s-i-1) +l+s+1} (s-i-1,1, l+ i) =  21\cdot2^{3k+ 3s-2i+ l-4}\quad for\quad 0\leq i\leq s-3 \\
  & \Delta _{k}\omega_{l+s+3} (1,1, l+ s-2) = {1\over 2}\cdot\Delta _{k}\omega_{l+s+4} (2,0, l+ s-2)\\
   &   = {1\over 2}\cdot21\cdot2^{3k+ l+s+1} = 21\cdot2^{3k+ l+s}  \quad for \quad i = s-2
 \end{align*}
We then deduce :\\
\begin{align*}
 & \sum_{i = 0}^{s-2}2^{4i+2}\cdot \Delta _{k}\omega_{2(s-i-1)+ l+s+1} (s-i-1,1, l+ i) = 
 \sum_{i = 0}^{s-2}2^{4i+2}\cdot21\cdot2^{3k+ 3s-2i+ l-4}\\
 &  = 7\cdot2^{3k+ 5s+l-6} - 7\cdot2^{3k+ 3s+l-4}
 & \\
 \end{align*}
\end{proof}
  \subsection{\textbf{Computation of $ \Delta _{k}\Gamma_{2s +m +j}^{\left[s\atop{ s +m\atop s+m} \right]} -   8\cdot\Delta _{k}\Gamma_{2s+ m+ j -1}^{\left[s\atop{ s +m-1\atop s +m } \right]}  \quad for\; 0\leq j \leq m $}}
\label{subsec 9} 
  \begin{lem}
\label{lem 11.10}
 We have for $ m \geq 3 $ the following reduction formula holding for $0\leq j\leq s+m.$\\
 \begin{align}
& \Delta _{k}\Gamma_{2s+m+j}^{\left[s\atop{ s+m  \atop s+m  } \right]}
    -8\cdot\Delta _{k}\Gamma_{2s+m+j-1}^{\left[s\atop{ s +m -1\atop s +m  } \right]} \label{eq 11.32}\\ 
    & = 2^{2m}\cdot\left[\Delta _{k}\Gamma_{2s + j}^{\left[s \atop{ s \atop s +m} \right]}
  -8\cdot \Delta _{k}\Gamma_{2s +j-1}^{\left[s \atop{ s  \atop s +m -1} \right]}\right]
   + \sum_{i = 0}^{m-1}2^{2i}\cdot \Delta _{k}\omega_{2s +m+j-i} (s,m-i,  i), \nonumber 
 \end{align}
 where 
\begin{equation}
\label{eq 11.33}
 \Delta _{k}\Gamma_{2s+j}^{\left[s\atop{ s \atop s+m} \right]} -   8\cdot\Delta _{k}\Gamma_{2s+j -1}^{\left[s\atop{ s \atop s+m-1} \right]}    =
 \begin{cases}
 3\cdot2^{2k+2s-2} - 33\cdot2^{k+4s-6} & \text{if  } j = 0, \\
  -15\cdot2^{2k+2s-2}  +3\cdot2^{k+4s-3}   & \text{if  } j = 1, \\
   -9\cdot2^{2k+2s +2j-4}      & \text{if  } 2\leq j\leq m-1,  \\
     -81\cdot2^{2k+2s +2m-4}      & \text{if  } j = m,  
    \end{cases}
\end{equation}
\begin{align}
\label{eq 11.34}
\sum_{i = 0}^{m-1}2^{2i}\cdot \Delta _{k}\omega_{2s+m+j-i} (s,m-i, i) =
 \begin{cases}
    -3\cdot 2^{2k +2s+m-2} -3 \cdot 2^{2k +2s +m-2}   & \text{if  } j = 0, \\
     15\cdot 2^{2k +2s+2m-2} - 69\cdot 2^{2k +2s +m-2} & \text{if  } j = 1, \\
    - 63\cdot 2^{2k +2s+m+3j-5}+9\cdot 2^{2k +2s +2m+2j-4}   & \text{if  } 2\leq j\leq m-1,  \\
   3\cdot 2^{2k +2s+4m-5}   & \text{if  } j = m.  
    \end{cases}
\end{align}
Combining \eqref{eq 11.32}, \eqref{eq 11.33} and \eqref{eq 11.34}  we obtain \\

\begin{equation}
\label{eq 11.35}
 \Delta _{k}\Gamma_{2s+m+j}^{\left[s\atop{ s +m\atop s+m} \right]}
  -   8\cdot\Delta _{k}\Gamma_{2s+m+j -1}^{\left[s\atop{ s +m-1\atop s+m} \right]}    =
 \begin{cases}
 -3\cdot2^{2k+2s+m-2} - 33\cdot2^{k+4s+2m-6} & \text{if  } j = 0 \\
  -69\cdot2^{2k+2s+m-2}  +3\cdot2^{k+4s+2m-3}   & \text{if  } j = 1 \\
   -63\cdot2^{2k+2s +m+3j-5}      & \text{if  } 2\leq j\leq m-1  \\
     -159\cdot2^{2k+2s +4m-5}      & \text{if  } j = m   
    \end{cases}
\end{equation}
 \end{lem}
 \begin{proof}
  \underline{Proof of \eqref{eq 11.32}} \\
 Let q be a rational integer between zero and  m-1.
 From \eqref{eq 8.3} we obtain successively the following equations. \\
    \begin{align*}
 &  \Delta _{k}\Gamma_{2s+m+j}^{\left[s\atop{ s +m\atop s+m } \right]}
   - 4\cdot \Delta _{k}\Gamma_{2s+m+j -1}^{\left[s\atop{ s +m-1\atop s+m } \right]}  = 
   8\cdot\left[  \Delta _{k}\Gamma_{2s+m+j -1}^{\left[s\atop{ s+m -1 \atop s +m } \right]}
   - 4\cdot \Delta _{k}\Gamma_{2s+m+j -2}^{\left[s\atop{ s +m-1\atop s +m-1} \right]}\right] + \Delta _{k}\omega_{2s+m+j} (s,m,0) \\
  &  4\cdot\left[  \Delta _{k}\Gamma_{2s+m+j-1}^{\left[s \atop{ s +m-1\atop s+m } \right]}
   - 4\cdot \Delta _{k}\Gamma_{2s+m+j -2}^{\left[s\atop{ s+m-2 \atop s +m } \right]}\right]  = 
4\cdot\left[   8\cdot\Big[  \Delta _{k}\Gamma_{2s+m+j -2}^{\left[s\atop{ s+m-1 \atop s +m-1 } \right]}
   - 4\cdot \Delta _{k}\Gamma_{2s+m+j -3}^{\left[s\atop{ s +m-2\atop s +m-1 } \right]}\Big] + \Delta _{k}\omega_{2s+m+j-1} (s,m-1, 1) \right]\\ 
     &  4^{2}\cdot\left[   \Delta _{k}\Gamma_{2s+m+j-2}^{\left[s \atop{ s +m-2 \atop s+m } \right]}
   - 4\cdot \Delta _{k}\Gamma_{2s+m+j -3}^{\left[s\atop{ s +m-3 \atop s +m} \right]}\right]  = 
 4^{2}\cdot\left[  8\cdot\Big[  \Delta _{k}\Gamma_{2s+m+j -3}^{\left[s \atop{ s +m-2 \atop s +m-1} \right]}
   - 4\cdot \Delta _{k}\Gamma_{2s+m+j - 4}^{\left[s \atop{ s +m-3 \atop s +m-1 } \right]}\Big] + \Delta _{k}\omega_{2s+m+j-2} (s,m-2,2)\right] \\
     &  \vdots       \\
   &  4^{q}\cdot\left[   \Delta _{k}\Gamma_{2s+m+j - q}^{\left[s \atop{ s +m - q) \atop s+m } \right]}
   - 4\cdot \Delta _{k}\Gamma_{2s+m+j - (q+1)}^{\left[s\atop{ s +m - (q+1) \atop s +m} \right]}\right]\\
    & =  4^{q}\cdot\left[  8\cdot\Big[  \Delta _{k}\Gamma_{2s+m+j -(q+1)}^{\left[s \atop{ s +m- q \atop s +m-1} \right]}
   - 4\cdot \Delta _{k}\Gamma_{2s+m+j - (q+2)}^{\left[s \atop{ s +m- (q+1) \atop s +m -1 } \right]}\Big] + \Delta _{k}\omega_{2s+m +j-q} (s,m- q, q)\right]  
 \end{align*}
    By summing  the above equations we obtain after reduction the following formula  \\ 
     \begin{align*}
    &  \Delta _{k}\Gamma_{2s+m+j}^{\left[s\atop{ s+m  \atop s+m  } \right]}
    -8\cdot\Delta _{k}\Gamma_{2s+m+j-1}^{\left[s\atop{ s +m -1\atop s +m } \right]}\\
     & =   4^{q+1}\cdot \Delta _{k}\Gamma_{2s +m+j-(q+1)}^{\left[s \atop{ s +m-(q+1)\atop s +m} \right]}
   - 8\cdot 4^{q+1}\cdot \Delta _{k}\Gamma_{2s + m+j -(q+2)}^{\left[s \atop{ s  +m -(q+1)\atop s + m -1} \right]}
      +\sum_{i = 0}^{q}4^{i}\cdot \Delta _{k}\omega_{2s +m+j-i} (s,m-i,  i) \\
      & \text{Setting q = m-1 in the above formula we get \eqref{eq 11.32}} 
      \end{align*}
     \underline{Proof of \eqref{eq 11.33}} \\[0.5 cm]
     \eqref{eq 11.33} is obtained by \eqref{eq 11.13} with $l\longmapsto m$.\\[0.5 cm]
     \underline{Proof of \eqref{eq 11.34}} \\[0.5 cm]
  The proof is based on results in  Lemma \ref{lem 7.4} \\
\textbf{The case $ j = 0 $}, we have \\
\begin{align*}
 &  \Delta _{k}\omega_{2s+(m-i)} (s,m-i,i) = -3\cdot2^{2k+2s+(m-i)-2}\quad  \text{if $ 1\leq i\leq s -3 $},\\
 &  \Delta _{k}\omega_{2s+m} (s,m,0) = -9\cdot2^{2k+2s+m-2}  \quad  \text{if  i = 0 }.\\
 & \text{It follows that}\\
 &  \sum_{i = 0}^{m-1}2^{2i}\cdot \Delta _{k}\omega_{2s+(m-i)} (s,m-i,i) = (-9\cdot2^{2k+2s+m-2}) + \sum_{i = 1}^{m-1}2^{2i}\cdot( -3\cdot2^{2k+2s+(m-i)-2}) \\
 & =   -3\cdot 2^{2k +2s+m-2} -3 \cdot 2^{2k +2s +m-2}
\end{align*}

\textbf{The case $ j = 1 $}, we obtain similarly \\
\begin{align*}
 &  \Delta _{k}\omega_{2s+(m-i)+1} (s,m-i,i) = 15\cdot2^{2k+2s+(m-i)-2}\quad  \text{if $ 2\leq i\leq m-1 $},\\
  &  \Delta _{k}\omega_{2s+(m-1)+1} (s,m-1,1) = -9\cdot2^{2k+2s+(m-1)-2}\quad  \text{if $ i=1 $},\\
   &  \Delta _{k}\omega_{2s+m+1} (s,m,0) = 9\cdot2^{2k+2s+m-2}\quad  \text{if $ i=0 $},\\
 & \text{It follows that}\\
 &  \sum_{i = 0}^{m-1}2^{2i}\cdot \Delta _{k}\omega_{2s+(m-i)+1} (s,m-i, i) \\
 & =   9\cdot2^{2k+2s+m-2} -9\cdot2^{2k+2s+(m-1)} 
 + \sum_{i = 2}^{m-1}2^{2i}\cdot( 15\cdot2^{2k+2s+(m-i)-2}) \\
 &  =   15\cdot 2^{2k +2s+2m-2} - 69\cdot 2^{2k +2s +m-2}
\end{align*}

\textbf{The case $2\leq j\leq m-1$}, similarly \\
  \begin{align*}
&  \Delta _{k}\omega_{2s + (m-i)+j} (s,m-i,  i) =  0 \quad \text{if \; $ 0\leq i\leq j -2 $   }\\
&  \Delta _{k}\omega_{2s + (m- i) +j} (s,m -i,  i) = 3\cdot2^{2k +2s + m+j-3}\quad \text{if \; i = j -1   }\\
&  \Delta _{k}\omega_{2s + (m- i) +j} (s,m -i,  i) = -15\cdot2^{2k +2s + m+j-4}\quad \text{if \; i = j }\\
&  \Delta _{k}\omega_{2s + (m- i) +j} (s,m -i,  i) = 9\cdot2^{2k +2s + m + 2j-4 -i}\quad \text{if \;$ i\geq j+1 $ }\\
& \text{We then deduce}\\
&  \sum_{i = 0}^{m-1}2^{2i}\cdot \Delta _{k}\omega_{2s + (m-i)+ j} (s,m-i,  i) = 
=  \sum_{i = j-1}^{m-1}2^{2i}\cdot \Delta _{k}\omega_{2s + (m-i)+j} (s,m-i,  i) \\
    & = ( 3\cdot2^{2k +2s + m +j - 3})\cdot2^{2(j-1)} + ( -15\cdot2^{2k +2s + m +j-4})\cdot2^{2j} 
    + \sum_{i = j+1}^{m-1}2^{2i}\cdot 9\cdot2^{2k +2s + m + 2j-4 -i}\\
    & =   - 63\cdot 2^{2k +2s+m+3j-5}+9\cdot 2^{2k +2s +2m+2j-4} 
\end{align*}
\textbf{The case $j=m $}, similarly \\
 \begin{align*}
&  \Delta _{k}\omega_{2s + (m-i)+m} (s,m-i,  i) =  0 \quad \text{if \; $ 0\leq i\leq m -2 $ }\\
&  \Delta _{k}\omega_{2s + 1 + m} (s,1,  m-1) = 3\cdot2^{2k +2s + 2m-3}\quad \text{if \; i = m -1   }\\
& \text{We then deduce}\\
& \sum_{i = 0}^{m-1}2^{2i}\cdot \Delta _{k}\omega_{2s + (m-i)+m} (s,m-i,  i) 
 = \sum_{i = m-1}^{m-1}2^{2i}\cdot \Delta _{k}\omega_{2s + (m-i)+m} (s,m-i, i) \\
 & =  3\cdot2^{2k +2s + 4m -5}
 \end{align*}

  \end{proof}
  \subsection{\textbf{Computation of $ \Delta _{k}\Gamma_{2s +2m +j+1}^{\left[s\atop{ s +m\atop s+m} \right]} -   8\cdot\Delta _{k}\Gamma_{2s+ 2m+ j }^{\left[s\atop{ s +m-1\atop s +m } \right]}  \quad for\; 0\leq j \leq s-1 $}}
\label{subsec 10} 
  \begin{lem}
\label{lem 11.11}
We have for $m\geq 3$ the following reduction formula holding for $0\leq j\leq s-1.$\\
 \begin{align}
& \Delta _{k}\Gamma_{2s+2m+j+1}^{\left[s\atop{ s+m  \atop s+m  } \right]}
    -8\cdot\Delta _{k}\Gamma_{2s+2m+j}^{\left[s\atop{ s +m -1\atop s +m  } \right]} \label{eq 11.36}\\ 
    & = 2^{2m}\cdot\left[\Delta _{k}\Gamma_{2s +m+1+ j}^{\left[s \atop{ s \atop s +m} \right]}
  -8\cdot \Delta _{k}\Gamma_{2s+m+j}^{\left[s \atop{ s  \atop s +m -1} \right]}\right]
   + \sum_{i = 0}^{m-1}2^{2i}\cdot \Delta _{k}\omega_{2s +2m+j+1-i} (s,m-i,  i), \nonumber 
 \end{align}
 where 
\begin{equation}
\label{eq 11.37}
 \Delta _{k}\Gamma_{2s+m+1+j}^{\left[s\atop{ s \atop s+m} \right]} -   8\cdot\Delta _{k}\Gamma_{2s+m+j}^{\left[s\atop{ s \atop s+m-1} \right]}    =
 \begin{cases}
 -315\cdot2^{2k+2s+2m+4j-2}      & \text{if  } 0 \leq j \leq s-2,  \\
     -21\cdot2^{3k+3s+m-3} +21\cdot2^{2k+6s+2m-6}    & \text{if  } j=s-1,  
    \end{cases}
\end{equation}
\begin{align}
\label{eq 11.38}
\sum_{i = 0}^{m-1}2^{2i}\cdot \Delta _{k}\omega_{2s+2m+j+1-i} (s,m-i, i) =
 \begin{cases}
   0   & \text{if  } 2\leq j\leq s-2,  \\
   -21\cdot2^{3k+3s+2m-3} +21\cdot2^{2k+3s+3m-3}    & \text{if  } j = s-1.  
    \end{cases}
\end{align}
Combining \eqref{eq 11.36}, \eqref{eq 11.37} and \eqref{eq 11.38}  we obtain \\

\begin{equation}
\label{eq 11.39}
 \Delta _{k}\Gamma_{2s+2m+1+j}^{\left[s\atop{ s +m\atop s+m} \right]}
  -   8\cdot\Delta _{k}\Gamma_{2s+2m+j}^{\left[s\atop{ s +m-1\atop s+m} \right]}    =
 \begin{cases}
 -315\cdot2^{2k+2s +4m+4j-2}      & \text{if  } 0\leq j\leq s-2  \\
   -21\cdot2^{3k+3s+2m-3} +21\cdot2^{2k+6s+4m-6}  & \text{if  } j = s-1   
    \end{cases}
\end{equation}
 \end{lem}
   \begin{proof}
  \underline{Proof of \eqref{eq 11.36}} \\[0.5 cm]
  \eqref{eq 11.36} follows from \eqref{eq 11.32} with $j\longrightarrow m+j+1.$\\[0.5 cm]
  \underline{Proof of \eqref{eq 11.37}} \\[0.5 cm]
   \eqref{eq 11.37} is obtained by  \eqref{eq 11.20}, \eqref{eq 11.26} and \eqref{eq 11.31} with $l\longrightarrow m.$\\[0.5 cm]
    \underline{Proof of \eqref{eq 11.38}} \\[0.5 cm]
   \textbf{The case $0\leq j\leq s-2$}  \\[0.5 cm]
 \begin{align*}
 & \text{From \eqref{eq 7.51}}\\
  &  \Delta _{k}\omega_{2s + 2m +j +1-i} (s,m-i,  i) =  0 \\
  & \Updownarrow \\
  & m-i +2s +i+2   \leq  2s + 2m +j +1-i \leq  2(m-i) +3s + i -1\\
   & \Updownarrow \\
&  3\leq  i\leq \inf(m-1, m+j-1)=m-1 \\
 & \text{From \eqref{eq 7.50}}\\
  &  \Delta _{k}\omega_{2s + 2m +j -1} (s,m-2,  2) =  0 \\
  & \Updownarrow \\
  & m-2 +2s + 4   \leq  2s + 2m +j -1 \leq  2(m-2) +3s + 1\\
   & \Updownarrow \\
&\sup(0, -m+3) \leq j\leq s-2  \quad \text{recall that $m\geq 3$}\\
 & \text{From \eqref{eq 7.49}}\\
  &  \Delta _{k}\omega_{2s + 2m +j} (s,m-1,  1) =  0 \\
  & \Updownarrow \\
  & m-1 +2s + 3   \leq  2s + 2m +j  \leq  2(m-1) +3s \\
   & \Updownarrow \\
&\sup(0, -m+2) \leq j\leq s-2  \quad \text{recall that $m\geq 3$}\\
 & \text{From \eqref{eq 7.48}}\\
  &  \Delta _{k}\omega_{2s + 2m +j+1} (s,m,0) =  0 \\
  & \Updownarrow \\
  & m +2s +2   \leq  2s + 2m +j +1 \leq  2m +3s -1\\
   & \Updownarrow \\
&\sup(0, -m+1) \leq j\leq s-2  \quad \text{recall that $m\geq 3$}\\
& \text{We then deduce}\\
& \sum_{i = 0}^{m-1}2^{2i}\cdot \Delta _{k}\omega_{2s+2m+j+1-i} (s,m-i, i) = 0 \quad \text{if}\quad 0\leq j\leq s-2
 \end{align*}
  \textbf{The case $j=s-1$}  \\[0.5 cm]
  \begin{align*}
    & \text{Similar to the proof of \eqref{eq 11.38} in the case $0\leq j\leq s-2$}\\
&  \Delta _{k}\omega_{3s + 2m -i} (s,m-i,  i) =  \Delta _{k}\omega_{3s + 2(m-i) +i} (s,m-i,  i) = 21\cdot 2^{3k+3s +2m -i-3} \quad \text{if \; $ 0\leq i\leq m - 1 $   }\\
& \text{We then deduce}\\ 
& \sum_{i = 0}^{m-1}2^{2i}\cdot \Delta _{k}\omega_{3s+2m-i} (s,m-i, i) = \sum_{i = 0}^{m-1}2^{2i}\cdot  21\cdot 2^{3k+3s +2m -i-3}\\
& =  21\cdot2^{2k+3s +3m -3} -21\cdot2^{3k+3s +2m -3}
\end{align*}

  \end{proof}

   \subsection{\textbf{Formulas summary}}
  \label{subsec 11}
   \begin{lem}
\label{lem 11.12} 
We have the following formulas :\\
\begin{equation}
\label{eq 11.40}
 \Delta _{k}\Gamma_{s+j}^{\left[s\atop{ s \atop s+l} \right]} -   8\cdot\Delta _{k}\Gamma_{s+j -1}^{\left[s\atop{ s \atop s+l-1} \right]} 
= \begin{cases}
3\cdot2^{k+s-1} & \text{if   } j = 0,\; l\geq 0. \\
15\cdot2^{k+s-1} & \text{if   } j = 1,\;  l\geq 0.\\
 63\cdot2^{k+s+3j-6} &   \text{if  }\; 2\leq j\leq s-1, \: l\geq 0.
 \end{cases}
 \end{equation}
 We have for $l\geq 3:$
\begin{equation}
\label{eq 11.41}
 \Delta _{k}\Gamma_{2s+j}^{\left[s\atop{ s \atop s+l} \right]} -   8\cdot\Delta _{k}\Gamma_{2s+j -1}^{\left[s\atop{ s \atop s+l-1} \right]}    =
 \begin{cases}
 3\cdot2^{2k+2s-2} - 33\cdot2^{k+4s-6} & \text{if  } j = 0 \\
  -15\cdot2^{2k+2s-2}  +3\cdot2^{k+4s-3}   & \text{if  } j = 1 \\
   -9\cdot2^{2k+2s +2j-4}      & \text{if  } 2\leq j\leq l-1  \\
     -81\cdot2^{2k+2s +2j-4}      & \text{if  } j = l   
    \end{cases}
\end{equation}  
 \begin{equation}
\label{eq 11.42}
 \Delta _{k}\Gamma_{2s+l+1+j}^{\left[s\atop{ s  \atop s +l } \right]}
    -8\cdot\Delta _{k}\Gamma_{2s+ l+j}^{\left[s\atop{ s  \atop s +l -1 } \right]}=
   \begin{cases}
  -315\cdot2^{2k+2s+2l+4j-2} & \text{if  } 0\leq j\leq s-2  \\
  -21\cdot2^{3k+3s+l-3} + 21\cdot2^{2k+6s+2l-6}  & \text{if  } j = s-1 
   \end{cases}  
   \end{equation} 
\begin{equation}
\label{eq 11.43}
  \Delta _{k}\Gamma_{s+j}^{\left[s\atop{ s+m\atop s+m+l} \right]} =
 \begin{cases}
 2^{k+s-1} & \text{if   } j = 0 \\
  \gamma _{j}\cdot2^{k+s-1} = (21\cdot2^{3j-4} - 3\cdot2^{2j-3})\cdot2^{k+s-1} & \text{if } 1\leq j\leq m-1,\; l\geq 0\\
    21\cdot2^{k+s+3m-5} + 13\cdot2^{k+s+2m-4} & \text{if }j = m\geq 1,\;l \geq 1\\
    21\cdot2^{k+s+3m-5} + 45\cdot2^{k+s+2m-4} & \text{if }j =  m\geq 1,\;l = 0
  \end{cases}
\end{equation}
\begin{equation}
\label{eq 11.44}
  \Delta _{k}\Gamma_{s+m+1}^{\left[s\atop{ s+m\atop s+m} \right]} =
 \begin{cases}
147\cdot 2^{k+s-1} & \text{if   } m = 0 \\
315\cdot 2^{k+s} & \text{if   } m = 1 \\
  21\cdot2^{k+s+3m -2} + 273\cdot2^{k+s+2m- 2} & \text{if }  m\geq 2\;(holds \; for \; m\geq 0 )\\
  \end{cases}
  \end{equation}
\begin{equation}
\label{eq 11.45}
\Delta _{k}\Gamma_{s+m+j}^{\left[s\atop{ s+m \atop s+m+l} \right]} -   8\cdot\Delta _{k}\Gamma_{s+m+ j -1}^{\left[s\atop{ s+m \atop s+m+l-1} \right]}    =
 \begin{cases}
9\cdot2^{k+s-1} & \text{if   } j = 0,\; m = 1,\; l\geq 0 \\
19\cdot2^{k+s +2m-4} & \text{if   } j = 0,\; m \geq 2,\; l\geq 0 \\
31\cdot2^{k+s +2m -2} & \text{if   } j = 1, \; m\geq 1,\; l\geq 0 \\
 63\cdot2^{k+s+2m +3j-6} &   \text{if  }\;  2\leq j\leq s-1,\; m\geq 0, \; l\geq 0.
 \end{cases}
\end{equation}
We have for $m\geq 3:$
 \begin{equation}
\label{eq 11.46}
 \Delta _{k}\Gamma_{2s+m+j}^{\left[s\atop{ s +m\atop s+m} \right]}
  -   8\cdot\Delta _{k}\Gamma_{2s+m+j -1}^{\left[s\atop{ s +m-1\atop s+m} \right]}    =
 \begin{cases}
 -3\cdot2^{2k+2s+m-2} - 33\cdot2^{k+4s+2m-6} & \text{if  } j = 0 \\
  -69\cdot2^{2k+2s+m-2}  +3\cdot2^{k+4s+2m-3}   & \text{if  } j = 1 \\
   -63\cdot2^{2k+2s +m+3j-5}      & \text{if  } 2\leq j\leq m-1  \\
     -159\cdot2^{2k+2s +4m-5}      & \text{if  } j = m   
    \end{cases}
\end{equation} 
  \begin{equation}
\label{eq 11.47}
 \Delta _{k}\Gamma_{2s+2m+1+j}^{\left[s\atop{ s +m\atop s+m} \right]}
  -   8\cdot\Delta _{k}\Gamma_{2s+2m+j}^{\left[s\atop{ s +m-1\atop s+m} \right]}    =
 \begin{cases}
 -315\cdot2^{2k+2s +4m+4j-2}      & \text{if  } 0\leq j\leq s-2  \\
   -21\cdot2^{3k+3s+2m-3} +21\cdot2^{2k+6s+4m-6}  & \text{if  } j = s-1   
    \end{cases}
\end{equation}
  \end{lem}
 \begin{proof}
 We need only to summarize the principal  formulas obtained in the subsections \\
  11.3, 11.4,$\ldots, 11.10$
  \end{proof}

      \subsection{\textbf{Computation of $ \Delta _{k}\Gamma_{s+m+j}^{\left[s\atop{ s +m\atop s+m} \right]}\; for\;m\geq 1,\; 1\leq j\leq s-1  $}}
\label{subsec 12} 
\begin{lem}
\label{lem 11.12} 
We have the following reduction formula holding for $m\geq 1,\; 2\leq j\leq s-1$ :\\
   \begin{equation}
   \label{eq 11.48}
 \Delta _{k}\Gamma_{s+m+j}^{\left[s\atop{ s+m \atop s+m} \right]}
     = 2^{6}\cdot\Delta _{k}\Gamma_{s+(m-1)+ (j-1)}^{\left[s\atop{ s+ (m-1)\atop s+ (m-1)} \right]} +  189\cdot2^{k+s+2m+3j-6} 
\end{equation}
By induction on j we deduce from \eqref{eq 11.48}\\[0.1 cm]
 \begin{align}
 & (H_{j})\quad \Delta _{k}\Gamma_{s+m+j}^{\left[s\atop{ s+m \atop s+m} \right]}=
  21\cdot[2^{k+s+3m-2}\cdot2^{3(j-1)} + a_{j}\cdot2^{k+s+2m-3}\cdot2^{3(j-1)}] & \text{where\quad $a_{j} = 35\cdot2^{j-1} -9$}\nonumber \\
  & \Updownarrow  \nonumber\\
   & \Delta _{k}\Gamma_{s+m+j}^{\left[s\atop{ s+m \atop s+m} \right]} \nonumber\\
   & =  21\cdot[2^{k+s+3m +3j -5} +  35\cdot2^{k+s+2m +4j -7} - 9\cdot2^{k+s+2m +3j -6}] & \text{for\; $m\geq 1,\;1\leq j\leq s-1$}  \label{eq 11.49} 
\end{align}
\end{lem}

\begin{proof}
\textbf{proof of \eqref{eq 11.48}}\\[0.2 cm]
We have from  \eqref{eq 7.48} and \eqref{eq 8.3}, \eqref{eq 11.45} with $l=0$\\[0.2 cm]
\begin{align*}
 & \Delta _{k}\Gamma_{s+m+j}^{\left[s\atop{ s+m \atop s+m} \right]} -   8\cdot\Delta _{k}\Gamma_{s+m+ j -1}^{\left[s\atop{ s+m -1\atop s+m} \right]}
 = 63\cdot2^{k+s+2m+3j-6} \\
 & \Updownarrow \\
 &   4\cdot\Delta _{k}\Gamma_{s+m+j-1}^{\left[s\atop{ s+m-1 \atop s+m} \right]} -   32\cdot\Delta _{k}\Gamma_{s+m+ j -2}^{\left[s\atop{ s+m -1\atop s+m-1} \right]}+\Delta _{k}\omega_{s+m+j } (s,m, 0) =
  63\cdot2^{k+s+2m+3j-6} \\
  & \Updownarrow \\
   &   \Delta _{k}\Gamma_{s+m+j-1}^{\left[s\atop{ s+m-1 \atop s+m} \right]} -   8\cdot\Delta _{k}\Gamma_{s+m+ j -2}^{\left[s\atop{ s+m -1\atop s+m-1} \right]} =
  63\cdot2^{k+s+2m+3j-8} \\[0.1 cm]
   & \text{We then deduce  }\\
    & \Delta _{k}\Gamma_{s+m+j}^{\left[s\atop{ s+m \atop s+m} \right]}= 8\cdot\left[ 8\cdot\Delta _{k}\Gamma_{s+m+ j -2}^{\left[s\atop{ s+m -1\atop s+m-1} \right]} +  63\cdot2^{k+s+2m+3j-8}\right] 
    +  63\cdot2^{k+s+2m+3j-6}\\
    &  = 2^{6}\cdot\Delta _{k}\Gamma_{s+(m-1)+ (j-1)}^{\left[s\atop{ s+ (m-1)\atop s+ (m-1)} \right]} +  189\cdot2^{k+s+2m+3j-6} 
\end{align*}
\textbf{proof of \eqref{eq 11.49}}\\[0.2 cm]
Proof by induction on j\\[0.2 cm]
\begin{align*}
& \text{We obtain from \eqref{eq 11.44}}:\quad
  \Delta _{k}\Gamma_{s+m+1}^{\left[s\atop{ s+m \atop s+m} \right]}= 21\cdot[2^{k+s+3m-2} +26\cdot2^{k+s+2m-3}]
\quad \text{and so $ (H_{j}) $ holds for j=1}\\
& \text{We get by \eqref{eq 11.48} with j=2 using $(H_{1})$}:
  \Delta _{k}\Gamma_{s+m+2}^{\left[s\atop{ s+m \atop s+m} \right]}
   = 2^{6}\cdot\Delta _{k}\Gamma_{s+(m-1)+1}^{\left[s\atop{ s+ (m-1)\atop s+ (m-1)} \right]} +  189\cdot2^{k+s+2m} \\[0.2 cm]
 & = 2^{6}\cdot\big( 21\cdot[2^{k+s+3(m-1)-2} +26\cdot2^{k+s+2(m-1)-3}]\big) +189\cdot2^{k+s+2m} \\
 & = 21\cdot\big[2^{k+s+3m-2}\cdot2^{3(2-1)} + a_{2}\cdot2^{k+s+2m-3}\cdot2^{3(2-1)}\big]\\
  & =  21\cdot\big[2^{k+s+3m +1} +  35\cdot2^{k+s+2m +1} - 9\cdot2^{k+s+2m}\big] 
\quad \text{ and so $ (H_{j}) $ holds for j=2}\\
& \text{Assume now that $(H_{j-1})$ holds, we get by  \eqref{eq 11.48}}\\[0.2 cm]
 & \Delta _{k}\Gamma_{s+m+j}^{\left[s\atop{ s+m \atop s+m} \right]}
  = 2^{6}\cdot\Delta _{k}\Gamma_{s+(m-1)+ (j-1)}^{\left[s\atop{ s+ (m-1)\atop s+ (m-1)} \right]} +  189\cdot2^{k+s+2m+3j-6}\\
& =  2^{6}\cdot21\cdot[2^{k+s+3(m-1) +3(j-1) -5} +  35\cdot2^{k+s+2(m-1) +4(j-1) -7} - 9\cdot2^{k+s+2(m-1) +3(j-1) -6}]\\
 &  +  189\cdot2^{k+s+2m+3j-6} =  21\cdot\big[2^{k+s+3m +3j -5} +  35\cdot2^{k+s+2m +4j -7} - 9\cdot2^{k+s+2m +3j -6}\big]
 \end{align*}
\end{proof}
 \subsection{\textbf{Computation of $ \Delta _{k}\Gamma_{2s+m}^{\left[s\atop{ s +m\atop s+m} \right]}\; for\;m\geq 1 $}}
\label{subsec 13} 
\begin{lem}
\label{lem 11.13} 
We have : \\
\begin{equation}
\label{eq 11.50}
\Delta _{k}\Gamma_{2s+ m}^{\left[s\atop{ s+m\atop s+m} \right]} = 9\cdot2^{2k+2s+m-2} +21\cdot2^{k+4s+3m- 5}
    + 735\cdot2^{k+5s+2m- 7} -477\cdot2^{k+4s+2m- 6} 
\end{equation}
\end{lem}
\begin{proof}
\begin{align*}
 & \text{We obtain using  \eqref{eq 11.46} and \eqref{eq 7.48}}\\
 & \Delta _{k}\Gamma_{2s+m}^{\left[s\atop{ s+m \atop s+m} \right]} -   8\cdot\Delta _{k}\Gamma_{2s+m -1}^{\left[s\atop{ s+m -1\atop s+m} \right]}
 = -3\cdot2^{2k+2s+m-2}  -33\cdot2^{k+4s+2m-6}  \\
& \Updownarrow \\
&  4\cdot\Delta _{k}\Gamma_{2s+m-1}^{\left[s\atop{ s+m-1 \atop s+m} \right]} -   32\cdot\Delta _{k}\Gamma_{2s+m -2}^{\left[s\atop{ s+m -1\atop s+m-1} \right]}+\Delta _{k}\omega_{2s+m } (s,m, 0) =
 -3\cdot2^{2k+2s+m-2}  -33\cdot2^{k+4s+2m-6} \\
 & \Updownarrow \\
 &  \Delta _{k}\Gamma_{2s+m-1}^{\left[s\atop{ s+m-1 \atop s+m} \right]} -   8\cdot\Delta _{k}\Gamma_{2s+m -2}^{\left[s\atop{ s+m -1\atop s+m-1} \right]}\\
 & ={1\over 4}\cdot[ 9\cdot2^{2k+2s+m-2}   -3\cdot2^{2k+2s+m-2}  -33\cdot2^{k+4s+2m-6}]  =    3\cdot2^{2k+2s+m-3}  -33\cdot2^{k+4s+2m-8} \\ 
   & \text{We then deduce from \eqref{eq 11.49} with $m\longrightarrow m-1 (\geq 1)\; and \; j\longrightarrow  s-1$ }\\
    & \Delta _{k}\Gamma_{2s+m}^{\left[s\atop{ s+m \atop s+m} \right]}= 8\cdot\left[ 8\cdot\Delta _{k}\Gamma_{2s+m -2}^{\left[s\atop{ s+m -1\atop s+m-1} \right]} +    3\cdot2^{2k+2s+m-3}  -33\cdot2^{k+4s+2m-8}   \right] 
    -3\cdot2^{2k+2s+m-2}  -33\cdot2^{k+4s+2m-6}\\
     & \Leftrightarrow \Delta _{k}\Gamma_{2s+m}^{\left[s\atop{ s+m \atop s+m} \right]} = 2^{6}\cdot\Delta _{k}\Gamma_{s+(m-1)+ (s-1)}^{\left[s\atop{ s+ (m-1)\atop s+ (m-1)} \right]} +  9\cdot2^{2k+2s+m-2}  -99\cdot2^{k+4s+2m-6} \\
     & = 2^{6}\cdot21\cdot[ 2^{k+s+3(m-1) +3(s-1) -5} +  35\cdot2^{k+s+2(m-1) +4(s-1) -7} - 9\cdot2^{k+s+2(m-1) +3(s-1) -6}]\\
     &  + 9\cdot2^{2k+2s+m-2}  -99\cdot2^{k+4s+2m-6}\\
     & =   9\cdot2^{2k+2s+m-2} +21\cdot2^{k+4s+3m- 5} + 735\cdot2^{k+5s+2m- 7} -477\cdot2^{k+4s+2m- 6}
   \end{align*}  
\end{proof}
   \subsection{\textbf{Computation of $ \Delta _{k}\Gamma_{2s+m + j}^{\left[s\atop{ s +m\atop s+m} \right]} \; for \;m\geq 2,\; 1\leq j\leq m $}}
\label{subsec 14} 
\begin{lem} 
\label{lem 11.15} 
 We have  the following reduction formulas for $2\leq j\leq m-1$\\
\begin{align}
 &  \Delta _{k}\Gamma_{2s+m + j}^{\left[s\atop{ s+m \atop s+m} \right]} = 2^{6}\cdot\Delta _{k}\Gamma_{2s +(m-1)+(j-1)}^{\left[s\atop{ s+ (m-1)\atop s+ (m-1)} \right]} -  189\cdot2^{2k+2s+m +3j-5} \label{eq 11.51} \\
  &     (H_{j})\quad   \Delta _{k}\Gamma_{2s+m +j}^{\left[s\atop{ s+m \atop s+m} \right]}
 = 16^{2(j-1)}\cdot\Delta _{k-2(j-1)}\Gamma_{2s+ (m-j+1)+1}^{\left[s\atop{ s+ (m-j+1) \atop s+ (m-j+1)} \right]} \quad holds \; for \; 2\leq j\leq m \label{eq 11.52}
  \end{align}
  Applying \eqref{eq 11.51} and \eqref{eq 11.52} we obtain by  induction on j  the following formula 
\begin{equation}
\label{eq 11.53}
  \Delta _{k}\Gamma_{2s+ m +j}^{\left[s\atop{ s+m\atop s+m} \right]} 
 =  \begin{cases}
   63\cdot2^{2k+2s+ m +3j -5}  + 21\cdot2^{k+4s+3m +3j -5}\\
    +735\cdot2^{k+5s+2m +4j-7} - 945\cdot2^{k+4s+2m+ 4j -7}& \text{if } 1\leq j\leq m-1\\
    159\cdot2^{2k +2s +4m -5} +  735\cdot2^{k +5s +6m-7} -1629\cdot2^{k+4s +6m-7}& \text{if } j = m 
    \end{cases}
  \end{equation}
\end{lem}
\begin{proof}
\textbf{Proof of \eqref{eq 11.53} in the case $j=1$}
  \begin{align*}
 & \text{We have from  \eqref{eq 11.46} and \eqref{eq 7.48}}\\
 & \Delta _{k}\Gamma_{2s+m +1}^{\left[s\atop{ s+m \atop s+m} \right]} -   8\cdot\Delta _{k}\Gamma_{2s+m}^{\left[s\atop{ s+m -1\atop s+m} \right]}
 = -69\cdot2^{2k+2s+m-2}  + 3\cdot2^{k+4s+2m - 3}  \\
& \Updownarrow \\
&  4\cdot\Delta _{k}\Gamma_{2s+m}^{\left[s\atop{ s+m-1 \atop s+m} \right]} -   32\cdot\Delta _{k}\Gamma_{2s+m -1}^{\left[s\atop{ s+m -1\atop s+m-1} \right]}+\Delta _{k}\omega_{2s+m+1 } (s,m, 0) =
 -69\cdot2^{2k+2s+m-2}  + 3\cdot2^{k+4s+2m - 3} \\
 & \Updownarrow \\
 &   \Delta _{k}\Gamma_{2s+m}^{\left[s\atop{ s+m-1 \atop s+m} \right]} -   8\cdot\Delta _{k}\Gamma_{2s+m - 1}^{\left[s\atop{ s+m -1\atop s+m-1} \right]}\\
 & ={1\over 4}\cdot[- 9\cdot2^{2k+2s+m-2}   -  69\cdot2^{2k+2s+m-2}  + 3\cdot2^{k+4s+2m - 3}]  =    -39\cdot2^{2k+2s+m-3} + 3\cdot2^{k+4s+2m-5} \\ 
   & \text{We then deduce from \eqref{eq 11.50}with $m\longrightarrow m-1 (\geq 1)\; and \;j=s$ }\\
    & \Delta _{k}\Gamma_{2s+m+1}^{\left[s\atop{ s+m \atop s+m} \right]}= 8\cdot\left[ 8\cdot\Delta _{k}\Gamma_{2s +(m-1)}^{\left[s\atop{ s+m -1\atop s+m-1} \right]}  -  39\cdot2^{2k+2s+m-3} + 3\cdot2^{k+4s+2m-5} \right] 
 -69\cdot2^{2k+2s+m-2}  + 3\cdot2^{k+4s+2m - 3}   \\
  & \Updownarrow \\
  & \Delta _{k}\Gamma_{2s+m +1}^{\left[s\atop{ s+m \atop s+m} \right]} = 2^{6}\cdot\Delta _{k}\Gamma_{2s +(m-1)}^{\left[s\atop{ s+ (m-1)\atop s+ (m-1)} \right]} -  225\cdot2^{2k+2s+m-2}   + 9\cdot2^{k+4s+2m-3} \\
     & = 2^{6}\cdot[9\cdot2^{2k+2s+(m-1) -2} +  21\cdot2^{k+4s+3(m-1) - 5} +  735\cdot2^{k+5s+2(m-1) - 7}   - 477\cdot2^{k+4s+2(m-1)  -6}]\\
     & -  225\cdot2^{2k+2s+m-2}   + 9\cdot2^{k+4s+2m-3}  \\
     & =   63\cdot2^{2k+2s+m-2} +21\cdot2^{k+4s+3m- 2} + 735\cdot2^{k+5s+2m- 3} - 945\cdot2^{k+4s+2m- 3}
   \end{align*} 
\textbf{Proof of \eqref{eq 11.53} in the case $2\leq j\leq m-1$}    
 \begin{align*}
 & \text{We have from \eqref{eq 11.46} and \eqref{eq 8.3}}\\
 & \Delta _{k}\Gamma_{2s+m +j}^{\left[s\atop{ s+m \atop s+m} \right]} -   8\cdot\Delta _{k}\Gamma_{2s+m +j-1}^{\left[s\atop{ s+m -1\atop s+m} \right]}
 = -63\cdot2^{2k+2s+m +3j-5}   \\
& \Leftrightarrow   4\cdot\Delta _{k}\Gamma_{2s+m +j-1}^{\left[s\atop{ s+m-1 \atop s+m} \right]} -   32\cdot\Delta _{k}\Gamma_{2s+m +j-2}^{\left[s\atop{ s+m -1\atop s+m-1} \right]}+\Delta _{k}\omega_{2s+m+j } (s,m, 0) =
 -63\cdot2^{2k+2s+m +3j -5}   \\
  & \text{We get from \eqref{eq 7.48} }\\
  & \Delta _{k}\omega_{2s+m+j } (s,m, 0) = 0 \Leftrightarrow m+2s+2\leq 2s+m+j\leq 2m+3s-1\Leftrightarrow 2\leq j\leq m-1\\
  & \text{Then} \quad \Delta _{k}\Gamma_{2s+m +j-1}^{\left[s\atop{ s+m-1 \atop s+m} \right]} -   8\cdot\Delta _{k}\Gamma_{2s+m +j-2}^{\left[s\atop{ s+m -1\atop s+m-1} \right]}
  ={1\over 4}\cdot[- 63\cdot2^{2k+2s+m +3j-5}]  =    -63\cdot2^{2k+2s+m +3j-7}  \\ 
   & \text{We then deduce \eqref{eq 11.51} }\\
    & \Delta _{k}\Gamma_{2s+m+ j}^{\left[s\atop{ s+m \atop s+m} \right]}= 8\cdot\left[ 8\cdot\Delta _{k}\Gamma_{2s + m+j-2}^{\left[s\atop{ s+m -1\atop s+m-1} \right]}  -  63\cdot2^{2k+2s+m +3j-7}  \right] 
 -63\cdot2^{2k+2s+m +3j-5}    \\
  & \Updownarrow \\
 &  \Delta _{k}\Gamma_{2s+m + j}^{\left[s\atop{ s+m \atop s+m} \right]} = 2^{6}\cdot\Delta _{k}\Gamma_{2s +(m-1)+(j-1)}^{\left[s\atop{ s+ (m-1)\atop s+ (m-1)} \right]} -  189\cdot2^{2k+2s+m +3j-5}  \\
  \end{align*}   
    We shall now prove by induction on j that $(H_{j})$ holds for $2\leq j\leq m-1.$  
    \begin{equation*}
     (H_{j})\quad   \Delta _{k}\Gamma_{2s+m +j}^{\left[s\atop{ s+m \atop s+m} \right]}
 = 16^{2(j-1)}\cdot\Delta _{k-2(j-1)}\Gamma_{2s+ (m-j+1)+1}^{\left[s\atop{ s+ (m-j+1) \atop s+ (m-j+1)} \right]}\quad   \text{if }\; 2\leq j\leq m-1
\end{equation*}
\underline{$(H_{j}) $ holds for j = 2}\\[0.2 cm]

From \eqref{eq 11.51} with $j=2$ and \eqref{eq 11.53} with $j=1,\;m\longrightarrow m-1$ we obtain:\\
\begin{align*}
 &  \Delta _{k}\Gamma_{2s+m + 2}^{\left[s\atop{ s+m \atop s+m} \right]} =
  2^{6}\cdot\Delta _{k}\Gamma_{2s +(m-1)+ 1}^{\left[s\atop{ s+ (m-1)\atop s+ (m-1)} \right]} -  189\cdot2^{2k+2s+m + 1}  \\
& =  2^{6}\cdot[63\cdot2^{2k+2s+(m-1) -2} +  21\cdot2^{k+4s+3(m-1) - 2} +  735\cdot2^{k+5s+2(m-1) - 3}
 - 945\cdot2^{k+4s+2(m-1) - 3}] -  189\cdot2^{2k+2s+m + 1} \\
 & = 63\cdot2^{2k+2s+ m+1}  + 21\cdot2^{k+4s+3m +1} +735\cdot2^{k+5s+2m +1} - 945\cdot2^{k+4s+2m+1}\\
 & 16^{2}\cdot\Delta _{k-2}\Gamma_{2s+ (m- 1)+1}^{\left[s\atop{ s+ (m -1) \atop s+ (m -1)} \right]}\\
 & = 2^{8}\cdot[ 63\cdot2^{2(k-2)+2s+ (m-1) -2}  + 21\cdot2^{(k-2)+4s+3(m-1) -2} +735\cdot2^{(k-2)+5s+2(m-1) -3} - 945\cdot2^{(k-2)+4s+2(m-1) -3}]\\
 & =  63\cdot2^{2k+2s+ m+1}  + 21\cdot2^{k+4s+3m +1} +735\cdot2^{k+5s+2m +1} - 945\cdot2^{k+4s+2m+1}
\end{align*} 
 
 \underline{Assume  $(H_{j}) $ holds for j = r} \\[0.2 cm]
 We shall show that  $(H_{j}) $ holds for j = r +1\\[0.1 cm]
  \begin{align*}
  & \text{Indeed we have from \eqref{eq 11.53} in the case $m\rightarrow m-r,\;k\rightarrow k-2(r-1),\;j=1$}\\
  &  \text{ $(H_{j}) $ holds for j = r with $m \rightarrow m-1$}\\
  & \Updownarrow \\
  &  \Delta _{k}\Gamma_{2s+(m-1) +r}^{\left[s\atop{ s+(m-1) \atop s+(m-1)} \right]}
 = 16^{2(r-1)}\cdot\Delta _{k-2(r-1)}\Gamma_{2s+ ((m-1)-r+1)+1}^{\left[s\atop{ s+ ((m-1)-r+1) \atop s+ ((m-1)-r+1)} \right]}\\
 & = 2^{8r-8}\cdot[ 63\cdot2^{2(k-2(r-1))+2s+ ((m-1)-r+1) -2}  + 21\cdot2^{(k-2(r-1))+4s+3((m-1)-r+1) -2}\\
 &  +735\cdot2^{(k-2(r-1))+5s+2((m-1)-r+1) -3} - 945\cdot2^{(k-r+1)+4s+2((m-1)-r+1)-3}] \\
 & =   63\cdot2^{2k+2s+ m +3r -6}  + 21\cdot2^{k+4s+3m +3r -8} +735\cdot2^{k+5s+2m +4r-9} - 945\cdot2^{k+4s+2m+ 4r -9}\\
 & \text{Now from \eqref{eq 11.51} with $j=r+1$}\\
  & \Delta _{k}\Gamma_{2s+m + r+1}^{\left[s\atop{ s+m \atop s+m} \right]} = 2^{6}\cdot\Delta _{k}\Gamma_{2s +(m-1)+ r}^{\left[s\atop{ s+ (m-1)\atop s+ (m-1)} \right]} -  189\cdot2^{2k+2s+m +3(r+1)-5}  \\
  & =  2^{6}\cdot[ 63\cdot2^{2k+2s+ m +3r -6}  + 21\cdot2^{k+4s+3m +3r -8}\\
  &  +735\cdot2^{k+5s+2m +4r-9} - 945\cdot2^{k+4s+2m+ 4r -9}] -  189\cdot2^{2k+2s+m +3(r+1)-5} \\
  & =   63\cdot2^{2k+2s+ m +3r -2}  + 21\cdot2^{k+4s+3m +3r -2} +735\cdot2^{k+5s+2m +4r-3} - 945\cdot2^{k+4s+2m+ 4r -3}.\\
   & \text{On the other hand}\\
  &  16^{2((r+1)-1)}\cdot\Delta _{k-2((r+1)-1)}\Gamma_{2s+ (m-(r+1)+1)+1}^{\left[s\atop{ s+ (m-(r+1)+1) \atop s+ (m-(r+1)+1)} \right]}
    = 2^{8r +4}\cdot\Delta _{k-2r}\Gamma_{2s+ (m- r)+1}^{\left[s\atop{ s+ (m- r) \atop s+ (m- r)} \right]}\\
    & =  2^{8r}\cdot[  63\cdot2^{2(k-2r)+2s+(m-r)-2} +21\cdot2^{(k-2r)+4s+3(m-r)- 2} + 735\cdot2^{(k-2r)+5s+2(m-r)- 3} - 945\cdot2^{(k-2r)+4s+2(m-r)- 3}]\\
    & =  63\cdot2^{2k+2s+ m +3r -2}  + 21\cdot2^{k+4s+3m +3r -2} +735\cdot2^{k+5s+2m +4r-3} - 945\cdot2^{k+4s+2m+ 4r -3}\\
    & \text{Then we have for $ 2\leq j\leq m-1 $}\\
    & \Delta _{k}\Gamma_{2s+m +j}^{\left[s\atop{ s+m \atop s+m} \right]}\\
    & =  63\cdot2^{2k+2s+ m +3j -2}  + 21\cdot2^{k+4s+3m +3j -2} +735\cdot2^{k+5s+2m +4j-3} - 945\cdot2^{k+4s+2m+ 4j -3}
   \end{align*}   
   \textbf{Proof of \eqref{eq 11.53} in the case $j=m$}    
       \begin{align*}
 & \text{We have from \eqref{eq 11.46} with $j=m$}\\
 & \Delta _{k}\Gamma_{2s+2m}^{\left[s\atop{ s+m \atop s+m} \right]} -   8\cdot\Delta _{k}\Gamma_{2s+2m-1}^{\left[s\atop{ s+m -1\atop s+m} \right]}
 = -159\cdot2^{2k+2s+4m -5}   \\
& \Leftrightarrow   4\cdot\Delta _{k}\Gamma_{2s+2m -1}^{\left[s\atop{ s+m-1 \atop s+m} \right]} -   32\cdot\Delta _{k}\Gamma_{2s+2m -2}^{\left[s\atop{ s+m -1\atop s+m-1} \right]}+\Delta _{k}\omega_{2s+m+j } (s,m, 0) =
 -159\cdot2^{2k+2s+4m  -5}   \\
  & \text{Using \eqref{eq 7.48} we obtain}\\
  & \Delta _{k}\omega_{2s+2m } (s,m, 0) = 0 \Leftrightarrow m+2s+2\leq 2s+ 2m\leq 2m+3s-1\Leftrightarrow 2\leq m\\
  & \text{Then} \quad \Delta _{k}\Gamma_{2s+2m -1}^{\left[s\atop{ s+m-1 \atop s+m} \right]} -   8\cdot\Delta _{k}\Gamma_{2s+2m -2}^{\left[s\atop{ s+m -1\atop s+m-1} \right]}
  ={1\over 4}\cdot[- 159\cdot2^{2k+2s+4m -5}]  =    -159\cdot2^{2k+2s+4m -7}  \\ 
   & \text{We then deduce  }\\
    & \Delta _{k}\Gamma_{2s+2m}^{\left[s\atop{ s+m \atop s+m} \right]}= 8\cdot\left[ 8\cdot\Delta _{k}\Gamma_{2s +2m-2}^{\left[s\atop{ s+m -1\atop s+m-1} \right]}  -  159\cdot2^{2k+2s +4m -7}  \right] 
 -159\cdot2^{2k+2s+4m -5}    \\
  & \Leftrightarrow \Delta _{k}\Gamma_{2s+2m }^{\left[s\atop{ s+m \atop s+m} \right]} = 2^{6}\cdot\Delta _{k}\Gamma_{2s +2(m-1)}^{\left[s\atop{ s+ (m-1)\atop s+ (m-1)} \right]} -  477\cdot2^{2k+2s+4m -5}  
  \end{align*}
  \begin{align*}
  & \text{We obtain successively}\\
   &  \Delta _{k}\Gamma_{2s+2m }^{\left[s\atop{ s+m \atop s+m} \right]}
    = 2^{6}\cdot\Delta _{k}\Gamma_{2s +2(m-1)}^{\left[s\atop{ s+ (m-1)\atop s+ (m-1)} \right]} -  477\cdot2^{2k+2s+4m -5}\\
     &  2^{6}\cdot \Delta _{k}\Gamma_{2s+2(m-1) }^{\left[s\atop{ s+(m-1) \atop s+(m-1)} \right]}
    =    2^{6}\cdot\left[2^{6}\cdot\Delta _{k}\Gamma_{2s +2(m-2)}^{\left[s\atop{ s+ (m-2)\atop s+ (m-2)} \right]} -  477\cdot2^{2k+2s+4(m-1) -5}\right]\\  
     & ( 2^{6})^{2}\cdot \Delta _{k}\Gamma_{2s+2(m-2) }^{\left[s\atop{ s+(m-2) \atop s+(m-2)} \right]}
    =   ( 2^{6})^{2}\cdot\left[2^{6}\cdot\Delta _{k}\Gamma_{2s +2(m-3)}^{\left[s\atop{ s+ (m-3)\atop s+ (m-3)} \right]} -  477\cdot2^{2k+2s+4(m-2) -5}\right]\\  
    & \vdots \\
     & ( 2^{6})^{m-2}\cdot \Delta _{k}\Gamma_{2s+2(m-(m-2)) }^{\left[s\atop{ s+(m-(m-2)) \atop s+(m-(m-2))} \right]}
    =   ( 2^{6})^{m-2}\cdot\left[2^{6}\cdot\Delta _{k}\Gamma_{2s +2(m-(m-1))}^{\left[s\atop{ s+ (m-(m-1))\atop s+ (m-(m-1))} \right]} -  477\cdot2^{2k+2s+4(m-(m-2)) -5}\right]\\   
    & \text{By summing the above equations we get}\\
    &  \Delta _{k}\Gamma_{2s+2m }^{\left[s\atop{ s+m \atop s+m} \right]}= 
     ( 2^{6})^{m-1}\cdot\Delta _{k}\Gamma_{2(s+1)}^{\left[s\atop{ s+ 1\atop s+ 1} \right]}
      -477\cdot2^{2k+2s+4m-5}\cdot\sum_{i = 0}^{m-2}2^{2i}.\\
      & \text{From \eqref{eq 9.13}  with $j=s+1$ and $s\rightarrow s+1$ we obtain}\\
       & \Delta _{k}\Gamma_{2s+2m }^{\left[s\atop{ s+m \atop s+m} \right]} = 2^{6m-6}\cdot[159\cdot2^{2k +2s -1} +  735\cdot2^{k +5s -1} -1629\cdot2^{k+4s-1}] -  159\cdot2^{2k+2s+4m-5}\cdot(2^{2m-2} -1)\\
      & = 159\cdot2^{2k +2s +4m -5} +  735\cdot2^{k +5s +6m-7} -1629\cdot2^{k+4s +6m-7}\\
        & \text{$(H_{j})$ holds for j = m}\\
        & \text{Indeed we have}\\
        & 16^{2(m-1)}\cdot\Delta _{k-2(m-1)}\Gamma_{2s+ 2}^{\left[s\atop{ s+ 1 \atop s+ 1} \right]}\\
        & = 2^{8m-8}\cdot[ 159\cdot2^{2(k-2(m-1)) +2s -1} +  735\cdot2^{(k-2(m-1)) +5s -1} -1629\cdot2^{(k-2(m-1))+4s -1}]\\
        & =  159\cdot2^{2k +2s +4m -5} +  735\cdot2^{k +5s +6m-7} -1629\cdot2^{k+4s +6m-7}
  \end{align*}

  \end{proof}  

 \subsection{\textbf{Computation of $ \Delta _{k}\Gamma_{2s+2m+1+j}^{\left[s\atop{ s +m\atop s+m} \right]}\; for \; m\geq 2, \quad  0\leq j\leq s-1 $}}
\label{subsec 15}
\begin{lem}
\label{lem 11.16}
We have the following reduction formula for $0\leq j\leq s-1$\\[0.1 cm]
\begin{equation}
\label{eq 11.54}
(H_{j})\qquad   \Delta _{k}\Gamma_{2s+2m +j+1}^{\left[s\atop{ s+m \atop s+m} \right]}
 = 16^{2m+3j}\cdot\Delta _{k- (2m+3j)}\Gamma_{2(s-j)+1}^{\left[s -j\atop{ s -j \atop s -j} \right]},
\end{equation}
furthermore \\[0.1 cm]
\begin{align}
\label{eq 11.55}
 &  \Delta _{k}\Gamma_{2s+ 2m + 1+j}^{\left[s\atop{ s+m\atop s+m} \right]} \\
 &  = \begin{cases}
  315\cdot2^{2k+2s+4m +4j-2} + 105\cdot( 7\cdot2^{k +5s +6m +4j -3} -31\cdot2^{k+4s +6m +5j -3})    & \text{if } 0\leq j\leq s-2\\
     7\cdot2^{3k +3s +2m -3} - 21\cdot2^{2k+6s+4m -6} + 7\cdot2^{k+9s +6m-8}   & \text{if } j = s-1 
    \end{cases}\nonumber
 \end{align}
\end{lem}
\begin{proof}
\textbf{The case  $ 0\leq j\leq s-2 $}\\[0.2 cm]
 \begin{align*}
 & \text{We have from \eqref{eq 11.47} and  \eqref{eq 8.3} with $l=0$}\\
 & \Delta _{k}\Gamma_{2s+2m+1+j}^{\left[s\atop{ s+m \atop s+m} \right]} -   8\cdot\Delta _{k}\Gamma_{2s+2m+j}^{\left[s\atop{ s+m -1\atop s+m} \right]}
 = -315\cdot2^{2k+2s+4m +4j-2}   \\
& \Leftrightarrow   4\cdot\Delta _{k}\Gamma_{2s+2m +j}^{\left[s\atop{ s+m-1 \atop s+m} \right]} -  
 32\cdot\Delta _{k}\Gamma_{2s+2m +j-1}^{\left[s\atop{ s+m -1\atop s+m-1} \right]}+\Delta _{k}\omega_{2s+2m+1+j } (s,m, 0) =
 -315\cdot2^{2k+2s+4m +4j -2}   \\[0.1 cm]
  & \text{Now \eqref{eq 7.48} yields}\\[0.1 cm]
  & \Delta _{k}\omega_{2s+2m +1+j} (s,m, 0) = 0 \Leftrightarrow m+2s+2\leq 2s+ 2m +1+j\leq 2m+3s-1\\
  & \Leftrightarrow  \sup(1-m,0)\leq j\leq s-2 \\
  & \text{Then} \quad \Delta _{k}\Gamma_{2s+2m +j}^{\left[s\atop{ s+m-1 \atop s+m} \right]} -   8\cdot\Delta _{k}\Gamma_{2s+2m +j-1}^{\left[s\atop{ s+m -1\atop s+m-1} \right]}
  ={1\over 4}\cdot[- 315\cdot2^{2k+2s+4m +4j -2}]  =    -315\cdot2^{2k+2s+4m +4j-4}  \\[0.2 cm] 
   & \text{We then deduce  }\\[0.2 cm]
    & \Delta _{k}\Gamma_{2s+2m +1+j}^{\left[s\atop{ s+m \atop s+m} \right]}= 8\cdot\left[ 8\cdot\Delta _{k}\Gamma_{2s +2m +j-1}^{\left[s\atop{ s+m -1\atop s+m-1} \right]}  -  315\cdot2^{2k+2s +4m +4j-4}  \right] 
 -315\cdot2^{2k+2s+4m +4j -2}    \\
  & \Leftrightarrow \Delta _{k}\Gamma_{2s+2m +1+j}^{\left[s\atop{ s+m \atop s+m} \right]}
   = 2^{6}\cdot\Delta _{k}\Gamma_{2s +2(m-1)+j+1}^{\left[s\atop{ s+ (m-1)\atop s+ (m-1)} \right]} -  945\cdot2^{2k+2s+4m +4j-2}  
  \end{align*}
  \begin{align*}
  & \text{We obtain successively}\\
   &  \Delta _{k}\Gamma_{2s+2m +1+j}^{\left[s\atop{ s+m \atop s+m} \right]}
    = 2^{6}\cdot\Delta _{k}\Gamma_{2s +2(m-1) +j+1}^{\left[s\atop{ s+ (m-1)\atop s+ (m-1)} \right]} -  945\cdot2^{2k+2s+4m +4j-2}\\
     &  2^{6}\cdot \Delta _{k}\Gamma_{2s+2(m-1)+j+1 }^{\left[s\atop{ s+(m-1) \atop s+(m-1)} \right]}
    =    2^{6}\cdot\left[2^{6}\cdot\Delta _{k}\Gamma_{2s +2(m-2)+j+1}^{\left[s\atop{ s+ (m-2)\atop s+ (m-2)} \right]} -  945\cdot2^{2k+2s+4(m-1) +4j -2}\right]\\  
     & ( 2^{6})^{2}\cdot \Delta _{k}\Gamma_{2s+2(m-2) +j+1}^{\left[s\atop{ s+(m-2) \atop s+(m-2)} \right]}
    =   ( 2^{6})^{2}\cdot\left[2^{6}\cdot\Delta _{k}\Gamma_{2s +2(m-3)+j+1}^{\left[s\atop{ s+ (m-3)\atop s+ (m-3)} \right]} -  945\cdot2^{2k+2s+4(m-2) +4j -2}\right]\\  
    & \vdots \\
     & ( 2^{6})^{m-2}\cdot \Delta _{k}\Gamma_{2s+2(m-(m-2))+j+1 }^{\left[s\atop{ s+(m-(m-2)) \atop s+(m-(m-2))} \right]}\\
      &  =   ( 2^{6})^{m-2}\cdot\left[2^{6}\cdot\Delta _{k}\Gamma_{2s +2(m-(m-1))+j+1}^{\left[s\atop{ s+ (m-(m-1))\atop s+ (m-(m-1))} \right]} -  945\cdot2^{2k+2s+4(m-(m-2)) +4j -2}\right]\\   
    & \text{By summing the above equations we get}\\
    &  \Delta _{k}\Gamma_{2s+2m +j+1}^{\left[s\atop{ s+m \atop s+m} \right]}= 
     ( 2^{6})^{m-1}\cdot\Delta _{k}\Gamma_{2(s+1)+j+1}^{\left[s\atop{ s+ 1)\atop s+ 1)} \right]}
      -945\cdot2^{2k+2s+4m +4j-2}\cdot\sum_{i = 0}^{m-2}2^{2i}\\
       & = 2^{6m-6}\cdot\big[105\cdot(3\cdot2^{2k +2s +4j+2} +  7\cdot2^{k +5s +4j+3} -31\cdot2^{k+4s +5j+3}\big]\\
       & - 315\cdot2^{2k+2s+4m +4j-2}\cdot(2^{2m-2} -1)\\
       & = 315\cdot2^{2k+2s+4m +4j-2} + 105\cdot( 7\cdot2^{k +5s +6m +4j -3} -31\cdot2^{k+4s +6m +5j -3})
      \end{align*}
      On the other hand we have by \eqref{eq 9.7} with $k\longrightarrow k-(2m+3j),\;s\longrightarrow s-j,\; j=1$\\
    \begin{align*}
   & 16^{2m+3j}\cdot\Delta _{k- (2m+3j)}\Gamma_{2(s-j)+1}^{\left[s -j\atop{ s -j \atop s -j} \right]} =
    2^{8m+12j}\cdot\Delta _{k- (2m+3j)}\Gamma_{2(s-j)+1}^{\left[s -j\atop{ s -j \atop s -j} \right]}\\
    & =   2^{8m+12j}\cdot\big[105\cdot(3\cdot2^{2(k- (2m+3j))+2(s-j) -2} + 7\cdot2^{(k- (2m+3j))+ 5(s-j) -3}\\
   &  - 31\cdot2^{(k- (2m+3j))+4(s-j) -3})\big] \\
    & = 315\cdot2^{2k+2s+4m +4j-2} + 105\cdot( 7\cdot2^{k +5s +6m +4j -3} -31\cdot2^{k+4s +6m +5j -3})
    \end{align*}
    Thus we get \eqref{eq 11.54} and \eqref{eq 11.55} in the case $0\leq j\leq s-2$\\[0.2 cm]
    
   \textbf{The case  $ j=s-1 $}\\[0.2 cm] 
     \begin{align*}
 & \text{We have by \eqref{eq 11.47} and  \eqref{eq 8.3} with $l=0$}\\
 & \Delta _{k}\Gamma_{3s+2m}^{\left[s\atop{ s+m \atop s+m} \right]} -   8\cdot\Delta _{k}\Gamma_{3s+2m -1}^{\left[s\atop{ s+m -1\atop s+m} \right]}
 = -21\cdot2^{3k+3s+2m -3} +   21\cdot2^{2k+6s+4m -6} \\
 & \Leftrightarrow   4\cdot\Delta _{k}\Gamma_{3s+2m -1}^{\left[s\atop{ s+m-1 \atop s+m} \right]} -  
 32\cdot\Delta _{k}\Gamma_{3s+2m -2}^{\left[s\atop{ s+m -1\atop s+m-1} \right]}+\Delta _{k}\omega_{3s+2m } (s,m, 0) =
  -21\cdot2^{3k+3s+2m -3} +   21\cdot2^{2k+6s+4m -6}  \\
  & \text{We get by \eqref{eq 7.48} }\\
  & \Delta _{k}\omega_{3s+2m } (s,m, 0) =  21\cdot2^{3k+3s+2m-3}\\
   & \text{Then} \\
   & \Delta _{k}\Gamma_{3s+2m -1}^{\left[s\atop{ s+m-1 \atop s+m} \right]} -   8\cdot\Delta _{k}\Gamma_{3s+2m -2}^{\left[s\atop{ s+m -1\atop s+m-1} \right]}
  ={1\over 4}\cdot[- 21\cdot2^{3k+3s+2m-3} -21\cdot2^{3k+3s+2m -3} +   21\cdot2^{2k+6s+4m -6} ]\\
  &  =  -21\cdot2^{3k+3s+2m -4} +   21\cdot2^{2k+6s+4m -8} \\ 
   & \text{We then deduce  }\\
    & \Delta _{k}\Gamma_{3s+2m}^{\left[s\atop{ s+m \atop s+m} \right]}= 8\cdot\left[ 8\cdot\Delta _{k}\Gamma_{3s +2m -2}^{\left[s\atop{ s+m -1\atop s+m-1} \right]}  -21\cdot2^{3k+3s+2m -4} +   21\cdot2^{2k+6s+4m -8}  \right] \\
    &  -21\cdot2^{3k+3s+2m -3} +   21\cdot2^{2k+6s+4m -6}   \\
  & \Leftrightarrow \Delta _{k}\Gamma_{3s+2m }^{\left[s\atop{ s+m \atop s+m} \right]}
   = 2^{6}\cdot\Delta _{k}\Gamma_{3s +2(m-1)}^{\left[s\atop{ s+ (m-1)\atop s+ (m-1)} \right]} - 105\cdot2^{3k+3s+2m -3} +   63\cdot2^{2k+6s+4m -6} 
  \end{align*}
    \begin{align*}
  & \text{We obtain successively}\\
   &  \Delta _{k}\Gamma_{3s+2m }^{\left[s\atop{ s+m \atop s+m} \right]}
    = 2^{6}\cdot\Delta _{k}\Gamma_{3s +2(m-1)}^{\left[s\atop{ s+ (m-1)\atop s+ (m-1)} \right]} - 105\cdot2^{3k+3s+2m -3} +   63\cdot2^{2k+6s+4m -6} \\
     &  2^{6}\cdot \Delta _{k}\Gamma_{2s+2(m-1) }^{\left[s\atop{ s+(m-1) \atop s+(m-1)} \right]}
    =    2^{6}\cdot\left[2^{6}\cdot\Delta _{k}\Gamma_{2s +2(m-2)}^{\left[s\atop{ s+ (m-2)\atop s+ (m-2)} \right]}  - 105\cdot2^{3k+3s+2(m-1) -3} +   63\cdot2^{2k+6s+4(m-1) -6}  \right]\\  
     & ( 2^{6})^{2}\cdot \Delta _{k}\Gamma_{3s+2(m-2) }^{\left[s\atop{ s+(m-2) \atop s+(m-2)} \right]}
    =   ( 2^{6})^{2}\cdot\left[2^{6}\cdot\Delta _{k}\Gamma_{3s +2(m-3)}^{\left[s\atop{ s+ (m-3)\atop s+ (m-3)} \right]}   - 105\cdot2^{3k+3s+2(m-1) -3} +   63\cdot2^{2k+6s+4(m-1) -6}\right]\\  
    & \vdots \\
     & ( 2^{6})^{m-2}\cdot \Delta _{k}\Gamma_{3s+2(m-(m-2)) }^{\left[s\atop{ s+(m-(m-2)) \atop s+(m-(m-2))} \right]}\\
      &  =   ( 2^{6})^{m-2}\cdot\left[2^{6}\cdot\Delta _{k}\Gamma_{3s +2(m-(m-1))}^{\left[s\atop{ s+ (m-(m-1))\atop s+ (m-(m-1))} \right]} 
      - 105\cdot2^{3k+3s+2(m-(m-2)) -3} +   63\cdot2^{2k+6s+4(m-(m-2)) -6} \right]\\   
    & \text{By summing the above equations we get}\\
    &  \Delta _{k}\Gamma_{3s+2m }^{\left[s\atop{ s+m \atop s+m} \right]}= 
     ( 2^{6})^{m-1}\cdot\Delta _{k}\Gamma_{3s+2}^{\left[s\atop{ s+ 1)\atop s+ 1)} \right]} +
      \sum_{i = 0}^{m-2}[   - 105\cdot2^{3k+3s+2(m- i) -3} +   63\cdot2^{2k+6s+4(m- i) -6} ]\cdot2^{6i} \\
      & = 2^{6m-6}\cdot\Delta _{k}\Gamma_{3s+2}^{\left[s\atop{ s+ 1)\atop s+ 1)} \right]} - 105\cdot2^{3k +3s +2m -3}\cdot \sum_{i = 0}^{m-2}2^{4i}
      +  63\cdot2^{2k+6s+4m -6}\cdot\sum_{i = 0}^{m-2}2^{2i}\\
      & =  2^{6m-6}\cdot\big[7\cdot2^{3k+3s -1} -21\cdot2^{2k+6s -2} + 7\cdot2^{k+9s -2}\big]\\
      &  - 7\cdot2^{3k +3s +2m -3}\cdot(2^{4m-4} -1)  + 21\cdot2^{2k+6s+4m -6}\cdot(2^{2m-2} -1) \\
      & =  7\cdot2^{3k +3s +2m -3} - 21\cdot2^{2k+6s+4m -6} + 7\cdot2^{k+9s +6m-8}
      \end{align*}
    \begin{align*}
&  \text{$ (H_{s-1})$  holds } \Longleftrightarrow \Delta _{k}\Gamma_{3s+2m }^{\left[s\atop{ s+m \atop s+m} \right]}
 = 16^{2m+3(s-1)}\cdot\Delta _{k- (2m+3(s-1))}\Gamma_{3}^{\left[1\atop{ 1 \atop 1} \right]} \\
 & \text{Indeed we deduce  from  Lemma \ref{lem 1.23} in the case $m=0$} \\
  & 16^{2m+3(s-1)}\cdot\Delta _{k- (2m+3(s-1))}\Gamma_{3}^{\left[ 1\atop{ 1 \atop 1} \right]} =
    2^{8m+ 12s-12}\cdot\Delta _{k- (2m+3(s-1))}\Gamma_{3}^{\left[1\atop{ 1 \atop 1} \right]}\\[0.1 cm]
    & =   2^{8m+12s -12}\cdot\big[7\cdot2^{3( k-2m -3s+3)} -21\cdot2^{2( k-2m -3s+3)} +   7\cdot2^{ k-2m -3s+3 }\big ] \\[0.1 cm]
    & = 7\cdot2^{3k +3s+2m-3)} -21\cdot2^{2k +6s +4m -6} +   7\cdot2^{ k +9s +6m -8}
\end{align*}  
        \end{proof}
  \subsection{\textbf{Summary of principal results in section 11}}
\label{subsec 16} 
\begin{lem}
  \label{lem 11.17} 
  We have for $m\geq 2$:\\
   \begin{equation}
\label{eq 11.56}
   \Delta _{k}\Gamma_{s+j}^{\left[s\atop{ s+m\atop s+m} \right]} =
 \begin{cases}
 2^{k+s-1} & \text{if   } j = 0 \\
  (21\cdot2^{3j-4} - 3\cdot2^{2j-3})\cdot2^{k+s-1} & \text{if } 1\leq j\leq m-1 \\
  21\cdot2^{k+s+3m-5} + 45\cdot2^{k+s+2m-4} & \text{if }j =  m 
     \end{cases}
  \end{equation}
   \begin{equation}
\label{eq 11.57}
  \Delta _{k}\Gamma_{s+ m +j}^{\left[s\atop{ s+m\atop s+m} \right]} =
 \begin{cases}
 21\cdot[2^{k+s+3m +3j -5} +  35\cdot2^{k+s+2m +4j -7} - 9\cdot2^{k+s+2m +3j -6}]  & \text{if } 1\leq j\leq s-1\\
   9\cdot2^{2k+2s+m-2} +21\cdot2^{k+4s+3m- 5}\\
    + 735\cdot2^{k+5s+2m- 7} -477\cdot2^{k+4s+2m- 6} & \text{if } j = s
  \end{cases}
 \end{equation}
\begin{equation}
\label{eq 11.58}
  \Delta _{k}\Gamma_{2s+ m +1+j}^{\left[s\atop{ s+m\atop s+m} \right]} =
 \begin{cases}
   63\cdot2^{2k+2s+ m +3j -2}  + 21\cdot2^{k+4s+3m +3j -2}\\
   \ +735\cdot2^{k+5s+2m +4j-3} - 945\cdot2^{k+4s+2m+ 4j -3}& \text{if } 0\leq j\leq m-2\\
    159\cdot2^{2k +2s +4m -5} +  735\cdot2^{k +5s +6m-7} -1629\cdot2^{k+4s +6m-7}& \text{if } j = m-1 
    \end{cases}
  \end{equation}
 \begin{equation}
\label{eq 11.59}
   \Delta _{k}\Gamma_{2s+ 2m + 1+j}^{\left[s\atop{ s+m\atop s+m} \right]} =
 \begin{cases}
  315\cdot2^{2k+2s+4m +4j-2}\\
   + 105\cdot( 7\cdot2^{k +5s +6m +4j -3} -31\cdot2^{k+4s +6m +5j -3})    & \text{if } 0\leq j\leq s-2\\
     7\cdot2^{3k +3s +2m -3} - 21\cdot2^{2k+6s+4m -6} + 7\cdot2^{k+9s +6m-8}   & \text{if } j = s-1 
    \end{cases}
  \end{equation}
 Further we have the following reduction formulas \\[0.2 cm]
 \begin{align}
   &  \Delta _{k}\Gamma_{2s+m +j}^{\left[s\atop{ s+m \atop s+m} \right]}
 = 16^{2(j-1)}\cdot\Delta _{k-2(j-1)}\Gamma_{2s+ (m-j+1)+1}^{\left[s\atop{ s+ (m-j+1) \atop s+ (m-j+1)} \right]}
 \quad   \text{if }\; 2\leq j\leq m \label{eq 11.60}\\[0.1 cm]
  &   \Delta _{k}\Gamma_{2s+2m +j+1}^{\left[s\atop{ s+m \atop s+m} \right]}
 = 16^{2m+3j}\cdot\Delta _{k- (2m+3j)}\Gamma_{2(s-j)+1}^{\left[s -j\atop{ s -j \atop s -j} \right]}
 \quad   \text{if }\; 0 \leq j\leq s-1  \label{eq 11.61}\\
 & \nonumber
  \end{align}
  
 \textbf{The case $m=1$}\\[0.2 cm]
\begin{equation}
\label{eq 11.62}
 \Delta _{k}\Gamma_{s+j}^{\left[s\atop{ s+1\atop s+1} \right]}=
 \begin{cases}
2^{k+s- 1}  &\text{if  }  j = 0 \\
33\cdot2^{k+s-1}& \text{if   } j = 1\\
7\cdot(105\cdot2^{j-2} - 15)\cdot2^{k+s+3j-7}  & \text{if   } 2\leq j\leq s\\
9\cdot2^{2k+2s-1} + 21\cdot35\cdot2^{k+5s-5} - 393\cdot2^{k+4s-4} & \text{if   } j = s +1\\
159\cdot2^{2k+2s-1} + 735\cdot2^{k+5s- 1} - 1629\cdot2^{k+4s-1} & \text{if   } j = s +2
\end{cases}
\end{equation}

 \begin{equation}
\label{eq 11.63}
 \Delta _{k}\Gamma_{2s+1+j}^{\left[s\atop{ s+1\atop s+1} \right]}=
 \begin{cases}
105\cdot(3\cdot2^{2k+2s +4j-6} + 7\cdot2^{k+5s+4j-5} - 31\cdot2^{k+4s +5j-7}) & \text{if   } 2\leq j\leq s\\
7\cdot2^{3k+3s-1} -21\cdot2^{2k+6s- 2} + 7\cdot2^{k+9s-2} & \text{if   } j = s +1
\end{cases}
\end{equation}
  \end{lem}  
 \begin{proof} 
 We obtain \eqref{eq 11.62} and \eqref{eq 11.63} by combining 
 \eqref{eq 8.9},  \eqref{eq 8.10} and \eqref{eq 9.5}, \eqref{eq 9.7}.
  
 \end{proof}    
 \section{\textbf{An explicit formula for $ \Gamma_{i}^{\left[s\atop{ s +m\atop s+m} \right]\times k} , \quad for\; s \leq i \leq 3s+2m,\; k\geq i $}}
\label{sec 12} 
    \subsection{Notation}
  \label{subsec 1}
  \begin{defn}
  \label{defn 12.1}
We recall that  $  \Gamma_{i}^{\left[s\atop{ s+m \atop s+m} \right]\times k} $ denotes the number  of rank i  matrices of the form 
   $\left[{A\over{B \over C}}\right] $   where A is  a  $ s \times k $ persymmetric matrix,  B a  $ (s+m) \times k $ persymmetric matrix and C is a  $ (s+m) \times k $ persymmetric matrix.

  \end{defn}
   \subsection{Introduction}
  \label{subsec 2}
We adapt the method in Section \ref{sec 10} to compute explicitly  the number $  \Gamma_{i}^{\left[s\atop{ s+m \atop s+m} \right]\times k} .$

 \subsection{\textbf{Computation of $ \Gamma_{s+i}^{\left[s\atop{ s +m\atop s+m} \right]\times k} -\Gamma_{s+i}^{\left[s\atop{ s +m\atop s+m} \right]\times (s+i+1)} , \quad for\; 0 \leq i \leq 2s+2m,\; k\geq s+i+1 $}}
  \label{subsec 3}
 \begin{lem}
 \label{lem 12.2}
 We have  for $m\geq 2,\;k\geq s+j+1$ :\\[0.2 cm]
  \begin{align}
& \Gamma_{s+j}^{\left[s\atop{ s +m\atop s+m} \right]\times k} -\Gamma_{s+j}^{\left[s\atop{ s +m\atop s+m} \right]\times (s+j+1)}\label{eq 12.1} \\
& =  \begin{cases}
2^{s-1}\cdot\big(2^{k}-2^{s+1}\big)  &\text{if  }  j = 0 \\
\big(21\cdot2^{s+3j-5}-3\cdot2^{s+2j-4}\big)\cdot\big(2^{k}-2^{s+j+1}\big)& \text{if   }\; 1\leq j\leq m-1\\
\big(21\cdot2^{s+3m-5} +45\cdot2^{s+2m-4}\big)\cdot\big(2^{k}-2^{s+m+1}\big)  & \text{if   } \; j=m, \nonumber\\
\end{cases} \\
& \nonumber
\end{align}
  for $m\geq 2,\;k\geq s+m+j+1$ :\\[0.2 cm]
  \begin{align}
& \Gamma_{s+m+j}^{\left[s\atop{ s +m\atop s+m} \right]\times k} -\Gamma_{s+m+j}^{\left[s\atop{ s +m\atop s+m} \right]\times (s+m+j+1)}\label{eq 12.2} \\
& =  \begin{cases}
 \big(2^{k}-2^{s+m+j+1}\big)\cdot\big(21\cdot2^{s+3m+3j-5}+735\cdot2^{s+2m+4j-7}-189\cdot2^{s+2m+3j-6}\big) &\text{if  }1\leq j\leq s-1\\
 3\cdot2^{2s+m-2}\cdot\big(2^{2k}-2^{4s+2m+2}\big) \\
+   \big(2^{k}-2^{2s+m+1}\big)\cdot\big(21\cdot2^{4s+3m-5}+735\cdot2^{5s+2m-7}-477\cdot2^{4s+2m-6}\big)   &\text{if  }  j=s, \nonumber\\
\end{cases}\\
& \nonumber
\end{align}
  for $m\geq 2,\;k\geq 2s+m+j+1$ :\\[0.2 cm]
  \begin{align}
& \Gamma_{2s+m+1+j}^{\left[s\atop{ s +m\atop s+m} \right]\times k} -\Gamma_{2s+m+1+j}^{\left[s\atop{ s +m\atop s+m} \right]\times (2s+m+2+j)}\label{eq 12.3} \\
& =  \begin{cases}
 \big(2^{k}-2^{2s+m+2+j}\big)\cdot\big(21\cdot2^{4s+3m+3j-2}+735\cdot2^{5s+2m+4j-3}-945\cdot2^{4s+2m+4j-3}\big) \\
+ 21\cdot2^{2s+m+3j-2}\cdot\big(2^{2k}-2^{4s+2m+2j+4}\big)  &\text{if  }0\leq j\leq m-2\\
  53\cdot2^{2s+4m-5}\cdot\big(2^{2k}-2^{4s+4m+2}\big) \\
    +   \big(2^{k}-2^{2s+2m+1}\big)\cdot\big(735\cdot2^{5s+6m-7}-1629\cdot2^{4s+6m-7}\big)   &\text{if  }  j=m-1, \nonumber\\
\end{cases}\\
& \nonumber
\end{align}

  for $m\geq 2,\;k\geq 2s+2m+j+1$ :\\[0.2 cm]
  \begin{align}
& \Gamma_{2s+2m+1+j}^{\left[s\atop{ s +m\atop s+m} \right]\times k} -\Gamma_{2s+2m+1+j}^{\left[s\atop{ s +m\atop s+m} \right]\times (2s+2m+2+j)}\label{eq 12.4} \\
& =  \begin{cases}
 \big(2^{k}-2^{2s+2m+2+j}\big)\cdot\big(735\cdot2^{5s+6m+4j-3}-3255\cdot2^{4s+6m+5j-3}\big) \\
+ 105\cdot2^{2s+4m+4j-2}\cdot\big(2^{2k}-2^{4s+4m+2j+4}\big)  &\text{if  }0\leq j\leq s-2\\
2^{3s+2m-3}\cdot\big(2^{3k}-2^{9s+6m+3}\big) - 7\cdot2^{6s+4m-6}\cdot\big(2^{2k}-2^{6s+4m+2}\big)  \\
+ 7\cdot2^{9s+6m-8}\cdot\big(2^{k}-2^{3s+2m+1}\big) &\text{if  } j=s-1\\
\end{cases}
& \nonumber
\end{align}
We have for $ m=1,\; k\geq s+j+1  $\\[0.2 cm]
  \begin{align}
& \Gamma_{s+j}^{\left[s\atop{ s +1\atop s+1} \right]\times k} -\Gamma_{s+j}^{\left[s\atop{ s +1\atop s+1} \right]\times (s+j+1)}\label{eq 12.5} \\
& =  \begin{cases}
2^{s-1}\cdot\big(2^{k}-2^{s+1}\big)  &\text{if  }  j = 0 \\
33\cdot2^{s-1}\cdot\big(2^{k}-2^{s+2}\big)  &\text{if  }  j = 1 \\
\big(735\cdot2^{s+4j-9}-105\cdot2^{s+3j-7}\big)\cdot\big(2^{k}-2^{s+j+1}\big)& \text{if   }\; 2\leq j\leq s\\
3\cdot2^{2s-1}\cdot(2^{2k}-2^{4s+4}) +
\big(735\cdot2^{5s-5}-393\cdot2^{4s-4}\big)\cdot\big(2^{k}-2^{2s+2}\big)  & \text{if   } \; j=s+1, \\
53\cdot2^{2s-1}\cdot(2^{2k}-2^{4s+6}) +
\big(735\cdot2^{5s-1}-1629\cdot2^{4s-1}\big)\cdot\big(2^{k}-2^{2s+3}\big)  & \text{if   } \; j=s+2, \\
\end{cases}\nonumber \\
& \nonumber
\end{align}

We have for $ m=1,\; k\geq 2s+2+j  $\\[0.2 cm]
  \begin{align}
& \Gamma_{2s+1+j}^{\left[s\atop{ s +1\atop s+1} \right]\times k} -\Gamma_{2s+1+j}^{\left[s\atop{ s +1\atop s+1} \right]\times (2s+2+j)}\label{eq 12.6} \\
& =  \begin{cases}
105\cdot2^{2s+4j-6}\cdot(2^{2k}-2^{4s+2j+4}) +
\big(735\cdot2^{5s+4j-5}-3255\cdot2^{4s+5j-7}\big)\cdot\big(2^{k}-2^{2s+j+2}\big)  & \text{if   } \; 2\leq j\leq s, \\
2^{3s-1}\cdot(2^{3k}-2^{9s+6}) -7\cdot2^{6s-2}\cdot\big(2^{2k}-2^{6s+6}\big)  
+ 7\cdot2^{9s-2}\cdot\big(2^{k}-2^{3s+3}\big) &\text{if  } j=s+1\\
 \end{cases}\nonumber \\
& \nonumber
\end{align}
 \end{lem}
\begin{proof}
By way of example we prove \eqref{eq 12.2} in the case $j=s$.\\[0.1 cm]
 We obtain from \eqref{eq 11.57} in the case $j=s$ \\[0.1 cm]
  \begin{align*}
&  \Delta _{k}\Gamma_{2s+m}^{\left[s\atop{ s+m\atop s+m} \right]} = 
  \Gamma_{2s+m}^{\left[s\atop{ s +m\atop s+m} \right]\times (k+1)} -\Gamma_{2s+m}^{\left[s\atop{ s +m\atop s+m} \right]\times k} \\[0.1 cm]
  &  =  9\cdot2^{2k+2s+m-2} +21\cdot2^{k+4s+3m- 5} + 735\cdot2^{k+5s+2m- 7} -477\cdot2^{k+4s+2m- 6}, \\[0.1 cm]
  & \text{hence the above equality gives}\\[0.1 cm]
  & \sum_{i=2s+m+1}^{k}\Big(\Gamma_{2s+m}^{\left[s\atop{ s +m\atop s+m} \right]\times (i+1)} -\Gamma_{2s+m}^{\left[s\atop{ s +m\atop s+m} \right]\times i}\Big )
  =  9\cdot2^{2s+m-2}\cdot\sum_{i=2s+m+1}^{k}2^{2i}\\
   &   + \Big(21\cdot2^{4s+3m- 5} + 735\cdot2^{5s+2m- 7} -477\cdot2^{4s+2m- 6}\Big)\cdot\sum_{i=2s+m+1}^{k}2^{i}\\[0.1 cm]
   & \Updownarrow \\[0.1 cm]
&   \Gamma_{2s+m}^{\left[s\atop{ s +m\atop s+m} \right]\times (k+1)} -\Gamma_{2s+m}^{\left[s\atop{ s +m\atop s+m} \right]\times (2s+m+1)}
  =  9\cdot2^{2s+m-2}\cdot2^{4s+2m+2}\cdot\bigg({(2^2)^{k-2s-m}-1\over 2^2-1} \bigg)\\[0.1 cm]
  & + \Big(21\cdot2^{4s+3m- 5} + 735\cdot2^{5s+2m- 7} -477\cdot2^{4s+2m- 6}\Big)\cdot2^{2s+m+1}\cdot\big( 2^{k-2s-m}-1\big)
   \end{align*}
\end{proof}

 \subsection{\textbf{Computation of $ \Gamma_{s+j}^{\left[s\atop{ s +m\atop s+m} \right]\times k}  , \quad for\;  m\geq 2,\;  0 \leq j \leq m,\; k\geq s+m \; and \; j=m+1,\;k=s+m+1
 $}}
  \label{subsec 4}
 \begin{lem}
 \label{lem 12.3}
    We have  for $m\geq 2$ :
  \begin{equation}
 \label{eq 12.7}
   \Gamma_{s+j}^{\left[s\atop{ s +m\atop s+m} \right]\times k} 
 =  \begin{cases}
2^{6s+2m-3}-7\cdot2^{4s-6}+3\cdot2^{3s-5}  &\text{if  }  j = 0 ,\;k=s\\
2^{k+s-1} -2^{2s}+21\cdot\big(5\cdot2^{4s-6}-2^{3s-5}\big)  &\text{if  }  j = 0 ,\;k\geq s+1\\
2^{6s+2m}-7\cdot2^{4s-2}+3\cdot2^{3s-2}  &\text{if  }  j = 1 ,\;k=s+1\\
9\cdot2^{k+s-1} +21\cdot\big(5\cdot2^{4s-2}-2^{3s-2}-2^{2s}\big)  &\text{if  }  j = 1 ,\;k\geq s+2\\
2^{6s+2m+3}-7\cdot2^{4s+2}+3\cdot2^{3s+1} +2^{2s+1} &\text{if  }  j = 2 ,\;k=s+2\\
39\cdot2^{k+s} +21\cdot\big(5\cdot2^{4s+2}-2^{3s+1}-9\cdot2^{2s+1}\big)  &\text{if  }  j = 2 ,\;k\geq s+3\\
2^{6s+2m+6}-7\cdot2^{4s+6}+3\cdot2^{3s+4} +3\cdot2^{2s+4} &\text{if  }  j = 3 ,\;k=s+3\\
81\cdot2^{k+s+2} +21\cdot\big(5\cdot2^{4s+6}-2^{3s+4}-19\cdot2^{2s+4}\big)  &\text{if  }  j = 3 ,\;k\geq s+4\\
2^{6s+2m+9}-7\cdot2^{4s+10}+3\cdot2^{3s+7} +7\cdot2^{2s+7} &\text{if  }  j = 4 ,\;k=s+4\\
165\cdot2^{k+s+4} +21\cdot\big(5\cdot2^{4s+10}-2^{3s+7}-39\cdot2^{2s+7}\big)  &\text{if  }  j = 4 ,\;k\geq s+5\\
2^{6s+2m+12}-7\cdot2^{4s+14}+3\cdot2^{3s+10} +15\cdot2^{2s+10} &\text{if  }  j = 5 ,\;k=s+5\\
(21\cdot2^{j-1}-3)\cdot2^{k+s+2j-4}\\
 +21\cdot\big(5\cdot2^{4s+4j-6}-2^{3s+3j-5}-5\cdot2^{2s+4j-6}+2^{2s+3j-5}\big) 
  &\text{if  }  1\leq j\leq m-1 ,\;k\geq s+j+1\\
 2^{6s+2m+3j-3}-7\cdot2^{4s+4j-6}+3\cdot2^{3s+3j-5}\\
  +(2^{j-1}-1)\cdot2^{2s+3j-5} &\text{if  } 1\leq j\leq m+1 ,\;k=s+j  \\
 \big(21\cdot2^{s+3m-5}+45\cdot2^{s+2m-4}\big)\cdot\big(2^{k}-2^{s+m+1}\big)\\
 + 105\cdot2^{4s+4m-6}-21\cdot2^{3s+3m-5}-21\cdot2^{2s+4m-6}+9\cdot2^{2s+3m-5} 
   &\text{if  }  j=m ,\;k\geq s+m+1
   \end{cases}
\end{equation}
 Let  $ (H_{q-1}) $ denote the following statement :\\
\begin{align*}
& \Gamma_{s+j}^{\left[s\atop{ s +m\atop s+m} \right]\times k}
= (21\cdot2^{j-1}-3)\cdot2^{k+s+2j-4}
  +21\cdot\big(5\cdot2^{4s+4j-6}-2^{3s+3j-5}-5\cdot2^{2s+4j-6}+2^{2s+3j-5}\big) \\
   &\text{for  }  1\leq j\leq q-1\leq  m-2 ,\;k\geq s+j+1
\end{align*}
 \end{lem}
\begin{proof}
Proof by strong induction \\[0.1 cm]
We start by prove that $(H_{q-1})$ holds for $q=2,3,4,5$\\[0.1 cm]
Next we assume that  $(H_{q-1})$ holds. then we prove that this can be used to show that  $(H_{q})$ holds, that is :\\[0.1 cm]
\begin{align*}
& \Gamma_{s+q}^{\left[s\atop{ s +m\atop s+m} \right]\times k}
= (21\cdot2^{q-1}-3)\cdot2^{k+s+2q-4}
  +21\cdot\big(5\cdot2^{4s+4q-6}-2^{3s+3q-5}-5\cdot2^{2s+4q-6}+2^{2s+3q-5}\big) \\
   &\text{for  }  2\leq q\leq m-1,\;k\geq s+q+1
\end{align*}
\textbf{The case $j=0,\;k=s$}\\[0.1 cm]
Follows from \eqref{eq 6.41} with $l=0,\;i=s$\\[0.1 cm]
\textbf{The case $j=0,\;k\geq s+1$}\\[0.1 cm]
Follows from \eqref{eq 12.1} with $j=0$ using \eqref{eq 6.42} with $l=0$\\[0.1 cm]
\textbf{The case $j=1,\;k= s+1$}\\[0.1 cm]
Follows from \eqref{eq 6.43} with $l=0$\\[0.1 cm]
\textbf{The case $j=1,\;k\geq  s+2$ and the case $j=2,\;k=s+2$}\\[0.1 cm]
From \eqref{eq 6.1} with $k=s+2,\;l=0$ we obtain \\[0.1 cm]
\begin{align*}
& \sum_{i=0}^{s+2}\Gamma_{i}^{\left[s\atop{ s +m\atop s+m} \right]\times (s+2)}= 2^{6s+2m+3}\\
&\Updownarrow \\
&  \sum_{i=0}^{s-1}\Gamma_{i}^{\left[s\atop{ s +m\atop s+m} \right]\times (s+2)} + \Gamma_{s}^{\left[s\atop{ s +m\atop s+m} \right]\times (s+2)}
+ \Gamma_{s+1}^{\left[s\atop{ s +m\atop s+m} \right]\times (s+2)} + \Gamma_{s+2}^{\left[s\atop{ s +m\atop s+m} \right]\times (s+2)}= 2^{6s+2m+3}\\
& \Updownarrow \\
& \text{using \eqref{eq 6.39} and \eqref{eq 12.7} in the case $j=0,\;k\geq s+1$ with $k=s+2$ we obtain}\\[0.1 cm]
& \big(7\cdot2^{4s-6}-3\cdot2^{3s-5}\big) + \big(2^{2s}+ 21\cdot(5\cdot2^{4s-6}-2^{3s-5})\big)+\Gamma_{s+1}^{\left[s\atop{ s +m\atop s+m} \right]\times (s+2)} + \Gamma_{s+2}^{\left[s\atop{ s +m\atop s+m} \right]\times (s+2)}= 2^{6s+2m+3}\\
& \Updownarrow \\
& \Gamma_{s+1}^{\left[s\atop{ s +m\atop s+m} \right]\times (s+2)} + \Gamma_{s+2}^{\left[s\atop{ s +m\atop s+m} \right]\times (s+2)}= 2^{6s+2m+3}
-7\cdot2^{4s-2}+3\cdot2^{3s-2}-2^{2s}
\end{align*}
From \eqref{eq 6.2} with $k=s+2,\;l=0$ we get \\[0.1 cm]
\begin{align*}
& \sum_{i=0}^{s+2}\Gamma_{i}^{\left[s\atop{ s +m\atop s+m} \right]\times (s+2)}\cdot2^{-i}= 2^{5s+2m+1}+2^{3s+3}-2^{2s+1}\\
&\Updownarrow \\
&  \sum_{i=0}^{s-1}\Gamma_{i}^{\left[s\atop{ s +m\atop s+m} \right]\times (s+2)}\cdot2^{-i} + \Gamma_{s}^{\left[s\atop{ s +m\atop s+m} \right]\times (s+2)}\cdot2^{-s}
+ \Gamma_{s+1}^{\left[s\atop{ s +m\atop s+m} \right]\times (s+2)}\cdot2^{-(s+1)} + \Gamma_{s+2}^{\left[s\atop{ s +m\atop s+m} \right]\times (s+2)}\cdot2^{-(s+2)}\\
& = 2^{5s+2m+1}+2^{3s+3}-2^{2s+1} \\
& \Updownarrow  \text{using \eqref{eq 6.40} and \eqref{eq 12.7} in the case $j=0,\;k\geq s+1$ with $k=s+2$}\\[0.2 cm]
& \big(15\cdot2^{3s-6}-7\cdot2^{2s-5}\big) + \big(2^{s}+ 21\cdot(5\cdot2^{3s-6}-2^{2s-5})\big)\\
& +\Gamma_{s+1}^{\left[s\atop{ s +m\atop s+m} \right]\times (s+2)}\cdot2^{-(s+1)} + \Gamma_{s+2}^{\left[s\atop{ s +m\atop s+m} \right]\times (s+2)}\cdot2^{-(s+2)}= 2^{5s+2m+1}+2^{3s+3}-2^{2s+1}\\
& \Updownarrow \\
& 2\cdot \Gamma_{s+1}^{\left[s\atop{ s +m\atop s+m} \right]\times (s+2)} + \Gamma_{s+2}^{\left[s\atop{ s +m\atop s+m} \right]\times (s+2)}
= 2^{6s+2m+3} +49\cdot2^{4s-1}-9\cdot2^{3s-1}-2^{2s+2} \\
& \text{hence we have (see above)}\\
& \Gamma_{s+1}^{\left[s\atop{ s +m\atop s+m} \right]\times (s+2)} = 105\cdot2^{4s-2}-21\cdot2^{3s-2}-3\cdot2^{2s}\\[0.1 cm]
& \Gamma_{s+2}^{\left[s\atop{ s +m\atop s+m} \right]\times (s+2)} =  2^{6s+2m+3} -7\cdot2^{4s+2}+3\cdot2^{3s+1}+2^{2s+1}
\end{align*}
From \eqref{eq 12.1} with $j=1$ we get \\[0.1 cm]
  \begin{align*}
& \Gamma_{s+1}^{\left[s\atop{ s +m\atop s+m} \right]\times k} -\Gamma_{s+1}^{\left[s\atop{ s +m\atop s+m} \right]\times (s+2)} 
= \big(21\cdot2^{s-2}-3\cdot2^{s-2}\big)\cdot\big(2^{k}-2^{s+2}\big)\\
& \Updownarrow \\
& \Gamma_{s+1}^{\left[s\atop{ s +m\atop s+m} \right]\times k}=  105\cdot2^{4s-2}-21\cdot2^{3s-2}-3\cdot2^{2s} +\big(21\cdot2^{s-2}-3\cdot2^{s-2}\big)\cdot\big(2^{k}-2^{s+2}\big)\\
& \Updownarrow \\
& \Gamma_{s+1}^{\left[s\atop{ s +m\atop s+m} \right]\times k} = 9\cdot2^{k+s-1} +21\cdot\big(5\cdot2^{4s-2}-2^{3s-2}-2^{2s}\big)  
\end{align*}
\textbf{The case $j=2,\;k\geq  s+3$ and the case $j=3,\;k=s+3$}\\[0.1 cm]
From \eqref{eq 6.1} with $k=s+3,\;l=0$ we obtain \\[0.1 cm]
\begin{align*}
& \sum_{i=0}^{s+3}\Gamma_{i}^{\left[s\atop{ s +m\atop s+m} \right]\times (s+3)}= 2^{6s+2m+6}\\
&\Updownarrow \\
&  \sum_{i=0}^{s-1}\Gamma_{i}^{\left[s\atop{ s +m\atop s+m} \right]\times (s+3)} + \Gamma_{s}^{\left[s\atop{ s +m\atop s+m} \right]\times (s+3)}
+ \Gamma_{s+1}^{\left[s\atop{ s +m\atop s+m} \right]\times (s+3)} + \Gamma_{s+2}^{\left[s\atop{ s +m\atop s+m} \right]\times (s+3)}
 + \Gamma_{s+3}^{\left[s\atop{ s +m\atop s+m} \right]\times (s+3)}= 2^{6s+2m+6}
 \end{align*}
 Using \eqref{eq 6.39} and \eqref{eq 12.7} in the cases $j=0,\;k=s+3$ and $j=1,\;k=s+3$ we have\\[0.1 cm]
 \begin{align*}
& \big(7\cdot2^{4s-6}-3\cdot2^{3s-5}\big) + \big(3\cdot2^{2s}+ 21\cdot(5\cdot2^{4s-6}-2^{3s-5})\big)+\big(15\cdot2^{2s}+ 21\cdot(5\cdot2^{4s-2}-2^{3s-2})\big) \\
& + \Gamma_{s+2}^{\left[s\atop{ s +m\atop s+m} \right]\times (s+3)} + \Gamma_{s+3}^{\left[s\atop{ s +m\atop s+m} \right]\times (s+3)}= 2^{6s+2m+6}\\
& \Updownarrow \\
& \Gamma_{s+2}^{\left[s\atop{ s +m\atop s+m} \right]\times (s+3)} + \Gamma_{s+3}^{\left[s\atop{ s +m\atop s+m} \right]\times (s+3)}
= 2^{6s+2m+6}-7\cdot2^{4s+2}+3\cdot2^{3s+1}-9\cdot2^{2s+1}
\end{align*}
From \eqref{eq 6.2} with $k=s+3,\;l=0$ we obtain \\[0.1 cm]
\begin{align*}
& \sum_{i=0}^{s+3}\Gamma_{i}^{\left[s\atop{ s +m\atop s+m} \right]\times (s+2)}\cdot2^{-i}= 2^{5s+2m+3}+2^{3s+6}-2^{2s+3}\\
&\Updownarrow \\
&  \sum_{i=0}^{s-1}\Gamma_{i}^{\left[s\atop{ s +m\atop s+m} \right]\times (s+3)}\cdot2^{-i} + \Gamma_{s}^{\left[s\atop{ s +m\atop s+m} \right]\times (s+3)}\cdot2^{-s}
+ \Gamma_{s+1}^{\left[s\atop{ s +m\atop s+m} \right]\times (s+3)}\cdot2^{-(s+1)} + \Gamma_{s+2}^{\left[s\atop{ s +m\atop s+m} \right]\times (s+3)}\cdot2^{-(s+2)}\\
& + \Gamma_{s+3}^{\left[s\atop{ s +m\atop s+m} \right]\times (s+3)}\cdot2^{-(s+3)}  = 2^{5s+2m+3}+2^{3s+6}-2^{2s+3} 
\end{align*}
 Using \eqref{eq 6.40} and \eqref{eq 12.7} in the cases $j=0,\;k=s+3$ and $j=1,\;k=s+3$ we have\\[0.1 cm]
\begin{align*} 
& \big(15\cdot2^{3s-6}-7\cdot2^{2s-5}\big) + \big(3\cdot2^{s}+ 21\cdot(5\cdot2^{3s-6}-2^{2s-5})\big) + \big(15\cdot2^{s-1}+ 21\cdot(5\cdot2^{3s-3}-2^{2s-3})\big) \\
& +\Gamma_{s+2}^{\left[s\atop{ s +m\atop s+m} \right]\times (s+3)}\cdot2^{-(s+2)} + \Gamma_{s+3}^{\left[s\atop{ s +m\atop s+m} \right]\times (s+3)}\cdot2^{-(s+3)}= 2^{5s+2m+3}+2^{3s+6}-2^{2s+3}\\
& \Updownarrow \\
& 2\cdot \Gamma_{s+2}^{\left[s\atop{ s +m\atop s+m} \right]\times (s+3)} + \Gamma_{s+3}^{\left[s\atop{ s +m\atop s+m} \right]\times (s+3)}
= 2^{6s+2m+6} +49\cdot2^{4s+3}-9\cdot2^{3s+2}-21\cdot2^{2s+2}, \\
& \text{Hence we have (see above)}\\
& \Gamma_{s+2}^{\left[s\atop{ s +m\atop s+m} \right]\times (s+3)} = 105\cdot2^{4s+2}-21\cdot2^{3s+1}-33\cdot2^{2s+1}\\[0.1 cm]
& \Gamma_{s+3}^{\left[s\atop{ s +m\atop s+m} \right]\times (s+3)} =  2^{6s+2m+6} -7\cdot2^{4s+6}+3\cdot2^{3s+4}+3\cdot2^{2s+4}
\end{align*}
From \eqref{eq 12.1} with $j=2$ we get \\[0.1 cm]
  \begin{align*}
& \Gamma_{s+2}^{\left[s\atop{ s +m\atop s+m} \right]\times k} -\Gamma_{s+2}^{\left[s\atop{ s +m\atop s+m} \right]\times (s+3)} 
= \big(21\cdot2^{s+1}-3\cdot2^{s}\big)\cdot\big(2^{k}-2^{s+3}\big)\\
& \Updownarrow \\
& \Gamma_{s+2}^{\left[s\atop{ s +m\atop s+m} \right]\times k}=  105\cdot2^{4s+2}-21\cdot2^{3s+1}-33\cdot2^{2s+1} +39\cdot2^{s}\cdot\big(2^{k}-2^{s+3}\big)\\
& \Updownarrow \\
& \Gamma_{s+2}^{\left[s\atop{ s +m\atop s+m} \right]\times k} = 39\cdot2^{k+s} +21\cdot\big(5\cdot2^{4s+2}-2^{3s+1}-9\cdot2^{2s+1}\big)  
\end{align*}
From the above results we now deduce that $(H_{q-1})$ holds for $q=2,3,4,5.$\\[0.1 cm]
Besides we get the formulas for $\Gamma_{s+j}^{\left[s\atop{ s +m\atop s+m} \right]\times (s+j)}$ for $j=0,1,2,3,4,5$\\[0.1 cm]
\textbf{Assume that $(H_{q-1})$ holds.}\\[0.1 cm]
From \eqref{eq 6.1} with $k=s+q+1,\;l=0$ we obtain \\[0.1 cm]
\begin{align*}
& \sum_{i=0}^{s+q+1}\Gamma_{i}^{\left[s\atop{ s +m\atop s+m} \right]\times (s+q+1)}= 2^{6s+2m+3q}\\
&\Updownarrow \\
&  \sum_{i=0}^{s-1}\Gamma_{i}^{\left[s\atop{ s +m\atop s+m} \right]\times (s+q+1)} + \Gamma_{s}^{\left[s\atop{ s +m\atop s+m} \right]\times (s+q+1)}
+\sum_{j=1}^{q-1} \Gamma_{s+j}^{\left[s\atop{ s +m\atop s+m} \right]\times (s+q+1)}\\
& + \Gamma_{s+q}^{\left[s\atop{ s +m\atop s+m} \right]\times (s+q+1)}
 + \Gamma_{s+q+1}^{\left[s\atop{ s +m\atop s+m} \right]\times (s+q+1)}= 2^{6s+2m+3q}\\
 & \Updownarrow \\
  & \big(7\cdot2^{4s-6}-3\cdot2^{3s-5}\big) + \big(2^{2s+q}-2^{2s}+ 21\cdot(5\cdot2^{4s-6}-2^{3s-5})\big) \\
   & +\sum_{j=1}^{q-1}(21\cdot2^{j-1}-3)\cdot2^{2s+q+2j-3}
  +21\cdot\big(5\cdot2^{4s+4j-6}-2^{3s+3j-5}-5\cdot2^{2s+4j-6}+2^{2s+3j-5}\big) \\
  & + \Gamma_{s+q}^{\left[s\atop{ s +m\atop s+m} \right]\times (s+q+1)} + \Gamma_{s+q+1}^{\left[s\atop{ s +m\atop s+m} \right]\times (s+q+1)}= 2^{6s+2m+3q}\\
& \Updownarrow \\
& \Gamma_{s+q}^{\left[s\atop{ s +m\atop s+m} \right]\times (s+q+1)} + \Gamma_{s+q+1}^{\left[s\atop{ s +m\atop s+m} \right]\times (s+q+1)}\\[0.1 cm]
& = 2^{6s+2m+3q}-7\cdot2^{4s+4q-6}+3\cdot2^{3s+3q-5}-5\cdot2^{2s+4q-6}+2^{2s+3q-5}
\end{align*}
Using \eqref{eq 6.2} with $k=s+q+1,\;l=0$ we obtain \\[0.1 cm]
\begin{align*}
& \sum_{i=0}^{s+q+1}\Gamma_{i}^{\left[s\atop{ s +m\atop s+m} \right]\times (s+q+1)}\cdot2^{-i}= 2^{5s+2m+2q-1}+2^{3s+3q}-2^{2s+2q-1}\\
&\Updownarrow \\
&  \sum_{i=0}^{s-1}\Gamma_{i}^{\left[s\atop{ s +m\atop s+m} \right]\times (s+q+1)}\cdot2^{-i} + \Gamma_{s}^{\left[s\atop{ s +m\atop s+m} \right]\times (s+q+1)}\cdot2^{-s}
+\sum_{j=1}^{q-1} \Gamma_{s+j}^{\left[s\atop{ s +m\atop s+m} \right]\times (s+q+1)}\cdot2^{-(s+j)}\\
& + \Gamma_{s+q}^{\left[s\atop{ s +m\atop s+m} \right]\times (s+q+1)}\cdot2^{-(s+q)}
 + \Gamma_{s+q+1}^{\left[s\atop{ s +m\atop s+m} \right]\times (s+q+1)}\cdot2^{-(s+q+1)}= 2^{5s+2m+2q-1}+2^{3s+3q}-2^{2s+2q-1} \\
 & \Updownarrow \\
  & \big(15\cdot2^{3s-6}-7\cdot2^{2s-5}\big) + \big(2^{s+q}-2^{s}+ 21\cdot(5\cdot2^{3s-6}-2^{2s-5})\big) \\
   & +\sum_{j=1}^{q-1}(21\cdot2^{j-1}-3)\cdot2^{s+q+j-3}
  +21\cdot\big(5\cdot2^{3s+3j-6}-2^{2s+2j-5}-5\cdot2^{s+3j-6}+2^{s+2j-5}\big) \\
  & + \Gamma_{s+q}^{\left[s\atop{ s +m\atop s+m} \right]\times (s+q+1)}\cdot2^{-(s+q)} + \Gamma_{s+q+1}^{\left[s\atop{ s +m\atop s+m} \right]\times (s+q+1)}\cdot2^{-(s+q+1)}= 2^{5s+2m+2q-1}+2^{3s+3q}-2^{2s+2q-1}\\
& \Updownarrow \\
& 2\cdot\Gamma_{s+q}^{\left[s\atop{ s +m\atop s+m} \right]\times (s+q+1)} + \Gamma_{s+q+1}^{\left[s\atop{ s +m\atop s+m} \right]\times (s+q+1)}\\[0.1 cm]
& = 2^{6s+2m+3q}+49\cdot2^{4s+4q-5}-9\cdot2^{3s+3q-4}+5\cdot2^{2s+3q-4}-13\cdot2^{2s+4q-5}\\
& \text{ Hence we have (see above)}\\
& \Gamma_{s+q}^{\left[s\atop{ s +m\atop s+m} \right]\times (s+q+1)} = 105\cdot2^{4s+4q-6}-21\cdot2^{3s+3q-5}-21\cdot2^{2s+4q-6}+9\cdot2^{2s+3q-5}\\[0.1 cm]
& \Gamma_{s+q+1}^{\left[s\atop{ s +m\atop s+m} \right]\times (s+q+1)} =  2^{6s+2m+3q} -7\cdot2^{4s+4q-2}+3\cdot2^{3s+3q-2}+2^{2s+4q-2}-2^{2s+3q-2}
\end{align*}

From \eqref{eq 12.1} with $j=q$ we get \\[0.1 cm]
  \begin{align*}
& \Gamma_{s+q}^{\left[s\atop{ s +m\atop s+m} \right]\times k} -\Gamma_{s+q}^{\left[s\atop{ s +m\atop s+m} \right]\times (s+q+1)} 
= \big(21\cdot2^{s+3q-5}-3\cdot2^{s+2q-4}\big)\cdot\big(2^{k}-2^{s+q+1}\big)\\
& \Updownarrow \\
& \Gamma_{s+q}^{\left[s\atop{ s +m\atop s+m} \right]\times k}= 
  105\cdot2^{4s+4q-6}-21\cdot2^{3s+3q-5}-21\cdot2^{2s+4q-6}+9\cdot2^{2s+3q-5}\\
  & +\big(21\cdot2^{s+3q-5}-3\cdot2^{s+2q-4}\big)\cdot\big(2^{k}-2^{s+q+1}\big)\\
& \Updownarrow \\
& \Gamma_{s+q}^{\left[s\atop{ s +m\atop s+m} \right]\times k} = (21\cdot2^{q-1}-3)\cdot2^{k+s+2q-4}
 +21\cdot\big(5\cdot2^{4s+4q-6}-2^{3s+3q-5}-5\cdot2^{2s+4q-6}+2^{2s+3q-5}\big)  
\end{align*}
\textbf{The case $j=m,\;k\geq s+m+1$}\\[0.1 cm]
We obtain from the proof by induction with $q=m$ (see above)\\[0.1 cm]
\begin{align*}
& \Gamma_{s+m}^{\left[s\atop{ s +m\atop s+m} \right]\times (s+m+1)} = 105\cdot2^{4s+4m-6}-21\cdot2^{3s+3m-5}-21\cdot2^{2s+4m-6}+9\cdot2^{2s+3m-5}\\[0.1 cm]
& \Gamma_{s+m+1}^{\left[s\atop{ s +m\atop s+m} \right]\times (s+m+1)} =  2^{6s+5m} -7\cdot2^{4s+4m-2}+3\cdot2^{3s+3m-2}+2^{2s+4m-2}-2^{2s+3m-2}
\end{align*}
From \eqref{eq 12.1} with $j=m$ we get \\[0.1 cm]
  \begin{align*}
& \Gamma_{s+m}^{\left[s\atop{ s +m\atop s+m} \right]\times k} -\Gamma_{s+m}^{\left[s\atop{ s +m\atop s+m} \right]\times (s+m+1)} 
= \big(21\cdot2^{s+3m-5}+45\cdot2^{s+2m-4}\big)\cdot\big(2^{k}-2^{s+m+1}\big)\\
& \Updownarrow \\
& \Gamma_{s+m}^{\left[s\atop{ s +m\atop s+m} \right]\times k}= 
  105\cdot2^{4s+4m-6}-21\cdot2^{3s+3m-5}-21\cdot2^{2s+4m-6}+9\cdot2^{2s+3m-5}\\
  & +\big(21\cdot2^{s+3m-5}+45\cdot2^{s+2m-4}\big)\cdot\big(2^{k}-2^{s+m+1}\big)\\
& \Updownarrow \\
& \Gamma_{s+m}^{\left[s\atop{ s +m\atop s+m} \right]\times k} = \big(21\cdot2^{s+3m-5}+45\cdot2^{s+2m-4}\big)\cdot\big(2^{k}-2^{s+m+1}\big)\\
& + 105\cdot2^{4s+4m-6}-21\cdot2^{3s+3m-5}-21\cdot2^{2s+4m-6}+9\cdot2^{2s+3m-5}
\end{align*}
\end{proof}
 
  \subsection{\textbf{Computation of $ \Gamma_{s+m+j}^{\left[s\atop{ s +m\atop s+m} \right]\times k} 
   , \quad for\; 1 \leq j \leq s+2,\;k\geq s+m+j,\; m\geq 2,\;j=s+3,\;k=2s+m+3. $}}
  \label{subsec 5}
 \begin{lem}
 \label{lem 12.4}
   We have  for $m\geq 2,\;k\geq s+m $ :\\[0.1 cm]
  \begin{align}
& \sum_{j=1}^{m-1}\Gamma_{s+j}^{\left[s\atop{ s +m\atop s+m} \right]\times k} \label{eq 12.8} \\
& =3\cdot2^{k+s+3m-5}- 3\cdot2^{3s+3m-5}+3\cdot2^{2s+3m-5} -2^{k+s-1} \nonumber \\
& +3\cdot2^{3s-2}+2^{2s}+7\cdot2^{4s+4m-6}-7\cdot2^{2s+4m-6}-7\cdot2^{4s-2}-2^{k+s+2m-4}\nonumber\\
& \nonumber
 \end{align}
  \begin{align}
& \sum_{j=1}^{m-1}\Gamma_{s+j}^{\left[s\atop{ s +m\atop s+m} \right]\times k}\cdot2^{-(s+j} \label{eq 12.9} \\
& =7\cdot2^{k+2m-5}- 7\cdot2^{2s+2m-5}+7\cdot2^{s+2m-5} -2^{k-1} \nonumber \\
& +7\cdot2^{2s-3}+2^{s}+15\cdot2^{3s+3m-6}-15\cdot2^{s+3m-6}-15\cdot2^{3s-3}-3\cdot2^{k+m-4}\nonumber\\
& \nonumber
 \end{align}
 \end{lem}
 \begin{proof}
 \textbf{Proof of \eqref{eq 12.8}} \\[0.1 cm]
 From Lemma \ref{lem 12.3} in the case $1\leq  j \leq  m-1,\;k \geq s+j+1$ we obtain :\\[0.1 cm]
 \begin{align*}
 & \sum_{j=1}^{m-1}\Gamma_{s+j}^{\left[s\atop{ s +m\atop s+m} \right]\times k}\\
 & = \big(21\cdot2^{k+s-5}-21\cdot2^{3s-5}+21\cdot2^{2s-5} \big)\cdot\sum_{j=1}^{m-1}2^{3j}\\
 & + \big(105\cdot2^{4s-6}-105\cdot2^{2s-6} \big)\cdot\sum_{j=1}^{m-1}2^{4j}
  -\big(3\cdot2^{k+s-4} \big)\cdot\sum_{j=1}^{m-1}2^{2j}\\
  & =  \big(3\cdot2^{k+s-5}-3\cdot2^{3s-5}+3\cdot2^{2s-5} \big)\cdot(2^{3m} -2^3)\\
 & + \big(7\cdot2^{4s-6}-7\cdot2^{2s-6} \big)\cdot(2^{4m} -2^4)
  -\big(2^{k+s-4} \big)\cdot(2^{2m} -2^2)
   \end{align*}
 \textbf{Proof of \eqref{eq 12.9}} \\[0.1 cm]
Similar to the proof of  \eqref{eq 12.8} : \\[0.1 cm]
 \begin{align*}
 & \sum_{j=1}^{m-1}\Gamma_{s+j}^{\left[s\atop{ s +m\atop s+m} \right]\times k}\cdot2^{-(s+j)}\\
 & = \big(21\cdot2^{k-5}-21\cdot2^{2s-5}+21\cdot2^{s-5} \big)\cdot\sum_{j=1}^{m-1}2^{2j}\\
 & + \big(105\cdot2^{3s-6}-105\cdot2^{s-6} \big)\cdot\sum_{j=1}^{m-1}2^{3j}
  -\big(3\cdot2^{k-4} \big)\cdot\sum_{j=1}^{m-1}2^{j}\\
  & =  \big(7\cdot2^{k-5}-7\cdot2^{2s-5}+7\cdot2^{s-5} \big)\cdot(2^{2m} -2^2)\\
 & + \big(15\cdot2^{3s-6}-15\cdot2^{s-6} \big)\cdot(2^{3m} -2^3)
  -\big(3\cdot2^{k-4} \big)\cdot(2^{m} -2)
   \end{align*}
 \end{proof}
 
  \begin{lem}
 \label{lem 12.5}
   We have  for $m\geq 2,\;1\leq j\leq s-1$ :
  \begin{align}
& \Gamma_{s+m+j}^{\left[s\atop{ s +m\atop s+m} \right]\times k} \nonumber \\
& =  \begin{cases} \label{eq 12.10}
 105\cdot2^{4s+4m-2}-21\cdot2^{3s+3m-2}-21\cdot2^{2s+4m-2}-63\cdot2^{2s+3m-2}  &\text{if  }  j = 1 ,\;k= s+m+2\\
  2^{6s+5m+3}-7\cdot2^{4s+4m+2}+3\cdot2^{3s+3m+1}+2^{2s+4m+2}+5\cdot2^{2s+3m+1} &\text{if  } j=2 ,\;k=s+m+2\\
 21\cdot\big[ 2^{k+s+3m-2}+13\cdot2^{k+s+2m-2}+ 5\cdot2^{4s+4m-2}\\
  -2^{3s+3m-2}- 5\cdot2^{2s+4m-2}-55\cdot2^{2s+3m-2}\big ] &\text{if  } j=1 ,\;k\geq s+m+2\\
  \\
   105\cdot2^{4s+4m+2}-21\cdot2^{3s+3m+1}-21\cdot2^{2s+4m+2}-441\cdot2^{2s+3m+1}  &\text{if  }  j = 2 ,\;k= s+m+3\\
   2^{6s+5m+6}-7\cdot2^{4s+4m+6}+3\cdot2^{3s+3m+4}+2^{2s+4m+6}+19\cdot2^{2s+3m+5} &\text{if  } j=3 ,\;k=s+m+3\\
   21\cdot\big[ 2^{k+s+3m+1}+61\cdot2^{k+s+2m}+ 5\cdot2^{4s+4m+2}\\
  -2^{3s+3m+1}- 5\cdot2^{2s+4m+2}-265\cdot2^{2s+3m+1}\big ] &\text{if  } j=2 ,\;k\geq s+m+3\\
  \\
   105\cdot2^{4s+4m+6}-21\cdot2^{3s+3m+4}-21\cdot2^{2s+4m+6}-1071\cdot2^{2s+3m+5}  &\text{if  }  j = 3 ,\;k= s+m+4\\
 2^{6s+5m+9}-7\cdot2^{4s+4m+10}+3\cdot2^{3s+3m+7}+2^{2s+4m+10}+47\cdot2^{2s+3m+9} &\text{if  } j=4 ,\;k=s+m+4\\
  21\cdot\big[ 2^{k+s+3m+4}+131\cdot2^{k+s+2m+3}+ 5\cdot2^{4s+4m+6}\\
  -2^{3s+3m+4}- 5\cdot2^{2s+4m+6}-575\cdot2^{2s+3m+5}\big ] &\text{if  } j=3 ,\;k\geq s+m+4\\
  \\
   105\cdot2^{4s+4m+10}-21\cdot2^{3s+3m+7}-21\cdot2^{2s+4m+10}-2331\cdot2^{2s+3m+9}  &\text{if  }  j = 4 ,\;k= s+m+5\\
 2^{6s+5m+12}-7\cdot2^{4s+4m+14}+3\cdot2^{3s+3m+10}+2^{2s+4m+14}+103\cdot2^{2s+3m+13} &\text{if  } j=5 ,\;k=s+m+5\\
  21\cdot\big[ 2^{k+s+3m+7}+271\cdot2^{k+s+2m+6}+ 5\cdot2^{4s+4m+10}\\
    -2^{3s+3m+7}- 5\cdot2^{2s+4m+10}-1195\cdot2^{2s+3m+9}\big ] &\text{if  } j=4 ,\;k\geq s+m+5\\
    \\
  105\cdot2^{4s+4m+4j-6}-21\cdot2^{3s+3m+3j-5}-21\cdot2^{2s+4m+4j-6}\\
  -315\cdot2^{2s+3m+5j-8}+189\cdot2^{2s+3m+4j-7}  &\text{if  }  k= s+m+j+1\\
 2^{6s+5m+3j-3}-7\cdot2^{4s+4m+4j-6}+3\cdot2^{3s+3m+3j-5}\\
+ 2^{2s+4m+4j-6}+7\cdot2^{2s+3m+5j-8}- 9\cdot2^{2s+3m+4j-7}&\text{if  } k=s+m+j\\
  21\cdot\big[ 2^{k+s+3m+3j-5}+35\cdot2^{k+s+2m+4j-7}-9\cdot2^{k+s+2m+3j-6}  + 5\cdot2^{4s+4m+4j-6}\\
  -2^{3s+3m+3j-5}- 5\cdot2^{2s+4m+4j-6}-155\cdot2^{2s+3m+5j-8}+45\cdot2^{2s+3m+4j-7}\big ] &\text{if  }k\geq s+m+j+1
   \end{cases} 
\end{align}
 Let  $ (H_{q-1}) $ denote the following statement :
\begin{align*}
& \Gamma_{s+m+j}^{\left[s\atop{ s +m\atop s+m} \right]\times k}
=   21\cdot\big[ 2^{k+s+3m+3j-5}+35\cdot2^{k+s+2m+4j-7}-9\cdot2^{k+s+2m+3j-6}  + 5\cdot2^{4s+4m+4j-6}\\
&  -2^{3s+3m+3j-5}- 5\cdot2^{2s+4m+4j-6}-155\cdot2^{2s+3m+5j-8}+45\cdot2^{2s+3m+4j-7}\big ] \\
&   \text{for  }  1\leq j\leq q-1\leq  s-2 ,\;k\geq s+m+j+1
\end{align*}
 \end{lem} 
 \begin{proof}
 Similar to the proof of Lemma \ref{lem 12.3}.\\[0.1 cm]
  We start by computing  successively 
$ \Gamma_{s+m+j}^{\left[s\atop{ s +m\atop s+m} \right]\times k} $ for
 $j=1,2,3,4\; and \;k\geq s+m+j+1.$\\[0.1 cm]
 The statement $ (H_{q-1}) $ is now obtained by observing their  structures. 

 \end{proof}
 
  \begin{lem}
   \label{lem 12.6}
   We have  for $m\geq 2,\;k\geq 2s+m $ :\\[0.2 cm]
  \begin{align}
& \sum_{j=1}^{s-1}\Gamma_{s+m+j}^{\left[s\atop{ s +m\atop s+m} \right]\times k} \label{eq 12.11} \\
& =3\cdot2^{k+4s+3m-5}- 27\cdot2^{k+4s+2m-6}-3\cdot2^{k+s+3m-2} -11\cdot2^{k+s+2m-2} \nonumber \\
& +49\cdot2^{k+5s+2m-7}+51\cdot2^{6s+3m-7}+3\cdot2^{3s+3m-2}+7\cdot2^{8s+4m-6}-7\cdot2^{6s+4m-6}\nonumber\\
& -7\cdot2^{4s+4m-2} + 7\cdot2^{2s+4m-2} -105\cdot2^{7s+3m-8}+21\cdot2^{2s+3m-2} \nonumber
 \end{align}
Setting  $k=2s+m$ in \eqref{eq 12.11} we get \\[0.1 cm]
 \begin{align}
& \sum_{j=1}^{s-1}\Gamma_{s+m+j}^{\left[s\atop{ s +m\atop s+m} \right]\times (2s+m)} \label{eq 12.12} \\
& =7\cdot2^{8s+4m-6}+7\cdot2^{2s+4m-2}+21\cdot2^{2s+3m-2} -7\cdot2^{4s+4m-2} \nonumber \\
& -7\cdot2^{7s+3m-8}-3\cdot2^{3s+4m-2}-3\cdot2^{6s+3m-7}-2^{6s+4m-6}-2^{3s+3m+1}\nonumber
 \end{align}
  \begin{align}
& \sum_{j=1}^{s-1}\Gamma_{s+m+j}^{\left[s\atop{ s +m\atop s+m} \right]\times k}\cdot2^{-(s+m+j)} \label{eq 12.13} \\
& =7\cdot2^{k+2s+2m-5}- 63\cdot2^{k+2s+m-6}-7\cdot2^{k+2m-3} -21\cdot2^{k+m-3} \nonumber \\
& +107\cdot2^{4s+2m-7}+7\cdot2^{2s+2m-3}+105\cdot2^{k+3s+m-7}+15\cdot2^{6s+3m-6}\nonumber \\
& -15\cdot2^{4s+3m-6}-15\cdot2^{3s+3m-3} + 15\cdot2^{s+3m-3}-217\cdot2^{5s+2m-8}+41\cdot2^{s+2m-3}\nonumber
& \nonumber
 \end{align}
 Setting  $k=2s+m$ in \eqref{eq 12.13} we get \\[0.1 cm]
  \begin{align}
& \sum_{j=1}^{s-1}\Gamma_{s+m+j}^{\left[s\atop{ s +m\atop s+m} \right]\times (2s+m)}\cdot2^{-(s+m+j)} \label{eq 12.14} \\
& =-2^{4s+3m-6}-19\cdot2^{4s+2m-7}-7\cdot2^{2s+3m-3} -7\cdot2^{2s+2m-2} \nonumber \\
& -7\cdot2^{5s+2m-8}-15\cdot2^{3s+3m-3}+15\cdot2^{6s+3m-6}+15\cdot2^{s+3m-3}+41\cdot2^{s+2m-3}\nonumber
 \end{align}
  \end{lem}
 \begin{proof}
 The assertions follow from Lemma \ref{lem 12.5}
 
 \end{proof}
 \newpage
   \begin{lem}
 \label{lem 12.7}
   We have  for $m\geq 2,\; 0\leq j\leq m$ :
     \begin{align}
& \Gamma_{2s+m+j}^{\left[s\atop{ s +m\atop s+m} \right]\times k}\nonumber \\
& =  \begin{cases}\label{eq 12.15}
 2^{9s+5m-3}-7\cdot2^{8s+4m-6}+7\cdot2^{7s+3m-8}+3\cdot2^{6s+3m-7}+2^{6s+4m-6} &\text{if  } j=0 ,\;k=2s+m\\
 21\cdot\big[ 5\cdot2^{8s+4m-6 }+ 5\cdot2^{6s+3m-7}
  -2^{6s+4m-6}- 15\cdot2^{7s+3m-8}\big ] &\text{if  } j=0 ,\;k=2s+m+1\\
 2^{9s+5m}-7\cdot2^{8s+4m-2}+7\cdot2^{7s+3m-3}-3\cdot2^{6s+3m-3}+2^{6s+4m-2} &\text{if  } j=1 ,\;k=2s+m+1\\ 
  3\cdot2^{2k+2s+m-2}+21\cdot2^{k+4s+3m-5}+735\cdot2^{k+5s+2m-7} - 477\cdot2^{k+4s+2m-6}\\
 +105\cdot2^{8s+4m-6} -105\cdot2^{6s+4m-6}-3255\cdot2^{7s+3m-8}+1629\cdot2^{6s+3m-7}  &\text{if  } j=0 ,\;k\geq 2s+m+1\\ 
  21\cdot\big[ 5\cdot2^{8s+4m-2 }+ 7\cdot2^{6s+3m-3}
  -2^{6s+4m-2}- 15\cdot2^{7s+3m-3}\big ] &\text{if  } j=1,\;k=2s+m+2\\
  2^{9s+5m+3}-7\cdot2^{8s+4m+2}+7\cdot2^{7s+3m+2}-3\cdot2^{6s+3m+2}+2^{6s+4m+2} &\text{if  } j=2 ,\;k=2s+m+2\\ 
  21\cdot2^{2k+2s+m-2}+21\cdot2^{k+4s+3m-2}+735\cdot2^{k+5s+2m-3} - 945\cdot2^{k+4s+2m-3}\\
 +105\cdot2^{8s+4m-2} -105\cdot2^{6s+4m-2}-3255\cdot2^{7s+3m-3}+3255\cdot2^{6s+3m-3}  &\text{if  } j=1,\;k\geq 2s+m+2\\
   21\cdot\big[ 5\cdot2^{8s+4m+2 }+ 7\cdot2^{6s+3m+2}
  -2^{6s+4m+2}- 15\cdot2^{7s+3m+2}\big ] &\text{if  } j=2,\;k=2s+m+3\\
  2^{9s+6m+6}-7\cdot2^{8s+4m+6}+7\cdot2^{7s+3m+7}-3\cdot2^{6s+3m+7}+2^{6s+4m+6} &\text{if  } j=3 ,\;k=2s+m+3\\ 
  21\cdot2^{2k+2s+m+1}+21\cdot2^{k+4s+3m+1}+735\cdot2^{k+5s+2m+1} - 945\cdot2^{k+4s+2m+1}\\
 +105\cdot2^{8s+4m+2} -105\cdot2^{6s+4m+2}-3255\cdot2^{7s+3m+2}+3255\cdot2^{6s+3m+2}  &\text{if  } j=2,\;k\geq 2s+m+3
  \end{cases}
  \end{align}
 \end{lem}
 \begin{proof}
 
 \textbf{The case $j=0,\;k=2s+m$}\\[0.1 cm]
From \eqref{eq 6.1} with $k=2s+m,\;l=0$ we obtain using \eqref{eq 12.7},\eqref{eq 12.8} and \eqref{eq 12.12}\\[0.1 cm]
\begin{align*}
& \sum_{i=0}^{2s+m}\Gamma_{i}^{\left[s\atop{ s +m\atop s+m} \right]\times (2s+m)}= 2^{9s+5m-3}\\
&\Updownarrow \\
&  \sum_{i=0}^{s-1}\Gamma_{i}^{\left[s\atop{ s +m\atop s+m} \right]\times (2s+m)} + \Gamma_{s}^{\left[s\atop{ s +m\atop s+m} \right]\times (2s+m)}
+\sum_{j=1}^{m-1} \Gamma_{s+j}^{\left[s\atop{ s +m\atop s+m} \right]\times (2s+m)}\\
& + \Gamma_{s+m}^{\left[s\atop{ s +m\atop s+m} \right]\times (2s+m)}
 +  \sum_{j=1}^{s-1} \Gamma_{s+m+j}^{\left[s\atop{ s +m\atop s+m} \right]\times (2s+m)} +  \Gamma_{2s+m}^{\left[s\atop{ s +m\atop s+m} \right]\times (2s+m)}= 2^{9s+6m-3}\\
 & \Updownarrow \\
 & \big(7\cdot2^{4s-6}-3\cdot2^{3s-5}\big) + \big(2^{3s+m-1}-2^{2s}+ 21\cdot(5\cdot2^{4s-6}-2^{3s-5})\big) \\
   & +\big( 3\cdot2^{3s+4m-2}-3\cdot2^{3s+3m-5}+3\cdot2^{2s+3m-5}-2^{3s+m-1} \\
   & +3\cdot2^{3s-2}+2^{2s}+7\cdot2^{4s+4m-6}-7\cdot2^{2s+4m-6}-7\cdot2^{4s-2}-2^{3s+3m-4} \big)\\
    & +\big(  21\cdot2^{3s+4m-5}+45\cdot2^{3s+3m-4}-21\cdot2^{2s+4m-4}-45\cdot2^{2s+3m-3} \\
   & +105\cdot2^{4s+4m-6}-21\cdot2^{3s+3m-5}-21\cdot2^{2s+4m-6}+9\cdot2^{2s+3m-5} \big) \\
    & +\big( 7\cdot2^{8s+4m-6}+7\cdot2^{2s+4m-2}+21\cdot2^{2s+3m-2}-7\cdot2^{4s+4m-2} \\
   & -7\cdot2^{7s+3m-8}-3\cdot2^{3s+4m-2}-3\cdot2^{6s+3m-7}-2^{6s+4m-6}-2^{3s+3m+1} \big)\\
   & + \Gamma_{2s+m}^{\left[s\atop{ s +m\atop s+m} \right]\times (2s+m)} = 2^{9s+5m-3}\\
& \Updownarrow \\
& \Gamma_{2s+m}^{\left[s\atop{ s +m\atop s+m} \right]\times (2s+m)} 
 = 2^{9s+5m-3}-7\cdot2^{8s+4m-6}+3\cdot2^{6s+3m-7}+7\cdot2^{7s+3m-8}+2^{6s+4m-6}
\end{align*}
 \textbf{The case $j=0,\;k=2s+m+1$ and the case  $j=1,\;k=2s+m+1.$}\\[0.1 cm]
  From \eqref{eq 6.1} with $k=2s+m+1,\;l=0$ we obtain using \eqref{eq 12.7} and \eqref{eq 12.8},\eqref{eq 12.11} \\[0.1 cm]
\begin{align*}
& \sum_{i=0}^{2s+m+1}\Gamma_{i}^{\left[s\atop{ s +m\atop s+m} \right]\times (2s+m+1)}= 2^{9s+5m}\\
&\Updownarrow \\
&  \sum_{i=0}^{s-1}\Gamma_{i}^{\left[s\atop{ s +m\atop s+m} \right]\times (2s+m+1)} + \Gamma_{s}^{\left[s\atop{ s +m\atop s+m} \right]\times (2s+m+1)}
+\sum_{j=1}^{m-1} \Gamma_{s+j}^{\left[s\atop{ s +m\atop s+m} \right]\times (2s+m+1)}\\
& + \Gamma_{s+m}^{\left[s\atop{ s +m\atop s+m} \right]\times (2s+m+1)}
+ \sum_{j=1}^{s-1} \Gamma_{s+m+j}^{\left[s\atop{ s +m\atop s+m} \right]\times (2s+m+1)}
 + \Gamma_{2s+m}^{\left[s\atop{ s +m\atop s+m} \right]\times (2s+m+1)} +  \Gamma_{2s+m+1}^{\left[s\atop{ s +m\atop s+m} \right]\times (2s+m+1)}= 2^{9s+5m}\\
 & \Updownarrow \\
 & \big(7\cdot2^{4s-6}-3\cdot2^{3s-5}\big) + \big(2^{3s+m}-2^{2s}+ 21\cdot(5\cdot2^{4s-6}-2^{3s-5})\big) \\
  & +\big( 3\cdot2^{3s+4m-4}-3\cdot2^{3s+3m-5}+3\cdot2^{2s+3m-5}-2^{3s+m} \\
   & +3\cdot2^{3s-2}+2^{2s}+7\cdot2^{4s+4m-6}-7\cdot2^{2s+4m-6}-7\cdot2^{4s-2}-2^{3s+3m-3} \big)\\
   & +\big(  21\cdot2^{3s+4m-4}+45\cdot2^{3s+3m-3}-21\cdot2^{2s+4m-4}-45\cdot2^{2s+3m-3} \\
   & +105\cdot2^{4s+4m-6}-21\cdot2^{3s+3m-5}-21\cdot2^{2s+4m-6}+9\cdot2^{2s+3m-5} \big) \\
 & + \big(3\cdot2^{6s+4m-4}- 27\cdot2^{6s+3m-5}-3\cdot2^{3s+4m-1} -11\cdot2^{3s+3m-1}  \\
& +49\cdot2^{7s+3m-6}+51\cdot2^{6s+3m-7}+3\cdot2^{3s+3m-2}+7\cdot2^{8s+4m-6}-7\cdot2^{6s+4m-6}\\
& -7\cdot2^{4s+4m-2} + 7\cdot2^{2s+4m-2} -105\cdot2^{7s+3m-8}+21\cdot2^{2s+3m-2} \big)\\
 & + \Gamma_{2s+m}^{\left[s\atop{ s +m\atop s+m} \right]\times (2s+m+1)} + \Gamma_{2s+m+1}^{\left[s\atop{ s +m\atop s+m} \right]\times (2s+m+1)}= 2^{9s+5m}\\
& \Updownarrow \\
& \Gamma_{2s+m}^{\left[s\atop{ s +m\atop s+m} \right]\times (2s+m+1)} + \Gamma_{2s+m+1}^{\left[s\atop{ s +m\atop s+m} \right]\times (2s+m+1)}\\[0.1 cm]
& = 2^{9s+5m}-7\cdot2^{8s+4m-6}-5\cdot2^{6s+4m-6}+ 57\cdot2^{6s+3m-7}-91\cdot2^{7s+3m-8}
\end{align*}
   From \eqref{eq 6.2} with $k=2s+m+1,\;l=0$ we obtain using \eqref{eq 12.7} and \eqref{eq 12.9},\eqref{eq 12.13}\\[0.1 cm]
\begin{align*}
& \sum_{i=0}^{2s+m+1}\Gamma_{i}^{\left[s\atop{ s +m\atop s+m} \right]\times (2s+m+1)}\cdot2^{-i}= 2^{7s+4m-1}+ 2^{6s+3m}- 2^{4s+2m-1}\\
&\Updownarrow \\
&  \sum_{i=0}^{s-1}\Gamma_{i}^{\left[s\atop{ s +m\atop s+m} \right]\times (2s+m+1)}\cdot2^{-i} + \Gamma_{s}^{\left[s\atop{ s +m\atop s+m} \right]\times (2s+m+1)}\cdot2^{-s}
+\sum_{j=1}^{m-1} \Gamma_{s+j}^{\left[s\atop{ s +m\atop s+m} \right]\times (2s+m+1)}\cdot2^{-(s+j)}\\
& + \Gamma_{s+m}^{\left[s\atop{ s +m\atop s+m} \right]\times (2s+m+1)}\cdot2^{-(s+m)}
+ \sum_{j=1}^{s-1} \Gamma_{s+m+j}^{\left[s\atop{ s +m\atop s+m} \right]\times (2s+m+1)}\cdot2^{-(s+m+j)}\\
& + \Gamma_{2s+m}^{\left[s\atop{ s +m\atop s+m} \right]\times (2s+m+1)}\cdot2^{-(2s+m)} +  \Gamma_{2s+m+1}^{\left[s\atop{ s +m\atop s+m} \right]\times (2s+m+1)}\cdot2^{-(2s+m+1)}= 2^{7s+4m-1}+ 2^{6s+3m}- 2^{4s+2m-1}\\
 & \Updownarrow \\
  & \big(15\cdot2^{3s-6}-7\cdot2^{2s-5}\big) + \big(2^{2s+m}-2^{s}+ 21\cdot(5\cdot2^{3s-6}-2^{2s-5})\big) \\
   & = \big(7\cdot2^{2s+3m-4}- 7\cdot2^{2s+2m-5}+7\cdot2^{s+2m-5} -2^{2s+m}  \\
& +7\cdot2^{2s-3}+2^{s}+15\cdot2^{3s+3m-6}-15\cdot2^{s+3m-6}-15\cdot2^{3s-3}-3\cdot2^{2s+2m-3}\big)\\
 & +\big(  21\cdot2^{2s+3m-4}+45\cdot2^{2s+2m-3}-21\cdot2^{s+3m-4}-45\cdot2^{s+2m-3} \\
   & +105\cdot2^{3s+3m-6}-21\cdot2^{2s+2m-5}-21\cdot2^{s+3m-6}+9\cdot2^{s+2m-5} \big) \\
   & + \big(7\cdot2^{4s+3m-4}- 63\cdot2^{4s+2m-5}-7\cdot2^{2s+3m-2} -21\cdot2^{2s+2m-2}  \\
& +107\cdot2^{4s+2m-7}+7\cdot2^{2s+2m-3}+105\cdot2^{5s+2m-6}+15\cdot2^{6s+3m-6}\\
& -15\cdot2^{4s+3m-6}-15\cdot2^{3s+3m-3} + 15\cdot2^{s+3m-3}-217\cdot2^{5s+2m-8}+41\cdot2^{s+2m-3}\big)\\
  & + \Gamma_{2s+m}^{\left[s\atop{ s +m\atop s+m} \right]\times (2s+m+1)}\cdot2^{-(2s+m)} + \Gamma_{2s+m+1}^{\left[s\atop{ s +m\atop s+m} \right]\times (2s+m+1)}\cdot2^{-(2s+m+1)}= 2^{7s+4m-1}+ 2^{6s+3m}- 2^{4s+2m-1}\\
& \Updownarrow \\
& 2\cdot\Gamma_{2s+m}^{\left[s\atop{ s +m\atop s+m} \right]\times (2s+m+1)} + \Gamma_{2s+m+1}^{\left[s\atop{ s +m\atop s+m} \right]\times (2s+m+1)}\\[0.1 cm]
& = 2^{9s+5m}+49\cdot2^{8s+4m-5}+81\cdot2^{6s+3m-6}-13\cdot2^{6s+4m-5}-203\cdot2^{7s+3m-7}
\end{align*}
    \textbf{The case $j=0,\;k\geq 2s+m+1$}\\[0.1 cm]
  From \eqref{eq 12.2} with $j=s$ we get \\[0.1 cm]
  \begin{align*}
& \Gamma_{2s+m}^{\left[s\atop{ s +m\atop s+m} \right]\times k} -\Gamma_{2s+m}^{\left[s\atop{ s +m\atop s+m} \right]\times (2s+m+1)} \\
& = 3\cdot2^{2s+m-2}\cdot\big(2^{2k}-2^{4s+2m+2}\big) 
+   \big(2^{k}-2^{2s+m+1}\big)\cdot\big(21\cdot2^{4s+3m-5}+735\cdot2^{5s+2m-7}-477\cdot2^{4s+2m-6}\big)   \\
& \Updownarrow \\
& \Gamma_{2s+m}^{\left[s\atop{ s +m\atop s+m} \right]\times k}= 
   21\cdot\big[ 5\cdot2^{8s+4m-6 }+ 5\cdot2^{6s+3m-7} -2^{6s+4m-6}- 15\cdot2^{7s+3m-8}\big ] \\
 & +  3\cdot2^{2s+m-2}\cdot\big(2^{2k}-2^{4s+2m+2}\big) 
+   \big(2^{k}-2^{2s+m+1}\big)\cdot\big(21\cdot2^{4s+3m-5}+735\cdot2^{5s+2m-7}-477\cdot2^{4s+2m-6}\big)  \\
& \Updownarrow \\
& \Gamma_{2s+m}^{\left[s\atop{ s +m\atop s+m} \right]\times k} =  3\cdot2^{2k+2s+m-2}+21\cdot2^{k+4s+3m-5}+735\cdot2^{k+5s+2m-7} - 477\cdot2^{k+4s+2m-6}\\
& +105\cdot2^{8s+4m-6} -105\cdot2^{6s+4m-6}-3255\cdot2^{7s+3m-8}+1629\cdot2^{6s+3m-7} 
  \end{align*}    
    \textbf{The case $j=1,\;k=2s+m+2$ and the case  $j=2,\;k=2s+m+1$}\\[0.1 cm]  
   From \eqref{eq 6.1} with $k=2s+m+2,\;l=0$ we obtain using \eqref{eq 12.7} and \eqref{eq 12.8},\eqref{eq 12.11} \\[0.1 cm]
\begin{align*}
& \sum_{i=0}^{2s+m+2}\Gamma_{i}^{\left[s\atop{ s +m\atop s+m} \right]\times (2s+m+2)}= 2^{9s+5m+3}\\
&\Updownarrow \\
&  \sum_{i=0}^{s-1}\Gamma_{i}^{\left[s\atop{ s +m\atop s+m} \right]\times (2s+m+2)} + \Gamma_{s}^{\left[s\atop{ s +m\atop s+m} \right]\times (2s+m+2)}
+\sum_{j=1}^{m-1} \Gamma_{s+j}^{\left[s\atop{ s +m\atop s+m} \right]\times (2s+m+2)}\\
& + \Gamma_{s+m}^{\left[s\atop{ s +m\atop s+m} \right]\times (2s+m+2)}
+ \sum_{j=1}^{s-1} \Gamma_{s+m+j}^{\left[s\atop{ s +m\atop s+m} \right]\times (2s+m+2)}
 + \Gamma_{2s+m}^{\left[s\atop{ s +m\atop s+m} \right]\times (2s+m+2)}\\
&  +  \Gamma_{2s+m+1}^{\left[s\atop{ s +m\atop s+m} \right]\times (2s+m+2)}+ \Gamma_{2s+m+2}^{\left[s\atop{ s +m\atop s+m} \right]\times (2s+m+2)}= 2^{9s+5m+3}\\
 & \Updownarrow \\ 
 & \big(7\cdot2^{4s-6}-3\cdot2^{3s-5}\big) + \big(2^{3s+m}-2^{2s}+ 21\cdot(5\cdot2^{4s-6}-2^{3s-5})\big) \\
  & +\big( 3\cdot2^{3s+4m-4}-3\cdot2^{3s+3m-5}+3\cdot2^{2s+3m-5}-2^{3s+m} \\
   & +3\cdot2^{3s-2}+2^{2s}+7\cdot2^{4s+4m-6}-7\cdot2^{2s+4m-6}-7\cdot2^{4s-2}-2^{3s+3m-3} \big)\\
   & +\big(  21\cdot2^{3s+4m-3}+45\cdot2^{3s+3m-2}-21\cdot2^{2s+4m-4}-45\cdot2^{2s+3m-3} \\
   & +105\cdot2^{4s+4m-6}-21\cdot2^{3s+3m-5}-21\cdot2^{2s+4m-6}+9\cdot2^{2s+3m-5} \big) \\
 &  \big(3\cdot2^{6s+4m-3}- 27\cdot2^{6s+3m-4}-3\cdot2^{3s+4m} -11\cdot2^{3s+3m}  \\
& +49\cdot2^{7s+3m-5}+51\cdot2^{6s+3m-7}+3\cdot2^{3s+3m-2}+7\cdot2^{8s+4m-6}-7\cdot2^{6s+4m-6}\\
& -7\cdot2^{4s+4m-2} + 7\cdot2^{2s+4m-2} -105\cdot2^{7s+3m-8}+21\cdot2^{2s+3m-2} \big)\\
& +  3\cdot2^{6s+3m+2}+21\cdot2^{6s+4m-3}+735\cdot2^{7s+3m-5} - 477\cdot2^{6s+3m-4}\\
& +105\cdot2^{8s+4m-6} -105\cdot2^{6s+4m-6}-3255\cdot2^{7s+3m-8}+1629\cdot2^{6s+3m-7} \\
 & + \Gamma_{2s+m+1}^{\left[s\atop{ s +m\atop s+m} \right]\times (2s+m+2)} + \Gamma_{2s+m+2}^{\left[s\atop{ s +m\atop s+m} \right]\times (2s+m+2)}= 2^{9s+5m+3}\\
& \Updownarrow \\
& \Gamma_{2s+m+1}^{\left[s\atop{ s +m\atop s+m} \right]\times (2s+m+2)} + \Gamma_{2s+m+2}^{\left[s\atop{ s +m\atop s+m} \right]\times (2s+m+2)}\\[0.1 cm]
& = 2^{9s+5m+3}-7\cdot2^{8s+4m-2}-5\cdot2^{6s+4m-2}+ 51\cdot2^{6s+3m-3}-91\cdot2^{7s+3m-3}
\end{align*}
 
    From \eqref{eq 6.2} with $k=2s+m+2,\;l=0$ we obtain using \eqref{eq 12.7} and \eqref{eq 12.9},\eqref{eq 12.13}\\[0.1 cm]
\begin{align*}
& \sum_{i=0}^{2s+m+2}\Gamma_{i}^{\left[s\atop{ s +m\atop s+m} \right]\times (2s+m+2)}\cdot2^{-i}= 2^{7s+4m+1}+ 2^{6s+3m+3}- 2^{4s+2m+1}\\
&\Updownarrow \\
&  \sum_{i=0}^{s-1}\Gamma_{i}^{\left[s\atop{ s +m\atop s+m} \right]\times (2s+m+2)}\cdot2^{-i} + \Gamma_{s}^{\left[s\atop{ s +m\atop s+m} \right]\times (2s+m+2)}\cdot2^{-s}
+\sum_{j=1}^{m-1} \Gamma_{s+j}^{\left[s\atop{ s +m\atop s+m} \right]\times (2s+m+2)}\cdot2^{-(s+j)}\\
& + \Gamma_{s+m}^{\left[s\atop{ s +m\atop s+m} \right]\times (2s+m+2)}\cdot2^{-(s+m)}
+ \sum_{j=1}^{s-1} \Gamma_{s+m+j}^{\left[s\atop{ s +m\atop s+m} \right]\times (2s+m+2)}\cdot2^{-(s+m+j)}
 + \Gamma_{2s+m}^{\left[s\atop{ s +m\atop s+m} \right]\times (2s+m+2)}\cdot2^{-(2s+m)}\\
& +  \Gamma_{2s+m+1}^{\left[s\atop{ s +m\atop s+m} \right]\times (2s+m+2)}\cdot2^{-(2s+m+1)}+\Gamma_{2s+m+2}^{\left[s\atop{ s +m\atop s+m} \right]\times (2s+m+2)}\cdot2^{-(2s+m+2)} = 2^{7s+4m-1}+ 2^{6s+3m}- 2^{4s+2m-1}\\
 & \Updownarrow \\
  & \big(15\cdot2^{3s-6}-7\cdot2^{2s-5}\big) + \big(2^{2s+m+1}-2^{s}+ 21\cdot(5\cdot2^{3s-6}-2^{2s-5})\big) \\
   & = \big(7\cdot2^{2s+3m-3}- 7\cdot2^{2s+2m-5}+7\cdot2^{s+2m-5} -2^{2s+m+1}  \\
& +7\cdot2^{2s-3}+2^{s}+15\cdot2^{3s+3m-6}-15\cdot2^{s+3m-6}-15\cdot2^{3s-3}-3\cdot2^{2s+2m-2}\big)\\
& +\big(  21\cdot2^{2s+3m-3}+45\cdot2^{2s+2m-2}-21\cdot2^{s+3m-4}-45\cdot2^{s+2m-3} \\
   & +105\cdot2^{3s+3m-6}-21\cdot2^{2s+2m-5}-21\cdot2^{s+3m-6}+9\cdot2^{s+2m-5} \big) \\
 & + \big(7\cdot2^{4s+3m-3}- 63\cdot2^{4s+2m-4}-7\cdot2^{2s+3m-1} -21\cdot2^{2s+2m-1}  \\
& +107\cdot2^{4s+2m-7}+7\cdot2^{2s+2m-3}+105\cdot2^{5s+2m-5}+15\cdot2^{6s+3m-6}\\
& -15\cdot2^{4s+3m-6}-15\cdot2^{3s+3m-3} + 15\cdot2^{s+3m-3}-217\cdot2^{5s+2m-8}+41\cdot2^{s+2m-3}\big)\\
& \\
& +  3\cdot2^{4s+2m+2}+21\cdot2^{4s+3m-3}+735\cdot2^{5s+2m-5} - 477\cdot2^{4s+2m-4}\\
& +105\cdot2^{6s+3m-6} -105\cdot2^{4s+3m-6}-3255\cdot2^{5s+2m-8}+1629\cdot2^{4s+2m-7} \\
 & + \Gamma_{2s+m+1}^{\left[s\atop{ s +m\atop s+m} \right]\times (2s+m+2)}\cdot2^{-(2s+m+1)} + \Gamma_{2s+m+2}^{\left[s\atop{ s +m\atop s+m} \right]\times (2s+m+2)}\cdot2^{-(2s+m+2)}= 2^{7s+4m+1}+ 2^{6s+3m+3}- 2^{4s+2m+1}\\
& \Updownarrow \\
& 2\cdot\Gamma_{2s+m+1}^{\left[s\atop{ s +m\atop s+m} \right]\times (2s+m+2)} + \Gamma_{2s+m+2}^{\left[s\atop{ s +m\atop s+m} \right]\times (2s+m+2)}\\[0.1 cm]
& = 2^{9s+5m+3}+49\cdot2^{8s+4m-1}+99\cdot2^{6s+3m-2}-13\cdot2^{6s+4m-1}-203\cdot2^{7s+3m-2},\\
& \text{hence we have (see above)}\\
& \Gamma_{2s+m+1}^{\left[s\atop{ s +m\atop s+m} \right]\times (2s+m+2)}=  21\cdot\big[ 5\cdot2^{8s+4m-2 }+ 7\cdot2^{6s+3m-3}
  -2^{6s+4m-2}- 15\cdot2^{7s+3m-3}\big ] \\
&  \Gamma_{2s+m+2}^{\left[s\atop{ s +m\atop s+m} \right]\times (2s+m+2)}=  2^{9s+6m}-7\cdot2^{8s+4m+2}+7\cdot2^{7s+3m+2}-3\cdot2^{6s+3m+2}+2^{6s+4m+2}
\end{align*}
  \textbf{The case $j=1,\;k\geq 2s+m+2$}\\[0.1 cm]
  From \eqref{eq 12.3} with $j=0$ we get \\[0.1 cm]
  \begin{align*}
& \Gamma_{2s+m+1}^{\left[s\atop{ s +m\atop s+m} \right]\times k} -\Gamma_{2s+m+1}^{\left[s\atop{ s +m\atop s+m} \right]\times (2s+m+2)} \\
& = 21\cdot2^{2s+m-2}\cdot\big(2^{2k}-2^{4s+2m+4}\big) 
+   \big(2^{k}-2^{2s+m+2}\big)\cdot\big(21\cdot2^{4s+3m-2}+735\cdot2^{5s+2m-3}-945\cdot2^{4s+2m-3}\big)   \\
& \Updownarrow \\
& \Gamma_{2s+m+1}^{\left[s\atop{ s +m\atop s+m} \right]\times k}= 
   21\cdot\big[ 5\cdot2^{8s+4m-2 }+ 7\cdot2^{6s+3m-3} -2^{6s+4m-2}- 15\cdot2^{7s+3m-3}\big ] \\
 & + 21\cdot2^{2s+m-2}\cdot\big(2^{2k}-2^{4s+2m+4}\big) 
+   \big(2^{k}-2^{2s+m+2}\big)\cdot\big(21\cdot2^{4s+3m-2}+735\cdot2^{5s+2m-3}-945\cdot2^{4s+2m-3}\big)   \\
& \Updownarrow \\
& \Gamma_{2s+m+1}^{\left[s\atop{ s +m\atop s+m} \right]\times k} =   21\cdot2^{2k+2s+m-2}+21\cdot2^{k+4s+3m-2}+735\cdot2^{k+5s+2m-3} - 945\cdot2^{k+4s+2m-3}\\
& +105\cdot2^{8s+4m-2} -105\cdot2^{6s+4m-2}-3255\cdot2^{7s+3m-3}+3255\cdot2^{6s+3m-3}
  \end{align*}
   \textbf{The case $j=2,\;k=2s+m+3$ and the case  $j=3,\;k=2s+m+3$}\\[0.1 cm]  
   From \eqref{eq 6.1} with $k=2s+m+3,\;l=0$ we obtain using \eqref{eq 12.7} and \eqref{eq 12.8},\eqref{eq 12.11} \\[0.1 cm]
\begin{align*}
& \sum_{i=0}^{2s+m+3}\Gamma_{i}^{\left[s\atop{ s +m\atop s+m} \right]\times (2s+m+3)}= 2^{9s+5m+6}\\
&\Updownarrow \\
&  \sum_{i=0}^{s-1}\Gamma_{i}^{\left[s\atop{ s+m\atop s+m} \right]\times (2s+m+3)} + \Gamma_{s}^{\left[s\atop{ s +m\atop s+m} \right]\times (2s+m+3)}
+\sum_{j=1}^{m-1} \Gamma_{s+j}^{\left[s\atop{ s +m\atop s+m} \right]\times (2s+m+3)}\\
& + \Gamma_{s+m}^{\left[s\atop{ s +m\atop s+m} \right]\times (2s+m+3)}
+ \sum_{j=1}^{s-1} \Gamma_{s+m+j}^{\left[s\atop{ s +m\atop s+m} \right]\times (2s+m+3)}
 + \Gamma_{2s+m}^{\left[s\atop{ s +m\atop s+m} \right]\times (2s+m+3)}
 +  \Gamma_{2s+m+1}^{\left[s\atop{ s +m\atop s+m} \right]\times (2s+m+3)}\\
 & + \Gamma_{2s+m+2}^{\left[s\atop{ s +m\atop s+m} \right]\times (2s+m+3)}+\Gamma_{2s+m+3}^{\left[s\atop{ s +m\atop s+m} \right]\times (2s+m+3)} = 2^{9s+5m+6}\\
 & \Updownarrow \\ 
 & \big(7\cdot2^{4s-6}-3\cdot2^{3s-5}\big) + \big(2^{3s+m+2}-2^{2s}+ 21\cdot(5\cdot2^{4s-6}-2^{3s-5})\big) \\
  & +\big( 3\cdot2^{3s+4m-2}-3\cdot2^{3s+3m-5}+3\cdot2^{2s+3m-5}-2^{3s+m+2} \\
   & +3\cdot2^{3s-2}+2^{2s}+7\cdot2^{4s+4m-6}-7\cdot2^{2s+4m-6}-7\cdot2^{4s-2}-2^{3s+3m-1} \big)\\
    & +\big(  21\cdot2^{3s+4m-2}+45\cdot2^{3s+3m-1}-21\cdot2^{2s+4m-4}-45\cdot2^{2s+3m-3} \\
   & +105\cdot2^{4s+4m-6}-21\cdot2^{3s+3m-5}-21\cdot2^{2s+4m-6}+9\cdot2^{2s+3m-5} \big) \\
    & + \big(3\cdot2^{6s+4m-2}- 27\cdot2^{6s+3m-3}-3\cdot2^{3s+4m+1} -11\cdot2^{3s+3m+1}  \\
& +49\cdot2^{7s+3m-4}+51\cdot2^{6s+3m-7}+3\cdot2^{3s+3m-2}+7\cdot2^{8s+4m-6}-7\cdot2^{6s+4m-6}\\
& -7\cdot2^{4s+4m-2} + 7\cdot2^{2s+4m-2} -105\cdot2^{7s+3m-8}+21\cdot2^{2s+3m-2} \big)\\
& +  3\cdot2^{6s+3m+4}+21\cdot2^{6s+4m-2}+735\cdot2^{7s+3m-4} - 477\cdot2^{6s+3m-3}\\
& +105\cdot2^{8s+4m-6} -105\cdot2^{6s+4m-6}-3255\cdot2^{7s+3m-8}+1629\cdot2^{6s+3m-7} \\
& +   21\cdot2^{6s+3m+4}+21\cdot2^{6m+4m+1}+735\cdot2^{7s+3m} - 945\cdot2^{6s+3m}\\
& +105\cdot2^{8s+4m-2} -105\cdot2^{6s+4m-2}-3255\cdot2^{7s+3m-3}+3255\cdot2^{6s+3m-3}  \\
 & + \Gamma_{2s+m+1}^{\left[s\atop{ s +m\atop s+m} \right]\times (2s+m+2)} + \Gamma_{2s+m+2}^{\left[s\atop{ s +m\atop s+m} \right]\times (2s+m+2)}= 2^{9s+5m+3}\\
& \Updownarrow \\
& \Gamma_{2s+m+1}^{\left[s\atop{ s +m\atop s+m} \right]\times (2s+m+2)} + \Gamma_{2s+m+2}^{\left[s\atop{ s +m\atop s+m} \right]\times (2s+m+2)}\\[0.1 cm]
& = 2^{9s+5m+6}-7\cdot2^{8s+4m+2}-5\cdot2^{6s+4m+2}+ 51\cdot2^{6s+3m+2}-91\cdot2^{7s+3m+2}
\end{align*}

    From \eqref{eq 6.2} with $k=2s+m+3,\;l=0$ we obtain using \eqref{eq 12.7} and \eqref{eq 12.9},\eqref{eq 12.13}\\[0.1 cm]
\begin{align*}
& \sum_{i=0}^{2s+m+3}\Gamma_{i}^{\left[s\atop{ s +m\atop s+m} \right]\times (2s+m+3)}\cdot2^{-i}= 2^{7s+4m+3}+ 2^{6s+3m+6}- 2^{4s+2m+3}\\
&\Updownarrow \\
&  \sum_{i=0}^{s-1}\Gamma_{i}^{\left[s\atop{ s +m\atop s+m} \right]\times (2s+m+3)}\cdot2^{-i} + \Gamma_{s}^{\left[s\atop{ s +m\atop s+m} \right]\times (2s+m+3)}\cdot2^{-s}
+\sum_{j=1}^{m-1} \Gamma_{s+j}^{\left[s\atop{ s +m\atop s+m} \right]\times (2s+m+3)}\cdot2^{-(s+j)}\\
& + \Gamma_{s+m}^{\left[s\atop{ s +m\atop s+m} \right]\times (2s+m+3)}\cdot2^{-(s+m)}
+ \sum_{j=1}^{s-1} \Gamma_{s+m+j}^{\left[s\atop{ s +m\atop s+m} \right]\times (2s+m+3)}\cdot2^{-(s+m+j)}
 + \Gamma_{2s+m}^{\left[s\atop{ s +m\atop s+m} \right]\times (2s+m+3)}\cdot2^{-(2s+m)}\\
& +  \Gamma_{2s+m+1}^{\left[s\atop{ s +m\atop s+m} \right]\times (2s+m+3)}\cdot2^{-(2s+m+1)}
+\Gamma_{2s+m+2}^{\left[s\atop{ s +m\atop s+m} \right]\times (2s+m+3)}\cdot2^{-(2s+m+2)}\\
 & +\Gamma_{2s+m+3}^{\left[s\atop{ s +m\atop s+m} \right]\times (2s+m+3)}\cdot2^{-(2s+m+3)} = 2^{7s+4m+3}+ 2^{6s+3m+6}- 2^{4s+2m+3}\\
 & \Updownarrow \\
  & \big(15\cdot2^{3s-6}-7\cdot2^{2s-5}\big) + \big(2^{2s+m+2}-2^{s}+ 21\cdot(5\cdot2^{3s-6}-2^{2s-5})\big) \\
   & + \big(7\cdot2^{2s+3m-2}- 7\cdot2^{2s+2m-5}+7\cdot2^{s+2m-5} -2^{2s+m+2}  \\
& +7\cdot2^{2s-3}+2^{s}+15\cdot2^{3s+3m-6}-15\cdot2^{s+3m-6}-15\cdot2^{3s-3}-3\cdot2^{2s+2m-1}\big)\\
& +\big(  21\cdot2^{2s+3m-2}+45\cdot2^{2s+2m-1}-21\cdot2^{s+3m-4}-45\cdot2^{s+2m-3} \\
   & +105\cdot2^{3s+3m-6}-21\cdot2^{2s+2m-5}-21\cdot2^{s+3m-6}+9\cdot2^{s+2m-5} \big) \\
 & + \big(7\cdot2^{4s+3m-2}- 63\cdot2^{4s+2m-3}-7\cdot2^{2s+3m} -21\cdot2^{2s+2m}  \\
& +107\cdot2^{4s+2m-7}+7\cdot2^{2s+2m-3}+105\cdot2^{5s+2m-4}+15\cdot2^{6s+3m-6}\\
& -15\cdot2^{4s+3m-6}-15\cdot2^{3s+3m-3} + 15\cdot2^{s+3m-3}-217\cdot2^{5s+2m-8}+41\cdot2^{s+2m-3}\big)\\
& +  3\cdot2^{4s+2m+4}+21\cdot2^{4s+3m-2}+735\cdot2^{5s+2m-4} - 477\cdot2^{4s+2m-3}\\
& +105\cdot2^{6s+3m-6} -105\cdot2^{4s+3m-6}-3255\cdot2^{5s+2m-8}+1629\cdot2^{4s+2m-7} \\
& + 21\cdot2^{4s+2m+3}+21\cdot2^{4s+3m}+735\cdot2^{5s+2m-1} - 945\cdot2^{4s+2m-1}\\
& +105\cdot2^{6s+3m-3} -105\cdot2^{4s+3m-3}-3255\cdot2^{5s+2m-4}+3255\cdot2^{4s+2m-4}\\
 & + \Gamma_{2s+m+2}^{\left[s\atop{ s +m\atop s+m} \right]\times (2s+m+3)}\cdot2^{-(2s+m+2)} + \Gamma_{2s+m+3}^{\left[s\atop{ s +m\atop s+m} \right]\times (2s+m+3)}\cdot2^{-(2s+m+3)}= 2^{7s+4m+3}+ 2^{6s+3m+6}- 2^{4s+2m+3}\\
& \Updownarrow \\
& 2\cdot\Gamma_{2s+m+2}^{\left[s\atop{ s +m\atop s+m} \right]\times (2s+m+3)} + \Gamma_{2s+m+3}^{\left[s\atop{ s +m\atop s+m} \right]\times (2s+m+3)}\\[0.1 cm]
& = 2^{9s+5m+6}+49\cdot2^{8s+4m+3}+99\cdot2^{6s+3m+3}-13\cdot2^{6s+4m+3}-203\cdot2^{7s+3m+3},\\
& \text{hence we have (see above)}\\
& \Gamma_{2s+m+2}^{\left[s\atop{ s +m\atop s+m} \right]\times (2s+m+3)}=  21\cdot\big[ 5\cdot2^{8s+4m+2 }+ 7\cdot2^{6s+3m+2}
  -2^{6s+4m+2}- 15\cdot2^{7s+3m+2}\big ] \\
&  \Gamma_{2s+m+3}^{\left[s\atop{ s +m\atop s+m} \right]\times (2s+m+3)}=  2^{9s+6m+6}-7\cdot2^{8s+4m+6}+7\cdot2^{7s+3m+7}-3\cdot2^{6s+3m+7}+2^{6s+4m+6}
\end{align*}

   \textbf{The case $j=2,\;k\geq 2s+m+3$}\\[0.1 cm]
  From \eqref{eq 12.3} with $j=1$ we get \\[0.1 cm]
  \begin{align*}
& \Gamma_{2s+m+2}^{\left[s\atop{ s +m\atop s+m} \right]\times k} -\Gamma_{2s+m+2}^{\left[s\atop{ s +m\atop s+m} \right]\times (2s+m+3)} \\
& = 21\cdot2^{2s+m+1}\cdot\big(2^{2k}-2^{4s+2m+6}\big) 
+   \big(2^{k}-2^{2s+m+3}\big)\cdot\big(21\cdot2^{4s+3m+1}+735\cdot2^{5s+2m+1}-945\cdot2^{4s+2m+1}\big)   \\
& \Updownarrow \\
& \Gamma_{2s+m+2}^{\left[s\atop{ s +m\atop s+m} \right]\times k}= 
 21\cdot\big[ 5\cdot2^{8s+4m+2 }+ 7\cdot2^{6s+3m+2} -2^{6s+4m+2}- 15\cdot2^{7s+3m+2}\big ] \\  
 & +   21\cdot2^{2s+m+1}\cdot\big(2^{2k}-2^{4s+2m+6}\big) 
+   \big(2^{k}-2^{2s+m+3}\big)\cdot\big(21\cdot2^{4s+3m+1}+735\cdot2^{5s+2m+1}-945\cdot2^{4s+2m+1}\big)  \\
& \Updownarrow \\
& \Gamma_{2s+m+2}^{\left[s\atop{ s +m\atop s+m} \right]\times k} =  21\cdot2^{2k+2s+m+1}+21\cdot2^{k+4s+3m+1}+735\cdot2^{k+5s+2m+1} - 945\cdot2^{k+4s+2m+1}\\
& +105\cdot2^{8s+4m+2} -105\cdot2^{6s+4m+2}-3255\cdot2^{7s+3m+2}+3255\cdot2^{6s+3m+2} 
  \end{align*}
   \end{proof}
   
  \subsection{\textbf{An explicit formula for $ \Gamma_{i}^{\left[s\atop{ s+1\atop s+1} \right]\times k} , \quad for\; s \leq i \leq 3s+2,\; k\geq i $}}
 \label{subsec 6}

We shall use the following results obtained in Lemma \ref{lem 12.2} 

 \begin{lem}
 \label{lem 12.8}

We have for $ m=1,\; k\geq s+j+1  $\\[0.05 cm]
  \begin{align}
& \Gamma_{s+j}^{\left[s\atop{ s +1\atop s+1} \right]\times k} -\Gamma_{s+j}^{\left[s\atop{ s +1\atop s+1} \right]\times (s+j+1)}\label{eq 12.16} \\
& =  \begin{cases}
2^{s-1}\cdot\big(2^{k}-2^{s+1}\big)  &\text{if  }  j = 0 \\
33\cdot2^{s-1}\cdot\big(2^{k}-2^{s+2}\big)  &\text{if  }  j = 1 \\
\big(735\cdot2^{s+4j-9}-105\cdot2^{s+3j-7}\big)\cdot\big(2^{k}-2^{s+j+1}\big)& \text{if   }\; 2\leq j\leq s\\
3\cdot2^{2s-1}\cdot(2^{2k}-2^{4s+4}) +
\big(735\cdot2^{5s-5}-393\cdot2^{4s-4}\big)\cdot\big(2^{k}-2^{2s+2}\big)  & \text{if   } \; j=s+1, \\
53\cdot2^{2s-1}\cdot(2^{2k}-2^{4s+6}) +
\big(735\cdot2^{5s-1}-1629\cdot2^{4s-1}\big)\cdot\big(2^{k}-2^{2s+3}\big)  & \text{if   } \; j=s+2, \\
\end{cases}\nonumber \\
& \nonumber
\end{align}

We have for $ m=1,\; k\geq 2s+2+j  $\\[0.05 cm]
  \begin{align}
& \Gamma_{2s+1+j}^{\left[s\atop{ s +1\atop s+1} \right]\times k} -\Gamma_{2s+1+j}^{\left[s\atop{ s +1\atop s+1} \right]\times (2s+2+j)}\label{eq 12.17} \\
& =  \begin{cases}
105\cdot2^{2s+4j-6}\cdot(2^{2k}-2^{4s+2j+4}) +
\big(735\cdot2^{5s+4j-5}-3255\cdot2^{4s+5j-7}\big)\cdot\big(2^{k}-2^{2s+j+2}\big)  & \text{if   } \; 2\leq j\leq s, \\
2^{3s-1}\cdot(2^{3k}-2^{9s+6}) -7\cdot2^{6s-2}\cdot\big(2^{2k}-2^{6s+6}\big)  
+ 7\cdot2^{9s-2}\cdot\big(2^{k}-2^{3s+3}\big) &\text{if  } j=s+1\\
 \end{cases}\nonumber \\
& \nonumber
\end{align}
 \end{lem}
  \begin{thm}
 \label{thm 12.9}
 We have \\
 
  \begin{equation}
\label{eq 12.18}
\Gamma_{s+j}^{\left[s\atop{ s+1\atop s+1} \right]\times (s+j+1)} 
 = \begin{cases}
  105\cdot2^{4s-6} - 21\cdot2^{3s -5}  &  \text{if  }  j = 0, \\
 105\cdot2^{4s-2} - 21\cdot2^{3s -2}-3\cdot2^{2s} &    \text{if }  j = 1,  \\
  21\cdot[5\cdot2^{4s+2} -2^{3s+1}-5\cdot2^{2s+1} ]  &    \text{if }  j = 2,  \\
  21\cdot[5\cdot2^{4s+6} -2^{3s+4}-25\cdot2^{2s+4} ]  &    \text{if }  j = 3,  \\ 
    21\cdot[5\cdot2^{4s+10} -2^{3s+7}-55\cdot2^{2s+8} ]  &    \text{if }  j = 4,  \\ 
   21\cdot\big[5\cdot2^{4s+4j-6} - 2^{3s+3j-5} - 15\cdot2^{2s+5j-10}+ 5\cdot2^{2s+4j-8}\big]    &    \text{if }  2\leq j\leq s+2,  \\  
    21\cdot[5\cdot2^{8s-2} +2^{6s-4}-15\cdot2^{7s-5} ]  &    \text{if }  j = s+1,  \\ 
     21\cdot[5\cdot2^{8s+2} +3\cdot2^{6s}-15\cdot2^{7s} ]  &    \text{if }  j = s+2
     \end{cases} 
     \end{equation}
     
   \begin{equation}
   \label{eq 12.19}
\Gamma_{2s+3+j}^{\left[s\atop{ s+1\atop s+1} \right]\times (2s+4+j)} 
 = \begin{cases}
  21\cdot[5\cdot2^{8s+4j+6} +5\cdot2^{6s+6j+5}-15\cdot2^{7s+5j+5} ]  &    \text{if } 0\leq j\leq s-2,  \\ 
 315\cdot2^{12s-1}   &    \text{if } j = s-1                      
    \end{cases} 
\end{equation}
  
    We have \\
 \begin{equation}
\label{eq 12.20}
\Gamma_{i}^{\left[s\atop{ s+1\atop s+1} \right]\times k} 
 = \begin{cases}
 1  & \text{if  } i = 0, \\
 105\cdot2^{4i-6} - 21\cdot2^{3i -5}  &  \text{if  }  1\leq i\leq s-1,\;k\geq i+1, \\
  2^{k+s-1} -2^{2s} +105\cdot2^{4s-6} -21\cdot2^{3s-5}   &    \text{if }  i = s, \; k\geq s+1,\; s\geq 1,  \\
   33\cdot2^{k+s-1} + [105\cdot2^{4s-2} -21\cdot2^{3s-2}-69\cdot2^{2s} ]  &    \text{if }  i = s+1, \; k\geq s+2,  \\
 630\cdot2^{k+s-1} + 21\cdot[5\cdot2^{4s+2} -2^{3s+1}-65\cdot2^{2s+1} ]  &    \text{if }  i = s+2, \; k\geq s+3,  \\
  1365\cdot2^{k+s+2} + 21\cdot[5\cdot2^{4s+6} -2^{3s+4}-285\cdot2^{2s+4} ]  &    \text{if }  i = s+3, \; k\geq s+4,  \\
  2835\cdot2^{k+s+5} + 21\cdot[5\cdot2^{4s+10} -2^{3s+7}-595\cdot2^{2s+8} ]  &    \text{if }  i = s+4, \; k\geq s+5,  \\
   105\cdot(7\cdot2^{j-2} - 1)\cdot2^{k+s+3j-7} \\
    + 21\cdot\big[5\cdot2^{4s+4j-6} - 2^{3s+3j-5} - 155\cdot2^{2s+5j-10} +25\cdot2^{2s+4j-8} \big]
     &  \text{if  }  i = s+j, \\
     & \; 2\leq j\leq s,\; k\geq s+j+1,\\
    3\cdot2^{2s-1}\cdot(2^{2k}-2^{4s+4}) 
+ \big(735\cdot2^{5s-5}-393\cdot2^{4s-4}\big)\cdot\big(2^{k}-2^{2s+2}\big)  \\
+   21\cdot[5\cdot2^{8s-2} +2^{6s-4}-15\cdot2^{7s-5} ]  &    \text{if }  i = 2s+1, \; k\geq 2s+2,  \\
    53\cdot2^{2s-1}\cdot(2^{2k}-2^{4s+6}) 
+ \big(735\cdot2^{5s-1}-1629\cdot2^{4s-1}\big)\cdot\big(2^{k}-2^{2s+3}\big)  \\
+   21\cdot[5\cdot2^{8s+2} +3\cdot2^{6s}-15\cdot2^{7s} ]  &    \text{if }  i = 2s+2, \; k\geq 2s+3,  \\
 105\cdot\big(2^{2k+2s+4j +2}+ 7\cdot2^{k+5s+4j+3} - 31\cdot2^{k+4s+5j+3} \big)  \\
     + 105\cdot\big(2^{8s+4j+6} - 31\cdot2^{7s+5j+5} +93\cdot2^{6s+6j+5} \big) 
     &  \text{if   }\quad  i=2s +3+j, \; k\geq 2s+4+j, \\
     &  0 \leq j\leq s-2, \\
    2^{3k+3s-1} -7\cdot2^{2k+6s-2} +7\cdot2^{k+9s-2} - 2^{12s-1}  &  \text{if  } i=3s+2, \; k\geq 3s +2 
     \end{cases}
\end{equation}  
  \begin{equation}
\label{eq 12.21}
\Gamma_{i}^{\left[s\atop{ s+1\atop s+1} \right]\times i} 
 = \begin{cases}
  2^{3s+3i-1} -7\cdot2^{4i-6} +3\cdot2^{3i-5}  &  \text{if  }\; 1\leq i\leq s+1, \\
   2^{6s+3j-1} -7\cdot2^{4s+4j-6} +3\cdot2^{3s+3j-5} \\
    + 7\cdot2^{2s+5j-10}- 5\cdot2^{2s+4j-8}&  \text{if  }\;i = s+j, \;2\leq j \leq s+3, \\
   2^{9s+2} +7\cdot2^{7s-5}-7\cdot2^{8s-2}+7\cdot2^{6s-4} &  \text{if  }\;i = 2s+1,\\
    2^{9s+5} +7\cdot2^{7s}-7\cdot2^{8s+2}+2^{6s} &  \text{if  }\;i = 2s+2,\\
      2^{9s+8} +7\cdot2^{7s+5}-7\cdot2^{8s+6}-2^{6s+5} &  \text{if  }\;i = 2s+3,\\
     2^{9s+3j+8} -7\cdot2^{8s+4j+6} +7\cdot2^{7s+5j+5} - 2^{6s+6j+5}  &  \text{if  }\;i = 2s+3+j, \;0 \leq j \leq s-1
  \end{cases}
\end{equation} 
We have the following reduction formulas\\
 
\begin{align}
\label{eq 12.22 }
\Gamma_{2s+2+j}^{\left[s\atop{ s+1\atop s+1} \right]\times k} &
 = 16^{2j}\Gamma_{2s+2-j}^{\left[s\atop{ s+(1-j)\atop s+(1-j)} \right]\times (k-2j)} & \text{if  }  0\leq j\leq 1 \\
\Gamma_{2s+3+j}^{\left[s\atop{ s+1\atop s+1} \right]\times k} &
 = 16^{2+3j}\Gamma_{2(s-j)+1}^{\left[s -j\atop{ s -j\atop s -j} \right]\times (k-2-3j)} & \text{if  }  0\leq j\leq s-1
 \label{eq 12.23}
\end{align}
\begin{example}
We have for $s=3:$\\
 \begin{equation*}
\Gamma_{i}^{\left[3\atop{ 3+1\atop 3+1} \right]\times k} 
 = \begin{cases}
 1  & \text{if  } i = 0,\; k\geq 1, \\
 21  &  \text{if  }  i=1,\;  k\geq 2,\\
 378  &  \text{if  }  i=2, \; k\geq 3,\\ 
 2^{k+2} + 6320 &    \text{if }  i = 3, \; k\geq 4, \\
  33\cdot 2^{k+2} + 100416 &    \text{if }  i = 4, \; k\geq 5, \\ 
  630\cdot2^{k+2} + 1524096 &    \text{if }  i = 5,\;  k\geq 6, \\
   1365\cdot2^{k+5} + 21224448  &    \text{if }  i = 6, \; k\geq 7, \\
   96\cdot2^{2k} +163008\cdot2^{k+2}+1029\cdot2^{18} &    \text{if }  i = 7,\; k\geq 8,  \\
   1696\cdot2^{2k} +2176512\cdot2^{k+2}+5723\cdot2^{18} &    \text{if }  i = 8, \; k\geq 9, \\
   105\cdot2^{2k+8} +2625\cdot2^{k+15}-90720\cdot2^{18} &    \text{if }  i = 9,\; k\geq 10,  \\
    105\cdot2^{2k+12} -315\cdot2^{k+20}+215040\cdot2^{18} &    \text{if }  i = 10,\;  k\geq 11, \\
    2^{3k+8} -7\cdot2^{2k+16} +7\cdot2^{k+25} - 2^{35}  &  \text{if  } i=11. \; k\geq 11.
     \end{cases}
\end{equation*}

\end{example}
\end{thm} 
\begin{proof}
  To prove Theorem \ref{thm 12.9} we proceed as in Theorem \ref{thm 10.9}.

 \end{proof}
 
   \subsection{\textbf{Computation of $ \Gamma_{i}^{\left[s\atop{ s +m\atop s+m} \right]\times k} 
    \quad for\; 2s+m+1 \leq i \leq 3s+2m,\;k\geq i,\; m\geq 2 $}}
  \label{subsec 7}
      \begin{lem}
      \label{lem 12.10}
   We have  for $m\geq 2,\; 0\leq j\leq m-1$ :
     \begin{align}
     \label{eq 12.24}
& \Gamma_{2s+m+1+j}^{\left[s\atop{ s +m\atop s+m} \right]\times k}\nonumber \\
& =  \begin{cases}
 21\cdot2^{2k+2s+m-2}+21\cdot2^{k+4s+3m-2}+735\cdot2^{k+5s+2m-3} - 945\cdot2^{k+4s+2m-3}\\
 +105\cdot2^{8s+4m-2} -105\cdot2^{6s+4m-2}-3255\cdot2^{7s+3m-3}+3255\cdot2^{6s+3m-3}  &\text{if  } j=0,\\
 & k\geq 2s+m+2\\
   21\cdot2^{2k+2s+m+3j-2}+21\cdot2^{k+4s+3m+3j-2}+735\cdot2^{k+5s+2m+4j-3} - 945\cdot2^{k+4s+2m+4j-3}\\
 +105\cdot2^{8s+4m+4j-2} -105\cdot2^{6s+4m+4j-2}-3255\cdot2^{7s+3m+5j-3}+3255\cdot2^{6s+3m+5j-3}
   &\text{if  } 0\leq j\leq m-2, \\
 &  k\geq 2s+m+2+j \\
   53\cdot2^{2s-1}\cdot(2^{2k+4m-4}-2^{4s+8m-2}) 
+ \big(735\cdot2^{5s-1}-1629\cdot2^{4s-1}\big)\cdot\big(2^{k+6m-6}-2^{2s+8m-5}\big)  \\
+   21\cdot[5\cdot2^{8s+8m-6} +3\cdot2^{6s+8m-8}-15\cdot2^{7s+8m-8} ]  &    \text{if }  j = m-1,\\
&  k\geq 2s+2m+1
  \end{cases}
  \end{align}

     \begin{align}
     \label{eq 12.25}
& \Gamma_{2s+m+1+j}^{\left[s\atop{ s +m\atop s+m} \right]\times (2s+m+1+j)} \nonumber \\
& =  \begin{cases}
 2^{9s+6m}-7\cdot2^{8s+4m-2}+7\cdot2^{7s+3m-3}-3\cdot2^{6s+3m-3}+2^{6s+4m-2} &\text{if  } j=0\\ 
  2^{9s+6m+2j}-7\cdot2^{8s+4m+4j-2}+7\cdot2^{7s+3m+5j-3}-3\cdot2^{6s+3m+5j-3}+2^{6s+4m+4j-2} &\text{if  } 0\leq j\leq m-2\\
    2^{9s+8m-3}-7\cdot2^{8s+8m-6}+7\cdot2^{7s+8m-8}-2^{6s+8m-8} &\text{if  } j=m-1\\
    & \\
     \end{cases}
  \end{align}
     We have  for $m\geq 2,\; 0\leq j\leq s-1$ :
     \begin{align}
     \label{eq 12.26}
& \Gamma_{2s+2m+1+j}^{\left[s\atop{ s +m\atop s+m} \right]\times k}\nonumber \\
& =  \begin{cases}
 105\cdot\big(2^{2k+2s+4m+4j-2}+ 7\cdot2^{k+5s+6m+4j-3} - 31\cdot2^{k+4s+6m+5j-3} \big)  \\
    + 105\cdot\big(2^{8s+8m+4j-2} - 31\cdot2^{7s+8m+5j-3} +93\cdot2^{6s+8m+6j-3} \big)  
 &\text{if  } 0\leq j\leq s-2,\; k\geq 2s+2m+2+j\\
  2^{3k+2m+3s-3}-7\cdot2^{2k+4m+6s-6}+7\cdot2^{k+6m+9s-8}-2^{8m+12s-9} &\text{if  } j=s-1,\; k\geq 3s+2m
   \end{cases}
  \end{align}

     \begin{align}
     \label{eq 12.27}
& \Gamma_{2s+2m+1+j}^{\left[s\atop{ s +m\atop s+m} \right]\times (2s+2m+1+j)} \nonumber \\
& =  \begin{cases}
 2^{9s+6m}-7\cdot2^{8s+8m-2}+7\cdot2^{7s+8m-3}-2^{6s+8m-3} &\text{if  } j=0\\ 
  2^{9s+8m+3j}-7\cdot2^{8s+8m+4j-2}+7\cdot2^{7s+8m+5j-3}-2^{6s+8m+6j-3} &\text{if  } 0\leq j\leq s-2 \\
    21\cdot2^{8m+12s-9} &\text{if  } j=s-1
   \end{cases}
  \end{align}
  We have the following reduction formulas\\
 
\begin{align}
\label{eq 12.28}
\Gamma_{2s+m+1+j}^{\left[s\atop{ s+m\atop s+m} \right]\times k} &
 = 16^{2j}\Gamma_{2s+1+(m-j)}^{\left[s\atop{ s+(m-j)\atop s+(m-j)} \right]\times (k-2j)} & \text{if  }  0\leq j\leq m-1 \\
\Gamma_{2s+2m+1+j}^{\left[s\atop{ s+m\atop s+m} \right]\times k} &
 = 16^{2m+3j}\Gamma_{2(s-j)+1}^{\left[s -j\atop{ s -j\atop s -j} \right]\times (k-2m-3j)} & \text{if  }  0\leq j\leq s-1
 \label{eq 12.29}
\end{align}
    \end{lem}
\begin{proof}    
    
 \end{proof}   
  \begin{example}
$ s = m =2 $ \\
     \begin{align*}
& \Gamma_{i}^{\left[2\atop{ 2+2\atop 2+2} \right]\times k} \\
& =  \begin{cases}
     1       &\text{if  } i =0,k\geq 1\\
   21      &\text{if  } i =1,k\geq 2\\
  2^{k+1} +362  &\text{if  } i =2,k\geq 3\\
   9\cdot2^{k+1} +6048 &\text{if  } i =3,k\geq 4\\
   348\cdot2^{k} +92640  &\text{if  } i =4,k\geq 5\\
     5712\cdot2^{k} +1295616  &\text{if  } i =5,k\geq 6\\ 
  3\cdot2^{2k+4} +1161\cdot2^{k+6} +15808512 &\text{if  } i =6,k\geq 7\\ 
   21\cdot2^{2k+4} +2163\cdot2^{k+9} +92897280  &\text{if  } i =7,k\geq 8\\ 
  424\cdot2^{2k+4} +167808\cdot2^{k+6} -1485832192  &\text{if  } i =8,k\geq 9\\ 
   105\cdot2^{2k+10} -315\cdot2^{k+17} +3523215360  &\text{if  } i =9,k\geq 10\\ 
     2^{3k+7}-7\cdot2^{2k+14}+7\cdot2^{k+22}-2^{31}&\text{if  } i =10,k\geq 10
      \end{cases}
  \end{align*}
 \end{example}
 \newpage
  \subsection{\textbf{An explicit formula for $ \Gamma_{i}^{\left[s\atop{ s+m\atop s+m} \right]\times k} , \quad for\; 1 \leq i \leq 3s+2m,\; k\geq i $}}
 \label{subsec 8} 
  \begin{thm}
 \label{thm 12.11}
 We have \\
  \begin{equation}
   \label{eq 12.30} 
   \Gamma_{i}^{\left[s\atop{ s +m\atop s+m} \right]\times k}  \\
  =  \begin{cases}
 1  & \text{if  } i = 0,\;k\geq 1 \\
 105\cdot2^{4i-6} - 21\cdot2^{3i -5}  &  \text{if  }  1\leq i\leq s-1,\;k\geq i+1, \\ 
2^{k+s-1} -2^{2s}+21\cdot\big(5\cdot2^{4s-6}-2^{3s-5}\big)  &\text{if  }  i = s,\;k\geq s+1\\
& \\
(21\cdot2^{j-1}-3)\cdot2^{k+s+2j-4}\\
 +21\cdot\big(5\cdot2^{4s+4j-6}-2^{3s+3j-5}-5\cdot2^{2s+4j-6}+2^{2s+3j-5}\big) 
  &\text{if  } i=s+j,\;  1\leq j\leq m-1 ,\\
  & k\geq s+j+1\\
  & \\
   \big(21\cdot2^{s+3m-5}+45\cdot2^{s+2m-4}\big)\cdot\big(2^{k}-2^{s+m+1}\big)\\
 + 105\cdot2^{4s+4m-6}-21\cdot2^{3s+3m-5}-21\cdot2^{2s+4m-6}+9\cdot2^{2s+3m-5} 
   &\text{if  }  i=s+m ,\;k\geq s+m+1\\
   & \\
    21\cdot\big[ 2^{k+s+3m+3j-5}+35\cdot2^{k+s+2m+4j-7}-9\cdot2^{k+s+2m+3j-6}\\
       + 5\cdot2^{4s+4m+4j-6} -2^{3s+3m+3j-5}- 5\cdot2^{2s+4m+4j-6}\\
       -155\cdot2^{2s+3m+5j-8}+45\cdot2^{2s+3m+4j-7}\big ]
   &\text{if  } i=s+m+j,\;1\leq j\leq s-1,\\
   & k\geq s+m+j+1\\
   & \\
   3\cdot2^{2k+2s+m-2}+21\cdot2^{k+4s+3m-5}+735\cdot2^{k+5s+2m-7} - 477\cdot2^{k+4s+2m-6}\\
 +105\cdot2^{8s+4m-6} -105\cdot2^{6s+4m-6}-3255\cdot2^{7s+3m-8}+1629\cdot2^{6s+3m-7} 
  &\text{if  } i=2s+m ,\;k\geq 2s+m+1\\ 
  & \\
    21\cdot2^{2k+2s+m+3j-2}+21\cdot2^{k+4s+3m+3j-2}+735\cdot2^{k+5s+2m+4j-3}\\
      - 945\cdot2^{k+4s+2m+4j-3} +105\cdot2^{8s+4m+4j-2} -105\cdot2^{6s+4m+4j-2}\\
      -3255\cdot2^{7s+3m+5j-3}+3255\cdot2^{6s+3m+5j-3}
   &\text{if  }i=2s+m+1+j,\\
   & 0\leq j\leq m-2, \;  k\geq 2s+m+2+j \\
   & \\
    53\cdot2^{2s-1}\cdot(2^{2k+4m-4}-2^{4s+8m-2}) \\
 + \big(735\cdot2^{5s-1}-1629\cdot2^{4s-1}\big)\cdot\big(2^{k+6m-6}-2^{2s+8m-5}\big)  \\
+   21\cdot[5\cdot2^{8s+8m-6} +3\cdot2^{6s+8m-8}-15\cdot2^{7s+8m-8} ]  &  \text{if } i=2s+2m,\;
  k\geq 2s+2m+1 \\
 & \\
  105\cdot\big(2^{2k+2s+4m+4j-2}+ 7\cdot2^{k+5s+6m+4j-3} - 31\cdot2^{k+4s+6m+5j-3} \big)  \\
    + 105\cdot\big(2^{8s+8m+4j-2} - 31\cdot2^{7s+8m+5j-3} +93\cdot2^{6s+8m+6j-3} \big)  
 &\text{if  }i=2s+2m+1+j, \\
 &  0\leq j\leq s-2,\; k\geq 2s+2m+2+j\\
 & \\
  2^{3k+2m+3s-3}-7\cdot2^{2k+4m+6s-6}+7\cdot2^{k+6m+9s-8}-2^{8m+12s-9} &\text{if  } i=3s+2m,\; k\geq 3s+2m
 \end{cases}
\end{equation}
 
 \begin{equation}
 \label{eq 12.31}
   \Gamma_{i}^{\left[s\atop{ s +m\atop s+m} \right]\times i} \\
  =  \begin{cases}
 2^{3s+2m+3i-3}-7\cdot2^{4i-6}+3\cdot2^{3i-5}  &\text{if  }  1\leq i\leq s+1\\
 & \\
  2^{6s+2m+3j-3}-7\cdot2^{4s+4j-6}+3\cdot2^{3s+3j-5}\\
  +(2^{j-1}-1)\cdot2^{2s+3j-5} &\text{if  }i=s+j,\; 1\leq j\leq m+1  \\
  & \\
   2^{6s+5m+3j-3}-7\cdot2^{4s+4m+4j-6}+3\cdot2^{3s+3m+3j-5}\\
+ 2^{2s+4m+4j-6}+7\cdot2^{2s+3m+5j-8}- 9\cdot2^{2s+3m+4j-7}&\text{if  } i=s+m+j,\;1\leq j\leq s-1\\
& \\
 2^{9s+5m-3}-7\cdot2^{8s+4m-6}+7\cdot2^{7s+3m-8}+3\cdot2^{6s+3m-7}+2^{6s+4m-6} &\text{if  } i=2s+m\\
  & \\
 2^{9s+5m+2j}-7\cdot2^{8s+4m+4j-2}+7\cdot2^{7s+3m+5j-3}-3\cdot2^{6s+3m+5j-3}+2^{6s+4m+4j-2} &\text{if  }i=2s+m+1+j,\\
 &  0\leq j\leq m-2\\
 & \\
  2^{9s+8m-3}-7\cdot2^{8s+8m-6}+7\cdot2^{7s+8m-8}-2^{6s+8m-8} &\text{if  } i=2s+2m\\ 
  & \\
 2^{9s+8m+3j}-7\cdot2^{8s+8m+4j-2}+7\cdot2^{7s+8m+5j-3}-2^{6s+8m+6j-3} &\text{if  }i=2s+2m+1+j,\\
 &  0\leq j\leq s-2 \\
 & \\
  21\cdot2^{8m+12s-9} &\text{if  } i=3s+2m
 \end{cases}
\end{equation}
  We have for $ 0\leq j\leq m-2,\;k \geq 2s+m+2+j$\\[0.1 cm]
  
  \begin{align}
 \Gamma_{2s+m+1+j}^{\left[s\atop{ s+m\atop s+m} \right]\times k}
 = 16^{2j}\Gamma_{2s+1+(m-j)}^{\left[s\atop{ s+(m-j)\atop s+(m-j)} \right]\times (k-2j)} \label{eq 12.32}\\
 =  2^{8j}\cdot\big[  21\cdot2^{2(k-2j)+2s+(m-j)-2}+21\cdot2^{(k-2j)+4s+3(m-j)-2}+735\cdot2^{k-2j+5s+2(m-j)-3}\nonumber\\
      - 945\cdot2^{k-2j+4s+2(m-j)-3} +105\cdot2^{8s+4(m-j)-2} -105\cdot2^{6s+4(m-j)-2}\nonumber\\
      -3255\cdot2^{7s+3(m-j)-3}+3255\cdot2^{6s+3(m-j)-3}\big]\nonumber\\
      =    21\cdot2^{2k+2s+m+3j-2}+21\cdot2^{k+4s+3m+3j-2}+735\cdot2^{k+5s+2m+4j-3}\nonumber\\
      - 945\cdot2^{k+4s+2m+4j-3} +105\cdot2^{8s+4m+4j-2} -105\cdot2^{6s+4m+4j-2}\nonumber\\
      -3255\cdot2^{7s+3m+5j-3}+3255\cdot2^{6s+3m+5j-3}\nonumber\\
       & \nonumber 
       \end{align}
   We have for $ j=m-1,\;k \geq 2s+2m+1$\\[0.1 cm]
   
    \begin{align}
 \Gamma_{2s+2m}^{\left[s\atop{ s+m\atop s+m} \right]\times k}
 = 16^{2m-2}\Gamma_{2s+2}^{\left[s\atop{ s+1\atop s+1} \right]\times (k-2(m-1))} \label{eq 12.33}\\
= 2^{8m-8}\cdot\bigg[
  53\cdot2^{2s-1}\cdot(2^{2(k-2m+2)}-2^{4s+6}) 
+ \big(735\cdot2^{5s-1}-1629\cdot2^{4s-1}\big)\cdot\big(2^{k-2m+2}-2^{2s+3}\big) \nonumber \\
+   21\cdot[5\cdot2^{8s+2} +3\cdot2^{6s}-15\cdot2^{7s} ] \bigg] \nonumber \\
=  53\cdot2^{2s-1}\cdot(2^{2k+4m-4}-2^{4s+8m-2}) 
 + \big(735\cdot2^{5s-1}-1629\cdot2^{4s-1}\big)\cdot\big(2^{k+6m-6}-2^{2s+8m-5}\big)\nonumber  \\
+   21\cdot[5\cdot2^{8s+8m-6} +3\cdot2^{6s+8m-8}-15\cdot2^{7s+8m-8} ]  \nonumber \\
& \nonumber
 \end{align}

 We have for $ 0\leq j\leq s-2,\;k \geq 2s+2m+2+j$
  \begin{align}
  \Gamma_{2s+2m+1+j}^{\left[s\atop{ s+m\atop s+m} \right]\times k} 
 = 16^{2m+3j}\Gamma_{2(s-j)+1}^{\left[s -j\atop{ s -j\atop s -j} \right]\times (k-2m-3j)} \label{eq 12.34} \\
 =  2^{8m+12j}\cdot105\cdot\big(2^{2(k-2m-3j)+2(s-j) -2}+ 7\cdot2^{k-2m-3j+5(s-j)-3} - 31\cdot2^{k-2m-3j+4(s-j)-3} \nonumber \\
     + 2^{8(s-j)-2} - 31\cdot2^{7(s-j)-3} +93\cdot2^{6(s-j)-3} \big)\nonumber \\
    =  105\cdot\big(2^{2k+2s+4m+4j-2}+ 7\cdot2^{k+5s+6m+4j-3} - 31\cdot2^{k+4s+6m+5j-3} \big) \nonumber \\
    + 105\cdot\big(2^{8s+8m+4j-2} - 31\cdot2^{7s+8m+5j-3} +93\cdot2^{6s+8m+6j-3} \big)  \nonumber \\
    & \nonumber
  \end{align}
We have for $j=s-1,\;k\geq 3s+2m$\\[0.1 cm]

   \begin{align}
  \Gamma_{3s+2m}^{\left[s\atop{ s+m\atop s+m} \right]\times k} 
 = 16^{2m+3s-3}\Gamma_{3}^{\left[1\atop{ 1\atop 1} \right]\times (k-2m-3s+3)}  \label{eq 12.35}\\
= 2^{8m+12s-12}\cdot\big[2^{3(k-2m-3s+3)}-7\cdot2^{2(k-2m-3s+3)}+7\cdot2^{k-2m-3s+3+1} -2^{3}\big]\nonumber\\
= 2^{3k+2m+3s-3}-7\cdot2^{2k+4m+6s-6}+7\cdot2^{k+6m+9s-8}-2^{8m+12s-9}\nonumber\\
 & \nonumber
  \end{align}

 \begin{align}
 & \text{ Let q be a rational integer $ \geq 1,$ then }\label{eq 12.36} \\
  &  g_{k,s,m}(t,\eta,\xi ) =  g(t,\eta,\xi )  = \sum_{deg Y\leq k-1}\sum_{deg Z \leq  s-1}E(tYZ)\sum_{deg U \leq s+m-1}E(\eta YU)\sum_{deg V \leq s+m-1}E(\eta YV) \nonumber \\
 & =  2^{3s+2m+k- r(D^{\left[s\atop{ s+m\atop s+m} \right]\times k}(t,\eta,\xi)) }\nonumber \\
 &\int_{\mathbb{P}^{3}} g^{q}(t,\eta,\xi  )d t d \eta d\xi = \nonumber\\
& = \sum_{(t,\eta,\xi  )\in \mathbb{P}/\mathbb{P}_{k+s-1}\times \mathbb{P}/\mathbb{P}_{k+s+m-1}\times \mathbb{P}/\mathbb{P}_{k+s+m-1}\atop { r( D^{\left[\stackrel{s}{\stackrel{s+m}{s+m}}\right] \times k }(t,\eta ,\xi )  =  i}}
 2^{\big(k+3s+2m-  r( D^{\left[\stackrel{s}{\stackrel{s+m}{s+m}}\right] \times k }(t,\eta ,\xi ) )\big)q  } \int_{\mathbb{P}_{k+s-1}}dt \int_{\mathbb{P}_{k+s+m-1}}d\eta \int_{\mathbb{P}_{k+s+m-1}}d\xi \nonumber \\
 & = \sum_{i = 0}^{\inf(k,3s+2m)}\sum_{(t,\eta,\xi  )\in \mathbb{P}/\mathbb{P}_{k+s-1}\times \mathbb{P}/\mathbb{P}_{k+s+m-1}\times \mathbb{P}/\mathbb{P}_{k+s+m-1}}
 2^{(k+3s+2m-i)q  } \int_{\mathbb{P}_{k+s-1}}dt \int_{\mathbb{P}_{k+s+m-1}}d\eta \int_{\mathbb{P}_{k+s+m-1}}d\xi \nonumber \\
  & = 2^{(k+3s+2m)q -(3k +3s+2m-3)}
 \sum_{i=0}^{\inf(k,3s+2m)}\Gamma_{i}^{\left[s\atop{ s+m \atop s+m }\right]\times k}\cdot2^{-iq} \nonumber
 & \nonumber \\
 & \nonumber
 \end{align}  

 We denote by  $ R_{q}(k,s,m) $ the number of solutions \\
 $(Y_1,Z_1,U_{1},V_{1}, \ldots,Y_q,Z_q,U_{q},V_{q}) $  of the polynomial equations
   \[\left\{\begin{array}{c}
 Y_{1}Z_{1} +Y_{2}Z_{2}+ \ldots + Y_{q}Z_{q} = 0,  \\
   Y_{1}U_{1} + Y_{2}U_{2} + \ldots  + Y_{q}U_{q} = 0,\\
    Y_{1}V_{1} + Y_{2}V_{2} + \ldots  + Y_{q}V_{q} = 0,\\  
 \end{array}\right.\]
  satisfying the degree conditions \\
                   $$  degY_i \leq k-1 , \quad degZ_i \leq s-1 ,\quad degU_{i}\leq s+m-1,\quad degV_{i}\leq s+m-1 \quad for \quad 1\leq i \leq q. $$ \\                           
Then \\
\begin{align}
  R_{q}(k,s,m) & =  \int_{\mathbb{P}\times \mathbb{P}\times \mathbb{P}} g_{k,s,m}^{q}(t,\eta,\xi  )dtd\eta d\xi  
  = 2^{(k+3s+2m)q -(3k +3s+2m-3)}
 \sum_{i=0}^{\inf(k,3s+2m)}\Gamma_{i}^{\left[s\atop{ s+m \atop s+m }\right]\times k}\cdot2^{-iq} \label{eq 12.37}
\end{align}

  \begin{example}
We have for $s=3,\;m=4,\;k=10:$\\
 \begin{equation*}
\Gamma_{i}^{\left[3\atop{ 3+4\atop 3+4} \right]\times 10} 
 = \begin{cases}
 1  & \text{if  } i = 0 \\
 21  &  \text{if  }  i=1\\
 378  &  \text{if  }  i=2\\ 
 10416 &    \text{if }  i = 3\\
  140352 &    \text{if }  i = 4,  \\ 
  1994112 &    \text{if }  i = 5, \\
   29598720 &    \text{if }  i = 6 \\
   458661888 &    \text{if }  i = 7, \\
   109389\cdot2^{16} &    \text{if }  i = 8, \\
   213759\cdot2^{19} &    \text{if }  i = 9, \\
  2^{44} -14273\cdot2^{23}  &  \text{if  } i=10
     \end{cases}
\end{equation*}

\end{example}

\begin{example} s = 3, m=4, k = 7, q = 3\\

 The number  $\Gamma_{i}^{\left[3\atop{ 3+4\atop 3+4} \right]\times 7}$ of rank i matrices of the form \\
  $$  \left ( \begin{array} {ccccccc}
\alpha _{1} & \alpha _{2} & \alpha _{3} & \alpha _{4} & \alpha _{5} & \alpha _{6} & \alpha _{7} \\
\alpha _{2 } & \alpha _{3} & \alpha _{4} & \alpha _{5} & \alpha _{6} & \alpha _{7}& \alpha _{8} \\
\alpha _{3 } & \alpha _{4} & \alpha _{5} & \alpha _{6} & \alpha _{7} & \alpha _{8}& \alpha _{9} \\
\beta  _{1} & \beta  _{2} & \beta  _{3}  & \beta  _{4} & \beta  _{5} & \beta  _{6} & \beta  _{7}\\
\beta  _{2} & \beta  _{3} & \beta  _{4}  & \beta  _{5} & \beta  _{6} & \beta  _{7} & \beta  _{8}\\
\beta  _{3} & \beta  _{4} & \beta  _{5}  & \beta  _{6} & \beta  _{7} & \beta  _{8} & \beta  _{9}\\
\beta  _{4} & \beta  _{5} & \beta  _{6}  & \beta  _{7} & \beta  _{8} & \beta  _{9} & \beta  _{10}\\
\beta  _{5} & \beta  _{6} & \beta  _{7}  & \beta  _{8} & \beta  _{9} & \beta  _{10} & \beta  _{11}\\
\beta  _{6} & \beta  _{7} & \beta  _{8}  & \beta  _{9} & \beta  _{10} & \beta  _{11} & \beta  _{12}\\
\beta  _{7} & \beta  _{8} & \beta  _{9}  & \beta  _{10} & \beta  _{11} & \beta  _{12} & \beta  _{13}\\
\gamma  _{1} & \gamma  _{2} & \gamma  _{3}  & \gamma  _{4} & \gamma  _{5} & \gamma  _{6}  & \gamma  _{7}\\
\gamma  _{2} & \gamma  _{3} & \gamma  _{4}  & \gamma  _{5} & \gamma  _{6} & \gamma  _{7} & \gamma  _{8}\\
\gamma  _{3} & \gamma  _{4} & \gamma  _{5}  & \gamma  _{6} & \gamma  _{7} & \gamma  _{8}  & \gamma  _{9}\\
\gamma  _{4} & \gamma  _{5} & \gamma  _{6}  & \gamma  _{7} & \gamma  _{8} & \gamma  _{9} & \gamma  _{10}\\
\gamma  _{5} & \gamma  _{6} & \gamma  _{7}  & \gamma  _{8} & \gamma  _{9} & \gamma  _{10}  & \gamma  _{11}\\
\gamma  _{6} & \gamma  _{7} & \gamma  _{8}  & \gamma  _{9} & \gamma  _{10} & \gamma  _{11} & \gamma  _{12}\\
\gamma  _{7} & \gamma  _{8} & \gamma  _{9}  & \gamma  _{10} & \gamma  _{11} & \gamma  _{12} & \gamma  _{13}
 \end{array}  \right) $$
is equal to

\[ \begin{cases}
 1  & \text{if  } i = 0 \\
 21  &  \text{if  }  i=1\\
 378  &  \text{if  }  i=2\\ 
 6832 &    \text{if }  i = 3\\
  108096 &    \text{if }  i = 4,  \\ 
  1714560 &    \text{if }  i = 5, \\
   27276288 &    \text{if }  i = 6 \\
  2^{35} -3553\cdot2^{13} &    \text{if }  i = 7, \\
      \end{cases}\]

 The number of solutions \\
 $(Y_1,Z_1,U_{1},V_{1},Y_2,Z_2,U_{2},V_{2} ,Y_3,Z_3,U_{3},V_{3}) $  of the polynomial equations
   \[\left\{\begin{array}{c}
 Y_{1}Z_{1} +Y_{2}Z_{2} + Y_{3}Z_{3} = 0,  \\
   Y_{1}U_{1} + Y_{2}U_{2}  + Y_{3}U_{3} = 0,\\
    Y_{1}V_{1} + Y_{2}V_{2} + Y_{3}V_{3} = 0,\\  
 \end{array}\right.\]
  satisfying the degree conditions \\
                   $$  degY_i \leq 6 , \quad degZ_i \leq 2 ,\quad degU_{i}\leq 6,\quad degV_{i}\leq 6 \quad for \quad 1\leq i \leq 3. $$ \\                           

is equal to 
\begin{align*}
&  R_{3}(7,3,4)  =  \int_{\mathbb{P}\times \mathbb{P}\times \mathbb{P}} g_{7,3,4}^{3}(t,\eta,\xi  )dtd\eta d\xi  
 = 2^{37}\cdot \sum_{i = 0}^{7} \Gamma_{i}^{\left[3\atop{ 3+4\atop 3+4} \right]\times 7}  \cdot2^{- 3i}\\
 & = 2^{37}\cdot\big(1 + 21\cdot2^{-3} + 378\cdot2^{-6} + 6832\cdot2^{-9} + 108096 \cdot 2^{-12} + 1714560 \cdot2^{-15}\\
& +  27276288 \cdot 2^{-18}  +  (2^{35} -3553\cdot2^{13}) \cdot 2^{ -21}  \big )  = 4243395\cdot2^{29}
\end{align*}
\end{example}

 \begin{example} The fraction of square triple persymmetric   $ \left[s\atop{ s+m\atop s+m} \right]\times (3s+2m) $   matrices which are  invertible is equal to 
   $ \displaystyle \frac{\Gamma_{3s+2m}^{\left[s\atop{ s+m\atop s+m} \right]\times (3s+2m)}}
 {\sum_{i = 0}^{3s+2m}\Gamma_{i}^{\left[s\atop{ s+m\atop s+m} \right]\times (3s+2m)} } = \frac{21}{64}. $
   \end{example}

 \end{thm}

 \end{document}